\documentclass[a4paper,twoside,12pt]{memoir}

\title{Multitriangulations, pseudotriangulations and some problems of realization of polytopes}
\author{Vincent Pilaud}

\usepackage[latin1]{inputenc} 
\usepackage[T1]{fontenc} 
\usepackage{xcolor} 
\usepackage[frenchb,spanish,english]{babel} 
\usepackage{aeguill} 
\usepackage{xspace} 
\usepackage{microtype} 
\usepackage{amsfonts,amsmath,amssymb,amsthm} 
\usepackage{stmaryrd} 
\usepackage{wasysym} 
\usepackage{bbm} 

\usepackage{hyperref} 
\hypersetup{pdftex, bookmarksnumbered, a4paper, colorlinks=true, linkcolor=black, citecolor=darkblue, urlcolor=blue} 
\usepackage{memhfixc} 
\hypersetup{pdftitle={\thetitle}, pdfauthor={\theauthor}} 
\usepackage{hypcap} 
\usepackage{url} 

\usepackage{graphicx} 
\graphicspath{{introduction/figures/}{introduction_multitriangulations/figures/}{stars/figures/}{mpt/figures/}{multiassociahedron/figures/}{appendix/figures/}{introduction_polytopes/figures/}{nonpolytopal/figures/}{psn/figures/}{miscellaneous/figures/}} 
\usepackage{rotating} 
\usepackage{fancybox} 
\usepackage{algorithmic} 
\usepackage{verbatim} 
\usepackage{import} 
\usepackage{paralist} 
\tightlists 
\usepackage{a4wide} 
\usepackage{multicol} 

\definecolor{red}{rgb}{1,0,0}
\definecolor{blue}{rgb}{0,0,1}
\definecolor{darkblue}{rgb}{0,0,0.7} 
\definecolor{gris}{gray}{0.4} 
\definecolor{green}{rgb}{0,1,0}
\definecolor{gray25}{gray}{0.75}
\definecolor{gray50}{gray}{0.5}
\definecolor{gray75}{gray}{0.25}
\definecolor{black}{gray}{0}

\newcommand{\darkblue}{\color{darkblue}}

\newcommand{\black}{\color{black}}

\newcommand{\ligne}{{\black\hrule}} 
\newcommand{\demiligne}{{\black\hrule width 5cm}} 
\newcommand{\carre}{\;\;{\darkblue\rule{1ex}{1ex}}\;\;} 
\newenvironment{lignes}{\vspace{.4cm}\noindent\begin{minipage}{\textwidth}{\black\hrule}\vspace*{.3cm}}{\vspace*{-.2cm}{\black\hrule}\end{minipage}\vspace{.4cm}} 

\makeatletter 

\newcommand{\partfigure}{}
\newcommand{\namepart}{}
\renewcommand*{\toclevel@part}{0} 
	\def\@partapp{Part}
	\renewcommand*{\partnumfont}{\black\fontseries{m}\fontshape{sl}\fontsize{15}{15}\selectfont\centering} 
	\renewcommand*{\printpartnum}{\partnumfont \@partapp \space \numtoName{\c@part}} 

\makechapterstyle{chapitre}{%
	\openright 
	\setlength{\beforechapskip}{0pt} 
	\setlength{\midchapskip}{18pt} 
	\setlength{\afterchapskip}{30pt} 
	\renewcommand*{\chapnumfont}{\black\fontseries{m}\fontshape{sl}\fontsize{15}{15}\selectfont\flushright} 
	\renewcommand*{\chaptitlefont}{\darkblue\fontseries{m}\fontshape{sc}\fontsize{25}{25}\selectfont\centering} 
	\renewcommand*{\printchapternum}{\chapnumfont \@chapapp \space \numtoName{\c@chapter}} 
	\renewcommand*{\printchaptertitle}[1]{\chaptitlefont \ligne \begin{center} ##1 \end{center} \ligne%
										\ifnum \value{compteurlof}>\z@ \addcontentsline{lof}{chapter}{\ifnum \c@chapter>\z@ \thechapter~~~ \fi ##1} \fi%
										\setcounter{figure}{0}%
										\thispagestyle{baspage}} 
}
\chapterstyle{chapitre}

\maxsecnumdepth{subsubsection} 
\setsecheadstyle{\fontseries{m}\fontshape{sc}\fontsize{18}{18}\selectfont \noindent \textcolor{darkblue}} 
\setsubsecheadstyle{\fontseries{b}\fontshape{n}\fontsize{15}{15}\selectfont \noindent \textcolor{darkblue}} 
\setsubsubsecheadstyle{\fontseries{b}\fontshape{n}\fontsize{13}{13}\selectfont \noindent \textcolor{darkblue}} 
\setparaheadstyle{\fontseries{b}\fontshape{n}\fontsize{12}{12}\selectfont \noindent \textcolor{darkblue}} 

\settocdepth{subsection} 

\makepagestyle{hautpage} 
	\makeevenhead{hautpage}{\textbf \thepage}{\leftmark}{} 
	\makeoddhead{hautpage}{}{\rightmark}{\textbf \thepage} 
	\makeheadrule{hautpage}{\textwidth}{.5pt} 
\makepagestyle{hautpageintro} 
	\makeevenhead{hautpageintro}{\textbf \thepage}{\namepart}{} 
	\makeoddhead{hautpageintro}{}{\namepart}{\textbf \thepage} 
	\makeheadrule{hautpageintro}{\textwidth}{.5pt} 
\makepagestyle{baspage} 
	\makeoddfoot{baspage}{\namepart\carre\theauthor}{}{\textbf \thepage} 
	\makefootrule{baspage}{\textwidth}{.4pt}{1pt} 
\renewcommand*{\chaptermark}[1]{\markboth{\@chapapp\ \thechapter\ -\ #1}{\@chapapp\ \thechapter\ -\ #1}} 
\makeindex

\setverbatimfont{\fontencoding{T1}\fontfamily{phv}\fontseries{n}\fontsize{10.5}{13}\selectfont} 
\tabson[4] 
\wrappingon 
\bvbox 
\makeatother

\newtheorem{theorem}{Theorem}[chapter]
\newtheorem{proposition}[theorem]{Proposition}
\newtheorem{lemma}[theorem]{Lemma}
\newtheorem{corollary}[theorem]{Corollary}
\newtheorem{conjecture}[theorem]{Conjecture}
\newtheorem{definition}[theorem]{Definition}

\theoremstyle{definition} 
\newtheorem{example}[theorem]{Example}

\newtheorem{remark}[theorem]{Remark}

\newtheorem{observation}[theorem]{Observation}
\newtheorem{question}[theorem]{Question}

\newcommand{\defn}[1]{\emph{\darkblue #1}} 

\newcommand{\ie}{\textit{i.e.}~} 
\newcommand{\etc}{etc.} 
\newcommand{\eg}{\textit{e.g.}~} 
\newcommand{\viceversa}{\textit{vice~versa}} 
\newcommand{\apriori}{\textit{a~priori}\xspace} 
\newcommand{\vs}{\textit{vs.}~} 

\newcommand{\gon}[1]{\mbox{$#1$-gon}}
\newcommand{\kcross}[1]{\mbox{$#1$-crossing}}
\newcommand{\ktri}[1]{\mbox{$#1$-trian}\-gu\-la\-tion}
\newcommand{\krel}[1]{\mbox{$#1$-rele}\-vant}
\newcommand{\kbound}[1]{\mbox{$#1$-boun}\-dary}
\newcommand{\kirrel}[1]{\mbox{$#1$-irre}\-le\-vant}
\newcommand{\kear}[1]{\mbox{$#1$-ear}}
\newcommand{\kstar}[1]{\mbox{$#1$-star}}
\newcommand{\kcolorable}[1]{\mbox{$#1$-colorable}}
\newcommand{\kcoloring}[1]{\mbox{$#1$-coloring}}
\newcommand{\kpath}[1]{\mbox{$#1$-path}}
\newcommand{\kdiag}[1]{\mbox{$#1$-diagonal}}
\newcommand{\kzz}[1]{\mbox{$#1$-zigzag}}
\newcommand{\kaccordion}[1]{\mbox{$#1$-accor}\-dion}
\newcommand{\kvalid}[1]{\mbox{$#1$-valid}}
\newcommand{\kdepth}[1]{\mbox{$#1$-depth}}
\newcommand{\kcompatible}[1]{\mbox{$#1$-com}\-pa\-ti\-ble}
\newcommand{\kkernel}[1]{\mbox{$#1$-kernel}}
\newcommand{\kpointed}[1]{\mbox{$#1$-pointed}}
\newcommand{\kalter}[1]{\mbox{$#1$-alter}nation}
\newcommand{\pt}[1]{\mbox{$#1$-pseu}\-do\-trian\-gu\-la\-tion}

\newcommand{\equivelar}[1]{\mbox{$#1$-equi}\-velar}
\newcommand{\vfarther}[1]{\mbox{$#1$-farther}}
\newcommand{\vmax}[1]{\mbox{$#1$-maximal}}
\newcommand{\dimensional}[1]{\mbox{$#1$-dimen}\-sio\-nal}
\newcommand{\zoseq}{\mbox{$0/1$-sequence}}
\newcommand{\neighborly}[1]{\mbox{$#1$-neigh}\-bor\-ly}
\newcommand{\tuple}[1]{\mbox{$#1$-tuple}}
\newcommand{\spars}[1]{\mbox{$#1$-sparse}}
\newcommand{\tight}[1]{\mbox{$#1$-tight}}
\newcommand{\arbores}[1]{\mbox{$#1$-arbo}rescence}
\newcommand{\dfsvector}{\mbox{\dfs-vector}}
\newcommand{\dfslabeling}{\mbox{\dfs-labeling}}
\newcommand{\dfspolytope}{\mbox{\dfs-polytope}}
\newcommand{\mpt}{multi\-pseudo\-trian\-gu\-la\-tion\xspace}
\newcommand{\piantiperiodic}{\mbox{$\pi$-anti}\-periodic}
\newcommand{\increasing}[1]{\mbox{$#1$-increa}sing}
\newcommand{\decreasing}[1]{\mbox{$#1$-decrea}sing}
\newcommand{\greedy}[1]{\mbox{$#1$-greedy}}
\newcommand{\poly}[1]{\mbox{$#1$-poly}}
\newcommand{\face}[1]{\mbox{$#1$-face}}
\newcommand{\simp}[1]{\mbox{$#1$-simplex}}
\newcommand{\cycle}[1]{\mbox{$#1$-cycle}}
\newcommand{\clique}[1]{\mbox{$#1$-clique}}
\newcommand{\fvector}{\mbox{$f$-vector}}
\newcommand{\skeleton}[1]{\mbox{$#1$-skele}ton}
\newcommand{\connected}[1]{\mbox{$#1$-connec}ted}
\newcommand{\complex}[1]{\mbox{$#1$-complex}}
\newcommand{\regular}[1]{\mbox{$#1$-regular}}
\newcommand{\xpsn}[1]{\mbox{$#1$-\psn}}
\newcommand{\xppsn}[1]{\mbox{$#1$-\ppsn}}
\newcommand{\domino}[1]{\mbox{$#1$-domino}}

\newcommand{\R}{\ensuremath{\mathbb{R}}} 
\newcommand{\N}{\ensuremath{\mathbb{N}}} 
\newcommand{\Z}{\ensuremath{\mathbb{Z}}} 
\newcommand{\D}{\ensuremath{\mathbb{D}}} 
\newcommand{\PP}{\ensuremath{\mathbb{P}}} 
\newcommand{\LL}{\ensuremath{\mathbb{L}}} 
\newcommand{\UU}{\ensuremath{\mathbb{U}}} 
\newcommand{\GG}{\ensuremath{\mathbb{G}}} 
\newcommand{\cA}{\ensuremath{\mathcal{A}}} 
\newcommand{\cC}{\ensuremath{\mathcal{C}}} 
\newcommand{\cE}{\ensuremath{\mathcal{E}}} 
\newcommand{\cF}{\ensuremath{\mathcal{F}}} 
\newcommand{\cK}{\ensuremath{\mathcal{K}}} 
\newcommand{\cM}{\ensuremath{\mathcal{M}}} 
\newcommand{\cL}{\ensuremath{\mathcal{L}}} 
\newcommand{\cP}{\ensuremath{\mathcal{P}}} 
\newcommand{\cS}{\ensuremath{\mathcal{S}}} 
\newcommand{\cT}{\ensuremath{\mathcal{T}}} 
\newcommand{\cV}{\ensuremath{\mathcal{V}}} 
\newcommand{\cX}{\ensuremath{\mathcal{X}}} 
\newcommand{\cY}{\ensuremath{\mathcal{Y}}} 
\newcommand{\cZ}{\ensuremath{\mathcal{Z}}} 
\newcommand{\ra}{\ensuremath{\mathrm{a}}} 
\newcommand{\rb}{\ensuremath{\mathrm{b}}} 
\newcommand{\rc}{\ensuremath{\mathrm{c}}} 
\newcommand{\rd}{\ensuremath{\mathrm{d}}} 
\newcommand{\re}{\ensuremath{\mathrm{e}}} 
\newcommand{\rf}{\ensuremath{\mathrm{f}}} 
\newcommand{\rg}{\ensuremath{\mathrm{g}}} 
\newcommand{\rG}{\ensuremath{\mathrm{G}}} 
\newcommand{\rh}{\ensuremath{\mathrm{h}}} 
\newcommand{\rH}{\ensuremath{\mathrm{H}}} 
\newcommand{\rI}{\ensuremath{\mathrm{I}}} 
\newcommand{\rJ}{\ensuremath{\mathrm{J}}} 
\newcommand{\rK}{\ensuremath{\mathrm{K}}} 
\newcommand{\rL}{\ensuremath{\mathrm{L}}} 
\newcommand{\rx}{\ensuremath{\mathrm{x}}} 
\newcommand{\ry}{\ensuremath{\mathrm{y}}} 
\newcommand{\rz}{\ensuremath{\mathrm{z}}} 
\newcommand{\rt}{\ensuremath{\mathrm{t}}} 
\newcommand{\subf}{\ensuremath{\mathit{sf}}} 
\newcommand{\tu}{\ensuremath{\mathit{tu}}} 
\newcommand{\hl}{\ensuremath{\mathit{hl}}} 
\newcommand{\VS}{\ensuremath{\textsc{vs}}} %

\newcommand{\cl}{\ensuremath{\prec}} 
\newcommand{\cle}{\ensuremath{\preccurlyeq}} 
\newcommand{\tle}{\ensuremath{\trianglelefteqslant}} 
\newcommand{\ssm}{\ensuremath{\smallsetminus}} 
\newcommand{\diffsym}{\ensuremath{\triangle}} 

\newcommand{\ens}[2]{\ensuremath{\left\{#1\,\middle|\,#2\right\}}} 
\newcommand{\Floor}[1]{\ensuremath{\left\lfloor#1\right\rfloor}} 
\newcommand{\Fracfloor}[2]{\ensuremath{\left\lfloor\frac{#1}{#2}\right\rfloor}} 
\newcommand{\flattening}[2]{\ensuremath{#1\,\rotatebox{90}{\tiny $\ll$}\, #2}} 
\newcommand{\inflating}[2]{\ensuremath{#1\,\rotatebox{90}{\tiny $\gg$}\, #2}} 

\DeclareMathOperator{\conv}{conv} 
\DeclareMathOperator{\prism}{prism} 
\DeclareMathOperator{\cm}{\textsc{cm}} 
\DeclareMathOperator{\sign}{sign} 
\DeclareMathOperator{\curry}{curry} 
\DeclareMathOperator{\subarr}{sub} 
\DeclareMathOperator{\gr}{gr} 

\newcommand{\rrb}{\rrbracket} 
\newcommand{\llb}{\llbracket} 
\newcommand{\simplex}{\ensuremath{\triangle}} 
\newcommand{\polar}{\ensuremath{\diamond}} 
\newcommand{\dual}{\#} 
\newcommand{\one}{\ensuremath{\mathbbm{1}}} 
\newcommand{\dotprod}[2]{\langle#1 \,|\, #2\rangle} 
\newcommand{\kdeg}[3]{\ensuremath{\deg_{#1}(#2,#3)}} 
\newcommand{\sym}{\ensuremath{\mathfrak{S}}} 
\newcommand{\vartriangledown}{\ensuremath{\rotatebox[origin=c]{180}{$\triangle$}}} 
\newcommand{\sub}[1]{\ensuremath{\underline{#1}}} 
\newcommand{\bn}{{\overline{n}}}
\newcommand{\bm}{{\overline{m}}}
\newcommand{\bchi}{{\overline{\chi}}}
\newcommand{\defA}{A^{\sim}}
\newcommand{\defb}{b^{\sim}}
\newcommand{\defP}{P^{\sim}}
\newcommand{\KG}{\mathrm{KG}}
\newcommand{\eqdef}{\mbox{~\raisebox{0.2ex}{\scriptsize\ensuremath{\mathrm:}}\ensuremath{=} }}
\newcommand{\eqfed}{\mbox{~\ensuremath{=}\raisebox{0.15ex}{\scriptsize\ensuremath{\mathrm:}} }}

\newcommand{\svs}{\vspace{.2cm}} 
\newcommand{\mvs}{\vspace{.5cm}} 
\newcommand{\bvs}{\vspace{.8cm}} 

\newcommand{\dfs}{\textsc{dfs}} 
\newcommand{\haskell}{\textsc{Haskell}\xspace} 
\newcommand{\GP}{Grassmann-Pl\"ucker\xspace} 
\newcommand{\psn}{\textsc{psn}\xspace} 
\newcommand{\ppsn}{\textsc{ppsn}\xspace} 
\newcommand{\edots}{\rotatebox[origin=c]{60}{$\ddots$}} 
\newcommand{\siecle}[1]{\textsc{#1}\textsuperscript{e}~siècle} 

\begin{document}

\frontmatter 
\newcounter{compteurlof}
\setcounter{compteurlof}{0}


\aliaspagestyle{plain}{empty}

\thispagestyle{empty}

\vspace*{-2cm}
\enlargethispage{2cm}

\noindent
\includegraphics[height=2.5cm]{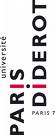}
\hspace{0pt plus 1 fill}
\begin{minipage}[b][2.5cm][c]{10cm}
\begin{center}
\textsc{Université Paris Diderot -- Paris 7} \\
\svs
\textsc{Universidad de Cantabria -- Santander}
\end{center}
\end{minipage}
\hspace{0pt plus 1 fill}
\includegraphics[height=2.5cm]{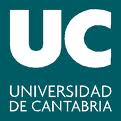}


\vspace{0pt plus 1 fill}

\begin{center}

{\fontseries{m}\fontshape{n}\fontsize{17}{17}\selectfont
THÈSE
}

\mvs

présentée et soutenue publiquement le 31 mai 2010 pour l'obtention du

\mvs

{\fontseries{m}\fontshape{n}\fontsize{17}{17}\selectfont
DOCTORAT DE L'UNIVERSITÉ PARIS DIDEROT -- PARIS~7
}

\mvs

(spécialité Informatique)

\vspace{1cm}

{\fontseries{m}\fontshape{n}\fontsize{23}{28}\selectfont
Multitriangulations, pseudotriangulations \\
et quelques problèmes de réalisation de polytopes

}

\bvs

{\fontseries{m}\fontshape{n}\fontsize{23}{23}\selectfont
Vincent \textsc{Pilaud}
}

\vspace{0pt plus 1 fill}

\begin{minipage}[t]{14cm}
Membres du jury~: \\
\hspace*{.5cm} Mireille \textsc{Bousquet-Melou} (CNRS, Univ. Bordeaux 1) \hfill \textit{Examinatrice} \\
\hspace*{.5cm}  Christian \textsc{Choffrut} (Univ. Paris 7) \hfill \textit{Examinateur} \\
\hspace*{.5cm}  Jakob \textsc{Jonsson} (KTH Stockholm) \hfill \textit{Rapporteur} \\
\hspace*{.5cm}  Marc \textsc{Noy} (UPC Barcelona) \hfill \textit{Examinateur} \\
\hspace*{.5cm}  Michel \textsc{Pocchiola} (Univ. Paris 6) \hfill \textit{Directeur de thèse} \\
\hspace*{.5cm}  Jorge \textsc{Ramirez Alfonsin} (Univ. Montpellier 2) \hfill \textit{Examinateur} \\
\hspace*{.5cm}  Francisco \textsc{Santos Leal} (Univ. Cantabria) \hfill \textit{Directeur de thèse}

\svs
Autre rapporteur~:\\
\hspace*{.5cm}  Stefan \textsc{Felsner} (TU Berlin) \hfill \textit{Rapporteur}
\end{minipage}

\vspace{0pt plus 1 fill}

\begin{minipage}[t]{5.3cm}
\begin{center}
\includegraphics[height=1.5cm]{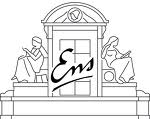}\\
Laboratoire d'Informatique de l'École Normale Supérieure
\end{center}
\end{minipage}
\begin{minipage}[t]{5.3cm}
\begin{center}
\includegraphics[height=1.5cm]{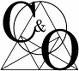}\\
Équipe Combinatoire et Optimisation de Paris 6
\end{center}
\end{minipage}
\begin{minipage}[t]{5.3cm}
\begin{center}
\includegraphics[height=1.5cm]{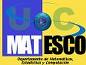}\\
Departamento de Matemáticas, Estadística y Computación, Universidad de Cantabria
\end{center}
\end{minipage}

\end{center}

\cleardoublepage

\thispagestyle{empty}

\vspace*{-8.6cm}\hspace*{12.5cm}\includegraphics[scale=1]{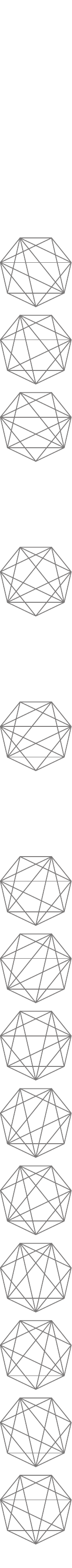}
\vspace*{-25cm}


\begin{center}

\vspace*{.3cm}

{\fontseries{m}\fontshape{n}\fontsize{23}{23}\selectfont
\darkblue
Multitriangulations, pseudotriangulations \\ et quelques problèmes de réalisation de polytopes

\vspace{1.7cm}

Multitriangulaciones, pseudotriangulaciones \\ y algunos problemas de realización de politopos

\vspace{1.7cm}

Multitriangulations, pseudotriangulations \\ and some problems of realization of polytopes

}

\vspace{1.5cm}
{\black\rule{7cm}{.5pt}}
\vspace{1.5cm}

{\fontseries{m}\fontshape{sc}\fontsize{22}{22}\selectfont
Vincent Pilaud

}

\vspace{1.8cm}

\begin{minipage}[c]{4.3cm}
\begin{flushright}
sous la direction de \\ bajo la dirección de \\ under the supervision of

\end{flushright}
\end{minipage}
\quad
\begin{minipage}[c]{10.8cm}
\fontseries{m}\fontshape{n}\fontsize{17}{13}\selectfont
\textsc{Michel Pocchiola} \\ {\normalsize École Normale Supérieure, Paris}  \\ {\normalsize Université Pierre et Marie Curie, Paris} \\ \\
\textsc{Francisco Santos Leal} \\ {\normalsize Universidad de Cantabria, Santander}
\end{minipage}

\end{center}

\cleardoublepage

\thispagestyle{empty}
\vspace*{0 pt plus 1 fill}
\hspace{0 pt plus 3 fill} À Claude et Bernard \hspace{0 pt plus 1 fill}
\vspace*{0 pt plus 3 fill}

\cleardoublepage

\chapter*{Remerciements}
\thispagestyle{empty}
\addcontentsline{toc}{chapter}{Remerciements}

Je tiens à remercier en premier lieu mes deux directeurs de thèse, qui m'ont guidé chacun à leur manière pour mes premiers pas dans le monde de la recherche. Qu'ils sachent tous les deux ma reconnaissance et mon admiration.

À Paris, \defn{Michel Pocchiola} m'a orienté et conseillé à bien des égards depuis mon master. Sur le plan scientifique, il m'a fait découvrir ses thèmes de recherche avec patience et m'a proposé des sujets riches et stimulants. Il m'a transmis une noble vision de la recherche, que j'espère pouvoir un jour défendre aussi bien que lui. Sur le plan personnel, il a fait preuve d'attention et d'amitié tout au long de cette thèse, et a toujours su m'encourager dans les moments compliqués.

À Santander, \defn{Francisco Santos} m'a accueilli chaleureusement au sein de l'Université de Cantabrie. Toujours entrain à me faire partager avec beaucoup de pédagogie et d'enthousiasme ses connaissances intarissables et ses idées originales, il a aussi su m'écouter et rendre intéressantes les fantaisies que je lui exposais. Je veux en particulier le remercier pour l'attention qu'il a portée à la rédaction de ce manuscrit. Par ailleurs, j'ai eu l'occasion de découvrir et d'apprécier ses qualités humaines,  son altruisme, ainsi que son hospitalité lors de mes séjours à Santander.

\svs
J'adresse mes sincères remerciements à \defn{Stefan Felsner} et \defn{Jakob Jonsson} pour avoir accepté d'être rapporteurs de cette thèse. Je suis très honoré de l'intérêt qu'ils portent à mon travail, au regard de leurs diverses compétences. Merci également aux examinateurs qui ont accepté de siéger dans le jury de ma soutenance, malgré les emplois du temps chargés des uns et l'éloignement géographique des autres.

\svs
Au cours de cette thèse, j'ai eu l'occasion de faire la connaissance de nombreuses personnes passionnantes et bienveillantes. En particulier, j'ai beaucoup apprécié les relations scientifiques et personnelles que j'ai eues avec mes coauteurs. Outre mes directeurs de thèse, je pense parti\-culièrement à \defn{Julian Pfeifle} dont la passion pour la géométrie est communicative et contagieuse, et qui m'a accueilli à plusieurs reprises chez lui (j'ai ainsi pu constater que ses enfants ont su construire des polytopes avant d'apprendre à faire du vélo). Un grand merci aussi à \defn{Benjamin Matschke} pour un séjour à Berlin aussi agréable que productif, et pour ses nombreuses idées lancées entre deux balles de baby-foot. Je suis par ailleurs très reconnaissant à \defn{J\"urgen Bokowski} pour notre travail commun qui a beaucoup alimenté ma réflexion.

Je tiens à remercier les membres des laboratoires qui m'ont accueilli durant ma thèse. À l'\textsc{éns}, merci à \defn{Éric Colin de Verdière} pour son soutien et ses conseils, mais aussi à l'équipe administrative pour son efficacité et sa disponibilité. À Paris 6, j'ai apprécié la bonne humeur des membres de l'équipe C\&O, et les échanges que j'ai eus avec les uns et les autres. À Santander, je remercie les thésards avec qui j'ai partagé un bureau et mes premiers mots en espagnol.

Je veux aussi exprimer ma gratitude envers les organisateurs et les participants du \emph{DocCourse Combinatorics and Geometry 2009} de Barcelone, qu'ils ont su rendre aussi intéressant que convivial. Je n'oublierai pas le soutien qu'ils m'ont apporté dans un moment personnel difficile.

\svs
J'aimerais pouvoir dire tout ce que je dois à mes amis, en qui je puise chaque jour l'énergie de mener à bien mes projets. 
En premier lieu, \defn{Lionel} et \defn{Anthony} que je tiens à remercier (le mot est faible) pour leurs conseils, leur soutien et leurs sourires. Je veux aussi faire un clin d'oeil à mes amis \textsc{hx2} de Rhône-Alpes et d'ailleurs, qu'aucune distance ne semble pouvoir séparer. Merci encore à \defn{Pascal} et \defn{Pierre}, mes deux compères de la prépa agreg, qui m'ont toujours tiré vers le haut, et dont l'amitié me touche. Enfin, je ne sais pas dans quelle langue remercier \defn{Cassandra}~: en français parce que c'est sa langue de coeur, in English for all that she taught me, o en espagnol por todo lo que compartimos.

\svs
Je souhaite exprimer toute mon affection à ma famille. Même si tous connaissent les sentiments que j'ai pour eux, il y a des évidences qui sont parfois bonnes à dire. Merci donc tout d'abord à \defn{Claude} pour son soutien indéfectible et ses attentions permanentes. Qu'elle trouve dans cette thèse les échos de tout ce qu'elle m'a appris, et dans ces mots le témoignage de ma gratitude. Merci à \defn{Christophe} et \defn{Thomas}~: ils sont comme les deux autres côtés d'un triangle que chaque épreuve rend plus solide, et il faudrait une thèse entière pour en étudier les propriétés. Je pense aussi à \defn{DD} et \defn{Tatou} et je suis touché de lire dans leurs yeux l'admiration objective des grands-parents. Enfin, je n'oublie pas \defn{Christiane}, \defn{Yvette}, \defn{Daniel} et \defn{Gérard} qui m'ont toujours traité plus comme un fils que comme un neveu, et dont j'ai su apprécier la présence dans les jours sombres.

\svs
En terminant cette thèse, j'ai une pensée pour \defn{Bernard}, pour notre complicité, pour tout ce qu'il était pour moi. Je me souviens avec beaucoup d'émotion de nos \emph{Saturday Maths fever}, durant lesquels il m'a transmis sa passion pour les mathématiques, son goût de l'esthétique et de la simplicité des démonstrations. J'aurais aimé qu'il puisse feuilleter cette thèse comme il aurait déambulé dans les dédales d'un musée, en s'émerveillant devant chaque figure comme face à un tableau de Rothko.

\cleardoublepage

\aliaspagestyle{plain}{baspage}


\pagestyle{hautpageintro}


\renewcommand{\namepart}{Contents}
\tableofcontents

\newpage

\renewcommand{\namepart}{List of figures}
\listoffigures



\mainmatter 
\setcounter{compteurlof}{1}


\counterwithout{section}{chapter}
\counterwithout{subsection}{chapter}
\counterwithin{subsection}{section}
\counterwithout{figure}{chapter}

\setcounter{chapter}{-1}
\addtocontents{toc}{ \vspace{.5cm} \demiligne \vspace{.1cm} }
\renewcommand*{\figurename}{Figure}
\makeatletter\renewcommand{\thefigure}{\@arabic\c@figure~}\makeatother
\renewcommand{\namepart}{Résumé}

\theoremstyle{theorem} 
\newtheorem{theorem_F}{Théorème}
\newtheorem{proposition_F}[theorem_F]{Proposition}
\newtheorem{definition_F}[theorem_F]{Définition}

\theoremstyle{definition} 
\newtheorem{example_F}[theorem_F]{Exemple}
\newtheorem{remark_F}[theorem_F]{Remarque}
\newtheorem{observation_F}[theorem_F]{Observation}

\newcommand{\kcrois}[1]{\mbox{$#1$-croi}\-se\-ment}
\newcommand{\kpert}[1]{\mbox{$#1$-perti}\-nen\-te}
\newcommand{\kbord}[1]{\mbox{$#1$-bord}}
\newcommand{\ketoile}[1]{\mbox{$#1$-étoile}}
\newcommand{\squelette}[1]{\mbox{$#1$-squelette}}

\chapter*{Résumé}\label{chap:resumeF}
\phantomsection
\addcontentsline{toc}{chapter}{Résumé}


\section{Introduction}

Les thèmes abordés dans cette thèse s'inscrivent dans le champ de la \defn{géométrie discrète et \mbox{algorithmique}}~\cite{hdcg-04}~: les problèmes rencontrés font en général intervenir des ensembles finis d'objets géométriques élémentaires (points, droites, demi-espaces, \etc), et les questions portent sur la manière dont ils sont reliés, arrangés, placés les uns par rapport aux autres (\mbox{comment} s'inter\-sectent-ils~? comment se voient-ils~? que délimitent-ils~? \etc). Cette thèse traite de deux sujets particuliers~: les \defn{multitriangulations} et les \defn{réalisations polytopales de produits}. Leurs diverses connexions avec la géométrie discrète nous en ont fait découvrir de nombreuses facettes, dont nous donnons quelques exemples dans cette introduction. Nous commençons par présenter leur problématique commune, la recherche d'une réalisation polytopale d'une structure donnée, qui a orienté nos recherches sur ces deux sujets.

\svs
Un polytope (convexe) est l'enveloppe convexe d'un ensemble fini de points d'un espace euclidien. Si l'intérêt porté à certains polytopes remonte à l'Antiquité (solides de Platon), leur étude systématique est relativement récente et les principaux résultats datent du \siecle{xx} (voir~\cite{g-cp-03,z-lp-95} et leurs références). L'étude des polytopes porte non seulement sur leurs propriétés géométriques mais surtout sur leurs aspects plus combinatoires. Il s'agit notamment de comprendre leurs faces (leurs intersections avec un hyperplan support) et le treillis qu'elles composent (c'est-à-dire les relations d'inclusion entre ces faces).

Les questions de \defn{réalisation polytopale} forment en quelque sorte le problème inverse~: elles portent sur l'existence et la construction de polytopes à partir d'une structure combinatoire donnée. Par \mbox{exemple}, étant donné un graphe, on voudrait déterminer s'il s'agit du graphe d'un polytope~: on dira alors que le graphe est polytopal. Le résultat fondateur dans ce domaine est le Théorème de Steinitz~\cite{s-pr-22} qui caractérise les graphes des \poly{3}topes. Dès la dimension~$4$, la situation est nettement moins satisfaisante~: malgré certaines conditions nécessaires~\cite{b-gscps-61,k-ppg-64,b-ncp-67}, les graphes de polytopes n'admettent pas de caractérisation locale en dimension générale~\cite{rg-rsp-96}. Lorsqu'un graphe est polytopal, on s'interroge sur les propriétés des éventuelles réalisations, par exemple leur nombre de faces, leur dimension, \etc{} On cherche à construire des exemples qui optimisent certaines de ces propriétés~: typiquement, on veut cons\-truire un polytope de dimension minimale ayant un graphe donné~\cite{g-ncp-63,jz-ncp-00,sz-capdp}. Ces questions de réalisations polytopales sont intéressantes pour des graphes qui proviennent soit de graphes de transformation sur des ensembles combinatoires ou géométriques (graphe des transpositions adjacentes sur les permutations, graphe des flips sur les triangulations, \etc), soit d'opérations sur des graphes (qui peuvent être locales comme la transformation~$\Delta Y$, ou globales comme le produit cartésien). Ces questions ne se posent pas uniquement pour un graphe, mais plus généralement pour n'importe quel sous-ensemble de treillis.


\paragraph{Polytopalité de graphes de flips.} L'existence de réalisations polytopales est d'abord étudiée pour des graphes de transformation sur des structures combinatoires ou géométriques. Citons ici le permutoèdre dont les sommets correspondent aux permutations de~$[n]$ et dans lequel deux sommets sont reliés par une arête si les permutations correspondantes diffèrent d'une transposition de deux positions adjacentes. D'autres exemples, ainsi que des classes de polytopes permettant de réaliser certaines structures sont exposés dans~\cite[Lecture~9]{z-lp-95}. D'une manière générale, les questions de polytopalité de structures combinatoires sont intéressantes non seulement pour leurs résultats, mais aussi parce que leur étude force à comprendre la combinatoire des objets et amène à développer des méthodes nouvelles.

Deux exemples particuliers de structures combinatoires polytopales ont un rôle important dans le cadre de cette thèse. Nous rencontrons d'abord l'associaèdre, dont le bord réalise le dual du complexe simplicial formé par tous les ensembles de cordes du \gon{n}e qui ne se croisent pas. Le graphe de l'associaèdre correspond au graphe des flips sur les triangulations du \gon{n}e. L'associaèdre apparaît dans divers contextes et plusieurs réalisations polytopales ont été proposées~\cite{l-atg-89,bfs-ccsp-90,gkz-drmd-94,l-rsp-04,hl-rac-07}. Nous rencontrons ensuite le polytope des pseudotriangulations d'un ensemble de points du plan euclidien~\cite{rss-empppt-03}. Introduites pour l'étude du complexe de visibilité d'obstacles convexes disjoints du plan~\cite{pv-tsvcp-96,pv-vc-96}, les pseudotriangulations ont été utilisées dans divers contextes géométriques~\cite{rss-pt-06}. Leurs propriétés de rigidité~\cite{s-ptrmp-05} ont permis d'établir la polytopalité de leur graphe de flips~\cite{rss-empppt-03}.

\begin{figure}
	\capstart
	\centerline{\includegraphics[width=\textwidth]{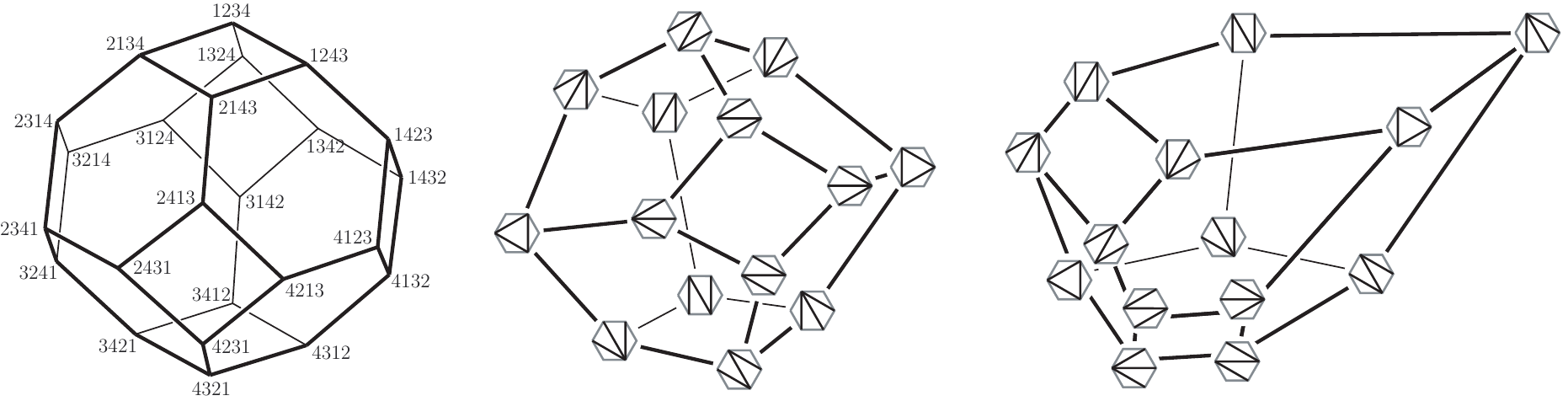}}
	\caption[Le permutoèdre et deux réalisations de l'associaèdre]{Le permutoèdre et deux réalisations de l'associaèdre.}
	\label{F:fig:permutohedronassociahedra}
\end{figure}

\svs
\enlargethispage{.1cm}
Dans la première partie de cette thèse, nous nous intéressons à la polytopalité du graphe des flips sur les multitriangulations. Ces objets, apparus de manière relativement contingente~\cite{cp-tttccp-92,dkm-lahp-02,dkkm-2kn..-01}, possèdent une riche structure combinatoire~\cite{n-gdfcp-00,j-gt-03,j-gtdfssp-05}. Une \defn{\ktri{k}} est un ensemble maximal de cordes du \gon{n}e ne contenant pas de sous-ensemble de $k+1$ cordes qui se croisent deux à deux. Nous consi\-dérons le graphe des flips dans lequel~deux multitriangulations sont reliées si elles diffèrent d'une corde. Comme pour les triangulations que l'on retrouve lorsque~$k=1$, ce graphe est régulier et connexe, et nous nous interrogeons sur sa polytopalité. Jakob Jonsson~\cite{j-gt-03} a fait un premier pas dans cette direction en montrant que le complexe simplicial des ensembles de cordes ne contenant pas de sous-ensemble de $k+1$ cordes qui se croisent deux à deux est une sphère topologique. Bien que nous n'ayons que des réponses partielles à cette question, elle a entraîné des résultats inattendus que nous exposons dans cette thèse.

\svs
Plusieurs constructions de l'associaèdre~\cite{bfs-ccsp-90,l-rsp-04} sont basées, directement ou indirectement, sur les triangles des triangulations. Pour les multitriangulations, aucun objet élémentaire similaire n'apparaît dans les travaux antérieurs. Nous avons donc cherché dans un premier temps à comprendre ce que deviennent les triangles dans les multitriangulations. Les \defn{étoiles} que nous introduisons au Chapitre~\ref{chap:stars} répondent à cette question. Au même titre que les \mbox{triangles} dans les triangulations, nous estimons que ces étoiles donnent le bon point de vue pour comprendre les multitriangulations. Pour preuve, nous commençons par retrouver à l'aide de ces étoiles toutes les propriétés combinatoires de base connues jusqu'alors sur les multitriangulations. D'abord, nous étudions les relations d'incidence entre les étoiles et les cordes (chaque corde interne est contenue dans deux étoiles) qui nous permettent de retrouver que toutes les \ktri{k}s du \gon{n}e ont même cardinal. Ensuite, en considérant les bissectrices des étoiles, nous donnons une interprétation locale de l'opération de flip (une corde interne est remplacée par l'unique bissectrice commune des deux étoiles qui lui sont adjacentes), ce qui éclaire l'étude du graphe des flips et de son diamètre. Nous redéfinissons par ailleurs en termes d'étoiles des opérations inductives sur les multitriangulations qui permettent d'ajouter ou de supprimer un sommet au \gon{n}e. Enfin, nous utilisons la décomposition d'une multitriangulation en étoiles pour l'interpréter comme une décomposition polygonale d'une surface, et nous appliquons cette interprétation à la construction de décompositions régulières de surfaces.

\begin{figure}
	\capstart
	\centerline{\includegraphics[scale=1]{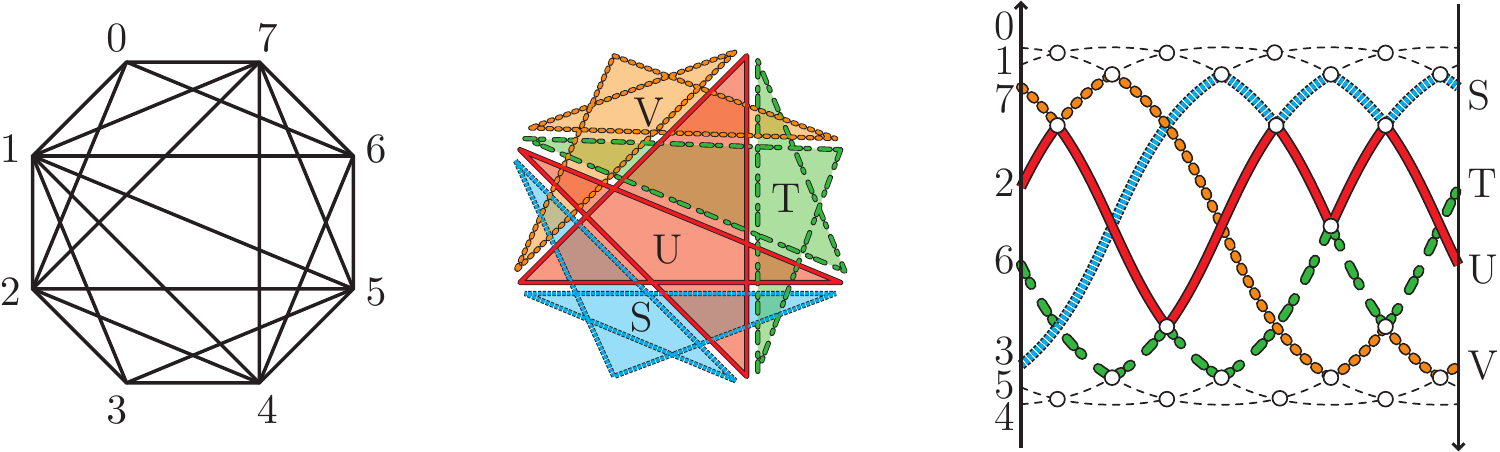}}
	\caption[Une \ktri{2}, sa décomposition en étoiles, et son arrangement dual]{Une \ktri{2}, sa décomposition en étoiles, et son arrangement dual.}
	\label{F:fig:ktristardual}
\end{figure}

\svs
Dans un deuxième temps, nous avons cherché à comprendre les multitriangulations par dua\-lité. L'espace des droites du plan est un ruban de M\"obius~; l'ensemble des droites qui passent par un point du plan est une pseudodroite du ruban de M\"obius~; et les pseudodroites duales d'une configuration de points forment un arrangement de pseudodroites~\cite{g-pa-97}.

Notre point de départ est la dualité entre les pseudotriangulations d'un ensemble~$P$ de points et les arrangements de pseudodroites dont le support est l'arrangement dual de~$P$ privé de son premier niveau~\cite{pv-ot-94,pv-ptta-96}. Nous établissons une dualité similaire pour les multitriangulations. D'une part, l'ensemble des bissectrices d'une étoile est une pseudodroite du ruban de M\"obius. D'autre part, les pseudodroites duales aux étoiles d'une \ktri{k} du \gon{n}e forment un arrangement de pseudodroites avec points de contact, supporté par l'arrangement dual des sommets du \gon{n}e privé de ses $k$ premiers niveaux. Nous montrons que tout arrangement avec points de contact ayant ce support est bien l'arrangement dual d'une multitriangulation. Cette dualité fait le lien entre les pseudotriangulations et les multitriangulations et explique ainsi leurs propriétés communes (nombre de cordes, flip, \etc).

Plus généralement, nous étudions au Chapitre~\ref{chap:mpt} les \defn{arrangements de pseudodroites avec points de contact} qui partagent un même support. Nous définissons une opération de flip qui correspond au flip des multitriangulations et des pseudotriangulations, et nous étudions le graphe des flips. Les propriétés de certains arrangements gloutons, définis comme les sources de certaines orientations acycliques de ce graphe, nous permettent en particulier d'énumérer ce graphe de flips en gardant un espace polynomial. Notre travail éclaire ainsi l'algorithme existant d'énumération des pseudotriangulations~\cite{bkps-ceppgfa-06} et en donne une preuve complémentaire.

\svs
Pour finir, nous donnons au Chapitre~\ref{chap:multiassociahedron} des pistes pour l'étude de trois problèmes ouverts qui reflètent la richesse combinatoire et géométrique des multitriangulations.

Le premier concerne le comptage des multitriangulations. Jakob Jonsson~\cite{j-gtdfssp-05} a prouvé (grâce à des considérations sur des $0/1$-remplissages de polyominos qui évitent certains motifs) que les multitriangulations sont comptées par un certain déterminant de Hankel de nombres de Catalan, qui compte aussi certaines familles de $k$-uplets de chemins de Dyck. Cependant, mis à part des résultats partiels~\cite{e-btdp-07,n-abtdp-09}, on ne connaît pas de preuve bijective de ce résultat. On entre ici dans le domaine de la \defn{combinatoire bijective} dont le but est de  construire des bijections entre des familles combinatoires qui conservent des paramètres caractéristiques~: on voudrait ici exhiber une bijection qui permette de lire les étoiles sur les $k$-uplets de chemins de Dyck.

Le second problème est celui de la \defn{rigidité}. Alors que la rigidité des graphes est bien comprise en dimension~$2$~\cite{l-grpss-70,g-cf-01,gss-cr-93}, aucune caractérisation satisfaisante n'est connue à partir de la dimension~$3$. Nous montrons que certaines propriétés typiques des graphes rigides en dimension~$2k$ sont vérifiées par les \ktri{k}s. Ceci amène naturellement à conjecturer que les \ktri{k}s sont rigides en dimension~$2k$, ce que nous prouvons lorsque~$k=2$. Une réponse positive à cette conjecture permettrait de se rapprocher de la polytopalité du graphes des flips sur les multitriangulations, de la même façon que le polytope des pseudotriangulations~\cite{rss-empppt-03} a pu être construit en s'appuyant sur leurs propriétés de rigidité.

Enfin, nous revenons sur la \defn{réalisation polytopale du graphe des flips} sur les multitriangulations. Dans un premier temps, nous étudions le premier exemple non trivial~: nous montrons que le graphe des flips sur les \ktri{2}s de l'octogone est bien le graphe d'un polytope de dimension~$6$. Pour trouver un tel polytope, nous décrivons complètement l'espace des réalisations polytopales symétriques de ce graphe en dimension~$6$, en étudiant d'abord tous les matroides orientés~\cite{bvswz-om-99,b-com-06} symétriques qui peuvent le réaliser. Dans un deuxième temps, en généralisant une construction de l'associaèdre due à Jean-Louis Loday~\cite{l-rsp-04}, nous construi\-sons un polytope qui réalise le graphe des flips restreint aux multitriangulations dont le graphe dual orienté est acyclique.

\svs
Par ailleurs, nous présentons dans l'Annexe~\ref{app:implementations} les résultats d'un algorithme d'énumération des petits arrangements de pseudodroites et de \defn{double pseudodroites}. À l'instar des arrangements de pseudodroites qui offrent un modèle combinatoire des configurations de points, les arrangements de double pseudodroites ont été introduits comme modèles de configurations de convexes disjoints~\cite{hp-adp-08}. Notre travail d'implémentation nous a permis de manipuler ces objets et de nous familiariser avec leurs propriétés, ce qui s'est révélé utile lors de notre étude de la dualité.


\paragraph{Polytopalité de produits cartésiens.} Le \defn{produit cartésien} de graphes est défini de sorte que le graphe d'un produit de polytopes est le produit des graphes de ses facteurs. Ainsi, un produit de graphes polytopaux est automatiquement polytopal. Nous étudions dans un premier temps la réciproque~: la polytopalité d'un produit de graphes implique-t-elle celle de ses facteurs~? Le Chapitre~\ref{chap:nonpolytopal} apporte des éléments de réponse en prêtant une attention particulière aux graphes réguliers et à leur réalisation comme polytopes simples (les polytopes simples ont des propriétés très particulières~: entre autres, ils sont déterminés par leur graphe~\cite{bm-ppi-87,k-swtsp-88}). Nous discutons en particulier la question de la polytopalité du produit de deux graphes de \mbox{Petersen}, qui a été soulevée par G\"unter Ziegler~\cite{crm}. Ce travail sur la polytopalité des produits nous a conduit à étudier des exemples de graphes qui se révèlent non-polytopaux bien qu'ils satisfassent les conditions nécessaires connues pour être polytopaux~\cite{b-gscps-61,k-ppg-64,b-ncp-67}.

\begin{figure}
	\capstart
	\centerline{\includegraphics[scale=1]{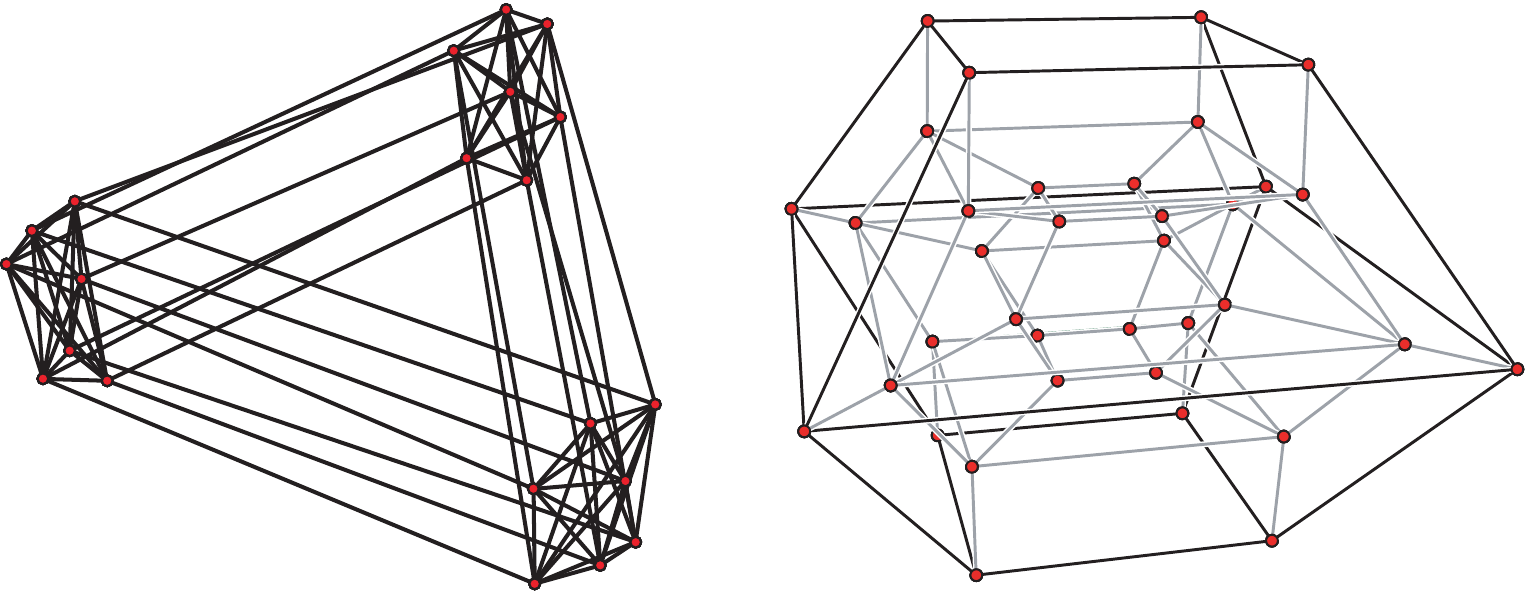}}
	\caption[Exemples de produits de graphes]{(Gauche) Un produit de graphes complets. Ce graphe est celui d'un produit de simplexes, mais aussi celui de polytopes de dimensions inférieures. (Droite) Un produit polytopal de graphes non-polytopaux.}
	\label{F:fig:products}
\end{figure}

\svs
Nous recherchons ensuite la dimension minimale d'un polytope dont le graphe est isomorphe à celui d'un produit fixé de polytopes. Cette question s'inscrit dans la recherche systématique de polytopes de dimension extrémale ayant des propriétés fixées~: par exemple, les \poly{d}topes \emph{neighborly} \cite{g-ncp-63} sont ceux dont le \squelette{\Fracfloor{d}{2}} est complet (et qui atteignent le maximum de faces autorisé par le Théorème de la Borne Supérieure~\cite{m-mnfcp-70})~; les \poly{d}topes \emph{neighborly} cubiques~\cite{jz-ncp-00,z-ppp-04,js-ncps-07} sont ceux dont le \squelette{\Fracfloor{d}{2}} est celui d'un $n$-cube,~\etc{} Pour construire des polytopes de petite dimension ayant une propriété fixée, il est naturel et souvent efficace de partir d'une dimension suffisante pour s'assurer de l'existence de tels polytopes, et de projeter ces polytopes sur des sous-espaces de dimension inférieure en conservant la propriété recherchée. Ces techniques et leurs limites ont été largement étudiées dans la littérature~\cite{az-dpmsp-99,z-ppp-04,sz-capdp,s-tovnms-09,sanyal-phd}. Nous les appliquons au Chapitre~\ref{chap:psn} pour construire des polytopes \defn{$(k,\sub{n})$-neighborly prodsimpliciaux}, qui ont le même \squelette{k} qu'un certain produit de simplexes~$\simplex_{\sub{n}} \eqdef \simplex_{n_1}\times\cdots\times\simplex_{n_r}$. Dans ce chapitre, nous donnons aussi des coordonnées entières pour de tels polytopes à l'aide de sommes de Minkowski explicites de polytopes~cycliques.


\section{Multitriangulations}

Fixons~$n$ sommets sur le cercle unité et considérons les cordes entre ces sommets. On dit que deux cordes se croisent lorsque les segments ouverts s'intersectent. Un \defn{\kcrois{\ell}} est un ensemble de~$\ell$ cordes qui se croisent deux à deux. À l'instar des triangulations, nous nous intéressons aux ensembles maximaux qui évitent ces motifs~:

\begin{definition_F}
Une \defn{\ktri{k}} du \gon{n}e est un ensemble maximal de cordes du \gon{n}e sans \kcrois{(k+1)}.
\end{definition_F}

Par définition, une corde ne peut apparaître dans un \kcrois{(k+1)} que si elle a au moins~$k$ sommets de chaque côté. Nous disons qu'une telle corde est \defn{\kpert{k}}. Par maximalité, toute \ktri{k} contient toutes les cordes qui ne sont pas \kpert{k}s. Parmi ces cordes, celles qui séparent~$k-1$ sommets de tous les autres sommets jouent un rôle particulier~: nous les appelons \defn{cordes du \kbord{k}}.

\begin{example_F}\label{F:exm:petitscas}
Pour certaines valeurs de~$k$ et~$n$, il est facile de décrire les \ktri{k}s~du~\gon{n}e~:
\begin{itemize}
\item[~~~~\fbox{$k=1$}] Les \ktri{1}s sont juste les triangulations du \gon{n}e.
\item[~~~~\fbox{$n=2k+1$}] Le graphe complet~$K_{2k+1}$ est l'unique \ktri{k} du \gon{(2k+1)}e car aucune de ses cordes n'est \kpert{k}.
\item[~~~~\fbox{$n=2k+2$}] L'ensemble des cordes du \gon{(2k+2)}e contient exactement un \kcrois{(k+1)} formé par les~$k+1$ diagonales du \gon{(2k+2)}e reliant deux points antipodaux. Il y a donc~$k+1$ \ktri{k}s du \gon{(2k+2)}e, obtenues en retirant au graphe complet l'une de ses longues diagonales.
\item[~~~~\fbox{$n=2k+3$}] Les cordes \kpert{k}s du \gon{(2k+3)}e forment un cycle polygonal dans lequel deux cordes se croisent toujours sauf si elles sont consécutives. Par conséquent, les \ktri{k}s du \gon{(2k+3)}e sont les unions disjointes de~$k$ paires de cordes \kpert{k}s consécutives (auxquelles on ajoute toutes les cordes qui ne sont pas \kpert{k}s).
\end{itemize}
\end{example_F}


\subsection*{Résultats antérieurs}

Les ensembles de cordes du \gon{n}e sans \kcrois{(k+1)} apparaissent dans le cadre de la théorie extrémale des graphes géométriques (voir~\cite[Chapitre~14]{pa-cg-95}, \cite[Chapitre~1]{f-gga-04} et la discussion dans~\cite{cp-tttccp-92}). Dans le graphe d'intersection des cordes du \gon{n}e, ces ensembles induisent en effet des sous-graphes sans \clique{k}, et sont donc à rapprocher du résultat classique de Tur\'an qui borne le nombre d'arêtes d'un graphe sans \clique{k}. Vasilis~Capoyleas et Janos Pach~\cite{cp-tttccp-92} ont montré que ces ensembles ne peuvent pas avoir plus de $k(2n-2k-1)$ cordes. Tomoki \mbox{Nakamigawa}~\cite{n-gdfcp-00}, et indépendamment Andreas Dress, Jacobus Koolen et Vincent Moulton~\cite{dkm-lahp-02}, ont prouvé que toutes les \ktri{k}s atteignent cette borne. Ces deux preuves sont basées sur l'opération de flip qui transforme une \ktri{k} en une autre en changeant la position d'une  seule corde. Tomoki~Nakamigawa~\cite{n-gdfcp-00} a montré que toute corde \kpert{k} dans une \ktri{k} peut être flippée, et que le graphe des flips est connexe. Nous tenons à observer que dans tous ces travaux antérieurs, les résultats sont obtenus de ``manière \mbox{indirecte}''~: d'abord, ils introduisent une opération de contraction (similaire à la contraction d'une arête du bord dans une triangulation), puis ils l'utilisent pour démontrer l'existence du flip (en toute généra\-lité~\cite{n-gdfcp-00}, ou seulement dans des cas particuliers~\cite{dkm-lahp-02}) et la connexité du graphe des flips, dont découle la formule sur le nombre de cordes. Pour résumer~:

\begin{theorem_F}[\cite{cp-tttccp-92,n-gdfcp-00,dkm-lahp-02}]\label{F:theo:fundamental}
\begin{enumerate}[(i)]
\item Il existe une opération inductive qui transforme les \ktri{k}s du \gon{(n+1)}e en \ktri{k}s du \gon{n}e et \viceversa~\cite{n-gdfcp-00}.
\item Toute corde \kpert{k} d'une \ktri{k} du \gon{n}e peut être flippée 
et le graphe des flips est régulier et connexe~\cite{n-gdfcp-00,dkm-lahp-02}.
\item Toute \ktri{k} du \gon{n}e a~$k(2n-2k-1)$ cordes~\cite{cp-tttccp-92,n-gdfcp-00,dkm-lahp-02}.\qed
\end{enumerate}
\end{theorem_F}

Jakob~Jonsson~\cite{j-gt-03,j-gtdfssp-05} a complété ces résultats dans deux directions. D'une part, il a étudié les propriétés énumératives des multitriangulations~:

\begin{theorem_F}[\cite{j-gtdfssp-05}]\label{F:theo:enumeration}
Le nombre de \ktri{k}s du \gon{n}e est donné par~:
$$\det(C_{n-i-j})_{1\le i,j\le k}=\det\begin{pmatrix} C_{n-2} & C_{n-3} & \edots & C_{n-k-1} \\ C_{n-3} & \edots & \edots & \edots \\ \edots & \edots & \edots & C_{n-2k+1} \\ C_{n-k-1} & \edots & C_{n-2k+1} & C_{n-2k} \end{pmatrix},$$
où~$C_p \eqdef \frac{1}{p+1}{2p \choose p}$ dénote le $p$-ième nombre de Catalan.\qed
\end{theorem_F}

La preuve de ce théorème repose sur des résultats plus généraux concernant les $0/1$-remplissa\-ges de polyominos maximaux pour certaines restrictions sur leurs suites diagonales de~$1$.

Lorsque~$k=1$, on retrouve simplement les nombres de Catalan qui comptent non seulement les triangulations, mais aussi les chemins de Dyck. Il se trouve que le déterminant de Hankel qui apparaît dans le théorème précédent compte aussi certaines familles de $k$-uplets de chemins de Dyck qui ne se croisent pas~\cite{gv-bdphlf-85}. L'égalité des cardinaux de ces deux familles combinatoires motive la recherche d'une preuve bijective qui pourrait éclairer le résultat du Théorème~\ref{F:theo:enumeration}. Même si Sergi~Elizalde~\cite{e-btdp-07} et Carlos Nicolas~\cite{n-abtdp-09} ont proposé deux bijections différentes dans le cas~$k=2$, cette question reste ouverte dans le cas général. Nous discutons ce problème au Chapitre~\ref{chap:multiassociahedron}.

\svs
D'autre part, Jakob~Jonsson~\cite{j-gt-03} a étudié le complexe simplicial~$\Delta_{n,k}$ formé par les ensembles de cordes  \kpert{k}s du \gon{n}e sans \kcrois{(k+1)}. Comme tous les éléments maximaux de ce complexe ont~$k(n-2k-1)$ éléments, il est pur de dimension~$k(n-2k-1)-1$. Jakob~Jonsson a en fait démontré le théorème suivant~:

\begin{theorem_F}[\cite{j-gt-03}]
Le complexe simplicial~$\Delta_{n,k}$ est une sphère épluchable de dimension $k(n-2k-1)-1$.\qed
\end{theorem_F}

Volkmar Welker a conjecturé que ce complexe simplicial est même polytopal. Cette conjecture est vraie lorsque~$k=1$ (le complexe simplicial des ensembles sans croisement de cordes internes du \gon{n}e est isomorphe au complexe de bord du dual de l'associaèdre) et dans les cas de l'Exemple~\ref{F:exm:petitscas} (où l'on obtient respectivement un point, un simplexe et un polytope cyclique). Nous discutons cette question au Chapitre~\ref{chap:multiassociahedron}.


\subsection*{Décompositions des multitriangulations en étoiles}

Notre contribution à l'étude des \ktri{k}s est basée sur leurs étoiles, qui généralisent les triangles des triangulations~:

\begin{definition_F}
Une \defn{\ketoile{k}} est un polygone (non-simple) avec $2k+1$~sommets $s_0,\dots,s_{2k}$ circulairement ordonnés sur le cercle unité, et~$2k+1$ cordes $[s_0,s_k],[s_1,s_{k+1}],\dots,[s_{2k},s_{k-1}]$.
\end{definition_F}

Les étoiles jouent pour les multitriangulations exactement le même rôle que les triangles pour les triangulations~: elles les décomposent en entités géométriques plus simples et permettent d'en comprendre la combinatoire. Pour preuve, l'étude des propriétés d'incidence des étoiles dans les multitriangulations nous amène à retrouver au Chapitre~\ref{chap:stars} toutes les propriétés combinatoires des multitriangulations connues jusqu'alors. Notre premier résultat structurel est le suivant~:

\begin{theorem_F}\label{F:theo:incidences}
Soit~$T$ une \ktri{k} du \gon{n}e (avec~$n\ge 2k+1$).
\begin{enumerate}[(i)]
\item Une corde \kpert{k} est contenue dans deux \ketoile{k}s de~$T$, une de chaque côté~; une corde du \kbord{k} est contenue dans une \ketoile{k} de~$T$~; une corde qui n'est ni \kpert{k} ni du \kbord{k} n'est contenue dans aucune \ketoile{k} de~$T$.
\item $T$~a exactement~$n-2k$ \ketoile{k}s et~$k(2n-2k-1)$ cordes.\qed
\end{enumerate}
\end{theorem_F}

\begin{figure}[h]
	\centerline{\includegraphics[scale=1]{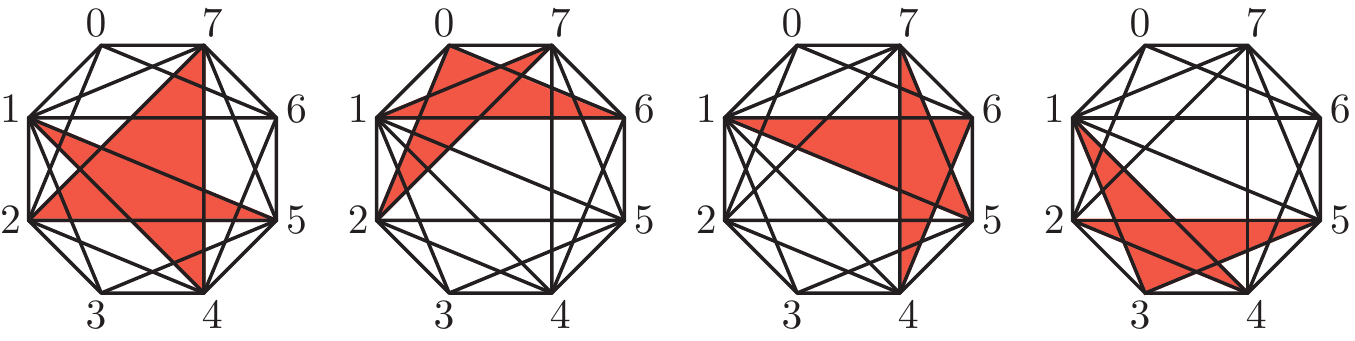}}
	\caption[Les \ketoile{2}s de la \ktri{2} de l'octogone de la~\fref{F:fig:ktristardual}]{Les \ketoile{2}s de la \ktri{2} de l'octogone de la~\fref{F:fig:ktristardual}\hspace{-.1cm}.}
	\label{F:fig:2triang8pointsstars}
\end{figure}

Nous prouvons le point~(ii) de ce théorème par un argument direct de double comptage, en exhibant deux relations indépendantes entre le nombre de cordes et le nombre de \ketoile{k}s d'une \ktri{k}. Notre première relation provient du point~(i) du théorème, tandis que la seconde repose sur les propriétés des bissectrices des \ketoile{k}s. Une \defn{bissectrice} d'une \ketoile{k} est une bissectrice de l'un de ses angles, c'est-à-dire, une droite qui passe par l'un de ses sommets et sépare ses autres sommets en deux ensembles de cardinal~$k$. Notre seconde relation est une conséquence directe  de la correspondance entre les paires de \ketoile{k}s de~$T$ et les cordes du \gon{n}e qui ne sont pas dans~$T$~:

\begin{theorem_F}
Soit~$T$ une \ktri{k} du \gon{n}e.
\begin{enumerate}[(i)]
\item Toute paire de \ketoile{k}s de~$T$ a une unique bissectrice commune, qui n'est pas dans~$T$.
\item Réciproquement, toute corde qui n'est pas dans~$T$ est la bissectrice commune d'une unique paire de \ketoile{k}s de~$T$.\qed
\end{enumerate}
\end{theorem_F}

Nous utilisons par ailleurs les étoiles et leurs bissectrices pour éclairer l'\defn{opération de flip} et pour en donner une interprétation locale. De la même façon qu'un flip dans une triangulation remplace une diagonale par une autre dans un quadrangle formé par deux triangles adjacents, un flip dans une multitriangulation s'interprète comme une transformation ne faisant intervenir que deux étoiles adjacentes~:

\begin{theorem_F}\label{F:theo:flip}
Soit~$T$ une \ktri{k} du \gon{n}e, soit~$e$ une corde \kpert{k} de~$T$ et soit~$f$ la bissectrice commune aux deux \ketoile{k}s de~$T$ qui contiennent~$e$. Alors~$T\diffsym\{e,f\}$ est une \ktri{k} du \gon{n}e et c'est la seule, avec~$T$, qui contienne~$T\ssm\{e\}$.~\qed
\end{theorem_F}

\begin{figure}[h]
	\centerline{\includegraphics[scale=1]{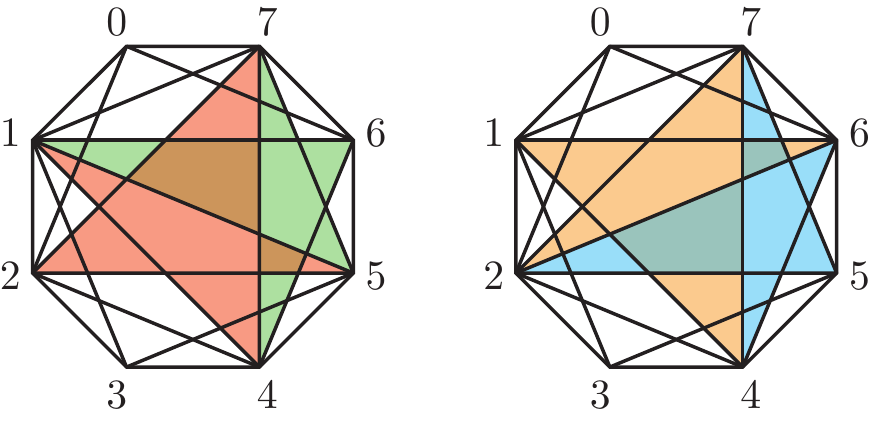}}
	\caption[Un flip dans la \ktri{2} de l'octogone de la~\fref{F:fig:ktristardual}]{Un flip dans la \ktri{2} de l'octogone de la~\fref{F:fig:ktristardual}\hspace{-.1cm}.}
	\label{F:fig:2triang8pointsflip}
\end{figure}

Cette interprétation simplifie l'étude du graphe des flips (dont les sommets sont les \ktri{k}s du \gon{n}e et dont les arêtes sont les flips entre elles) et fournit de nouvelles preuves et des extensions partielles des résultats de~\cite{n-gdfcp-00}~:

\begin{theorem_F}
Le graphe des flips sur les \ktri{k}s du \gon{n}e est $k(n-2k-1)$-régulier, connexe, et pour tout~$n>4k^2(2k+1)$, son diamètre~$\delta_{n,k}$ est borné par $$2\Fracfloor{n}{2}\left(k+\frac{1}{2}\right)-k(2k+3) \le \delta_{n,k} \le 2k(n-4k-1).$$
\end{theorem_F}
\vspace{-1cm}\qed
\vspace{1cm}

Nous utilisons aussi les étoiles pour étudier les \defn{$k$-oreilles} dans les multitriangulations, c'est-à-dire les cordes qui séparent~$k$ sommets des autres sommets. Nous proposons une preuve simple du fait que toute \ktri{k} a au moins~$2k$ $k$-oreilles~\cite{n-gdfcp-00}, puis nous donnons plusieurs caractérisations des \ktri{k}s qui atteignent cette borne~:

\begin{theorem_F}
Le nombre de $k$-oreilles d'une \ktri{k} est égal à~$2k$ plus le nombre de \defn{\ketoile{k}s internes}, \ie qui ne contiennent pas de cordes du \kbord{k}.

Si~$k>1$ et $T$ est une \ktri{k}, les assertions suivantes sont équivalentes~:
\begin{enumerate}[(i)]
\item $T$~a exactement $2k$~$k$-oreilles~;
\item $T$~n'a aucune \ketoile{k} interne~;
\item $T$~est \defn{$k$-coloriable}, \ie il existe une $k$-coloration de ses cordes \kpert{k}s telle qu'aucun croisement n'est monochromatique~;
\item l'ensemble des cordes \kpert{k}s de~$T$ est l'union disjointe de $k$ accordéons (un \defn{accordéon} est une suite de cordes~$[a_i,b_i]$ telle que pour tout~$i$, soit $a_{i+1}=a_i$ et $b_{i+1}=b_i+1$ soit $a_{i+1}=a_i-1$ et $b_{i+1}=b_i$).\qed
\end{enumerate}
\end{theorem_F}

Ensuite, nous réinterprétons en termes d'étoiles l'\defn{opération inductive} du Théorème~\ref{F:theo:fundamental}(i) qui transforme les \ktri{k}s du \gon{(n+1)}e en \ktri{k}s du \gon{n}e et \viceversa. Dans un sens on écrase une \ketoile{k} (contenant une corde du \kbord{k}) en un \kcrois{k}, et dans l'autre sens on enfle un \kcrois{k} en une \ketoile{k}. Volontairement, nous ne présentons cette opération qu'à la fin du Chapitre~\ref{chap:stars} pour souligner qu'aucune de nos preuves précédentes n'utilise cette transformation inductive.

\begin{figure}[h]
	\centerline{\includegraphics[scale=.8]{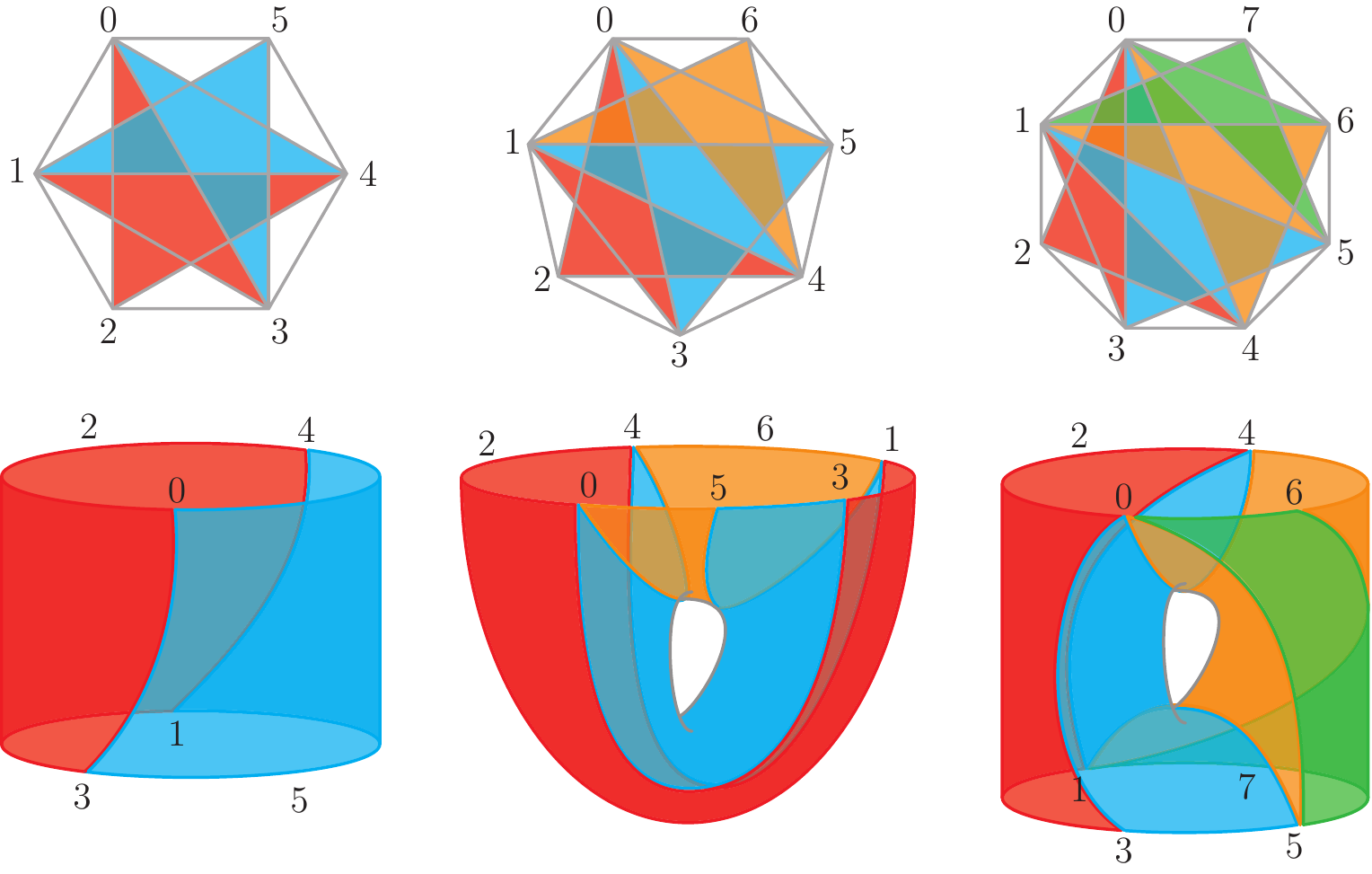}}
	\caption[Décompositions de surfaces associées aux multitriangulations]{Décompositions de surfaces associées aux multitriangulations.}
	\label{F:fig:surfaces}
\end{figure}

\svs
Pour finir le Chapitre~\ref{chap:stars}, nous utilisons le résultat du Théorème~\ref{F:theo:incidences} pour interpréter une \ktri{k} comme une \defn{décomposition de surface} en $n-2k$ \gon{k}es. Des exemples sont présentés sur la \fref{F:fig:surfaces}\hspace{-.1cm}. Nous exploitons cette interprétation pour construire, à l'aide de multitriangulations, des décompositions très régulières d'une famille infinie de surfaces.

\subsection*{Multitriangulations, pseudotriangulations et dualité}

Dans le Chapitre~\ref{chap:mpt}, nous nous intéressons à l'interprétation des multitriangulations dans l'espace des droites du plan, c'est-à-dire dans le ruban de M\"obius. Le cadre est celui de la dualité classique (voir \cite[Chapitres~5,6]{f-gga-04} et \cite{g-pa-97}) entre configurations de points et arrangements de pseudodroites~: l'ensemble~$p^*$ des droites qui passent par un point~$p$ du plan est une \defn{pseudodroite} de l'espace des droites (une courbe fermée simple non-séparatrice), et l'ensemble ${P^* \eqdef \ens{p^*}{p\in P}}$ des pseudodroites duales à un ensemble fini~$P$ de points est un \defn{arrangement de pseudodroites} (deux pseudodroites s'intersectent exactement une fois).

\begin{figure}
	\capstart
	\centerline{\includegraphics[scale=1]{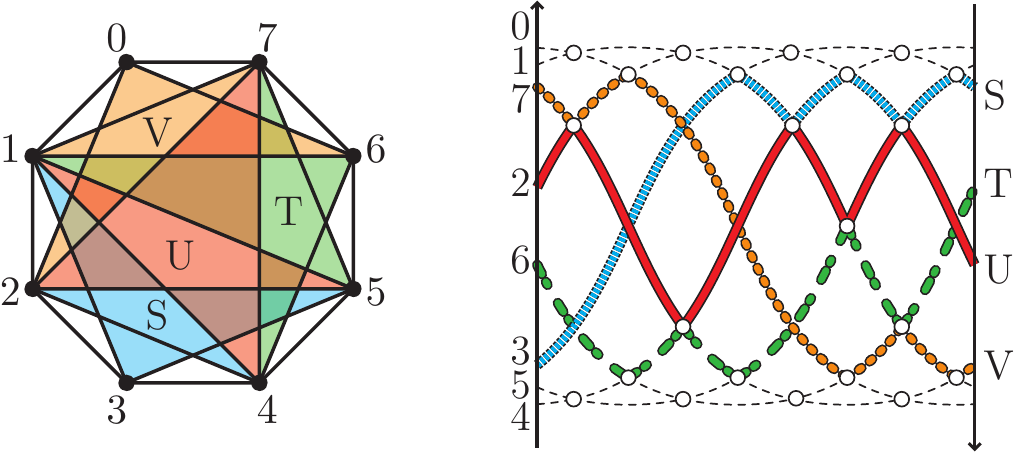}}
	\caption[Arrangement dual d'une multitriangulation]{Une \ktri{2} de l'octogone et son arrangement dual.}
	\label{F:fig:dualmulti}
\end{figure}

Le point de départ du Chapitre~\ref{chap:mpt} est l'observation suivante (voir \fref{F:fig:dualmulti}\hspace{-.1cm})~:

\begin{observation_F}
Soit $T$ une \ktri{k} du \gon{n}e. Alors~:
\begin{enumerate}[(i)]
\item l'ensemble~$S^*$ des bissectrices d'une \ketoile{k}~$S$ de~$T$ est une pseudodroite~;
\item la bissectrice commune à deux \ketoile{k}s~$R$ et~$S$ est un point de croisement des pseudodroites $R^*$ et~$S^*$, tandis qu'une corde commune à~$R$ et~$S$ est un point de contact~entre~$R^*$~et~$S^*$~;
\item l'ensemble~$T^* \eqdef \ens{S^*}{S \text{ \ketoile{k} de } T}$ de toutes les pseudodroites duales aux \ketoile{k}s de~$T$ est un \defn{arrangement de pseudodroites avec points de contact} (deux pseudodroites se croisent exactement une fois, mais peuvent se toucher en un nombre fini de contacts)~;
\item le support de cet arrangement couvre exactement le support de l'arrangement dual~$V_n^*$ des sommets du \gon{n}e privé de ses $k$~premiers niveaux.
\end{enumerate}
\end{observation_F}

Nous montrons la réciproque de cette observation~:

\begin{theorem_F}
Tout arrangement de pseudodroites avec points de contact dont le support couvre exactement le support de l'arrangement dual~$V_n^*$ des sommets du \gon{n}e privé de ses $k$~premiers niveaux est l'arrangement dual d'une \ktri{k} du \gon{n}e.\qed
\end{theorem_F}

Dans~\cite{pv-ot-94,pv-ptta-96}, une observation similaire avait été faite pour les pseudotriangulations pointées d'un ensemble de points en position générale. Introduites pour l'étude du complexe de visibilité d'obstacles convexes disjoints du plan~\cite{pv-tsvcp-96,pv-vc-96}, les pseudotriangulations ont été utilisées dans divers contextes géométriques (tels que la planification de mouvements et la rigidité \cite{s-ptrmp-05,horsssssw-pmrgpt-05}) et leurs propriétés ont été largement étudiées (nombre de pseudotriangulations~\cite{aaks-cmpt-04,aoss-nptcps-08}, polytope des pseudotriangulations~\cite{rss-empppt-03}, considérations algorithmiques~\cite{b-eptp-05,bkps-ceppgfa-06,hp-cpcpb-07}, \etc). Nous renvoyons à~\cite{rss-pt-06} pour plus de détails.

\begin{definition_F}
Un \defn{pseudotriangle} est un polygone~$\Delta$ avec exactement trois angles convexes reliés par trois chaînes polygonales concaves. Une droite est \defn{tangente} à~$\Delta$ si elle passe par un angle convexe de~$\Delta$ et sépare ses deux arêtes adjacentes, ou si elle passe par un angle concave de~$\Delta$ et ne sépare pas ses deux arêtes adjacentes. Une \defn{pseudotriangulation} d'un ensemble~$P$ de points en position générale est un ensemble d'arêtes qui décompose l'enveloppe convexe de~$P$ en pseudotriangles.
\end{definition_F}

Dans cette thèse, nous ne considérons que des pseudotriangulations \defn{pointées}, c'est-à-dire les pseudotriangulations telles que par chaque point~$p\in P$ passe une droite qui définit un demi-plan fermé ne contenant aucune des arêtes issues de~$p$. Pour une pseudotriangulation pointée~$T$ de~$P$, il a été observé~\cite{pv-ot-94,pv-ptta-96} que~:
\begin{enumerate}[(i)]
\item l'ensemble~$\Delta^*$ des droites tangentes à un pseudotriangle~$\Delta$ est une pseudodroite~; et
\item l'ensemble~$T^* \eqdef \ens{\Delta^*}{\Delta \text{ pseudotriangle de } T}$ est un arrangement de pseudodroites avec points de contact supporté par l'arrangement~$P^*$ dual de~$P$ privé de son premier niveau.
\end{enumerate}

\begin{figure}
	\capstart
	\centerline{\includegraphics[scale=1]{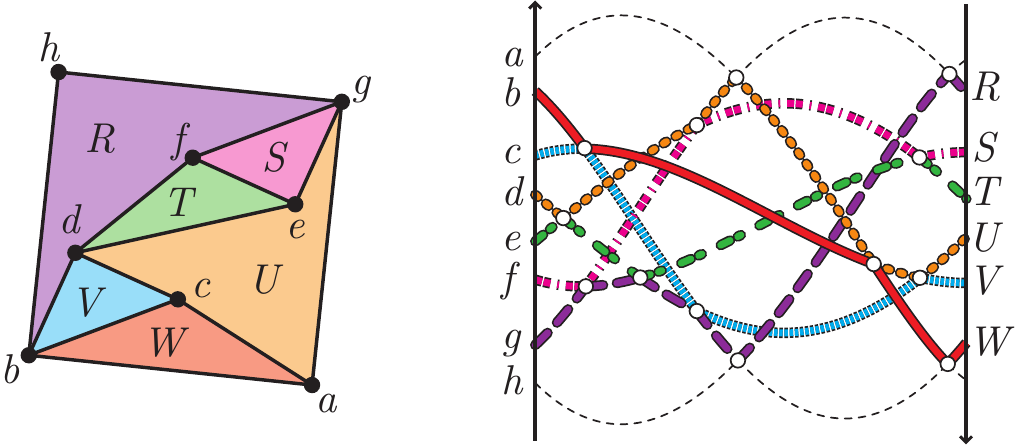}}
	\caption[Arrangement dual d'une pseudotriangulation]{Une pseudotriangulation et son arrangement dual.}
	\label{F:fig:dualpseudo}
\end{figure}

Nous donnons différentes preuves de la réciproque~:

\begin{theorem_F}
Soit~$P$ un ensemble de points en position générale et~$P^*$ son arrangement dual. Tout arrangement de pseudodroites avec points de contact dont le support couvre précisément le support de l'arrangement~$P^*$ privé de son premier niveau est l'arrangement dual d'une pseudotriangulation de~$P$.\qed
\end{theorem_F}

Motivés par ces deux théorèmes, nous considérons alors les arrangements de pseudodroites qui partagent un même support. Nous définissons une opération de flip qui correspond au flip dans les multitriangulations et les pseudotriangulations~: 

\begin{definition_F}
Deux arrangements de pseudodroites avec le même support sont reliés par un flip si la différence symétrique de leurs ensembles de points de contact est réduite à une paire~$\{u,v\}$. Dans ce cas, dans l'un des arrangements, $u$ est un point de contact de deux pseudodroites qui se croisent en~$v$, tandis que dans l'autre arrangement, $v$ est un point de contact de deux pseudodroites qui se croisent en~$u$.
\end{definition_F}

Nous étudions le graphe~$G(\cS)$ des flips sur les arrangements supportés par~$\cS$. Par exemple, lorsque~$\cS$ est le support d'un arrangement de deux pseudodroites avec~$p$ points de contact, le graphe~$G(\cS)$ est le graphe complet à~$p+1$ sommets. Nous nous intéressons à certaines orien\-tations acycliques du graphe~$G(\cS)$, données par des coupes verticales du support~$\cS$. Une telle orientation a une unique source que nous appelons \defn{arrangement glouton} et que nous caractérisons en terme de réseaux de tri. L'étude de ces arrangements gloutons et de leur transformation lorsque l'orientation de~$G(\cS)$ évolue fournit un algorithme d'énumération des arrangements de pseudodroites avec points de contact partageant un même support, dont l'espace de travail reste polynomial. Ainsi, nous éclairons et donnons une preuve complémentaire de l'algorithme similaire existant pour énumérer les pseudotriangulations d'un ensemble de points~\cite{bkps-ceppgfa-06}.

\svs
Nous revenons ensuite dans un cadre plus particulier. Au regard de la dualité entre les multitriangulations (resp.~les pseudotriangulations) et les arrangements de pseudodroites avec points de contact, nous proposons la généralisation suivante des multitriangulations lorsque les points ne sont pas en position convexe~:

\begin{definition_F}
Une \defn{\pt{k}} d'un arrangement de pseudodroites~$L$ est un arrangement de pseudodroites avec points de contact dont le support est l'arrangement~$L$ privé de ses~$k$ premiers niveaux. Une \pt{k} d'un ensemble de points~$P$ en position générale est un ensemble d'arêtes~$T$ qui correspond par dualité aux points de contact d'une \pt{k}~$T^*$ de~$P^*$.
\end{definition_F}

Nous montrons que toutes les \pt{k}s d'un ensemble de points~$P$ ont exactement ${k(2|P|-2k-1)}$ arêtes et ne peuvent pas contenir de configuration de~$2k+1$ arêtes alternantes, mais qu'elles peuvent en revanche éventuellement contenir un \kcrois{(k+1)}. Nous étudions ensuite leurs étoiles~: une étoile d'une \pt{k}~$T$ de~$P$ est un polygone formé par l'ensemble des arêtes correspondantes aux points de contact sur une pseudodroite fixée de~$T^*$. Nous discutons leur nombre possible de coins, et nous montrons que pour tout point~$q$ du plan, la somme des indices des étoiles de~$T$ autour de~$q$ est indépendante de~$T$.

\svs
Nous terminons le Chapitre~\ref{chap:mpt} avec trois questions reliées aux \mbox{\mpt{}s~:}
\begin{enumerate}[(i)]
\item Nous étudions d'abord les \defn{\mpt{}s itérées}~: une \pt{k} d'une \pt{m} d'un arrangement de pseudodroites~$L$ est une \pt{(k+m)} de~$L$. Nous donnons cependant un exemple de \ktri{2} qui ne contient pas de triangulation. Nous montrons en revanche que les \mpt{}s gloutonnes d'un arrangement sont des itérées de la pseudotriangulation gloutonne.
\item Nous donnons ensuite une caractérisation des arêtes de la \pt{k} gloutonne d'un ensemble de points~$P$ en terme de \defn{$k$-arbres d'horizon} de~$P$. Cette caractérisation généralise une observation de Michel Pocchiola~\cite{p-htvpt-97} pour les \mbox{pseudotriangulations}.
\item Finalement, nous définissons les \mpt{}s de configurations de convexes disjoints du plan et étudions leurs propriétés élémentaires. Les arrangements duaux des configurations de convexes sont les \defn{arrangements de double pseudodroites}, introduits par Luc Habert et Michel Pocchiola~\cite{hp-adp-08}. Dans l'Appendice~\ref{app:implementations}, nous énumérons à isomorphisme près les arrangements d'au plus~$5$ double pseudodroites.
\end{enumerate}


\subsection*{Trois problèmes ouverts}

Finalement, nous discutons au Chapitre~\ref{chap:multiassociahedron} trois problèmes ouverts qui illustrent la richesse combinatoire et géométrique des multitriangulations. Notre but est de présenter des idées naturelles basées sur les étoiles qui peuvent être fertiles bien qu'elles n'apportent pour l'heure que des réponses partielles à ces problèmes.

\svs
Le premier problème que nous discutons est celui de trouver une \defn{bijection explicite} entre l'ensemble des multitriangulations et l'ensemble des $k$-uplets de chemins de Dyck sans croisement, qui sont tous deux comptés par le déterminant de Hankel du Théorème~\ref{F:theo:enumeration}. Si dans une triangulation~$T$, on note~$\delta_i(T)$ le nombre de triangles de~$T$ dont le premier sommet est~$i$, alors l'application~$T\mapsto N^{\delta_0(T)}EN^{\delta_1(T)}E\dots N^{\delta_{n-3}(T)}E$ est une bijection des triangulations du \gon{n}e dans les chemins de Dyck de demi-longueur~$n-2$ (où~$N$ et~$E$ désignent les pas nord et est). Jakob Jonsson~\cite{j-gt-03} a généralisé cette remarque en comparant la répartition des séquences de degrés entrants des \ktri{k}s avec les signatures des $k$-uplets de chemins de Dyck (sans donner les définitions précises, signalons que ces deux $k$-uplets de séquences généralisent la séquence~$(\delta_i(T))$ d'une triangulation~$T$ et la séquence des puissances de~$N$ dans un chemin de Dyck). Motivés par ce résultat, nous avons cherché à définir des $k$-colorations des cordes \kpert{k}s d'une \ktri{k} de sorte que les séquences respectives de degrés entrant de chaque couleur définissent des $k$-uplets de chemins de Dyck qui ne se croisent pas. Nous présentons une coloration basée sur les étoiles, qui vérifie cette propriété, mais pour laquelle l'application des multitriangulations dans les $k$-uplets de chemins de Dyck n'est malheureusement pas bijective.

\svs
Notre second problème concerne les propriétés de \defn{rigidité} des multitriangulations. Une trian\-gulation est \defn{minimalement rigide} dans le plan~: les seuls mouvements de ses sommets qui préservent les longueurs de ses arêtes sont des isométries du plan, et retirer n'importe quelle arête rend la structure flexible. De manière équivalente, une triangulation vérifie la condition de Laman~: elle a~$2n-3$ arêtes et tout sous-graphe à~$m$ sommets a au plus~$2m-3$ arêtes. Nous observons deux connexions intéressantes entre les multitriangulations et la théorie de la rigidité~:
\begin{enumerate}[(i)]
\item Tout d'abord, une \ktri{k} est $\left(2k,{2k+1 \choose 2}\right)$-\defn{raide}~: elle contient~$2kn-{2k+1 \choose 2}$ arêtes et tout sous-graphe à~$m$ sommets a au plus~$2km-{2k+1 \choose 2}$ arêtes. Cette propriété combinatoire fait d'une \ktri{k} un candidat raisonnable pour être un graphe génériquement minimalement rigide en dimension~$2k$. Nous prouvons cette conjecture lorsque~$k=2$.
\item Ensuite, nous montrons que le graphe dual d'une $k$-triangulation est~$(k,k)$-raide. En particulier, il peut être décomposé en~$k$ arbres couvrants arêtes-disjoints, ce qui doit être rapproché du fait que le dual d'une triangulation est un arbre.
\end{enumerate}

\svs
Finalement, nous revenons à la question de la \defn{réalisation polytopale} du complexe simplicial~$\Delta_{n,k}$ formé par les ensembles de cordes \kpert{k}s du \gon{n}e sans \kcrois{(k+1)}. Nous présentons deux contributions modestes à cette question~:
\begin{enumerate}[(i)]
\item Nous répondons d'une part au premier cas non-trivial en décrivant l'espace des réalisations symétriques de~$\Delta_{8,2}$. Ce résultat est obtenu en deux étapes~: d'abord nous énumérons par ordinateur tous les matroides orientés symétriques qui réalisent notre complexe simplicial~; ensuite nous étudions les réalisations polytopales symétriques de ces matroides orientés.
\item Nous nous intéressons d'autre part à une construction de l'associaèdre due à Jean-Louis \mbox{Loday}~\cite{l-rsp-04}. Nous proposons une interprétation de sa construction en termes d'arrange\-ments duaux des triangulations, qui se généralise naturellement aux multitriangulations. Nous obtenons un polytope dont les facettes sont définies par des inégalités simples et qui réalise le graphe des flips restreint à certaines multitriangulations (celles dont le graphe dual est acyclique). Ce polytope aurait pu être \apriori une projection d'un polytope réalisant~$\Delta_{n,k}$, mais nous montrons toutefois qu'une telle projection est impossible.
\end{enumerate}


\section{Polytopalité de produits}

Dans la deuxième partie de cette thèse, nous nous intéressons à des questions de réalisation polytopale de graphes (ou de squelettes) obtenus par produits cartésiens d'autres graphes (ou squelettes).

Le \defn{produit cartésien} de deux polytopes~$P,Q$ est le polytope~$P\times Q \eqdef \ens{(p,q)}{p\in P,q\in Q}$. Sa structure combinatoire ne dépend que de celle de ses facteurs~: la dimension de~$P\times Q$ est la somme des dimensions de~$P$ et~$Q$ et les faces non-vides de~$P\times Q$ sont précisément les produits d'une face non-vide de~$P$ par une face non-vide de~$Q$.


\subsection*{Polytopalité de produits de graphes non-polytopaux}

Dans le Chapitre~\ref{chap:nonpolytopal}, nous nous intéressons au produit de graphes. Le \defn{produit cartésien} de deux graphes~$G$ et~$H$ est le graphe~$G\times H$ dont les sommets sont~$V(G\times H) \eqdef V(G)\times V(H)$ et dont les arêtes sont~$E(G\times H) \eqdef \big(E(G)\times V(H)\big)\cup\big(V(G)\times E(H)\big)$. Autrement dit, pour tous sommets $a,c\in V(G)$ et~$b,d\in V(H)$, les sommets~$(a,b)$ et~$(c,d)$ de~$G\times H$ sont adjacents si~$a=c$ et $\{b,d\}\in E(H)$, ou~$b=d$ et~$\{a,c\}\in E(G)$. Ce produit est en accord avec le produit de polytopes~: le graphe d'un produit de polytopes est le produit de leurs graphes. En particulier, le produit de deux graphes polytopaux est automatiquement polytopal. Dans ce chapitre, nous étudions la question réciproque~: étant donnés deux graphes~$G$ et~$H$, la polytopalité du produit~$G\times H$ implique-t-elle celle de ses facteurs~$G$ et~$H$~?

Le produit d'un triangle par un chemin de longueur~$2$ est un contre-exemple simple à cette question~: bien que le chemin ne soit pas polytopal, le produit est le graphe d'un \poly{3}tope obtenu en collant deux prismes triangulaires par une face triangulaire. Nous éliminons de tels exemples en exigeant que chaque facteur soit un graphe régulier. Si~$G$ et~$H$ sont réguliers de degré respectifs~$d$ et~$e$, alors le produit~$G\times H$ est $(d+e)$-régulier et il est naturel de se demander s'il est le graphe d'un \poly{(d+e)}tope simple. La réponse est donnée par le théorème suivant~:

\begin{theorem_F}\label{F:theo:simpleproduct}
Un produit~$G\times H$ est le graphe d'un polytope simple si et seulement si ses deux facteurs~$G$ et~$H$ sont des graphes de polytopes simples. Dans ce cas, il existe un unique polytope simple dont le graphe est~$G\times H$~: c'est précisément le produit des uniques polytopes simples dont les graphes respectifs sont~$G$ et~$H$.\qed
\end{theorem_F}

Dans ce théorème, l'unicité du polytope simple réalisant~$G\times H$ découle du fait qu'un polytope simple est complètement déterminé par son graphe~\cite{bm-ppi-87,k-swtsp-88}. Ces résultats reposent sur la propriété suivante~: tout ensemble de~$k+1$ arêtes adjacentes à un même sommet d'un polytope simple~$P$ définit une \face{k} de~$P$.

En application du Théorème~\ref{F:theo:simpleproduct} nous obtenons une famille infinie de graphes $4$-réguliers non-polytopaux~: le produit d'un graphe $3$-régulier non-polytopal par un segment est non-polyto\-pal et $4$-régulier.

\begin{figure}[h]
	\centerline{\includegraphics[width=.9\textwidth]{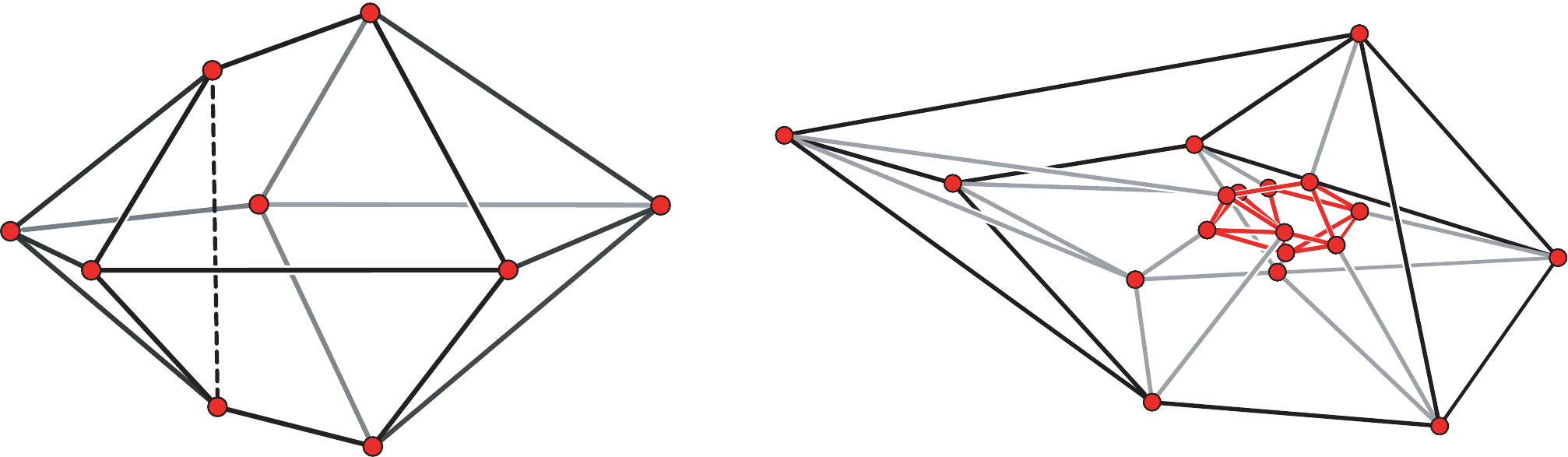}}
	\caption[Un graphe~$H$ $4$-régulier et non-polytopal et le diagramme de Schlegel d'un \poly{4}tope dont le graphe est le produit de~$H$ par un segment]{Un graphe~$H$ $4$-régulier et non-polytopal qui est le graphe d'une subdivision régulière d'un \poly{3}tope (gauche) et le diagramme de Schlegel d'un \poly{4}tope dont le graphe est le produit de~$H$ par un segment (droite).}
	\label{F:fig:truncatedoctahedron}
\end{figure}

Nous cherchons ensuite à savoir si un produit de graphes réguliers non-polytopaux peut être polytopal dans une dimension plus petite que son degré. Les exemples suivants répondent à cette question~:

\begin{theorem_F}
\begin{enumerate}[(i)]
\item Lorsque~$n\ge 3$, le produit~$K_{n,n}\times K_2$ d'un graphe complet bipartite par un segment n'est pas polytopal.
\item Le produit d'un graphe \poly{d}topal par le graphe d'une subdivision régulière d'un \poly{e}\-tope est \poly{(d+e)}topal. Ceci fournit des produits polytopaux de graphes réguliers non-polytopaux (par exemple le produit de deux dominos de la \fref{F:fig:products}et le graphe de la \fref{F:fig:truncatedoctahedron}\hspace{-.1cm}).\qed
\end{enumerate}
\end{theorem_F}

\subsection*{Polytopes \emph{neighborly} prodsimpliciaux}

Au Chapitre~\ref{chap:psn}, nous considérons des polytopes dont les squelettes sont ceux d'un produit de simplexes~:

\begin{definition_F}
Soient~$k\ge0$ et~$\sub{n} \eqdef (n_1,\dots,n_r)$, avec~$r\ge1$ et~$n_i\ge 1$ pour tout~$i$. Un polytope est \defn{$(k,\sub{n})$-neighborly prodsimplicial}~---~ou en abrégé \defn{\xpsn{(k,\sub{n})}}~---~si son \squelette{k} est combinatoirement équivalent à celui du produit de simplexes~$\simplex_{\sub{n}} \eqdef \simplex_{n_1}\times\cdots\times\simplex_{n_r}$.
\end{definition_F}

Cette définition généralise deux classes particulières de polytopes~:
\begin{enumerate}[(i)]
\item On retrouve les polytopes \defn{neighborly} lorsque~$r=1$. Dans la littérature, un polytope est \defn{\neighborly{k}} si tout sous-ensemble de~$k$ de ses sommets forme une face. Avec notre definition, un tel polytope est \xpsn{(k-1,n)}.
\item On retrouve les polytopes \defn{neighborly cubiques}~\cite{jz-ncp-00,js-ncps-07,sz-capdp} lorsque~$\sub{n}=(1,1,\dots,1)$.
\end{enumerate}

\svs
Le produit~$\simplex_{\sub{n}}$ est un polytope \xpsn{(k,\sub{n})} de dimension~$\sum n_i$. Nous cherchons naturellement des polytopes \xpsn{(k,\sub{n})} en dimension inférieure. Par exemple, le polytope cyclique~$C_{2k+2}(n+1)$ est un polytope \xpsn{(k,n)} en dimension~$2k+2$. Nous notons~$\delta(k,\sub{n})$ la dimension minimale que peut avoir un polytope \xpsn{(k,\sub{n})}.

Certains polytopes \psn s'obtiennent en projetant le produit~$\simplex_{\sub{n}}$, ou tout polytope combinatoirement équivalent, sur un sous-espace de dimension inférieure. Par exemple, le polytope cyclique~$C_{2k+2}(n+1)$, comme n'importe quel polytope avec~$n+1$ sommets, est une projection du simplexe~$\simplex_n$ sur~$\R^{2k+2}$.

\begin{definition_F}
Un polytope \xpsn{(k,\sub{n})} est \defn{$(k,\sub{n})$-neighborly prodsimplicial projeté} ---~ou en abrégé \defn{\xppsn{(k,\sub{n})}}~---~si c'est une projection d'un polytope combinatoirement équivalent à~$\simplex_{\sub{n}}$.
\end{definition_F}

Nous notons~$\delta_{pr}(k,\sub{n})$ la dimension minimale d'un polytope \xppsn{(k,\sub{n})}.

\mvs
Notre Chapitre~\ref{chap:psn} se divise en deux parties. Dans la première, nous présentons trois méthodes pour construire des polytopes \ppsn en petite dimension~:
\begin{enumerate}[(i)]
\item par des réflexions de polytopes cycliques~;
\item par des sommes de Minkowski de polytopes cycliques~;
\item par des ``projections de produits déformés'', dans l'esprit des constructions de Raman Sanyal et G\"unter Ziegler~\cite{z-ppp-04,sz-capdp}.
\end{enumerate}
Dans la seconde partie, nous obtenons des obstructions topologiques à l'existence de tels objets, en utilisant des techniques développées par Raman Sanyal~\cite{s-tovnms-09} pour borner le nombre de sommets d'une somme de Minkowski. Au regard de ces obstructions, nos constructions de la première partie se révèlent optimales pour un large spectre de paramètres.

\paragraph{Constructions.}
Notre premier exemple non-trivial est un polytope \xpsn{(k,(1,n))} en dimension~$2k+2$, obtenu par réflexion du polytope cyclique~$C_{2k+2}(n+1)$ par rapport à un hyperplan bien choisi~:

\begin{proposition_F}
Pour tous~$k\ge0$,~$n\ge2k+2$ et~$\lambda\in\R$ suffisamment grand, le polytope
$$\conv\left(\ens{(t_i,\dots,t_i^{2k+2})^T}{i\in[n+1]} \cup \ens{(t_i,\dots,t_i^{2k+1},\lambda-t_i^{2k+2})^T}{i\in[n+1]}\right)$$
est un polytope \xpsn{(k,(1,n))} de dimension~$2k+2$.\qed
\end{proposition_F}

Par exemple, cette construction fournit un \poly{4}tope dont le graphe est~$K_2\times K_n$~($n\ge3$).

\svs
Ensuite, à l'aide de sommes de Minkowski bien choisies de polytopes cycliques, nous obtenons des coordonnées explicites de polytopes \xppsn{(k,\sub{n})}~:

\begin{theorem_F}\label{F:theo:UBminkowskiCyclic}
Soient~$k\ge0$ et~$\sub{n} \eqdef (n_1,\dots,n_r)$ avec~$r\ge1$ et~$n_i\ge1$ pour tout~$i$. Il existe des ensembles~$I_1,\dots,I_r\subset\R$, avec~$|I_i|=n_i$ pour tout~$i$, tels que le polytope
$$\conv\ens{w_{a_1,\dots,a_r}}{(a_1,\dots,a_r)\in I_1\times\cdots\times I_r} \subset \R^{2k+r+1}$$
soit \xppsn{(k,\sub{n})}, où $w_{a_1,\dots,a_r} \eqdef \big(a_1,\dots,a_r,\sum_{i\in[r]}a_i^2,\dots,\sum_{i\in[r]} a_i^{2k+2}\big)^T$. Par conséquent,
$$\delta(k,\sub{n}) \le \delta_{pr}(k,\sub{n}) \le 2k+r+1.$$
\end{theorem_F}
\vspace{-.75cm}\qed
\vspace{.35cm}

Lorsque~$r=1$ nous retrouvons les polytopes \emph{neighborly}. 

\svs
Finalement, nous étendons la technique de ``projection de produits déformés de polygones'' de Raman Sanyal et G\"unter Ziegler~\cite{z-ppp-04,sz-capdp} aux produits de polytopes simples arbitraires~: nous projetons un polytope bien choisi, combinatoirement équivalent à un produit de polytopes simples, de sorte à préserver son \squelette{k} complet. Plus concrètement, nous décrivons comment utiliser des colorations des graphes des polytopes polaires des facteurs du produit pour augmenter la dimension du squelette préservé. La version basique de cette technique fournit le résultat suivant~:

\begin{proposition_F}
Soient~$P_1,\dots,P_r$ des polytopes simples. Pour chaque polytope~$P_i$, on note~$n_i$ sa dimension,~$m_i$ son nombre de facettes, et~$\chi_i \eqdef \chi(\gr(P_i^\polar))$ le nombre chromatique du graphe du polytope polaire~$P_i^\polar$. Pour un entier fixé~$d\le n$, soit~$t$ l'entier maximal tel que~${\sum_{i=1}^t n_i\le d}$. Alors il existe un \poly{d}tope dont le \squelette{k} est combinatoirement équivalent à celui du produit ${P_1\times\cdots\times P_r}$ dès lors que
$$0 \le k \le \sum_{i=1}^r (n_i-m_i) + \sum_{i=1}^t (m_i-\chi_i) + \Floor{\frac{1}{2}\left(d-1+\sum_{i=1}^t(\chi_i-n_i)\right)}.$$
\end{proposition_F}
\vspace{-1.3cm}\qed
\vspace{1.1cm}

\svs
En spécialisant cette proposition à un produit de simplexes, nous obtenons une autre construction de polytopes \ppsn. Lorsque certains des simplexes sont petits comparés à~$k$, cette technique produit en fait nos meilleurs exemples de polytopes \ppsn~:

\begin{theorem_F}
Pour tous~$k\ge0$ et~$\sub{n} \eqdef (n_1,\dots,n_r)$ avec~${1=n_1=\cdots=n_s<n_{s+1}\le\cdots\le n_r}$,
$$\delta_{pr}(k,\sub{n}) \le
\begin{cases}
     2(k+r)-s-t & \text{si } 3s \le 2k+2r, \\
     2(k+r-s)+1 & \text{si } 3s = 2k+2r+1, \\
     2(k+r-s+1) & \text{si } 3s \ge 2k+2r+2,
\end{cases}$$
où~$t\in\{s,\dots,r\}$ est maximal tel que~$3s+\sum_{i=s+1}^{t}(n_i+1) \le 2k+2r$.\qed
\end{theorem_F}

Si~$n_i=1$ pour tout~$i$, nous retrouvons les polytopes \emph{neighborly} cubiques de~\cite{sz-capdp}.

\paragraph{Obstructions.}
Pour obtenir des bornes inférieures sur la dimension minimale~$\delta_{pr}(k,\sub{n})$ que peut avoir un polytope \xppsn{(k,\sub{n})}, nous appliquons une méthode due à Raman Sanyal~\cite{s-tovnms-09}. À toute projection qui préserve le \squelette{k} de~$\simplex_{\sub{n}}$, nous associons par dualité de Gale un complexe simplicial qui doit être plongeable dans un espace d'une certaine dimension. L'argument repose ensuite sur une obstruction topologique dérivée du critère de Sarkaria pour le plongement d'un complexe simplicial en termes de colorations de graphes de Kneser~\cite{m-ubut-03}. Nous obtenons le résultat suivant~:

\begin{theorem_F}\label{F:theo:topObstr}
Soit~$\sub{n} \eqdef (n_1,\dots,n_r)$ avec~$1=n_1=\cdots=n_s<n_{s+1}\le\cdots\le
n_r$.
\begin{enumerate}[(i)]
\item Si
$$0 \le k \le \sum_{i=s+1}^r \Fracfloor{n_i-2}{2} + \max\left\{0,\Fracfloor{s-1}{2}\right\},$$
alors~$\delta_{pr}(k,\sub{n}) \ge 2k+r-s+1$.
\item Si~$k\ge \Floor{\frac{1}{2} \sum_i n_i}$ alors~$\delta_{pr}(k,\sub{n}) \ge \sum_i n_i$.\qed
\end{enumerate}
\end{theorem_F}

En particulier, les bornes inférieure et supérieure des Théorèmes~\ref{F:theo:UBminkowskiCyclic} et~\ref{F:theo:topObstr} se rejoignent sur un large champ de paramètres~:

\begin{theorem_F}\label{F:theo:mainResult}
Pour tout~$\sub{n} \eqdef (n_1,\dots,n_r)$ avec~$r\ge1$ et~$n_i\ge2$ pour tout~$i$, et pour tout~$k$ tel que~$0\le k\le \sum_{i\in [r]}\Fracfloor{n_i-2}{2}$, le plus petit polytope \xppsn{(k,\sub{n})} est de dimension exactement $2k+r+1$. Autrement dit~:
$$\delta_{pr}(k,\sub{n}) = 2k+r+1.$$
\end{theorem_F}
\vspace{-.8cm}\qed

\begin{remark_F}
Les techniques de projection de polytopes et les obstructions que nous utilisons ont été développées par Raman Sanyal et G\"unter Ziegler~\cite{z-ppp-04,sz-capdp,s-tovnms-09}. Nous avons décidé de les présenter dans cette thèse parce que leur application aux produits de simplexes fournit des résultats nouveaux qui complètent notre étude sur la polytopalité des produits. Par ailleurs, après avoir appliqué ces méthodes aux produits de simplexes, nous avons découvert que Thilo~R\"orig et Raman~Sanyal avaient un travail similaire en cours sur le même sujet~\cite{rs-npps} (voir aussi~\cite{sanyal-phd,rorig-phd}).
\end{remark_F}


\cleardoublepage

\setcounter{chapter}{-2}
\setcounter{section}{0}
\selectlanguage{spanish}
\renewcommand*{\figurename}{Figura}
\makeatletter\renewcommand{\thefigure}{\@arabic\c@figure}\makeatother
\renewcommand{\namepart}{Resumen}

\theoremstyle{theorem} 
\newtheorem{theorem_E}{Teorema}
\newtheorem{proposition_E}[theorem_E]{Proposición}
\newtheorem{definition_E}[theorem_E]{Definición}

\theoremstyle{definition} 
\newtheorem{example_E}[theorem_E]{Ejemplo}
\newtheorem{remark_E}[theorem_E]{Comentario}
\newtheorem{observación_E}[theorem_E]{Observación}

\newcommand{\ktric}[1]{\mbox{$#1$-trian}\-gu\-la\-ción}
\newcommand{\ptc}[1]{\mbox{$#1$-pseu}\-do\-trian\-gu\-la\-ción}
\newcommand{\ptcs}[1]{\mbox{$#1$-pseu}\-do\-trian\-gu\-la\-ciones}
\newcommand{\mptcs}{multipseudotriangulaciones}
\newcommand{\ktrics}[1]{\mbox{$#1$-trian}\-gu\-la\-ciones}
\newcommand{\kcruce}[1]{\mbox{$#1$-cruce}}
\newcommand{\kborde}[1]{\mbox{$#1$-borde}}
\newcommand{\kestrella}[1]{\mbox{$#1$-estrella}}
\newcommand{\esqueleto}[1]{\mbox{$#1$-esqueleto}}
\newcommand{\poli}[1]{\mbox{$#1$-poli}}

\chapter*{Resumen}\label{chap:resumeE}
\phantomsection
\addcontentsline{toc}{chapter}{Resumen}



\section{Introducción}

Los temas tratados en esta memoria se insertan en el campo de la \defn{geometría discreta y algorítmica}~\cite{hdcg-04}: los problemas encontrados en general involucran conjuntos finitos de objetos geométricos elementales (como puntos, rectas, semi-espacios, \etc), y las preguntas tratan de la manera como están relacionados, colocados, situados unos respecto a otros (cómo se intersecan, cómo se ven, qué delimitan, \etc). Esta memoria versa sobre dos sujetos particulares: las \defn{\mbox{multitriangulaciones}} y las \defn{realizaciones politopales de productos}. Sus varias conexiones con la geometría discreta nos permitirán descubrir algunas de sus numerosas facetas. Presentamos algunos ejemplos en esta introducción. Empezamos por presentar su problemática común, la búsqueda de una realización politopal de una estructura, que ha dirigido nuestra investigación sobre estos sujetos.

\svs
Un politopo (convexo) es la envoltura convexa de un conjunto finito de puntos de un espacio euclidiano. Aunque el interés despertado por ciertos politopos se remonta a la Antigüedad (sólidos Platónicos), su estudio sistemático es relativamente reciente y los principales resultados datan del siglo \textsc{xx} (ver~\cite{g-cp-03,z-lp-95} y sus referencias).
El estudio de los politopos versa no solamente sobre sus propiedades geométricas pero también sobre sus aspectos más combinatorios. Se trata especialmente de comprender sus caras (sus intersecciones con un hiperplano soporte) y el retículo que componen (es decir las relaciones de inclusión entre estas caras).

Las preguntas de \defn{realización politopal} forman en cierto modo el problema recíproco: versan sobre la existencia y la construcción de politopos a partir de una estructura combinatoria dada. Por ejemplo, dado un grafo, querríamos determinar si es el grafo de un politopo: diremos entonces que el grafo es politopal. El resultado fundador en este campo es el \mbox{Teorema} de Steinitz~\cite{s-pr-22} que caracteriza los grafos de \poli{3}topos. En cuanto pasamos a dimensión~$4$, la situación es mucho menos satisfactoria: a pesar de ciertas condiciones necesarias~\cite{b-gscps-61,k-ppg-64,b-ncp-67}, los grafos de politopos no admiten ninguna caracterización local en dimensión general~\cite{rg-rsp-96}. Cuando un grafo es politopal, nos preguntamos por las propiedades de sus realizaciones, por ejemplo su número de caras, su dimensión, \etc{} Intentamos construir ejemplos que optimizan ciertas de estas propiedades: típicamente, un politopo de dimensión mínima con un grafo dado~\cite{g-ncp-63,jz-ncp-00,sz-capdp}. Estas preguntas de realizaciones politopales son interesantes para grafos que provienen bien de grafos de transformación sobre conjuntos combinatorios o geométricos (grafo de transposiciones adyacentes sobre las permutaciones, grafo de flips sobre las triangulaciones, \etc), bien de operaciones sobre otros grafos (que pueden ser locales como la transformación~$\Delta Y$, o globales como el producto cartesiano). Estas preguntas son interesantes no solamente para un grafo, sino más generalmente para cualquier subconjunto de un retículo.


\paragraph{Politopalidad de grafos de flips.}
La existencia de realizaciones politopales es primero estudiada para grafos de transformación sobre estructuras combinatorias o geométricas. Podemos mencionar aquí el permutoedro cuyos vértices corresponden a las permutaciones de~$[n]$ y donde dos vértices están relacionados por una arista si las permutaciones correspondientes difieren por una transposición de dos posiciones adyacentes. Otros ejemplos, así como clases de politopos que permitan realizar ciertas estructuras combinatorias están expuestos en~\cite[Lectura~9]{z-lp-95}. De manera general, las preguntas de politopalidad de estructuras combinatorias son interesantes no solamente por sus resultados, sino también porque su estudio obliga a entender la combinatoria de los objetos y a desarrollar nuevas métodos.

Dos ejemplos particulares de estructuras combinatorias politopales desempeñan un papel importante en esta memoria. Encontramos primero el asociaedro cuyo borde realiza el dual del complejo simplicial formado por todos los conjuntos sin cruce de cuerdas del \gon{n}o. El grafo del asociaedro corresponde al grafo de flips de las triangulaciones del \gon{n}o. El asociaedro aparece en diversos contextos y varias realizaciones politopales han sido propuestas~\cite{l-atg-89,bfs-ccsp-90,gkz-drmd-94,l-rsp-04,hl-rac-07}. Encontramos después el politopo de las pseudotriangulaciones de un conjunto de puntos del plano euclideano~\cite{rss-empppt-03}. Introducidas para el estudio del complejo de visibilidad de obstáculos convexos disjuntos del plano~\cite{pv-tsvcp-96,pv-vc-96}, las pseudotriangulaciones han sido utilizadas en diversos contextos geométricos~\cite{rss-pt-06}. Sus propiedades de rigidez~\cite{s-ptrmp-05} permitían establecer la politopalidad de su grafo de flips~\cite{rss-empppt-03}.

\begin{figure}
	\capstart
	\centerline{\includegraphics[width=\textwidth]{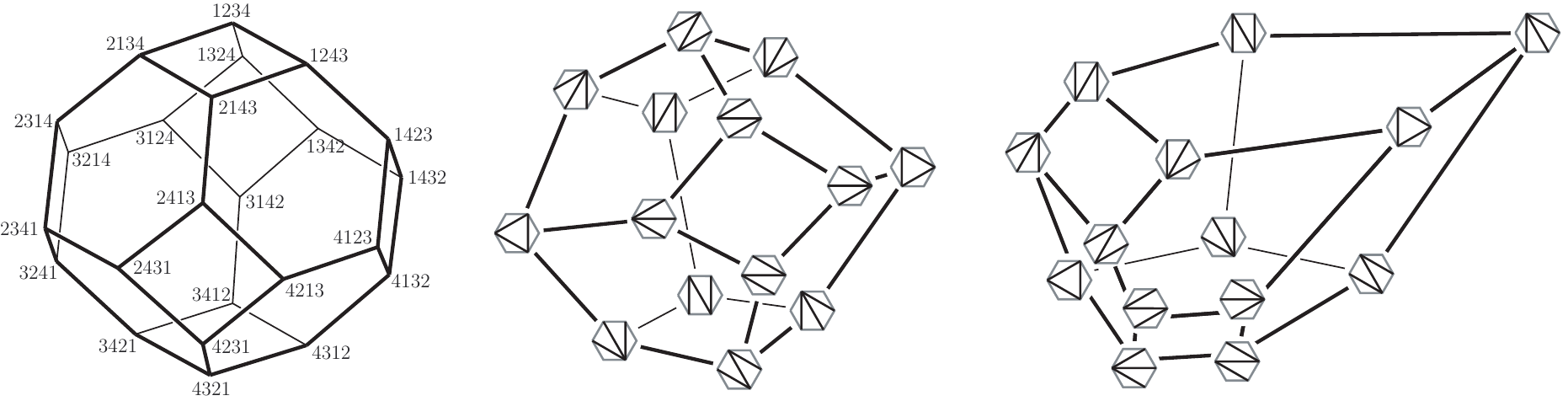}}
	\caption[El permutoedro y dos realizaciones del asociaedro]{El permutoedro y dos realizaciones del asociaedro.}
	\label{E:fig:permutohedronassociahedra}
\end{figure}

\svs
En la primera parte de esta memoria, nos interesamos por la politopalidad del grafo de flips de las multitriangulaciones. Estos objetos, aparecidos de manera relativamente contingente~\cite{cp-tttccp-92,dkm-lahp-02,dkkm-2kn..-01}, tienen una rica estructura combinatoria~\cite{n-gdfcp-00,j-gt-03,j-gtdfssp-05}. Una \defn{\ktric{k}} es un conjunto maximal de cuerdas del \gon{n}o que no contiene ningún subconjunto de $k+1$ cuerdas que se cruzan mutuamente. Consideramos el grafo de flips en el cual dos multitriangulaciones están relacionadas si difieren por una cuerda. Como para las triangulaciones, que ocurren cuando~$k=1$, este grafo es regular y conexo, y nos preguntamos por su politopalidad. Jakob Jonsson~\cite{j-gt-03} hizo un primer paso en esta dirección demostrando que el complejo simplicial formado por los conjuntos de cuerdas que no contienen ningún subconjunto de $k+1$ cuerdas que se cruzan mutuamente es una esfera topológica. Aunque solo tenemos respuestas parciales, esta pregunta reveló resultados interesantes que exponemos en esta memoria.

\svs
Algunas construcciones del asociaedro~\cite{bfs-ccsp-90,l-rsp-04} están basadas, directamente o indirectamente, en los triángulos de las triangulaciones. Para las multitriangulaciones, ningún objeto elemental similar aparece en los trabajos anteriores. En primer lugar, hemos intentado entender qué ha sido de los triángulos en las multitriangulaciones. Las \defn{estrellas} que introducimos en el Capítulo~\ref{chap:stars} responden a esta pregunta. Tal como los triángulos en las triangulaciones, estimamos que estas estrellas dan el punto de vista pertinente para entender las multitriangulaciones. Como prueba, empezamos por recobrar con ayuda de estas estrellas todas las propiedades combinatorias elementales conocidas hasta entonces sobre las multitriangulaciones. Primero, estudiamos las relaciones de incidencia entre las estrellas y las cuerdas (cada cuerda interna está contenida en dos estrellas) que nos permiten recobrar que todas las \ktrics{k} del \gon{n}o tienen el mismo cardinal. Después, considerando las bisectrices de las estrellas, damos una interpretación local de la operación de flip (una cuerda interna es substituida por la única bisectriz común de las dos estrellas que la contienen), lo que aclara el estudio del grafo de flips y de su diámetro. Redefinimos además en términos de estrellas operaciones inductivas sobre multitriangulaciones que permiten añadir o suprimir un vértice al \gon{n}o. En fin, utilizamos la descomposición de una multitriangulación en estrellas para interpretarla como una descomposición poligonal de superficie, y aplicamos esta interpretación a la construcción de descomposiciones regulares de superficies.

\begin{figure}
	\capstart
	\centerline{\includegraphics[scale=1]{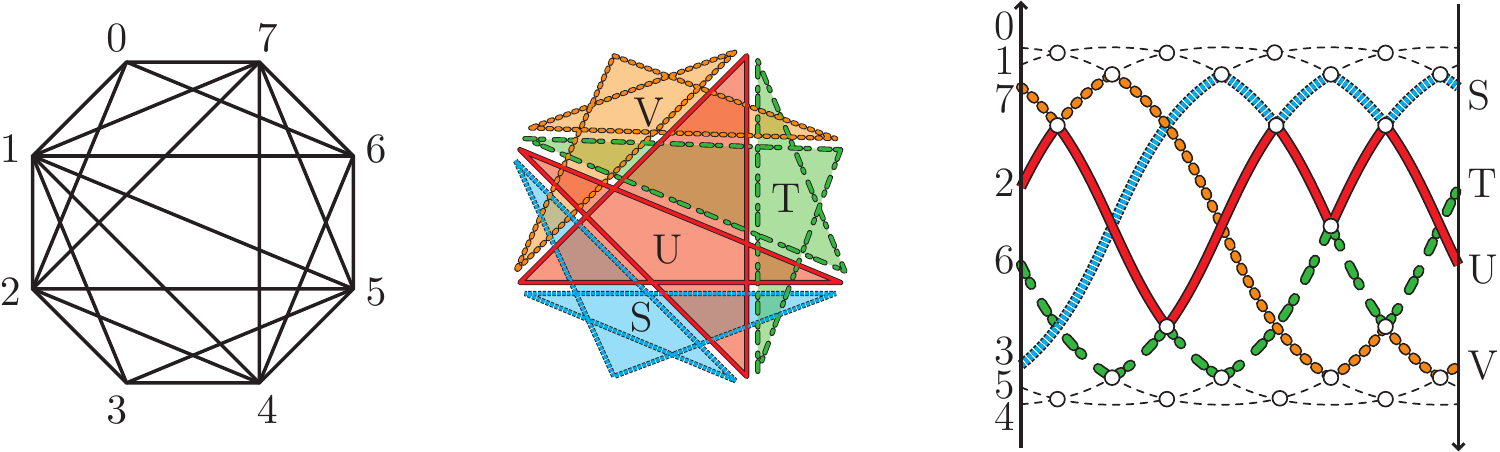}}
	\caption[Una \ktric{2}, su descomposición en estrellas, y su arreglo dual]{Una \ktric{2}, su descomposición en estrellas, y su arreglo dual.}
	\label{E:fig:ktristardual}
\end{figure}

\svs
En segundo lugar, hemos intentado entender las multitriangulaciones por dualidad. El espacio de las rectas del plano es una banda de M\"obius; el conjunto de las rectas que pasan por un punto del plano es una pseudorecta de la banda de M\"obius; y las pseudorectas duales de una configuración de puntos forman un arreglo de pseudorectas~\cite{g-pa-97}.

Nuestro punto de partida es la dualidad entre las pseudotriangulaciones de un conjunto~$P$ de puntos y los arreglos de pseudorectas cuyo soporte es el arreglo dual de~$P$ salvo su primer nivel. Establecemos una dualidad similar para multitriangulaciones. Por una parte, el conjunto de las bisectrices de una estrella es una pseudorecta de la banda de M\"obius. Por otra parte, las pseudorectas duales de las estrellas de una \ktric{k} del \gon{n}o forman un arreglo de pseudorectas con puntos de contacto, cuyo soporte es exactamente el soporte del arreglo dual de los vértices del \gon{n}o salvo sus $k$ primeros niveles. Mostramos que cualquier arreglo con puntos de contacto que tiene este soporte es efectivamente el arreglo dual de una multitriangulación. Esta dualidad relaciona las pseudotriangulaciones y las multitriangulaciones y explica así sus propiedades comunes (número de cuerdas, flip, \etc).

Más generalmente, estudiamos en el Capítulo~\ref{chap:mpt} los \defn{arreglos de pseudorectas con puntos de contacto} que comparten el mismo soporte. Definimos una operación de flip que corresponde al flip en las multitriangulaciones, y estudiamos el grafo de flips. Las propiedades de ciertos arreglos ``glotones'', definidos como las fuentes de ciertas orientaciones acíclicas de este grafo, nos permiten en particular enumerar este grafo de flips manteniendo un espacio polinomial. Nuestro trabajo aclara así el algoritmo existente para enumerar las pseudotriangulaciones~\cite{bkps-ceppgfa-06} y ofrece una prueba complementaria.

\svs
Por último, discutimos en el Capítulo~\ref{chap:multiassociahedron} tres problemas abiertos que reflejan la riqueza combinatoria y geométrica de las multitriangulaciones.

El primero concierne el conteo de las multitriangulaciones. Jakob Jonsson~\cite{j-gtdfssp-05} probó (considerando $0/1$-rellenos de poliominos que evitan ciertos patrones) que las multitriangulaciones están contadas por un cierto determinante de Hankel de números de Catalan, que cuentan también ciertas familias de $k$-tuplas de caminos de Dyck. Sin embargo, excepto dos resultados parciales~\cite{e-btdp-07,n-abtdp-09}, no existe ninguna prueba biyectiva de este resultado. Entramos aquí el dominio de la \defn{combinatoria biyectiva} cuyo objetivo es construir biyecciones entre familias combinatorias que conserven parámetros característicos: nos gustaría aquí una biyección que permitiera leer las estrellas sobre las $k$-tuplas de caminos de Dyck.

Nuestro segundo problema es el de la \defn{rigidez}. Aunque la rigidez de grafos es bien entendida en dimensión~$2$~\cite{l-grpss-70,g-cf-01,gss-cr-93}, ninguna caracterización satisfactoria es conocida a partir de dimensión~$3$. Mostramos que las \ktrics{k} satisfacen las propiedades típicas de los grafos rígidos de dimensión~$2k$. Esto incita naturalmente a conjeturar que las \ktrics{k} son rígidas en dimensión~$2k$, lo que demostramos cuando~$k=2$. Una respuesta positiva a esta conjetura permitiría acercarse de la politopalidad del grafos de flips de las multitriangulaciones, tal como el politopo de las pseudotriangulaciones~\cite{rss-empppt-03} ha sido construido apoyándose en sus propiedades de rigidez.

Por fin, volvemos sobre la \defn{realización politopal del grafo de flips} de las multitriangulaciones. En primer lugar, estudiamos el primer ejemplo no-trivial: mostramos que el grafo de flips sobre las \ktric{2}s del octógono es efectivamente el grafo de un politopo de dimensión~$6$. Para encontrar un tal politopo, describimos completamente el espacio de las realizaciones politopales simétricas de este grafo en dimensión~$6$, estudiando primero todos los matroides orientados~\cite{bvswz-om-99,b-com-06} simétricos que pueden realizarlo. En segundo lugar, generalizando la construcción del asociaedro de Jean-Louis Loday~\cite{l-rsp-04}, construimos un politopo que realiza el grafo de los flips restringido a las multitriangulaciones cuyo grafo dual es acíclico.

\svs
Además, presentamos en el Apéndice~\ref{app:implementations} los resultados de un algoritmo de enumeración de pequeños arreglos de pseudorectas y de \defn{doble pseudorectas}. A semejanza de los arreglos de pseudorectas que ofrecen un modelo combinatorio de las configuraciones de puntos, los arreglos de doble pseudorectas han sido introducidos como modelos de configuraciones de convexos disjuntos~\cite{hp-adp-08}. Nuestro trabajo de implementación nos permitía manipular estos objetos y familiarizarnos con sus propiedades, lo que resultó útil para nuestro estudio de la dualidad.


\paragraph{Politopalidad de productos cartesianos.} El \defn{producto cartesiano} de grafos está definido de tal manera que el grafo de un producto de politopos es el producto de los grafos de sus factores. Así, un producto de grafos politopales es automáticamente politopal. Estudiamos en primer lugar la recíproca: ¿Son los factores de un producto politopal necesariamente politopales? El Capítulo~\ref{chap:nonpolytopal} responde a esta pregunta con una atención particular a los grafos regulares y a sus realizaciones como politopos simples (los politopos simples tienen propiedades muy particulares: entre otras, están determinados por sus grafos~\cite{bm-ppi-87,k-swtsp-88}). Discutimos en particular la pregunta de la politopalidad del producto de dos grafos de \mbox{Petersen}, planteada por G\"unter Ziegler~\cite{crm}. Este trabajo sobre la politopalidad de productos nos llevaba además a estudiar ejemplos de grafos no-politopales aunque satisfagan las condiciones necesarias conocidas para ser politopales~\cite{b-gscps-61,k-ppg-64,b-ncp-67}.

\begin{figure}
	\capstart
	\centerline{\includegraphics[scale=.92]{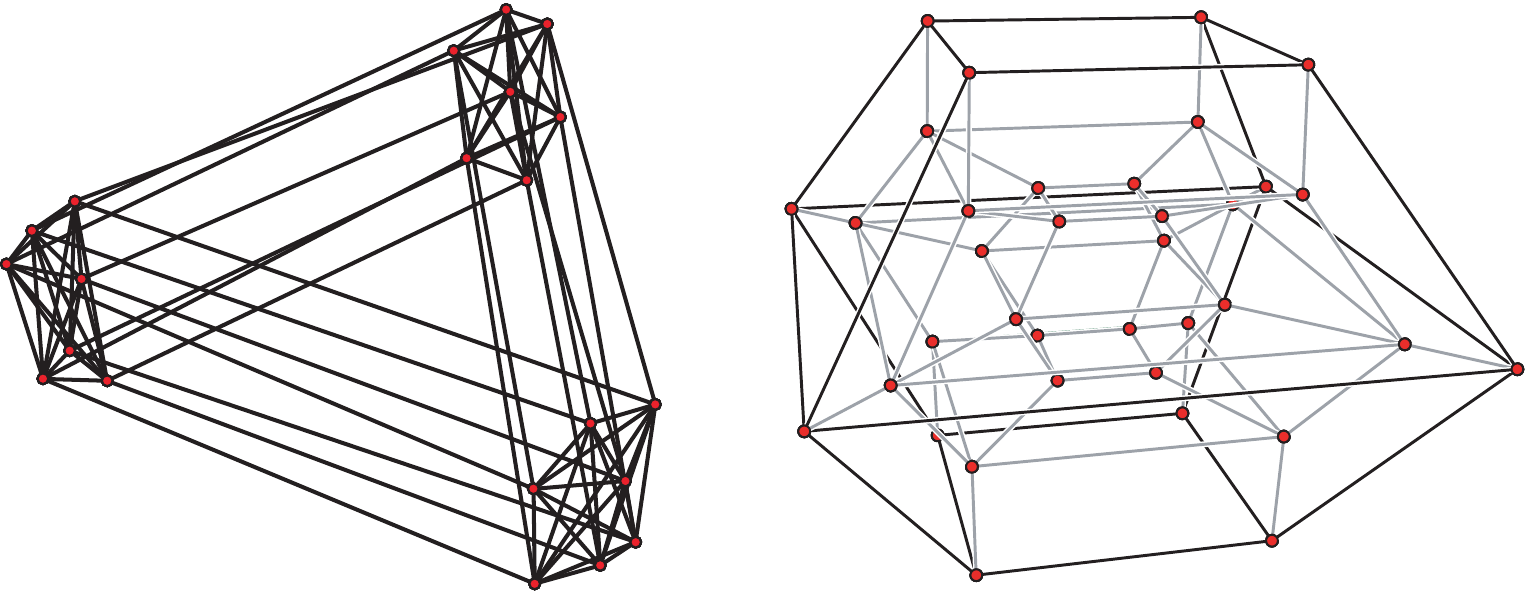}}
	\caption[Ejemplos de productos de grafos]{(Izquierda) Un producto de grafos completos. Este grafo es el de un producto de símplices, pero también el de politopos de dimensiónes inferiores. (Derecha) Un producto politopal de grafos no-politopales.}
	\label{E:fig:products}
\end{figure}

\svs
Buscamos después la dimensión mínima de un politopo cuyo grafo es isomorfo al grafo de un producto dado de politopos. Esta pregunta se inserta en la búsqueda sistemática de politopos de dimensión extrema con propiedades fijadas: por ejemplo, los \poli{d}topos \emph{neighborly} \cite{g-ncp-63} son aquéllos cuyo \esqueleto{\Fracfloor{d}{2}} es completo (y consiguen el máximo número de caras permitido por el Teorema de la Cota Superior~\cite{m-mnfcp-70}); los \poli{d}topos cúbicos \emph{neighborly}~\cite{jz-ncp-00,z-ppp-04,js-ncps-07} son aquéllos cuyo \esqueleto{\Fracfloor{d}{2}} es el del $n$-cubo,~\etc{} Para construir politopos de pequeña dimensión con una propiedad dada, es natural y a menudo eficaz empezar en una dimensión suficiente para asegurarse de la existencia de tales politopos, y proyectar estos politopos sobre subespacios de dimensión inferior conservando la propiedad buscada. Estas técnicas y sus límites han sido ampliamente estudiadas en la literatura~\cite{az-dpmsp-99,z-ppp-04,sz-capdp,s-tovnms-09,sanyal-phd}. Las aplicamos en el Capítulo~\ref{chap:psn} para construir politopos \defn{$(k,\sub{n})$-neighborly prodsimpliciales}, que tienen el mismo \esqueleto{k} que cierto producto de símplices~$\simplex_{\sub{n}} \eqdef \simplex_{n_1}\times\cdots\times\simplex_{n_r}$. En este capítulo, también damos coordenadas enteras para tales politopos utilizando sumas de Minkowski explícitas de politopos cíclicos.


\section{Multitriangulaciones}

Fijamos~$n$ vértices sobre un círculo y consideramos las cuerdas entre ellos. Dos cuerdas se cruzan cuando los segmentos abiertos se intersecan. Un \defn{\kcruce{\ell}} es un conjunto de~$\ell$ cuerdas que se cruzan mutuamente. 
A semejanza de las triangulaciones, nos interesamos por los conjuntos maximales que evitan estas configuraciones:

\begin{definition_E}
Una \defn{\ktric{k}} del \gon{n}o es un conjunto maximal de cuerdas del \gon{n}o sin \kcruce{(k+1)}s.
\end{definition_E}

Por definición, una cuerda puede aparecer en un \kcruce{(k+1)} solamente si deja al menos~$k$ vértices de cada lado.  Decimos que una tal cuerda es \defn{\krel{k}e}. Por maximalidad, cualquier \ktric{k} contiene todas las cuerdas que no son \krel{k}es. Entre estas cuerdas, llamamos \defn{cuerdas del \kborde{k}} a las que separan~$k-1$ vértices de todos los otros vértices.

\begin{example_E}\label{E:exm:petitscas}
Para ciertas valores de~$k$ y~$n$, se describen fácilmente las \ktrics{k}~del~\gon{n}o:
\begin{itemize}
\item[~~~~\fbox{$k=1$}] Las \ktrics{1} son las triangulaciones del \gon{n}o.
\item[~~~~\fbox{$n=2k+1$}] El grafo completo~$K_{2k+1}$ es la única \ktric{k} del \gon{(2k+1)}o porque ningúna de sus cuerdas es \krel{k}e.
\item[~~~~\fbox{$n=2k+2$}] El conjunto de las cuerdas del \gon{(2k+2)}o contiene un \kcruce{(k+1)} formado por las~$k+1$ diagonales largas del \gon{(2k+2)}o. Así existen~$k+1$ \ktrics{k} del \gon{(2k+2)}o, obtenidas retirando del grafo completo una de sus diagonales largas.
\item[~~~~\fbox{$n=2k+3$}] Las cuerdas \krel{k}es del \gon{(2k+3)}o forman un ciclo poligonal en el cual dos cuerdas se cruzan siempre salvo si son consecutivas. Por consiguiente, las \ktrics{k} del \gon{(2k+3)}o son las uniones disjuntas de~$k$ pares de cuerdas \krel{k}es consecutivas (además de todas las cuerdas que no son \krel{k}es).
\end{itemize}
\end{example_E}


\subsection*{Resultados anteriores}

Los conjuntos de cuerdas del \gon{n}o sin \kcruce{(k+1)}s aparecen en el contexto de la teoría extrema de los grafos geométricos (ver~\cite[Capítulo~14]{pa-cg-95}, \cite[Capítulo~1]{f-gga-04} y la discusión en~\cite{cp-tttccp-92}). En el grafo de intersección de las cuerdas del \gon{n}o, estos conjuntos inducen en efecto subgrafos sin $k$-clan, y se pueden así acertar del resultado clásico de Turán (que acota el  número de aristas de un grafo sin $k$-clan). Vasilis~Capoyleas y Janos Pach~\cite{cp-tttccp-92} demostraron que estos conjuntos no pueden tener más que $k(2n-2k-1)$ cuerdas. Después, \mbox{Tomoki} \mbox{Nakamigawa}~\cite{n-gdfcp-00}, e independientemente Andreas Dress, Jacobus Koolen y Vincent \mbox{Moulton}~\cite{dkm-lahp-02}, probaron que todas las \ktrics{k} alcanzan este cota. Estas dos pruebas están basadas en la operación de flip que transforma una \ktric{k} en otra cambiando la posición de una sola cuerda. Tomoki Nakamigawa~\cite{n-gdfcp-00} demostró que cualquiera cuerda \krel{k}e de una \ktric{k} puede ser flipada, y que el grafo de flips es conexo. Queremos observar que en todos estos trabajos anteriores, los resultados están obtenidos de ``manera \mbox{indirecta}'': primero, una operación de contracción (similar a la contracción de una arista del borde de una triangulación) es introducida y utilizada para demostrar la existencia del flip (en toda generalidad~\cite{n-gdfcp-00}, o solamente en casos particulares~\cite{dkm-lahp-02}) y la conexidad del grafo de flips, de la cual es deducida la formula del número de cuerdas. Para resumir:

\begin{theorem_E}[\cite{cp-tttccp-92,n-gdfcp-00,dkm-lahp-02}]\label{E:theo:fundamental}
\begin{enumerate}[(i)]
\item Existe una operación inductiva que transforma las \ktrics{k} del \gon{(n+1)}o en \ktrics{k} del \gon{n}o y \viceversa~\cite{n-gdfcp-00}.
\item Cualquiera cuerda \krel{k}e de una \ktric{k} del \gon{n}o puede ser flipada
y el grafo de flips es regular y conexo~\cite{n-gdfcp-00,dkm-lahp-02}.
\item Todas las \ktrics{k} del \gon{n}o tienen exactamente~$k(2n-2k-1)$ cuerdas~\cite{cp-tttccp-92,n-gdfcp-00,dkm-lahp-02}.\qed
\end{enumerate}
\end{theorem_E}

Jakob~Jonsson~\cite{j-gt-03,j-gtdfssp-05} completó estos resultados en dos direcciones. Por una parte, estudió las propiedades enumerativas de las multitriangulaciones:

\begin{theorem_E}[\cite{j-gtdfssp-05}]\label{E:theo:enumeration}
El número de \ktrics{k} del \gon{n}o esta dado por:
$$\det(C_{n-i-j})_{1\le i,j\le k}=\det\begin{pmatrix} C_{n-2} & C_{n-3} & \edots & C_{n-k-1} \\ C_{n-3} & \edots & \edots & \edots \\ \edots & \edots & \edots & C_{n-2k+1} \\ C_{n-k-1} & \edots & C_{n-2k+1} & C_{n-2k} \end{pmatrix},$$
donde~$C_p \eqdef \frac{1}{p+1}{2p \choose p}$ denota el $p$-ésimo número de Catalan.\qed
\end{theorem_E}

La prueba de este teorema se apoya en resultados más generales concerniendo a los $0/1$-relle\-nos de poliominos maximales para ciertas restricciones sobre sus secuencias diagonales de~$1$.

Cuando~$k=1$, encontramos simplemente los números de Catalan que cuentan no solamente las triangulaciones, sino también los caminos de Dyck. El determinante de Hankel que aparece en el teorema precedente también cuenta ciertas familias de $k$-tuplas de caminos de Dyck sin cruce~\cite{gv-bdphlf-85}. La igualdad de los cardinales de estas dos familias combinatorias motiva la búsqueda de una prueba biyectiva que podría aclarar el resultado del Teorema~\ref{E:theo:enumeration}. Aunque Sergi \mbox{Elizalde}~\cite{e-btdp-07} y Carlos Nicolas~\cite{n-abtdp-09} propusieran dos biyecccións diferentes en el caso~$k=2$, esta pregunta queda abierta en el caso general. Discutimos este problema en el Capítulo~\ref{chap:multiassociahedron}.

\svs
Por otra parte, Jakob~Jonsson~\cite{j-gt-03} estudió el complejo simplicial~$\Delta_{n,k}$ formado por los conjuntos de cuerdas \krel{k}es del \gon{n}o sin \kcruce{(k+1)}s. Como todos sus elementos maximales tienen~$k(n-2k-1)$ cuerdas, este complejo es puro de dimensión~$k(n-2k-1)-1$. Jakob~Jonsson demostró el teorema siguiente:

\begin{theorem_E}[\cite{j-gt-03}]
El complejo simplicial~$\Delta_{n,k}$ es una esfera de dimensión ${k(n-2k-1)-1}$.\qed 
\end{theorem_E}

Volkmar Welker conjeturó de que este complejo simplicial es además politopal. Esta conjetura es cierta cuando~${k=1}$ (el complejo simplicial de los conjuntos sin cruce de cuerdas internas del \gon{n}o es isomorfo al complejo de borde del dual del asociaedro) y en los casos del Ejemplo~\ref{E:exm:petitscas} (donde obtenemos respectivamente un punto, un símplex y un politopo cíclico). Discutimos esta pregunta en el Capítulo~\ref{chap:multiassociahedron}.


\subsection*{Descomposicones de multitriangulaciones en estrellas}

Nuestra contribución al estudio de las \ktrics{k} está basada en las estrellas, que generalizan los triángulos de las triangulaciones:

\begin{definition_E}
Una \defn{\kestrella{k}} es un polígono (no-simple) con $2k+1$~vértices $s_0,\dots,s_{2k}$ en orden sobre el círculo, y~$2k+1$ cuerdas $[s_0,s_k],[s_1,s_{k+1}],\dots,[s_{2k},s_{k-1}]$.
\end{definition_E}

Las estrellas desempeñan para las multitriangulaciones exactamente el mismo papel que los triángulos por las triangulaciones: las descomponen en entidades geométricas más simples y permiten entender su combinatoria. Como prueba, el estudio de las propiedades de incidencia de las estrellas en las multitriangulaciones lleva a nuevas pruebas de todas las  propiedades combinatorias de las multitriangulaciones conocidas, que presentamos en el Capítulo~\ref{chap:stars}. Nuestro primer resultado estructural es el siguiente:

\begin{theorem_E}\label{E:theo:incidences}
Sea~$T$ una \ktric{k} del \gon{n}o (con~$n\ge 2k+1$).
\begin{enumerate}[(i)]
\item Una cuerda \krel{k}e está contenida en dos \kestrella{k}s de~$T$, una de cada lado; una cuerda del \kborde{k} está contenida en una \kestrella{k} de~$T$; una cuerda que no es \krel{k}e ni del \kborde{k} no está contenida en ninguna \kestrella{k} de~$T$.
\item $T$~tiene exactamente~$n-2k$ \kestrella{k}s y~$k(2n-2k-1)$ cuerdas.\qed
\end{enumerate}
\end{theorem_E}

\begin{figure}[h]
	\centerline{\includegraphics[scale=1]{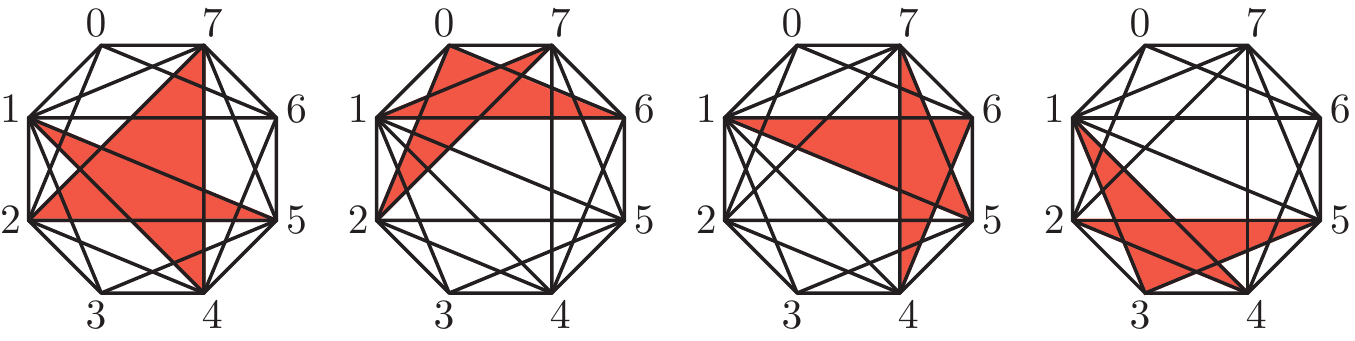}}
	\caption[Las \kestrella{2}s de la \ktric{2} del octógono de la Figura~\ref{E:fig:ktristardual}]{Las \kestrella{2}s de la \ktric{2} del octógono de la Figura~\ref{E:fig:ktristardual}.}
	\label{E:fig:2triang8pointsstars}
\end{figure}

Probamos el punto~(ii) de este teorema con un argumento directo de doble cuenta, exhibiendo dos relaciones independientes entre el número de cuerdas y el número de \kestrella{k}s de una \ktric{k}. Nuestra primera relación proviene del punto~(i) del teorema, mientras que la segunda se apoya en las propiedades de las bisectrices de las \kestrella{k}s. Una \defn{bisectriz} de una \kestrella{k} es una bisectriz de uno de sus ángulos, es decir, una recta que pasa por uno de sus vértices y que separa sus otros vértices en dos conjuntos de cardinal~$k$. Nuestra segunda relación es una consecuencia directa de la correspondencia entre los pares de \kestrella{k}s de~$T$ y las cuerdas del \gon{n}o que no están en~$T$:

\begin{theorem_E}
Sea~$T$ una \ktric{k} del \gon{n}o.
\begin{enumerate}[(i)]
\item Cada par de \kestrella{k}s de~$T$ tiene una única bisectriz común, que no está en~$T$.
\item Recíprocamente, cada cuerda que no está en~$T$ es la bisectriz común de un único par de \kestrella{k}s de~$T$.\qed
\end{enumerate}
\end{theorem_E}

Utilizamos además las estrellas y sus bisectrices para aclarar la \defn{operación de flip} y para proponer una interpretación local. De la misma manera que un flip en una triangulación substituye una diagonal por la otra en un cuadrángulo formado por dos triángulos adyacentes, un flip en una multitriangulación se interpreta como una transformación donde solo intervienen  dos estrellas adyacentes:

\begin{theorem_E}\label{E:theo:flip}
Sea~$T$ una \ktric{k} del \gon{n}o, sea~$e$ una cuerda \krel{k}e de~$T$ y sea~$f$ la bisectriz común a las \kestrella{k}s de~$T$ que contienen~$e$. Entonces~$T\diffsym\{e,f\}$ es una \ktric{k} del \gon{n}o y es la única, salvo~$T$ misma, que contiene~$T\ssm\{e\}$.~\qed
\end{theorem_E}

\begin{figure}[h]
	\centerline{\includegraphics[scale=1]{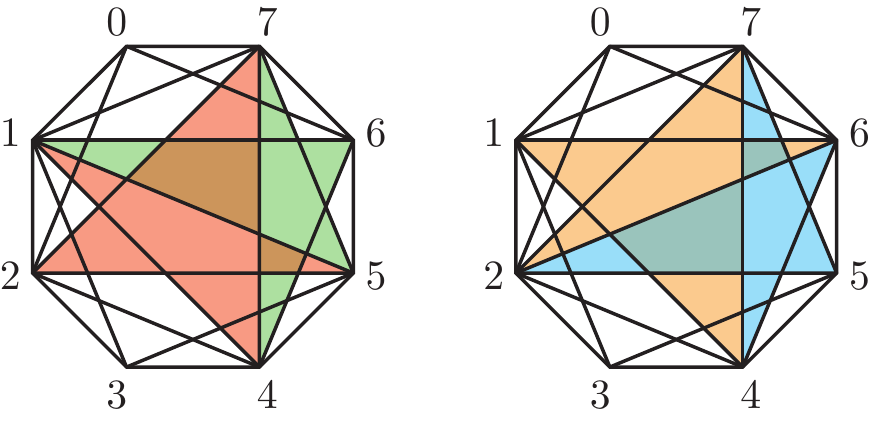}}
	\caption[Un flip en la \ktric{2} del octógono de la Figura~\ref{E:fig:ktristardual}]{Un flip en la \ktric{2} del octógono de la Figura~\ref{E:fig:ktristardual}.}
	\label{E:fig:2triang8pointsflip}
\end{figure}

Esta interpretación simplifica el estudio del grafo de flips (cuyos vértices son las \ktrics{k} del \gon{n}o y cuyas aristas son los flips entre ellas) y proporciona nuevas pruebas y extensiones parciales de los resultados de~\cite{n-gdfcp-00}:

\begin{theorem_E}
El grafo de flips de las \ktrics{k} del \gon{n}o es \regular{k(n-2k-1)}, conexo, y por todo~$n>4k^2(2k+1)$, su diámetro~$\delta_{n,k}$ está acotado por $$2\Fracfloor{n}{2}\left(k+\frac{1}{2}\right)-k(2k+3) \le \delta_{n,k} \le 2k(n-4k-1).$$
\end{theorem_E}
\vspace{-1cm}\qed
\vspace{1cm}

Utilizamos también las estrellas para estudiar las \defn{$k$-orejas} de las multitriangulaciones, es decir las cuerdas que separan~$k$ vértices de los otros vértices. Proponemos una prueba simple del hecho que cada \ktric{k} tiene al menos~$2k$ $k$-orejas~\cite{n-gdfcp-00}, y damos luego varias caracterizaciones de las \ktrics{k} que alcanzan esta cota:

\begin{theorem_E}
El número de $k$-orejas de una \ktric{k} es igual a~$2k$ más el número de \defn{\kestrella{k}s internas}, \ie las que no contienen ninguna cuerdas del \kborde{k}.

Si~$k>1$ y $T$ es una \ktric{k}, las afirmaciones siguientes son equivalentes:
\begin{enumerate}[(i)]
\item $T$~tiene exactamente $2k$~$k$-orejas;
\item $T$ no tiene ninguna \kestrella{k} interna;
\item $T$~es \defn{$k$-colorable}, \ie existe una $k$-coloración de sus cuerdas \krel{k}es tal que ningun cruce es monocromático;
\item el conjunto de cuerdas \krel{k}es de~$T$ es la unión disjunta de $k$ acordeones (un \defn{acordeón} es una secuencia de cuerdas~$[a_i,b_i]$ tal que para todo~$i$, o bien $a_{i+1}=a_i$ y $b_{i+1}=b_i+1$ o bien $a_{i+1}=a_i-1$ y $b_{i+1}=b_i$).\qed
\end{enumerate}
\end{theorem_E}

Después, reinterpretamos en términos de estrellas la \defn{operación inductiva} del Teorema~\ref{E:theo:fundamental}(i) que transforma las \ktrics{k} del \gon{(n+1)}o en \ktrics{k} del \gon{n}o y \textit{vice versa}. En un sentido aplastamos una \kestrella{k} (que contiene una cuerda del \kborde{k}) hasta un \kcruce{k}, y en el otro sentido inflamos un \kcruce{k} para obtener una \kestrella{k}. Intencionadamente, solo presentamos esta operación al final del Capítulo~\ref{chap:stars} para destacar que ninguna de nuestras pruebas precedentes usa este transformación inductiva.

\begin{figure}[h]
	\centerline{\includegraphics[scale=.8]{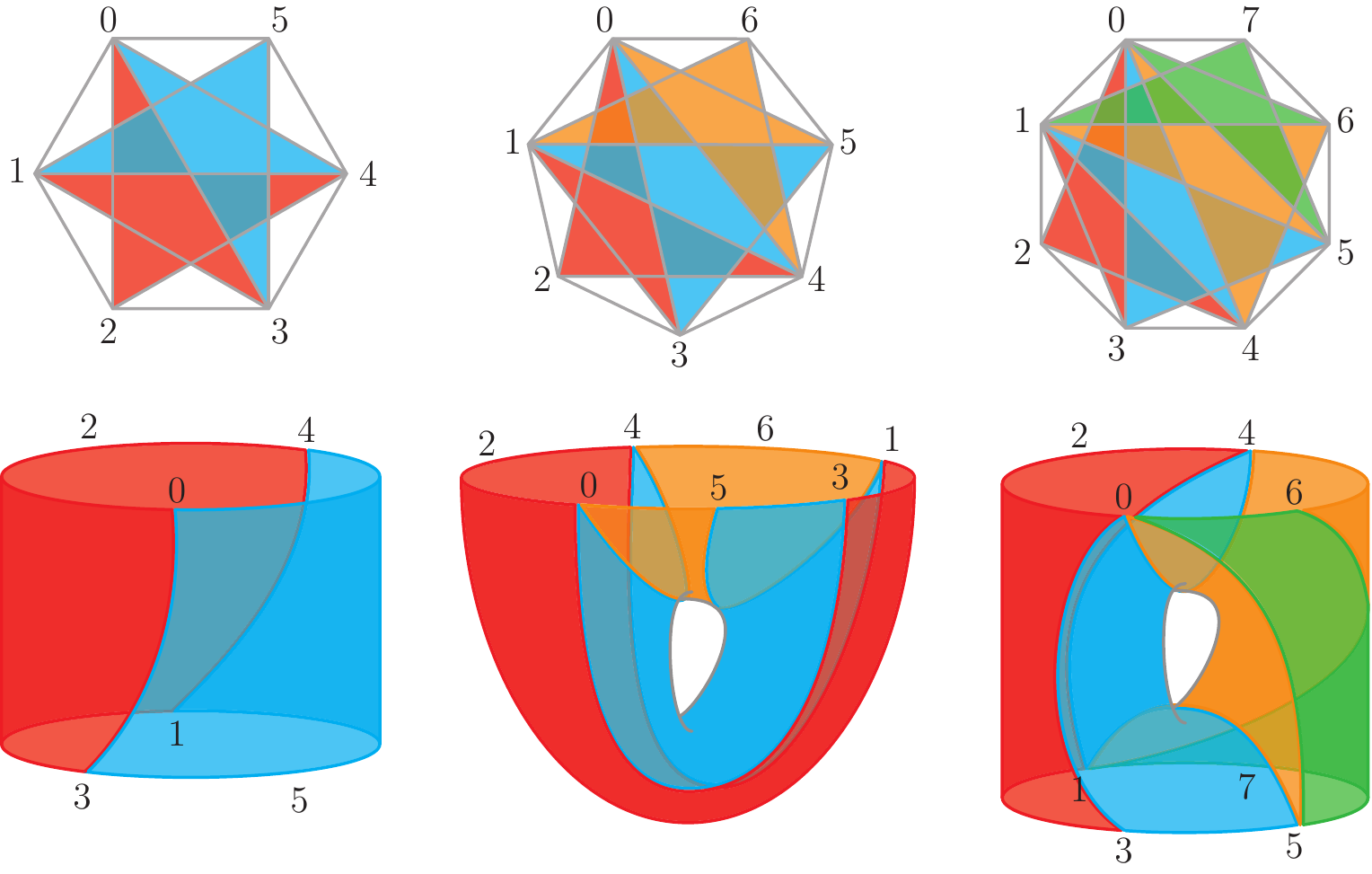}}
	\caption[Descomposiciones de superficies asociadas a las multitriangulaciones]{Descomposiciones de superficies associadas a multitriangulaciones.}
	\label{E:fig:surfaces}
\end{figure}

\svs
Para acabar el Capítulo~\ref{chap:stars}, utilizamos el resultado del Teorema~\ref{E:theo:incidences} para interpretar una \ktric{k} como una \defn{descomposición de superficie} en $n-2k$ \gon{k}os. Algunos ejemplos son presentados en la Figura~\ref{E:fig:surfaces}. Explotamos esta interpretación para construir, con ayuda de multitriangulaciones, unas descomposiciones muy regulares de una familia infinita de superficies.

\subsection*{Multitriangulaciones, pseudotriangulaciones y dualidad}

En el Capítulo~\ref{chap:mpt}, nos interesamos por la interpretación de las multitriangulaciones en el espacio de rectas del plano, es decir en la banda de M\"obius. Estamos en el contexto de la dualidad clásica (ver~\cite[Capítulos~5,6]{f-gga-04} y~\cite{g-pa-97}) entre configuraciones de puntos y arreglos de pseudorectas: el conjunto~$p^*$ de las rectas que pasan por un punto~$p$ del plano es una \defn{pseudorecta} del espacio de rectas (una curva cerrada simple no-separando), y el conjunto ${P^* \eqdef \ens{p^*}{p\in P}}$ de las pseudorectas duales de un conjunto finito~$P$ de puntos es un \defn{arreglo de pseudorectas} (dos pseudorectas se intersecan exactamente una vez).

\begin{figure}
	\capstart
	\centerline{\includegraphics[scale=1]{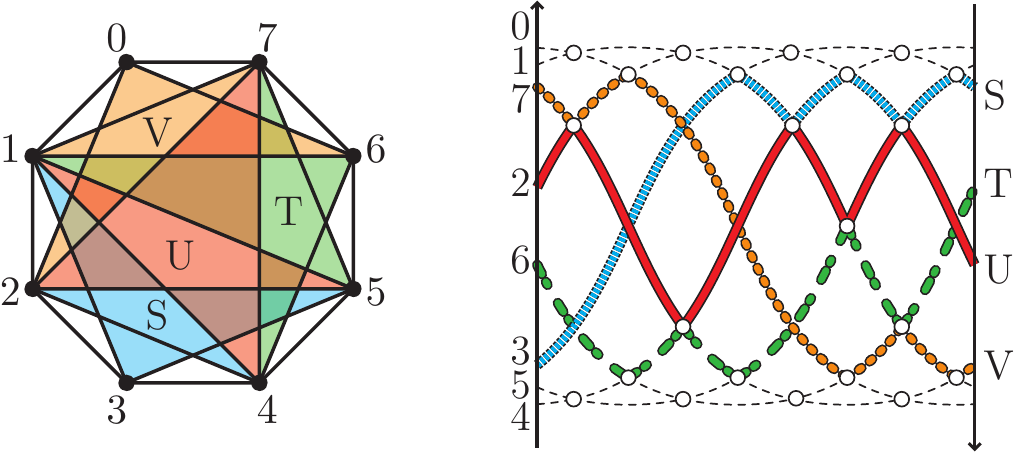}}
	\caption[Arreglo dual de una multitriangulación]{Una \ktric{2} del octógono y su arreglo dual.}
	\label{E:fig:dualmulti}
\end{figure}

El punto de partida del Capítulo~\ref{chap:mpt} es la observación siguiente (ver Figura~\ref{E:fig:dualmulti}):

\begin{observación_E}
Sea $T$ una \ktric{k} del \gon{n}o. Entonces:
\begin{enumerate}[(i)]
\item el conjunto~$S^*$ de las bisectrices de una \kestrella{k}~$S$ de~$T$ es una pseudorecta;
\item la bisectriz común de dos \kestrella{k}s~$R$ y~$S$ es un punto de cruce de las pseudorectas $R^*$ y~$S^*$, mientras que una cuerda común de~$R$ y~$S$ es un punto de contacto~entre~$R^*$~et~$S^*$;
\item el conjunto~$T^* \eqdef \ens{S^*}{S \text{ \kestrella{k} de } T}$ de todas las pseudorectas duales de las \kestrella{k}s de~$T$ es un \defn{arreglo de pseudorectas con puntos de contacto} (dos pseudorectas se cruzan exactamente una vez, pero pueden tocarse en un número finito de puntos de contacto);
\item el soporte de este arreglo cubre exactamente el soporte del arreglo dual~$V_n^*$ de los vértices del \gon{n}o salvo sus $k$~primeros niveles.
\end{enumerate}
\end{observación_E}

Demostramos la recíproca de esta observación:

\begin{theorem_E}
Cualquier arreglo de pseudorectas con puntos de contacto cuyo soporte cubre exactamente el soporte del arreglo dual~$V_n^*$ de los vértices del \gon{n}o salvo sus $k$~primeros niveles es el arreglo dual de una \ktric{k} del \gon{n}o.\qed
\end{theorem_E}

En~\cite{pv-ot-94,pv-ptta-96}, una observación similar ha sido hecha para las pseudotriangulaciones puntiagudas de un conjunto de puntos en posición general. Introducidas para el estudio del complejo de visibilidad de obstáculos convexos disjuntos del plano~\cite{pv-tsvcp-96,pv-vc-96}, las pseudotriangulaciones han sido utilizadas en varios contextos geométricos (como la planificación de movimiento y la rigidez \cite{s-ptrmp-05,horsssssw-pmrgpt-05}) y sus propiedades han sido ampliamente estudiadas (número de pseudotriangulaciones~\cite{aaks-cmpt-04,aoss-nptcps-08}, politopo de las pseudotriangulaciones~\cite{rss-empppt-03}, consideraciones algorítmicas~\cite{b-eptp-05,bkps-ceppgfa-06,hp-cpcpb-07}, \etc). Remitimos a~\cite{rss-pt-06} para más detalles.

\begin{definition_E}
Un \defn{pseudotriángulo} es un polígono~$\Delta$ con exactamente tres ángulos convexos, relacionados por tres cadenas poligonales cóncavas. Una recta es \defn{tangente} a~$\Delta$ si pasa por un ángulo convexo de~$\Delta$ y separa sus dos aristas adyacentes, o si pasa por un ángulo cóncavo de~$\Delta$ y no separa sus dos aristas adyacentes. Una \defn{pseudotriangulación} de un conjunto~$P$ de puntos en posición general es un conjunto de aristas que descompone la envoltura convexa de~$P$ en pseudotriángulos.
\end{definition_E}

En esta memoria,  consideramos sólo pseudotriangulaciones \defn{puntiagudas}, es decir pseudotriangulaciones tales que por cada punto~$p\in P$ pasa una recta que define un semi-plano cerrado que no contiene ninguna arista saliendo de~$p$. Para una pseudotriangulación puntiaguda~$T$ de~$P$, ha sido observado~\cite{pv-ot-94,pv-ptta-96} que:
\begin{enumerate}[(i)]
\item el conjunto~$\Delta^*$ de las rectas tangentes a un pseudotriángulo~$\Delta$ es una pseudorecta; y
\item el conjunto~$T^* \eqdef \ens{\Delta^*}{\Delta \text{ pseudotriángulo de } T}$ es un arreglo de pseudorectas con puntos de contacto soportada por el arreglo~$P^*$ dual de~$P$ salvo su primer nivel.
\end{enumerate}

\begin{figure}
	\capstart
	\centerline{\includegraphics[scale=1]{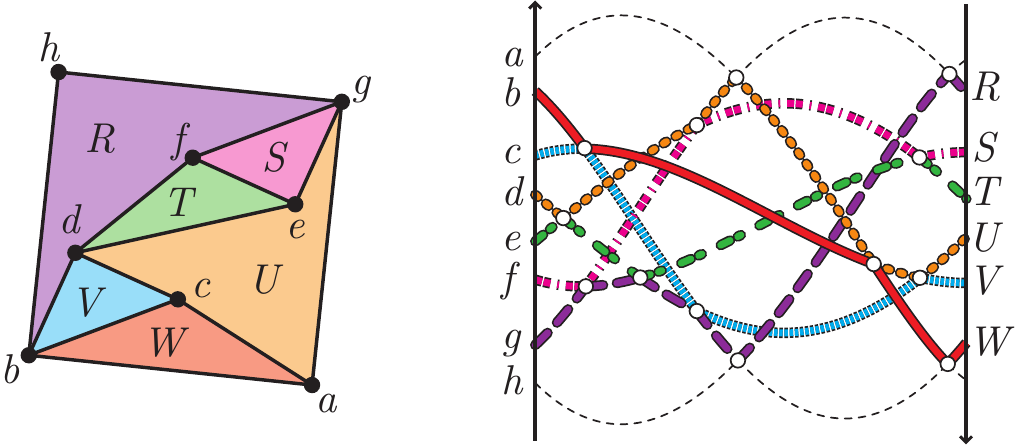}}
	\caption[Arreglo dual de una pseudotriangulación]{Una pseudotriangulación y su arreglo dual.}
	\label{E:fig:dualpseudo}
\end{figure}

Damos diferentes pruebas de la recíproca:

\begin{theorem_E}
Sea~$P$ un conjunto de puntos en posición general y~$P^*$ su arreglo dual. Cualquier arreglo de pseudorectas con puntos de contacto cuyo soporte es precisamente el soporte del arreglo~$P^*$ salvo su primer nivel es el arreglo dual de una pseudotriangulación~de~$P$.\qed
\end{theorem_E}

Motivados por estos dos teoremas, consideramos luego los arreglos de pseudorectas que comparten un mismo soporte. Definimos una operación de flip que corresponde al flip en las multitriangulaciones y las pseudotriangulaciones: 

\begin{definition_E}
Dos arreglos de pseudorectas con el mismo soporte están relacionados por un flip si la diferencia simétrica de sus conjuntos de puntos de contacto se reduce a un par~$\{u,v\}$. En este caso, en uno de los arreglos, $u$ es un punto de contacto de dos pseudorectas que se cruzan en~$v$, mientras que en el otro arreglo, $v$ es un punto de contacto de dos pseudorectas que se cruzan en~$u$.
\end{definition_E}

Estudiamos el grafo~$G(\cS)$ de flips de los arreglos soportados por~$\cS$. Por ejemplo, cuando~$\cS$ es el soporte de un arreglo de dos pseudorectas con~$p$ puntos de contacto, el grafo~$G(\cS)$ es el grafo completo con~$p+1$ vértices. Nos interesamos por ciertas orientaciones acíclicas del grafo~$G(\cS)$, dadas por cortes verticales del soporte~$\cS$. Una tal orientación tiene una única fuente que llamamos \defn{arreglo glotón} y que caracterizamos en términos de redes de selección. El estudio de estos arreglos glotones y de sus transformaciones cuando la orientación de~$G(\cS)$ evoluciona proporciona un algoritmo de enumeración de los arreglos de pseudorectas con puntos de contacto que comparten un mismo soporte, cuyo espacio de trabajo queda polinomial. Así, aclaramos y damos una prueba complementaría del algoritmo similar que existía para la enumeración de las pseudotriangulaciones de un conjunto de puntos~\cite{bkps-ceppgfa-06}.

\svs
Volvemos después a un contexto más particular. A la luz de la dualidad entre las multitriangulaciones (resp.~las pseudotriangulaciones) y los arreglos de pseudorectas con puntos de contacto, proponemos la generalización siguiente de las multitriangulaciones cuando los puntos no están en posición convexa:

\begin{definition_E}
Una \defn{\ptc{k}} de un arreglo de pseudorectas~$L$ es un arreglo de pseudorectas con puntos de contacto cuyo soporte es el arreglo~$L$ salvo sus~$k$ primer niveles. Una \ptc{k} de un conjunto de puntos~$P$ en posición general es un conjunto de aristas~$T$ que corresponden por dualidad a los puntos de contacto de una \ptc{k}~$T^*$ de~$P^*$.
\end{definition_E}

Mostramos que todas las \ptcs{k} de un conjunto de puntos~$P$ tienen exactamente ${k(2|P|-2k-1)}$ aristas, y que no pueden contener configuraciones de~$2k+1$ aristas alternativas, pero que pueden eventualmente contener un \kcruce{(k+1)}. Estudiamos después sus estrellas: una estrella de una \ptc{k}~$T$ de~$P$ es un polígono formado por el conjunto de las aristas correspondientes a los puntos de contacto de una pseudorecta fijada de~$T^*$. Discutimos su número possible de esquinas, y mostramos que para todo punto~$q$ del plano, la suma de los índices de las estrellas de~$T$ alrededor de~$q$ es independiente de~$T$.

\svs
Acabamos el Capítulo~\ref{chap:mpt} con tres preguntas relacionadas con las \mbox{\mptcs{}:}
\begin{enumerate}[(i)]
\item Estudiamos primero las \defn{\mptcs{} iteradas}: una \ptc{k} de una \ptc{m} de un arreglo de pseudorectas~$L$ es un \ptc{(k+m)} de~$L$. Presentamos sin embargo un ejemplo de \ktric{2} que no contiene ninguna triangulación. Mostramos en cambio que las \mptcs{} glotonas de un arreglo son iteradas de la pseudotriangulación glotona.
\item Damos después una caracterización de las aristas de la \ptc{k} glotona de un conjunto de puntos~$P$ en términos de \defn{$k$-árboles de horizonte} de~$P$. Esta caracterización generaliza una observación de Michel Pocchiola~\cite{p-htvpt-97} para las \mbox{pseudotriangulaciones}.
\item Finalmente, definimos las \mptcs{} de configuraciones de convexos disjuntos del plano y estudiamos sus propiedades elementales. Los arreglos duales de las configuraciones de convexos son los \defn{arreglos de doble pseudorectas}, introducidos por Luc Habert y Michel Pocchiola~\cite{hp-adp-08}. En el Apéndice~\ref{app:implementations}, manipulamos estos arreglos para enumerar los arreglos con a lo más~$5$ doble pseudorectas.
\end{enumerate}


\subsection*{Tres problemas abiertos}

Finalmente, discutimos en el Capítulo~\ref{chap:multiassociahedron} tres problemas abiertos que ilustran la riqueza combinatoria y geométrica de las multitriangulaciones. Nuestro objectivo no es únicamente presentar algunos resultados positivos parciales, sino también destacar ideas naturales basadas en las estrellas que pueden ser fértiles aunque suspendieran para resolver estos problemas.

\svs
El primer problema que discutimos es encontrar una \defn{biyección explícita} entre el conjunto de las multitriangulaciones y el conjunto de las $k$-tuplas de caminos de Dyck sin cruce, que están ambos enumerados por el determinante de Hankel del Teorema~\ref{E:theo:enumeration}. Si en una triangulación~$T$, denotamos por~$\delta_i(T)$ el número de triángulos de~$T$ cuyo primer vértice es~$i$, entonces la aplicación~$T\mapsto N^{\delta_0(T)}EN^{\delta_1(T)}E\dots N^{\delta_{n-3}(T)}E$ es una biyección de las triangulaciones del \gon{n}o a los caminos de Dyck de semi-longitud~$n-2$ (donde~$N$ y~$E$ denotan los pasos norte y este). Jakob Jonsson~\cite{j-gt-03} generalizó esta observación comparando la repartición de las secuencias de grados entrantes de las \ktrics{k} con las signaturas de las $k$-tuplas de caminos de Dyck (sin dar las definiciones concisas, señalamos que estas dos $k$-tuplas de secuencias generalizan la secuencia~$(\delta_i(T))$ de una triangulación~$T$ y la secuencia de los potencias de~$N$ en un camino de Dyck). Motivados por este resultado, buscamos como definir $k$-coloraciones de las cuerdas \krel{k}es de una \ktric{k} de tal modo que las secuencias respectivas de grado entrante de cada color defina $k$-tuplas de caminos de Dyck sin cruces. Presentamos una coloración basada en las estrellas, que satisface esta propiedad, pero con la cual la función de las multitriangulaciones en las $k$-tuplas de caminos de Dyck no es desgraciadamente biyectiva.

\svs
Nuestro segundo problema concierne las propiedades de \defn{rigidez} de las multitriangulaciones. Una triangulación es \defn{mínimalmente rígida} en el plano: los únicos movimientos de sus vértices que preserven las longitudes de sus aristas son las isometrías del plano, y retirar cualquiera arista vuelve la estructura flexible. De manera equivalente, una triangulación satisface la condición de Laman: tiene~$2n-3$ aristas y cada subgrafo con~$m$ vértices tiene a lo más~$2m-3$ aristas. Observamos dos conexiones interesantes entre las multitriangulaciones y la teoría de la rigidez:
\begin{enumerate}[(i)]
\item Primero, una \ktric{k} es $\left(2k,{2k+1 \choose 2}\right)$-\defn{tenso}: contiene~$2kn-{2k+1 \choose 2}$ aristas y cada subgrafo con~$m$ vértices tiene a lo más~$2km-{2k+1 \choose 2}$ aristas. Esta propiedad combinatoria hace de una \ktric{k} un candidato razonable para ser un grafo genéricamente mínimalmente rígido en dimensión~$2k$. Probamos esta conjectura cuando~$k=2$.
\item Después, mostramos que el grafo dual de una $k$-triangulación es~$(k,k)$-tenso. En particular, puede ser descompuesto en~$k$ árboles de expansión con aristas disjuntas (en el dual de una triangulación, solo hay un árbol).
\end{enumerate}

\svs
Finalmente, volvemos a la pregunta de la \defn{realización politopal} del complejo simplicial~$\Delta_{n,k}$ formado por los conjuntos de cuerdas \krel{k}es del \gon{n}o sin \kcruce{(k+1)}s. Presentamos dos contribuciones modestas a esta pregunta:
\begin{enumerate}[(i)]
\item Respondemos por una parte al caso no-trivial mínimo describiendo el espacio de las realizaciones simétricas de~$\Delta_{8,2}$. Obtenemos este resultado en dos etapas: primero enumeramos por ordenador todos los matroides orientados simétricos que realizan nuestro complejo simplicial; después buscamos las realizaciones politopales simétricas de estos matroides~orientados.
\item Consideramos por otra parte la construcción del asociaedro de Jean-Louis Loday~\cite{l-rsp-04}. Interpretamos su construcción con los arreglos duales de las triangulaciones, lo que se generaliza naturalmente a las multitriangulaciones. Obtenemos un politopo cuyas facetas están definidas por desigualdades simples, y que realiza el grafo de flips restringido a ciertas multitriangulaciones (aquéllas con un grafo dual acíclico). Aunque este politopo habría podido ser \apriori una proyección de una realización de~$\Delta_{n,k}$, mostramos que no~lo~es.
\end{enumerate}


\section{Politopalidad de productos}

En la segunda parte de esta tesis, nos interesamos por preguntas de realización politopal de grafos (o de esqueletos) obtenidos como productos cartesianos de otros grafos (o esqueletos).

El \defn{producto cartesiano} de dos politopos~$P,Q$ es el politopo~$P\times Q \eqdef \ens{(p,q)}{p\in P,q\in Q}$. Su estructura combinatoria solo depende de la de sus factores: la dimensión de~$P\times Q$ es la suma de las dimensiones de~$P$ y~$Q$ y las caras no-vacías de~$P\times Q$ son precisamente los productos de una cara no-vacía de~$P$ por una cara no-vacía de~$Q$.


\subsection*{Politopalidad de productos de grafos no-politopales}

En el Capítulo~\ref{chap:nonpolytopal}, nos interesamos por el producto de grafos. El \defn{producto cartesiano} de dos grafos~$G$ y~$H$ es el grafo~$G\times H$ cuyos vértices son~$V(G\times H) \eqdef V(G)\times V(H)$ y cuyas aristas son $E(G\times H) \eqdef \big(E(G)\times V(H)\big)\cup\big(V(G)\times E(H)\big)$. Dicho de otra manera, para todos vértices $a,c\in V(G)$ y~$b,d\in V(H)$, los vértices~$(a,b)$ y~$(c,d)$ de~$G\times H$ son adyacentes si~$a=c$ y $\{b,d\}\in E(H)$, ó~$b=d$ y~$\{a,c\}\in E(G)$. Este producto es coherente con el producto de politopos: el grafo de un producto de politopos es el producto de sus grafos. En particular, el producto de dos grafos politopales es automáticamente politopal. En este capítulo, estudiamos la pregunta recíproca: ¿~Si un producto~$G\times H$ es politopal, son necesariamente sus factores~$G$ y~$H$ politopales~?

El producto de un triángulo por un camino de longitud~$2$ es un contraejemplo simple a esta pregunta: aunque el camino no sea politopal, el producto es el grafo de un \poli{3}topo obtenido pegando dos prismos triangulares por una cara triangular. Eliminamos tales ejemplos exigiendo que cada factor sea un grafo regular. Si~$G$ y~$H$ son regulares de grado respectivos~$d$ y~$e$, entonces el producto~$G\times H$ es \regular{(d+e)} y es natural preguntarse si es el grafo de un \poli{(d+e)}topo simple. El teorema siguiente responde a esta pregunta:

\begin{theorem_E}\label{E:theo:simpleproduct}
Un producto~$G\times H$ es el grafo de un politopo simple si y solo si sus dos factores~$G$ y~$H$ son grafos de politopos simples. En este caso, existe un único politopo simple cuyo grafo es~$G\times H$: es precisamente el producto de los únicos politopos simples cuyos grafos respectivos son~$G$ y~$H$.\qed
\end{theorem_E}

En este teorema, la unicidad del politopo simple realizando~$G\times H$ es una aplicación directa del hecho de que un politopo simple es completamente determinado por su grafo~\cite{bm-ppi-87,k-swtsp-88}. Estos resultados se apoyan en la propiedad siguiente de los politopos simples: cada conjunto de~$k+1$ aristas adyacentes a un mismo vértice de un politopo simple~$P$ define una $k$-cara de~$P$.

Como aplicación del Teorema~\ref{E:theo:simpleproduct} obtenemos una familia infinita de grafos \regular{4}es no-politopales: el producto de un grafo \regular{3} no-politopal por un segmento es no-polito\-pal y \regular{4}.

\begin{figure}[h]
	\centerline{\includegraphics[width=.9\textwidth]{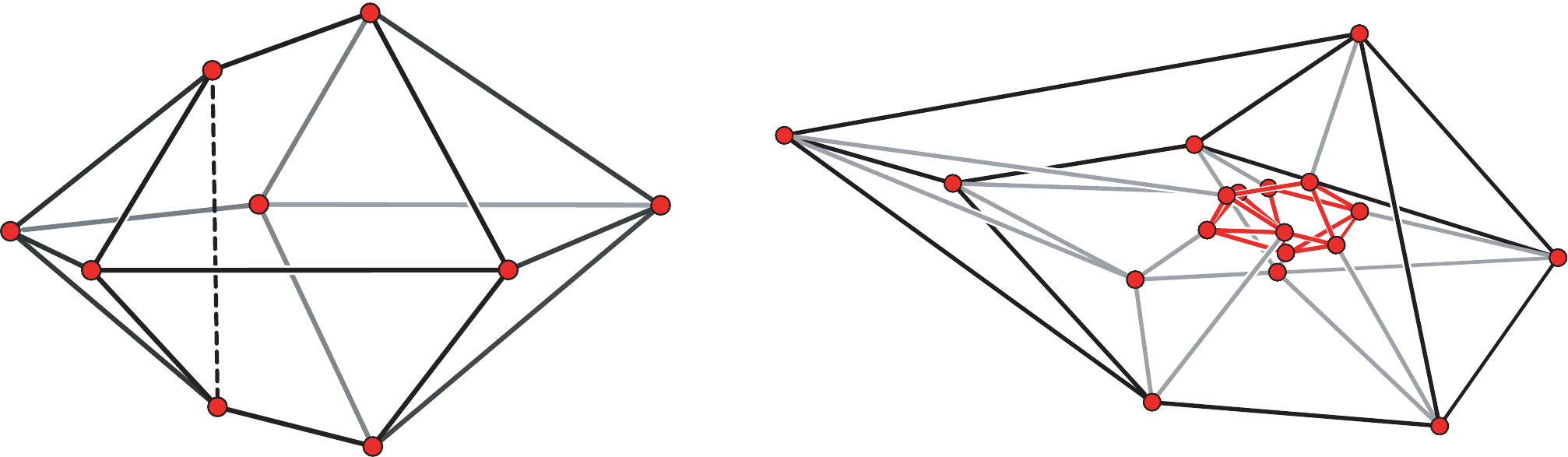}}
	\caption[Un grafo~$H$ \regular{4} y no-politopal y el diagrama de Schlegel de un \poli{4}topo cuyo grafo es el producto de~$H$ por un segmento]{Un grafo~$H$ \regular{4} y no-politopal que es el grafo de una subdivisión regular de un \poli{3}topo (izquierda) y el diagrama de Schlegel de un \poli{4}topo cuyo grafo es el producto de~$H$ por un segmento (derecha).}
	\label{E:fig:truncatedoctahedron}
\end{figure}

\enlargethispage{-.3cm}
Preguntamos después cuándo un producto de grafos regulares no-politopales puede ser politopal en una dimensión más pequeña que su grado. Los ejemplos siguientes responden parcialmente a esta pregunta:

\begin{theorem_E}
\begin{enumerate}[(i)]
\item Cuando~$n\ge 3$, el producto~$K_{n,n}\times K_2$ de un grafo completo bipartito por un segmento no es politopal.
\item El producto de un grafo \poli{d}topal por el grafo de una subdivisión regular de un \poli{e}\-topo es \poli{(d+e)}topal. Esto proporciona productos politopales de grafos regulares no-politopales (por ejemplo el producto de dos dominós de la Figura~\ref{E:fig:products} y el grafo de la Figura~\ref{E:fig:truncatedoctahedron}).\qed
\end{enumerate}
\end{theorem_E}

\subsection*{Politopos \emph{neighborly} prodsimpliciales}

En el Capítulo~\ref{chap:psn}, consideramos politopos cuyos esqueletos son los de un producto de símplices:

\begin{definition_E}
Sean~$k\ge0$ y~$\sub{n} \eqdef (n_1,\dots,n_r)$, con~$r\ge1$ y~$n_i\ge 1$ por cada~$i$. Un politopo es \defn{$(k,\sub{n})$-neighborly prodsimplicial} ---~o para acotar \defn{\xpsn{(k,\sub{n})}} ---~si su \esqueleto{k} es combinatorialmente equivalente al del producto de símplices~$\simplex_{\sub{n}} \eqdef \simplex_{n_1}\times\cdots\times\simplex_{n_r}$.
\end{definition_E}

Esta definición generaliza dos clases particulares de politopos:
\begin{enumerate}[(i)]
\item Encontramos los politopos \defn{neighborly} cuando~$r=1$. En la literatura, un politopo es \defn{\neighborly{k}} si cada subconjunto de~$k$ vértices forma una cara. Con nuestra definición, un tal politopo es \xpsn{(k-1,n)}.
\item Encontramos los politopos \defn{neighborly cúbicos}~\cite{jz-ncp-00,js-ncps-07,sz-capdp} cuando~$\sub{n}=(1,\dots,1)$.
\end{enumerate}

\svs
El producto~$\simplex_{\sub{n}}$ es un politopo \xpsn{(k,\sub{n})} de dimensión~$\sum n_i$. Buscamos naturalmente politopos \xpsn{(k,\sub{n})} en dimensión inferior. Por ejemplo, el politopo cíclico~$C_{2k+2}(n+1)$ es un politopo \xpsn{(k,n)} en dimensión~$2k+2$. Denotamos~$\delta(k,\sub{n})$ la dimensión mínima que puede tener un politopo \xpsn{(k,\sub{n})}.

Ciertos politopos \psn se obtienen proyectando el producto~$\simplex_{\sub{n}}$, o cualquier politopo combinatorialmente equivalente, sobre un subespacio de dimensión inferior. Por ejemplo, el politopo cíclico~$C_{2k+2}(n+1)$, como cualquier politopo con~$n+1$ vértices, es una proyección del símplice~$\simplex_n$ sobre~$\R^{2k+2}$.

\begin{definition_E}
Un politopo \xpsn{(k,\sub{n})} es \defn{$(k,\sub{n})$-neighborly prodsimplicial proyectado} ---~o para acotar \defn{\xppsn{(k,\sub{n})}} ---~si es una proyección de un politopo combinatorialmente equivalente a~$\simplex_{\sub{n}}$.
\end{definition_E}

Denotamos~$\delta_{pr}(k,\sub{n})$ la dimensión mínima de un politopo \xppsn{(k,\sub{n})}.

\mvs
Nuestro Capítulo~\ref{chap:psn} se divide en dos partes. En la primera, presentamos tres métodos para construir politopos \ppsn en pequeña dimensión:
\begin{enumerate}[(i)]
\item con reflexiones de politopos cíclicos;
\item con sumas de Minkowski de politopos cíclicos;
\item con ``proyecciones de productos deformados'', en el espíritu de las construcciones de Raman Sanyal y G\"unter Ziegler~\cite{z-ppp-04,sz-capdp}.
\end{enumerate}
En la segunda parte, obtenemos obstrucciones topológicas a la existencia de tal objetos, utilizando técnicas desarrolladas por Raman Sanyal~\cite{s-tovnms-09} para acotar el número de vértices de una suma de Minkowski. A la luz de estas obstrucciones, nuestras construcciones de la primera parte se revelan optimas sobre un ancho espectro de parámetros.

\paragraph{Construcciones.}
Nuestro primer ejemplo no-trivial es un politopo \xpsn{(k,(1,n))} en dimensión~$2k+2$, obtenido por reflexión del politopo cíclico~$C_{2k+2}(n+1)$ respecto a un hiperplano bien elegido:

\begin{proposition_E}
Para todos~$k\ge0$,~$n\ge2k+2$ y~$\lambda\in\R$ suficientemente grande, el politopo
$$\conv\left(\ens{(t_i,\dots,t_i^{2k+2})^T}{i\in[n+1]} \cup \ens{(t_i,\dots,t_i^{2k+1},\lambda-t_i^{2k+2})^T}{i\in[n+1]}\right)$$
es un politopo \xpsn{(k,(1,n))} de dimensión~$2k+2$.\qed
\end{proposition_E}

Por ejemplo, esta construcción produce un \poli{4}topo cuyo grafo es~$K_2\times K_n$~($n\ge3$).

\svs
Después, con ayuda de sumas de Minkowski bien elegidas de politopos cíclicos, obtenemos coordenadas explícitas de politopos \xppsn{(k,\sub{n})}:

\begin{theorem_E}\label{E:theo:UBminkowskiCyclic}
Sean~$k\ge0$ y~$\sub{n} \eqdef (n_1,\dots,n_r)$ con~$r\ge1$ y~$n_i\ge1$ para todo~$i$. Existen conjuntos~$I_1,\dots,I_r\subset\R$, con~$|I_i|=n_i$ para todo~$i$, tales que el politopo
$$\conv\ens{w_{a_1,\dots,a_r}}{(a_1,\dots,a_r)\in I_1\times\cdots\times I_r} \subset \R^{2k+r+1}$$
sea \xppsn{(k,\sub{n})}, donde $w_{a_1,\dots,a_r} \eqdef \big(a_1,\dots,a_r,\sum_{i\in[r]}a_i^2,\dots,\sum_{i\in[r]} a_i^{2k+2}\big)^T$. Por consiguiente,
$$\delta(k,\sub{n}) \le \delta_{pr}(k,\sub{n}) \le 2k+r+1.$$
\end{theorem_E}
\vspace{-.8cm}\qed
\vspace{.4cm}

Cuando~$r=1$ encontramos los politopos \emph{neighborly}. 

\svs
Finalmente, extendemos la técnica de ``proyección de productos deformados de polígonos'' de Raman Sanyal y G\"unter Ziegler~\cite{z-ppp-04,sz-capdp} a productos de politopos simples arbitrarios: proyectamos un politopo bien elegido, combinatorialmente equivalente a un producto de politopos simples, de tal manera que preservamos su \esqueleto{k} completo. Más concretamente, describimos como usar coloraciones de los grafos de los politopos polares de los factores del producto para aumentar la dimensión del esqueleto preservado. La versión básica de esta técnica produce el resultado siguiente:

\begin{proposition_E}
Sean~$P_1,\dots,P_r$ unos politopos simples. Para cada politopo~$P_i$, denotamos~$n_i$ su dimensión,~$m_i$ su número de facetas, y~$\chi_i \eqdef \chi(\gr(P_i^\polar))$ el número cromático del grafo del politopo polar~$P_i^\polar$. Para un entero fijo~$d\le n$, sea~$t$ el entero máximo tal que~${\sum_{i=1}^t n_i\le d}$. Entonces existe un \poli{d}topo cuyo \esqueleto{k} es combinatorialmente equivalente al del producto ${P_1\times\cdots\times P_r}$ en cuanto
$$0 \le k \le \sum_{i=1}^r (n_i-m_i) + \sum_{i=1}^t (m_i-\chi_i) + \Floor{\frac{1}{2}\left(d-1+\sum_{i=1}^t(\chi_i-n_i)\right)}.$$
\end{proposition_E}
\vspace{-1.2cm}\qed
\vspace{1cm}

\svs
Especializando esta proposición a un producto de símplices, obtenemos otra construcción de politopos \ppsn. Cuando algunos de los símplices son pequeños respecto a~$k$, esta técnica produce de hecho nuestro mejores ejemplos de politopos \ppsn:

\begin{theorem_E}
Para todos~$k\ge0$ y~$\sub{n} \eqdef (n_1,\dots,n_r)$ con~${1=n_1=\cdots=n_s<n_{s+1}\le\cdots\le n_r}$,
$$\delta_{pr}(k,\sub{n}) \le
\begin{cases}
     2(k+r)-s-t & \text{si } 3s \le 2k+2r, \\
     2(k+r-s)+1 & \text{si } 3s = 2k+2r+1, \\
     2(k+r-s+1) & \text{si } 3s \ge 2k+2r+2,
\end{cases}$$
donde~$t\in\{s,\dots,r\}$ es máximo tal que~$3s+\sum_{i=s+1}^{t}(n_i+1) \le 2k+2r$.\qed
\end{theorem_E}

Si~$n_i=1$ para todo~$i$, encontramos los politopos \emph{neighborly} cúbicos de~\cite{sz-capdp}.

\paragraph{Obstrucciones.}
Para obtener cotas inferiores sobre la dimensión mínima~$\delta_{pr}(k,\sub{n})$ que puede tener un politopo \xppsn{(k,\sub{n})}, aplicamos una método de Raman Sanyal~\cite{s-tovnms-09}. A cada proyección que preserve el \esqueleto{k} de~$\simplex_{\sub{n}}$, asociamos por dualidad de Gale un complejo simplicial que debe ser encajable en un espacio de una cierta dimensión. El argumento se apoya después en una obstrucción topológica derivada del criterio de Sarkaria para el encaje de un complejo simplicial en términos de coloraciones de grafos de Kneser~\cite{m-ubut-03}. Obtenemos el resultado siguiente:

\begin{theorem_E}\label{E:theo:topObstr}
Sea~$\sub{n} \eqdef (n_1,\dots,n_r)$ con~$1=n_1=\cdots=n_s<n_{s+1}\le\cdots\le
n_r$.
\begin{enumerate}[(i)]
\item Si
$$0 \le k \le \sum_{i=s+1}^r \Fracfloor{n_i-2}{2} + \max\left\{0,\Fracfloor{s-1}{2}\right\},$$
entonces~$\delta_{pr}(k,\sub{n}) \ge 2k+r-s+1$.
\item Si~$k\ge \Floor{\frac{1}{2} \sum_i n_i}$ entonces~$\delta_{pr}(k,\sub{n}) \ge \sum_i n_i$.\qed
\end{enumerate}
\end{theorem_E}

En particular, las cotas inferior y superior de los Teoremas~\ref{E:theo:UBminkowskiCyclic} y~\ref{E:theo:topObstr} se juntan sobre un largo campo de parámetros:

\begin{theorem_E}\label{E:theo:mainResult}
Para todos~$\sub{n} \eqdef (n_1,\dots,n_r)$ con~$r\ge1$ y~$n_i\ge2$ para cada~$i$, y para cada~$k$ tal que~$0\le k\le \sum_{i\in [r]}\Fracfloor{n_i-2}{2}$, el politopo \xppsn{(k,\sub{n})} mínimo es de dimensión exactamente $2k+r+1$. Dicho de otra manera:
$$\delta_{pr}(k,\sub{n}) = 2k+r+1.$$
\end{theorem_E}
\vspace{-.8cm}\qed

\begin{remark_E}
Las técnicas de proyección de politopos y las obstrucciones que utilizamos han sido desarrolladas por Raman Sanyal y G\"unter Ziegler \cite{z-ppp-04,sz-capdp,s-tovnms-09}. Decidimos presentarlas en esta memoria porque sus aplicaciones a productos de símplices producen resultados nuevos que completan nuestro estudio sobre la politopalidad de productos. Por otra parte, después de que acabamos de estudiar estos métodos a productos de símplices, descubrimos que Thilo~R\"orig y Raman~Sanyal trabajaron sobre el mismo asunto~\cite{rs-npps} (ver también~\cite{sanyal-phd,rorig-phd}).
\end{remark_E}



\cleardoublepage
\addtocontents{toc}{ \vspace{.5cm} \demiligne \vspace{.1cm} }
\addtocontents{lof}{ \vspace{.5cm} \demiligne \vspace{.1cm} }


\vspace*{0 pt plus 1 fill}
\ligne 
\begin{center}
\fontseries{m}\fontshape{sc}\fontsize{35}{40}\selectfont
\darkblue
Multitriangulations, pseudotriangulations \\ and some problems of realization of polytopes
\end{center}
\ligne
\vspace*{0 pt plus 1 fill}
\begin{center}\includegraphics[width=\textwidth]{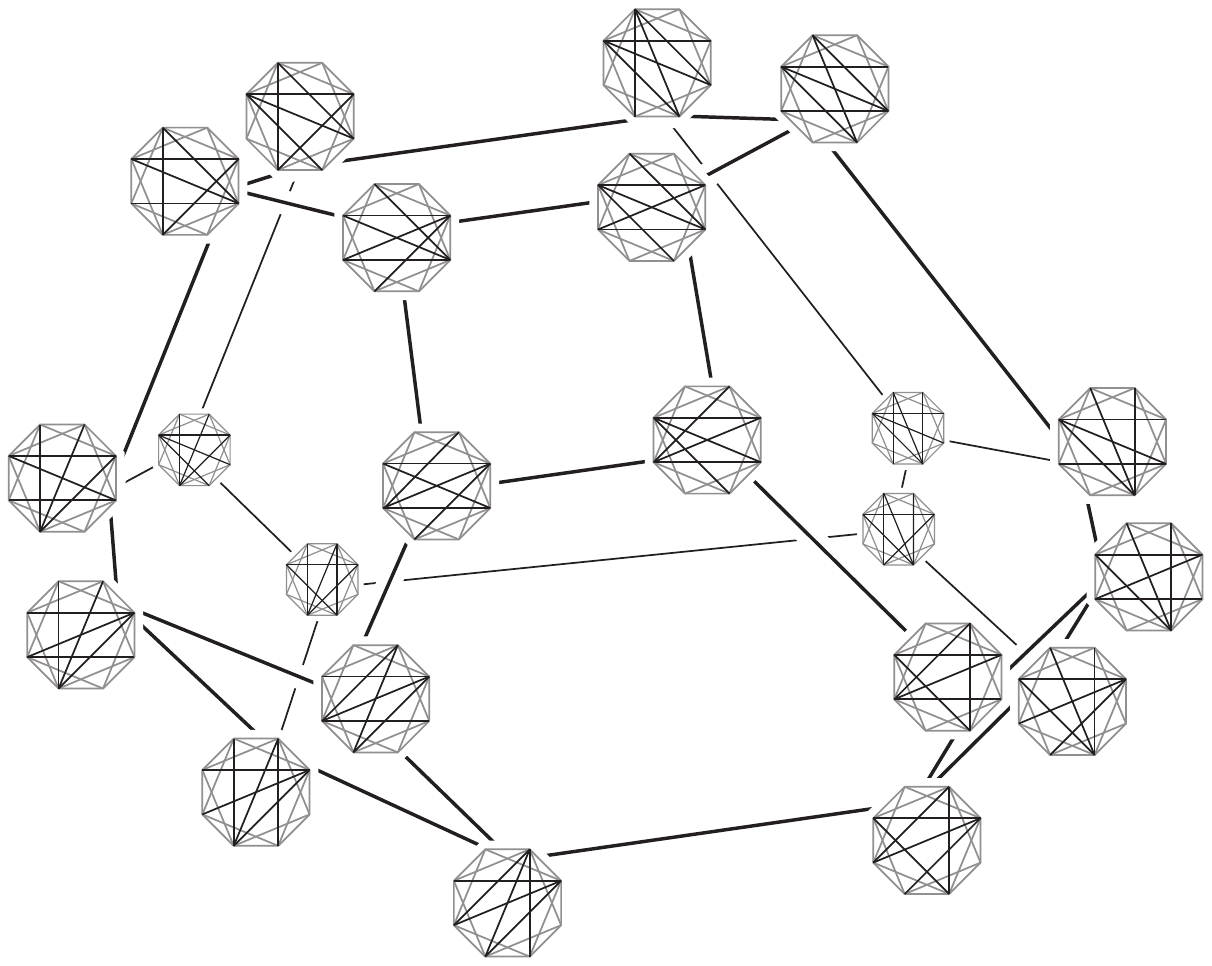}\end{center}

\thispagestyle{empty}
\addtocontents{toc}{\vspace*{.5cm}\centerline{\darkblue\Large\textsc{Multitriangulations, pseudotriangulations}}\par}
\addtocontents{toc}{\centerline{\darkblue\Large\textsc{and some problems of realization of polytopes}}\par}


\setcounter{chapter}{-3}
\selectlanguage{english}
\renewcommand{\namepart}{Introduction}
\chapter*{Introduction}\label{chap:introduction}
\addcontentsline{toc}{chapter}{Introduction}

This thesis explores subjects in the field of \defn{discrete and computational geometry}~\cite{hdcg-04}; it typically studies the ways in which simple geometric objects (such as points, lines, half-spaces,~\etc) are related, arranged, or placed with respect to each other (whether they intersect, how they see each other, \etc). The two specific topics of \defn{multitriangulations} and \defn{polytopal \mbox{realizations} of products} are closely examined. Their various connections with discrete geometry enable us to uncover some of its multiple facets, examples of which are provided in this introduction. We studied these two topics with respect to a common problem, the study of the existence of a polytopal realization of a given structure.

\svs
\index{polytope}
A (convex) polytope is the convex hull of a finite point set of an Euclidean space. While interest in certain polytopes dates back to antiquity (Platonic solids), their systematic study is relatively recent and the main results only appeared last century (see~\cite{g-cp-03,z-lp-95} and the references therein). The study of polytopes deals not only with their geometric properties, but mainly with more combinatorial aspects. The goal is to understand their faces (their intersections with a supporting hyperplane) and the lattice they form (\ie the inclusion relations between these~faces).

\index{polytopality|hbf}
\index{polytopal|hbf}
\index{realization!polytopal ---|hbf}
\index{Steinitz' Theorem}
The questions of \defn{polytopal realization} are in a way the inverse problem: they deal with the existence and the construction of polytopes with a prescribed combinatorial structure. For example, given a graph, we would like to determine if it is the graph of a polytope: we then say that it is polytopal. Steinitz' Theorem~\cite{s-pr-22}, which characterizes the graphs of \poly{3}topes, provides the fundamental result on this question. Already in dimension~$4$, the situation is much less satisfactory: except certain necessary conditions~\cite{b-gscps-61,k-ppg-64,b-ncp-67}, polytopal graphs do not allow for local characterization in general dimension~\cite{rg-rsp-96}. When a graph is polytopal, we are interested in the properties of its realizations, for example the number of faces, the dimension,~\etc{} We often try to construct examples which optimize some of these properties; typically, we want to construct a polytope of minimal dimension for a given graph~\cite{g-ncp-63,jz-ncp-00,sz-capdp}. These questions of polytopal realizations are interesting either for transformation graphs on combinatorial or geometric objects (for example, the graph of adjacent transpositions on permutations, or the graph of flips on triangulations, \etc) or for graphs derived from operations on other graphs (which may be local such as the $\Delta Y$ transformation, or global such as the Cartesian product). These questions are not only interesting for a graph, but more generally for any subset of a lattice.


\paragraph{Polytopality of flip graphs.} The existence of polytopal realizations is studied for transformation graphs on combinatorial or geometric structures. We can cite as an example the permutahedron, whose vertices correspond to permutations of~$[n]$, and whose edges correspond to pairs of permutations which differ by an adjacent transposition. Other examples, as well as classes of polytopes providing realizations of certain general structures are presented in~\cite[Lecture~9]{z-lp-95}. In general, polytopality questions on combinatorial structures are not only interesting in and of themself, but also because studying them forces us to understand the combinatorics of the objects and often leads to the development of new methods and results.

Two particular examples of polytopal combinatorial structures play an important role in this dissertation. We first encounter the associahedron whose boundary complex realizes the dual of the simplicial complex formed by all crossing-free sets of chords of the \gon{n}. The graph of the associahedron corresponds to the graph of flips on triangulations of the \gon{n}. The \mbox{associahedron} appears in several contexts, and various polytopal realizations have been developed \cite{l-atg-89,bfs-ccsp-90,gkz-drmd-94,l-rsp-04,hl-rac-07}. We then meet the polytope of pseudotriangulations of an Euclidean point set~\cite{rss-empppt-03}. Originally introduced in the study of the visibility complex of disjoint obstacles in the plane~\cite{pv-tsvcp-96,pv-vc-96}, pseudotriangulations have since been used in different geometric contexts~\cite{rss-pt-06}. In particular, their interesting rigidity properties~\cite{s-ptrmp-05} yield the construction of a polytope which realizes their flip graph.

\begin{figure}
	\capstart
	\centerline{\includegraphics[width=\textwidth]{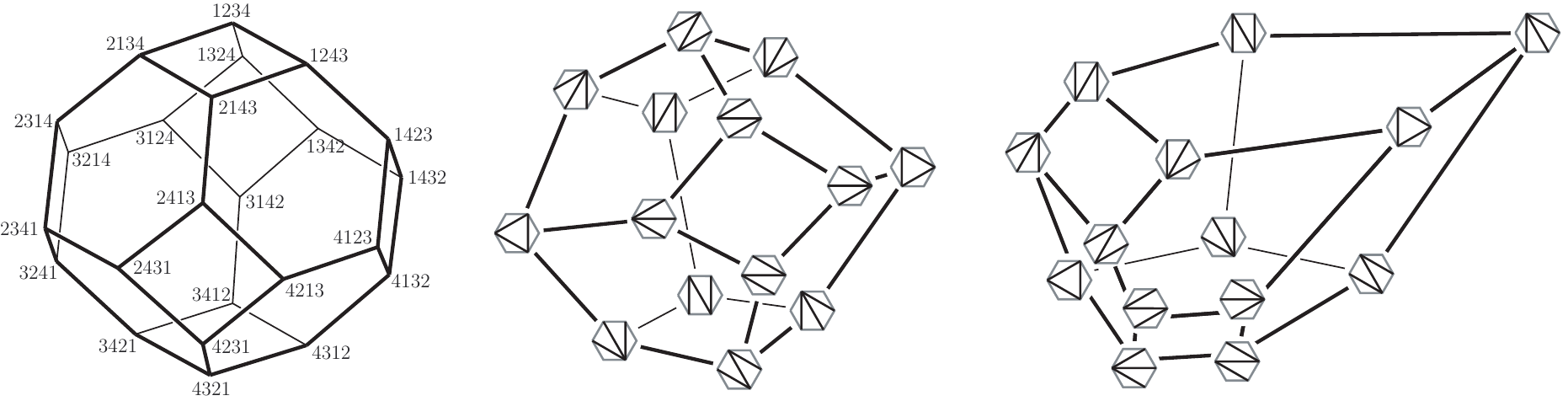}}
	\caption[The permutahedron and two realizations of the associahedron]{The permutahedron and two realizations of the associahedron.}
	\label{I:fig:permutohedronassociahedra}
\end{figure}

\svs
In the first part of this dissertation, we focus on the flip graph on multitriangulations. These objects appeared under different motivations~\cite{cp-tttccp-92,dkm-lahp-02,dkkm-2kn..-01} and revealed a rich combinatorial structure~\cite{n-gdfcp-00,j-gt-03,j-gtdfssp-05}. A \defn{\ktri{k}} is a maximal set of chords of the \gon{n} such that no~$k+1$ of them mutually cross. We consider the flip graph where two multitriangulations are related if they differ by one chord. As for triangulations which occur when~$k=1$, this graph is regular and connected, and we consider the question of its polytopality. Jakob Jonsson~\cite{j-gt-03} made an important step towards answering this question in proving that the simplicial complex formed by sets of chords with no~$k+1$ mutually crossing subsets is a topological sphere. Although we have only partial answers to this question, our study yields surprising results which we present in this dissertation. 

\svs
Several constructions of the associahedron~\cite{bfs-ccsp-90,l-rsp-04} are directly or indirectly based on the triangles of the triangulations. For multitriangulations, no similar elementary object has appeared in previous works. We thus started with the question: what are triangles for multitriangulations? The \defn{stars} that we introduce in Chapter~\ref{chap:stars} address this question. We consider stars to be the right way of approaching multitriangulations, in much the same way that triangles, rather than chords, are the right way of approaching triangulations. As evidence for this, we use them to arrive at new proofs of the basic properties of multitriangulations known to date. First, we study the incidence relations between stars and chords (each internal chord is contained in two stars) which provide a new proof showing that all \ktri{k}s of the \gon{n} have the same cardinality. Then, considering the bisectors of the stars, we provide a local interpretation of the flip operation (an internal chord is replaced by the unique bisector of the two stars adjacent to it), which sheds light upon the study of the flip graph and of its diameter. We also redefine inductive operations on multitriangulations which enable us to insert or suppress a vertex to the \gon{n}. Finally, we use the decomposition of a multitriangulation into stars to interpret it as a polyhedral decomposition of a surface, and we apply this interpretation to the construction of regular decompositions of surfaces.

\begin{figure}
	\capstart
	\centerline{\includegraphics[scale=1]{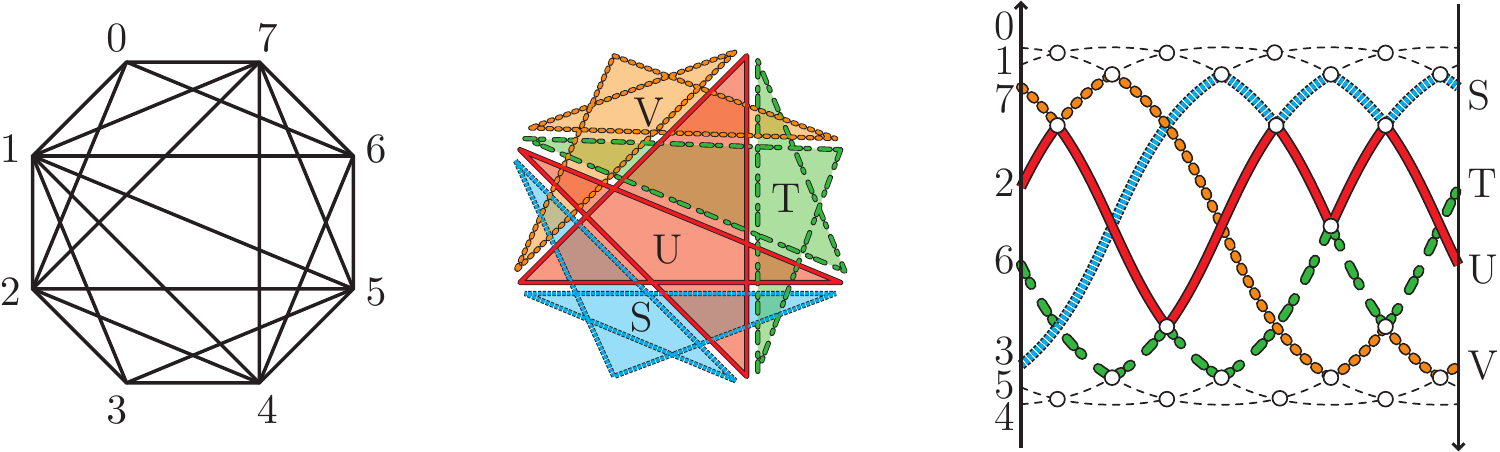}}
	\caption[A \ktri{2}, its decomposition into stars, and its dual arrangement]{A \ktri{2}, its decomposition into stars, and its dual arrangement.}
	\label{I:fig:ktristardual}
\end{figure}

\svs
In Chapter~\ref{chap:mpt}, we consider multitriangulations by duality. The line space of the plane is a M\"obius strip; the set of lines passing through a point of the plane is a pseudoline of the M\"obius strip; and the dual pseudolines to a point  configuration form a pseudoline arrangement~\cite{g-pa-97}.

Our starting point is the duality between the pseudotriangulations of a point set~$P$ and the pseudoline arrangements whose support is the dual arrangement of~$P$ minus its first level. We introduce a similar duality for multitriangulations. On the one hand, the set of bisectors of a star is a pseudoline of the M\"obius strip. On the other hand, the dual pseudolines to the stars of a \ktri{k} of the \gon{n} form a pseudoline arrangement with contact points supported by the dual arrangement of the vertex set of the \gon{n} minus its first~$k$ levels. We prove that any arrangement with contact points with this support is effectively the dual arrangement of a multitriangulation. This duality relates pseudotriangulations and multitriangulations and thus provides an explanation for their common properties (number of chords, flips, \etc).

More generally, we consider \defn{pseudoline arrangements with contact points} which share a common support. We define a flip operation which corresponds to flips in multitriangulations and pseudotriangulations, and we study the flip graph. The properties of certain greedy arrangements, defined as sources of certain acyclic orientations on these graphs, enable us in particular to enumerate this flip graph with a polynomial working space. Our work elucidates the existing enumeration algorithm for pseudotriangulations~\cite{bkps-ceppgfa-06} and provides a complementary proof for it.

\svs
To finish, we discuss in Chapter~\ref{chap:multiassociahedron} three open problems which reflect the combinatorial and geometric richness of multitriangulations.

The first one looks at the enumeration of multitriangulations. Jakob Jonsson~\cite{j-gtdfssp-05} proved (considering $0/1$-fillings of polyominoes avoiding certain patterns) that the multitriangulations are counted by a Hankel deteminant of Catalan numbers, which also counts certain families of \tuple{k}s of Dyck paths. However, apart from partial results~\cite{e-btdp-07,n-abtdp-09}, no bijective proof of this result is known. We enter here the field of \defn{bijective combinatorics} which aims to construct bijections between combinatorial families which preserve characteristic parameters; here, we would like to find a bijection which would enable us to read the stars on the \tuple{k}s of Dyck paths.

The second problem is that of \defn{rigidity}. Although the rigidity of \dimensional{2} graphs is well understood~\cite{l-grpss-70,g-cf-01,gss-cr-93}, no satisfactory characterization is known in dimension~$3$ or higher. We prove that \ktri{k}s satisfy some typical properties of rigid graphs in dimension~$2k$. This leads to conjecture that \ktri{k}s are rigid in dimension~$2k$, which we prove when~$k=2$. A positive answer to this conjecture would provide insight to the polytopality question of the flip graph on multitriangulations, in much the same way as the polytope of pseudotriangulations~\cite{rss-empppt-03} relies on their rigidity properties.

Finally, we come back to the \defn{polytopal realization of the flip graph} on multitriangulations. We study the first non-trivial example and we prove that the graph of flips on the \ktri{2}s of the octagon is the graph of a \dimensional{6} polytope. To find such a polytope, we completely describe its symmetric polytopal realizations in dimension~$6$, studying first all the symmetric oriented matroids~\cite{bvswz-om-99,b-com-06} which can realize it. Then generalizing Jean-Louis Loday's construction of the associahedron~\cite{l-rsp-04}, we construct a polytope which realizes the graph of flips restricted to multitriangulations whose oriented dual graph is acyclic.

\svs
In addition, we present in Appendix~\ref{app:implementations} the results of an enumeration algorithm for small arrangements of pseudolines and \defn{double pseudolines}. Much like pseudoline arrangements which model point configurations, double pseudoline arrangements have been introduced as a combinatorial model of configurations of disjoint convex bodies~\cite{hp-adp-08}. Our work of implementation enabled us to manipulate these objects and to become familiar with their properties, which proved to be useful for our work on duality.


\paragraph{Polytopality of Cartesian products.} The \defn{Cartesian product} of graphs is defined in such a way to be coherent with products of polytopes: the graph of a product of polytopes is the product of their graphs. Thus, a product of polytopal graphs is automatically polytopal. We first study the reciprocal question: does the polytopality of a product of graphs imply the polytopality of its factors? Chapter~\ref{chap:nonpolytopal} addresses this question, paying particular attention to regular graphs and their realization as simple polytopes (simple polytopes have  special properties, including the fact that they are determined by their graph~\cite{bm-ppi-87,k-swtsp-88}). We discuss in particular the polytopality of the product of two Petersen graphs, which was raised by G\"unter Ziegler~\cite{crm}. This work also led us to study examples of non-polytopal graphs which satisfy all known necessary conditions to be polytopal~\cite{b-gscps-61,k-ppg-64,b-ncp-67}.

\begin{figure}
	\capstart
	\centerline{\includegraphics[scale=.92]{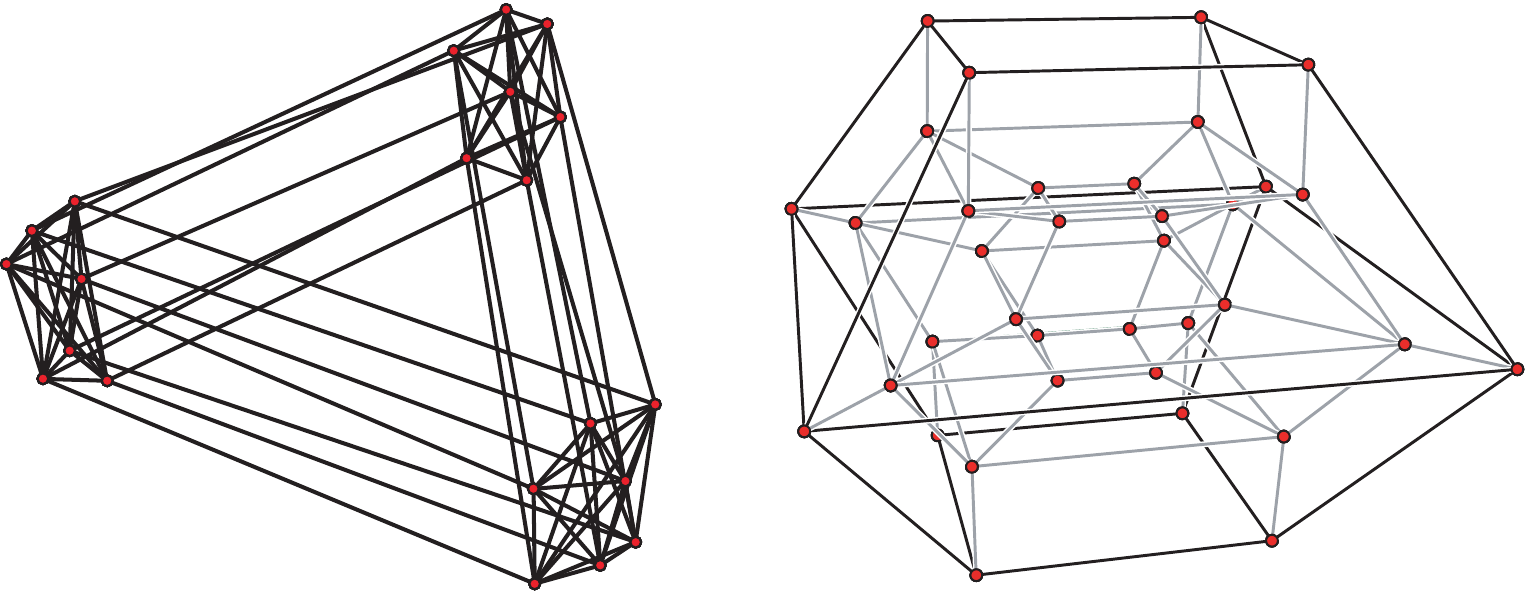}}
	\caption[Examples of products of graphs]{(Left) A product of complete graphs. This graph is that of a product of simplices, but also that of polytopes in smaller dimension. (Right) A polytopal product of non-polytopal graphs.}
	\label{I:fig:products}
\end{figure}

\svs
Second, we look for the minimal dimension of a polytope whose graph is isomorhpic to that of a fixed product of polytopes. This question is an example of the systematic research of extremal dimensions of polytopes with given properties. For example, neighborly polytopes~\cite{g-ncp-63} are those whose \skeleton{\Fracfloor{d}{2}} is complete (and they reach the maximal number of faces allowed by the Upper Bound Theorem~\cite{m-mnfcp-70}); neighborly cubical polytopes~\cite{jz-ncp-00,z-ppp-04,js-ncps-07} are those whose \skeleton{\Fracfloor{d}{2}} is that of the $n$-cube, \etc{} To construct polytopes in small dimension with a given property, it is natural and often efficient to start from a sufficiently large dimension to guarantee the existence of such polytopes, and to project these polytopes on low dimensional subspaces in such a way as to preserve the desired property. These techniques and their limits have been widely studied in the literature~\cite{az-dpmsp-99,z-ppp-04,sz-capdp,s-tovnms-09,sanyal-phd}. We apply them in Chapter~\ref{chap:psn} to construct and study the minimal possible dimension of \defn{\mbox{$(k,\sub{n})$-prodsimplicial neighborly}} polytopes whose \skeleton{k} is that of the product of simplices ${\simplex_{\sub{n}} \eqdef \simplex_{n_1}\times\cdots\times\simplex_{n_r}}$. In this chapter, we also provide explicit integer coordinates of such polytopes using suitable Minkowski sums of cyclic polytopes.

\newpage

\setcounter{chapter}{0}
\renewcommand*{\figurename}{Figure}
\makeatletter\renewcommand{\thefigure}{\thechapter.\@arabic\c@figure}\makeatother
\pagestyle{hautpage}
\counterwithin{figure}{chapter}
\counterwithin{section}{chapter}
\counterwithin{subsection}{chapter}
\counterwithin{subsection}{section}


\renewcommand{\partfigure}{irreducible}
\renewcommand{\namepart}{Multitriangulations}
\part{Multitriangulations}
	\chapter{Introduction}\label{intro:chap:introduction_multitriangulations}


\section{On triangulations\ldots~of a convex polygon}\label{intro:sec:triangulations}

Fix~$n$ points on the unit circle (the vertices of a convex \gon{n}), and consider the set of~${n \choose 2}$ straight chords between them. We say that two such chords \defn{cross}\index{crossing} if they have an intersection point which is not one of their vertices, and we consider sets of non-crossing chords of the \gon{n}. Among these sets, we focus on those that are maximal, which of course contain all the boundary edges of the \gon{n} and decompose it into triangles:

\begin{definition}
\index{triangulation}
A \defn{triangulation} of the convex \gon{n} is:
\begin{enumerate}[(i)]
\item a maximal crossing-free subset of edges of the \gon{n}; or (equivalently)
\item a set of non-overlapping triangles covering the \gon{n}.
\end{enumerate}
\end{definition}

\begin{figure}
	\capstart
	\centerline{\includegraphics[scale=1]{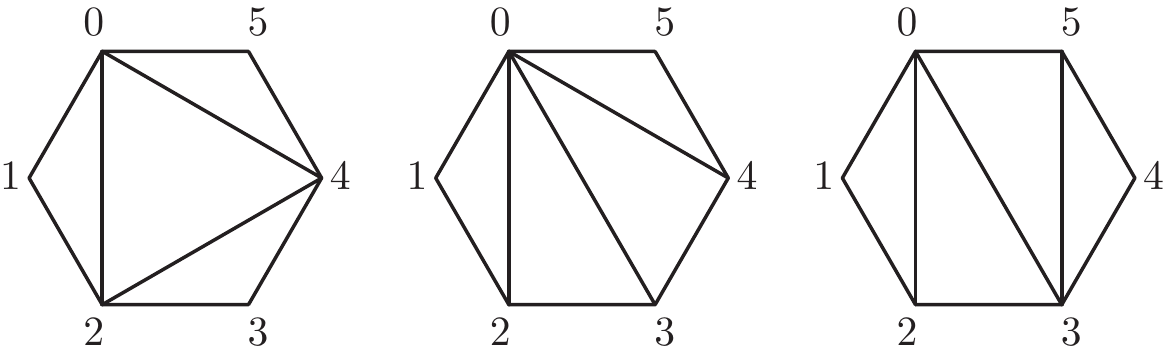}}
	\caption[Three triangulations of the hexagon]{Three triangulations of the hexagon.}
	\label{intro:fig:trianghexagon}
\end{figure}

For example, \fref{intro:fig:trianghexagon} represents some triangulations of a convex hexagon (in fact, all of them modulo rotations and reflections). Observe that they all have the same number of triangles (namely~$4$) and of edges (namely~$9$), which can be generalized for any~$n$:

\begin{lemma}\label{intro:lem:2n-3}
A triangulation of the \gon{n} has exactly~$n-2$ triangles and~$2n-3$ edges.
\end{lemma}

This simple fact can be proved with many different approaches: we use it as an excuse to briefly recall some of the multiple facets of triangulations and to give some flavor of the surprisingly rich combinatorial structure of these elementary objects.

\paragraph{Double counting.}

Our first proof is probably the most direct one since it only involves a single triangulation. Indeed, Lemma~\ref{intro:lem:2n-3} is easily derived from any two of the three following equalities between the number~$e$ of edges and the number~$t$ of triangles of a triangulation~$T$ of the \gon{n}:
\begin{enumerate}[(i)]
\item The first relation is derived from a \defn{double count} of the number of incidences between edges and triangles: a triangle contains three edges and an edge is contained in two triangles (except for boundary edges), which yields the equality~$3t=2e-n$.
\item The second one is an application of \defn{Euler's formula} (which relates the number of faces, edges and vertices of a polygonal decomposition of the sphere): $(t+1)-e+n=2$.
\item The last one is based on bisectors: a \defn{bisector}\index{bisector} of a triangle is any edge which bisects one of its angles. Any two triangles of~$T$ have a unique common bisector, while any edge not in~$T$ is the common bisector of a unique pair of triangles of~$T$. This provides a bijection between pairs of distinct triangles of~$T$ and edges not in~$T$ and yields~${t \choose 2} = {n \choose 2} - e$.
\end{enumerate}

\paragraph{Induction.}

Our second proof is a simple recurrence on the number~$n$ of vertices: in a triangulation~$T$ of the \gon{(n+1)}, when we  flatten a triangle which contains a boundary edge, we obtain a triangulation of the \gon{n} with one less triangle (the flattened one) and two less edges (in the flattened triangle, the boundary edge is contracted and the two other ones are merged). Lemma~\ref{intro:lem:2n-3} immediately follows.

This flattening operation can also be used to compute the number~$\theta(n)$ of triangulations of the \gon{n}, where the vertices of the \gon{n} are labeled, say from~$0$ to~$n-1$. Indeed, the contraction of a distinguished boundary edge, say~$[0,n]$, defines a surjective map from triangulations of the \gon{(n+1)} to triangulations of the \gon{n}, and the number of preimages of a triangulation of the \gon{n} is simply the degree of its merged vertex (call it~$\bar 0$). Thus, the quotient of~$\theta(n+1)$ by~$\theta(n)$ equals the average degree of vertex~$\bar 0$, which yields the formula~$n\theta(n+1)=(4n-6)\theta(n)$.

There is another way to understand triangulations inductively: any triangulation of the \gon{n} can be decomposed into its root triangle~$\{0,m,n-1\}$ (with $1\le m\le n-2$) together with one triangulation of the \gon{(m+1)}, formed by all the edges on the left of~$[0,m]$, and one triangulation of the \gon{(n-m)}, formed by all the edges on the right of~$[m,n-1]$. Since the left and the right triangulations are independent, the number of triangulations of the \gon{n} can also be computed by the induction formula $\theta(n)=\sum_{m=1}^{n-2} \theta(m+1)\theta(n-m)$.

From any of these two inductions, one can deduce the following closed formula:

\begin{proposition}\label{intro:prop:catalan}
\index{Catalan number}
The number~$\theta(n)$ of triangulations of a convex \gon{n} is the \defn{Catalan number}:
$$\theta(n)=C_{n-2} \eqdef \frac{1}{n-1}{2n-4 \choose n-2}.$$
\end{proposition}
\vspace{-1.05cm}\qed

\mvs
Catalan numbers count several combinatorial ``Catalan families'' such as triangulations, \defn{rooted binary trees}\index{dual!--- binary tree of a triangulation} (trees whose~$n-2$ nodes have either two or no children), or \defn{Dyck paths}\index{Dyck path} (lattice paths which use north steps~$N \eqdef (0,1)$ and east steps~$E \eqdef (1,0)$, start from~$(0,0)$, end at~$(n-2,n-2)$ and never go below the diagonal~$x=y$)~---~see~\cite{s-ec-99} for many other examples. Of course there are bijections from triangulations to any of these families; these bijections can be made explicit either by the inductive decomposition described above, or even directly. For example, if~$T$ is a triangulation of the \gon{n}, then:
\begin{enumerate}[(i)]
\item the dual graph of~$T$ (with one vertex for each triangle and one edge between two adjacent triangles) is a binary tree with~$n-2$ nodes, rooted at the triangle containing the boundary edge~$[0,n-1]$; and
\item the lattice path~$N^{\delta_0(T)}EN^{\delta_1(T)}E\dots N^{\delta_{n-3}(T)}E$, where~$\delta_i(T)$ is the number of triangles of~$T$ whose first vertex is~$i$, is a Dyck path (of semilength~$n-2$).
\end{enumerate}
These two maps, from triangulations to rooted binary trees and to Dyck paths are bijective.

\begin{figure}
	\capstart
	\centerline{\includegraphics[scale=1]{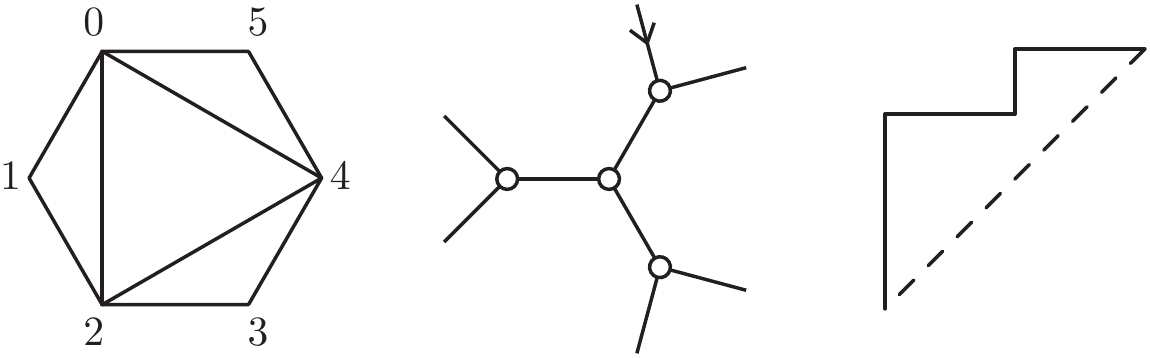}}
	\caption[Catalan families]{Three catalan families: a triangulation, its dual rooted binary tree, and the associated Dyck path.}
	\label{intro:fig:catalan}
\end{figure}

\paragraph{Flips.}

Our last proof of Lemma~\ref{intro:lem:2n-3} illustrates flips in triangulations. A \defn{flip}\index{flip} is a local operation which transforms a triangulation of the \gon{n} into another one: in any quadrangle formed by two adjacent triangles, we can replace one diagonal by the other one. Since a flip does not change neither the number of triangles, nor the number of edges, Lemma~\ref{intro:lem:2n-3} is a consequence of:

\begin{proposition}
\index{flip!graph of ---s}
The graph of flips, whose vertices are triangulations of the \gon{n} and whose edges are flips between them, is connected.\qed
\end{proposition}

This proposition is proved by observing that any triangulation can be transformed into the ``fan triangulation'' with apex~$0$ (\ie the triangulation obtained by relating vertex~$0$ to all the \mbox{others}): indeed, in any other triangulation, there is a flip which increases the degree of vertex~$0$. Observe that this even implies that we can reach the fan triangulation from any other triangulation with $n-1-\deg_0(T)\le n-3$ flips, and therefore, that the diameter of this graph is at most~$2n-6$.

This bound on the diameter can be slightly improved: we can relate any two triangulations of the \gon{n} through a fan triangulation whose apex maximizes the sum of the degrees (in our two triangulations). Since the average sum of degree is $4(2n-3)/n$, this bounds the diameter by~$2n-10-12/n$, and thus, by~$2n-10$ as soon as~$n>12$.

It is much less elementary to prove that this is in fact the exact diameter of the graph of flips:

\begin{theorem}[\cite{stt-rdthg-88}]\label{intro:theo:diameter}
\index{diameter (of the graph of flips)}
The graph of flips on triangulations of the \gon{n} has diameter~$2n-10$ for sufficiently large values of~$n$.\qed
\end{theorem}

We finish this section with a last important property of the graph of flips (see Section~\ref{ft:subsec:multiassociahedron:associahedron}):

\begin{theorem}[\cite{l-atg-89,bfs-ccsp-90,gkz-drmd-94,l-rsp-04,hl-rac-07}]\label{intro:theo:associahedron}
\index{associahedron}
The graph of flips on triangulations of the \gon{n} is the graph of a convex polytope of dimension~$n-3$, called the \defn{associahedron}.\qed
\end{theorem}

\begin{figure}
	\capstart
	\centerline{\includegraphics[scale=1]{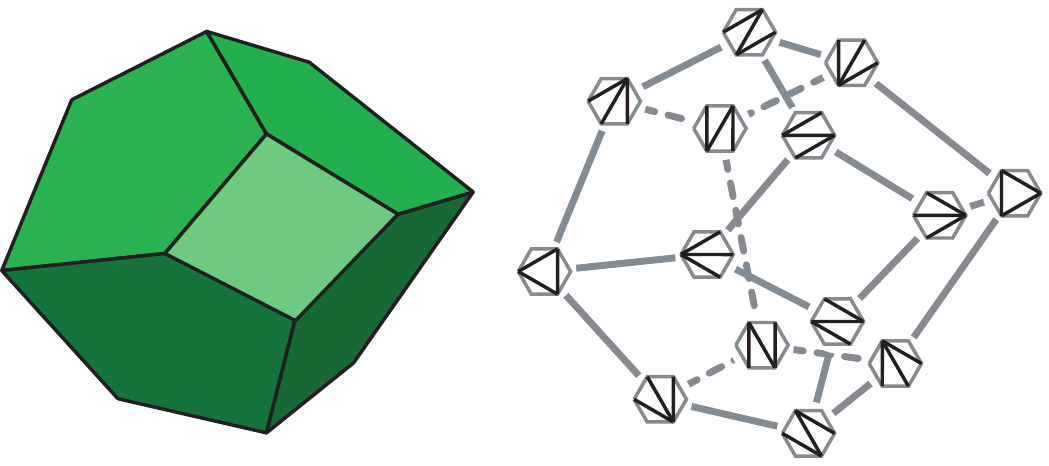}}
	\caption[The \dimensional{3} associahedron]{The \dimensional{3} associahedron, with its vertices corresponding to the~$14$ triangulations of the hexagon, and its edges corresponding to flips between them.}
	\label{intro:fig:associahedronintro}
\end{figure}


\section{Multitriangulations}\label{intro:sec:multitriangulations}

The first topic of this dissertation is a generalization of triangulations, in which some chords are allowed to cross as long as the cardinality of the maximal set of mutually crossing chords remains sufficiently small. To be more precise, for~$\ell\ge 2$, we call \defn{\kcross{\ell}}\index{crossing@\kcross{k}} any set of~$\ell$ mutually intersecting edges of the \gon{n}. As well as we were interested in non-crossing sets of edges, we now consider sets of edges which avoid \kcross{(k+1)}s. Among them, maximal such sets play a special role:

\begin{definition}
\index{triangulation@\ktri{k}}
\index{multitriangulation}
A \defn{\ktri{k}} is a maximal set of edges of the \gon{n} such that no~$k+1$ of them are mutually crossing.
\end{definition}

All throughout our presentation, we will use the \ktri{2} of the following picture to illustrate results on multitriangulations:

\begin{figure}[!h]
	\capstart
	\centerline{\includegraphics[scale=1]{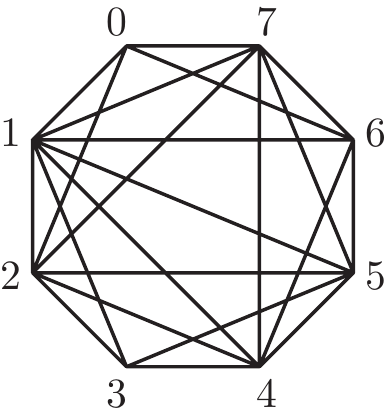}}
	\caption[A \ktri{2} of the octagon]{A \ktri{2} of the octagon.}
	\label{intro:fig:2triang8points}
\end{figure}

Let us insist on the fact that triangulations occur when the parameter~$k$ equals~$1$. Surprisingly, multitriangulations share many combinatorial and structural properties with triangulations. Before reviewing the existing results and summarizing those of this dissertation, we present several interpretations under which multitriangulations have been  motivated and studied.


\subsection{Four different interpretations}\label{intro:subsec:multitriangulations:interpretations}

\paragraph{Clique-free subgraphs of the intersection graph~{\normalfont\cite{cp-tttccp-92}}.}

\begin{figure}[b]
	\capstart
	\centerline{\includegraphics[scale=1]{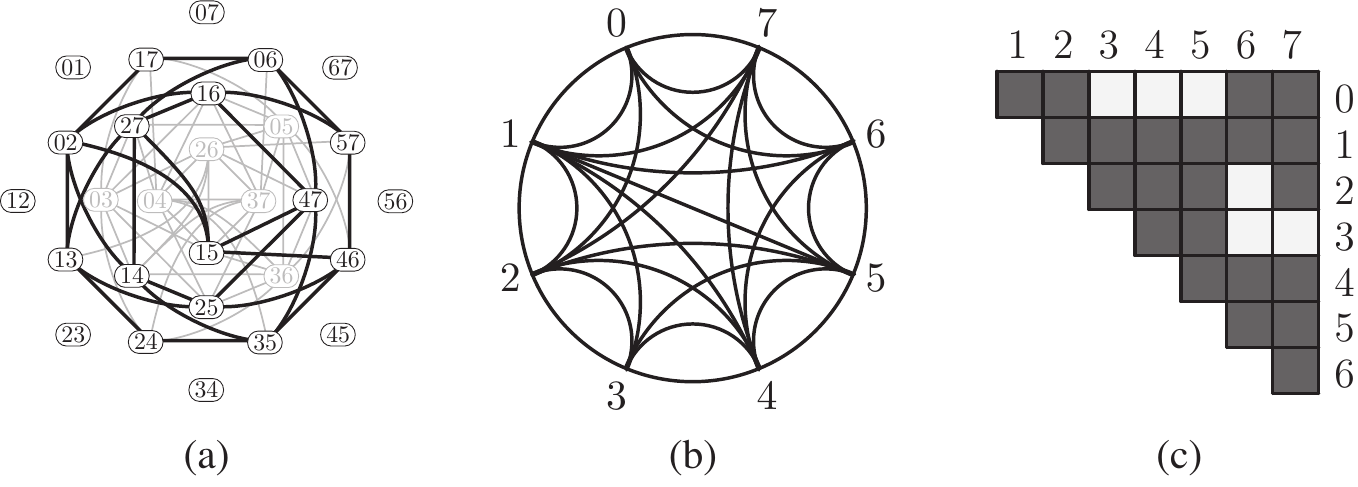}}
	\caption[Three different interpretations of multitriangulations]{Three interpretations of the \ktri{2} of \fref{intro:fig:2triang8points}: (a) A maximal \mbox{triangle-free} induced subgraph of the intersection graph; (b) A line arrangement in the hyperbolic plane (represented on the Poincaré Disk model); (c) A \kdiag{3}-free $0/1$-filling of a triangular polyomino.}
	\label{intro:fig:interpretations}
\end{figure}

Consider the \defn{intersection graph} of the complete geometric graph with vertices in convex position, that is, the graph whose vertices are the straight chords of the \gon{n} and whose edges relate any pair of crossing chords. By definition, an \kcross{\ell} of the \gon{n} translates into an $\ell$-clique in the intersection graph. Consequently, \ktri{k}s are those sets of chords which induce a subgraph of the intersection graph with no~$(k+1)$-cliques (see \fref{intro:fig:interpretations}(a)).

With this formulation~\cite{cp-tttccp-92}, multitriangulations are thought of as optimal graphs for a Tur\'an-type Theorem on chords of the \gon{n}. Remember indeed Tur\'an's famous Theorem:

\begin{theorem}[Tur\'an~\cite{t-tog-54}]
\index{Turan@Tur\'an's Theorem}
A graph on~$n$ vertices with no $p$-clique cannot have more than $\frac{n^2}{2}(1-\frac{1}{p-1})$ edges.\qed
\end{theorem}

\paragraph{Line arrangements in the hyperbolic plane~{\normalfont\cite{dkm-lahp-02}}.}

\index{hyperbolic plane}\index{plane!hyperbolic ---}
Consider a finite set~$\cL$ of lines in a (Euclidean or hyperbolic) plane. Denote by~$\kappa(\cL)$ the maximal number of mutually crossing lines of~$\cL$ and by~$\nu(\cL)$ the number of points at infinity that coincide with at least one line in~$\cL$. For the Euclidean plane these numbers satisfy~$1\le \kappa(\cL)=\nu(\cL)\le |\cL|$, since two non-crossing lines have the same point at infinity. Moreover, there is no further restriction: a line arrangement in general position satisfies~$\kappa(\cL)=\nu(\cL)=|\cL|$ and adding any line parallel to an existing line does not change~$\kappa(\cL)$ or~$\nu(\cL)$.

The situation for the hyperbolic plane is quite different: the number of lines cannot be arbitrarily large compared to the number of points at infinity. By definition, \ktri{k}s of the \gon{n} are exactly the representations in the Klein model of the hyperbolic plane of the maximal (for inclusion) hyperbolic line arrangements whose number of mutually crossing lines does not exceed~$k$ and whose number of points at infinity is~$n$. In other words, while understanding the relations between~$k$,~$n$ and the number of edges in \ktri{k}s of the \gon{n}, we will equivalently be dealing with~$\kappa(\cL)$ and~$\nu(\cL)$ in a hyperbolic line arrangement (see \fref{intro:fig:interpretations}(b)).

\paragraph{Diagonal-free fillings of polyominoes~{\normalfont\cite{j-gtdfssp-05}}.}

\index{polyomino}
Given a subset~$E$ of edges of the \gon{n}, consider the upper triangular part of its adjacency matrix (whose coefficient in the $i$th row and $j$th column is~$1$ if the edge~$[i,j]$ is in~$T$ and~$0$ otherwise, for all~$0\le i<j\le n-1$). Two edges~$[a,b]$ and~$[c,d]$ of~$E$ cross whenever the two corresponding $1$'s in the adjacency matrix form a strictly decreasing diagonal (one~$1$ is below and to the right of the other) whose enclosing rectangle is contained in the upper triangular part of the adjacency matrix. Similarly, an \kcross{\ell} in~$E$ translates in the adjacency matrix into a strictly decreasing \kdiag{\ell} whose enclosing rectangle is contained in the upper part of the adjacency matrix (see \fref{intro:fig:interpretations}(c)).

This replaces multitriangulations in the broader context of the study of $0/1$-fillings of polyominoes (\ie connected subsets of the integer lattice~$\Z^2$) with restrictions on the length on their longest increasing and decreasing diagonals. This interpretation provides in particular interesting enumeration results on multitriangulations (see the discussion in Section~\ref{ft:sec:Dyckpaths}).

\paragraph{Split systems~{\normalfont\cite{dkm-4n-10-04,dkkm-2kn..-01}}.}

A \defn{split} of a finite set~$X$ is a bipartition~$A\sqcup B=X$, and a \defn{split system} on~$X$ is a set of splits of~$X$. Call two splits~$A\sqcup B$ and~$A'\sqcup B'$ \defn{compatible} if at least one of the intersections $A\cap A'$, $A\cap B'$, $B\cap A'$, $B\cap B'$ is empty. Finally a split system is \defn{\kcompatible{k}} if it does not contain a subset of $k+1$ pairwise incompatible splits.

Multitriangulations correspond to the very special ``\kcompatible{k} and cyclic'' split systems: a split system~$S$ on~$X$ is \defn{cyclic} if there exists a bijection~$\phi:\{0,\dots,|X|-1\}\to X$ such that any split of~$S$ is of the form~$\phi([i,j))\sqcup\phi([i,j)^C)$ for some $0\le i<j\le n-1$. In this case, there is an obvious correspondence between a split~$\phi([i,j))\sqcup\phi([i,j)^C)$ and the chord~$[i,j]$ of the \gon{|X|}, in which incompatible splits correspond to crossings. Thus, \kcompatible{k} cyclic split systems on~$X$ are nothing else but \ktri{k}s of the \gon{|X|}.

In fact, when $k\le 2$, the only \kcompatible{k} split systems with maximal cardinality are cyclic split systems~\cite{dkkm-2kn..-01}. This is not true anymore as soon as~$k=3$ (for example the split system~$\ens{([i,j),[i,j)^C)}{0\le i<j\le 3}\cup\{(\{0,2\},\{1,3\})\}$ is \kcompatible{3}).


\subsection{Previous related results}\label{intro:subsec:multitriangulations:known}

As far as we know, multitriangulations first appeared in the work of Vasilis~Capoyleas and J\'anos~Pach in the context of extremal theory for geometric graphs (see~\cite[Chapter~14]{pa-cg-95}, \cite[Chapter~1]{f-gga-04} and the discussion in~\cite{cp-tttccp-92}). They proved that a \kcross{(k+1)}-free subset of edges of the \gon{n} cannot have more than~$k(2n-2k-1)$ edges. Tomoki \mbox{Nakamigawa}~\cite{n-gdfcp-00}, and independently Andreas Dress, Jacobus Koolen and Vincent~Moulton~\cite{dkm-lahp-02}, then proved that all \ktri{k}s of the \gon{n} reach this bound. Both proofs rely on the concept of flips in multitriangulations: as for triangulations, a flip creates a \ktri{k} from another one, removing and inserting a single edge. Tomoki~Nakamigawa~\cite{n-gdfcp-00} proved that every (sufficiently long) edge of a \ktri{k} can be flipped, and that the graph of flips is connected. Let us observe that in these previous works, the results are obtained in an ``indirect way'': first, a flattening operation, similar to the contraction of a boundary edge in triangulations, is introduced and used to prove the existence of flips in multitriangulations (either in full generality, like in~\cite{n-gdfcp-00}, or only for partial cases, like in~\cite{dkm-lahp-02}); then, the number of edges is deduced from the connectedness of the graph of flips. To summarize:

\begin{theorem}[\cite{cp-tttccp-92,n-gdfcp-00,dkm-lahp-02}]\label{intro:theo:fundamental}
\begin{enumerate}[(i)]
\item There is an inductive operation which transforms the \ktri{k}s of the \gon{(n+1)} into those of the \gon{n} and \viceversa~\cite{n-gdfcp-00}.
\item Any edge of length at least~$k+1$ in a \ktri{k} of the \gon{n} can be flipped
and the graph of flips is regular and connected~\cite{n-gdfcp-00,dkm-lahp-02}.
\item All \ktri{k}s of the \gon{n} have~$k(2n-2k-1)$ edges~\cite{cp-tttccp-92,n-gdfcp-00,dkm-lahp-02}.\qed
\end{enumerate}
\end{theorem}

In statement~(ii), the \defn{length} of an edge~$e$ of the \gon{n} is defined as the minimum between the numbers of remaining vertices of the \gon{n} located on each side of~$e$, augmented by~$1$ (boundary edges have length~$1$). Only edges of length greater than~$k$ are \defn{relevant} since the other ones cannot be part of a \kcross{(k+1)} (hence they appear in all \ktri{k}s).

\svs
Jakob~Jonsson~\cite{j-gt-03,j-gtdfssp-05} then completed these fundamental structural results in two directions. On the one hand, he studied enumerative properties of multitriangulations, generalizing in particular Proposition~\ref{intro:prop:catalan} as follows:

\begin{theorem}[\cite{j-gtdfssp-05}]\label{intro:theo:enumeration}
\index{Catalan number}
The number~$\theta(n,k)$ of \ktri{k}s of the \gon{n} equals:
$$\theta(n,k)=\det(C_{n-i-j})_{1\le i,j\le k}=\det\begin{pmatrix} C_{n-2} & C_{n-3} & \edots & \edots & C_{n-k-1} \\ C_{n-3} & \edots & \edots & C_{n-k-1} & \edots \\ \edots & \edots & \edots & \edots & \edots \\ \edots & C_{n-k-1} & \edots & \edots & C_{n-2k+1} \\ C_{n-k-1} & \edots & \edots & C_{n-2k+1} & C_{n-2k} \end{pmatrix}.$$
\end{theorem}
\vspace*{-.8cm}\qed

\mvs
The proof of this theorem relies on a more general enumeration result concerning \kdiag{k}-free $0/1$-fillings of polyominoes\index{polyomino}. Interestingly, the Hankel determinant in this statement is also known to count a certain family of non-intersecting \tuple{k}s of Dyck paths~\cite{gv-bdphlf-85}. The equality of the cardinalities of these two combinatorial families raises the question to find an explicit bijection (generalizing the bijection in \fref{intro:fig:catalan}). Even if Sergi~Elizalde~\cite{e-btdp-07} and \mbox{Carlos Nicolas}~\cite{n-abtdp-09} exhibit two different bijections for the case when~$k=2$, this question remains open for general~$k$. We discuss this problem in Section~\ref{ft:sec:Dyckpaths}.

\svs
On the other hand, Jakob~Jonsson studied in~\cite{j-gt-03} the simplicial complex~$\Delta_{n,k}$ formed by all \kcross{(k+1)}-free subsets of relevant edges of the \gon{n}. Since all its maximal elements have cardinality~$k(n-2k-1)$, this complex is pure of dimension~$k(n-2k-1)-1$. Jakob~Jonsson proved that this complex is even a~topological~sphere:

\begin{theorem}[\cite{j-gt-03}]
The simplicial complex~$\Delta_{n,k}$ is a vertex-decomposable piece-wise linear sphere of dimension~$k(n-2k-1)-1$.\qed
\end{theorem}

This theorem, together with small cases examples, leads to conjecture that this complex is even polytopal (\ie the boundary complex of a convex polytope) as happens when~$k=1$ (see Theorem~\ref{intro:theo:associahedron}). We discuss this open question in Section~\ref{ft:sec:multiassociahedron}.


\subsection{Overview of our results}\label{intro:subsec:multitriangulations:results}

Our contribution to the study of multitriangulations is based on their stars, which generalize triangles for triangulations:

\begin{definition}
\index{star@\kstar{k}}
A \defn{\kstar{k}} is the subset of edges of length~$k$ of a complete geometric graph on~$2k+1$ vertices in convex position. In other words, it is a (non-simple) polygon with~$2k+1$ vertices~$s_0,\dots,s_{2k}$ cyclically ordered, and~$2k+1$ edges~$[s_0,s_k],[s_1,s_{k+1}],\dots,[s_{2k},s_{k-1}]$.
\end{definition}

\begin{figure}[!h]
	\capstart
	\centerline{\includegraphics[scale=1]{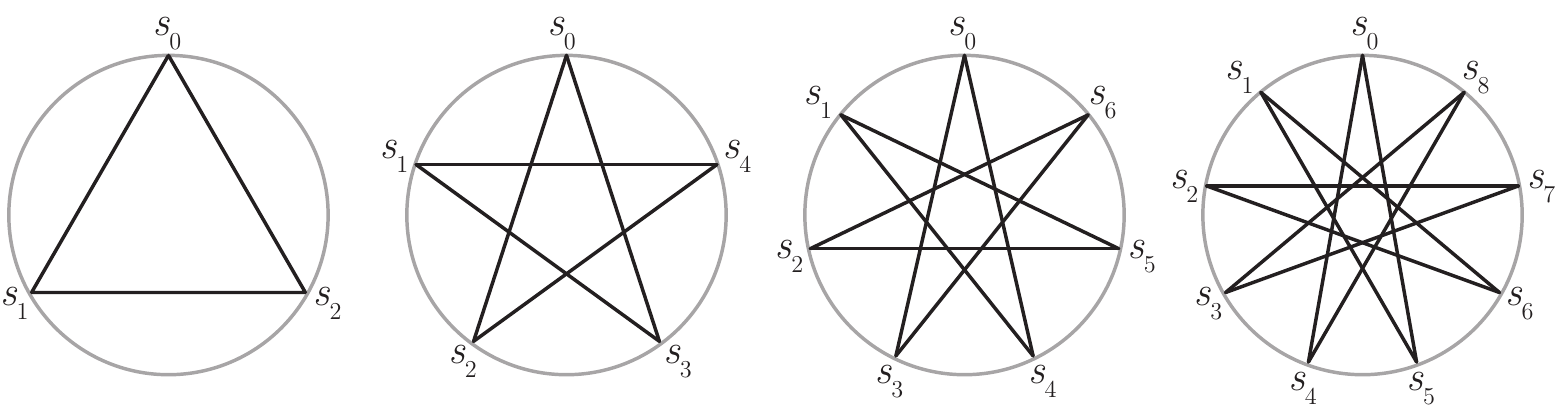}}
	\caption[Examples of stars]{Examples of \kstar{k}s, for $k\in[4]$.}
	\label{intro:fig:stars}
\end{figure}

Stars in multitriangulations play exactly the same role as triangles for triangulations: they decompose multitriangulations into simpler pieces and provide a useful tool to study their combinatorial properties. As evidence for this, we give in Chapter~\ref{chap:stars} direct and simple proofs of the fundamental properties of multitriangulations of Theorem~\ref{intro:theo:fundamental} (number of edges, definition and properties of the flip operation, \etc). Afterwards, we use stars in Chapter~\ref{chap:mpt} to relate multitriangulations with pseudotriangulations and pseudoline arrangements. This new interpretation provides a natural generalization of multitriangulations to any point set in general position in the plane, which keeps the rich combinatorial properties of multitriangulations (in particular the flip operation). Finally, in Chapter~\ref{chap:multiassociahedron}, we discuss in terms of stars the open questions presented above, and present some partial solutions as well as some natural failed attempts based on properties of stars. In the end of this section, we present in more details our contribution, summarizing and underlying our main results.

\subsubsection{Stars in multitriangulations}

In Chapter~\ref{chap:stars}, we study elementary properties of multitriangulations in terms of stars. In Section~\ref{stars:sec:angles}, we prove our main structural result concerning the incidence relations between the \kstar{k}s of a \ktri{k}:

\begin{theorem}\label{intro:theo:incidences}
Let~$T$ be a \ktri{k} of the \gon{n} (with~$n\ge 2k+1$).
\begin{enumerate}[(i)]
\item Each edge belongs to zero, one or two \kstar{k}s of~$T$, depending on whether its length is smaller, equal or greater than~$k$.
\item $T$~has exactly~$n-2k$ \kstar{k}s and~$k(2n-2k-1)$ edges.\qed
\end{enumerate}
\end{theorem}

Let us mention that we obtain Part~(ii) of Theorem~\ref{intro:theo:incidences} by a direct double counting argument similar to our first proof of Lemma~\ref{intro:lem:2n-3}: we use two simple independent relations between the number of edges and the number of \kstar{k}s in a \ktri{k}. The first relation is of course derived from Part~(i) of Theorem~\ref{intro:theo:incidences}, while the second one relies on the properties of the bisectors of \kstar{k}s. A \defn{bisector}\index{bisector} of a \kstar{k} is a bisector of one of its angles, that is, a line which passes through one of its vertices and separates its remaining vertices into two equal parts of cardinality~$k$. As for triangulations, the second relation is derived from the following identification between the pairs of \kstar{k}s of a \ktri{k}~$T$ and the edges of the \gon{n} which are not in~$T$:

\begin{theorem}\label{intro:theo:uniquebisector}
Let~$T$ be a \ktri{k} of the \gon{n}.
\begin{enumerate}[(i)]
\item Any pair of \kstar{k}s of~$T$ has a unique common bisector, which is not in~$T$.
\item Reciprocally, any edge not in~$T$ is the common bisector of a unique pair of \kstar{k}s~of~$T$.\qed
\end{enumerate}
\end{theorem}

We also use stars and their common bisectors in Section~\ref{stars:sec:flips} to understand the flip operation as a local transformation. In much the same way than a flip in a triangulation consists in replacing one diagonal by the other in a quadrangle formed by two adjacent triangles, flips in multitriangulations can be interpreted as a transformation involving only two adjacent stars:

\begin{theorem}\label{intro:theo:flip}
\index{flip}
Let~$T$ be a \ktri{k} of the \gon{n},~$e$ be an edge of length greater than~$k$, and~$f$ be the common bisector of the two \kstar{k}s of~$T$ containing~$e$. Then~$T\diffsym\{e,f\}$ is a \ktri{k} of the \gon{n} and it is the only one, except~$T$ itself, which contains~$T\ssm\{e\}$.\qed
\end{theorem}

\begin{figure}[h]
	\capstart
	\centerline{\includegraphics[scale=1]{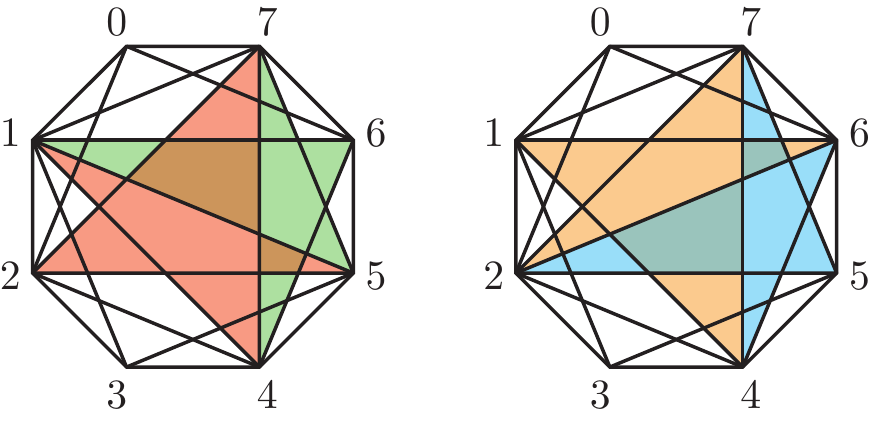}}
	\caption[A flip in the \ktri{2} of the octagon of \fref{intro:fig:2triang8points}]{A flip in the \ktri{2} of the octagon of \fref{intro:fig:2triang8points}.}
	\label{intro:fig:2triang8pointsflip}
\end{figure}

This interpretation simplifies the study of the graph of flips (whose vertices are the \ktri{k}s of the \gon{n} and whose edges are the flips between them) and yields to new proofs of the results in~\cite{n-gdfcp-00} as well as partial extensions of them:

\vspace{2.2cm}\qed
\vspace{-2.7cm}
\begin{theorem}
\index{flip!graph of ---s}
\index{diameter (of the graph of flips)}
The graph of flips on \ktri{k}s of the \gon{n} is \regular{k(n-2k-1)} and connected, and for all~$n>4k^2(2k+1)$, its diameter~$\delta_{n,k}$ is bounded by
$$2\Fracfloor{n}{2}\left(k+\frac{1}{2}\right)-k(2k+3) \le \delta_{n,k} \le 2k(n-4k-1).$$
\end{theorem}

The upper bound on the diameter in this theorem was already proved in~\cite{n-gdfcp-00} (and has to be compared with the exact value~$2n-10$ of the diameter of the graph of flips on triangulations of the \gon{n} mentioned in Theorem~\ref{intro:theo:diameter}). The lower bound is new, but only slightly improves the trivial bound $\delta_{n,k}\ge k(n-2k-1)$. 

\svs
We also use stars in Section~\ref{stars:sec:ears} to study the \defn{\kear{k}s}\index{ear@\kear{k}} (\ie the edges of length~$k+1$) of a \ktri{k}. We present a simple proof that any \ktri{k} contains at least~$2k$ \kear{k}s~\cite{n-gdfcp-00}, and we provide several characterizations of the \ktri{k}s which attain this bound:

\begin{theorem}
\index{star@\kstar{k}!internal ---}
\index{triangulation@\ktri{k}!colorable@\kcolorable{k} ---}
\index{accordion@\kaccordion{k}}
The number of \kear{k}s in a \ktri{k} equals~$2k$ plus its number of internal \kstar{k}s, \ie those that do not contain an edge of length~$k$.

Furthermore, when~$k>1$ and~$T$ is a \ktri{k}, the following conditions are equivalent:
\begin{enumerate}[(i)]
\item $T$~contains exactly~$2k$ \kear{k}s;
\item $T$~has no internal \kstar{k};
\item $T$~is \defn{\kcolorable{k}}, \ie it is possible to color its relevant edges with~$k$ colors such that no crossing is monochromatic;
\item the set of relevant edges of~$T$ is the union of $k$ disjoint accordions (an \defn{accordion} is a sequence of edges~$[a_i,b_i]$ such that for all~$i$, either $a_{i+1}=a_i$ and $b_{i+1}=b_i+1$ or ${a_{i+1}=a_i-1}$ and $b_{i+1}=b_i$).\qed
\end{enumerate}
\end{theorem}

\index{flattening (an external \kstar{k})}
\index{inflating (an external \kcross{k})}
In Section~\ref{stars:sec:flat-infl}, we reinterpret in terms of stars the inductive operation of Theorem~\ref{intro:theo:fundamental}(i) which transforms the \ktri{k}s of the \gon{(n+1)} into those of the \gon{n} and \viceversa. It consists in \defn{flattening} (or \defn{inflating}) a single \kstar{k} adjacent to an edge of length~$k$. We voluntarily only introduce this operation at the end of Chapter~\ref{chap:stars} to emphasize that none of our proofs so far make use of induction.

\svs
Finally, Section~\ref{stars:sec:surfaces} completes our study of the structural properties of multitriangulations in terms of stars. According to Theorem~\ref{intro:theo:incidences}, we see \ktri{k}s of the \gon{n} as polygonal decompositions of a certain surface into~$n-2k$ \gon{k}s. We exploit this interpretation to construct, via multitriangulations, very regular decompositions of an infinite family of surfaces.

\subsubsection{Multitriangulations, pseudotriangulations, and duality}

\index{dual!--- of a point}
\index{dual!--- of a point set}
Chapter~\ref{chap:mpt} is devoted to the interpretation of multitriangulations in the line space of the plane, \ie in the M\"obius strip. We consider the usual duality between configurations of points and pseudoline arrangements (see~\cite[Chapters~5,6]{f-gga-04} and~\cite{g-pa-97}): the set~$p^*$ of all lines passing through a point~$p$ of the plane is a \defn{pseudoline}\index{pseudoline} (a non-separating simple closed curve) of the M\"obius strip, and the set~${P^* \eqdef \ens{p^*}{p\in P}}$ of all pseudolines dual to the points of a finite set~$P$ is a \defn{pseudoline arrangement}\index{pseudoline!--- arrangement} (any two pseudolines have exactly one crossing point). The starting point of Chapter~\ref{chap:mpt} is the following observation:

\begin{observation}
\index{dual!--- of a \kstar{k}}
\index{dual!--- of a \ktri{k}}
Let~$T$ be \ktri{k} of the \gon{n}. Then:
\begin{enumerate}[(i)]
\item the set~$S^*$ of all bisectors to a \kstar{k}~$S$ of~$T$ is a pseudoline;
\item the common bisector of two \kstar{k}s~$R$ and~$S$ is a crossing point of the two pseudolines~$R^*$ and~$S^*$, while a common edge of~$R$ and~$S$ is a contact point between~$R^*$ and~$S^*$;
\item the set~$T^* \eqdef \ens{S^*}{S \text{ \kstar{k} of } T}$ of all pseudolines dual to the \kstar{k}s of~$T$ is a \defn{pseudoline arrangement with contact points} (two pseudolines cross exactly once, but can touch each other a finite number of times);
\item the support of this arrangement covers precisely the support of the dual arrangement~$V_n^*$ of the vertex set of the \gon{n}, except its first~$k$ levels.
\end{enumerate}
\end{observation}

We prove the reciprocal statement of this observation:

\begin{theorem}\label{intro:theo:dualmulti}
Any pseudoline arrangement with contact points whose support covers precisely the support of the dual arrangement~$V_n^*$ of the vertex set of the \gon{n} except its first~$k$ levels is the dual arrangement of a \ktri{k} of the \gon{n}.\qed
\end{theorem}

A similar observation was made in~\cite{pv-ot-94,pv-ptta-96} for pointed pseudotriangulations of a point set in general position. A \defn{pseudotriangle}\index{pseudotriangle} is a polygon with only three convex angles joined by three concave polygonal chains, and a \defn{pseudotriangulation}\index{pseudotriangulation} of a point set~$P$ in general position is a set of edges which decomposes the convex hull of~$P$ into pseudotriangles. It is \defn{pointed}\index{pseudotriangulation!pointed ---} if for any point~$p\in P$, there is a line passing through~$p$ which defines a closed half-plane containing no edge incident to~$p$. For a pointed pseudotriangulation~$T$ of~$P$, it was observed that:
\index{dual!--- of a pseudotriangle}
\index{dual!--- of a pseudotriangulation}
\begin{enumerate}[(i)]
\item the set~$\Delta^*$ of internal tangents to a pseudotriangle~$\Delta$ is a pseudoline of the M\"obius strip;
\item the set~$T^* \eqdef \ens{\Delta^*}{\Delta \text{ pseudotriangle de } T}$ is a pseudoline arrangement with contact points supported by the dual arrangement~$P^*$ minus its first level.
\end{enumerate}

\begin{figure}[h]
	\capstart
	\centerline{\includegraphics[scale=1]{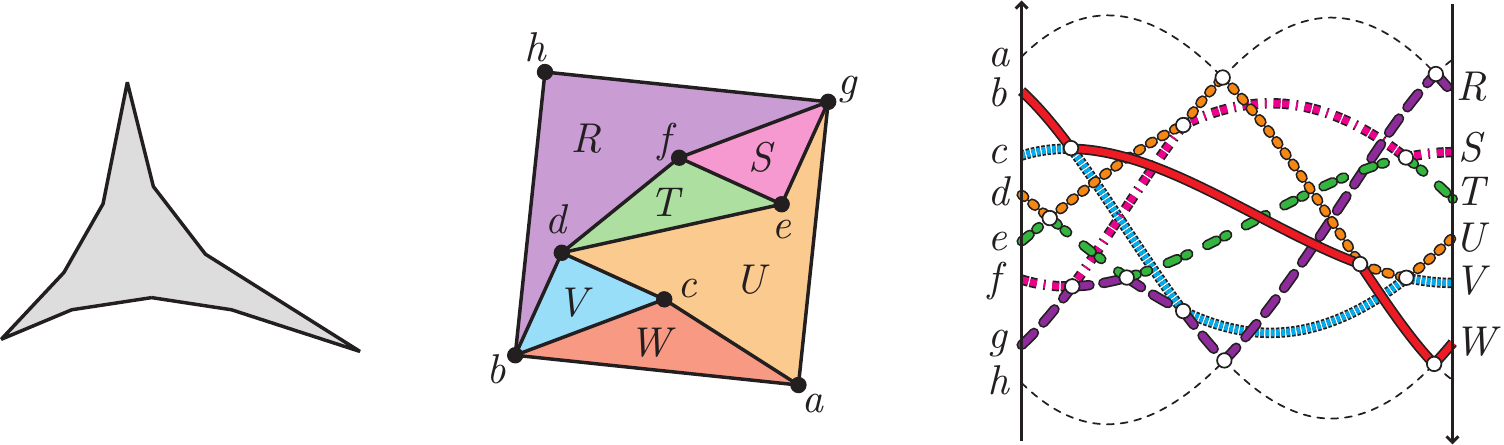}}
	\caption[A pseudotriangle, a pseudotriangulation, and its dual pseudoline arrangement]{A pseudotriangle, a pseudotriangulation, and its dual pseudoline arrangement.}
	\label{intro:fig:pseudotriangulationintro}
\end{figure}

We provide again different proofs of the reciprocal statement: 

\begin{theorem}\label{intro:theo:dualpseudo}
Let~$P$ be a point set in general position and~$P^*$ denote its dual pseudoline arrangement. Then any pseudoline arrangement with contact points whose support covers precisely the support of~$P^*$ except its first level is the dual arrangement of a pseudotriangulation of~$P$.\qed
\end{theorem}

\svs
\index{pseudoline!--- arrangement!--- --- with contact points}
Motivated by these two theorems, we then consider \defn{pseudoline arrangements with contact points} which share a common support. We define a flip operation corresponding to flips in multitriangulations and pseudotriangulations:

\begin{definition}
Two pseudoline arrangements with the same support are related by a flip if the symmetric difference of their sets of contact points is only a pair~$\{u,v\}$. In this case, in one of these arrangements $u$ is a contact point of two pseudolines which cross at~$v$, while in the other arrangement, $v$ is a contact point of two pseudolines which cross at~$u$.
\end{definition}

\index{greedy!--- pseudoline arrangement}
\index{pseudoline!--- arrangement!greedy --- ---}
We study in Section~\ref{mpt:sec:enumeration} the graph~$G(\cS)$ of flips on pseudoline arrangements with a common support~$\cS$. For example, when~$\cS$ is the support of an arrangement of two pseudolines with~$p$ contact points,~$G(\cS)$ is the complete graph on~$p+1$ vertices. We  consider certain acyclic orientations of the graph~$G(\cS)$, given by vertical cuts of the support~$\cS$. Each of these orientations has a unique source that we call the \defn{greedy pseudoline arrangement}, and that we characterize in terms of sorting networks. The properties of these greedy pseudoline arrangements provide an enumeration algorithm for pseudoline arrangements with the same support, whose working space is polynomial. We thus sheds light upon a similar algorithm existing for the enumeration of pseudotriangulations of point sets~\cite{bkps-ceppgfa-06}, providing a complementary proof of it.

\svs
We then come back in Sections~\ref{mpt:sec:duality} and~\ref{mpt:sec:mpt} to the specific cases of multitriangulations and pseudotriangulations. According to Theorems~\ref{intro:theo:dualmulti} and~\ref{intro:theo:dualpseudo}, we propose the following generalization of both multitriangulations and pseudotriangulations:

\begin{definition}
\index{pseudotriangulation@\pt{k}!--- of a pseudoline arrangement}
\index{multipseudotriangulation!--- of a pseudoline arrangement}
\index{pseudotriangulation@\pt{k}!--- of a point set}
\index{multipseudotriangulation!--- of a point set}
A \defn{\pt{k}} of a pseudoline arrangement~$L$ is a pseudoline arrangement whose support is~$L$ minus its first~$k$ levels. A \defn{\pt{k}} of a point set~$P$ in general position in the plane is a set~$T$ of edges which correspond by duality to the contact points of a \pt{k}~$T^*$ of~$P^*$.
\end{definition}

In contrast to other possible candidates for generalization of multitriangulations to non-convex position, this definition preserves the rich combinatorial structure of multitriangulations. For example, a \pt{k} of a point set~$P$ has automatically~$k(2|P|-2k-1)$ edges, and any (sufficiently internal) edge can be flipped to create a new \pt{k} of~$P$.

Following our directed line, we again focus on stars in \mpt{}s, which are now defined as follows: a \defn{star}\index{star} in a \pt{k}~$T$ is a set of edges corresponding to all the contact points of a single pseudoline of the dual pseudoline arrangement~$T^*$. As immediate consequences of the definition, the stars of a \mpt have analog properties to that of Theorems~\ref{intro:theo:incidences},~\ref{intro:theo:uniquebisector} and~\ref{intro:theo:flip}.

\svs
\enlargethispage{.5cm}
We conclude Chapter~\ref{chap:mpt} with three questions related to \mpt{}s:
\begin{enumerate}
\item We study \defn{iterated \mpt{}s}\index{multipseudotriangulation!iterated ---} in Section~\ref{mpt:sec:iterated}: a \pt{k} of an \pt{m} of a pseudoline arrangement~$L$ is a \pt{(k+m)} of~$L$. We give an example of \ktri{2} not containing any triangulation. We prove though that greedy \mpt{}s are iterations of greedy~\mbox{pseudotriangulations}.
\item In Section~\ref{mpt:subsec:furthertopics:horizon}, we characterize the edges of the greedy \pt{k}s of a point set~$P$ in terms of \defn{$k$-horizon trees} of~$P$. This characterization generalizes an observation of Michel Pocchiola~\cite{p-htvpt-97} in the context of pseudotriangulations.
\item Finally in Section~\ref{mpt:subsec:furthertopics:dpl}, we study \mpt{}s of configurations of disjoint convex bodies of the plane. The dual arrangements of configurations of disjoint convex bodies are \defn{double pseudoline arrangements}\index{double pseudoline!--- arrangement}, introduced in~\cite{hp-adp-08} by Luc Habert and Michel Pocchiola. Let us mention that we manipulate these arrangements in Appendix~\ref{app:sec:dpl} to enumerate all arrangements of at most~$5$ double pseudolines up to isomorphism.
\end{enumerate}

\subsubsection{Further topics}

Finally, we discuss in Chapter~\ref{chap:multiassociahedron} three remaining open problems on multritriangulations which illustrate the rich combinatorial and geometric structure of multitriangulations. Our goal is to present natural ideas based on stars which could be fertile, although they provide at the present time only partial results to these problems.

Our first problem (Section~\ref{ft:sec:Dyckpaths}) is that of finding an \defn{explicit bijection} between the set of multitriangulations and the set of non-crossing \tuple{k}s of Dyck paths, which are both counted by the Hankel determinant of Theorem~\ref{intro:theo:enumeration}. As mentioned in our introduction on triangulations, the indegree sequence of a triangulation defines a natural bijection to Dyck paths. Generalizing this simple remark, Jakob~Jonsson compared in~\cite{j-gt-03} the repartition of indegree sequences of \ktri{k}s with that of signatures of non-crossing \tuple{k}s of Dyck paths (see Definition~\ref{ft:def:signature} for the precise statement of the definition of these parameters). Following this result, we try to find a suitable way to decompose the relevant edges of a \ktri{k} into~$k$ disjoint sets such that the respective indegree sequences of these sets define~$k$ non-crossing Dyck paths. We point out a simple such decomposition based on stars, for which the resulting map from \ktri{k}s to non-crossing \tuple{k}s of Dyck paths is unfortunately not bijective.

\index{rigidity}
Our second open problem (Section~\ref{ft:sec:rigidity}) concerns \defn{rigidity} properties of multitriangulations. A triangulation is a \defn{minimally rigid} object in the plane: the only continuous motions of its vertices which preserve its edges' lengths are the isometries of the plane, and removing any of its edges makes it flexible. Equivalently, a triangulation of the \gon{n} satisfies Laman's property: it has~$2n-3$ edges and any subgraph on~$m$ vertices has at most~$2m-3$ edges. We observe two interesting connections of multitriangulations to rigidity (and related topics):
\begin{enumerate}[(i)]
\item First, a \ktri{k} of the \gon{n} is a $(2k,{2k+1 \choose 2})$-\defn{tight graph}\index{tight graph}: it has~$2kn-{2k+1 \choose 2}$ edges and any subgraph on~$m$ vertices has at most~$2km-{2k+1 \choose 2}$ edges. This yields the conjecture that any \ktri{k} is a (generically) minimally rigid graph in dimension~$2k$. We prove this conjecture when~$k=2$.
\item We then show that the dual graph of a \ktri{k} is \tight{(k,k)}: it can be decomposed into~$k$ edge-disjoint spanning trees. This has to be understood has a generalization of the dual tree of a triangulation.
\end{enumerate}

We finally focus (in Section~\ref{ft:sec:multiassociahedron}) on the \defn{polytopality question} of the simplicial complex~$\Delta_{n,k}$ formed by all \kcross{(k+1)}-free subsets of relevant edges of the \gon{n}. We propose two different contributions to this problem:
\begin{enumerate}[(i)]
\item On the one hand (Section~\ref{ft:subsec:multiassociahedron:delta82}), we answer the minimal non-trivial case: we describe the space of symmetric realizations of~$\Delta_{8,2}$. This result is obtained in two steps: first, we enumerate by computer all symmetric oriented matroids realizing our simplicial complex; then, we study the possible geometric realizations of the resulting oriented matroids. 
\item On the other hand (Section~\ref{ft:subsec:multiassociahedron:loday}), we discuss a natural generalization of a construction of the associahedron due to Jean-Louis~Loday~\cite{l-rsp-04}. We obtain a polytope with a particularly simple facet description which realizes the flip graph restricted to certain multitriangulations (those whose dual graph is acyclic). This polytope could have been \apriori a projection of a polytope realizing~$\Delta_{n,k}$, but we show that such a projection is impossible.
\end{enumerate}


\section{Sources of material}

The results presented in the first part of this dissertation have been partially published (or submitted) in the following papers:
\begin{enumerate}
\item Chapter~\ref{chap:stars}~---~where we introduce stars and use them to recover and simplify some structural properties of multitriangulations~---~is a joint work with Francisco Santos and appeared in~\cite{ps-mtcsp-09}. Only few minor results, including those concerning the diameter of the graph of flips (Lemmas~\ref{stars:conj:increasingdiametermin} and~\ref{stars:lem:lowerbounddiameter}) and the construction of equivelar surfaces (Section~\ref{stars:subsec:surfaces:equivelar}), have been added since the publication of this article. This paper also contains a detailed discussion of various open questions related to multitriangulations, some of which are explored in Chapter~\ref{chap:multiassociahedron}.
\item All the content of Chapter~\ref{chap:mpt}~---~the systematic study and enumeration of the graph of flips on pseudoline arrangements with the same support, the relationship between pseudotriangulations and multitriangulations, and the definition and properties of \mpt{}s~---~is the result of a joint work with Michel Pocchiola~\cite{pp-mpt-09}.
\item Chapter~\ref{chap:multiassociahedron} is the concatenation of various discussions I had in particular with Francisco Santos, and a large part of this chapter is inspired from~\cite[Section~8]{ps-mtcsp-09}. As far as the two main contributions of this chapter are concerned:
\begin{enumerate}
\item The description of the space of symmetric realizations of~$\Delta_{8,2}$ is a joint work with J\"urgen Bokowski~\cite{bp-srmt-09}. Some \haskell code developed to enumerate symmetric matroid realizations is presented and explained in Appendix~\ref{app:sec:enumerationmatroidpolytopes}.
\item The generalization of Loday's construction of the associahedron is a joint work with Francisco Santos~\cite{ps-gla-10}.
\end{enumerate}
\item Finally, Appendix~\ref{app:sec:dpl}~---~where we enumerate isomorphism classes of arrangements with few pseudolines and double pseudolines ~---~is a joint work with Julien Fert\'e and Michel Pocchiola~\cite{fpp-nsafdp-08}.
\end{enumerate}

	\chapter{Stars in multitriangulations}\label{chap:stars}

In this chapter, we study structural properties of multitriangulations. Our approach is to decompose a multitriangulation into a complex of stars, which generalize triangles for triangulations. In much the same way as triangles are central objects in triangulations, stars are our main tool to understand multitriangulations. In this chapter, we use them to:
\begin{enumerate}[(i)]
\item Give a direct proof that the number of edges is the same in all multitriangulations.
\item Define locally the ``flip operation'' which transforms a multitriangulation into another one by exchanging a single edge. We study extensively the graph of flips, with a particular attention to its diameter.
\item Characterize a particularly nice and simple family of multitriangulations obtained as superpositions of triangulations.
\item Understand the ``deletion operation'' which transforms multitriangulations of an \gon{(n+1)} into multitriangulations of an \gon{n}, and \viceversa. This transformation allows recursive proofs for some properties of multitriangulations.
\item Interpret multitriangulations as decompositions of a certain surface. As an application of multitriangulaitons, we obtain examples of very regular decompositions of surfaces.
\end{enumerate}


\section{Notations and examples}\label{stars:sec:notations}


\subsection{Multicrossings, multitriangulations, and relevant edges}\label{stars:subsec:notations:edges}

Let~$k$ and~$n$ be two integers such that~$k\ge 1$ and~$n\ge 2k+1$.

Let~$V_n$ be the \defn{set of vertices} of a convex \gon{n}, \ie any set of points on the unit circle, labeled counterclockwise by the cyclic set~$\Z_n$. We always refer to the points in~$V_n$ by their labels to simplify notation. For~$u,v,w\in V_n$, we write~$u\cl v\cl w$ meaning that~$u$,~$v$~and~$w$ are in counterclockwise order on the circle. For any~$u,v\in V_n$, let~$\llb u,v\rrb$ denote the \defn{cyclic interval}~$\ens{w\in V_n}{u\cle w\cle v}$. The intervals~$\rrb u,v\llb$, $\llb u,v\llb$~and~$\rrb u,v\rrb$ are defined similarly. Let~$|u-v| \eqdef \min(|\llb u, v\llb|,|\llb v, u\llb|)$ be the \defn{cyclic distance} between~$u$ and~$v$.

For~$u\ne v\in V_n$, let~$[u,v]$ denote the straight \defn{edge} connecting the vertices~$u$ and~$v$. We say that~$[u,v]$ is of \defn{length}~$|u-v|$. By definition, an edge of length~$\ell$ passes through two points of~$V_n$ and separates~$\ell-1$ points from the other~$n-\ell-1$. Let~$E_n \eqdef {V_n \choose 2}$ be the \defn{set of edges} of the complete graph on~$V_n$. Two edges~$[u,v]$ and~$[u',v']$ are said to \defn{cross}\index{crossing|hbf} when~$u$ and~$v$ lie one in each of the open cyclic intervals~$\rrb u',v'\llb$ and~$\rrb v',u'\llb$. In other words, the edges~$[u,v]$ and~$[u',v']$ cross when~$u\cl u'\cl v\cl v'\cl u$ or~$u'\cl u\cl v'\cl v\cl u'$. Observe that this definition is symmetric in the two edges involved and that it corresponds to an intersection of the straight open segments~$(u,v)$ and~$(u',v')$. Similarly to our discussion on triangulations in Section~\ref{intro:sec:triangulations}, where we considered subsets of $E_n$ avoiding crossings, we are now interested in subsets of $E_n$ avoiding the following configuration:

\begin{definition}\label{stars:def:kcrossing}
\index{crossing@\kcross{k}|hbf}
For~$\ell\in\N$, an \defn{\kcross{\ell}} is a set of~$\ell$~mutually intersecting edges.
\end{definition}

We will always order the edges of an \kcross{\ell}~$F \eqdef \ens{f_i}{i\in[\ell]}$ cyclically such that the endpoints~$x_i$ and~$y_i$ of~$f_i$ satisfy $x_1\cl x_2\cl\cdots\cl x_\ell\cl y_1\cl y_2\cl\cdots\cl y_\ell\cl x_1$ (see \fref{stars:fig:edges}).

As for triangulations, we have a particular attention to maximal sets of edges avoiding these forbidden configurations:

\begin{definition}\label{stars:def:ktriangulation}
\index{triangulation@\ktri{k}|hbf}
\index{multitriangulation|hbf}
A \defn{\ktri{k}} of the \gon{n} is a maximal \kcross{(k+1)}-free subset of~$E_n$.
\end{definition}

Obviously, an edge~$[u,v]$ of~$E_n$ can appear in a \kcross{(k+1)} only if it has at least~$k$~vertices on each side, \ie if~$|u-v|>k$. 
We distinguish three types of edges, according on whether their length is lower, equal or greater than~$k$~(see \fref{stars:fig:edges}):

\begin{figure}
	\capstart
	\centerline{\includegraphics[scale=1]{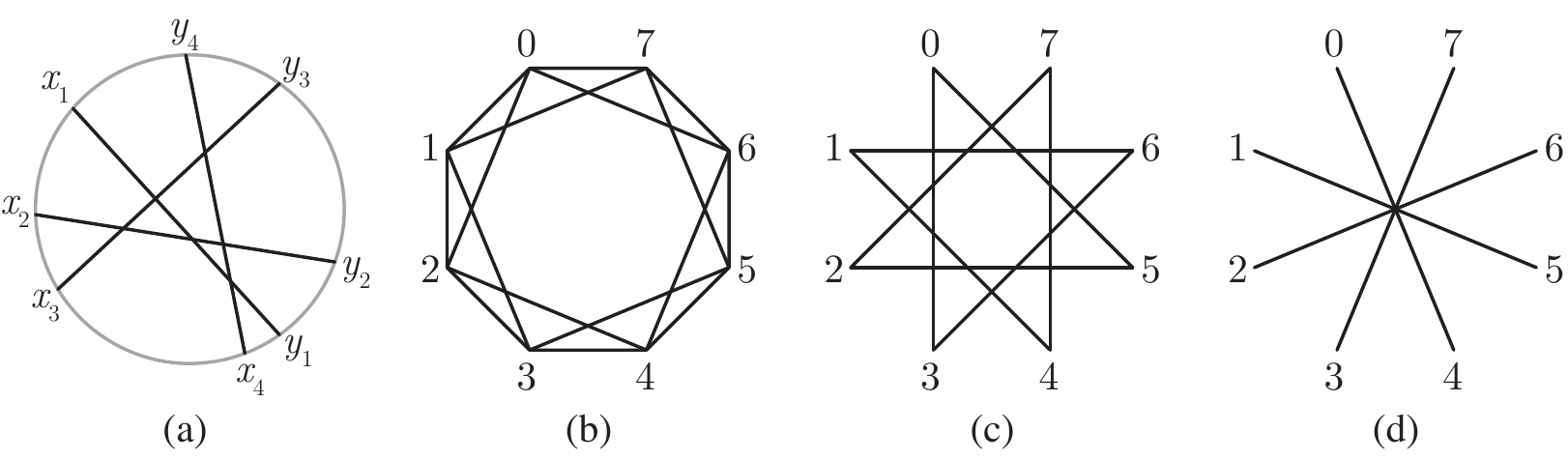}}
	\caption[A \kcross{4} and the \kirrel{3}, \kbound{3}, and \krel{3} edges of an octagon]{A \kcross{4}~(a) and the sets of \kirrel{3}~(b), \kbound{3}~(c), and \krel{3}~(d) edges of an octagon.}
	\label{stars:fig:edges}
\end{figure}

\begin{definition}\label{stars:def:krelevant}
\index{relevant@\krel{k}!--- edge|hbf}
\index{boundary@\kbound{k} edge|hbf}
\index{irrelevant@\kirrel{k} edge|hbf}
An edge~$[u,v]$ of~$E_n$ is:
\begin{enumerate}[(i)]
	\item a \defn{\krel{k}} edge if~$|u-v|>k$;
	\item a \defn{\kbound{k}} edge if~$|u-v|=k$;~and
	\item a \defn{\kirrel{k}} edge if~$|u-v|<k$.
\end{enumerate}
\end{definition}

By maximality, every \ktri{k} of the \gon{n} consists of all the~$kn$~\kirrel{k} plus \kbound{k} edges of~$E_n$ and some \krel{k} edges.


\subsection{Examples}\label{stars:subsec:notations:examples}

Before going further, we consider small values of~$k$ and~$n$ (namely~$k=1$ or~$n$ close to~$2k+1$), for which we can easily describe all \ktri{k}s of the \gon{n}.

\begin{example}[$k=1$]\label{stars:exm:k=1}
Since a \kcross{2} is just a crossing, the \ktri{1}s are just triangulations of the \gon{n}. All internal diagonals of the \gon{n} are \krel{1}, while the boundary edges of the \gon{n} are \kbound{1} (and there is no \kirrel{1} edges). See Section~\ref{intro:sec:triangulations} for a discussion on structural properties of triangulations of the \gon{n}.
\end{example}

\begin{example}[$n=2k+1$]\label{stars:exm:n=2k+1}
The complete graph~$K_{2k+1}$ on~$2k+1$ vertices does not contain~$k+1$ mutually intersecting edges, and thus it is the unique \ktri{k} of the \gon{(2k+1)}. The longest diagonals $[i,i+k]$ are \kbound{k}, while all the other edges are \kirrel{k}.
\end{example}

\begin{example}[$n=2k+2$]\label{stars:exm:n=2k+2}
The set of edges~$E_{2k+2}$ contains exactly one \kcross{(k+1)} formed by the~$k+1$~long diagonals~$[i,i+k+1]$ of the~\gon{(2k+2)}. Consequently, there are~$k+1$ \ktri{k}s of the \gon{(2k+2)}, each of them obtained from the complete graph~$K_{2k+2}$ by deleting one of these long diagonals (see Figures~\ref{stars:fig:edges}(b-c-d) and~\ref{stars:fig:2k+2}).

\begin{figure}[!h]
	\capstart
	\centerline{\includegraphics[scale=1]{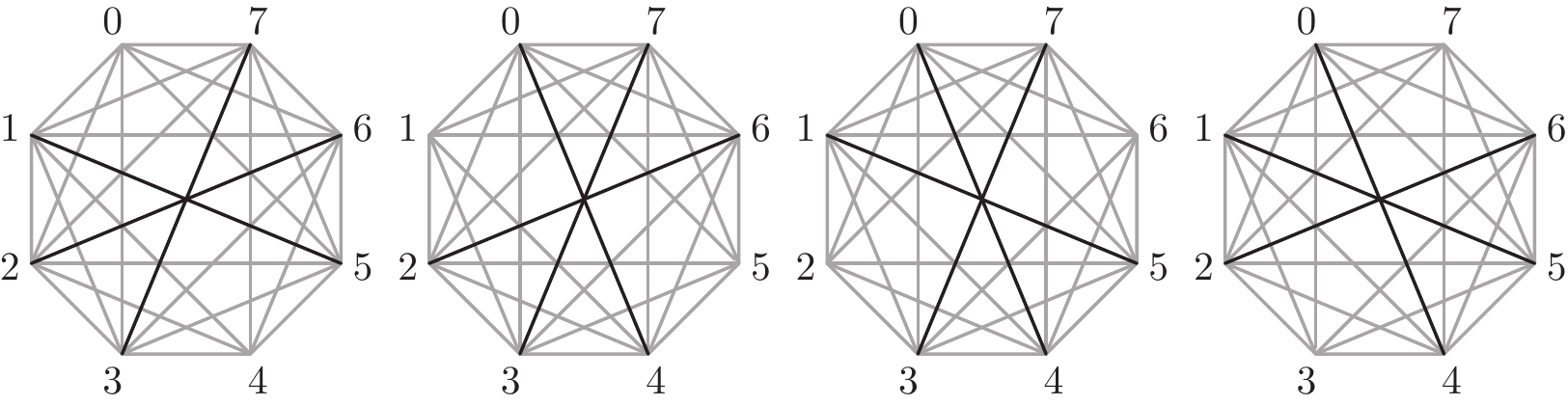}}
	\caption[The four \ktri{3}s of the octagon]{The four \ktri{3}s of the octagon. \krel{3} edges are black, while \kbound{3} and \kirrel{3} edges are gray.}
	\label{stars:fig:2k+2}
\end{figure}

\end{example}

\begin{example}[$n=2k+3$]\label{stars:exm:n=2k+3}
The \krel{k} edges of~$E_{2k+3}$ are exactly the longest diagonals ${e_i \eqdef [i,i+k+1]}$. They form a polygonal cycle\footnote{Remember that we have indexed the vertices of~$V_{2k+3}$ with the cyclic group~$\Z_{2k+3}$, \ie that the indices have to be understood modulo $2k+3$.}~$e_0,e_{k+1},e_{2(k+1)},\dots,e_{-(k+1)},e_0$ in which any two diagonals cross except if they are consecutive (see for example \fref{stars:fig:2k+3}(a)). Consequently, if~$E$ is a \kcross{(k+1)}-free subset of~$E_{2k+3}$, then:
\begin{enumerate}[(i)]
\item it cannot contain more than~$k+1$ non-consecutive long diagonals;~and
\item if $E$ contains $2\ell+1$~consecutive long diagonals~$e_i,e_{i+(k+1)},\dots,e_{i+2\ell(k+1)}$, but neither the diagonal $e \eqdef e_{i+(2\ell+1)(k+1)}$ nor the diagonal $f \eqdef e_{i-(k+1)}$, then both sets $E\cup\{e\}$ and ${E\cup\{f\}}$ are still \kcross{(k+1)}-free.
\end{enumerate}
We deduce from these two observations that the \ktri{k}s of the \gon{(2k+3)} are exactly the disjoint unions of~$k$~pairs of consecutive long diagonals of~$E_{2k+3}$ (plus of course the non-\krel{k} edges). For example, there are fourteen~\ktri{2}s of the heptagon, all obtained by rotation of one of the two \ktri{2}s of \fref{stars:fig:2k+3}(b-c).

\begin{figure}[!h]
	\capstart
	\centerline{\includegraphics[scale=1]{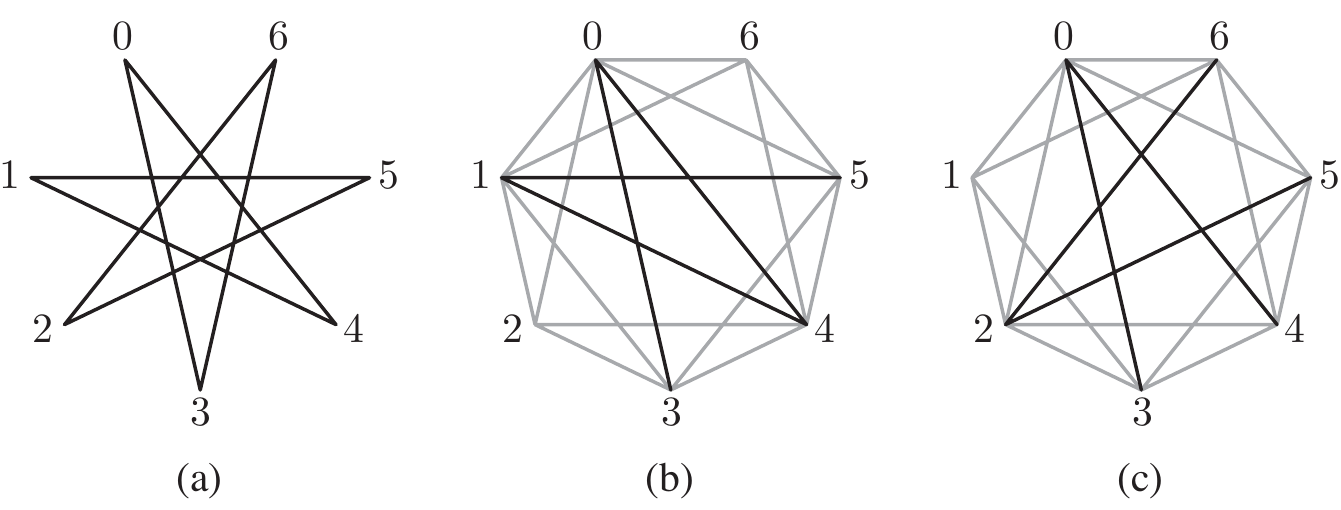}}
	\caption[The \krel{2} edges and two \ktri{2}s of the heptagon]{(a)~The \krel{2} edges of the heptagon and (b-c)~two \ktri{2}s of the heptagon (representing all \ktri{2}s of the heptagon up to rotation).}
	\label{stars:fig:2k+3}
\end{figure}

\end{example}




\subsection{Stars}\label{stars:subsec:notations:stars}

\begin{figure}[!h]
	\capstart
	\centerline{\includegraphics[scale=1]{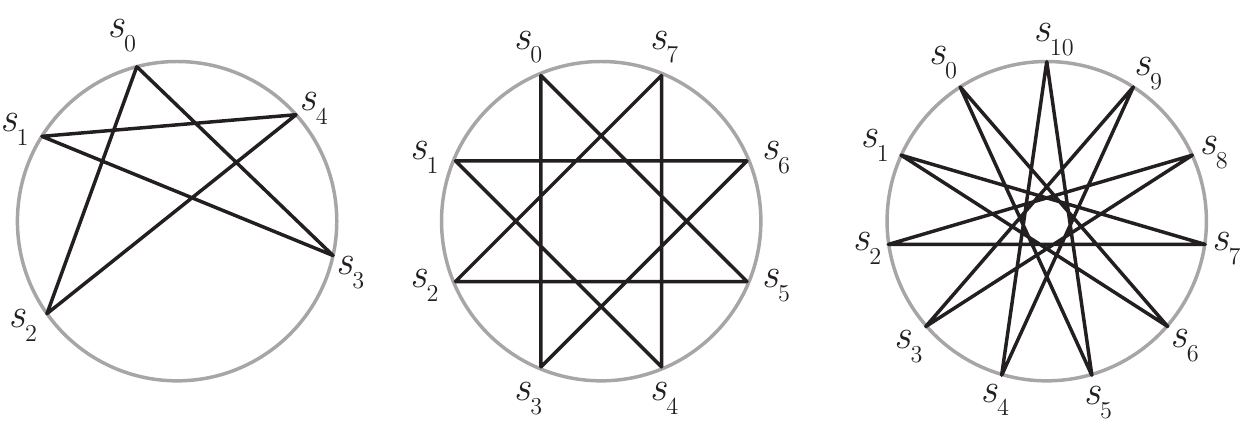}}
	\caption[Star polygons]{Star polygons of type~$\{5/2\}$,~$\{8/3\}$~and~$\{11/5\}$. The left one is a \kstar{2} and the right one is a \kstar{5}.}
	\label{stars:fig:starpolygons}
\end{figure}

\index{star!--- polygon}
As mentioned in the introduction, our study of \ktri{k}s is based on a natural generalization of the triangles for triangulations. If~$p$ and~$q$ are two coprime integers, a \defn{star polygon} of type~$\{p/q\}$ is a polygon formed by connecting a set~$V \eqdef \ens{s_j}{j\in\Z_p}$ of~$p$ points of the unit circle with the set~${E \eqdef \ens{[s_{j},s_{j+q}]}{j\in\Z_p}}$ of edges (see~\cite[\mbox{Chapter~2.3,~pp.~36-38}]{c-ig-69},~\cite[\mbox{Chapter~6,~pp.~93-95}]{c-rp-73} and~\fref{stars:fig:starpolygons}). The role of triangles for \ktri{k}s is played by some special star polygons:

\begin{definition}\label{stars:def:kstar}
\index{star@\kstar{k}|hbf}
A \defn{\kstar{k}} is a star polygon of type~$\left\{\frac{2k+1}{k}\right\}$: it is a (non-simple) polygon with~$2k+1$ vertices~$s_0\cl \cdots\cl s_{2k}$ cyclically ordered and with $2k+1$~edges~$[s_0,s_k],[s_1,s_{k+1}],\dots,[s_{2k},s_{k-1}]$.
\end{definition}

\index{star!--- order}
\index{cyclic order}
\index{order!star ---}
\index{order!cyclic ---}
Observe that there are two natural cyclic orders on the vertices of a \kstar{k}~$S$: the \defn{circle order}, defined as the cyclic order around the circle, and the \defn{star order}, defined as the cyclic order tracing the edges of~$S$. More precisely, if~$s_0,\dots,s_{2k}$ are the vertices of~$S$ cyclically ordered around the circle, we rename the vertices~$r_i \eqdef s_{ik}$ to obtain the star order~$r_0,\dots,r_{2k}$.


\section{Angles, bisectors, and star decomposition}\label{stars:sec:angles}

In this section, we present our main tool to study structural properties of multitriangulations: we show that \ktri{k}s can be decomposed into \kstar{k}s in the same way that triangulations can be thought of as decompositions into triangles. This structural result is based on a careful comparison of the angles of a \ktri{k} with the angles of its \kstar{k}s, where an angle is given by two edges adjacent to a common vertex and consecutive around it. More precisely:

\begin{definition}
\index{angle}
An \defn{angle}~$\angle(u,v,w)$ of a subset~$E$ of~$E_n$ is a pair~$\{[u,v],[v,w]\}$ of edges of~$E$ such that~$u\cl v\cl w$ and for all~$t\in\;\rrb w,u\llb$, the edge~$[v,t]$ is not in~$E$.
\end{definition}

We call~$v$ the \defn{vertex} of the angle~$\angle(u,v,w)$. If~$t\in\;\rrb w,u\llb$, we say that \defn{$t$ is contained} in the angle~$\angle(u,v,w)$. An angle is said to be \defn{\krel{k}}\index{relevant@\krel{k}!--- angle} if both its edges are either \krel{k} or \kbound{k} edges of~$E_n$.

To illustrate these definitions, let us describe the angles of a \kstar{k} (see \fref{stars:fig:starpolygons}(a-c)):

\begin{lemma}\label{stars:lem:anglekstar}
Let~$S$ be a \kstar{k} of $E_n$. Then:
\begin{enumerate}[(i)]
\item the angles of~$S$ are the $2k+1$~pairs of consecutive edges of~$S$ in star order;
\item all angles of~$S$ are \krel{k};~and
\item any vertex~$t$ not in~$S$ is contained in a unique angle~$\angle(u,v,w)$ of~$S$.
\end{enumerate}
\end{lemma}

\begin{proof}
Part~(i) is nothing else but the definition of angle.

Since any edge of~$S$ separates the other vertices of~$S$ into two parts of size~$k-1$ and~$k$ respectively, it is at least a \kbound{k} edge, which implies Part~(ii).

To prove Part~(iii), denote by~$\ens{s_i}{i\in\Z_{2k+1}}$ the vertices of~$S$ in cyclic order. Given a vertex~$t$ not in~$S$, there exists a unique~$i\in\Z_{2k+1}$ such that~$s_i\cl t\cl s_{i+1}$. The angle~$\angle(s_i,s_{i-k},s_{i+1})$ is the unique angle containing~$t$.
\end{proof}


\subsection{Angles of \kstar{k}s \vs angles of \ktri{k}s}\label{stars:subsec:angles:stardecomposition}

In order to compare the angles of a \ktri{k} with the angles of its \kstar{k}s, we start with the following simple observation:

\begin{lemma}\label{stars:lem:angles}
Let~$S$ be a \kstar{k} of a \kcross{(k+1)}-free subset~$E$ of~$E_n$. Then any angle of~$S$ is also a (\krel{k}) angle of~$E$.
\end{lemma}

\begin{proof}
Let~$\ens{s_j}{j\in\Z_{2k+1}}$ denote the vertices of~$S$ in star order. Suppose that~$E$ contains an edge~$[s_j,t]$ where~$j\in\Z_{2k+1}$ and~$t\in\;\rrb s_{j+1},s_{j-1}\llb$. Then the set of edges
$$\{[s_{j+1},s_{j+2}],[s_{j+3},s_{j+4}],\dots,[s_{j-2},s_{j-1}],[s_j,t]\}$$
forms a \kcross{(k+1)}. Consequently~$\angle(s_{j-1},s_j,s_{j+1})$ is an angle of~$E$.
\end{proof}

It is interesting to observe that the reciprocal statement of this lemma is true for triangulations of the \gon{n}: any angle of a triangulation belongs to exactly one of its triangles. The similar statement for multitriangulations is much harder to prove, and is the basis of our study of the structural properties of multitriangulations:

\begin{theorem}\label{stars:theo:angles}
Any \krel{k} angle of a \ktri{k}~$T$ belongs to a unique \kstar{k} of~$T$.
\end{theorem}

\begin{proof}
In this proof, we need the following definition (see \fref{stars:fig:farther}): let~$\angle(u,v,w)$ be a \krel{k} angle of~$T$ and let~$e$ and~$f$ be two edges of~$T$ that \defn{intersect}~$\angle(u,v,w)$ (\ie that intersect both~$[u,v]$ and~$[v,w]$). If~$a$,~$b$,~$c$ and~$d$ denote their vertices such that~$u\cl a\cl v\cl b\cl w$ and~$u\cl c\cl v\cl d\cl w$, then we say that~$e=[a,b]$ is \defn{\vfarther{v}} than~$f=[c,d]$ if~$u\cl a\cle c\cl v\cl d\cle b\cl w$. Let~$E$ and~$F$ be two \kcross{(k-1)}s that \defn{intersect}~$\angle(u,v,w)$. Let their edges be labeled~$e_1 \eqdef [a_1,b_1],e_2 \eqdef [a_2,b_2],\dots,e_{k-1} \eqdef [a_{k-1},b_{k-1}]$ and~$f_1 \eqdef [c_1,d_1]$, $f_2 \eqdef [c_2,d_2],\dots,f_{k-1} \eqdef [c_{k-1},d_{k-1}]$ such that~$u\cl a_1\cl a_2\cl \cdots\cl a_{k-1}\cl v\cl b_1\cl b_2\cl \cdots\cl b_{k-1}\cl w$ and~$u\cl c_1\cl c_2\cl \cdots\cl c_{k-1}\cl v\cl d_1\cl d_2\cl \cdots\cl d_{k-1}\cl w$. Then we say that~$E$ is \defn{\vfarther{v}} than~$F$ if~$e_i$ is \vfarther{v} than~$f_i$ for every~$i\in[k-1]$. We say that~$E$ is \defn{\vmax{v}} if it does not exist any \kcross{(k-1)} intersecting~$\angle(u,v,w)$ and \vfarther{v} than~$E$.

\begin{figure}[h]
	\capstart
	\centerline{\includegraphics[scale=1]{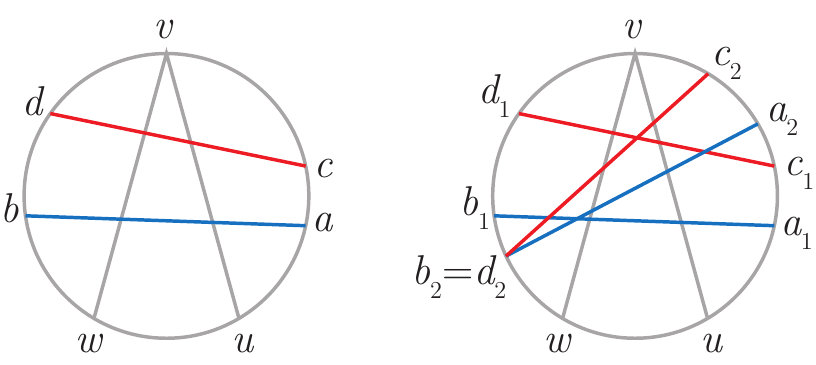}}
	\caption[The \vfarther{v} relation]{$[a,b]$ is \vfarther{v} than $[c,d]$ (left) and $\{[a_i,b_i]\}$ is \vfarther{v} than $\{[c_i,d_i]\}$ (right).}
	\label{stars:fig:farther}
\end{figure}

Let~$\angle(u,v,w)$ be a \krel{k} angle of~$T$. If the edge~$[u,v+1]$ is in~$T$, then the angle $\angle(v+1,u,v)$ is a \krel{k} angle of~$T$ and, if it is contained in a \kstar{k}~$S$ of~$T$, then so is~$\angle(u,v,w)$. Moreover, if~$n>2k+1$ ($n=2k+1$~is a trivial case), $T$~cannot contain all the edges~$\ens{[u+i,v+i]}{0\le i\le n-1}$ and~$\ens{[u+i,v+i+1]}{0\le i\le n-1}$. Consequently, we can assume that~$[u,v+1]$ is not in~$T$.

Thus we have a \kcross{k}~$E$ of the form~$e_1 \eqdef [a_1,b_1],\dots,e_k \eqdef [a_k,b_k]$ with~$u\cl a_1\cl \cdots\cl a_k\cl v+1$ and~$v+1\cl b_1\cl \cdots\cl b_k\cl u$. Since~$[u,v]\in T$, $a_k=v$~and since~$\angle(u,v,w)$ is an angle,~$v+1\cl b_k\cle w$. Consequently,~$\{e_1,\dots,e_{k-1}\}$ forms a \kcross{(k-1)} intersecting~$\angle(u,v,w)$ and we can assume that it is \vmax{v} (see \fref{stars:fig:twosteps}(a)). We will prove that the edges~$[u,b_1], [a_1,b_2],\dots, [a_{k-2},b_{k-1}],[a_{k-1},w]$ are in~$T$ such that the points~$u$, $a_1,\dots,a_{k-1}$, $v$, $b_1,\dots,b_{k-1}$, $w$ are the vertices of a \kstar{k} of~$T$ containing the angle~$\angle(u,v,w)$. To get this result, we use two steps: first we prove that~$\angle(a_1,b_1,u)$ is an angle of~$T$, and then we prove that the edges~$e_2,\dots,e_{k-1},[v,w]$ form a \kcross{(k-1)} intersecting~$\angle(a_1,b_1,u)$ and \vmax{b_1} (so that we can reiterate the argument).

\begin{figure}
	\capstart
	\centerline{\includegraphics[scale=1]{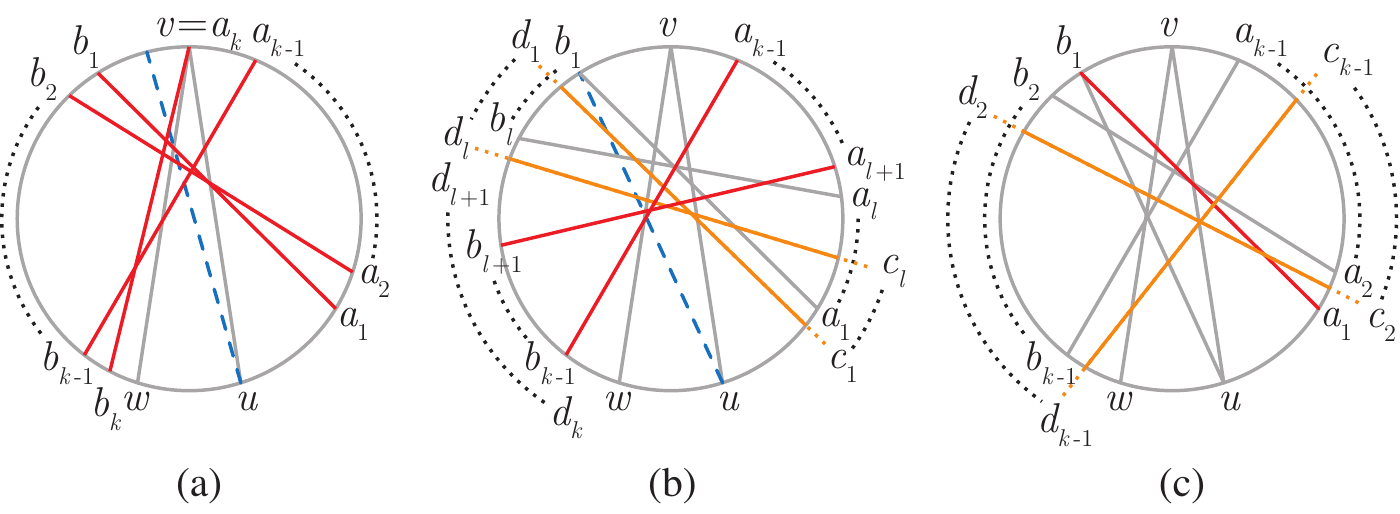}}
	\caption[Illustration of the proof of Theorem~\ref{stars:theo:angles}]{(a) The \kcross{k}~$E$; (b) $[u,b_1]\in T$; (c) $\{e_2,\dots,e_{k-1},[v,w]\}$ is \vmax{b_1}.}
	\label{stars:fig:twosteps}
\end{figure}

\paragraph{First step.} See \fref{stars:fig:twosteps}(b).

\noindent Suppose that~$[u,b_1]$ is not in~$T$. Thus we have a \kcross{k}~$F$ that prevents the edge~$[u,b_1]$. Let~$f_1 \eqdef [c_1,d_1],\dots,f_k \eqdef [c_k,d_k]$ denote its edges with~$u\cl c_1\cl \cdots\cl c_k\cl b_1$ and $b_1\cl d_1\cl \cdots\cl d_k\cl u$.

Note first that~$v\cl d_k\cle w$. Otherwise~$d_k\in\;\rrb w,u\llb$ and~$c_k\ne v$, because~$\angle(u,v,w)$ is an angle. Thus either~$c_k\in\;\rrb u,v\llb$ and then~$F\cup\{[u,v]\}$ forms a \kcross{(k+1)}, or~$c_k\in\;\rrb v,b_1\llb$ and then~$E\cup\{[c_k,d_k]\}$ forms a \kcross{(k+1)}. Consequently, we have~$b_1\cl d_1\cl \cdots\cl d_{k-1}\cl w$.

Let~$\ell \eqdef \max\ens{j\in [k-1]}{b_i\cl d_i\cl w \text{ for any } i \text{ with } i\in[j]}$. Then for any~$i\in[\ell]$, since $\{e_1,\dots,e_i\}\cup\{f_i,\dots,f_k\}$ does not form a \kcross{(k+1)}, we have~$u\cl c_i\cle a_i$. Thus for any~$i\in[\ell]$,~$u\cl c_i\cle a_i\cl v\cl b_i\cl d_i\cl w$, so that~$f_i$ is \vfarther{v} than~$e_i$. Furthermore, we have~$u\cl c_1\cl \cdots\cl c_\ell\cl a_{\ell+1}\cl\cdots\cl a_{k-1}\cl v\cl d_1\cl\cdots\cl d_\ell\cl b_{\ell+1}\cl \cdots\cl b_{k-1}\cl w$. Consequently, we get a \kcross{(k-1)}~$\{f_1,\dots,f_\ell,e_{\ell+1},\dots,e_{k-1}\}$ which is \vfarther{v} than~$\{e_1,\dots,e_{k-1}\}$; this contradicts the definition of~$\{e_1,\dots,e_{k-1}\}$. Thus we obtain that~$[u,b_1]\in T$.

Suppose now that~$\angle(a_1,b_1,u)$ is not an angle of~$T$. Then there exists~$a_0\in\;\rrb u,a_1\llb$ such that~$[b_1,a_0]\in T$. But then the \kcross{(k-1)}~$\{[a_0,b_1],e_2,\dots,e_{k-1}\}$ is \vfarther{v} than the \kcross{(k-1)}~$\{e_1,\dots,e_{k-1}\}$. This implies that~$\angle(a_1,b_1,u)$ is an angle of~$T$.

\paragraph{Second step.} See \fref{stars:fig:twosteps}(c).

\noindent Assume that there exists a \kcross{(k-1)}~$F$ intersecting~$\angle(a_1,b_1,u)$ and \vfarther{b_1} than the \kcross{(k-1)} $\{e_2,\dots,e_{k-1},[v,w]\}$. Let~$f_2 \eqdef [c_2,d_2],\dots,f_k \eqdef [c_k,d_k]$ denote its edges, with $a_1\cl c_2\cl \cdots\cl c_k\cl b_1\cl d_2\cl \cdots\cl d_k\cl u$. 

Note first that~$b_k\cle d_k\cle w$. Otherwise~$d_k\in\;\rrb w,u\llb$ and~$c_k\ne v$, because~$\angle(u,v,w)$ is an angle. Thus either~$c_k\in\;\rrb a_1,v\llb$ and then~$F\cup\{[u,v],e_1\}$ forms a \kcross{(k+1)}, or~$c_k\in\;\rrb v,b_1\llb$ and then~$E\cup\{[c_k,d_k]\}$ forms a \kcross{(k+1)}. Thus, we have~$b_1\cl d_2\cl\dots\cl d_{k-1}\cl w$.

Furthermore, for any~$2\le i\le k-1$, the edge $f_i$~is \vfarther{\angle(a_1,b_1,u)} than~$e_i$, so that $a_1\cl c_i\cle a_i\cl b_1\cl b_i\cle d_i\cl u$. In particular,~$a_1\cl c_{k-1}\cle a_{k-1}\cl v$ and we obtain that $u\cl a_1\cl c_2\cl\cdots\cl c_{k-1}\cl v\cl b_1\cl d_2\cl \cdots\cl d_{k-1}\cl w$. Consequently, the \kcross{(k-1)}~$\{e_1,f_2,\dots,f_{k-1}\}$ is \vfarther{v} than~$\{e_1,\dots,e_{k-1}\}$, which is a contradiction.
\end{proof}

As an immediate corollary of Theorem~\ref{stars:theo:angles}, we obtain the incidence relations between the \kstar{k}s and the edges of a \ktri{k}:

\begin{corollary}\label{stars:coro:incidences}
Let~$e$ be an edge of a \ktri{k}~$T$.
\begin{enumerate}[(i)]
\item If~$e$ is a \krel{k} edge, then it belongs to exactly two \kstar{k}s of~$T$ (one on each side).
\item If~$e$ is a \kbound{k} edge, then it belongs to exactly one \kstar{k} of~$T$ (on its ``inner'' side).
\item If~$e$ is a \kirrel{k} edge, then it does not belong to any \kstar{k} of~$T$.\qed
\end{enumerate}
\end{corollary}

We will see many repercussions of Corollary~\ref{stars:coro:incidences} throughout this chapter: in Corollary~\ref{stars:coro:starsenumeration}, we use it to obtain by double counting both the number of \kstar{k}s and of edges in a \ktri{k} of the \gon{n}; in Section~\ref{stars:sec:flips}, it takes part in the definition of flips for \ktri{k}s; and in Section~\ref{stars:sec:surfaces}, we reinterpret it as a polygonal decomposition on a surface. Using tools developed in Section~\ref{mpt:subsec:duality:mt}, we will also prove its reciprocal statement (see Theorem~\ref{mpt:theo:characterization}):

\begin{theorem}\label{stars:theo:characterization}
Let~$\Sigma$ be a set of \kstar{k}s of the \gon{n} such that:
\begin{enumerate}[(i)]
\item any \krel{k} edge of~$E_n$ is contained in zero or two \kstar{k}s of~$\Sigma$, one on each side;~and
\item any \kbound{k} edge of~$E_n$ is contained in exactly one \kstar{k} of~$\Sigma$.
\end{enumerate}
Then $\Sigma$ is the set of \kstar{k}s of a \ktri{k} of the \gon{n}.\qed
\end{theorem}

\begin{figure}
	\capstart
	\centerline{\includegraphics[scale=1]{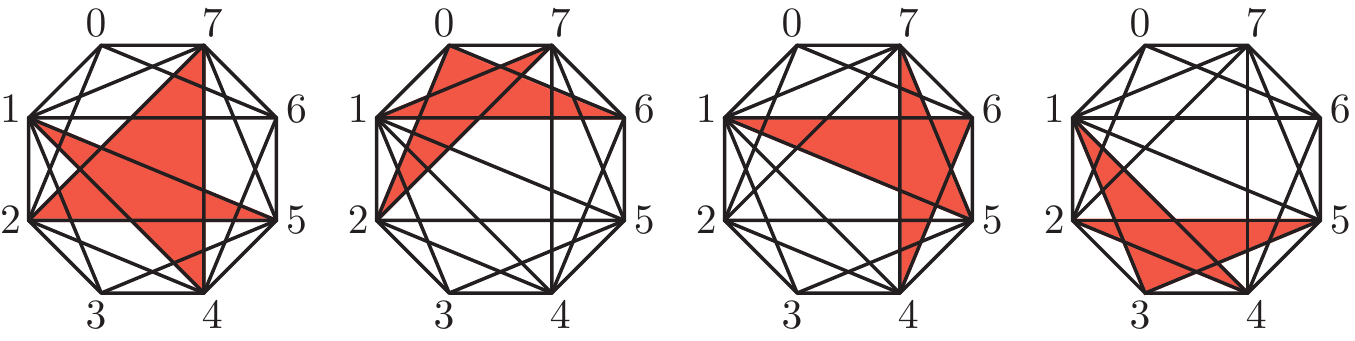}}
	\caption[The four \kstar{2}s in the \ktri{2} of the octagon of \fref{intro:fig:2triang8points}]{The four \kstar{2}s in the \ktri{2} of the octagon of \fref{intro:fig:2triang8points}.}
	\label{stars:fig:2triang8pointsstars}
\end{figure}

\begin{example}
To illustrate these results, \fref{stars:fig:2triang8pointsstars} shows the \kstar{2}s contained in the \ktri{2}~$T$ of the octagon of \fref{intro:fig:2triang8points}. Each of the twenty \krel{2} angles of $T$ is contained in exactly one \kstar{2} of~$T$ (Theorem~\ref{stars:theo:angles}) and every \krel{2} edge is contained in two \kstar{2}s of~$T$ (Corollary~\ref{stars:coro:incidences}).
\end{example}

To conclude, let us insist on the fact that the results of this section interpret \ktri{k}s as \defn{complexes of \kstar{k}s}: a \ktri{k} is just a way to glue together \kstar{k}s along \krel{k} edges. This motivates the title of~\cite{ps-mtcsp-09}.


\subsection{Common bisectors and the number of stars}\label{stars:subsec:angles:commonbisector}

We now want to obtain both the numbers of \kstar{k}s and of edges of a \ktri{k} of the \gon{n}. As in the case of triangulations (see the ``Double counting'' paragraph in Section~\ref{intro:sec:triangulations}), this requires two relations between \kstar{k}s and edges: one is given by Theorem~\ref{stars:theo:angles} and we derive the other one from the study of common bisectors of \kstar{k}s.

\begin{definition}\label{stars:def:bisector}
\index{bisector|hbf}
A \defn{bisector} of an angle~$\angle(u,v,w)$ is an edge~$[v,t]$ where $w\cl t\cl u$. A \defn{bisector} of a \kstar{k} is a bisector of one of its angles.
\end{definition}

Observe that the bisectors of a \kstar{k} are exactly the edges of~$E_n$ that pass through one of its vertices and split the other ones into two sets of size~$k$.

The following theorem will be used extensively in our study of multitriangulations, in particular for our local definition of flips in Section~\ref{stars:sec:flips}.

\begin{theorem}\label{stars:theo:bisector}
Every pair of \kstar{k}s whose union is \kcross{(k+1)}-free have a unique common bisector.
\end{theorem}

Using notions developed in Section~\ref{mpt:sec:duality}, we will obtain a very short proof of the existence of a common bisector. However, we propose here an alternative proof using only our current tools. We need the following two lemmas, in which $E$~denotes a \kcross{(k+1)}-free subset of edges of~$E_n$, and~$S$ denotes a \kstar{k} of~$E$.

\begin{lemma}
The numbers of vertices of~$S$ on each side of an edge~$e$ of~$E$ are different.
\end{lemma}

\begin{proof}
Suppose that~$S$ has the same number of vertices on both sides of an edge~$e$. Since the number of vertices of~$S$ is~$2k+1$, one of the two vertices of~$e$ is a vertex of~$S$, and~$e$ is a bisector of~$S$. Lemma~\ref{stars:lem:angles} then ensures that~$e$ is not in~$E$.
\end{proof}

Let~$V$ be the set of vertices of the \kstar{k}~$S$. If~$|\llb u,v\rrb\cap V|<|\llb v,u\rrb\cap V|$, then we say that~$S$ lies on the \defn{positive side} of the edge~$[u,v]$ oriented from~$u$ to~$v$ (otherwise we say that~$S$ lies on the \defn{negative side} of the edge~$[u,v]$ oriented from~$u$ to~$v$). The \kstar{k}~$S$ is said to be \defn{contained in an angle}~$\angle(u,v,w)$ if it lies on the positive side of both the edges~$[u,v]$ and~$[v,w]$ oriented from~$u$ to~$v$ and from~$v$ to~$w$ respectively.

\begin{lemma}
Let~$\angle(u,v,w)$ be an angle of~$E$ containing the \kstar{k}~$S$. Then:
\begin{enumerate}[(i)]
\item either~$v$ is a vertex of~$S$ and~$\angle(u,v,w)$ is an angle of~$S$;
\item or~$v$ is not a vertex of~$S$ and~$\angle(u,v,w)$ has a common bisector with an angle of~$S$.
\end{enumerate}
\end{lemma}

\begin{proof}
Suppose first that~$v$ is a vertex of~$S$. Let~$\angle(x,v,y)$ denote the angle of~$S$ at vertex~$v$. Since~$\angle(u,v,w)$ contains~$S$, we have~$w\cle y\cl x\cle u$. But since~$\angle(u,v,w)$ is an angle of~$E$, we have~$x=u$ and~$y=w$, so that~$\angle(u,v,w)$ is an angle of~$S$.

Suppose now that~$v$ is not a vertex of~$S$. Then, by Lemma~\ref{stars:lem:anglekstar}(iii) there exists a unique angle~$\angle(x,y,z)$ of~$S$ containing~$v$. If~$y\in\;\rrb u,v\llb$, then~$\rrb u,v\llb$ contains all the~$k+1$ vertices of~$S$ between~$y$ and~$z$, which is not possible (because~$S$ lies on the positive side of the edge~$[u,v]$, oriented from~$u$ to~$v$). For the same reason,~$y\not\in\;\rrb v,w \llb$.
If~$y=u$ or~$y=w$, then~$[u,v]$ or~$[v,w]$ is a bisector of~$\angle(x,y,z)$, which contradicts Lemma~\ref{stars:lem:angles}. Consequently, $\angle(u,v,w)$~contains~$y$, so that~$[v,y]$ is a common bisector of~$\angle(u,v,w)$ and~$\angle(x,y,z)$.
\end{proof}

\begin{proof}[Proof of Theorem~\ref{stars:theo:bisector}]
Let~$R$ and~$S$ be two \kstar{k}s whose union is \kcross{(k+1)}-free.

Let~$\ens{r_j}{j\in\Z_{2k+1}}$ denote the vertices of~$R$ in star order. Note that for any~$j\in\Z_{2k+1}$, if~$S$ lies on the negative side of the edge~$[r_{j-1},r_j]$ oriented from~$r_{j-1}$ to~$r_j$, then it lies on the positive side of the edge~$[r_j,r_{j+1}]$ oriented from~$r_j$ to~$r_{j+1}$. Since~$2k+1$ is odd, this implies that there exists~$j\in\Z_{2k+1}$ such that~$S$ lies on the positive side of both the edges~$[r_{j-1},r_j]$ and~$[r_j,r_{j+1}]$ oriented from~$r_{j-1}$ to~$r_j$ and from~$r_j$ to~$r_{j+1}$, respectively. That is, in the angle~$\angle(r_{j-1},r_j,r_{j+1})$. The previous lemma thus ensures the existence of a common bisector.

\begin{figure}
	\capstart
	\centerline{\includegraphics[scale=1]{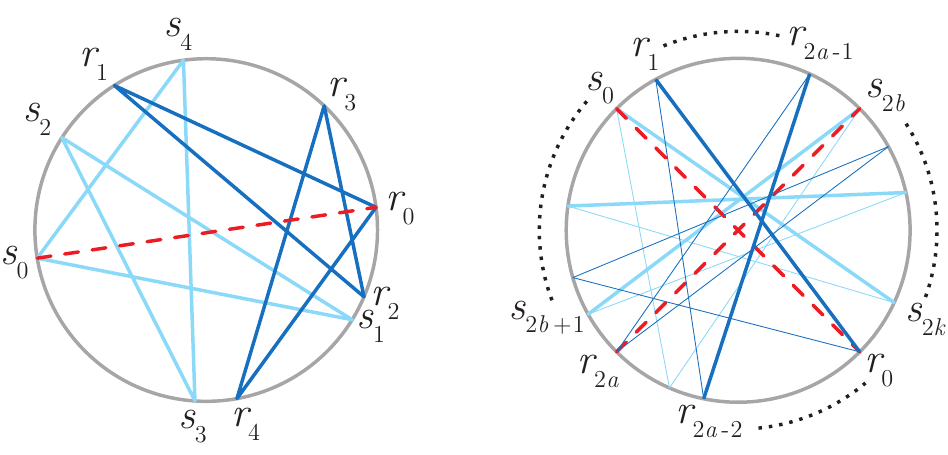}}
	\caption[Commmon bisectors of \kstar{k}s]{The unique common bisector of a \kcross{3}-free pair of \kstar{2}s (left) and a \kcross{(k+1)} in a pair of \kstar{k}s with two common bisectors (right).}
	\label{stars:fig:uniquecommonbisector}
\end{figure}

Suppose now that the two \kstar{k}s~$R$ and~$S$ have two different common bisectors~$e$ and~$f$, and let~$\ens{r_j}{j\in\Z_{2k+1}}$ and~$\ens{s_j}{j\in\Z_{2k+1}}$ denote the vertices of~$R$ and~$S$ in star order and labeled such that~$e=[r_0,s_0]$. Let~$\alpha,\beta\in\Z_{2k+1}$ be such that~$f=[r_\alpha,s_\beta]$. Note that certainly,~$\alpha\ne0$, ${\beta\ne0}$, and~$\alpha$ and~$\beta$ have the same parity. By symmetry, we can assume that~$\alpha=2a$,~$\beta=2b$ with $1\le b\le a\le k$. But then the set
$$\ens{[r_{2i},r_{2i+1}]}{0\le i\le a-1}\cup\ens{[s_{2j},s_{2j+1}]}{b\le j\le k}$$
forms a \kcross{(k+1+a-b)}, and~$k+1+a-b\ge k+1$. This proves uniqueness.
\end{proof}

\begin{proposition}\label{stars:prop:findingstars}
Let~$T$ be a \ktri{k}.
\begin{enumerate}[(i)]
\item For any \kstar{k}~$S$ in~$T$ and for any vertex~$r$ not in~$S$ there is a unique \kstar{k}~$R$ in~$T$ such that~$r$ is a vertex of the common bisector of~$R$ and~$S$.
\item Any \krel{k} edge not in~$T$ is the common bisector of a unique pair of \kstar{k}s of~$T$.
\end{enumerate}
\end{proposition}

\begin{proof}
Let~$\angle(u,s,v)$ be the unique angle of~$S$ which contains~$r$ (Lemma~\ref{stars:lem:anglekstar}(iii)). Let~$\angle(x,r,y)$ be the unique angle of~$T$ of vertex~$r$ which contains~$s$. According to Theorem~\ref{stars:theo:angles}, the angle $\angle(x,r,y)$ belongs to a unique \kstar{k}~$R$. The common bisector of~$R$ and~$S$ is~$[r,s]$ and~$R$ is the only such \kstar{k} of~$T$.

Let~$e \eqdef [r,s]$ be a \krel{k} edge, not in~$T$. Let~$\angle(x,r,y)$ (resp.~$\angle(u,s,v)$) denote the unique angle of~$T$ of vertex~$r$ (resp.~$s$) which contains~$s$ (resp.~$r$). According to Theorem~\ref{stars:theo:angles}, the angle~$\angle(x,r,y)$ (resp.~$\angle(u,s,v)$) belongs to a unique \kstar{k}~$R$ (resp.~$S$) of~$T$. The common bisector of~$R$ and~$S$ is~$[r,s]$ and~$(R,S)$ is the only such couple of \kstar{k}s of~$T$.
\end{proof}

Observe that Parts~(i)~and~(ii)~of this proposition give bijections between:
\begin{enumerate}[(i)]
\item ``vertices not used in the \kstar{k}~$S$ of~$T$'' and ``\kstar{k}s of~$T$ different from~$S$'';
\item ``\krel{k} edges not used in~$T$'' and ``pairs of \kstar{k}s of~$T$''.
\end{enumerate}
From any of these two bijections, and using Corollary~\ref{stars:coro:incidences} for double counting, it is easy to derive the number of \kstar{k}s and of edges in~$T$:

\begin{corollary}\label{stars:coro:starsenumeration}
\begin{enumerate}[(i)]
\item Any \ktri{k} of the \gon{n} contains~$n-2k$ \kstar{k}s, $k(n-2k-1)$ \krel{k} edges and~$k(2n-2k-1)$ edges.
\item \ktri{k}s  are exactly \kcross{(k+1)}-free subsets of~$k(2n-2k-1)$ edges of~$E_n$.
\end{enumerate}
\end{corollary}

\begin{proof}
Let~$T$ be a \ktri{k} of the \gon{n}. Let~$\sigma$ denote the number of \kstar{k}s of~$T$ and~$\rho$ denote its number of \krel{k} edges. Corollary~\ref{stars:coro:incidences} (or equivalently Theorem~\ref{stars:theo:angles}) ensures that $(2k+1)\sigma=2\rho+n$ and Proposition~\ref{stars:prop:findingstars}(ii) affirms that $nk+\rho = {n \choose 2}-{\sigma \choose 2}$. Thus $\sigma$ satisfies $\sigma^2 +2k\sigma-n(n-2k)=0$, which has $\sigma=n-2k$ as unique (positive) solution. Hence, $\rho=k(n-2k-1)$. Finally, the total number of edges of~$T$ is $nk+\rho=k(2n-2k-1)$.
\end{proof}


\section{Flips}\label{stars:sec:flips}

In this section, we introduce flips for multitriangulations: we transform a \ktri{k} into another one by just exchanging one \krel{k} edge. Even if this transformation was already known and used in~\cite{n-gdfcp-00,dkm-lahp-02,j-gt-03}, we provide a new definition based on \kstar{k}s, which involves only a substructure of the \ktri{k} during the transformation. In the second part of this section, we study the graph of flips and discuss upper and lower bounds on its diameter.


\subsection{Mutual position of two \kstar{k}s and the local definition of flips}\label{stars:subsec:flips:position}

To start this section, we consider two \kstar{k}s~$R$ and~$S$ of a \kcross{(k+1)}-free subset~$E$ of~$E_n$ and we study their mutual position.

\begin{lemma}\label{stars:lem:angledisjoint}
$R$~and~$S$ cannot share any angle.
\end{lemma}

\begin{proof}
By Lemma~\ref{stars:lem:angles}, the knowledge of one angle~$\angle(s_{j-1},s_j,s_{j+1})$ of~$S$ enables us to recover all the \kstar{k}~$S$: the vertex~$s_{j+2}$ is the unique vertex such that~$\angle(s_j,s_{j+1},s_{j+2})$ is an angle of~$E$ (\ie the first neighbour of~$s_{j+1}$ after~$s_j$ when moving clockwise), and so on.
\end{proof}

\begin{corollary}\label{stars:coro:commonedges}
$R$~and~$S$ cannot share more than~$k$ edges.\qed
\end{corollary}

Observe that, for example, the two \kstar{k}s of any \ktri{k} of a \gon{(2k+2)} share exactly~$k$ edges (see \fref{stars:fig:2k+2}).

By Theorem~\ref{stars:theo:bisector}, we know that $R$~and~$S$ have a unique common bisector $e$. In the following lemmas, we are interested in the position of their vertices relatively to this bisector. We denote by $\ens{r_j}{j\in\Z_{2k+1}}$ and $\ens{s_j}{j\in\Z_{2k+1}}$ the vertices of~$R$~and~$S$, in star order, and such that the edge~$e \eqdef [r_0,s_0]$ is the common bisector of~$R$~and~$S$.

\begin{lemma}\label{stars:lem:paralleledges}
For every $i\in[k]$, we have $r_{2i-1}\in\;\rrb r_0,s_{2i}\rrb$~and~$s_{2i-1}\in\;\rrb s_0,r_{2i}\rrb$. In particular, for every~$i\in\Z_{2k+1}$ the edges~$[r_i,r_{i+1}]$~and~$[s_i,s_{i+1}]$ do not cross.
\end{lemma}

\begin{proof}
Suppose that there exists $i\in[k]$ such that $r_{2i-1}\in\;\rrb s_{2i},s_0\llb$~or~$s_{2i-1}\in\;\rrb r_{2i},r_0\llb$. Let~$\gamma$ be the highest such integer, and assume for example that~$r_{2\gamma-1}\in\;\rrb s_{2\gamma},s_0\llb$. Then the definition of~$\gamma$ ensures that~$s_0\cl s_{2\gamma+1}\cle r_{2\gamma+2}\cl r_{2\gamma-2}\cl r_0\cl s_{2k}\cl s_{2\gamma}\cl r_{2\gamma-1}\cl r_1\cl s_0$~so that the set
$$\ens{[r_{2i},r_{2i+1}]}{0\le i\le\gamma-1}\cup\ens{[s_{2j},s_{2j+1}]}{\gamma\le j\le k}$$
forms a \kcross{(k+1)}.
\end{proof}

\begin{figure}
	\capstart
	\centerline{\includegraphics[scale=1]{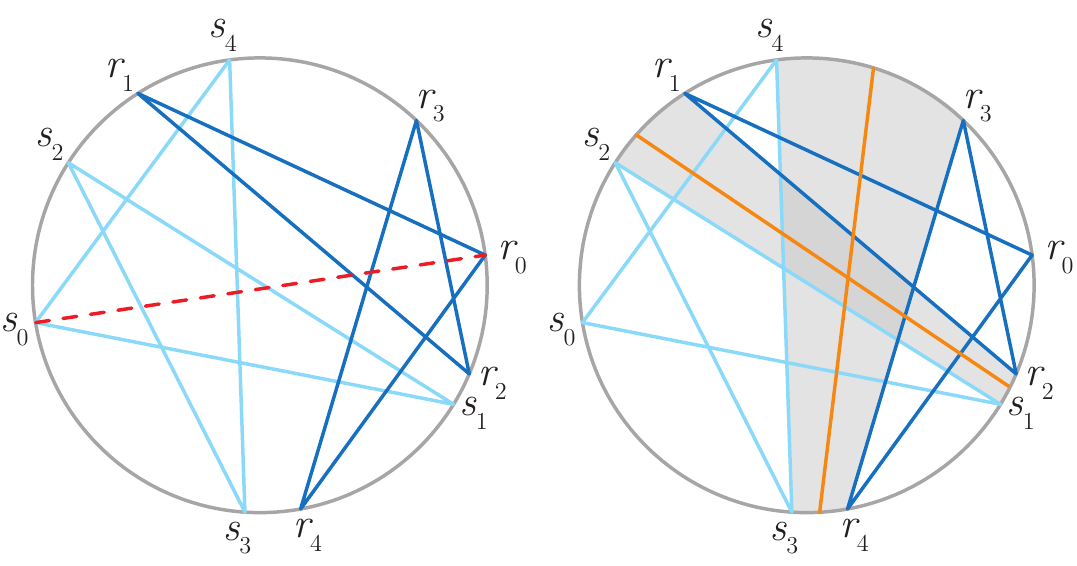}}
	\caption[Mutual position of two \kstar{k}s]{The common bisector of two \kstar{2}s (left) and a \kcross{2} that crosses it (right).}
	\label{stars:combisector}
\end{figure}

The previous lemma can be read as saying that corresponding edges of~$R$~and~$S$ are parallel. In fact,~$k$ of these~$2k+1$ pairs of parallel edges, the ones of the form~$([r_{2i-1},r_{2i}],[s_{2i-1},s_{2i}])$, with $i\in[k]$, separate~$R$~from~$S$, meaning that~$R$~and~$S$ lie on opposite sides of both. The next lemma says that any \kcross{k} that, in turn, crosses the common bisector~$e=[r_0,s_0]$, has one edge parallel to and between each such pair (see \fref{stars:combisector}).

\begin{lemma}\label{stars:lem:sandwich}
Let~$F$ be a \kcross{k} of~$E$ such that all its edges cross~$e=[r_0,s_0]$. Denote by $f_1 \eqdef [x_1,y_1],\dots,f_k \eqdef [x_k,y_k]$  the edges of~$F$, such that $r_0\cl x_1\cl \cdots\cl x_k\cl s_0\cl y_1\cl \cdots\cl y_k\cl r_0$. Then~$x_i\in\llb r_{2k-2i+1},s_{2k-2i+2}\rrb$ and~$y_i\in\llb s_{2k-2i+1},r_{2k-2i+2}\rrb$, for any $i\in[k]$.
\end{lemma}

\begin{proof}
Suppose that there exists~$i\in[k]$ such that~$r_0\cl x_i\cl r_{2k-2i+1}$ and let
$$\ell \eqdef \max\ens{i\in[k]}{r_0\cl x_i\cl r_{2k-2i+1}}.$$
If~$\ell=k$, then the set~$\{f_1,\dots,f_k,[r_0,r_1]\}$ is a \kcross{(k+1)} of~$E$, thus we assume that~$\ell<k$. In order for the set
$$\{f_1,\dots,f_\ell\}\cup\{[r_0,r_1],\dots,[r_{2k-2\ell},r_{2k-2\ell+1}]\}$$
not to be a \kcross{(k+1)}, we have~$r_{2k-2\ell}\cle y_\ell\cl r_0$, so that~$r_{2k-2\ell}\cl y_{\ell+1}\cl r_0$. But the definition of~$\ell$ implies that~$r_{2k-2\ell+1}\cl r_{2k-2\ell-1}\cle x_{\ell+1}\cl s_0$, so that the set
$$\{[r_{2k-2\ell},r_{2k-2\ell+1}],\dots,[r_{2k-2},r_{2k-1}],[r_{2k},r_0]\}\cup\{f_{\ell+1},\dots,f_k\}$$
is a \kcross{(k+1)} of~$E$.
By symmetry, the lemma is proved.
\end{proof}

\begin{figure}
	\capstart
	\centerline{\includegraphics[scale=1]{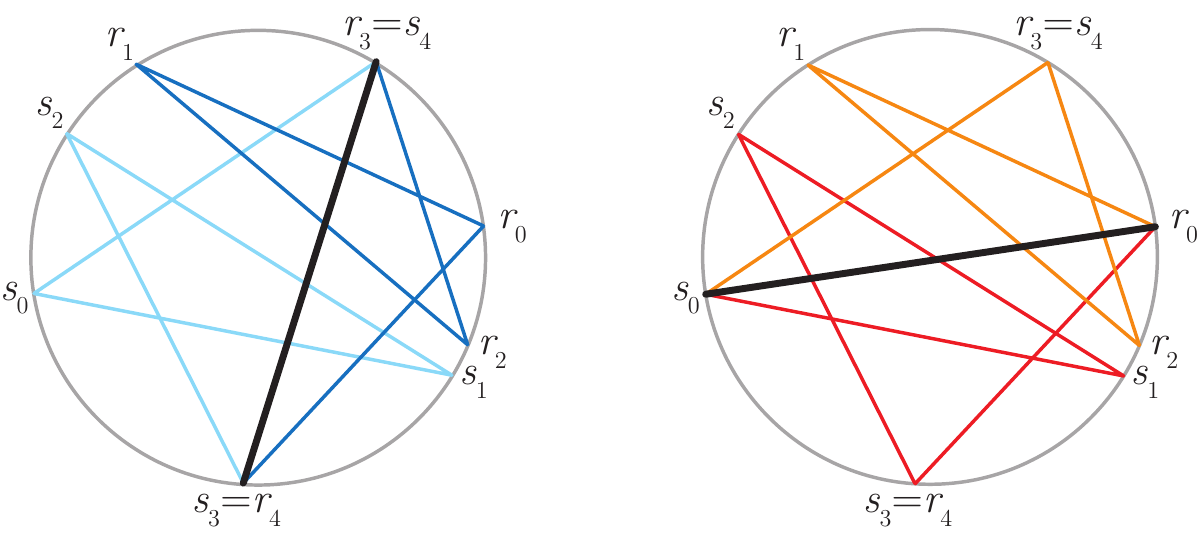}}
	\caption[The flip operation]{The flip of an edge.}
	\label{stars:fig:flip}
\end{figure}

In the following lemma, the symbol $\diffsym$ denotes the symmetric difference between two sets: $X\diffsym Y \eqdef (X\ssm Y)\cup(Y\ssm X)$. We refer to \fref{stars:fig:flip} for an illustration of this lemma.

\begin{lemma}\label{stars:lem:commonedge}
Let~$f$ be a common edge of~$R$~and~$S$. Then:
\begin{enumerate}[(i)]
\item there exists $i\in[k]$ such that~$f=[r_{2i-1},r_{2i}]=[s_{2i},s_{2i-1}]$;
\item $E\diffsym\{e,f\}$ is a \kcross{(k+1)}-free subset of~$E_n$;
\item the vertices $s_0,\dots,s_{2i-2},s_{2i-1}=r_{2i},r_{2i+1},\dots,r_{2k},r_0$ (resp.~$r_0,\dots,r_{2i-2},r_{2i-1}=s_{2i}$, $s_{2i+1},\dots,s_{2k},s_0$) are the vertices of a \kstar{k}~$X$ (resp.~$Y$) of~$E\diffsym\{e,f\}$, in star order;
\item $X$~and~$Y$ share the edge~$e$ and their common bisector is~$f$.
\end{enumerate}
\end{lemma}

\begin{proof}
Let~$u$~and~$v$ denote the vertices of~$f$. Lemma~\ref{stars:lem:paralleledges} ensures that~$\{r_0,s_0\}\cap\{u,v\}=\emptyset$ so that we can assume~$r_0\cl u\cl s_0\cl v\cl r_0$. Consequently, there exists $i,j\in[k]$ such that~$u=r_{2i-1}=s_{2j}$~and~$v=r_{2i}=s_{2j-1}$. Suppose that~$i>j$. Then according to Lemma~\ref{stars:lem:paralleledges}, we have~$r_0\cl r_{2i-1}\cle s_{2i}\cl s_{2j}=r_{2i-1}$ which is impossible. By symmetry, we obtain that~$i=j$ and $f=[r_{2i-1},r_{2i}]=[s_{2i},s_{2i-1}]$.

Lemma~\ref{stars:lem:sandwich} then proves that any \kcross{k} of~$E$ that prevents~$e$ being in~$E$ contains~$f$, so that~$E\diffsym\{e,f\}$ is \kcross{(k+1)}-free.

Let~$L$ be the list of vertices~$(s_0,\dots,s_{2i-2},r_{2i},\dots,r_{2k},r_0)$. Between two consecutive elements of~$L$ lie exactly~$k-1$ others points of~$L$ (for the circle order). This implies that~$L$ is in star order. Part (iii) thus follows from the fact that any edge connecting two consecutive points of~$L$ is in~$E\diffsym\{e,f\}$.

Finally, the edge~$e$ is clearly common to~$X$~and~$Y$. The edge~$f$ is a bisector of both angles $\angle(r_{2i-2},r_{2i-1},s_{2i+1})$~and~$\angle(s_{2i-2},s_{2i-1},r_{2i+1})$, so that it is the common bisector of~$X$~and~$Y$.
\end{proof}

This lemma is at the heart of the concept of flips between \ktri{k}s. Assume that~$T$ is a \ktri{k} of the \gon{n}, and let~$f$ be a \krel{k} edge of~$T$. Let~$R$~and~$S$ be the two \kstar{k}s of~$T$ containing~$f$ (Corollary~\ref{stars:coro:incidences}), and let~$e$ be the common bisector of~$R$~and~$S$ (Theorem~\ref{stars:theo:bisector}). Lemma~\ref{stars:lem:commonedge} affirms that~$T\diffsym\{e,f\}$ is a \kcross{(k+1)}-free subset of~$E_n$. Observe moreover that~$T\diffsym\{e,f\}$ is maximal: this follows from Corollary~\ref{stars:coro:starsenumeration}, but also from the fact that if~$T\diffsym\{e,f\}$ is properly contained in a \ktri{k}~$\widetilde{T}$, then~$\widetilde{T}\diffsym\{e,f\}$ is \kcross{(k+1)}-free and properly contains~$T$. Thus,~$T\diffsym\{e,f\}$ is a \ktri{k} of the \gon{n}.

\begin{lemma}\label{stars:lem:flip}
$T$~and~$T\diffsym\{e,f\}$ are the only \ktri{k}s of the \gon{n} containing~${T\ssm\{f\}}$.
\end{lemma}

\begin{proof}
Let~$e'$ be any edge of~$E_n\ssm T$ distinct from~$e$. Let~$R'$~and~$S'$ be the two \kstar{k}s with common bisector~$e'$ (Proposition~\ref{stars:prop:findingstars}(ii)). Either~$R'$ or~$S'$ (or both) does not contain~$f$, say~$R'$. Then~$R'\cup\{e'\}$ is contained in~$T\diffsym\{e',f\}$ and forms a \kcross{(k+1)}.
\end{proof}

We say that we obtain the \ktri{k}~$T\diffsym\{e,f\}$ from the \ktri{k}~$T$ by \defn{flipping} the edge~$f$. We insist on the definition using the common bisector of~$R$~and~$S$ (see \fref{stars:fig:flip}):

\begin{definition}\label{stars:def:flip}
\index{flip|hbf}
The \defn{flip} of a \krel{k} edge~$f$ in a \ktri{k}~$T$ of the \gon{n} is the transformation which replaces the edge~$f$ by the common bisector of the two \kstar{k}s of~$T$ containing~$f$.
\end{definition}

Observe that Lemma~\ref{stars:lem:flip} could be an alternative definition of flips (and it is the definition in all previous works involving flips in multitriangulations; see~\cite{n-gdfcp-00,dkm-lahp-02,j-gt-03,j-gtdfssp-05}). However, our definition is interesting theoretically as well as practically because it is a local definition: it only involves a subconfiguration formed by two adjacent \kstar{k}s. Given a \krel{k} edge~$f$ of a \ktri{k}~$T$, we need a constant time to flip~$f$ in~$T$:
\begin{enumerate}[(i)]
\item first, we compute the two \kstar{k}s of~$T$ containing~$f$ rotating around their angles (remember the proof of Lemma~\ref{stars:lem:angledisjoint});
\item then we compute the common bisector of these two \kstar{k}s and substitute it to~$f$.
\end{enumerate}
Observe the difference with the method consisting in testing, among all edges not in~$T$, which edge could be added to~$T\ssm\{f\}$ without creating a \kcross{(k+1)}.


\subsection{The graph of flips}\label{stars:subsec:flips:graph}

\index{flip!graph of ---s}
Let~$G_{n,k}$ be the \defn{graph of flips} on the set of \ktri{k}s of the \gon{n}, \ie the graph whose vertices are the \ktri{k}s of the \gon{n} and whose edges are the pairs of \ktri{k}s related by a flip.
It follows from Corollary~\ref{stars:coro:starsenumeration} and Lemma~\ref{stars:lem:flip} that~$G_{n,k}$ is regular of degree~$k(n-2k-1)$: every \krel{k} edge of a \ktri{k} can be flipped, in a unique way. In this section, we prove the connectedness of this graph and bound its diameter~$\delta_{n,k}$.

\begin{figure}
	\capstart
	\centerline{\includegraphics[scale=1]{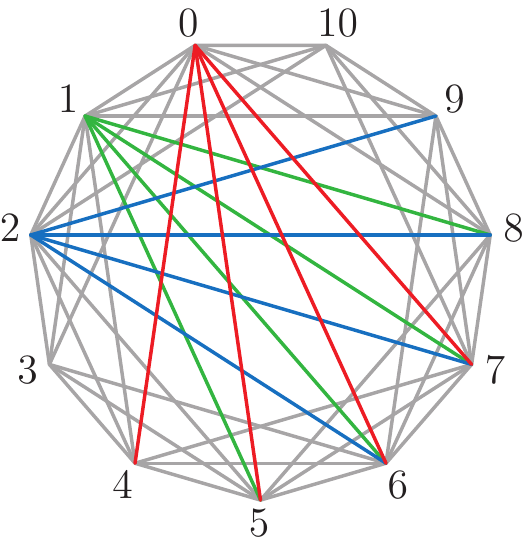}}
	\caption[The minimal triangulation $T_{11,3}^{\min}$]{The minimal triangulation $T_{11,3}^{\min}$.}
	\label{stars:fig:minimalktriangulation}
\end{figure}

Let~$T$ be a \ktri{k} of the \gon{n}, $f$~be a \krel{k} edge of~$T$ and~$e$~be the common bisector of the two \kstar{k}s of~$T$ containing~$f$. Then the edges~$e$~and~$f$ necessarily cross. In particular, if~$e \eqdef [\alpha,\beta]$ and~$f \eqdef [\gamma,\delta]$, with $0\cle\alpha\cl\beta\cle n-1$ and $0\cle\gamma\cl\delta\cle n-1$, then
\begin{enumerate}[(i)]
\item either $0\cle\alpha\cl\gamma\cl\beta\cl\delta\cle n-1$ and the flip is said to be \defn{slope-decreasing},
\item or $0\cle\gamma\cl\alpha\cl\delta\cl\beta\cle n-1$ and the flip is said to be \defn{slope-increasing}.
\end{enumerate}
\index{flip!increasing and decreasing ---s}
We define a partial order on the set of \ktri{k}s of the \gon{n} as follows: for two \ktri{k}s~$T$~and~$T'$, we write~$T<T'$ if and only if there exists a sequence of slope-increasing flips from~$T$~to~$T'$. That this is indeed a partial order follows from the fact that each slope increasing flip increases the \defn{total slope} of a \ktri{k}, where the slope of an edge~$[u,v]$ is defined as~$u+v$ (with the sum taken in~$\N$, not in~$\Z_{n}$)  and the total slope of a \ktri{k} is the sum of slopes of its edges. The following \ktri{k} of the \gon{n} (see \fref{stars:fig:minimalktriangulation}) is a minimal element for this order, and we will see that it is even the least element. 

\begin{lemma}\label{stars:lem:minimalktriangulation}
The set $F \eqdef \ens{[i,j]}{i\in\llb 0,k-1\rrb \text{ and } j\in\llb i+k+1,i-k-1\rrb}$ is the set of \krel{k} edges of a \ktri{k}~$T_{n,k}^{\min}$ of the \gon{n} which is minimal for the partial order $<$.
\end{lemma}

\begin{proof}
Since~$F$ is the union of~$k$~disjoint crossing-free fans $\ens{[i,j]}{j\in\llb i+k+1,i-k-1\rrb}$, it is \kcross{(k+1)}-free and has $k(n-2k-1)$~edges. Thus, it is the set of \krel{k} edges of a \ktri{k} by Corollary~\ref{stars:coro:starsenumeration}. Furthermore, since the slope of any \krel{k} edge~$e$ not in~$F$ is larger than the slope of any edge in~$F$ crossing $e$, any flip of $T_{n,k}^{\min}$ is increasing, and $T_{n,k}^{\min}$ is minimal.
\end{proof}

\begin{lemma}\label{stars:lem:increase}
For any \ktri{k} of the \gon{n}~$T$ different from~$T_{n,k}^{\min}$, there exists a \krel{k} edge~$f\in T\ssm T_{n,k}^{\min}$ such that the edge added by the flip of~$f$ is in~$T_{n,k}^{\min}$.
\end{lemma}

\begin{proof}
Let~$T$ be a \ktri{k} of the \gon{n} distinct from~$T_{n,k}^{min}$. Let
$$\ell \eqdef \max\ens{k+1\le i\le n-k-1}{\{[0,i],[1,i+1],\dots,[k-1,i+k-1]\}\nsubseteq T},$$
which exists because~$T\ne T_{n,k}^{min}$.

Let~$0\le j\le k-1$ such that the edge~$[j,\ell+j]$ is not in~$T$. Let~$\{[x_1,y_1],\dots,[x_k,y_k]\}$ denote a \kcross{k} which prevents~$[j,\ell+j]$ to be in~$T$, with the convention that:
\begin{enumerate}[(i)]
\item $x_1\cl \dots\cl x_k\cl y_1\cl \dots\cl y_k$;
\item if~$j>0$, then~$j\in\;\rrb x_j,x_{j+1}\llb$~and~$\ell+j\in\;\rrb y_j,y_{j+1}\llb$;~and
\item if~$j=0$, then~$0\in\;\rrb y_k,x_1\llb$~and~$\ell\in\;\rrb x_k,y_1\llb$.
\end{enumerate}
With this convention, we are sure that~$x_k\in\llb k,\ell-1\rrb$~and~$y_k\in\llb \ell+k,n-1\rrb$. If~$y_k\in\;\rrb \ell+k,n-1\rrb$, then the set~$\{[0,\ell+1],\dots,[k-1,\ell+k],[x_k,y_k]\}$ is a \kcross{(k+1)} of~$T$. Thus~$y_k=\ell+k$, and there exists an edge~$[x_k,\ell+k]$ with $x_k\in\llb k,\ell-1\rrb$.

Now let~$m \eqdef \min\ens{k\le i\le \ell-1}{[i,\ell+k]\in T}$. Let~$f$ be the edge~$[m,\ell+k]$, let~$S$ be the \kstar{k} containing the angle~$\angle(m,\ell+k,k-1)$ and let~$R$ be the other \kstar{k} containing~$f$. Let~$s_0, \dots,s_{k-2},s_{k-1}=k-1,s_k=m,s_{k+1},\dots,s_{2k-1},s_{2k}=\ell+k$ denote the vertices of the \kstar{k}~$S$ in circle order. Then~$s_0\in\llb \ell+k+1, 0\rrb$, and the only way not to get a \kcross{(k+1)} is to have~$s_0=0$. This implies that~$s_j=j$ for all~$0\le j\le k-1$. 

Let~$e$ denote the common bisector of~$R$~and~$S$ and let~$s$ denote its vertex in~$S$. Since ${f=[m,\ell+k]=[s_k,s_{2k}]}$ is a common edge of~$R$~and~$S$, it is sure that~$s\notin\{s_k,s_{2k}\}$. Moreover, since for any~$0\le j\le k-2$ the interval~$\rrb s_j,s_{j+1}\llb$ is empty,~$s\notin\{s_{k+1},s_{k+2},\dots,s_{2k-1}\}$. Consequently,~$s\in\{s_0,\dots,s_{k-1}\}=\llb 0,k-1\rrb$ and~$e\in T_{n,k}^{\min}\ssm T$.
\end{proof}

\begin{corollary}\label{stars:coro:graphflips}
\index{diameter (of the graph of flips)|(}
\begin{enumerate}[(i)]
\item The \ktri{k}~$T_{n,k}^{\min}$ is the unique least element of the set of \ktri{k}s of the \gon{n}, partially ordered by~$<$.
\item Any \ktri{k}~$T$ of the \gon{n} can be transformed into the minimal \ktri{k}~$T_{n,k}^{\min}$ by a sequence of $|T\ssm T_{n,k}^{\min}|$ flips.
\item The graph~$G_{n,k}$ is connected, regular of degree~$k(n-2k-1)$, and its diameter~$\delta_{n,k}$ is at most~$2k(n-2k-1)$.
\end{enumerate}
\end{corollary}

\begin{proof}
Parts~(i)~and~(ii) are immediate corollaries of the previous lemma: there is a sequence of $|T\ssm T_{n,k}^{\min}|$ slope-decreasing flips from~$T$ to~$T_{n,k}^{\min}$. As far as Part~(iii) is concerned, the regularity follows from the fact that any of the~$k(n-2k-1)$ \krel{k} edges of~$T$ can be flipped. Finally, the connectedness as well as the bound on the diameter are obtained by joining any two \ktri{k}s~$T$~and~$T'$ of the \gon{n} by a path of flips passing through~$T_{n,k}^{\min}$.
\end{proof}

\begin{figure}[!h]
	\capstart
	\centerline{\includegraphics[scale=1]{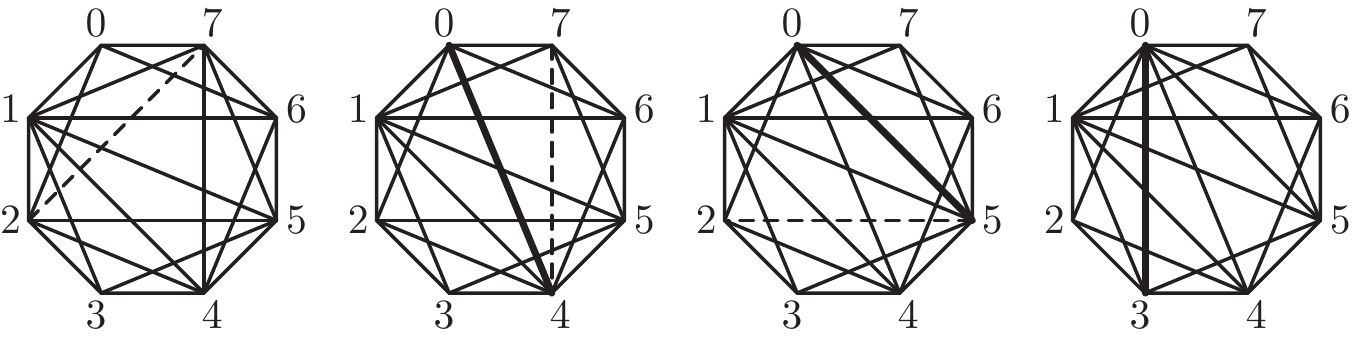}}
	\caption[A path of slope-decreasing flips from the \ktri{2} of \fref{intro:fig:2triang8points} to~$T_{8,2}^{\min}$]{A path of slope-decreasing flips from the \ktri{2} of \fref{intro:fig:2triang8points} to~$T_{8,2}^{\min}$. In each of the four pictures, the new edge is in bold and the dashed edge will be flipped.}
	\label{stars:fig:2triang8pointsflips}
\end{figure}

\begin{example}
\fref{stars:fig:2triang8pointsflips} shows a path of slope-decreasing flips from the \ktri{2} of \fref{intro:fig:2triang8points} to~$T_{8,2}^{\min}$.
\end{example}

\begin{remark}\label{stars:rem:flipandedgeenueration}
We want to observe that Corollary~\ref{stars:coro:graphflips}(iii) is the way how~\cite{n-gdfcp-00} and~\cite{dkm-lahp-02} proved that the number of edges in all \ktri{k}s of the \gon{n} is the same: indeed, any two \ktri{k}s are related by a sequence of flips which preserves the number of edges. Compare it to our direct proof in Corollary~\ref{stars:coro:starsenumeration}.
\end{remark}

Obviously, we would get symmetric results with the maximal \ktri{k}~$T_{n,k}^{\max}$ whose set of \krel{k} edges is $\ens{[i,j]}{i\in\llb n-k,n-1\rrb \text{ and } j\in\llb i+k+1,i-k-1\rrb}$. Observe also that we would obtain the same results with any rotation of the labeling of~$V_n$. Following~\cite{n-gdfcp-00}, this leads to the following improved bound for the diameter of the graph of flips:

\begin{corollary}\label{stars:coro:upperbounddiameter}
When $n>4k^2(2k+1)$, the diameter~$\delta_{n,k}$ of~$G_{n,k}$ is at most~$2k(n-4k-1)$.
\end{corollary}

\begin{proof}
Let~$\rho$ be the rotation~$t\mapsto t+1$. For any~$i\in\Z_n$, denote by~$T_i$ the \ktri{k}~$\rho^i(T_{n,k}^{\min})$. Using the same argument as before, any two \ktri{k}s~$T$~and~$T'$ of the \gon{n} are linked by a sequence of at most~$|T\ssm T_i|+|T'\ssm T_i|$ flips, for any~$i\in\Z_n$. In particular, there exists a path linking~$T$~and~$T'$ of length smaller than the average over~$i$ of~$|T\ssm T_i|+|T'\ssm T_i|$. Thus, if~$\kdeg{k}{i}{T}$ denotes the number of \krel{k} edges of~$T$ adjacent to the vertex~$i$, then the diameter of~$G_{n,k}$ is bounded by
\begin{align*}
\delta_{n,k} & \le \frac{1}{n}\sum_{i\in\Z_n} \big(|T\ssm T_i|+|T'\ssm T_i|\big) = 2k(n-2k-1) - \frac{2}{n}\sum_{i\in\Z_n}\sum_{j\in[k]} \kdeg{k}{i+j}{T} \\
& = 2k(n-2k-1) - \frac{4k^2(n-2k-1)}{n} = 2k(n-4k-1)+\frac{4k^2(2k+1)}{n}.
\end{align*}
Thus, if~$n>4k^2(2k+1)$, we obtain that $\delta_{n,k}\le 2k(n-4k-1)$.
\end{proof}

The following related result uses the same argument based on the average degree of the \ktri{k}s of the \gon{n}:

\begin{lemma}\label{stars:conj:increasingdiametermin}
The diameter~$\delta_{n,k}$ of the graph of flips satisfies $\delta_{n+1,k}\le\delta_{n,k}+4k-1$.
\end{lemma} 

\begin{proof}
Let~$T$ and~$T'$ be two \ktri{k}s of the \gon{(n+1)}, and~$i$ be a distinguished vertex of~$V_{n+1}$. Let~$T_i$ and~$T'_i$ denote the \ktri{k}s of the \gon{(n+1)} obtained from~$T$ and~$T'$ respectively by flipping all the \krel{k} edges adjacent to the vertex~$i$. Let~$\tau$ and~$\tau'$ denote the graphs obtained by forgetting from~$T_i$ and~$T'_i$ respectively the vertex~$i$ as well as the~$2k$ edges adjacent to it. Observe that Corollary~\ref{stars:coro:starsenumeration} implies that~$\tau$ and~$\tau'$ are \ktri{k}s of the \gon{n}. Consequently, there exists a path $\tau=\tau_0,\tau_1,\dots,\tau_p=\tau'$ of at most~$\delta_{n,k}$ flips which transform~$\tau$ into~$\tau'$. Remembering the vertex~$i$ and the~$2k$ edges adjacent to it, this translates into a path of at most~$\delta_{n,k}$ flips joining~$T_i$ to~$T'_i$. We deduce from this that there exists a path in~$G_{n+1,k}$ from~$T$~to~$T'$ with at most $\kdeg{k}{i}{T}+\delta_{n,k}+\kdeg{k}{i}{T'}$ flips (where~$\kdeg{k}{i}{T}$ denotes the number of \krel{k} edges of~$T$ adjacent to the vertex~$i$). Since the same argument works for any choice of the vertex $i\in\Z_{n+1}$, we can again use the average degree argument: the diameter of~$G_{n+1,k}$ is bounded by
\begin{align*}
\delta_{n+1,k} & \le \frac{1}{n+1}\sum_{i\in\Z_{n+1}}(\kdeg{k}{i}{T}+\delta_{n,k}+\kdeg{k}{i}{T'}) = \delta_{n,k}+\frac{2}{n+1}\sum_{i\in\Z_{n+1}}\kdeg{k}{i}{T} \\
& = \delta_{n,k}+\frac{4k(n+1-2k-1)}{n+1} = \delta_{n,k}+4k-\frac{4k(2k-1)}{n+1}.
\end{align*}
We finally obtain the result since~$\delta_{n+1,k}$ is an integer.
\end{proof}

Note that even if the improvement provided by Corollary~\ref{stars:coro:upperbounddiameter} is asymptotically not relevant, for the case~$k=1$, the improved bound of~$2n-10$ is actually the exact diameter of the associahedron for large values of~$n$~\cite{stt-rdthg-88}. The proof in~\cite{stt-rdthg-88} is based on the interpretation of a sequence of flips relating two triangulations~$T$~and~$T'$ of the \gon{n} as a triangulation of a \dimensional{3} simplicial polytope whose boundary is made glueing~$T$~and~$T'$ along their boundary edges. Finding such a polytope with no small triangulation is achieved using \dimensional{3} hyperbolic polytopes of large volume: indeed, since the volume of an hyperbolic tetrahedron is bounded by a constant, hyperbolic polytopes with large volume only admit triangulations with many tetrahedra.

When~$k\ne1$, the exact asymptotic value of the diameter~$\delta_{n,k}$ is unknown. As far as lower bounds are concerned, the following simple argument shows a lower bound of order~$kn$:

\begin{lemma}\label{stars:lem:diammin}
If~$n\ge 4k$, then the diameter~$\delta_{n,k}$ of~$G_{n,k}$ is at least~$k(n-2k-1)$.
\end{lemma}

\begin{proof}
The proof consists in displaying two \ktri{k}s with no \krel{k} edges in common. We call \defn{\kzz{k}}\index{zigzag@\kzz{k}|hbf} the following subset of \krel{k} edges of~$E_n$ (see \fref{stars:fig:zigzags}):
$$Z  \eqdef \ens{[q-1,-q-k]}{1\le q\le \Fracfloor{n-2k}{2}}\cup\ens{[q,-q-k]}{1\le q\le \Fracfloor{n-2k-1}{2}}.$$

\begin{figure}[!h]
	\capstart
	\centerline{\includegraphics[scale=1]{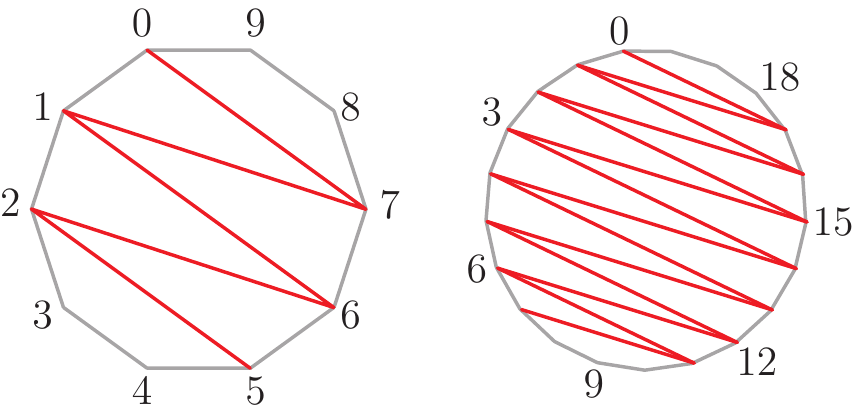}}
	\caption[Examples of \kzz{k}s of~$E_n$]{Examples of \kzz{k}s of~$E_n$ for~$(n,k)=(10,2)$~and~$(21,3)$.}
	\label{stars:fig:zigzags}
\end{figure}

Let~$\rho$ be the rotation~$t\mapsto t+1$. Observe that the~$2k$ \kzz{k}s~$Z,\rho(Z),\dots,\rho^{2k-1}(Z)$ are disjoint since~$n\ge 4k$. Moreover,
\begin{enumerate}[(i)]
\item there is no \kcross{2} in a zigzag, thus there is no \kcross{(k+1)} in the union of~$k$~of~them;
\item a \kzz{k} contains~$n-2k-1$~\krel{k} edges, so that the union of $k$ disjoint \kzz{k}s contains $k(n-2k-1)$ \krel{k} edges.
\end{enumerate}
According to Corollary~\ref{stars:coro:starsenumeration}, this proves that the union of~$k$~disjoint \kzz{k}s is the set of \krel{k} edges of a \ktri{k}. Thus, we obtain the two \ktri{k}s we were looking for with the sets of \krel{k} edges $\bigcup_{i=0}^{k-1} \rho^i(Z)$ and $\bigcup_{i=k}^{2k-1} \rho^i(Z)$.
\end{proof}

Another way to obtain Lemma~\ref{stars:lem:diammin} would be to prove the following conjecture:

\begin{conjecture}\label{stars:conj:increasingdiametermax}
The diameter $\delta_{n,k}$ of the graph of flips~$G_{n,k}$ satisfies  $\delta_{n+1,k} \ge \delta_{n,k}+k$.
\end{conjecture}

This conjecture is easy to prove when~$k=1$ (see~\cite[Chapter~1]{lrs-tri}). The proof uses an operation on triangulations generalized to \ktri{k}s in Section~\ref{stars:sec:flat-infl}, where we will discuss the reason why the natural generalization of this proof fails (see Remark~\ref{stars:rem:increasingdiameter}). Embarrassingly, we even have no proof that the diameter $\delta_{n,k}$ of the graph of flips~$G_{n,k}$ is an increasing function of~$n$ when~$k$ is fixed. This is closely related to the following conjecture (which is known to be true for triangulations~\cite[Lemma~3]{stt-rdthg-88}):

\begin{conjecture}
\begin{enumerate}[(i)]
\item The shortest path in the flip graph~$G_{n,k}$ between two \ktri{k}s~$T$ and~$T'$ of the \gon{n} never flips a common edge of~$T$ and~$T'$.
\item Let~$T$ and~$T'$ be two \ktri{k}s of the \gon{n},~$f$ be a \krel{k} edge of~$T$ and~$e$ be the unique bisector of the two \kstar{k}s of~$T$ containing~$f$. If~$e$ is  an edge of~$T'$, then there exists a shortest path in the graph of flips~$G_{n,k}$ from~$T$ to~$T'$ which first flips~$f$.
\end{enumerate}
\end{conjecture}

The lower bound of Lemma~\ref{stars:lem:diammin} can be slightly improved by a more careful choice of two disjoint \ktri{k}s. Our best lower bound is the following:

\begin{lemma}\label{stars:lem:lowerbounddiameter}
If~$m\ge 2k+1$, then both~$\delta_{2m,k}$ and~$\delta_{2m+1,k}$ are at least~$2m\left(k+\frac{1}{2}\right)-k(2k+3)$.
\end{lemma}

\begin{proof}
We first prove the result for the \gon{(2m)}. Choose an arbitrary \ktri{k}~$\tau$ of the \gon{m}, and define the following three subsets of~$E_{2m}$:
\begin{align*}
A &  \eqdef \ens{[2i+1,2i+1+j]}{i\in\Z_m \text{ and } j\in[k]},\\
B &  \eqdef \ens{[2i,2(i+j)-1]}{i\in\Z_m \text{ and } j\in[k]},\quad\text{and}\\
C &  \eqdef \ens{[2i,2j]}{i,j\in\Z_m,\;[i,j] \text{ edge of } \tau}.
\end{align*}

We consider the union~$T \eqdef A\cup B\cup C$ of these three subsets of~$E_{2m}$ and we denote by $T' \eqdef \rho(T)$~its image under the rotation~$\rho:t\mapsto t+1$ of the \gon{(2m)}. In a first step, we prove that~$T$ (and thus, $T'$ as well) is a \ktri{k} of the \gon{(2m)}, and in a second step, we prove that the number of flips needed to join $T$ to $T'$ is at least $2m\left(k+\frac{1}{2}\right)-k(2k+3)$.

To visualize an example, \fref{stars:fig:optdiameter} shows the two \ktri{2}s~$T$ and~$T'$ of the \gon{16} obtained from the \ktri{2}~$\tau$ of the \gon{8} represented in \fref{intro:fig:2triang8points}.

\begin{figure}
	\capstart
	\centerline{\includegraphics[scale=1]{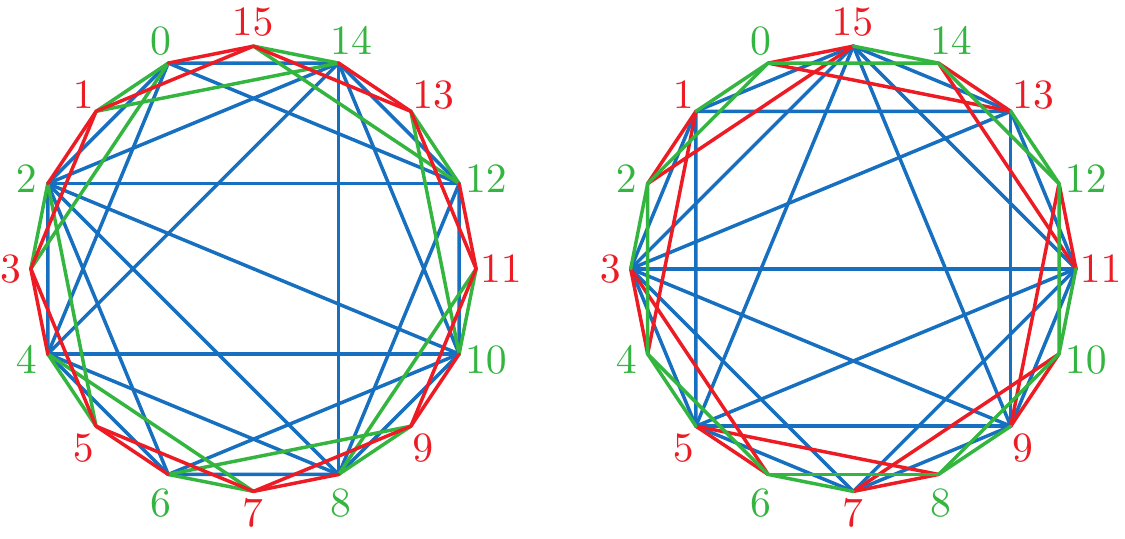}}
	\caption[Two \ktri{2}s far appart in the graph of flips]{The two \ktri{2}s~$T$ and~$T'$ of the \gon{16} obtained from the \ktri{2}~$\tau$ of the \gon{8} represented in \fref{intro:fig:2triang8points}.}
	\label{stars:fig:optdiameter}
\end{figure}

\paragraph{First step.} 
We prove first that the union~$T=A\cup B\cup C$ is disjoint and \kcross{(k+1)}-free. Since $|A|+|B|+|C|=mk+mk+k(2m-2k-1)=k(4m-2k-1)$, Corollary~\ref{stars:coro:starsenumeration} then ensures that~$T$ is a \ktri{k} of the \gon{(2m)}.

To see that $A$, $B$ and $C$ are disjoint, orient all their edges in such a way that each edge has less vertices of~$V_{2m}$ on its right than on its left. With this orientation an edge $[2i+1,2i+1+j]\in A$ has an odd tail, an edge $[2i,2(i+j)-1]\in B$ has an even tail and on odd head, while both tail and head of any edge $[2i,2j]\in C$ are even. This separates the three sets~$A$,~$B$~and~$C$.

Assume now that~$T$ contains a \kcross{(k+1)}~$F$. Since the edges of~$A$ are \kirrel{k}, $F$~is included in~$B\sqcup C$. We derive from~$F$ a subset~$F'$ of~$C$ as follows: 
\begin{enumerate}[(i)]
\item we replace each edge $[2i,2j+1]\in F\cap B$ by the edge $[2i,2j+2]\in C$;
\item we replace each edge $[2i,2j]\in F\cap C$, with $|i-j|<k$, by the edge $[2i,2j+2]$;~and
\item we keep all other edges of~$F$ (that is, the edges $[2i,2j]\in F\cap C$ with $|i-j|\ge k$).
\end{enumerate}
In other words, we shift to the next even vertex the head of all the edges of~$F$ with length strictly smaller than~$2k$. Observe that~$F'$ is indeed a subset of~$C$ since we only introduce edges of the form~$[2i,2j+2]$ where~$|i-j|< k$ which are necessarily in~$C$ (because~$\tau$, as any \ktri{k} of the \gon{m}, contains all non-\krel{k} edges). Also $F'$ has still~$k+1$ different edges since we only modified the head of some edges which keep all the tails different. Finally, we claim that~$F'$ is still a \kcross{(k+1)}. Indeed, choose any two edges in~$F$, and denote them~$e \eqdef [a,b]$~and~$f \eqdef [c,d]$ such that~$e$ (resp.~$f$) is oriented from~$a$ to~$b$ (resp.~$c$ to~$d$) and~$a\cl c\cl b$. The only possible way to uncross~$e$ and~$f$ would be that~$d=b+\varepsilon$, with~$\varepsilon\in\{1,2\}$ and that our transformation shifts the head of~$e$ by~$\varepsilon$, keeping the head of~$f$ at~$d$. But this would imply that~$e$ has length smaller than~$2k-\varepsilon$ while~$f$ has length greater or equal~$2k$, which is impossible if~$a\cl c\cl b$ and~$d=b+\varepsilon$. We obtain a contradiction with the \kcross{(k+1)}-freeness of~$\tau$; thus,~$T$ is \kcross{(k+1)}-free, and consequently is a \ktri{k}.

\paragraph{Second step.}
We shall now study the number of flips needed to transform~$T$ into~$T'=\rho(T)$. Observe first that~$T$ and~$T'$ have no \krel{k} edges in common, again because all the tails of the \krel{k} edges of~$T$ are even, while all the tails of the \krel{k} edges of~$T'$ are odd. Thus, as in the proof of the previous lemma, we need at least~$k(2m-2k-1)$ flips to destroy all the edges of~$T$. Now consider a \kstar{k}~$S$ of~$T$ with only even vertices. At the moment when we first destroy an edge of this star, we necessarily construct an edge of length at least~$2k+1$ and with one even endpoint (a vertex of~$S$). Such an edge does not exist in the \ktri{k}~$T'$ (the only edges of~$T'$ with an even endpoint have length strictly smaller than~$2k$) so that we will have to destroy it again. Since~$T$ contains $m-2k$ \kstar{k}s with all vertices even (the images of the \kstar{k}s of~$T$ under the map~$t\mapsto 2t$), the number of flips needed to join the \ktri{k}s~$T$ and~$T'$ is at least
$$k(2m-2k-1)+m-2k = 2m\left(k+\frac{1}{2}\right)-k(2k+3).$$

Finally, the proof for the \gon{(2m+1)} is obtained similarily. We use our two \ktri{k}s~$T$ and~$T'$ of the \gon{(2m)} to which we add a vertex~$2m$ (between vertices~$2m-1$ and~$0$), connected only with the vertices of~$\{2m-k,\dots,2m-1\}\cup\{0,\dots,k-1\}$. The analysis of the distance between the two resulting \ktri{k}s of the \gon{(2m+1)} is then exactly the same as in the case of even polygons.
\end{proof}

Observe that Conjecture~\ref{stars:conj:increasingdiametermax} would prove that:
$$\delta_{2m+1,k} \ge \delta_{2m,k}+k \ge \left(k+\frac{1}{2}\right)2m-k(2k+3)+k,$$
and consequently that, for any~$n\ge 2(2k+1)$,
$$\delta_{n,k} \ge \Floor{\left(k+\frac{1}{2}\right)n-k(2k+3)}.$$

\mvs
To summarize our results on the diameter of the graph of flips~$G_{n,k}$, we have proved that when~$n>4k^2(2k+1)$,
$$2\Fracfloor{n}{2}\left(k+\frac{1}{2}\right)-k(2k+3) \le \delta_{n,k} \le 2k(n-4k-1).$$
Furthermore, we know that the difference~$\delta_{n+1,k}-\delta_{n,k}$ between two consecutive values of this diameter is at most~$4k-1$, and we conjecture that it is at least~$k$ (for~$n$ sufficiently large).

\mvs
For completeness, we collect the couples~$(n,k)$ for which the exact value of~$\delta_{n,k}$ is known:
\begin{enumerate}
\item When~$k=1$, the exact value of the diameter of the associahedron is known~\cite{stt-rdthg-88} for the following little values of~$n$:
\svs
\begin{center}
\begin{tabular}{c|cccccccccccccccc}
$n$ & $3$ & $4$ & $5$ & $6$ & $7$ & $8$ & $9$ & $10$ & $11$ & $12$ & $13$ & $14$ & $15$ & $16$ & $17$ & $18$ \\
\hline
$\delta_{n,1}$ & $0$ & $1$ & $2$ & $4$ & $5$ & $7$ & $9$ & $11$ & $12$ & $15$ & $16$ & $18$ & $20$ & $22$ & $24$ & $26$
\end{tabular}
\end{center}
\svs
Furthermore, there exists an integer~$N$ such that for any~$n$ larger than~$N$, the diameter of the associahedron is exactly~$2n-10$~\cite{stt-rdthg-88}. The minimal value of~$N$ remains unknown (according to the previous table, it is conjectured to be~$13$).
\item When~$k=2$ and~$k=3$, a computer search provides the following tables:
\svs
\begin{center}
\begin{tabular}{c|ccccccc}
$n$ & $5$ & $6$ & $7$ & $8$ & $9$ & $10$ & $11$ \\
\hline
$\delta_{n,2}$ & $0$ & $1$ & $3$ & $6$ & $8$ & $11$ & $14$
\end{tabular}
\hspace{1.5cm}
\begin{tabular}{c|ccccc}
$n$ & $7$ & $8$ & $9$ & $10$ & $11$ \\
\hline
$\delta_{n,3}$ & $0$ & $1$ & $3$ & $6$ & $10$
\end{tabular}
\end{center}
\svs
\item When~$n=2k+1$ (see Example~\ref{stars:exm:n=2k+1}), there exist only one \ktri{k} and~$\delta_{2k+1,k}=0$.
\item When~$n=2k+2$ (see Example~\ref{stars:exm:n=2k+2}), the graph of flips~$G_{2k+2,k}$ is complete and~$\delta_{2k+2,k}=1$.
\item When~$n=2k+3$ (see Example~\ref{stars:exm:n=2k+3}), the graph of flips~$G_{2k+3,k}$ is the ridge graph of the \dimensional{2k} cyclic polytope with~$2k+3$ vertices (see Lemma~\ref{ft:lem:2k+3cyclicpolytope}), and its diameter is~$\delta_{2k+3,k}=3$ as soon as~$k\ge 2$ (we refer to~\cite{m-drgcp-09} for the diameter of the ridge graph of the cyclic polytope).
\end{enumerate}
\index{diameter (of the graph of flips)|)}


\section{\kear{k}s and \kcolorable{k} \ktri{k}s}\label{stars:sec:ears}

Throughout this section, we assume that~$n\ge 2k+3$. Generalizing similar notions for triangulations, we define:

\begin{definition}\label{stars:def:ear}
\index{ear@\kear{k}|hbf}
\index{star@\kstar{k}!external ---}
\index{star@\kstar{k}!internal ---|hbf}
A \defn{\kear{k}} of a \ktri{k} is any edge of length~$k+1$. A \kstar{k} of a \ktri{k} is \defn{external} if it contains at least one \kbound{k} edge, and \defn{internal} otherwise.
\end{definition}

It is well known and easy to prove that the number of ears in any triangulation equals its number of internal triangles plus~$2$. In this section, we prove that the number of \kear{k}s in any \ktri{k} equals its number of internal \kstar{k}s plus~$2k$. Then we characterize the triangulations that have no internal \kstar{k} in terms of colorability of their intersection complex.


\subsection{The maximal number of \kear{k}s}\label{stars:subsec:ears:2k}

Let~$S$ be a \kstar{k} and~$[u,v]$ be an edge of~$S$ such that~$S$ lies on the positive side of the oriented edge from~$u$ to~$v$. We say that~$[u,v]$ is a \defn{positive ear} of~$S$ if~$|\llb v,u\llb|=k+1$. Said differently,~$[u,v]$ is an ear, and~$S$ is the unique (by Corollary~\ref{stars:coro:incidences}) \kstar{k} on the ``outer'' side of it.

\begin{lemma}\label{stars:lem:starwithboundary}
Let~$b\ge 1$ be the number of \kbound{k} edges of an external \kstar{k}~$S$. Then,~$S$ has exactly~$b-1$ positive ears. Moreover, \kbound{k} edges and positive ears of~$S$ form an alternating path in the \kstar{k}.
\end{lemma}

\begin{proof}
Observe first that if~$[x,x+k]$ is a \kbound{k} edge of~$S$, then all~$k+1$ vertices of~$\llb x,x+k\rrb$ are vertices of~$S$. Since~$S$ has~$2k+1$ vertices, this implies that if~$[x,x+k]$ and~$[y,y+k]$ are two \kbound{k} edges of~$S$, then~$\llb x,x+k\rrb$ and~$\llb y,y+k\rrb$ intersect; that is, for example,~$y\in \llb x,x+k\rrb$. But then, all edges~$[i,i+k]$ with~$x\cle i\cle y$ are \kbound{k} edges of~$S$, and all edges $[i,i+k+1]$ with $x\cle i\cl y$ are positive ears of~$S$. Thus, \kbound{k} edges and positive ears of~$S$ form an alternating path in the \kstar{k}, beginning and ending by a \kbound{k} edge. In particular,~$S$ has exactly~$b-1$ positive ears.
\end{proof}

\begin{corollary}
\begin{enumerate}[(i)]
\item A \kstar{k} cannot have more than~$k+1$ \kbound{k} edges.
\item The \kbound{k} edges of a \kstar{k} are always consecutive around the circle.\qed
\end{enumerate}
\end{corollary}

Observe that any \kstar{k} with all vertices consecutive has~$k+1$ \kbound{k} edges (\ie attains the bound~(i) in the previous corollary). Such a \kstar{k} is an alternating path of~$k+1$ \kbound{k} edges and~$k$ \kear{k}s.

Lemma~\ref{stars:lem:starwithboundary} has also the following consequence:

\begin{corollary}\label{stars:coro:earsenumeration}
The number of \kear{k}s in a \ktri{k}~$T$ equals the number of internal \kstar{k}s of~$T$ plus~$2k$. In particular,~$T$~contains at least~$2k$ \kear{k}s.
\end{corollary}

\begin{proof}
For any~$0\le i\le 2k+1$, let~$\mu_i$ denote the number of \kstar{k}s of~$T$ with exactly~$i$ \kbound{k} edges. Let~$\nu$ denote the number of \kear{k}s of~$T$. Then
$$\sum_{i=0}^{2k+1} \mu_i=n-2k \quad \text{and}\quad \sum_{i=1}^{2k+1} i\mu_i=n.$$
Since any \kear{k} is a positive ear of exactly one \kstar{k}, Lemma~\ref{stars:lem:starwithboundary} ensures moreover that
$$\sum_{i=1}^{2k+1} (i-1)\mu_i=\nu.$$
Thus, we obtain~$\nu=n-(n-2k)+\mu_0=2k+\mu_0$.
\end{proof}

\begin{example}
The \ktri{2} in \fref{intro:fig:2triang8points} has five \kear{2}s and one internal \kstar{2}.
\end{example}


\subsection{\kcolorable{k} \ktri{k}s}\label{stars:subsec:ears:colorable}

We are now interested in a characterization of the \ktri{k}s that have exactly~$2k$ \kear{k}s, or equivalently that have no internal \kstar{k}. We need two additional definitions.

\begin{definition}\label{stars:def:kcolorable}
\index{triangulation@\ktri{k}!colorable@\kcolorable{k} ---|hbf}
We say that a \ktri{k} is \defn{\kcolorable{k}} if it is possible to color its \krel{k} edges with~$k$ colors such that there is no monochromatic \kcross{2}. Observe that, if this happens, then every \kcross{k} uses an edge of each color.
\end{definition}

\begin{definition}\label{stars:def:kaccordion}
\index{accordion@\kaccordion{k}|hbf}
A \defn{\kaccordion{k}} of~$E_n$ is a set~$Z \eqdef \ens{[a_i,b_i]}{i\in[n-2k-1]}$ of~$n-2k-1$ edges such that~$b_1=a_1+k+1$ and for any~$2\le i\le n-2k-1$, the edge~$[a_i,b_i]$ is either~$[a_{i-1},b_{i-1}+1]$ or~$[a_{i-1}-1,b_{i-1}]$.
\end{definition}

\begin{figure}[b]
	\capstart
	\centerline{\includegraphics[scale=1]{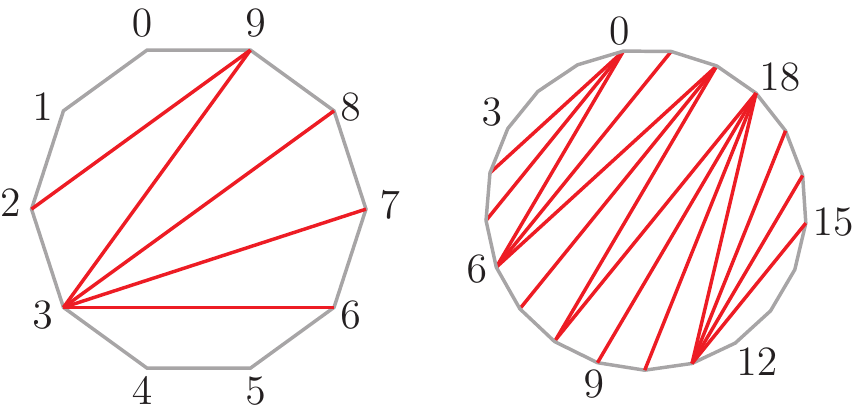}}
	\caption[Examples of \kaccordion{k}s of~$E_n$]{Examples of \kaccordion{k}s of~$E_n$ for~$(n,k)=(10,2)$ and~$(21,3)$.}
	\label{stars:fig:accordion}
\end{figure}

Note that, in this definition,~$b_i=a_i+k+i$ for any~$i\in[n-2k-1]$. When~$k\ge 2$, observe that the definition is equivalent to being a set of~$n-2k-1$ \krel{k} edges of~$E_n$ without \kcross{2}.

\begin{lemma}\label{stars:lem:unionzigzags}
The union of~$k$ disjoint \kaccordion{k}s of~$E_n$ is the set of \krel{k} edges of a \ktri{k} of the \gon{n}.
\end{lemma}

\enlargethispage{.4cm}
\begin{proof}
Observe that:
\begin{enumerate}[(i)]
\item A \kaccordion{k} has no \kcross{2}, so the union of~$k$ of them has no \kcross{(k+1)}.
\item A \kaccordion{k} contains~$n-2k-1$ \krel{k} edges; thus, the union of~$k$ disjoint \kaccordion{k}s contains~$k(n-2k-1)$ \krel{k} edges.
\end{enumerate}
This lemma thus follows from Corollary~\ref{stars:coro:starsenumeration}.
\end{proof}

We have already met some particular \kaccordion{k}s both when we constructed the triangulation $T_{n,k}^{\min}$ and in the proof of Lemma~\ref{stars:lem:diammin}. The two types of \kaccordion{k}s in these examples (the ``fan'' and the ``zigzag'') are somehow the two extremal examples of them: one has only alternating angles when the other has no of them.

The following theorem relates the two Definitions~\ref{stars:def:kcolorable} and~\ref{stars:def:kaccordion} and characterizes \ktri{k}s with exactly~$2k$ \kear{k}s:

\begin{theorem}\label{stars:theo:kcolorable}
Let~$T$ be a \ktri{k}, with~$k>1$. The following properties are equivalent:
\begin{enumerate}[(i)]
\item $T$~is \kcolorable{k}.
\item There exists a \kcoloring{k} of the \krel{k} edges of~$T$ such that no \kstar{k} of~$T$ contains three edges of the same color.
\item $T$~has no internal \kstar{k}.
\item $T$~contains exactly~$2k$ \kear{k}s.
\item $T$~contains exactly~$2k$ edges of each length $k+1,\dots, \Fracfloor{n-1}{2}$ (and~$k$ of length~$\frac{n}{2}$ if~$n$ is even).
\item The set of \krel{k} edges of $T$ is the union of $k$ disjoint \kaccordion{k}s.
\end{enumerate}
\end{theorem}

Note that for~$k=1$, only~(ii),~(iii),~(iv),~(v) and~(vi) are equivalent, while~(i) always holds.

\begin{proof}[Proof of Theorem~\ref{stars:theo:kcolorable}]
When~$k>1$, any three edges of a \kstar{k} form at least one \kcross{2}. Thus, any \kcoloring{k} of~$T$ without monochromatic \kcross{2} is such that no \kstar{k} of~$T$ contains three edges of the same color, and~(i)$\Rightarrow$(ii).

Let~$S$ be a \kstar{k} of a \ktri{k} whose \krel{k} edges are colored with~$k$ colors. If all edges of~$S$ are \krel{k}, then, by the pigeon-hole principle, there is a color that colors three edges of~$S$. Thus~(ii)$\Rightarrow$(iii).

Corollary~\ref{stars:coro:earsenumeration} ensures that~(iii)$\Leftrightarrow$(iv). Thus, since~(vi)$\Rightarrow$(v)$\Rightarrow$(iv)~and~(vi)$\Rightarrow$(i) are trivial, it only remains to prove that~(iv)$\Rightarrow$(vi). 

For this, we give an algorithm that finds the~$k$ disjoint \kaccordion{k}s in a \ktri{k}~$T$ with~$2k$ \kear{k}s. Recall that if~$S$ is an external \kstar{k} with~$b$ \kbound{k} edges, then~$S$ has~$b-1$ positive ears, and moreover \kbound{k} edges and positive ears of~$S$ form an alternating path in the \kstar{k}. Thus, the edges of~$S$ which are neither \kbound{k} nor positive ears of~$S$ form a path of even length. This defines naturally a ``pairing'' of them: we say that the first and the second (resp.~the third and the fourth, \etc) edges of this path form a \defn{pair} of edges in~$S$. Observe that such a pair of edges forms an angle containing no vertex.

Consider now a \kear{k}~$e_1$ of~$T$. Let~$S_0$ be the \kstar{k} of~$T$ for which~$e_1$ is a positive ear. Let~$S_1$ be the other \kstar{k} of~$T$ containing~$e_1$. Let~$e_2$ be the pair of~$e_1$ in~$S_1$. Let~$S_2$ be the other \kstar{k} of~$T$ containing~$e_2$ and let~$e_3$ be the pair of~$e_2$ in~$S_2$. Let continue so until we reach a \kear{k}. It is obvious that we get a \kaccordion{k} of~$E_n$. To get another \kaccordion{k}, we do the same with a \kear{k} which is neither~$e_1$ nor~$e_{n-2k-1}$. 

To prove the correctness of this algorithm, we only have to prove that the \kaccordion{k}s we construct are disjoint. Suppose that two of them~$\{e_1,\dots,e_{n-2k-1}\}$ and~$\{f_1,\dots,f_{n-2k-1}\}$ intersect. Let~$i,j$ be such that~$e_i=f_j$. Let~$S$ be a \kstar{k} containing~$e_i$. Then by construction, either~$e_{i+1}$ or~$e_{i-1}$ (resp.~either~$f_{j+1}$ or~$f_{j-1}$) is the pair of~$e_i$ in~$S$. Thus, either~$e_{i+1}=f_{j+1}$, or~$e_{i+1}=f_{j-1}$. By propagation, we get that~$\{e_1,\dots,e_{n-2k-1}\}=\{f_1,\dots,f_{n-2k-1}\}$.
\end{proof}

Observe that the above proof of~(iv)$\Rightarrow$(vi) also gives uniqueness (up to permutation of colors) of the \kcoloring{k} (decomposition into accordions) of a \kcolorable{k} \ktri{k}~$T$. Indeed, any \kcoloring{k} has to respect the pairing of edges in all \kstar{k}s of~$T$.

Let us also remark that Part~(v) implies that every \kcolorable{k} \ktri{k} contains exactly~$k(n-2p-1)$ \krel{p} edges, for any~$k\le p\le \Fracfloor{n-1}{2}$. It is proved by a flip method in~\cite{n-gdfcp-00} that any \ktri{k} of the \gon{n} contains at most this same number of \krel{p} edges. We will prove this result in the next section, as an application of the flattening of \kstar{k}s in \ktri{k}s.

Let us conclude this section with the following remark. An easy ``intuitive model'' for a \ktri{k} is just a superposition of~$k$ triangulations. Even if this model is sometimes useful, the results in this section say that it is also misleading:
\begin{itemize}
\item Theorem~\ref{stars:theo:kcolorable} says that the structure of \ktri{k}s obtained in this way is really particular: it is uniquely \kcolorable{k}, the number of edges of each length is~$2k$, all \kstar{k}s are external, \etc
\item The number of \kaccordion{k}s of~$E_n$ containing a given \kear{k} of~$E_n$ is~$2^{n-2k-2}$. In particular, the number of \kcolorable{k} \ktri{k}s of~$E_n$ is at most~${n \choose k}2^{k(n-2k-2)}\le 2^{(k+1)n}$. This is much smaller than the total number of \ktri{k}s, which for constant~$k$ equals~$4^{nk}$ modulo a polynomial factor in~$n$ (see Remark~\ref{ft:rem:enumeration}).
\item Let~$T$ be any triangulation with only two ears. Then it is easy to see that there exists a \kaccordion{2}~$Z$ disjoint from it, so that~$Z$ completes~$T$ to give the set of \krel{2} edges of a \kcolorable{2} \ktri{2}. Surprisingly, this property fails when~$k\ge3$: even if a \ktri{k}~$T$ is \kcolorable{k}, it is not always possible to find a \kaccordion{(k+1)} disjoint from~$T$ (see \fref{stars:fig:ctrexm}).

\begin{figure}[b]
	\capstart
	\centerline{\includegraphics[scale=1]{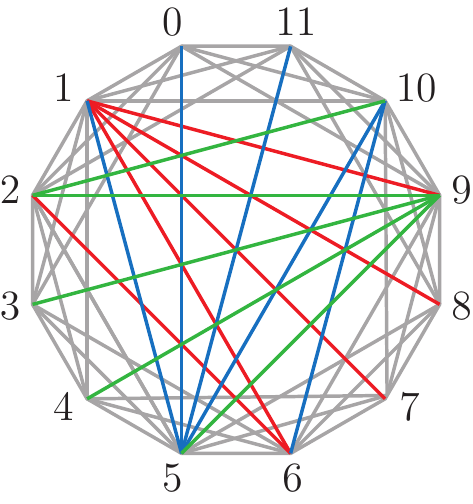}}
	\caption[A \kcolorable{3} \ktri{3} of the \gon{12} with no \kaccordion{4} disjoint from it]{A \kcolorable{3} \ktri{3} of the \gon{12} with no \kaccordion{4} disjoint from it.}
	\label{stars:fig:ctrexm}
\end{figure}

\item There even exist \ktri{k}s of the \gon{n} which do not contain a single triangulation of the \gon{n}. See Example~\ref{mpt:exm:irred} and \fref{mpt:fig:15gonctrexm}.
\end{itemize}


\section{Flattening a \kstar{k}, inflating a \kcross{k}}\label{stars:sec:flat-infl}

The goal of this section is to describe in terms of \kstar{k}s an operation that connects \ktri{k}s of~$n$ and of~$n+1$ vertices. This operation, already present in~\cite{n-gdfcp-00,j-gt-03}, and when $k=2$ in~\cite{e-btdp-07,n-abtdp-09}, is useful for recursive arguments and it was a step in all previous proofs of the flipability of \krel{k} edges (Lemma~\ref{stars:lem:flip}). We only present this operation here to emphasize that none of our proofs so far make use of induction. In the end of this section, we discuss the properties of this operation with respect to flips (namely, how does this property commutes with flips), and we present some applications.

Before giving the precise definitions, let us present a simplified picture of this operation:

\begin{figure}[!h]
	\capstart
	\centerline{\includegraphics[scale=1]{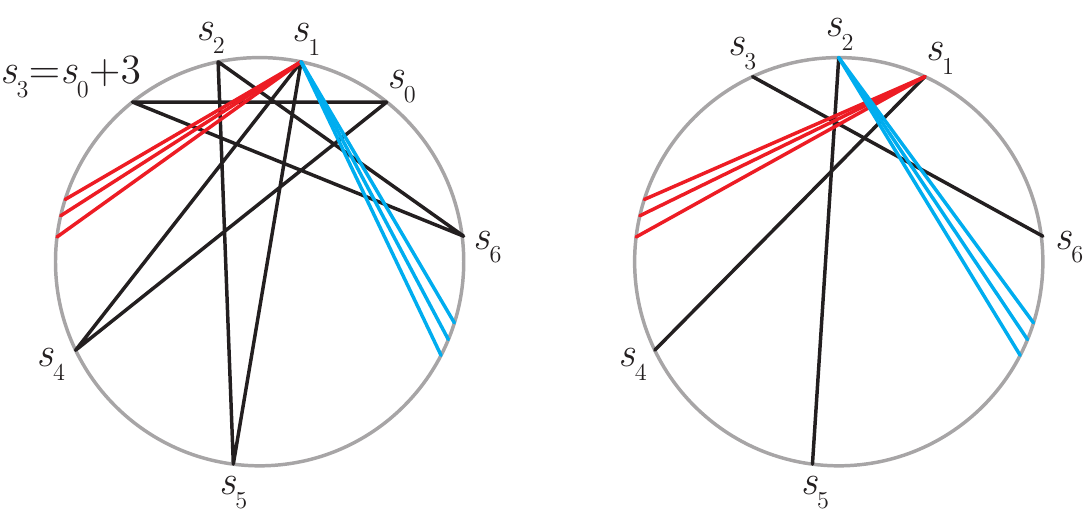}}
	\caption[Flattening a \kstar{k}, inflating a \kcross{k}]{Flattening a \kbound{3} edge, inflating a \kcross{3}.}
	\label{stars:fig:flat-infl}
\end{figure}


\subsection{Flattening a \kstar{k}}\label{stars:subsec:flat-infl:flattening}

\begin{definition}\label{stars:def:flattening}
\index{flattening (an external \kstar{k})|hbf}
Let~$T$ be a \ktri{k} of the \gon{(n+1)}, and~$b \eqdef [s_0,s_0+k]$ be a \kbound{k} edge of~$T$. Let~$s_0,s_1=s_0+1,\dots,s_k=s_0+k,s_{k+1},\dots,s_{2k}$ be the vertices of the unique \kstar{k}~$S$ of~$T$ containing~$b$, in their circle order.

We call \defn{flattening} of~$b$ in~$T$ the set of edges~$\flattening{T}{b}$ whose underlying set of points is~${V_{n+1}\ssm\{s_0\}}$ and which is constructed from~$T$ as follows (see \fref{stars:fig:flat-infl}):
\begin{enumerate}[(i)]
\item for any edge of~$T$ whose vertices are not in~$\{s_0,\dots,s_k\}$ just copy the edge;
\item forget all the edges~$[s_0,s_i]$, for~$i\in[k]$;
\item replace any edge of the form~$[s_i,t]$ with~$0\le i\le k$, and~$s_k\cl t \cle s_{k+i}$ (resp.~${s_{k+i+1}\cle t\cl s_0}$) by the edge~$[s_i,t]$ (resp.~$[s_{i+1},t]$).
\end{enumerate}
\end{definition}

\begin{remark}
We sometimes call this operation ``flattening of the \kstar{k}~$S$''. This is a more graphical description of the operation, but it is also a slight abuse of language: if~$S$ has more than one \kbound{k} edge, the result of the flattening  depends on the \kbound{k} edge we flatten, and not only on its adjacent \kstar{k}.
\end{remark}

We want to prove that~$\flattening{T}{b}$ is a \ktri{k} of the \gon{n}. Observe first that:
\begin{enumerate}[(1)]
\item If~$e$ is a \krel{k} edge of $\flattening{T}{b}$, then:
\begin{enumerate}[(i)]
\item either~$e$ is of the form~$[s_i,s_{k+i}]$, for some $i\in[k]$, and then it arises as the gluing of two edges~$e'=[s_{i-1},s_{k+i}]$~and~$e''=[s_i,s_{k+i}]$ of the initial \ktri{k}~$T$;
\item otherwise,~$e$~arises from a unique \krel{k} edge~$e'$ of~$T$.
\end{enumerate}

\item If~$e$ and~$f$ are two \krel{k} edges of~$\flattening{T}{b}$ arising from~$e'$ and~$f'$ respectively, then:
\begin{enumerate}[(i)]
\item either~$e$ and~$f$ do not cross, and then~$e'$ and~$f'$ do not cross; 
\item or~$e$~and~$f$ do cross, and then~$e'$~and~$f'$~do cross, unless the following happens: there exists~$i\in[k]$, and~$u,v$ two vertices such that~$s_k\cl u\cle s_{i+k}\cl s_{i+k+1}\cle v\cl s_0$ and~$e=[s_i,u]$,~$f=[s_{i+1},v]$,~$e'=[s_i,u]$ and~$f'=[s_i,v]$. Such a configuration is said to be a \defn{hidden configuration} (see \fref{stars:fig:flat-infl}).
\end{enumerate}
\end{enumerate}

It is easy to derive from this that~$|\flattening{T}{b}|=|T|-2k=k(2n-2k-1)$ and that any subset of~$E_n$ that properly contains~$\flattening{T}{b}$ contains a \kcross{(k+1)}. However, this is not sufficient to conclude that $\flattening{T}{b}$ is a \ktri{k}. 
Thus we have to prove that~$\flattening{T}{b}$ is \kcross{(k+1)}-free. Note that this provides a third proof (which this time is recursive on~$n$) of the formula for the number of edges of a \ktri{k} of the \gon{n}.


\begin{lemma}\label{stars:lem:flattening}
The set~$\flattening{T}{b}$ is \kcross{(k+1)}-free. Hence, it is a \ktri{k} of the \gon{n}.
\end{lemma}

\begin{proof}
Suppose that~$\flattening{T}{b}$ contains a \kcross{(k+1)}~$E$, and let~$e_0 \eqdef [x_0,y_0],\dots,e_k \eqdef [x_k,y_k]$ denote the edges of~$E$ ordered such that~$x_0\cl x_1\cl\cdots\cl x_k\cl y_0\cl y_1\cl \cdots\cl y_k$. Let~$e'_0 \eqdef [x'_0,y'_0],\dots,e'_k \eqdef [x'_k,y'_k]$ be~$k+1$ edges of~$T$ that give (when we flatten~$b$) the edges $e_0,\dots,e_k$ respectively.

It is clear that if there exists no~$0\le i\le k$ such that the four edges~$(e_i,e_{i+1},e'_i,e'_{i+1})$ form a hidden configuration, then the edges~$e'_0,\dots,e'_k$ form a \kcross{(k+1)} of~$T$, which is impossible. Thus we can suppose that the number of hidden configurations in~$\ens{(e_i,e_{i+1},e'_i,e'_{i+1})}{0\le i\le k}$ is at least~$1$. We can also assume that this number is minimal, that is that we cannot find a \kcross{(k+1)}~$F$ of~$\flattening{T}{b}$ arising from a set~$F'$ of edges of~$T$ such that there are strictly less hidden configurations in~$\ens{(f_i,f_{i+1},f'_i,f'_{i+1})}{0\le i\le k}$ than in~$\ens{(e_i,e_{i+1},e'_i,e'_{i+1})}{0\le i\le k}$. Here, we raise an absurdity by finding such sets~$F$ and~$F'$.

Let~$0\le i\le k$ be such that~$(e_i,e_{i+1},e'_i,e'_{i+1})$ is a hidden configuration. We can assume that~$x_i=s_i$ and~$x_{i+1}=s_{i+1}$ (if this is not the case, we renumber the edges of~$E$ such that this be true). Thus we know that~$y_i\cle s_{i+k}\cl s_{i+k+1}\cle y_{i+1}$. Let 
\begin{align*}
p &  \eqdef  \min\ens{0\le j< i}{\text{for any } j< \ell\le i, \; x_\ell=s_\ell \text{ and }  y_\ell\cle s_{\ell+k}}, \quad \text{and} \\
q &  \eqdef  \max\ens{i+1<j\le k+1}{\text{for any } i+1\le \ell< j, \; x_\ell=s_\ell\;\text{and}\;  s_{\ell+k}\cle y_\ell}.
\end{align*}
Let~$F$ be the set of~$k+1$ edges of~$\flattening{T}{b}$ deduced from~$E$ as follows:
\begin{itemize}
\item for all~$i$ with~$p<i<q$, let~$f_i$ be~$[s_i,s_{i+k}]$,
\item for all~$i$ with~$0\le i\le p$ and~$q\le i\le k$, let~$f_i$ be~$e_i$.
\end{itemize}
Let~$F'$ be the set of~$k+1$ edges of~$T$ constructed as follows:
\begin{itemize}
\item for all~$i$ with~$p<i<q$, let~$f'_i$ be~$[s_i,s_{i+k}]$,
\item for all~$i$ with~$0\le i\le p$ or~$q\le i\le k$, let~$f'_i$ be~$e'_i$.
\end{itemize}
It is quite clear that~$F$ is a \kcross{(k+1)} of~$\flattening{T}{b}$ arising from~$F'$. We just have to verify that the number of hidden configurations in the set~$\ens{(f_i,f_{i+1},f'_i,f'_{i+1})}{0\le i\le k}$ is less than in the set~$\ens{(e_i,e_{i+1},e'_i,e'_{i+1})}{0\le i\le k}$. But:
\begin{enumerate}[(i)]
\item the number of hidden configurations in~$\ens{(f_i,f_{i+1},f'_i,f'_{i+1})}{i<p \text{ or } q\le i}$ is exactly the same as in~$\ens{(e_i,e_{i+1},e'_i,e'_{i+1})}{i<p \text{ or } q\le i}$;
\item there is no hidden configuration in~$\ens{(f_i,f_{i+1},f'_i,f'_{i+1})}{p<i<q-1}$, whereas there is one in~$\ens{(e_i,e_{i+1},e'_i,e'_{i+1})}{p<i<q-1}$;~and
\item the edges~$(f_p,f_{p+1},f'_p,f'_{p+1})$ (resp.~$(f_{q-1},f_q,f'_{q-1},f'_q)$) do not form a hidden configuration.
\end{enumerate}
\end{proof}

\begin{figure}[!h]
	\capstart
	\centerline{\includegraphics[scale=1]{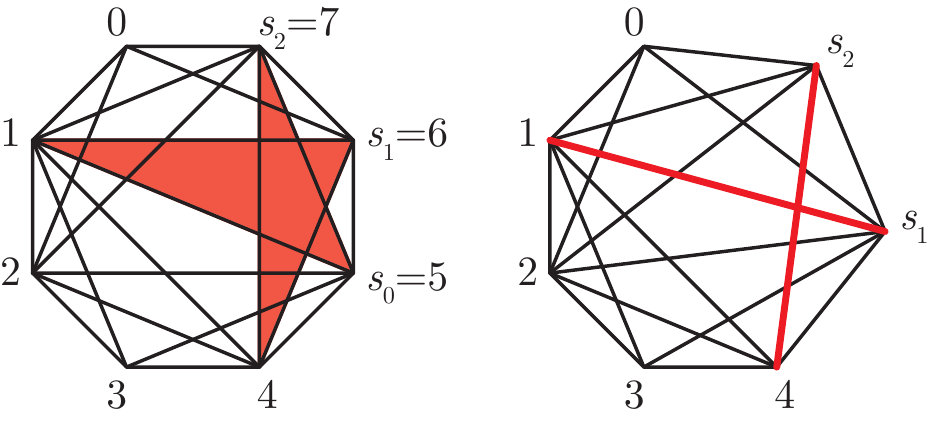}}
	\caption[An example of flattening]{Flattening the \krel{2} edge~$[5,7]$ in the \ktri{2} of \fref{intro:fig:2triang8points}.}
	\label{stars:fig:2triang8pointsflatten}
\end{figure}

\begin{example}
\fref{stars:fig:2triang8pointsflatten} shows the flattening of the \krel{2} edge~$[5,7]$ in the \ktri{2} of \fref{intro:fig:2triang8points}. The flattened \kstar{2} is colored, as well as the resulting \kcross{2}.
\end{example}

\subsection{Inflating a \kcross{k}}\label{stars:subsec:flat-infl:inflating}

\begin{definition}\label{stars:def:inflating}
\index{crossing@\kcross{k}!external ---}
\index{inflating (an external \kcross{k})|hbf}
Let~$T$ be a \ktri{k} of the \gon{n} and~$X$ be a \kcross{k}. Denote by~$[s_1,s_{1+k}],\dots,[s_k,s_{2k}]$ its edges, with~$s_1\cl s_2 \cl \dots \cl s_{2k}$. Assume further that~$s_1,\dots,s_k$ are consecutive in the cyclic order of the \gon{n} (that is~$s_k=s_1+k-1$). We call~$X$ an \defn{external} \kcross{k}.
Let~$s_0$ be a new vertex on the circle, between~$s_1-1$ and~$s_1$.

We call \defn{inflating} of~$X$ at~$\llb s_1,s_k\rrb$ in~$T$ the set of edges $\inflating{T}{X}$ whose underlying set of points is the set~$V_n\cup\{s_0\}$ and which is constructed from~$T$ as follows (see \fref{stars:fig:flat-infl}):
\begin{enumerate}[(i)]
\item for any edge of~$T$ whose vertices are not in~$\{s_1,\dots,s_k\}$ just copy the edge;
\item add all the edges~$[s_0,s_i]$, for~$i\in[k]$;~and
\item replace any edge of the form~$[s_i,t]$ with~$i\in[k]$, and~$s_k\cl t \cle s_{k+i}$ (resp.~$s_{k+i}\cle t\cl s_1$) by the edge~$[s_i,t]$ (resp.~$[s_{i-1},t]$).
\end{enumerate}
\end{definition}

\begin{remark}
We are abusing notation: if~$X$ contains more than~$k$ consecutive vertices, then the result of the inflating of~$X$ depends on the~$k$~consecutive points we choose. It would also be an abuse of language to say that we inflate the cyclic interval~$\llb s_1,s_k\rrb$ since several \kcross{k}s may be adjacent to these~$k$ points. For example, when~$k=1$, we have to specify both an edge and a vertex of this edge to define an inflating.
\end{remark}

Observe that:
\begin{enumerate}[(1)]
\item Any \krel{k} edge~$e$ of~$\inflating{T}{X}$ arises from a unique edge~$e'$ of~$T$.

\item If~$e$~and~$f$ are two \krel{k} edges of~$\inflating{T}{X}$ that cross, and~$e'$~and~$f'$ are the two edges of~$T$ that give~$e$~and~$f$ respectively, then~$e'$~and~$f'$ cross as well.
\end{enumerate}

Point~(2) ensures that~$\inflating{T}{X}$ is \kcross{(k+1)}-free. Moreover, it is easy to see that its cardinality is $|\inflating{T}{X}|=|T|+2k=k(2(n+1)-2k-1)$. Thus, by Corollary~\ref{stars:coro:starsenumeration}, we get:

\begin{lemma}\label{stars:lem:inflating}
$\inflating{T}{X}$ is a \ktri{k} of the \gon{(n+1)}.\qed
\end{lemma}

\begin{figure}[!h]
	\capstart
	\centerline{\includegraphics[scale=1]{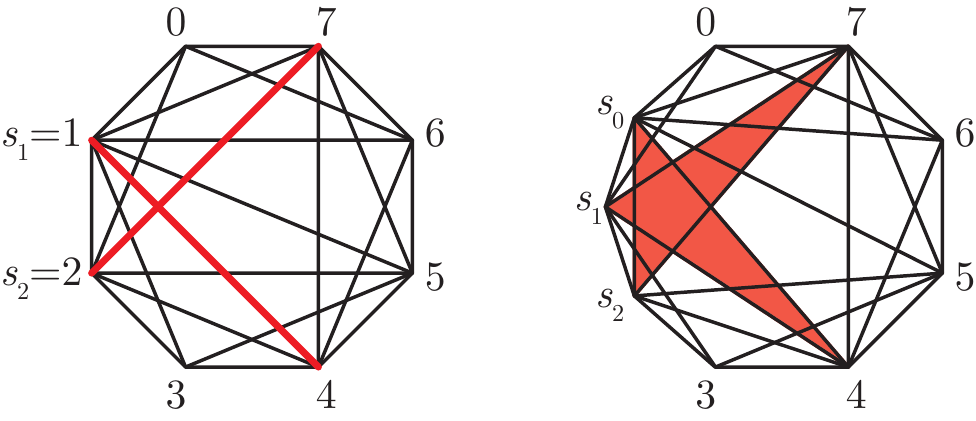}}
	\caption[An example of inflating]{Inflating the \kcross{2}~$\{[1,4],[2,7]\}$ at~$\llb 1,2\rrb$ in the \ktri{2} of \fref{intro:fig:2triang8points}.}
	\label{stars:fig:2triang8pointsinflat}
\end{figure}

\begin{example}
\fref{stars:fig:2triang8pointsinflat} shows the inflating of the \kcross{2}~$\{[1,4],[2,7]\}$ at~$\llb 1,2\rrb$ in the \ktri{2} of \fref{intro:fig:2triang8points}. The inflated \kcross{2} is colored, as well as the new \kstar{2}.
\end{example}

\subsection{Properties of the flattening \& inflating operations}\label{stars:subsec:flat-infl:properties}

First of all, the following theorem is an immediate consequence of the definitions:

\begin{theorem}\label{stars:theo:flat-infl}
Flattening and inflating are inverse operations. More precisely:
\begin{enumerate}[(i)]
\item if~$b$ is a \kbound{k} edge of a \ktri{k}~$T$, and~$X$ denotes the \kcross{k} of~$\flattening{T}{b}$ consisting of edges that arise by gluing two edges of~$T$, then $\inflating{(\flattening{T}{b})}{X}=T$;~and
\item if~$X$ is a \kcross{k} with~$k$ consecutive vertices~$s_1,\dots,s_k$, and~$b$ denotes the edge~$[s_0,s_k]$ of~$\inflating{T}{X}$, then~$\flattening{(\inflating{T}{X})}{b}=T$.\qed
\end{enumerate}
\end{theorem}

Another interesting property is the behavior of the flattening and inflating operations with respect to flips:

\begin{lemma}\label{stars:lem:flip-flat}
Let~$T$ be a \ktri{k} of the \gon{(n+1)},~$f$ be a \krel{k} edge of~$T$ and~$e$ be the common bisector of the two \kstar{k}s of~$T$ containing~$f$. Let~$b$ be a \kbound{k} edge of the \gon{(n+1)}, and~$f'$ (resp.~$e'$) be the image of the edge~$f$ (resp.~$e$) in~$\flattening{T}{b}$ (resp.~$\flattening{(T\diffsym\{e,f\})}{b}$) after flattening. If~$f$ is not an edge of the \kstar{k} of~$T$ containing~$b$ (or equivalently, if~$e$ is not an edge of the \kstar{k} of~$T\diffsym\{e,f\}$ containing~$b$), then flip and flattening commute: 
$$(\flattening{T}{b})\diffsym\{e',f'\}=\flattening{(T\diffsym\{e,f\})}{b}.$$
\end{lemma}

\begin{proof}
The \ktri{k}s~$T$ and~$T\diffsym\{e,f\}$ have the same \kstar{k} containing~$b$ (since~$f$ is not an edge of this \kstar{k}). Consequently, the image of any edge of $T\cup\{e\}$ is the same in~$\flattening{T}{b}$ and~$\flattening{(T\diffsym\{e,f\})}{b}$.
\end{proof}

\begin{remark}\label{stars:rem:increasingdiameter}

In the previous lemma, it is interesting (and surprising) to observe that if~$f$ is an edge of the \kstar{k} of~$T$ containing~$b$, then the two \ktri{k}s of the \gon{n}~$\flattening{T}{b}$ and~$\flattening{(T\diffsym\{e,f\})}{b}$ are not necessarily the same. Even worst, these two \ktri{k}s can  be relatively far appart in the flip graph~$G_{n,k}$. For example, let $T_n$ be the \ktri{2} of the \gon{n} formed by all the edges adjacent to vertices~$0$ or~$2$, plus all the non-\krel{2} edges ($T_n$~is indeed a \ktri{2} since its set of \krel{2} edges is the union of two disjoint fans of edges). Then the \ktri{2}~$\flattening{T_{n+1}}{[1,n]}$ equals the minimal \ktri{2}~$T_{n,2}^{\min}$ while the \ktri{2}~$\flattening{(T_{n+1}\diffsym\{[1,n-1],[2,n]\})}{[1,n]}$ is nothing else than~$T_n$ (see \fref{stars:fig:flipflat}). These two \ktri{2}s have only $(n-5)$ \krel{2} edges in common.

\begin{figure}[!h]
	\capstart
	\centerline{\includegraphics[scale=1]{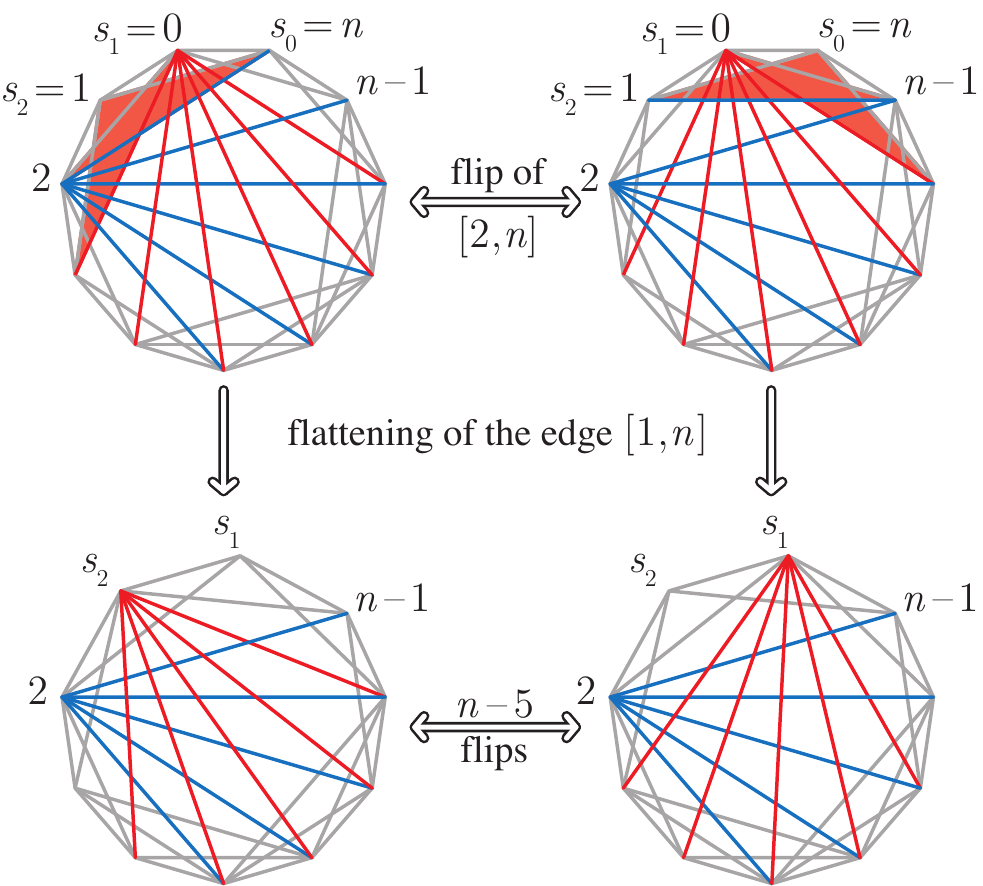}}
	\caption[When the flattening of a \kstar{2} takes away neighbors in the graph of flips]{When the flattening of a \kstar{2} takes away neighbors in the graph of flips.}
	\label{stars:fig:flipflat}
\end{figure}

Let us insist however on the fact that this property holds when~$k=1$: if~$T$ is a triangulation of the \gon{(n+1)}, with a distinguished boundary edge~$b$, if~$f$ is a \krel{1} edge of the triangle of~$T$ adjacent to~$b$, and if~$e$ is the common bisector of the two triangles of~$T$ adjacent to~$f$, then $\flattening{T}{b}=\flattening{(T\diffsym\{e,f\})}{b}$. A proof is sketched on \fref{stars:fig:flipflatk=1}. This property is interesting since it provides a simple proof of the fact that~$\delta_{n+1,1}\ge \delta_{n,1}+1$. If this property were also valid for \ktri{k}s, this proof would extend easily to the corresponding Conjecture~\ref{stars:conj:increasingdiametermax} stated in Section~\ref{stars:subsec:flips:graph}.

\begin{figure}
	\capstart
	\centerline{\includegraphics[scale=1]{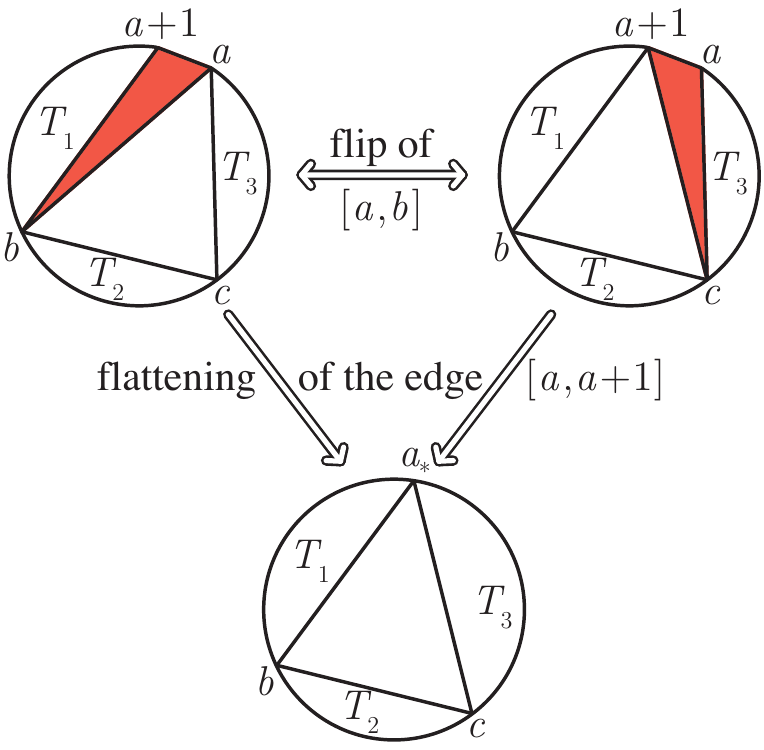}}
	\caption[Flip and flattening for a triangulation]{Flip and flattening for a triangulation.}
	\label{stars:fig:flipflatk=1}
\end{figure}
\end{remark}

\begin{corollary}\label{stars:coro:flip-infl}
Let~$T$ be a \ktri{k} of the \gon{n}, let~$f$ be a \krel{k} edge of~$T$ and let~$e$ be the common bisector of the two \kstar{k}s of~$T$ containing~$f$. Let~$X$ be an external \kcross{k} of~$T\ssm\{f\}$, and let~$f'$ (resp.~$e'$) be the image of the edge~$f$ (resp.~$e$) in~$\inflating{T}{X}$ (resp.~$\inflating{(T\diffsym\{e,f\})}{X}$) after inflating. Then flip and inflating commute:
$$(\inflating{T}{X})\diffsym\{e',f'\}=\inflating{(T\diffsym\{e,f\})}{X}.$$
\end{corollary}

To close this section, we want to observe that the result of flattening several \kbound{k} edges of a \ktri{k} is independent of the order:

\begin{lemma}\label{stars:lem:multiflatten}
Let~$b_1$ and~$b_2$ be two distinct \kbound{k} edges of a \ktri{k}~$T$. Let~$b_1'$ (resp.~$b_2'$) denote the edge of~$\flattening{T}{b_2}$ (resp.~$\flattening{T}{b_1}$) arising from~$b_1$ (resp.~$b_2$). Then~$b_1'$ (resp.~$b_2'$) is a \kbound{k} edge of~$\flattening{T}{b_2}$ (resp.~$\flattening{T}{b_1}$) and
$$\flattening{(\flattening{T}{b_2})}{b_1'}=\flattening{(\flattening{T}{b_1})}{b_2'}.$$
\end{lemma}
\vspace{-.8cm}\qed

\svs
In particular, one can define the flattening of a set of \kbound{k} edges. Similarly, it is possible to define the inflating of a set of edges-disjoint external \kcross{k}s of a \ktri{k}.


\subsection{Examples of application}\label{stars:subsec:flat-infl:examples}

We now present some applications using these inverse transformations to prove results on \ktri{k}s by induction. 

\begin{corollary}\label{stars:coro:enumerationinduction}
Let~$\theta(n,k)$ denote the number of \ktri{k}s of an \gon{n}. Then the quotient $\theta(n+1,k)/\theta(n,k)$ equals the average number, among all \ktri{k}s of the \gon{n}, of \kcross{k}s adjacent to the first~$k$ points.
\end{corollary}

\begin{proof}
Consider the set~$\cX$ of all couples~$(T,X)$ where~$T$ is a \ktri{k} of the \gon{n} and~$X$ is a \kcross{k} of~$T$ on the first~$k$ points, and the set~$\cY$ of all couples~$(T,[0,k])$ where~$T$ is a \ktri{k} of the \gon{(n+1)} where the \kbound{k} edge~$[0,k]$ is distinguished. Then the inflating \& flattening operations define opposite bijections between~$\cX$ and~$\cY$, which proves that:
$$\theta(n+1,k)=|\cY|=|\cX|=\sum_T |\{X\; k\text{-crossing of } T \text{ on the first } k \text{ points}\}|,$$
and we derive the result dividing by~$\theta(n,k)$.
\end{proof}

For example:
\begin{enumerate}[(i)]
\item For~$k=1$, we get that~$\theta(n+1,1)/\theta(n,1)$ equals the average degree of vertex~$0$ in triangulations of the \gon{n}, and we recover the well-known recursion for Catalan numbers:
$$C_{n-1}=\frac{4n-6}{n}C_{n-2}.$$
(Remember the proof of Proposition~\ref{intro:prop:catalan}.)
\item For~$n=2k+1$ we have that~$\theta(2k+1,k)=1$ (the unique \ktri{k} is the complete graph) and the number of \kcross{k}s using the first~$k$ vertices in this \ktri{k} is~$k+1$ (any choice of~$k$ of the last~$k+1$ vertices gives one \kcross{k}). In particular, we recover the fact that~$\theta(2k+2,k)=k+1$ (see Example~\ref{stars:exm:n=2k+2}).

\item Unfortunately, for~$n>2k+1$, the number of \kcross{k}s using~$k$ consecutive vertices is not independent of the \ktri{k}. Otherwise the quotient~$n\theta(n+1,k)/\theta(n,k)$ would be an integer, equal to that number (as happens in the case of triangulations).
\end{enumerate}

\mvs
The following lemma is another example of the use of a recursive argument. It is equivalent to Theorem~10~in~\cite{n-gdfcp-00}, where the proof uses flips. It is interesting to observe that \kcolorable{k} \ktri{k}s attain the following inequality (see Theorem~\ref{stars:theo:kcolorable}).

\begin{lemma}\label{stars:lem:countinternaledges}
Any \ktri{k} of the \gon{n} contains at most~$k(n-2p-1)$ \krel{p} edges, for any~$k\le p\le \frac{n-1}{2}$.
\end{lemma}

\begin{proof}
When~$n=2p+1$, it is true since there are no \krel{p} edges.

Suppose now that~$n>2p+1$. Let~$T$ be a \ktri{k} of the \gon{n}. Let~$b \eqdef [u,u+k]$ be a \kbound{k} edge of~$T$. It is easy to check that if~$g$ is a \krel{k} edge of~$T$ of length~$\ell$, then the corresponding edge~$g'$ in~$\flattening{T}{b}$ has length~$\ell$ or~$\ell-1$, and the latter is possible only if~$g=[v,v+\ell]$ with~$v\cle u\cl u+k\cle v+\ell$.

Let~$E$ be the \kcross{k} of~$\flattening{T}{b}$ that arises from the edges of the \kstar{k} of~$T$ containing~$b$. Let~$F$ be the set of edges of~$E$ of length at least~$p$. Let~$G$ denote the set of edges of length~$p$ of~$(\flattening{T}{b})\ssm E$ arising from an edge of length~$p+1$ of~$T$.

Any non-\krel{p} edge of~$(\flattening{T}{b})\ssm (E\cup G)$ (resp.~of~$E\ssm F$) arises from one (resp.~two) non-\krel{p} edge of~$T$. Thus we obtain by induction that the number of \krel{p} edges is at most~$k(n-1-2p-1)+|F|+|G|$.

To conclude, we only have to observe that any two edges of~$F\cup G$ are crossing, which implies that~$|F|+|G|=|F\cup G|\le k$.
\end{proof}

We finish this section by reinterpreting a construction we have done in Section~\ref{stars:sec:flips}:

\begin{remark}
Remember the construction in the proof of Lemma~\ref{stars:lem:lowerbounddiameter}: we started from a \ktri{k}~$\tau$ of the \gon{m} and constructed a \ktri{k}~$T$ of the \gon{(2m)} as the union of the sets $A=\ens{[2i,2(i+j)-1]}{i\in\Z_m \text{ and } j\in[k]}$, $B=\ens{[2i+1,2i+1+j]}{i\in\Z_m \text{ and } j\in[k]}$ and $C=\ens{[2i,2j]}{[i,j] \text{ edge of } \tau}$ (see \fref{stars:fig:optdiameter}).

For all $i\in\Z_n$, define the \kcross{k}~$X_i=\{[i,i-k],[i-1,i-k-1],\dots,[i-k+1,i-2k+1]\}$ of $\tau$. We just want to observe that the \ktri{k}~$T$ can be seen as the \ktri{k} obtained from~$\tau$ by sequentially inflating the \kcross{k}s $X_{m-1},X_{m-2},\dots,X_0$:
$$T=\inflating{(\inflating{(\inflating{\tau}{X_{m-1}})}{X_{m-2}})\dots)}{X_0}.$$
The inflated \kstar{k}s in the resulting \ktri{k} $T$ of the \gon{2n} are the \kstar{k}s with at least one odd vertex.
\end{remark}


\section{Decompositions of surfaces}\label{stars:sec:surfaces}

\index{polygonal decomposition of surface}
Regarding a \ktri{k}~$T$ of the \gon{n} as a complex of \kstar{k}s naturally defines a polygonal complex~$\cC(T)$ as follows:
\begin{enumerate}[(i)]
\item the vertices of~$\cC(T)$ are the vertices of the \gon{n};
\item the edges of~$\cC(T)$ are the \kbound{k} edges and \krel{k} edges of~$T$;
\item the facets of~$\cC(T)$ are the \kstar{k}s of $T$, considered as \gon{(2k+1)}s.
\end{enumerate}
Since every \krel{k} edge belongs to two \kstar{k}s and every \kbound{k} edge belongs to one, this complex is a surface with boundary, with $\gcd(n,k)$ boundary components. Also, it is orientable since the natural orientation of each \kstar{k} can be inherited on each polygon. Hence, from its Euler characteristic, we derive its genus:
$$
g_{n,k} \eqdef \frac{1}{2}(2-f+e-v-b)=\frac{1}{2}(2-n+k+kn-2k^2-\gcd(n,k)).
$$
That is, the surface does not depend on the \ktri{k}~$T$ but only on~$n$ and~$k$. We denote by~$\cS_{n,k}$ this surface. The polygonal complex~$\cC(T)$ of each \ktri{k}~$T$ provides a polygonal decomposition of~$\cS_{n,k}$. In other words~\cite[Chapter~3.1,~p.~78]{mt-gs-01}, we consider a \defn{\mbox{$2$-cell} embedding} of the \krel{k} plus \kbound{k} edges of~$T$ in the surface~$\cS_{n,k}$, where each cell is a \gon{(2k+1)}. 

\begin{figure}
	\capstart
	\centerline{\includegraphics[width=\textwidth]{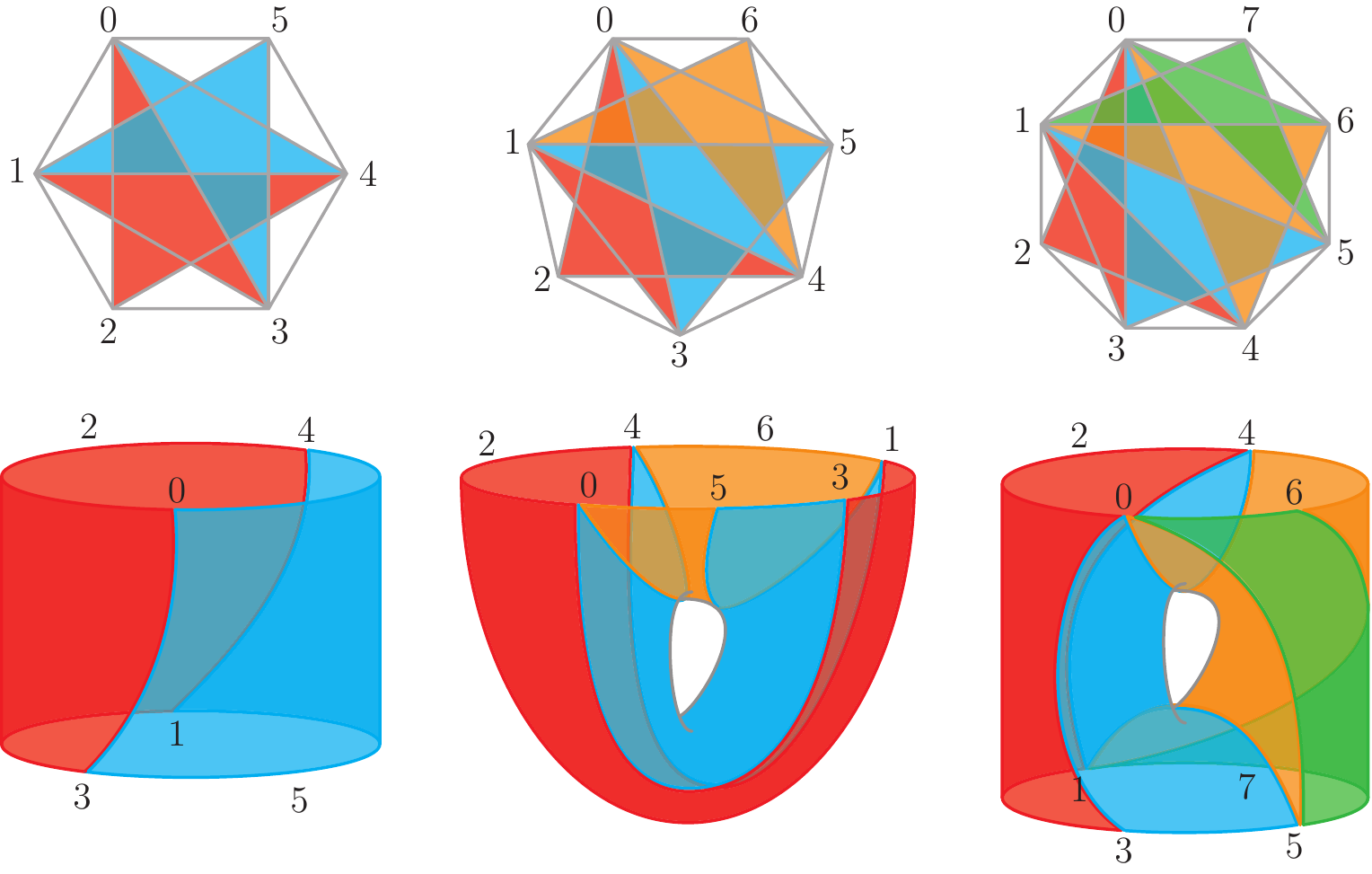}}
	\caption[Examples of decompositions of surfaces associated to \ktri{2}s]{Examples of decompositions of surfaces associated to the \ktri{2}s~$T_{6,2}^{\min}$, $T_{7,2}^{\min}$ and~$T_{8,2}^{\min}$.}
	\label{stars:fig:surfaces}
\end{figure}

\begin{example}
When~$k=2$, the genus of the surface~$\cS_{n,2}$ is~$g_{n,2}=\Fracfloor{n-5}{2}$. In \fref{stars:fig:surfaces}, we have represented the polygonal decompositions~$\cC(T_{n,2}^{\min})$ for~$n=6,7$ and~$8$.
\end{example}

\begin{remark}
Observe that two polygons in~$\cC(T)$ can intersect improperly, sharing for example two opposite edges (see for example \fref{stars:fig:surfaces}), but that each connected component of the intersection of any two polygons~$P,Q\in\cC(T)$ is either a common vertex, or a common edge of~$P$ and~$Q$.
\end{remark}


\subsection{Flips and twists}\label{stars:subsec:surfaces:flips}

We can now interpret the operation of flip in this new setting.

Let~$T$ be a \ktri{k} of the \gon{n}, and let~$f$ be a \krel{k} edge of~$T$. Let~$R$ and~$S$ be the two \kstar{k}s of~$T$ containing~$f$, and let~$e$ be the common bisector of~$R$ and~$S$. Let~$X$ and~$Y$ be the two \kstar{k}s of~$T\diffsym\{e,f\}$ containing~$e$.

Then~$T\ssm\{f\}$ can be viewed as a decomposition of~$\cS_{n,k}$ into ${n-2k-2}$ \gon{(2k+1)}s and one~\gon{4k}, obtained from~$\cC(T)$ by gluing the two~\gon{(2k+1)}s corresponding to~$R$ and~$S$ along the edge~$f$. And then~$T\diffsym\{e,f\}$ is obtained from~$T\ssm\{f\}$ by splitting the \gon{4k} into the two~\gon{(2k+1)}s corresponding to~$X$ and~$Y$.

\begin{example}
In \fref{stars:fig:monodromy}, we have drawn the decomposition of the cylinder corresponding to the \ktri{2} $T_{6,2}^{\min}$. The second \ktri{2} is obtained from $T_{6,2}^{\min}$ by flipping the edge $[1,4]$, and we have represented the decomposition of the cylinder obtained by flipping the edge on the surface. If we now flip $[0,3]$ and then $[2,5]$, we obtain again the triangulation $T_{6,2}^{\min}$. 

\begin{figure}
	\capstart
	\centerline{\includegraphics[scale=.9]{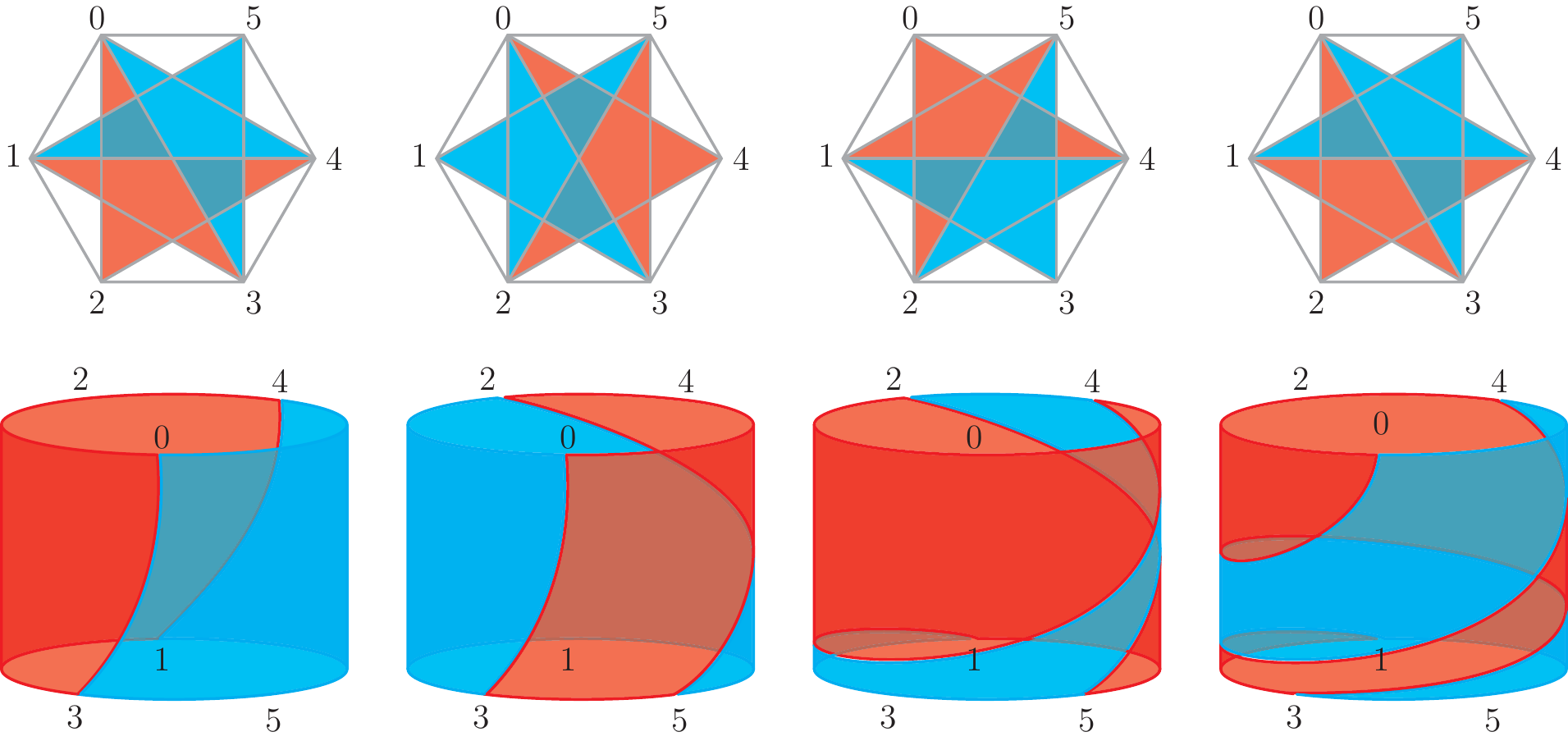}}
	\caption[Flips and twists]{A cycle of flips gives raise to a twist.}
	\label{stars:fig:monodromy}
\end{figure}

Observe that although the corresponding decomposition is combinatorially the same as in the first picture, it has been ``twisted'' on the surface. This phenomenon can be generalized to the following proposition:
\end{example}

\begin{proposition}\label{stars:prop:monodromy}
\index{mapping class group}
There is a homomorphism between:
\begin{enumerate}[(i)]
\item the \emph{fundamental group}~$\pi_{n,k}$ of the graph of flips~$G_{n,k}$, \ie the set of closed walks in~$G_{n,k}$, based at an initial \ktri{k}, considered up to homotopy, and
\item the \emph{mapping class group}~$\cM_{n,k}$ of the surface~$\cS_{n,k}$, \ie the set of diffeomorphisms of the surface~$\cS_{n,k}$ into itself that preserve the orientation and that fix the boundary of~$\cS_{n,k}$, considered up to isotopy~\cite{b-blmcg-74,i-mcg-02}.\qed
\end{enumerate}
\end{proposition}

Observe that we know the structure of the fundamental group~$\pi_{n,k}$: it is a free group, and its number of generators is the number of edges in~$G_{n,k}$ minus the number of vertices in~$G_{n,k}$ plus~$1$. This number grows as~$4^{nk}$ modulo a rational function of degree~$O(k^2)$ in~$n$ (see Remark~\ref{ft:rem:enumeration}). In contrast, the mapping class group~$\cM_{n,k}$ is generated by~$2g+1$ Dehn twists (see~\cite{i-mcg-02} for the definition). It may be interesting to understand the image and the kernel of this homomorphism. In particular, if this homomorphism is surjective, this interpretation provides a new combinatorial description of the mapping class group of~$\cS_{n,k}$.


\subsection{The dual graph}\label{stars:subsec:surfaces:dual}

\begin{figure}
	\capstart
	\centerline{\includegraphics[width=\textwidth]{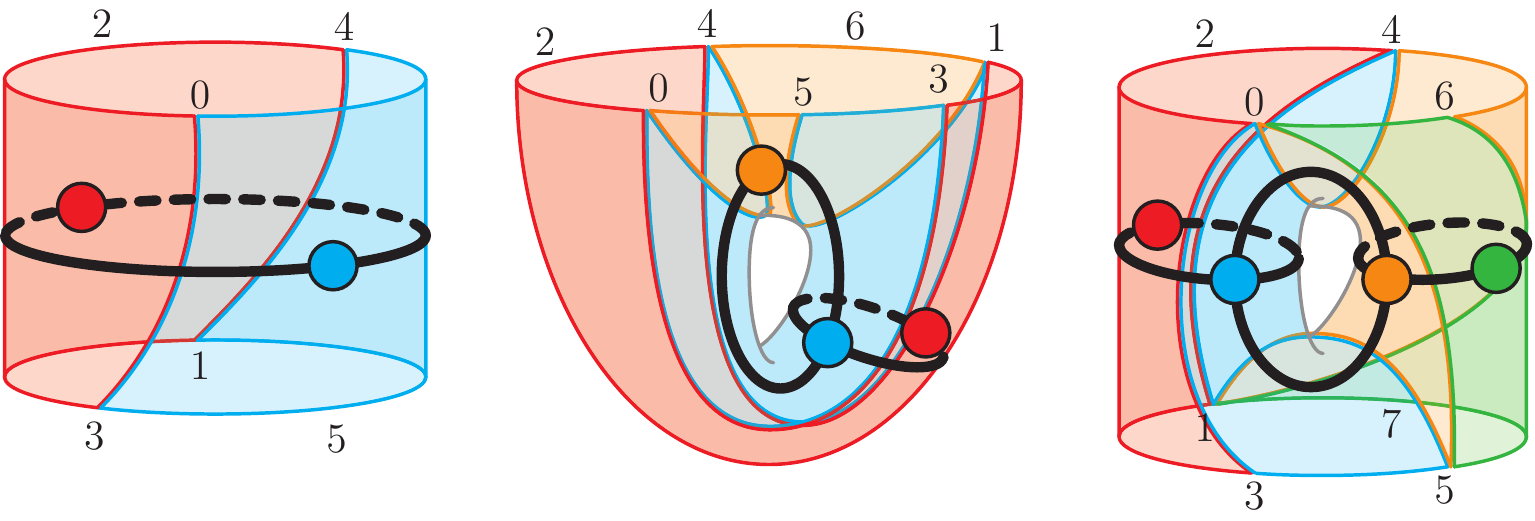}}
	\caption[The dual maps of the surface decompositions of \fref{stars:fig:surfaces}]{The dual maps of the surface decompositions of \fref{stars:fig:surfaces}.}
	\label{stars:fig:surfacesdual}
\end{figure}

\index{dual!--- graph of a \ktri{k}}
We now consider the dual map of the $2$-cell-embedding~$\cC(T)$ on the surface~$\cS_{n,k}$: we denote~$T^\dual$ the graph whose vertices are the \kstar{k}s of~$T$ and with an edge between two \kstar{k}s for each of their common \krel{k} edge. Remember that two \kstar{k}s can have until~$k$ common edges, which means that~$T^\dual$ is a multigraph. Observe also that it is naturally embedded on the surface~$\cS_{n,k}$ as the dual map of $\cC(T)$ (see \fref{stars:fig:surfacesdual}). Again, this embedding is ``naturally'' defined up to a boundary preserving diffeomorphism of~$\cS_{n,k}$.

By definition, this dual graph~$T^\dual$ has~$n-2k$ vertices and~$k(n-2k-1)$ edges. The degree of each vertex is between~$k$ and~$2k+1$. It is also interesting to observe what happens when we cut the surface~$\cS_{n,k}$ along~$T^\dual$:

\begin{lemma}\label{stars:lem:surfacesdual}
Cutting along the dual map~$T^\dual$ decomposes the surface~$\cS_{n,k}$ into~$\gcd(n,k)$ cylinders. In each of these cylinders, one boundary is a boundary of~$\cS_{n,k}$ and the other one is formed only by edges of~$T^\dual$.
\end{lemma}

\begin{proof}
When we have cut along all edges of~$T^\dual$, the polygons adjacent to a given boundary component~$B$ of~$\cS_{n,k}$ form a cycle of polygons, glued along edges. Since the surface is orientable, we obtain a cylinder, with~$B$ as one boundary, and a sequence of edges of~$T^\dual$ as the other boundary. Since all the vertices of~$\cC(T)$ are on the boundary of~$\cS_{n,k}$ there is no other pieces.
\end{proof}

For example, when~$k=1$, cutting along a tree included into a disk gives raise to a cylinder. It is intereseting to observe that duals of triangulations produce in fact all graphs which decompose the disk into a cylinder and whose vertices have degree between~$1$ and~$3$.


\subsection{Equivelar polygonal decompositions}\label{stars:subsec:surfaces:equivelar}

A polygonal decomposition of a surface is \defn{\equivelar{(p,q)}}\index{equivelar (decomposition of surface)} if all its faces are \gon{p}s and all its vertices have degree~$q$. Equivelarity is a local regularity condition: any regular polygonal decomposition (\ie for which the automorphism group acts transitively on the flags), is clearly equivelar. For example, the boundary complexes of the Platonic solids define equivelar polygonal decompositions of the \mbox{$2$-sphere}. Because of their symmetry, regular and equivelar decompositions of surfaces have received particular attention (see~\cite{bs-pm-97} and the references therein).

In this section, we use multitriangulations to derive equivelar polygonal decompositions of an infinite family of surfaces of high genus. The result is the following:

\begin{theorem}\label{stars:theo:decomp}
For any integers~$\ell\ge 2$ and~$m\ge 1$, there exists a \equivelar{\big(2\ell m+1,2(\ell m-m+1)\big)} polygonal decomposition of the closed surface of genus 
$$g \eqdef (\ell-1)\left(\ell^2m^2-\frac{\ell m}{2}-1\right),$$
in which no two polygons share an edge.
More precisely, we construct a proper polygonal decomposition of this surface with:
\begin{itemize}
\item $\ell(2\ell m+1)$~vertices, each of degree~$2(\ell m-m+1)$;
\item $\ell(\ell m-m+1)(2\ell m+1)$~edges;~and
\item $2\ell(\ell m-m+1)$~polygons, each of valence~$(2\ell m+1)$, and no two of which share~two~edges.
\end{itemize}
\end{theorem}

We prove this theorem by constructing a multitriangulation~$T$ whose associated polygonal decomposition~$\cC(T)$ is equivelar. Choose~$\ell\ge 2$ and~$m\ge 1$. Let~$k \eqdef \ell m$ and~${n \eqdef \ell(2\ell m+1)}$. Since~$\gcd(n,k)=\ell$, the surface~$\cS_{n,k}$ has $\ell$~boundary components, each of length~$2\ell m+1$, and has genus $g$. Thus, we obtain a closed surface~$\hat{\cS}_{n,k}$ of genus $g$ by gluing a \gon{(2\ell m+1)} on each boundary component of~$\cS_{n,k}$. Any \ktri{k} of the \gon{n} defines a decomposition of~$\hat{\cS}_{n,k}$ into $2\ell(\ell m-m+1)$ polygons of valence $(2\ell m+1)$, with $\ell(\ell m-m+1)(2\ell m+1)$~edges and $\ell(2\ell m+1)$~vertices. The end of this section is devoted to the construction of such a \ktri{k}~$T$ with the property that:
\begin{enumerate}[(i)]
\item the degree of any vertex in~$T$ is~$2m(2\ell-1)$ (thus, the degree in $\cC(T)$ is~$2(\ell m-m+1)$);~and
\item any \kstar{k} of~$T$ shares at most one edge with each boundary component and with each other \kstar{k} of~$T$.
\end{enumerate}

\mvs
We consider again (see \fref{stars:fig:zigzags}) the \kzz{k}~$Z$ of~$E_n$ defined by:
$$Z  \eqdef  \ens{[q-1,-q-k]}{1\le q\le \Fracfloor{n-2k}{2}}\cup\ens{[q,-q-k]}{1\le q\le \Fracfloor{n-2k-1}{2}}.$$
Let~$\rho$ be the rotation~$\rho:t\mapsto t+2k+1$, let~$\theta$ be the rotation~$\theta:t\mapsto t+2$ and let~$\sigma$ be the reflection~$\sigma:t\mapsto -t$. Let~$\Gamma$ denote the dihedral group of order~$2\ell$ generated by~$\rho$ and~$\sigma$.

\begin{lemma}
\index{zigzag@\kzz{k}}
The set
$$\bigcup_{i=0}^{m-1} \theta^i\left(\bigcup_{\gamma\in\Gamma} \gamma(Z)\right)$$
is the subset of the \krel{k} edges of a \ktri{k}~$T_{l,m}$ of the~\gon{n}.
\end{lemma}

\begin{proof}
Observe that~$Z$ has a stabilizer of cardinality~$2$ in~$\Gamma$, so that~$\bigcup_{\gamma\in\Gamma} \gamma(Z)$ consists of~$\ell$ (and not~$2\ell$) \kzz{k}s of~$E_n$. Indeed,
\begin{enumerate}[(i)]
\item if~$\ell$ is odd, then~$\sigma(Z)=Z$, and~$\bigcup_{\gamma\in\Gamma} \gamma(Z)$ is the disjoint union of the~$\ell$~\kzz{k}s~$\rho^i(Z)$, for~$0\le i\le \ell-1$;
\item if~$\ell$ is even, then~$\rho^{\ell/2}(Z)=Z$, and~$\bigcup_{\gamma\in\Gamma} \gamma(Z)$ is the disjoint union of the $\ell/2$~zigzags~$\rho^i(Z)$, for~$0\le i\le \ell/2-1$, and the $\ell/2$~zigzags~$\rho^i\circ\sigma(Z)$, for~$0\le i\le \ell/2-1$.
\end{enumerate}
Consequently, the set~$\bigcup_{i=0}^{m-1} \theta^i\left(\bigcup_{\gamma\in\Gamma} \gamma(Z)\right)$ is the union of~$\ell m=k$~disjoint \kzz{k}s, which proves the result by Lemma~\ref{stars:lem:unionzigzags}.
\end{proof}

In the end of this section, we always denote by~$\bar\Gamma$ the quotient of~$\Gamma$ by the stabilizer of~$Z$, so that each \kzz{k} of~$T_{l,m}$ is uniquely represented as the image of~$Z$ under a transformation of~$\bar\Gamma$.

\begin{figure}[!b]
	\capstart
	\centerline{\includegraphics[scale=1]{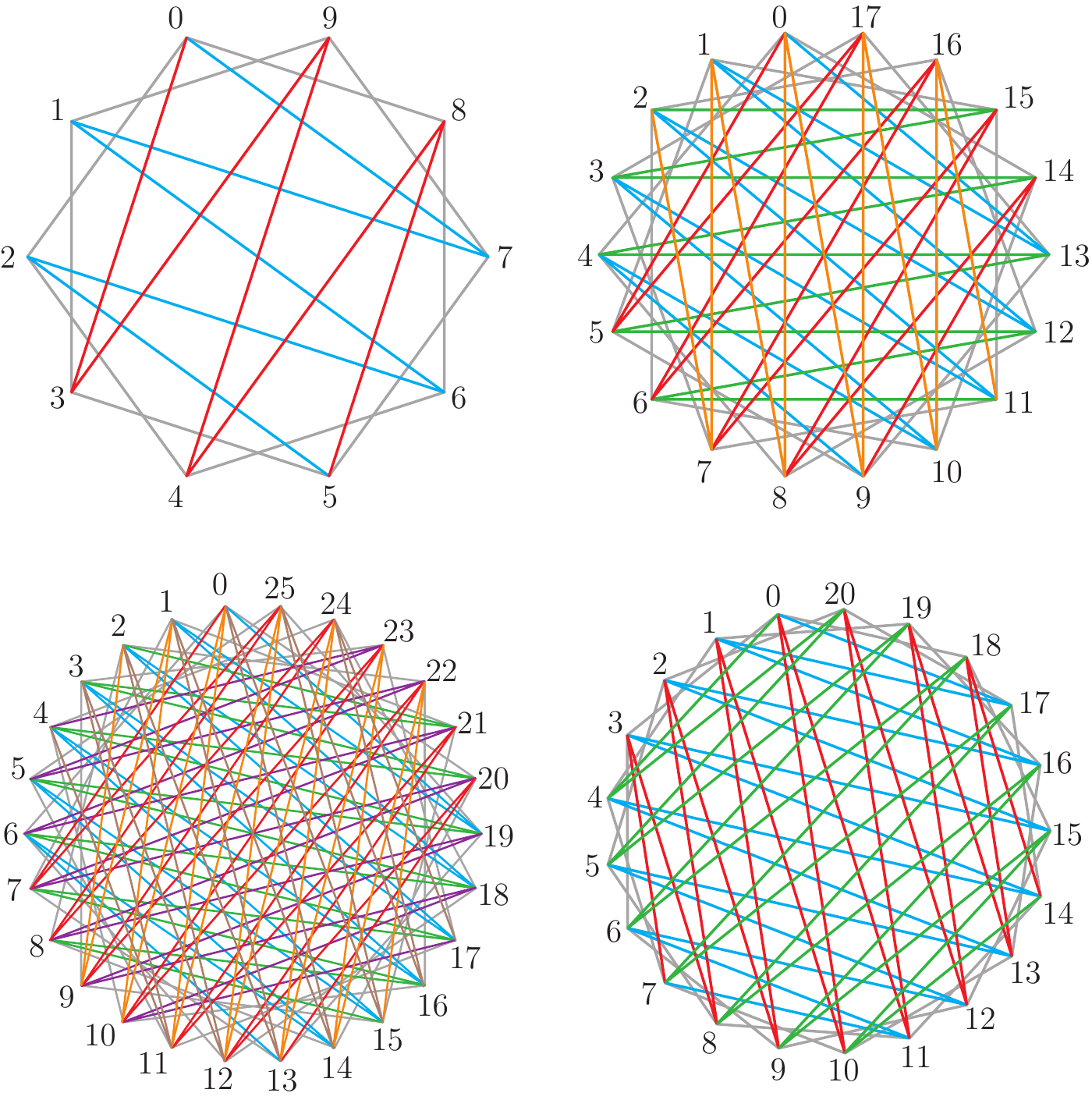}}
	\caption[Four very regular multitriangulations]{The very regular multitriangulation~$T_{l,m}$ for~$(\ell,m,n,k)=(2,1,10,2)$ (top left), $(2,2,18,4)$ (top right), $(2,3,26,6)$ (bottom left), and $(3,1,21,3)$ (bottom right).}
	\label{stars:fig:equivelar}
\end{figure}

\begin{example}
\fref{stars:fig:equivelar} shows these multitriangulations for the parameters
$$(\ell,m,n,k)\in\{(2,1,10,2),(2,2,18,4),(2,3,26,6),(3,1,21,3)\}.$$
\end{example}

\begin{lemma}\label{stars:lem:constantdegree}
The degree of any vertex in~$T_{l,m}$ is~$2m(2\ell-1)$.
\end{lemma}

\begin{proof}
Any vertex~$v$ of~$V_n$ is contained in~$2(\ell-1)$~edges of~$\bigcup_{\gamma\in\Gamma} \gamma(Z)$. Consequently, $v$~is contained in $2m(l-1)$~\krel{k} edges of~$T_{l,m}$. Since it is contained in $2k$~\kirrel{k} or \kbound{k} edges, the degree of~$v$ is~$2m(\ell-1)+2m\ell=2m(2\ell-1)$.
\end{proof}

\begin{lemma}\label{stars:lem:boundaryedges}
Any \kstar{k} of~$T_{l,m}$ contains at most one edge on each boundary component of~$\cS_{n,k}$.
\end{lemma}

\begin{proof}
Since $T_{l,m}$ does not contain any pair of consecutive \kear{k}s, any \kstar{k} of $T_{l,m}$ contains at most two \kbound{k} edges (Lemma~\ref{stars:lem:starwithboundary}). These \kbound{k} edges are consecutive (for the cyclic order), and therefore they are not in the same boundary component of~$\cS_{n,k}$.
\end{proof}

It only remains to show that two \kstar{k}s of~$T_{l,m}$ cannot share more than one edge. For this, we need to understand how \kstar{k}s and \kzz{k}s of~$T_{l,m}$ are arranged together. Let~$\bar Z$ be the unique \kzz{(k-1)} obtained by extending~$Z$ with to \kbound{k} edges.

\begin{lemma}\label{stars:lem:starcolor}
Any \kstar{k} of~$T_{l,m}$ contains  exactly one \kbound{k} edge of $T_{l,m}\ssm\cup_{\gamma\in\bar\Gamma} \gamma(\bar Z)$ and exactly one angle of the zigzag~$\gamma(\bar Z)$ for each~$\gamma\in\bar\Gamma$. Furthermore, the star order of these angles is given by the order of the set~$\ens{\gamma(0)}{\gamma\in\bar\Gamma}$ around the circle.
\end{lemma}

\begin{proof}
Since any two angles of a \kstar{k} define a crossing, a \kstar{k} can contain at most one angle of each zigzag. Since each zigzag $\gamma(\bar Z)$ has~$n-2k$ angles, each \kstar{k} of~$T_{l,m}$ has to contain exactly one angle of each zigzag. The last edge of each \kstar{k} is in no zigzags. Finally, the star order of the angles is given by the slope of the zigzags, which can be read on the positions of~$\gamma(0)$.
\end{proof}

\begin{lemma}\label{stars:lem:adjacency}
Any two \kstar{k}s of~$T_{l,m}$ share at most one edge. 
\end{lemma}

\begin{proof}
Assume that two \kstar{k}s~$R$~and~$S$ share two edges. By symmetry, we can assume that these edges belong to the \kzz{k}s~$Z$ and~$\gamma(Z)$, for a certain~$\gamma\in\bar\Gamma$. Let $[r,r+k]$ denote the unique edge of~$R$ not contained in $\lambda(\bar Z)$, for any~$\lambda\in\bar\Gamma$. Define similarly the edge~$[s,s+k]$ of~$S$. Denote by~$[a,b]$~and~$[c,d]$ the common edges of~$R$~and~$S$ with $r\cle a \cl b \cle r+k$~and $s\cle c\cl d\cle s+k$. The previous lemma ensures that $|b-a|=|d-c|$ is the number of angles appearing between~$Z$ and~$\gamma(Z)$ in star order around~$R$~and~$S$. We deduce from this that all the edges~$\ens{[a+i,c+i]}{i<|b-a|}$ are common to~$R$~and~$S$. This finally implies that the transformation~$\gamma$ satisfies~$\gamma(0)=1$ which is absurd.
\end{proof}

\begin{remark}
Observe however that two \kstar{k}s of~$T_{l,m}$ can share one edge plus some other vertices.
\end{remark}

\svs
Finally, we want to observe the following properties of this polygonal decomposition:
\begin{enumerate}
\item When~$\ell$ (or~$m$) is big, the number of edges is equivalent to twice the genus, which is optimal.
\item When~$m=1$ (see \fref{stars:fig:equivelar}), the degree of the vertices and the number of edges of the polygons differ just by~$1$.
\item When~$\ell=2$ (see \fref{stars:fig:equivelar}), we have an equivelar decomposition of a surface of genus~$g$ with~$4\sqrt{g}$ vertices of degree~$\sqrt{g}$, with $2g$~edges, and with $2\sqrt{g}$~polygons of degree $2\sqrt{g}$.
\end{enumerate}
Observe furthermore that the multitriangulation~$T_{l,m}$, and thus the corresponding surface decomposition, is invariant by the dihedral group~$\Gamma$.

	\chapter{Multipseudotriangulations}\label{chap:mpt}

In this chapter, we interpret multitriangulations in the line space of the plane. The set of all bisectors of a star is a pseudoline. Moreover, all pseudolines dual to the \kstar{k}s of a \ktri{k} of the \gon{n} form a \defn{pseudoline arrangement with contact points} whose support covers precisely the dual pseudoline arrangement of the vertex set of the \gon{n} except its first~$k$ levels. It relates multitriangulations with \defn{pseudotriangulations}, for which a similar interpretation is known, and thus provides an explanation to their common properties.

Motivated by these two examples, we investigate in Section~\ref{mpt:sec:enumeration} the set of all pseudoline arrangements with contact points which cover a given support. We define a notion of \defn{flip} between them, which extends the flip for multitriangulations and pseudotriangulations, and we study the graph of these flips. In particular, we provide an enumeration algorithm for pseudoline arrangements with a given support (similar to the enumeration algorithm of~\cite{bkps-ceppgfa-06} for pseudotriangulations), based on certain \defn{greedy pseudoline arrangements}. We define these greedy pseudoline arrangements as sources for certain orientations on the flip graph and we use extensively their relationship with \defn{primitive sorting networks}~\cite[Section~5.3.4]{k-acp-73}. Our algorithm requires a polynomial running time per multipseudotriangulation and its working space is polynomial.

The duality between pseudoline arrangements with contact points and multitriangulations (or pseudotriangulations) is only presented in Section~\ref{mpt:sec:duality}. This duality leads to a natural definition of \defn{\mpt{}s}, of which we study the elementary properties (number of edges, pointedness and crossings, stars, \etc) in Section~\ref{mpt:sec:mpt}.


\section{Pseudoline arrangements with the same support}\label{mpt:sec:enumeration}


\subsection{Pseudoline arrangements in the M\"obius strip}\label{mpt:subsec:enumeration:pa}

Let~$\cM$ denote the \defn{M\"obius strip}\index{Mobius strip@M\"obius strip} (without boundary), defined as the quotient set of the plane~$\R^2$ under the map $\tau: \R^2\to\R^2,\; (x,y)\mapsto(x+\pi,-y)$. The induced canonical projection will be denoted by~$\pi:\R^2\to\cM$.

A \defn{pseudoline}\index{pseudoline|hbf} is the image under the canonical projection~$\pi$ of the graph~$\ens{(x,f(x))}{x\in\R}$ of a continuous and \defn{\piantiperiodic} function~$f:\R\to\R$ (that is, which satisfies~$f(x+\pi)=-f(x)$ for all~$x\in\R$). We say that~$f$ represents the pseudoline.

\index{crossing and contact points}
When we consider two pseudolines, we always assume that they have a finite number of \defn{intersection points}. Thus, these intersection points can only be either \defn{crossing points} or \defn{contact points} (see \fref{mpt:fig:flip}(a)). Any two pseudolines always have an odd number of crossing points (in particular, at least one). When~$\lambda$ and~$\mu$ have exactly one crossing point, we denote it by~$\lambda\wedge\mu$.

\begin{figure}
	\capstart
	\centerline{\includegraphics[scale=1]{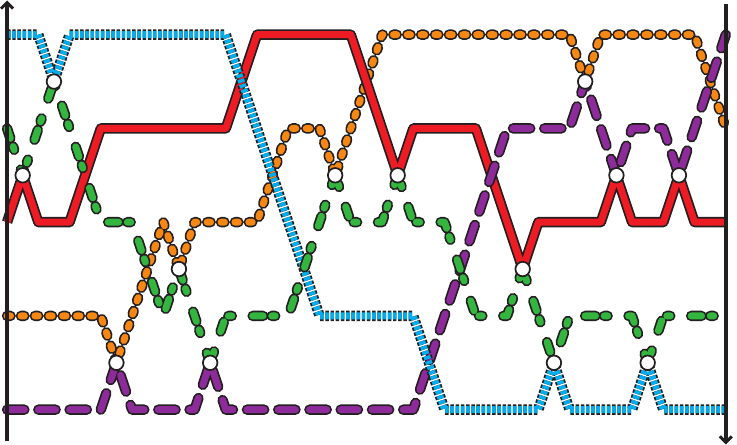}}
	\caption[A pseudoline arrangement in the M\"obius strip]{A pseudoline arrangement in the M\"obius strip. Its contact points are represented by white circles.}
	\label{mpt:fig:arrangement}
\end{figure}

A \defn{pseudoline arrangement}\index{pseudoline!--- arrangement!--- --- with contact points|hbf} is a finite set of pseudolines such that any two of them have exactly one crossing point and possibly some contact points (see \fref{mpt:fig:arrangement}). In this chapter, we are only interested in \defn{simple} arrangements, that is, where no three pseudolines meet in a common point. The \defn{support}\index{support (of an arrangement)} of a pseudoline arrangement is the union of its pseudolines. Observe that an arrangement is completely determined by its support together with its set of contact points.

\begin{remark}
\index{pseudoline!--- arrangement|hbf}
Usually, the definition of pseudoline arrangements does not allow contact points. In this chapter, they play a crucial role since we are interested in all pseudoline arrangements which share a common support, and which only differ by their sets of contact points. 

Pseudoline arrangements are also classically defined on the projective plane rather than the M\"obius strip. The projective plane is obtained from the M\"obius strip by adding a point at infinity.

For more details on pseudoline arrangements, we refer to the broad literature on the topic \cite{g-as-72,bvswz-om-99,k-ah-92} and in particular the survey paper~\cite{g-pa-97}. In Appendix~\ref{app:sec:dpl}, we discuss some classical features of arrangements (of pseudolines and double pseudolines) in the presentation of an enumeration algorithm for arrangements.
\end{remark}


\subsection{The flip graph and the greedy pseudoline arrangements}\label{mpt:subsec:enumeration:flipgraph}

\subsubsection{Flips}\label{mpt:subsubsec:flipgraph}

\begin{figure}
	\capstart
	\centerline{\includegraphics[scale=1]{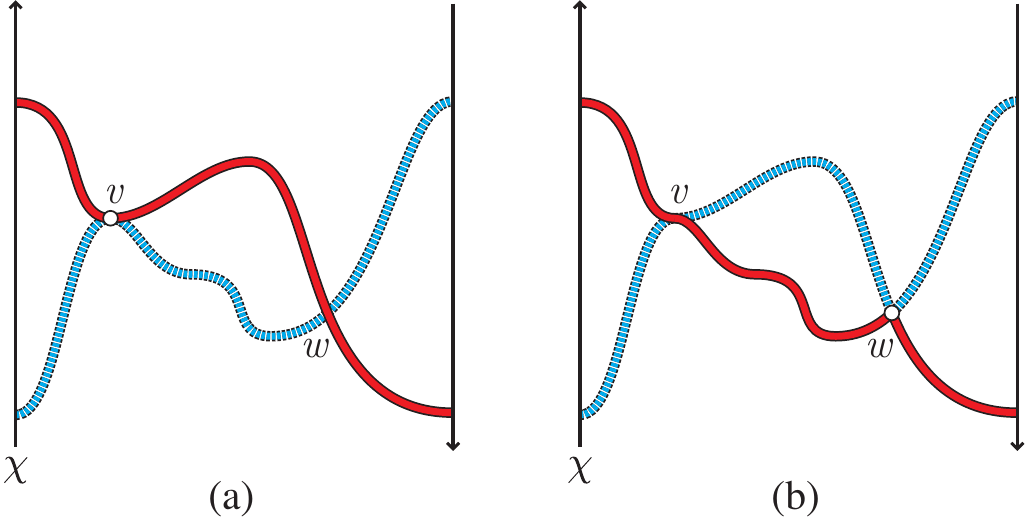}}
	\caption[Flipping a contact point in a pseudoline arrangement]{(a)~A pseudoline arrangement with one contact point~$v$ and one crossing point~$w$; (b)~Flipping~$v$ in the pseudoline arrangement of~(a).}
	\label{mpt:fig:flip}
\end{figure}

\index{flip}
\index{flip!graph of ---s}
The following lemma defines the notion of flips (see \fref{mpt:fig:flip}), to be studied in this chapter. As usual, we use the symbol~$\triangle$ for the symmetric difference:~$X\triangle Y \eqdef (X\ssm Y)\cup(Y\ssm X)$.

\begin{lemma}\label{mpt:lem:flip}
Let~$\Lambda$ be a pseudoline arrangement,~$\cS$ be its support, and~$V$ be the set of its contact points. Let~$v\in V$ be a contact point of two pseudolines of~$\Lambda$, and~$w$ denote their crossing point. Then~$V\triangle\{v,w\}$ is also the set of contact points of a pseudoline arrangement~$\Lambda'$ supported~by~$\cS$.
\end{lemma}

\begin{proof}
Let~$f$ and~$g$ represent the two pseudolines~$\lambda$ and~$\mu$ of~$\Lambda$ in contact at~$v$. Let~$x$~and~$y$ be such that~$v=\pi(x,f(x))$, $w=\pi(y,f(y))$~and~$x<y<x+\pi$. We define two functions~$f'$~and~$g'$~by:
$$ f' \eqdef 
	\begin{cases}
		f & \text{on } [x,y]+\Z\pi, \\
		g & \text{otherwise,}
	\end{cases}
\qquad \text{and} \qquad g' \eqdef 
	\begin{cases}
		g & \text{on } [x,y]+\Z\pi, \\
		f & \text{otherwise.}
	\end{cases}
$$
These two functions are continuous and \piantiperiodic, and thus define two pseudolines~$\lambda'$~and~$\mu'$. These two pseudolines have a contact point at~$w$ and a unique crossing point at~$v$, and they cross any pseudoline~$\nu$ of~$\Lambda\ssm\{\lambda,\mu\}$ exactly once (since~$\nu$ either crosses both $\lambda$~and~$\mu$ between $v$~and~$w$, or crosses both $\lambda$~and~$\mu$ between~$w$~and~$v$). Consequently, $\Lambda' \eqdef \Lambda\triangle\{\lambda,\mu,\lambda',\mu'\}$ is a pseudoline arrangement, with support~$\cS$, and contact points~$V\triangle\{v,w\}$.
\end{proof}

The pseudoline arrangements~$\Lambda$~and~$\Lambda'$ are the only two pseudoline arrangements supported by~$\cS$ whose sets of contact points contain~$V\ssm\{v\}$. We call \defn{flip} the transformation from $\Lambda$~to~$\Lambda'$, interchanging the roles of~$v$~and~$w$. The central object in this chapter is the flip~graph:

\begin{definition}
Let~$\cS$ be the support of a pseudoline arrangement. The \defn{flip graph} of~$\cS$, denoted by~$G(\cS)$, is the graph whose vertices are all the pseudoline arrangements supported by~$\cS$, and whose edges are flips between them.
\end{definition}

In other words, there is an edge in the graph~$G(\cS)$ between two pseudoline arrangements if and only if their sets of contact points differ by exactly two points. Observe that the graph~$G(\cS)$ is regular: there is one edge adjacent to a pseudoline arrangement~$\Lambda$ supported by~$\cS$ for each contact point of~$\Lambda$, and two pseudoline arrangements with the same support have the same number of contact points.

\begin{example}
The flip graph of the support of an arrangement of two pseudolines with~$p$ contact points is the complete graph on~$p+1$~vertices.
\end{example}

\subsubsection{Unique source orientations}\label{mpt:subsubsec:orientations}

Let~$\cS$ be the support of a pseudoline arrangement and~$\bar{\cS}$ denote its preimage under the projection~$\pi$. We orient the graph~$\bar{\cS}$ along the abscissa increasing direction, and the graph~$\cS$ by projecting the orientations of the edges of~$\bar{\cS}$. We denote by~$\cle$ the induced partial order on the vertex set of~$\bar{\cS}$ (defined by~$z\cle z'$ if there exists an oriented path on~$\bar{\cS}$ from~$z$~to~$z'$).

A \defn{filter} of~$\bar{\cS}$ is a proper set~$F$ of vertices of~$\bar{\cS}$ such that~$z\in F$ and~$z\cle z'$ implies~$z'\in F$. The corresponding \defn{antichain} is the set of all edges and faces of~$\bar{\cS}$ with one vertex in~$F$ and one vertex not in~$F$. This antichain has a linear structure, and thus, can be seen as the set of edges and faces that cross a vertical curve~$\bar \chi$ of~$\R^2$. The projection~$\chi \eqdef \pi(\bar \chi)$ of such a curve is called a \defn{cut}\index{cut} of~$\cS$. We see the fundamental domain located between the two curves~$\bar \chi$ and~$\tau(\bar \chi)$ as the result of cutting the M\"obius strip along the cut~$\chi$. For example, we use such a cut to represent pseudoline arrangements in all figures of this chapter. See for example \fref{mpt:fig:flip}.

The cut~$\chi$ defines a partial order~$\cle_\chi$ on the vertex set of~$\cS$: for all vertices~$v$~and~$w$ of~$\cS$, we write~$v\cle_\chi w$ if there is an oriented  path in~$\cS$ which does not cross~$\chi$. In other words,~$v\cle_\chi w$ if~$\bar v\cle\bar w$, where~$\bar v$ (resp.~$\bar w$) denotes the unique preimage of~$v$ (resp.~$w$) between~$\bar\chi$ and~$\tau(\bar\chi)$. For example, in the arrangements of \fref{mpt:fig:flip}, we have~$v\cl_\chi w$.

Let~$\Lambda$ be a pseudoline arrangement supported by~$\cS$, $v$~be a contact point between two pseudolines of~$\Lambda$ and~$w$~denote their crossing point. Since~$v$~and~$w$ lie on a same pseudoline on~$\cS$, they are comparable for~$\cle_\chi$. We say that the flip of~$v$ is \defn{\increasing{\chi}}\index{flip!increasing and decreasing ---s|hbf} if~$v\cl_\chi w$ and \mbox{\defn{\decreasing{\chi}}} otherwise. For example, the flip from~(a) to~(b) in \fref{mpt:fig:flip} is \increasing{\chi}. We denote by~$G_\chi(\cS)$ the directed graph of \increasing{\chi} flips on pseudoline arrangements supported by~$\cS$.

\begin{theorem}\label{mpt:theo:uniquesource}
The directed graph~$G_\chi(\cS)$ is an acyclic connected graph with a unique source.
\end{theorem}

\begin{definition}
\index{greedy!--- pseudoline arrangement|hbf}
\index{pseudoline!--- arrangement!greedy --- ---|hbf}
This unique source is denoted by~$\Gamma_\chi(\cS)$ and called the \defn{\greedy{\chi} pseudoline arrangement} on~$\cS$.
\end{definition}

\begin{proof}[Proof of Theorem~\ref{mpt:theo:uniquesource}]
If~$\Lambda$~and~$\Lambda'$ are two pseudoline arrangements supported by~$\cS$, we write $\Lambda\tle_\chi \Lambda'$ if there exists a bijection~$\phi$ between their sets of contact points such that~$v\cle_\chi \phi(v)$ for any contact point~$v$ of~$\Lambda$. It is easy to see that this relation is a partial order on the vertices of~$G_\chi(\cS)$. Since the edges of~$G_\chi(\cS)$ are oriented according to~$\tle_\chi$, the graph~$G_\chi(\cS)$ is acyclic. (In fact,~$G_\chi(\cS)$ is even the Hasse diagram of the partial order~$\tle_\chi$.)

That~$G_\chi(\cS)$ has a unique source will be proved next in Theorem~\ref{mpt:theo:sorting}, where we will provide a characterization of~$\Gamma_\chi(\cS)$ in terms of sorting networks.

Finally, an acyclic directed graph with only one source is necessarily connected.
\end{proof}

\subsubsection{Sorting networks}\label{mpt:subsubsec:sorting}

\begin{figure}
	\capstart
	\centerline{\includegraphics[scale=1]{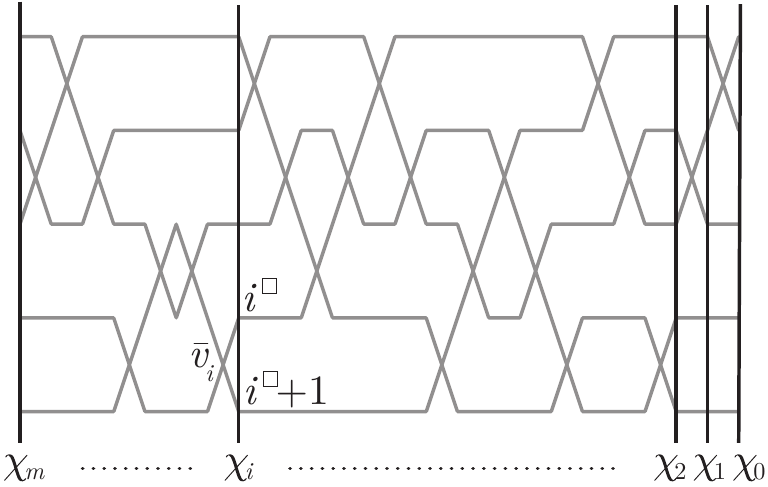}}
	\caption[Sweeping a pseudoline arrangement]{A (backward) sweep of the support of the pseudoline arrangement of \fref{mpt:fig:arrangement}.}
	\label{mpt:fig:sweep}
\end{figure}

Let~$n$ denote the number of pseudolines of the arrangements supported by~$\cS$ and $m\ge {n \choose 2}$ their number of intersection points (crossing points plus contact points). We consider a chain $F=F_m\supset F_{m-1}\supset \dots\supset F_1\supset F_0=\tau(F)$ of filters of~$\bar{\cS}$ such that two consecutive of them~$F_i$~and~$F_{i+1}$ only differ by a single element:~$\{\bar v_i\} \eqdef F_{i+1}\ssm F_i$. This corresponds to a (backward) \defn{sweep}\index{sweep (of an arrangement)}~$\chi=\chi_0,\chi_1,\dots,\chi_m=\chi$ of the M\"obius strip, where each cut~$\chi_{i+1}$ is obtained from the cut~$\chi_i$ by sweeping a maximal vertex~$v_i \eqdef \pi(\bar v_i)$ of~$\cS$ (for the partial order~$\cle_\chi$). For all~$i$, let~$e_i^1,e_i^2,\dots,e_i^n$ denote the sequence of edges of~$\bar{\cS}$ with exactly one vertex in~$F_i$, ordered from top to bottom, and let~$i^\square$ be the index such that~$\bar v_i$ is the common point of edges~$e_i^{i^\square}$ and~$e_i^{i^\square+1}$~(see \fref{mpt:fig:sweep}).

Let~$\Lambda \eqdef (\lambda_1,\dots,\lambda_n)$ be a pseudoline arrangement supported by~$\cS$. For all~$i$, we denote by~$\sigma_i^\Lambda$ the permutation of~$\{1,\dots,n\}$ whose $j$th entry~$\sigma_i^\Lambda(j)$ is the index of the pseudoline supporting~$e_i^j$, \ie such that~$\pi(e_i^j)\subset \lambda_{\sigma_i^\Lambda(j)}$. Up to reindexing the pseudolines of~$\Lambda$, we can assume that~$\sigma_0^\Lambda$ is the inverted permutation~$\sigma_0^\Lambda \eqdef (n,n-1,\dots,2,1)$, and consequently that~$\sigma_m^\Lambda$ is the identity permutation. Observe that for all~$i$:
\begin{enumerate}[(i)]
\item if~$v_i$ is a contact point of~$\Lambda$, then~$\sigma_i^\Lambda=\sigma_{i+1}^\Lambda$;
\item otherwise,~$\sigma_{i+1}^\Lambda$~is obtained from~$\sigma_i^\Lambda$ by inverting its $i^\square$th and~$(i^\square+1)$th entries.
\end{enumerate}

\begin{figure}
	\capstart
	\centerline{\includegraphics[scale=1]{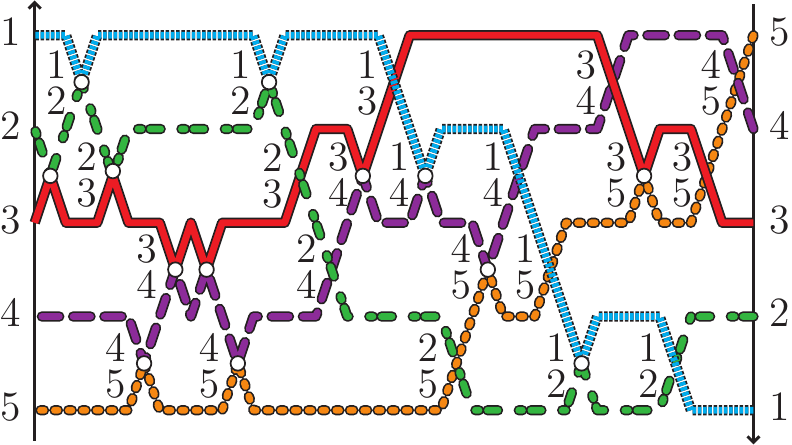}}
	\caption[The greedy pseudoline arrangement obtained by sorting]{The greedy pseudoline arrangement on the support of \fref{mpt:fig:arrangement}, obtained by sorting permutation~$(5,4,3,2,1)$ on the corresponding network. The result of each comparator is written after it (\ie to its left since we sweep backwards).}
	\label{mpt:fig:sorting}
\end{figure}

These permutations provide a characterization of the source of~$G_\chi(\cS)$ (see \fref{mpt:fig:sorting}):

\begin{theorem}\label{mpt:theo:sorting}
The unique source~$\Gamma_\chi(\cS)$ of the directed graph~$G_\chi(\cS)$ is characterized by the property that for all~$i$, the permutation~$\sigma_{i+1}^\Gamma$ is obtained from~$\sigma_i^\Gamma$ by sorting its $i^\square$th and $(i^\square+1)$th entries.
\end{theorem}

\begin{proof}
If~$\Gamma$ satisfies the above property, then it is obviously a source of the directed graph~$G_\chi(\cS)$: any flip of~$\Gamma$ is \increasing{\chi} since two of its pseudolines cannot touch before they cross.

Assume reciprocally that~$\Gamma$ is a source of~$G_\chi(\cS)$. Let~$a \eqdef \sigma_i^\Gamma(i^\square)$ and~$b \eqdef \sigma_i^\Gamma(i^\square+1)$. We have two possible situations:
\begin{enumerate}[(i)]
\item If~$a<b$, then the two pseudolines~$\lambda_a$~and~$\lambda_b$ of~$\Gamma$ already cross before~$v_i$. Consequently,~$v_i$ is necessarily a contact point of~$\Gamma$, which implies that~$\sigma_{i+1}^\Gamma(i^\square)=a$ and~$\sigma_{i+1}^\Gamma(i^\square+1)=b$.
\item If~$a>b$, then the two pseudolines~$\lambda_a$~and~$\lambda_b$ of~$\Gamma$ do not cross before~$v_i$. Since~$\Gamma$ is a source of~$G_\chi(\cS)$, $v_i$~is necessarily a crossing point of~$\Gamma$. Thus~$\sigma_{i+1}^\Gamma(i^\square)=b$ and~$\sigma_{i+1}^\Gamma(i^\square+1)=a$.
\end{enumerate}
In both cases,~$\sigma_{i+1}^\Gamma$ is obtained from~$\sigma_i^\Gamma$ by sorting its $i^\square$th and $(i^\square+1)$th entries.
\end{proof}

Let us reformulate Theorem~\ref{mpt:theo:sorting} in terms of sorting networks (see~\cite[Section~5.3.4]{k-acp-73} for a detailed presentation; see also~\cite{b-sms-74}). Let~$i<j$ be two integers. A \defn{comparator}\index{comparator|see{sorting network}}~$[i:j]$ transforms a sequence of numbers~$(a_1,\dots,a_n)$ by sorting~$(a_i,a_j)$, \ie replacing~$a_i$ by~$\min(a_i,a_j)$ and~$a_j$ by~$\max(a_i,a_j)$. A comparator~$[i:j]$ is \defn{primitive} if~$j=i+1$. A \defn{sorting network}\index{sorting network} is a sequence of comparators that sorts any sequence~$(a_1,\dots,a_n)$.

The support~$\cS$ of an arrangement of~$n$ pseudolines together with a sweep $F_m\supset\dots\supset F_0$ defines a primitive sorting network~$[1^\square:1^\square+1],\dots,[m^\square:m^\square+1]$ (see~\cite[Section 8]{k-ah-92}). Theorem~\ref{mpt:theo:sorting} affirms that sorting the permutation $(n,n-1,\dots,2,1)$ according to this sorting network provides a pseudoline arrangement supported~by~$\cS$, which depends only upon the support~$\cS$ and the filter~$F_0$ (not on the total order given by the sweep).

\subsubsection{Greedy set of contact points}\label{mpt:subsubsec:greedycontact}

In this paragraph, we give another construction of the greedy pseudoline arrangement.

Let~$\cS$ be the support of a pseudoline arrangement, and~$\chi$ be a cut of~$S$. We construct recursively a sequence~$v_1,\dots,v_q$ of vertices of~$\cS$ where for all~$i$, the vertex~$v_i$ is a minimal (for the partial order~$\cle_\chi$) remaining vertex of~$\cS$ such that~$\{v_1,\dots,v_i\}$ is a subset of the set of contact points of a pseudoline arrangement supported by~$\cS$. This procedure obviously finishes and provides a set~$V \eqdef \{v_1,\dots,v_q\}$ of vertices of~$\cS$, which turns out to be the set of contact points of the \greedy{\chi} pseudoline arrangement:

\begin{proposition}
The set~$V$ is exactly the set of contact points of~$\Gamma_\chi(\cS)$.
\end{proposition}

\begin{proof}
First of all,~$V$~is by construction the set of contact points of a pseudoline arrangement~$\Lambda$ supported by~$\cS$. If~$\Lambda$ is not the (unique) source~$\Gamma_\chi(\cS)$ of the oriented graph~$G_\chi(\cS)$, then there is a contact point~$v_i$ of~$\Lambda$ whose flip is \decreasing{\chi}. Let~$w$ denote the corresponding crossing point, and~$\Lambda'$ the pseudoline arrangement obtained from~$\Lambda$  by flipping $v_i$. This implies that~$\{v_1,\dots,v_{i-1},w\}$ is a subset of the set of contact points of~$\Lambda'$ and that~$w\cl_\chi v_i$, which contradicts the minimality of~$v_i$.
\end{proof}

\enlargethispage{.3cm}
Essentially, this proposition affirms that we obtain the same pseudoline arrangement when:
\begin{enumerate}[(i)]
\item sweeping~$\cS$ increasingly and place crossing points as long as possible; or
\item sweeping~$\cS$ decreasingly and place contact points as long as possible.
\end{enumerate}

\subsubsection{Constrained flip graph}\label{mpt:subsubsec:constrained}

Let~$\cS$ be the support of a pseudoline arrangement,~$V$~be a subset of vertices of~$\cS$, and~$\chi$ be a cut of~$\cS$. We denote by~$G(\cS\,|\,V)$ (resp.~$G_\chi(\cS\,|\,V)$) the subgraph of~$G(\cS)$ (resp.~$G_\chi(\cS)$) induced by the pseudoline arrangements with support~$\cS$ and whose set of contact points contains~$V$.

\begin{theorem}\label{mpt:theo:constrained}
The directed graph~$G_\chi(\cS\,|\,V)$ either is empty or is an acyclic connected graph with a unique source~$\Gamma_\chi(\cS\,|\,V)$ characterized by the property that for all~$i$:
\begin{enumerate}[(i)]
\item if~$v_i\in V$, then~$\sigma_{i+1}^\Gamma=\sigma_i^\Gamma$;
\item if~$v_i\notin V$, then~$\sigma_{i+1}^\Gamma$ is obtained from~$\sigma_i^\Gamma$ by sorting its $i^\square$th and $(i^\square+1)$th entries.
\end{enumerate}
\end{theorem}

\begin{proof}
We transform our support~$\cS$ into another one~$\cS'$ by opening all intersection points of~$V$ (the local picture of this transformation is \includegraphics[scale=1]{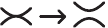}). If~$\cS'$ supports at least one pseudoline arrangement, we apply the result of Theorem~\ref{mpt:theo:sorting}: a pseudoline arrangement supported by~$\cS'$ corresponds to a pseudoline arrangement with support~$\cS$ whose set of contact points contains~$V$.
\end{proof}

In terms of sorting networks,~$\Gamma_\chi(\cS\,|\,V)$~is the result of the sorting of the inverted permutation $(n,n-1,\dots,2,1)$ by the restricted primitive network~$([i^\square:i^\square+1])_{i\in I}$, where~$I \eqdef \ens{i}{v_i\notin V}$.

Observe also that we can obtain, like in the previous subsection, the contact points of~$\Gamma_\chi(\cS\,|\, V)$ by an iterative procedure: we start from the set~$V$ and add recursively  a minimal (for the partial order~$\cle_\chi$) remaining vertex~$v_i$ of~$\cS$ such that~$V\cup\{v_1,\dots,v_i\}$ is a subset of the set of contact points of a pseudoline arrangement supported by~$\cS$. The vertex set produced by this procedure is the set of contact points of the \greedy{\chi} constrained pseudoline arrangement~$\Gamma_\chi(\cS\,|\, V)$.


\subsection{Greedy flip property and enumeration}\label{mpt:subsec:enumeration:gfp}

\subsubsection{The greedy flip property}\label{mpt:subsubsec:gfp}

\begin{figure}[!h]
	\capstart
	\centerline{\includegraphics[scale=1]{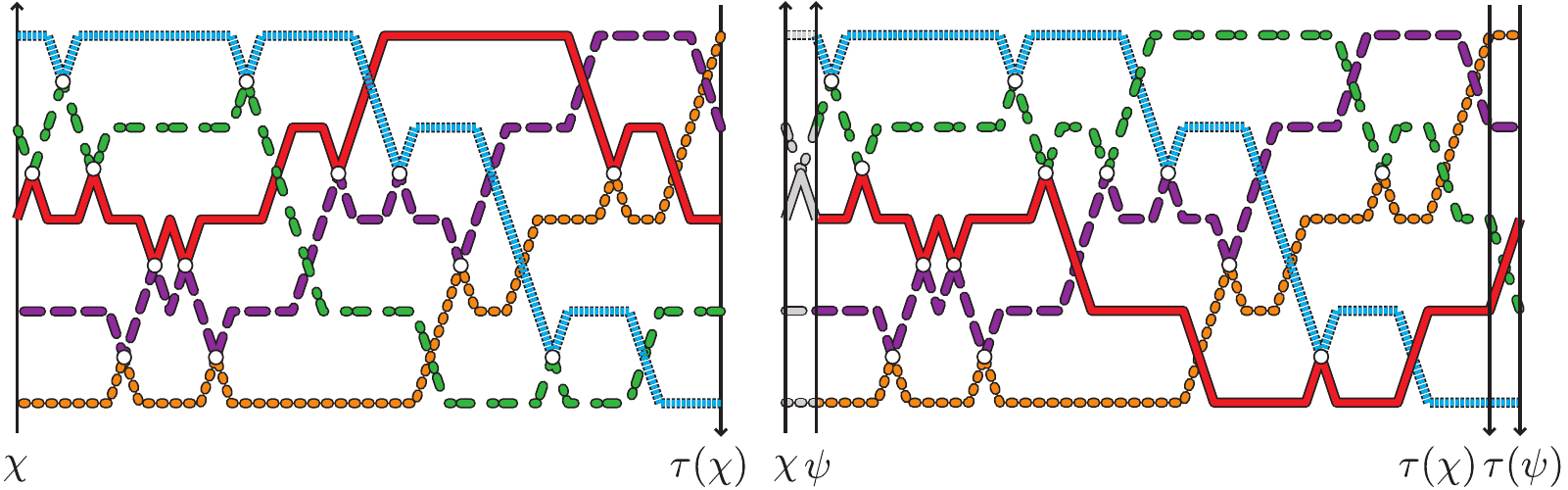}}
	\caption[The greedy flip property]{The greedy flip property.}
	\label{mpt:fig:gfp}
\end{figure}

\index{greedy!--- flip property}
We are now ready to state the \defn{greedy flip property} (see \fref{mpt:fig:gfp}) that says how to update the greedy pseudoline arrangement~$\Gamma_\chi(\cS\,|\,V)$ when either~$\chi$ or~$V$ are slightly perturbed.

\begin{theorem}[Greedy flip property]\label{mpt:theo:gfp}
Let~$\cS$ be the support of a pseudoline arrangement. Let~$\chi$ be a cut of~$\cS$, $v$~be a minimal (for the order~$\cle_\chi$) vertex of~$\cS$, and~$\psi$~denote the cut obtained from~$\chi$ by sweeping~$v$. Let~$V$ be a set of vertices of~$\cS$ (such that~$G(\cS\,|\,V)$ is not empty), and~${W \eqdef V\cup\{v\}}$. Then:
\begin{enumerate}
\item If~$v$ is a contact point of~$\Gamma_\chi(\cS\,|\,V)$ which is not in~$V$, then~$\Gamma_{\psi}(\cS\,|\,V)$ is obtained from $\Gamma_\chi(\cS\,|\,V)$ by flipping~$v$. Otherwise,~$\Gamma_{\psi}(\cS\,|\,V)=\Gamma_\chi(\cS\,|\,V)$.
\item If~$v$ is a contact point of~$\Gamma_\chi(\cS\,|\,V)$, then~$\Gamma_{\psi}(\cS\,|\,W)=\Gamma_\chi(\cS\,|\,V)$. Otherwise,~$G(\cS\,|\,W)$ is empty.
\end{enumerate}
\end{theorem}

\begin{proof}
We consider a sweep
$$F_{m+1}=F\supset F_m=F'\supset F_{m-1}\supset\dots\supset F_2\supset F_1=\tau(F)\supset F_0=\tau(F')$$
such that~$F$ (resp.~$F'$) is a filter corresponding to the cut~$\chi$ (resp.~$\psi$). Define~${\bar v_i \eqdef F_{i+1}\ssm F_i}$, ${v_i \eqdef \pi(\bar v_i)}$, and~$i^\square$ as previously.
Let~$\sigma_1,\dots,\sigma_{m+1}$ denote the sequence of permutations corresponding to~$\Gamma_\chi(\cS\,|\, V)$ on the sweep~$F_1\subset\dots\subset F_{m+1}$. In other words:
\begin{enumerate}[(i)]
\item $\sigma_1$~is the inverted permutation~$(n,n-1,\dots,2,1)$;
\item if~$v_i\in V$, then~$\sigma_{i+1}=\sigma_i$;
\item otherwise,~$\sigma_{i+1}$~is obtained from~$\sigma_i$ by sorting its $i^\square$th and $(i^\square+1)$th entries.
\end{enumerate}
Similarly, let~$\rho_0,\dots,\rho_m$ and $\omega_0,\dots,\omega_m$ denote the sequences of permutations corresponding to~$\Gamma_{\psi}(\cS\,|\ V)$ and~$\Gamma_{\psi}(\cS\,|\ W)$ respectively on the sweep~$F_0\subset\dots\subset F_m$.

Assume first that~$v$ is a contact point of~$\Gamma_\chi(\cS\,|\, V)$, but is not in~$V$. Let~$j$ denote the integer such that~$v_j$ is the crossing point of the two pseudolines of~$\Gamma_\chi(\cS\,|\, V)$ that are in contact at~$v$. We claim that in this case~$\Gamma_{\psi}(\cS\,|\, V)$ is obtained from~$\Gamma_\chi(\cS\,|\, V)$ by flipping~$v$, \ie that:
\begin{enumerate}[(i)]
\item for all~$1\le i\le j$, $\rho_i$ is obtained from~$\sigma_i$ by exchanging~$m^\square$ and~$m^\square+1$;
\item for all~$j<i\le m$,~$\rho_i=\sigma_i$.
\end{enumerate}
Indeed,~$\rho_1$ is obtained by exchanging~$m^\square$ and~$m^\square+1$ in~$\rho_0=(n,n-1,\dots,2,1)=\sigma_1$ (since~$m^\square$ and~$m^\square+1$ are the $(0^\square+1)$th and $0^\square$th entries of~$\rho_0$ respectively). Then, any comparison between two consecutive entries give the same result in~$\rho_i$ and in~$\sigma_i$, until~$m^\square$ and~$m^\square+1$ are compared again, \ie until~$i=j$. At this stage,~$m^\square$~and~$m^\square+1$ are already sorted in~$\rho_j$ but not in~$\sigma_j$. Consequently, we have to exchange~$m^\square$~and~$m^\square+1$ in~$\sigma_j$ and not in~$\rho_j$, and we obtain~${\sigma_{j+1}=\rho_{j+1}}$. After this, all the comparisons give the same result in~$\rho_i$~and~$\sigma_i$, and~$\rho_i=\sigma_i$ for all~${j<i\le m}$.

We prove similarly that:
\begin{itemize}
\item When~$v$ is not a crossing point of~$\Gamma_\chi(\cS\,|\, V)$, or is in~$V$, $\rho_i=\sigma_i$~for all~$i\in[m]$, and $\Gamma_{\psi}(\cS\,|\, V)=\Gamma_\chi(\cS\,|\, V)$.
\item When~$v$ is a contact point of~$\Gamma_\chi(\cS\,|\, V)$, $\omega_i=\sigma_i$~for all~$i\in[m]$, and ${\Gamma_{\psi}(\cS\,|\, W)=\Gamma_\chi(\cS\,|\, V)}$.
\end{itemize}

Finally, we prove that~$G(\cS\,|\, W)$ is empty when~$v$ is not a contact point of~$\Gamma_\chi(\cS\,|\, V)$. For this, assume that~$G(\cS\,|\, W)$ is not empty, and consider the greedy arrangement ${\Gamma=\Gamma_\chi(\cS\,|\, W)}$. The flip of any contact point of~$\Gamma$ not in~$W$ is \increasing{\chi}. Furthermore, since~$v$ is a minimal element for~$\cle_\chi$, the flip of~$v$ is also \increasing{\chi}. Consequently,~$\Gamma$ is a source in the graph~$G_\chi(\cS\,|\, V)$, which implies that~$\Gamma_\chi(\cS\,|\, V)=\Gamma$, and thus,~$v$ is a contact point of~$\Gamma_\chi(\cS\,|\, V)$.
\end{proof}

\subsubsection{Enumeration}\label{mpt:subsubsec:enumeration}

\index{enumeration algorithm!--- for pseudoline arrangements with contact points}
From this greedy flip property, we derive an algorithm to enumerate all pseudoline arrangements with a given support. In the next section, we will see the relation between this algorithm and the enumeration algorithm of~\cite{bkps-ceppgfa-06} for pseudotriangulations of a point~set.

Let~$\cS$ be the support of a pseudoline arrangement with a cut~$\chi$. We say that a pseudoline arrangement is \defn{colored} if its contact points are colored in blue, green or red. A green or red contact point is fixed for the end of the algorithm, while a blue one can be flipped. In the algorithm, we construct a binary tree~$\cT$, whose nodes are colored pseudoline arrangements supported by~$\cS$,~as~follows:
\begin{enumerate}[(i)]
\item The root of the tree is the \greedy{\chi} pseudoline arrangement on~$\cS$, entirely colored in blue.
\item Any node~$\Lambda$ of~$\cT$ is a leaf of~$\cT$ if either it contains a green contact point or it only contains red contact points.
\item Otherwise, we choose a minimal blue point~$v$ of~$\Lambda$. The right child of~$\Lambda$ is obtained by flipping~$v$ and coloring it in blue if the flip is \increasing{\chi} and in green if the flip is \decreasing{\chi}. The left child of~$\Lambda$ is obtained by changing the color of~$v$ into red.
\end{enumerate}

This algorithm depends upon what minimal blue point~$v$ we choose at each step. However, no matter what this choice is, we obtain all pseudoline arrangements supported by~$\cS$:

\begin{theorem}\label{mpt:theo:gfa}
The set of pseudoline arrangements supported by~$\cS$ is exactly the set of red-colored leafs of~$\cT$.
\end{theorem}

\begin{proof}
The proof is similar to that of Theorem~$9$ in~\cite{bkps-ceppgfa-06}.

We define inductively a cut~$\chi(\Lambda)$ for each node~$\Lambda$ of~$\cT$: the cut of the root is~$\chi$, and for each node~$\Lambda$ the cut of its children is obtained from~$\chi$ by sweeping the contact point~$v$. We also denote~$V(\Lambda)$ the set of red contact points of~$\Lambda$. With these notations, the greedy flip property (Theorem~\ref{mpt:theo:gfp}) ensures that~$\Lambda=\Gamma_{\chi(\Lambda)}(\cS\,|\,V(\Lambda))$, for each node~$\Lambda$ of~$\cT$.

The fact that any red-colored leaf of~$\cT$ is a pseudoline arrangement supported by~$\cS$ is obvious. Reciprocally, let us prove that any pseudoline arrangement supported by~$\cS$ is a red leaf of~$\cT$. Let~$\Lambda$ be a pseudoline arrangement supported by~$\cS$. We define inductively a path~$\Lambda_0,\dots,\Lambda_p$ in the tree~$\cT$ as follows:~$\Lambda_0$~is the root of~$\cT$ and for all~$i\ge 0$,~$\Lambda_{i+1}$~is the left child of~$\Lambda_i$ if the minimal blue contact point of~$\Lambda_i$ is a contact point of~$\Lambda$, and its right child otherwise (we stop when we reach a leaf). We claim that for all~$0\le i\le p$:
\begin{itemize}
\item the set~$V(\Lambda_i)$ is a subset of contact points of~$\Lambda$;
\item the contact points of~$\Lambda$ not in~$V(\Lambda_i)$ are not located between~$\chi(\Lambda)$~and~$\chi$;
\end{itemize}
from which we derive that~$\Lambda=\Lambda_p$ is a red-colored leaf.
\end{proof}

Let us briefly discuss the complexity of this algorithm. We assume that the input of the algorithm is a pseudoline arrangement and we consider a flip as an elementary operation. Then, this algorithm requires a polynomial running time per pseudoline arrangement supported by~$\cS$. As for many enumeration algorithms, the crucial point of this algorithm is that its working space is also polynomial (while the number of pseudoline arrangements supported by~$\cS$ is exponential).


\section{Dual pseudoline arrangements}\label{mpt:sec:duality}

In this section, we prove that both the graph of flips on ``(pointed) pseudotriangulations of a point set'' and the graph of flips on ``multitriangulations of a convex polygon'' can be interpreted as graphs of flips on ``pseudoline arrangements with a given support''. This interpretation is based on the classical duality that we briefly recall in the first subsection, and yields to a natural definition of ``\mpt{}s of a pseudoline arrangement'' that we present in Section~\ref{mpt:sec:mpt}.


\subsection{The dual pseudoline arrangement of a point set}\label{mpt:subsec:duality:points}

To a given oriented line in the Euclidean plane, we usually associate its angle~$\theta\in\R/2\pi\Z$ with the horizontal axis and its algebraic distance~$d\in\R$ to the origin (\ie the value of~$\dotprod{(-u,v)}{.}$ on the line, where~$(u,v)$ is its unitary direction vector). Since the same line oriented in the other direction gives an angle~$\theta+\pi$ and a distance~$-d$, this parametrization naturally associates a point of the M\"obius strip to each line of the Euclidean plane. In other words, the line space of the Euclidean plane is (isomorphic to) the M\"obius strip.

\begin{figure}[!h]
	\capstart
	\centerline{\includegraphics[scale=1]{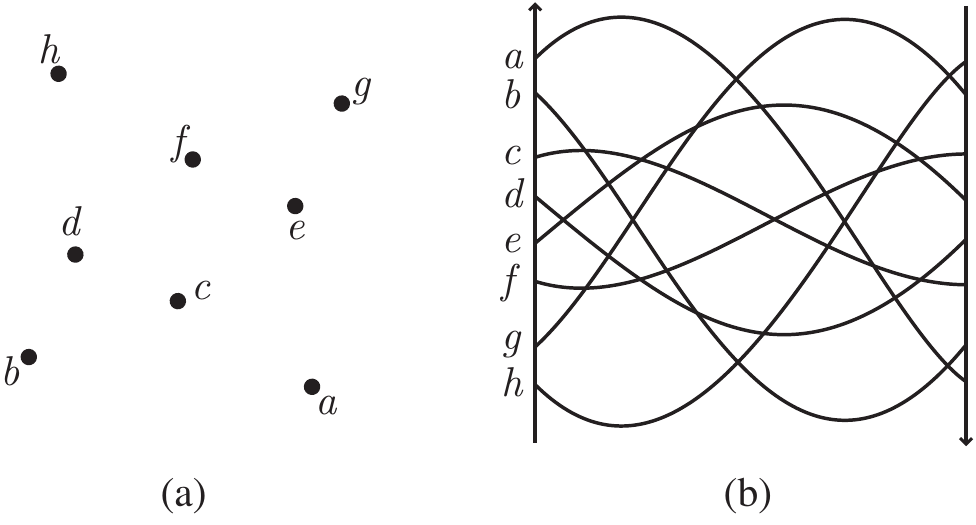}}
	\caption[Dual arrangement of a point set]{(a)~A point set~$P$ in general position and (b)~its dual pseudoline arrangement~$P^*$.}
	\label{mpt:fig:dual}
\end{figure}

\index{dual!--- of a point|hbf}
\index{dual!--- of a point set|hbf}
Via this parametrization, the set of lines passing through a point~$p$ forms a pseudoline~$p^*$. The pseudolines~$p^*$~and~$q^*$ dual to two distinct points~$p$~and~$q$ have a unique crossing point dual to the line~$(pq)$. Thus, for a finite point set~$P$ in the Euclidean plane, the set~$P^* \eqdef \ens{p^*}{p\in P}$ is a pseudoline arrangement without contact points (see \fref{mpt:fig:dual}). Again, we always assume that the point set~$P$ is in general position (no three points lie in a same line), so that the arrangement~$P^*$ is simple.

\index{plane!topological ---}
\index{topological plane}
This elementary duality also holds for any topological plane (or \mbox{$\R^2$-plane}, see~\cite{sbghls-cpp-95} for a definition), not only for the real Euclidean plane~$\R^2$. That is to say, the line space of a topological plane is (isomorphic to) the M\"obius strip and the dual of a finite set of points in a topological plane is a pseudoline arrangement without contact points. Let us also recall that any pseudoline arrangement of the M\"obius strip without contact points is the dual arrangement of a finite set of points in a certain topological plane~\cite{gpwz-atp-94}. Thus, in the rest of this chapter, we deal with sets of points and their duals without restriction to the Euclidean plane.


\subsection{The dual pseudoline arrangement of a pseudotriangulation}\label{mpt:subsec:duality:pt}

Introduced by Michel Pocchiola and Gert Vegter in their study of the visibility complex of a set of disjoint convex obstacles in the plane~\cite{pv-tsvcp-96,pv-vc-96}, pseudotriangulations were then used in different contexts such as motion planning and rigidity theory~\cite{s-ptrmp-05,horsssssw-pmrgpt-05}. Their combinatorial and geometric structure has been extensively studied in the last decade (number of pseudotriangulations~\cite{aaks-cmpt-04,aoss-nptcps-08}, polytope of pseudotriangulations~\cite{rss-empppt-03}, algorithmic issues~\cite{b-eptp-05,bkps-ceppgfa-06,hp-cpcpb-07}, \etc). We refer to~\cite{rss-pt-06} for a detailed survey on the subject, and just recall here the elementary definitions we need:

\begin{definition}\label{mpt:def:pseudotriangulation}
\index{pseudotriangle|hbf}
\index{corner}
\index{tangent|see{pseudotriangle}}
\index{pseudotriangulation|hbf}
A \defn{pseudotriangle} is a polygon~$\Delta$ with only three convex vertices (the \defn{corners} of~$\Delta$), joined by three concave polygonal chains (\fref{mpt:fig:pseudotriangle}). A line is said to be \defn{tangent} to~$\Delta$ if:
\begin{enumerate}[(i)]
\item either it passes through a corner of~$\Delta$ and separates the two edges incident to it;
\item or it passes through a concave vertex of~$\Delta$ and does not separate the two edges incident~to~it.
\end{enumerate}

\begin{figure}
	\capstart
	\centerline{\includegraphics[scale=1.3]{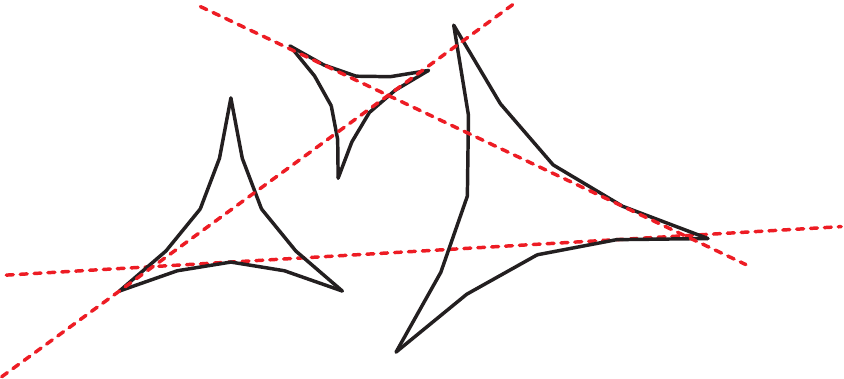}}
	\caption[Three pseudotriangles and their common tangents]{Three pseudotriangles and their common tangents.}
	\label{mpt:fig:pseudotriangle}
\end{figure}

A \defn{pseudotriangulation} of a point set~$P$ in general position is a set of edges of~$P$ which decomposes the convex hull of~$P$ into pseudotriangles.
\end{definition}

\index{pseudotriangulation!pointed ---|hbf}
\index{pointed|see{pseudotriangulation}}
In this dissertation, we only deal with \defn{pointed} pseudotriangulations, meaning that for any point~$p\in P$, there exists a line which passes through~$p$ and defines a half-plane containing all the edges incident to~$p$. Under this assumption, any two pseudotriangles of a pseudotriangulation have a unique common tangent. This leads to the following observation (first stated in~\cite{pv-ot-94,pv-ptta-96} in the context of pseudotriangulations of convex bodies):

\begin{observation}
\index{dual!--- of a pseudotriangle|hbf}
\index{dual!--- of a pseudotriangulation|hbf}
Let~$T$ be a pseudotriangulation of a point set~$P$ in general position.~Then:
\begin{enumerate}[(i)]
\item the set~$\Delta^*$ of all tangents to a pseudotriangle~$\Delta$ of~$T$ is a pseudoline in the M\"obius strip;
\item the set~$T^* \eqdef \ens{\Delta^*}{\Delta\text{ pseudotriangle of } T}$ of dual pseudolines of the pseudotriangles of~$T$ is a pseudoline arrangement (with contact points);~and
\item $T^*$ is supported by~$P^*$ minus its first level (see \fref{mpt:fig:pseudotriangulationpoints}(b)).
\end{enumerate}
\begin{figure}
	\capstart
	\centerline{\includegraphics[scale=1]{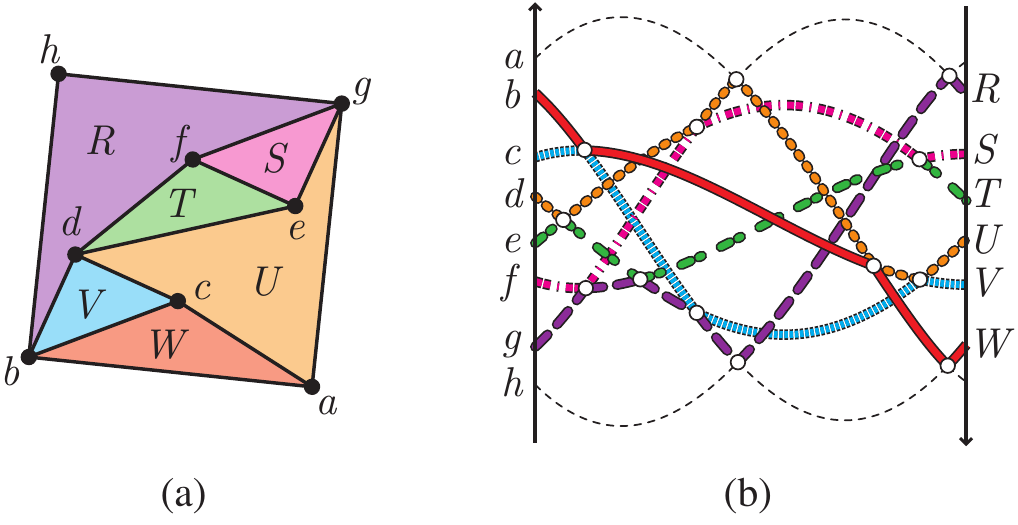}}
	\caption[A pseudotriangulation of a point set and its dual pseudoline arrangement]{A pseudotriangulation~(a) of the point set of \fref{mpt:fig:dual}(a) and its dual arrangement~(b). Each thick pseudoline in (b) corresponds to a pseudotriangle in (a); each contact point in (b) corresponds to an edge in (a); \etc.}
	\label{mpt:fig:pseudotriangulationpoints}
\end{figure}
\end{observation}

It turns out that this covering property characterizes pseudotriangulations:

\begin{theorem}\label{mpt:theo:dualitypt}
Let~$P$ be a finite point set (in general position), and~$P^{*1}$ denote the support of its dual pseudoline arrangement minus its first level. Then:
\begin{enumerate}[(i)]
\item The dual pseudoline arrangement~$T^* \eqdef \ens{\Delta^*}{\Delta\text{ pseudotriangle of }T}$ of a pseudotriangulation~$T$ of~$P$ is supported by~$P^{*1}$.
\item The primal set of edges~$E \eqdef \ens{[p,q]}{p,q\in P,\; p^*\wedge q^*\text{ is not a crossing point of } \Lambda}$ of a pseudoline arrangement~$\Lambda$ supported by~$P^{*1}$ is a pseudotriangulation of~$P$.
\end{enumerate}
\end{theorem}

\begin{figure}[h]
	\capstart
	\centerline{\includegraphics[scale=1]{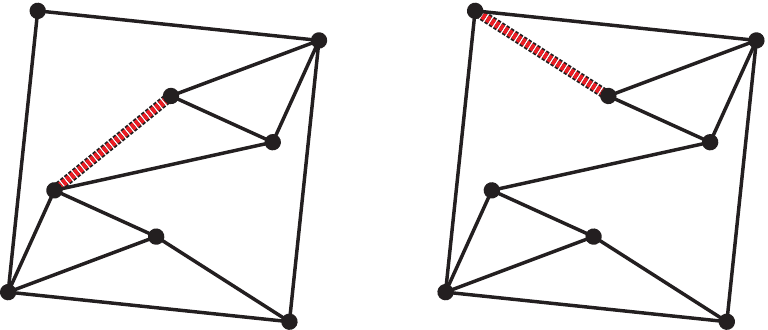}}
	\caption[A flip in a pseudotriangulation]{A flip in the pseudotriangulation of \fref{mpt:fig:pseudotriangulationpoints}(a).}
	\label{mpt:fig:primalflip}
\end{figure}

In this section, we provide two proofs of Part~(ii) of this result. The first one is based on flips. First, remember that there is also a simple flip operation on pseudotriangulations of~$P$: replacing any internal edge~$e$ in a pseudotriangulation of~$P$ by the common tangent of the two pseudotriangles containing~$e$ produces a new pseudotriangulation of~$P$. For example, \fref{mpt:fig:primalflip} shows two pseudotriangulations of the point set of \fref{mpt:fig:dual}(a), related by a flip. We denote by~$G(P)$ the graph of flips on pseudotriangulations of~$P$, whose vertices are pseudotriangulations of~$P$ and whose edges are flips between them. In other words, there is an edge in~$G(P)$ between two pseudotriangulations of~$P$ if and only if their symmetric difference is reduced to a pair.

\begin{proof}[Proof~$1$ of Theorem~\ref{mpt:theo:dualitypt}(ii)]
The two notions of flips (the primal notion on pseudotriangulations of~$P$ and the dual notion on pseudoline arrangements supported by~$P^{*1}$) coincide via duality: an internal edge~$e$ of a pseudotriangulation~$T$ of~$P$ corresponds to a contact point~$e^*$ of the dual pseudoline arrangement~$T^*$; the two pseudotriangles~$\Delta$~and~$\Delta'$ of~$T$ containing~$e$ correspond to the two pseudolines~$\Delta^*$~and~$\Delta'^*$ of~$T^*$ in contact at~$e^*$; and the common tangent~$f$ of~$\Delta$~and~$\Delta'$ corresponds to the crossing point~$f^*$ of~$\Delta^*$~and~$\Delta'^*$.

Thus, the graph~$G(P)$ is a subgraph of~$G(\cS)$. Since both are connected and regular of degree~$|P|-3$, they coincide. In particular, any pseudoline arrangement supported by~$P^{*1}$ is the dual of a pseudotriangulation of~$P$.
\end{proof}

\begin{remark}
\index{greedy!--- pseudotriangulation}
\index{enumeration algorithm!--- for pseudotriangulations}
Observe that this duality matches our greedy pseudoline arrangement supported by~$P^{*1}$ with the greedy pseudotriangulation of~\cite{bkps-ceppgfa-06} (originally of~\cite{pv-tsvcp-96} for convex bodies). In particular, the greedy flip property and the enumeration algorithm of Subsection~\ref{mpt:subsec:enumeration:gfp} are generalizations of results in~\cite{bkps-ceppgfa-06,pv-tsvcp-96}. It is interesting to see however that our new proof of the greedy flip property (with the interpretation as sorting networks) does not use the notion of visibility complex, which is essential in~\cite{bkps-ceppgfa-06,pv-tsvcp-96}.
\end{remark}

\index{witness pseudoline|hbf}
Our second proof of Theorem~\ref{mpt:theo:dualitypt} is a bit more complicated, but is more direct and introduces a ``witness method'' that we will repeatedly use throughout this chapter. It is based on the following characterization of pseudotriangulations:

\begin{lemma}[\cite{s-ptrmp-05}]\label{mpt:lem:streinu}
A graph~$T$ on a point set~$P$ is a pseudotriangulation of~$P$ if and only if it is non-crossing, pointed and has~$2|P|-3$ edges.\qed
\end{lemma}

\begin{proof}[Proof~$2$ of Theorem~\ref{mpt:theo:dualitypt}(ii)]
We check that~$E$ is non-crossing, pointed and has~$2|P|-3$ edges:

\vspace{-.4cm}
\paragraph{Cardinality.} First, the number of edges of~$E$ equals the difference between the number of crossing points of~$P^*$ and of~$\Lambda$:
$$|E|={\left|P^*\right| \choose 2}-{| \Lambda| \choose 2}={|P| \choose 2}-{|P|-2 \choose 2}=2|P|-3.$$

\begin{figure}
	\capstart
	\centerline{\includegraphics[scale=1]{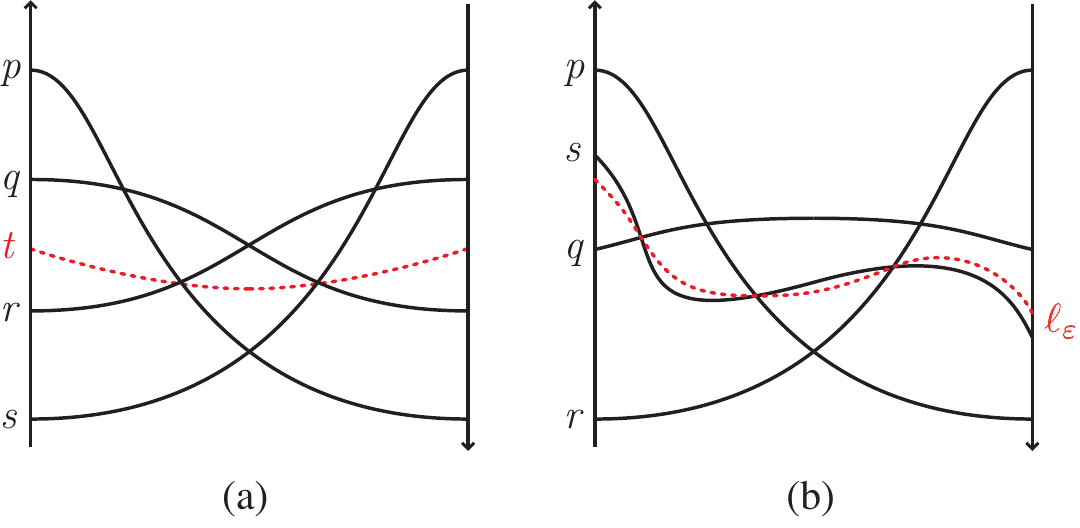}}
	\caption[Dual arrangements of forbidden configurations]{(a) Four points~$p,q,r,s$ in convex position with the intersection~$t$ of~$[p,r]$ and~$[q,s]$; (b)~A point~$s$ inside a triangle~$pqr$ with the witness pseudoline~$\ell_\varepsilon$.}
	\label{mpt:fig:pqrs}
\end{figure}

\vspace{-.6cm}
\paragraph{Crossing-free.} Let~$p,q,r,s$ be four points of~$P$ in convex position. Let~$t$ be the intersection of~$[p,r]$~and~$[q,s]$ (see \fref{mpt:fig:pqrs}(a)). We use the pseudoline~$t^*$ as a \defn{witness} to prove that~$[p,r]$ and~$[q,s]$ cannot both be in~$E$. For this, we count crossings of~$t^*$ with~$P^*$ and~$\Lambda$ respectively:
\begin{enumerate}[(i)]
\item Since~$P$ is in general position, $t$~is not in~$P$ and~$P^*\cup\{t^*\}=(P\cup\{t\})^*$ is a (non-simple) pseudoline arrangement. Thus,~$t^*$ crosses~$P^*$ exactly~$|P|$ times.
\item Since~$t^*$ is a pseudoline, it crosses each pseudoline of~$\Lambda$ at least once. Thus, it crosses~$\Lambda$ at least~$|\Lambda|=|P|-2$ times.
\item For each of the points~$p^*\wedge r^*$~and~$q^*\wedge s^*$, replacing the crossing point by a contact point removes two crossings with~$t^*$.
\end{enumerate}
Thus,~$[p,r]$~and~$[q,s]$ cannot both be in~$E$, and~$E$ is crossing-free.

\paragraph{Pointed.} Let~$p,q,r,s$ be four points of~$P$ such that~$s$ lies inside the triangle~$pqr$. We first construct a witness pseudoline (see \fref{mpt:fig:pqrs}(b)) that we use to prove that~$[p,s]$,~$[q,s]$~and~$[r,s]$ cannot all be in~$E$. Let~$f_p,f_q,f_r$~and~$f_s$ represent~$p^*,q^*,r^*$~and~$s^*$ respectively. Let~$x,y,z\in\R$ be such that~$f_p(x)=f_s(x)$, $f_q(y)=f_s(y)$~and~$f_r(z)=f_s(z)$. Let~$g$ be a continuous and \piantiperiodic{} function vanishing exactly on~$\{x,y,z\}+\Z\pi$ and changing sign each time it vanishes; say for example $g(t) \eqdef \sin(t-x)\sin(t-y)\sin(t-z)$. For all~$\varepsilon> 0$, we define the function $h_\varepsilon:\R\to\R$ by~$h_\varepsilon(t)=f_s(t)+\varepsilon g(t)$. The function~$h_\varepsilon$ is continuous and \piantiperiodic. The corresponding pseudoline~$\ell_\varepsilon$ crosses~$s^*$ three times. It is also easy to see that if~$\varepsilon$ is sufficently small, then~$\ell_\varepsilon$ crosses the pseudolines of~$(P\ssm\{s\})^*$ exactly as~$s^*$ does (see \fref{mpt:fig:pqrs}(b)). For such a small~$\varepsilon$, we count the crossings of~$\ell_\varepsilon$ with~$P^*$~and~$\Lambda$ respectively:
\begin{enumerate}[(i)]
\item $\ell_\varepsilon$~crosses~$P^*$ exactly~$|P|+2$ times (it crosses~$s^*$ three times and any other pseudoline of~$P^*$ exactly once).
\item Since~$\ell_\varepsilon$ is a pseudoline, it crosses~$\Lambda$ at least~$|\Lambda|=|P|-2$ times.
\item For each of the points~$p^*\wedge s^*$, $q^*\wedge s^*$~and~$r^*\wedge s^*$, replacing the crossing point by a contact point removes two crossings with~$\ell_\varepsilon$.
\end{enumerate}
Thus,~$[p,r]$,~$[q,s]$~and~$[r,s]$ cannot all be in~$E$, and~$E$ is pointed.
\end{proof}


\subsection{The dual pseudoline arrangement of a multitriangulation}\label{mpt:subsec:duality:mt}

Our study of stars (see Chapter~\ref{chap:stars}) yields to a similar observation for multitriangulations:

\begin{observation}
\index{dual!--- of a \kstar{k}|hbf}
\index{dual!--- of a \ktri{k}|hbf}
Let~$T$ be a \ktri{k} of a convex \gon{n}. Then:
\begin{enumerate}[(i)]
\item the set~$S^*$ of all bisectors of a \kstar{k}~$S$ of~$T$ is a pseudoline of the M\"obius strip;
\item the set~$T^* \eqdef \ens{S^*}{S\;k\text{-star of } T}$ of dual pseudolines of \kstar{k}s of~$T$ is a pseudoline arrangement (with contact points);~and
\item $T^*$ is supported by the dual pseudoline arrangement~$V_n^*$ of $V_n$ minus its first $k$~levels (see \fref{mpt:fig:2triang8points}(b)).
\end{enumerate}
\begin{figure}
	\capstart
	\centerline{\includegraphics[scale=1]{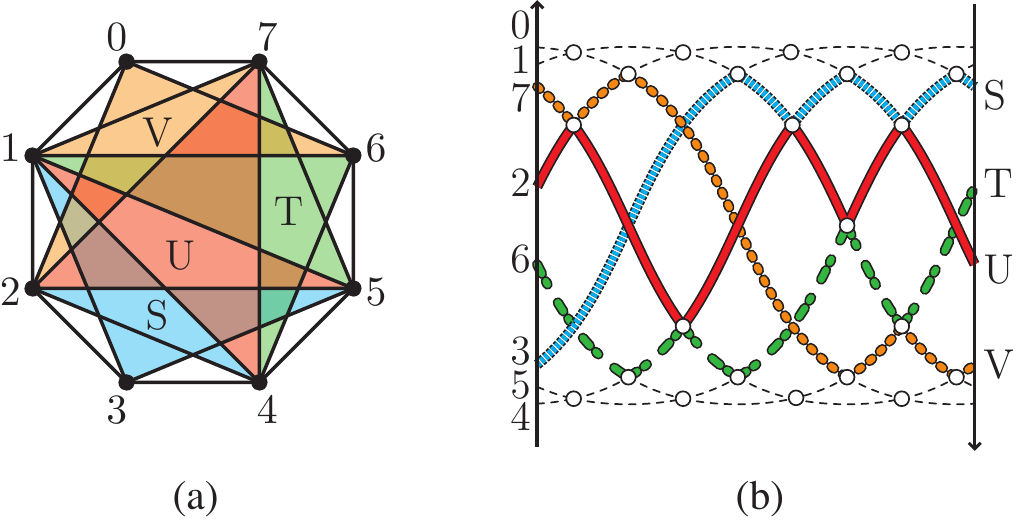}}
	\caption[A \ktri{2} of the octagon and its dual pseudoline arrangement]{A \ktri{2} of the octagon~(a) and its dual arrangement~(b). Each thick pseudoline in~(b) corresponds to a \kstar{2} in~(a); each contact point in~(b) corresponds to an edge in~(a);~\etc}
	\label{mpt:fig:2triang8points}
\end{figure}
\end{observation}

Again, it turns out that this observation provides a characterization of multitriangulations of a convex polygon (and we present again two different proofs of this result):

\begin{theorem}\label{mpt:theo:dualitymt}
Let~$V_n$ denote the set of vertices of a convex \gon{n}, and~$V_n^{*k}$ denote the dual pseudoline arrangement of~$V_n$ minus its first~$k$~levels. Then:
\begin{enumerate}[(i)]
\item The dual pseudoline arran\-gement~$T^* \eqdef \ens{S^*}{S\;k\text{-star of }T}$ of a \ktri{k}~$T$ of the \gon{n} is supported by~$V_n^{*k}$.
\item The primal set of edges~$E \eqdef \ens{[p,q]}{p,q\in V_n,\; p^*\wedge q^*\text{ is not a crossing point of } \Lambda}$ of a pseudoline arrangement~$\Lambda$ supported by~$V_n^{*k}$ is a \ktri{k} of the \gon{n}.
\end{enumerate}
\end{theorem}

\begin{proof}[Proof~$1$ of Theorem~\ref{mpt:theo:dualitymt}(ii)]
The two notions of flips (the primal notion on \ktri{k}s of the \gon{n} and the dual notion on pseudoline arrangements supported by~$V_n^{*k}$) coincide. Thus, the flip graph~$G_{n,k}$ on \ktri{k}s of the \gon{n} is a subgraph of~$G(V_n^{*k})$. Since they are both connected and regular of degree~$k(n-2k-1)$, these two graphs coincide. In particular, any pseudoline arrangement supported by~$V_n^{*k}$ is the dual of a \ktri{k} of the \gon{n}.
\end{proof}

\begin{proof}[Proof~$2$ of Theorem~\ref{mpt:theo:dualitymt}(ii)]
\index{witness pseudoline}
We follow the method of our second proof of Theorem~\ref{mpt:theo:dualitypt}(ii). Since~$E$ has the right number of edges (namely~$k(2n-2k-1)$), we only have to prove that it is \kcross{(k+1)}-free. We consider~$2k+2$ points~$p_0,\dots, p_k, q_0,\dots, q_k$ cyclically ordered around the unit circle. Since the definition of crossing (and thus, of \kcross{\ell}) is purely combinatorial, \ie depends only on the cyclic order of the points and not on their exact positions, we can move all the vertices of our \gon{n} on the unit circle while preserving their cyclic order. In particular, we can assume that the lines~$(p_iq_i)_{i\in\{0,\dots,k\}}$ all contain a common point~$t$. Its dual pseudoline~$t^*$ crosses~$V_n^*$ exactly~$n$ times and~$\Lambda$ at least~$|\Lambda|=n-2k$ times. Furthermore, for any point~${p_i^*\wedge q_i^*}$, replacing the crossing point by a contact point removes two crossings with~$t^*$. Thus, the pseudoline~$t^*$ provides a witness which proves that the edges~$[p_i,q_i]$,~$i\in\{0,\dots,k\}$, cannot be all in~$E$, and thus ensures that~$E$ is \kcross{(k+1)}-free.
\end{proof}

As an application of Theorem~\ref{mpt:theo:dualitymt}, we provide the promised proof of the Characterization Theorem~\ref{stars:theo:characterization} mentioned in Chapter~\ref{chap:stars}:

\begin{theorem}\label{mpt:theo:characterization}
Let~$\Sigma$ be a set of \kstar{k}s of the \gon{n} such that:
\begin{enumerate}[(i)]
\item any \krel{k} edge of~$E_n$ is contained in zero or two \kstar{k}s of~$\Sigma$, one on~each~side;~and
\item any \kbound{k} edge of~$E_n$ is contained in exactly one \kstar{k} of~$\Sigma$.
\end{enumerate}
Then~$\Sigma$ is the set of \kstar{k}s of a \ktri{k} of the \gon{n}.
\end{theorem}

\begin{proof}
By duality,~$\Sigma$ corresponds to a set of pseudolines~$\Lambda$ of the M\"obius strip. We just have to prove that~$\Lambda$ forms a pseudoline arrangement supported by~$V_n^{*k}$ (the dual arrangement of~$n$ pseudolines in convex position, except its first~$k$ levels).

Observe first that the pseudolines of~$\Lambda$ partition~$V_n^{*k}$ (meaning that each edge of~$V_n^{*k}$ is contained in precisely one pseudoline of~$\Lambda$). Indeed:
\begin{enumerate}
\item Consider the set~$A$ of all pseudosegments of a pseudoline~$\ell\in L$, which are located in~$V_n^{*k}$. The number of pseudolines of~$\Lambda$ containing a pseudosegment of~$A$ is constant on~$A$. Otherwise, when this number changes, the contact point would not satisfy condition (i).
\item By condition (ii), this number is~$1$ for each pseudoline.
\end{enumerate}

This immediately implies that~$|\Sigma|=n-2k$. Consequently (double counting) we know that we have~$kn$ non-\krel{k} edges and~$k(n-2k-1)$ \krel{k} ones. This only leaves ${n \choose 2}-nk-k(n-2k-1)={n-2k \choose 2}$ crossing points. Since two pseudolines of the M\"obius strip cross at least once, this implies that~$\Lambda$ is a pseudoline arrangement, and finishes the proof.
\end{proof}


\section{Multipseudotriangulations}\label{mpt:sec:mpt}

Motivated by Theorems~\ref{mpt:theo:dualitypt} and~\ref{mpt:theo:dualitymt}, we define in terms of pseudoline arrangements a natural generalization of both pseudotriangulations and multitriangulations. We then study elementary properties of the corresponding set of edges in the primal space.


\subsection{Definition}\label{mpt:subsec:mpt:definition}

\index{levels (of an arrangement)|hbf}
Let~$L$ be a pseudoline arrangement supported by~$\cS$. Define its \defn{\kkernel{k}}\index{kernel@\kkernel{k} (of an arrangement)}~$\cS^k$ to be its support minus its first~$k$ levels (which are the iterated external hulls of the support of the arrangement). Denote by~$V^k$ the set of contact points of~$L$ in~$\cS^k$.

\begin{definition}
\index{pseudotriangulation@\pt{k}!--- of a pseudoline arrangement|hbf}
\index{multipseudotriangulation!--- of a pseudoline arrangement|hbf}
A \defn{\pt{k}} of~$L$ is a pseudoline arrangement whose support is~$\cS^k$ and whose set of contact points contains~$V^k$. 
\end{definition}

Pseudotriangulations of a point set~$P$ correspond via duality to \pt{1}s of the dual pseudoline arrangement~$P^*$. Similarly, \ktri{k}s of the \gon{n} correspond to \pt{k}s of the pseudoline arrangement~$V_n^*$ in convex position.
If~$L$ is a pseudoline arrangement with no contact point, then any pseudoline arrangement supported by~$\cS^k$ is a \pt{k} of~$L$. In general, the condition that the contact points of~$L$ in its \kkernel{k} should be contact points of any \pt{k} of~$L$ is a natural assumption for iterating \mpt{}s (see Section~\ref{mpt:sec:iterated}).

Let~$\Lambda$ be a \pt{k} of~$L$. We denote by~$V(\Lambda)$ the union of the set of contact points of~$\Lambda$ with the set of intersection points of the first~$k$ levels of~$L$. In other words,~$V(\Lambda)$ is the set of intersection points of~$L$ which are not crossing points of~$\Lambda$. As for pseudoline arrangements, the set~$V(\Lambda)$ completely determines~$\Lambda$.

Flips for \mpt{}s are defined as in Lemma~\ref{mpt:lem:flip}, with the restriction that the contact points in~$V^k$ cannot be flipped. In other words, the flip graph on \pt{k}s of~$L$ is exactly the graph~$G(\cS^k\,|\, V^k)$. Section~\ref{mpt:sec:enumeration} asserts that the graph of flips is regular and connected, and provides an enumeration algorithm for \mpt{}s~of~$L$.

\begin{figure}
	\capstart
	\centerline{\includegraphics[scale=1]{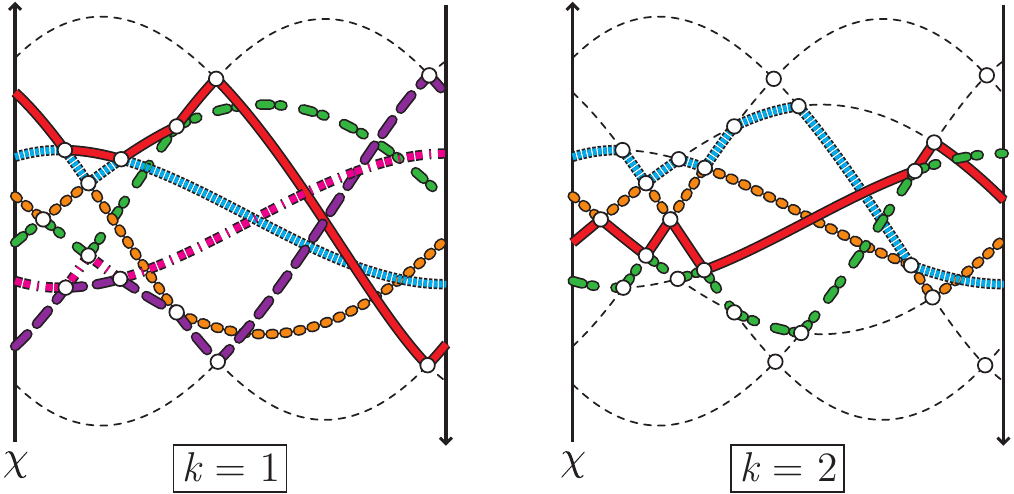}}
	\caption[Greedy \mpt{}s]{The \greedy{\chi} \pt{1} and the \greedy{\chi} \pt{2} of the pseudoline arrangement of \fref{mpt:fig:dual}(b).}
	\label{mpt:fig:greedy}
\end{figure}

Let~$\chi$ be a cut of (the support of)~$L$. It is also a cut of the \kkernel{k}~$\cS^k$ of~$L$. A particularly interesting example of \pt{k} of~$L$ is the source of the graph of \increasing{\chi} flips on \pt{k}s~of~$L$ (see \fref{mpt:fig:greedy} for an illustration):

\begin{definition}
\index{greedy!--- multipseudotriangulation}
\index{multipseudotriangulation!greedy ---}
\index{pseudotriangulation@\pt{k}!greedy ---}
The \defn{\greedy{\chi} \pt{k}} of~$L$, denoted~$\Gamma_\chi^k(L)$, is the greedy pseudoline arrangement~$\Gamma_\chi(\cS^k\,|\, V^k)$.
\end{definition}

\subsection{Pointedness and crossings}\label{mpt:subsec:mpt:pointedcrossing}

Let~$P$ be a point set in general position. Let~$\Lambda$ be a \pt{k} of~$P^*$ and~$V(\Lambda)$ be the set of crossing points of~$P^*$ which are not crossing points of~$\Lambda$. We call \defn{primal of~$\Lambda$} the set~$E \eqdef \ens{[p,q]}{p,q\in P,\; p^*\wedge q^*\in V(\Lambda)}$ of edges of~$P$ primal to~$V(\Lambda)$~---~see \fref{mpt:fig:2-pseudotriangulation}. Here, we discuss general properties of primals of \mpt{}s. We start with elementary properties that we already observed for the special cases of pseudotriangulations and multitriangulations in the proofs of Theorems~\ref{mpt:theo:dualitypt} and~\ref{mpt:theo:dualitymt}:

\begin{figure}[h]
	\capstart
	\centerline{\includegraphics[scale=1]{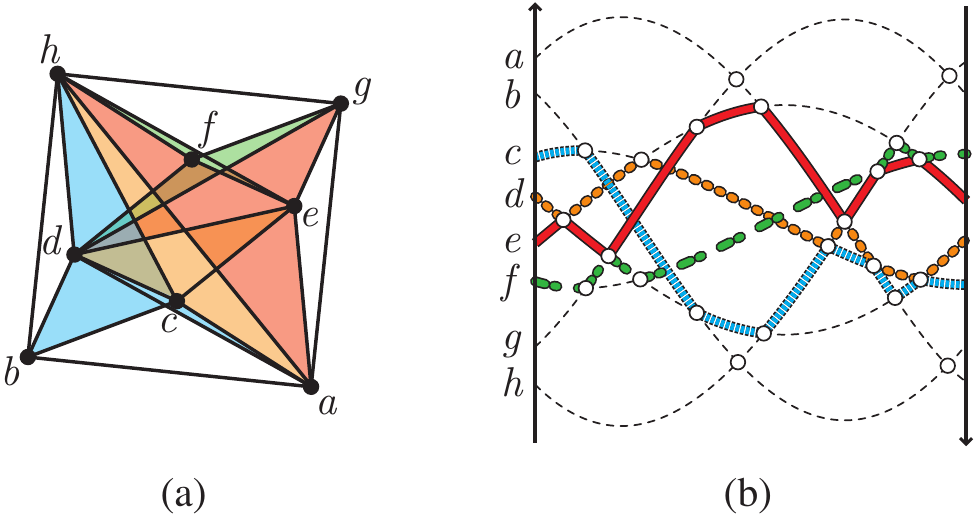}}
	\caption[The primal of a \pt{2}]{The primal set of edges~(a) of a \pt{2} (b)~of the dual pseudoline of the point set of~\fref{mpt:fig:dual}(a).}
	\label{mpt:fig:2-pseudotriangulation}
\end{figure}

\begin{lemma}
The set~$E$ has~$k(2|P|-2k-1)$ edges.
\end{lemma}

\proof
The number of edges of~$E$ is the difference between the number of crossing points in the pseudoline arrangements~$P^*$ and~$\Lambda$:
$$\hspace{2.4cm}|E|={\left|P^*\right| \choose 2}-{|\Lambda| \choose 2}={|P| \choose 2}-{|P|-2k \choose 2}=k(2|P|-2k-1).\hspace{1.8cm}\qed$$

We now discuss pointedness of~$E$. We call \defn{\kalter{k}}\index{alternation@\kalter{k}|hbf} any set~$\ens{f_i}{i\in\Z_{2k+1}}$ of~$2k+1$ edges all adjacent to a common vertex and whose cyclic order around it is given by $f_0\cl \bar f_{1+k} \cl f_1 \cl \bar f_{2+k} \cl \dots \cl f_{2k} \cl \bar f_k \cl f_0$, where~$\bar f_i$ denotes the opposite direction of the edge~$f_i$.

\begin{lemma}\label{mpt:lem:pointed}
The set~$E$ cannot contain a \kalter{k}.
\end{lemma}

\begin{proof}
We simply mimic the proof of pointedness in Theorem~\ref{mpt:theo:dualitypt}. Let~$p_0,\dots,p_{2k}$ and~$q$ be~${2k+2}$ points of~$P$ such that~$F \eqdef \ens{[p_i,q]}{i\in\Z_{2k+1}}$ is a \kalter{k}. We prove that~$F$ cannot be a subset of~$E$ by constructing a witness pseudoline~$\ell$ that separates all the crossing points~${p_i^*\wedge q^*}$ corresponding to~$F$, while crossing~$q^*$ exactly~$2k+1$ times and the other pseudolines of~$P^*$ exactly as~$q^*$ does. (We skip the precise construction, since it is exactly the same as in the proof of Theorem~\ref{mpt:theo:dualitypt}.) Counting the crossings of~$\ell$ with~$P^*$ and~$\Lambda$, we obtain:
\begin{enumerate}[(i)]
\item $\ell$~crosses~$P^*$ exactly~$|P|+2k$ times;
\item $\ell$~crosses~$\Lambda$ at least~$|\Lambda|=|P|-2k$ times;
\item for each of the points~$p_i^*\wedge q^*$, replacing the crossing point by a contact point removes two crossings with~$\ell$.\qedhere
\end{enumerate}
\end{proof}

\begin{figure}
	\capstart
	\centerline{\includegraphics[scale=1]{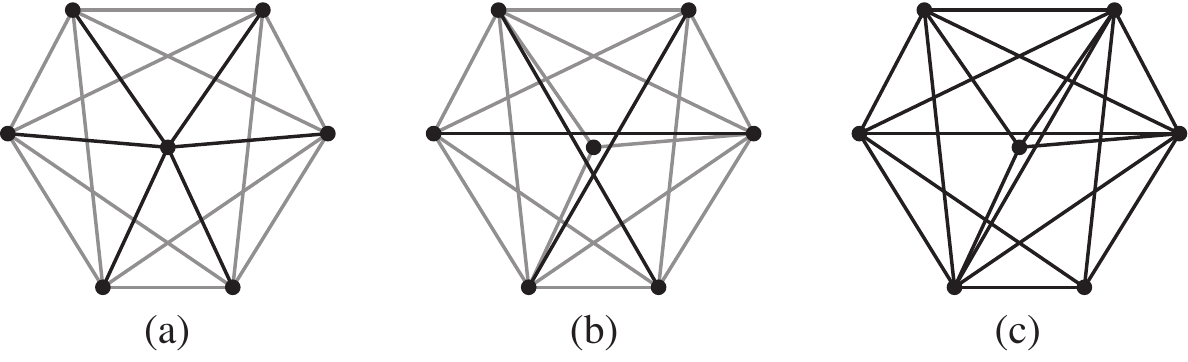}}
	\caption[A non-\kpointed{2} \pt{2}, a \pt{2} containing a \kcross{3} and a \kcross{3}-free \kalter{2}-free set not contained in a \pt{2}]{(a) A non-\kpointed{2} (but \kalter{2}-free) \pt{2}; (b)~A \pt{2} containing a \kcross{3}; (c) A \kcross{3}-free \kalter{2}-free set not contained in a \pt{2}.}
	\label{mpt:fig:3crossing}
\end{figure}

\begin{remark}
Observe that a set of edges is pointed if and only if it is \kalter{1}-free.
In contrast, we want to observe the difference between \kalter{k}-freeness and the following natural notion of \kpointed{k}ness: we say that a set~$F$ of edges with vertices in~$P$ is \defn{\kpointed{k}} if for all~$p$ in~$P$, there exists a line which passes through~$p$ and defines a half-plane that contains at most~$k-1$ segments of~$F$ adjacent to~$p$.
Observe that a \kpointed{k} set is automatically \kalter{k}-free but that the reciprocal statement does not hold (see \fref{mpt:fig:3crossing}(a)).
\end{remark}

Finally, contrarily to pseudotriangulations ($k=1$) and multitriangulations (convex position), the condition of avoiding \kcross{(k+1)}s does not hold for \pt{k}s in general:

\enlargethispage{.1cm}
\begin{remark}
There exist \pt{k}s with \kcross{(k+1)}s (see \fref{mpt:fig:3crossing}(b)) as well as \kcross{(k+1)}-free \kalter{k}-free sets of edges that are not subsets of \pt{k}s (see \fref{mpt:fig:3crossing}(c)).
\end{remark}

\subsection{Stars in \mpt{}s}\label{mpt:subsec:mpt:stars}

To complete our understanding of the primal of \mpt{}s, we need to generalize pseudotriangles of pseudotriangulations and \kstar{k}s of \ktri{k}s: both pseudotriangles and \kstar{k}s correspond to pseudolines of the covering pseudoline arrangement.

We keep the notations of the previous section:~$P$~is a point set in general position, $\Lambda$~is a \pt{k} of~$P^*$ and~$E$~is the primal set of edges of~$\Lambda$. Consider now a pseudoline~$\lambda$ of~$\Lambda$. We call \defn{star}\index{star|hbf} the set of edges ${S(\lambda) \eqdef \ens{[p,q]}{p,q\in P,\; p^*\wedge q^*\text{ contact point of }\lambda}}$ primal to the contact points of~$\lambda$.

\begin{lemma}
For any~$\lambda\in \Lambda$, the star~$S(\lambda)$ is non-empty.
\end{lemma}

\begin{proof}
We have to prove that any pseudoline~$\lambda$ of~$\Lambda$ supports at least one contact point. If it is not the case, then~$\lambda$ is also a pseudoline of~$P^*$, and all the~$|P|-1$ crossing points of~$\lambda$ with~$P^*\ssm\{\lambda\}$ should be crossing points of~$\lambda$ with~$\Lambda\ssm\{\lambda\}$. This is impossible since~$|\Lambda\ssm\{\lambda\}|=|P|-2k-1$.
\end{proof}

Similarly to the case of \ktri{k}s of the \gon{n}, we say that an edge~$[p,q]$ of~$E$ is a \defn{\krel{k}}\index{relevant@\krel{k}!--- edge} (resp.~\defn{\kbound{k}}\index{boundary@\kbound{k} edge}, resp.~\defn{\kirrel{k}}\index{irrelevant@\kirrel{k} edge}) edge if there remain strictly more than (resp.~exactly, resp.~strictly less than)~$k-1$~points of~$P$ on each side (resp. one side) of the line~$(pq)$. In other words,~$p^*\wedge q^*$ is located in the \kkernel{k} (resp.~in the intersection of the $k$th level and the \kkernel{k}, resp.~in the first $k$ levels) of the pseudoline arrangement~$P^*$. Thus, the edge~$[p,q]$ is contained in~$2$ (resp.~$1$, resp.~$0$) stars of~$\Lambda$.

\begin{figure}
	\capstart
	\centerline{\includegraphics[scale=1]{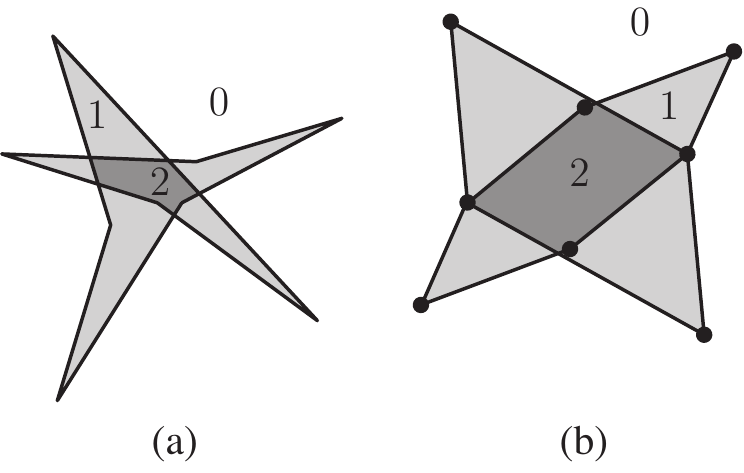}}
	\caption[Winding number of a star; depth in a point set]{(a) Winding number of a star; (b)~\kdepth{2} in the point set of \fref{mpt:fig:dual}(a).}
	\label{mpt:fig:depth_winding_number}
\end{figure}

The edges of a star~$S(\lambda)$ are cyclically ordered by the order of their dual contact points on~$\lambda$, and thus~$S(\lambda)$ forms a (not-necessarily simple) polygonal cycle. For any point~$q$ in the plane, let~$\sigma_\lambda(q)$ denote the \defn{winding number}\index{winding number (of a star)|hbf} of~$S(\lambda)$ around~$q$, that is, the number of rounds made by~$S(\lambda)$ around the point~$q$ (see \fref{mpt:fig:depth_winding_number}(a)). For example, the winding number of a point in the external face is~$0$.

We call \defn{\kdepth{k}}\index{depth@\kdepth{k}|hbf} of a point~$q$ the number~$\delta^k(q)$ of \kbound{k} edges of~$P$ crossed by any (generic continuous) path from~$q$ to the external face, counted positively when passing from the ``big'' side (the one containing at least~$k$ vertices of~$P$) to the ``small side'' (the one containing~$k-1$ vertices of~$P$), and negatively otherwise (see \fref{mpt:fig:depth_winding_number}(b)). That this number is independent from the path can be seen by mutation. For example,~$\delta^1(q)$ is~$1$ if~$q$ is in the convex hull of~$P$ and~$0$ otherwise.

\begin{proposition}\label{mpt:prop:decomposition}
Any point~$q$ of the plane is covered~$\delta^k(q)$ times by the stars~$S(\lambda)$, $\lambda\in \Lambda$, of the \pt{k}~$\Lambda$ of~$P^*$:
$$ \delta^k(q)=\sum_{\lambda\in\Lambda} \sigma_\ell(q).$$
\end{proposition}

The proposition is intuitively clear: let us walk on a continuous path from the external face to the point~$q$. Initially, the winding numbers of all stars of~$\Lambda$ are zero (we start outside all stars of~$\Lambda$). Then, each time we cross an edge~$e$:
\begin{enumerate}[(i)]
\item If~$e$ is \kirrel{k}, it is not contained in any star of~$\Lambda$, and we do not change the winding numbers of the stars of~$\Lambda$.
\item If~$e$ is a \kbound{k} edge, and if we cross it positively, we increase the winding number of the star~$S$ of~$\Lambda$ containing~$e$; if we cross~$e$ negatively, we decrease the winding~number~of~$S$.
\item If~$e$ is \krel{k}, then we decrease the winding number of one star of~$\Lambda$ containing~$e$ and increase the winding number of the other star of~$\Lambda$ containing~$e$.
\end{enumerate}
Let us give a formal proof in the dual:

\begin{proof}[Proof of Proposition~\ref{mpt:prop:decomposition}]
Both~$\sigma_\lambda(q)$ and~$\delta^k(q)$ can be read on the pseudoline~$q^*$:
\begin{enumerate}[(i)]
\item If~$\tau_\lambda(q)$ denotes the number of intersection points between~$q^*$ and~$\lambda$ (that is, the number of tangents to~$S(\lambda)$ passing through~$q$), then~$\sigma_\lambda(q)=(\tau_\lambda(q)-1)/2$.
\item If~$\gamma^k(q)$ denotes the number of intersection points between~$q^*$ and the first~$k$ levels of~$P^*$, then~$\delta^k(q)=k-\gamma^k(q)/2$.
\end{enumerate}
The pseudoline~$q^*$ has exactly~$|P|$ crossings with~$P^*$ (since~$P^*\cup\{q^*\}$ is an arrangement), which are crossings either with the pseudolines of~$\Lambda$ or with the first~$k$ levels of~$P^*$. Hence,
$$|P|=\gamma^k(q)+\sum_{\lambda\in \Lambda} \tau_\lambda(q)=2k-2\delta^k(q)+|\Lambda|+2\sum_{\lambda\in \Lambda}\sigma_\lambda(q),$$
and we get the aforementioned result since~$|\Lambda|=|P|-2k$.
\end{proof}

A \defn{corner}\index{corner|hbf} of the star~$S(\lambda)$ is an internal convex angle of it (see \fref{mpt:fig:star}(a)). In the following proposition, we are interested in the number of corners of~$S(\lambda)$:

\begin{figure}
	\capstart
	\centerline{\includegraphics[scale=1]{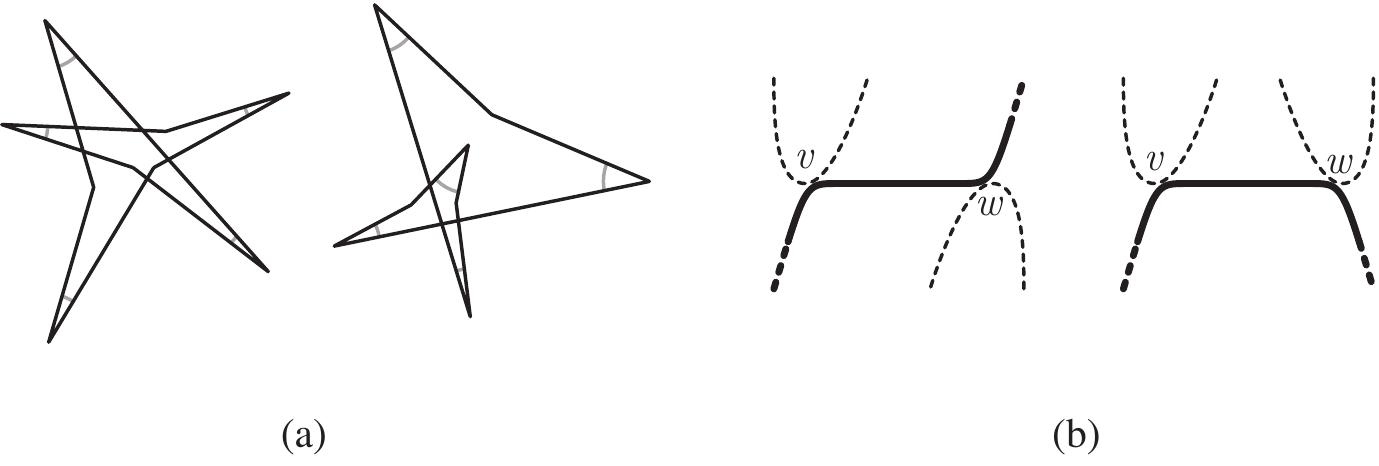}}
	\caption[Two stars with~$5$ corners; The two possible configurations of two consecutive contact points on~$\lambda$]{(a) Two stars with~$5$ corners; (b)~The two possible configurations of two consecutive contact points on~$\lambda$: convex (left) and concave (right).}
	\label{mpt:fig:star}
\end{figure}

\begin{proposition}\label{mpt:prop:corners}
The number of corners of a star~$S(\lambda)$ of a \pt{k} of~$P^*$ is odd and between~$2k+1$~and $2(k-1)|P|+2k+1$.
\end{proposition}

\begin{proof}
We read convexity of internal angles of~$S(\lambda)$ on the preimage~$\bar\lambda$ of the pseudoline~$\lambda$ under the projection~$\pi$. Let~$pqr$ be an internal angle, let~$v=p^*\wedge q^*$~and~$w=q^*\wedge r^*$ denote the contact points corresponding to the two edges~$[p,q]$~and~$[q,r]$ of this angle, and let~$\bar v$ and~$\bar w$ denote two consecutive preimages of~$v$~and~$w$ on~$\bar\lambda$ (meaning that~$\bar w$ is located between~$\bar v$~and~$\tau(\bar v)$). The angle~$pqr$ is a corner if and only if~$\bar v$~and~$\bar w$ lie on opposite sides of~$\bar\lambda$, meaning that the other curves touching~$\bar\lambda$ at~$\bar v$~and~$\bar w$ lie on opposite sides, one above and one below~$\bar\lambda$ (see \fref{mpt:fig:star}(b)).

\svs
In particular, the number~$c(\lambda)=c$ of corners of~$S(\lambda)$ is the number of opposite consecutive contact points on~$\bar\lambda$ between two versions~$\bar v$~and~$\tau(\bar v)$ of a contact point~$v$ of~$\lambda$. To see that~$c$ is odd, imagine that we are discovering the contact points of~$\lambda$ one by one. The first contact point~$v$ that we see corresponds to two opposite contact points~$\bar v$~and~$\tau(\bar v)$ on~$\bar\lambda$. Then, at each stage, we insert a new contact point~$\bar w$ between two old contact points that can be:
\begin{enumerate}[(i)]
\item either on opposite sides and then we are not changing~$c$;
\item or on the same side and we are adding to~$c$ either~$0$ (if~$\bar w$ is also on the same side) or~$2$ (if~$\bar w$ is on the opposite side).
\end{enumerate}
Thus,~$c$~remains odd in any case.

\medskip

To prove the lower bound, we use our witness method. We perturb a little bit~$\lambda$ to obtain a pseudoline~$\mu$ that passes on the opposite side of each contact point (this is possible since~$c$ is odd). This pseudoline~$\mu$ crosses~$\lambda$ between each pair of opposite contact points and crosses the other pseudolines of~$\Lambda$ exactly as~$\lambda$ does. Thus,~$\mu$ crosses~$\Lambda$ exactly~$|\Lambda|-1+c$ times. But since~$\mu$ is a pseudoline, it has to cross all the pseudolines of~$P^*$ at least once. Thus, we get~$|P| \le |\Lambda|-1+c=|P|-2k-1+c$ and~$c\ge 2k+1$.

From this lower bound, we derive automatically the upper bound. Indeed, we know that the number of corners around one point~$p$ is at most~$\deg(p)-1$. Consequently,
$$\sum_{p\in P} (\deg(p)-1) \ge \sum_{\nu\in \Lambda} c(\nu) = c(\lambda)+\sum_{\substack{\nu\in \Lambda\\ \nu\ne\lambda}} c(\nu).$$
The left sum equals~$2k(2|P|-2k-1)-|P|$ while, according to the previous lower bound, the right one is at least~$c+(|P|-2k-1)(2k+1)$. Thus we get~$c\le 2(k-1)|P|+2k+1$.
\end{proof}

\begin{figure}
	\capstart
	\centerline{\includegraphics[scale=1]{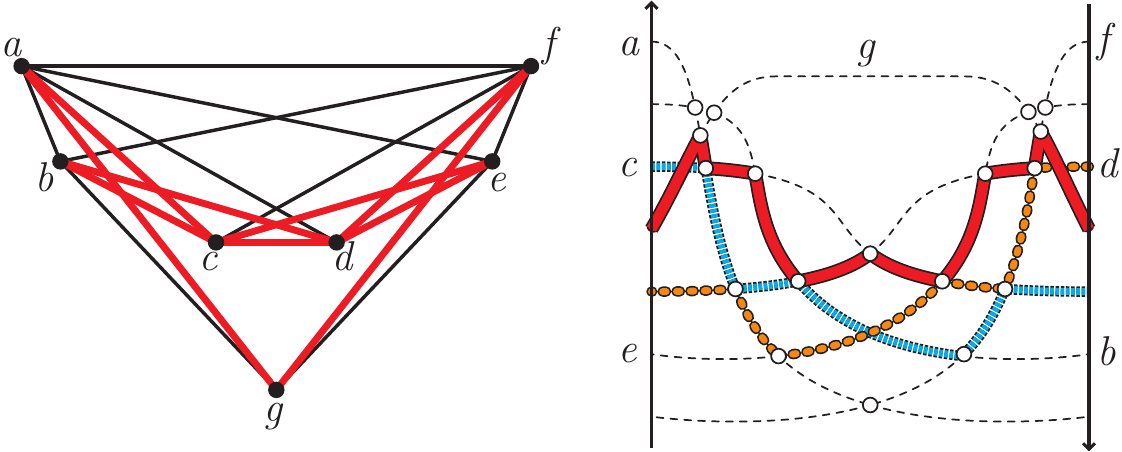}}
	\caption[A star of a \pt{2} with almost~$2|P|$ corners]{A star of a \pt{2} with almost~$2|P|$ corners.}
	\label{mpt:fig:2nstar}
\end{figure}

For general~$k$, contrary to the case of convex position (Section~\ref{mpt:subsec:duality:mt}), the number of corners of a star is not necessarily~$2k+1$ (\fref{mpt:fig:2nstar} presents a star of a \pt{2} with almost~$2|P|$ corners).

When~$k=1$ however, both bounds turn out to be~$3$ and this proposition affirms that the primal of a pseudoline of~$\Lambda$ is a pseudotriangle. Proposition~\ref{mpt:prop:decomposition} ensures that the convex hull of~$P$ is covered once by these~$|P|-2$ pseudotriangles, and hence provides another alternative proof of Theorem~\ref{mpt:theo:dualitypt}.

Finally, let us mention that the definitions developed in this section provide some alternative conditions in the Characterization Theorem~\ref{mpt:theo:characterization}:

\begin{lemma}\label{mpt:lem:boundary}
Let~$\Sigma$ be a set of \kstar{k}s of the convex \gon{n} such that each \krel{k} edge of~$E_n$ is contained either in zero or in two \kstar{k}s of~$\Sigma$, one on each side. Then the following properties are equivalent:
\begin{enumerate}
\item Every \kbound{k} edge of the \gon{n} is contained in exactly one \kstar{k} of~$\Sigma$.
\item On each of the~$\gcd(n,k)$ cycles of \kbound{k} edges of the \gon{n}, there is one \kbound{k} edge which is contained in exactly one \kstar{k} of~$\Sigma$.
\item The \kdepth{k} of every chamber of the \gon{n} equals the sum of the winding numbers of the \kstar{k}s of~$\Sigma$ around it.
\item For one chamber of the \gon{n} of \kdepth{k}~$k$, the sum of the winding numbers of the \kstar{k}s of~$\Sigma$ around it equals~$k$.
\item For every (generic) direction, each line parallel to this direction passing through one of the~$n-2k$ central vertices of the \gon{n} (for this direction) is the bisector of a unique \kstar{k} of~$\Sigma$.
\item There is one (generic) direction for which each line parallel to this direction passing through one of the~$n-2k$ central vertices of the \gon{n} (for this direction) is the bisector of a unique \kstar{k} of~$\Sigma$.\qed
\end{enumerate}
\end{lemma}


\section{Iterated \mpt{}s}\label{mpt:sec:iterated}

\index{multipseudotriangulation!iterated ---|hbf}
By definition, a \pt{k} of an \pt{m} of a pseudoline arrangement~$L$ is a \pt{(k+m)} of~$L$. In this section, we study these iterated sequences of \mpt{}s. In particular, we compare \mpt{}s with iterated sequences of \pt{1}s.


\subsection{Definition and examples}\label{mpt:subsec:iterated:definition}

Let~$L$ be a pseudoline arrangement. An \defn{iterated \mpt{}} of~$L$ is a sequence $\Lambda_1,\dots,\Lambda_r$ of pseudoline arrangements such that~$\Lambda_i$ is a \mpt of~$\Lambda_{i-1}$ for all~$i$ (by convention,~$\Lambda_0=L$). We call \defn{signature} of~$\Lambda_1,\dots,\Lambda_r$ the sequence~$k_1<\cdots<k_r$ of integers such that~$\Lambda_i$ is a \pt{k_i} of~$L$ for all~$i$. Observe that the assumption that contact points of a pseudoline arrangement~$L$ should be contact points of any \mpt of~$L$ is natural in this setting: iterated \mpt{}s correspond to decreasing sequences of sets of crossing points.

A \defn{decomposition} of a \mpt~$\Lambda$ of a pseudoline arrangement~$L$ is an iterated \mpt~$\Lambda_1,\dots,\Lambda_r$ of~$L$ such that~$\Lambda_r=\Lambda$ and~$r>1$. We say that~$\Lambda$ is \defn{decomposable} if such a decomposition exists, and \defn{irreducible}\index{multipseudotriangulation!irreducible ---} otherwise. The decomposition is \defn{complete} if its signature is~$1,2,\dots,r$.

\begin{example}[\kcolorable{k} \ktri{k}s of the convex \gon{n}]
Consider a \kcolorable{k} \ktri{k}~$T$ of the \gon{n} (see Section~\ref{stars:sec:ears}). It is easy to construct from the $k$~\kaccordion{k}s of~$T$ a deomposition~$\Lambda_1,\dots,\Lambda_k$ of~$T^*$ of~$V_n^*$, where each \mpt{}~$\Lambda_i$ is in convex position (\ie the dual pseudoline arrangement of points in convex position). In particular, the set of decomposable \ktri{k}s of the convex \gon{n} is much bigger than the set of \kcolorable{k} \ktri{k}s. However, as the following example shows, there are still irreducible multitriangulations.
\end{example}

\begin{example}[An irreducible \ktri{2} of the \gon{15}]\label{mpt:exm:irred}
We consider the geometric graph~$T$ of \fref{mpt:fig:15gonctrexm}. The edges are:
\begin{enumerate}[(i)]
\item all the \kirrel{2} and \kbound{2} edges of the \gon{15}, and
\item the five zigzags~$Z_a=\{[3a,3a+6], [3a+6,3a+1], [3a+1,3a+5], [3a+5,3a+2]\}$, for~$a\in\{0,1,2,3,4\}$.
\end{enumerate}
Thus,~$T$ has~$50$ edges and is \kcross{3}-free (since the only \krel{2} edges of~$T$ that cross a zigzag~$Z_a$ are edges of~$Z_{a-1}$ and~$Z_{a+1}$). Consequently,~$T$ is a \ktri{2} of the \gon{15}.

\begin{figure}
	\capstart
	\centerline{\includegraphics[scale=1]{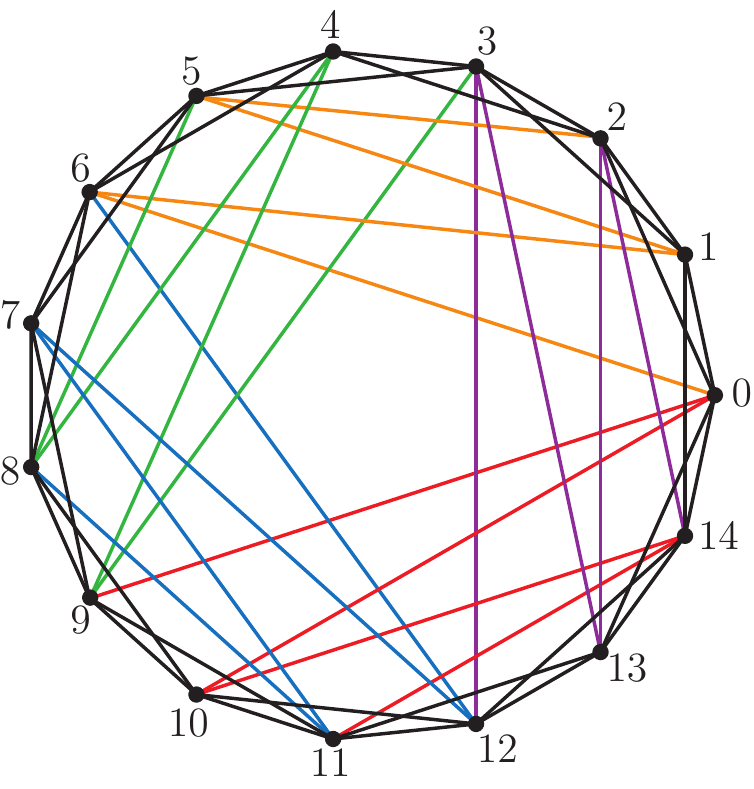}}
	\caption[An irreducible \ktri{2} of the \gon{15}]{An irreducible \ktri{2} of the \gon{15}: it contains no triangulation.}
	\label{mpt:fig:15gonctrexm}
\end{figure}

Let us now prove that~$T^*$ is irreducible, that is, that~$T$ contains no triangulation. Observe first that the edge~$[0,6]$ cannot be an edge of a triangulation contained in~$T$ since none of the triangles~$06i$,~${i\in\{7,\dots,14\}}$, is contained in~$T$. Thus, we are looking for a triangulation contained in~$T\ssm\{[0,6]\}$. Repeating the argument successively for the edges~$[1,6]$,~$[1,5]$~and~$[2,5]$, we prove that the zigzag~$Z_0$ is disjoint from any triangulation contained in~$T$. By symmetry, this prove the irreducibility of~$T^*$.
\end{example}


\subsection{Iterated greedy pseudotriangulations}\label{mpt:subsec:iterated:greedy}

Greedy \mpt{}s provide interesting examples of iteration of pseudotriangulations. Let~$L$ be a pseudoline arrangement, and~$\chi$ be a cut of~$L$. 

\begin{theorem}
For any positive integers~$a$ and~$b$, $\Gamma_\chi^{a+b}(L)=\Gamma_\chi^b(\Gamma_\chi^a(L))$. Consequently, for any integer~$k$, $\Gamma_\chi^k(L)=\Gamma_\chi^1\circ\Gamma_\chi^1\circ\cdots\circ\Gamma_\chi^1(L)$, where~$\Gamma_\chi^1(.)$ is iterated~$k$ times.
\end{theorem}

\begin{proof}
Since~$\chi$ is a cut of~$L$, it is also a cut of~$\Gamma^a_\chi(L)$ and thus~$\Gamma^b_\chi(\Gamma^a_\chi(L))$ is well-defined. Observe also that we can assume that~$L$ has no contact point (otherwise, we can open them). Let~$n \eqdef |L|$ and~$m \eqdef {n \choose 2}$.

Let~$\chi=\chi_0,\dots,\chi_m=\chi$ be a backward sweep of~$L$. For all~$i$, let~$v_i$ denote the vertex of~$L$ swept when passing from~$\chi_i$ to $\chi_{i+1}$, and~$i^\square$ denote the integer such that the pseudolines that cross at~$v_i$ are the $i^\square$th and $(i^\square+1)$th pseudolines of~$L$ on~$\chi_i$.

Let~$\sigma_0,\dots,\sigma_m$ denote the sequence of permutations corresponding to~$\Gamma^a_\chi(L)$ on the sweep $\chi_0,\dots,\chi_m$. In other words,~$\sigma_0$ is the permutation
$$(1,\dots,a,n-a,n-a-1,\dots,a+2,a+1,n-a+1,\dots,n),$$
whose first~$a$ and last~$a$ entries are preserved, while its~$n-2a$ intermediate entries are inverted. Then, for all~$i$, the permutation~$\sigma_{i+1}$ is obtained from~$\sigma_i$ by sorting its $i^\square$th and $(i^\square+1)$th entries.

Similarly, let~$\rho_0,\dots,\rho_m$ and~$\omega_0,\dots,\omega_m$ denote the sequences of permutations corresponding to~$\Gamma_\chi^{a+b}(L)$ and~$\Gamma^b_\chi(\Gamma^a_\chi(L))$ respectively: both~$\rho_0$ and~$\omega_0$ equal the permutation whose first and last~$a+b$ entries are preserved and whose~$n-2a-2b$ intermediate entries are inverted,~and~for~all~$i$:
\begin{itemize}
\item $\rho_{i+1}$ is obtained from~$\rho_i$ by sorting its $i^\square$th and $(i^\square+1)$th entries;
\item if~$v_i\notin\Gamma^a_\chi(L)$, then~$\omega_{i+1}$ is obtained from~$\omega_i$ by sorting its $i^\square$th and $(i^\square+1)$th entries; otherwise,~$\omega_{i+1}=\omega_i$.
\end{itemize}

We claim that for all~$i$,
\begin{enumerate}[(A)]
\item all the inversions of~$\rho_i$ are also inversions of~$\sigma_i$: $\rho_i(p)>\rho_i(q)$ implies~$\sigma_i(p)>\sigma_i(q)$ for all~$1\le p<q\le n$; and
\item $\rho_i=\omega_i$.
\end{enumerate}

We prove this claim by induction on~$i$. It is clear for~$i=0$. Assume that it is true for~$i$ and let us prove it for~$i+1$. We have two possible situations:

\begin{enumerate}
\item \textbf{First case:}~$\sigma_i(i^\square)<\sigma_i(i^\square+1)$. Then,~$\sigma_{i+1}=\sigma_i$ and~$v_i\in\Gamma^a_\chi(L)$. Thus,~${\omega_{i+1}=\omega_i}$. Furthermore, using Property~(A) at rank~$i$, we know that~$\rho_i(i^\square)<\rho_i(i^\square+1)$, and thus ${\rho_{i+1}=\rho_i}$. To summarize, $\sigma_{i+1}=\sigma_i$, $\omega_{i+1}=\omega_i$, and $\rho_{i+1}=\rho_i$, which trivially implies that Properties (A) and (B) remain true.

\item \textbf{Second case:}~$\sigma_i(i^\square)>\sigma_i(i^\square+1)$. Then,~$\sigma_{i+1}$ is obtained from~$\sigma_i$ by exchanging the $i^\square$th and $(i^\square+1)$th entries, and~$v_i\notin\Gamma^a_\chi(L)$. Consequently, $\rho_{i+1}$~and~$\omega_{i+1}$ are both obtained from~$\rho_i$~and~$\omega_i$ respectively by sorting their $i^\square$th and $(i^\square+1)$th entries. Thus, Property~(B) obviously remains true. As far as Property~(A) is concerned, the result is obvious if~$p$~and~$q$ are different from~$i^\square$~and~$i^\square+1$. By symmetry, it is enough to prove that for any~$p<i^\square$, $\rho_{i+1}(p)>\rho_{i+1}(i^\square)$~implies~$\sigma_{i+1}(p)>\sigma_{i+1}(i^\square)$. We have to consider two subcases:
\begin{enumerate}
\item \textbf{First subcase:}~$\rho_i(i^\square)<\rho_i(i^\square+1)$. Then~$\rho_{i+1}=\rho_i$. Thus, if~$p<i^\square$ is such that~$\rho_{i+1}(p)>\rho_{i+1}(i^\square)$, then we have~$\rho_i(p)>\rho_i(i^\square)$. Consequently, we obtain that ${\sigma_{i+1}(p)=\sigma_i(p)>\sigma_i(i^\square)>\sigma_i(i^\square+1)=\sigma_{i+1}(i^\square)}$.
\item \textbf{Second subcase:}~$\rho_i(i^\square)>\rho_i(i^\square+1)$. Then~$\rho_{i+1}$ is obtained from~$\rho_i$ by exchanging its $i^\square$th and $(i^\square+1)$th entries. If~$p<i^\square$ is such that~$\rho_{i+1}(p)>\rho_{i+1}(i^\square)$, then ${\rho_i(p)>\rho_i(i^\square+1)}$. Consequently,~$\sigma_{i+1}(p)=\sigma_i(p)>\sigma_i(i^\square+1)=\sigma_{i+1}(i^\square)$.
\end{enumerate}
\end{enumerate}

Obviously, Property~(B) of our claim proves the theorem.
\end{proof}


\subsection{Flips in iterated \mpt{}s}\label{mpt:subsec:iterated:flips}

Let~$\Lambda_1,\dots,\Lambda_r$  be an iterated \mpt of a pseudoline arrangement~$L$, with signature~$k_1<\dots<k_r$.
Let~$v$ be a contact point of~$\Lambda_r$ (which is not a contact point of~$L$), and let~$i$ denote the first integer for which~$v$ is a contact point of~$\Lambda_i$ (thus,~$v$ is a contact point of~$\Lambda_j$ if and only if~$i\le j\le r$). For all~$i\le j\le r$, let~$\Lambda'_j$ denote the pseudoline arrangement obtained from~$\Lambda_j$ by flipping~$v$, and let~$w_j$ denote the new contact point of~$\Lambda'_j$. Let~$j$ denote the biggest integer such that~$w_j=w_i$.
There are three possibilities:
\begin{enumerate}[(i)]
\item If~$j=r$, then~$\Lambda_1,\dots,\Lambda_{i-1},\Lambda'_i,\dots,\Lambda'_r$ is an iterated \mpt of~$L$. We say that it is obtained from~$\Lambda_1,\dots,\Lambda_r$ by a \defn{complete flip} of~$v$ .
\item If~$j<r$, and~$w_i=w_j$ is a contact point of~$\Lambda_{j+1}$, then~$\Lambda_1,\dots,\Lambda_{i-1},\Lambda'_i,\dots,\Lambda'_j,\Lambda_{j+1},\dots,\Lambda_r$ is an iterated \mpt of~$L$. We say that it is obtained from~$\Lambda_1,\dots,\Lambda_r$ by a \defn{partial flip} of~$v$.
\item If~$j<r$, and~$w_i=w_j$ is a crossing point in~$\Lambda_{j+1}$, then we cannot flip~$v$ in~$\Lambda_i$ maintaining an iterated \mpt of~$L$.
\end{enumerate}

\begin{figure}
	\capstart
	\centerline{\includegraphics[scale=1]{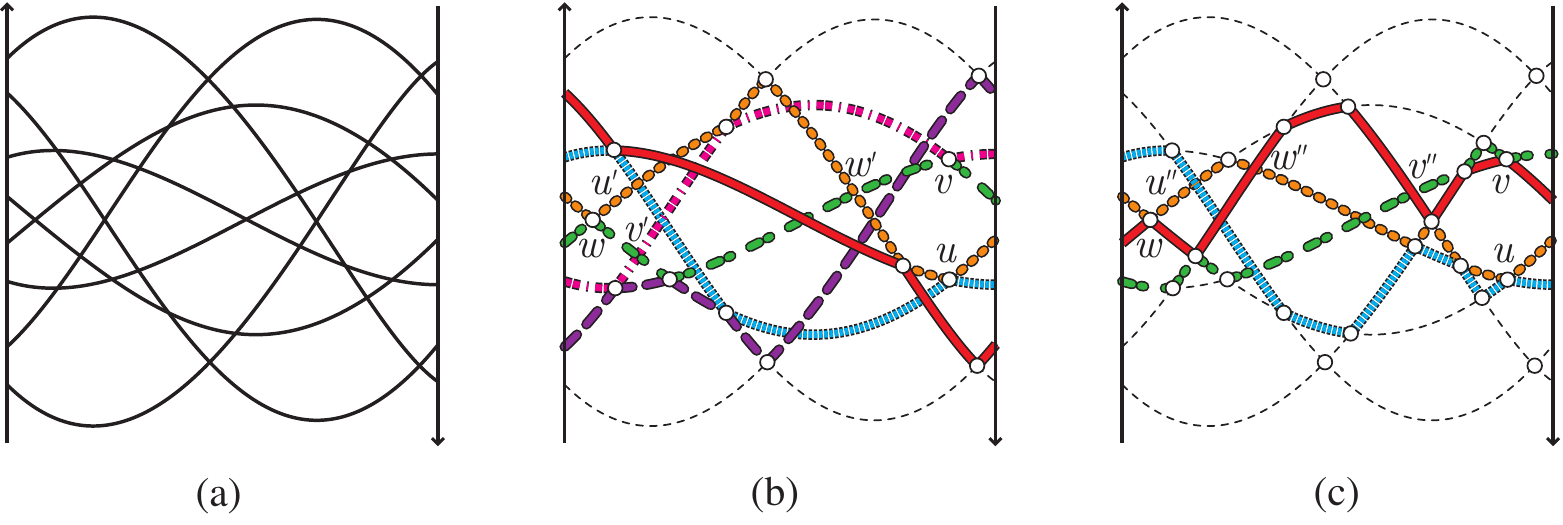}}
	\caption[The three possible situations for flipping a contact point in an iterated \mpt{}]{The three possible situations for flipping a contact point in an iterated \mpt{}.}
	\label{mpt:fig:iteratedflip}
\end{figure}

To illustrate these three possible cases, we have labeled on \fref{mpt:fig:iteratedflip} some intersection points of an iterated pseudotriangulation. We have choosen three contact points~$u,v,w$ in~(b). For~$z\in\{u,v,w\}$, we label $z'$ (resp.~$z''$) the crossing point corresponding to~$z$ in~(b) (resp. in~(c)). Observe that:
\begin{enumerate}[(i)]
\item points~$u'$~and~$u''$ coincide. Thus we can flip simultaneously point~$u$ in~(b) and~(c) (complete flip);
\item points~$v'$ is different from~$v''$ but is a contact point in~(c). Thus, we can just flip~$v$ in~(b), without changing~(c) and we preserve an iterated pseudotriangulation (partial flip);
\item point~$w'$ is a crossing point in~(c), different from~$w''$. Thus, we cannot flip~$w$ in~(b) maintaining an iterated pseudotriangulation.
\end{enumerate}

Let~$G^{k_1,\dots,k_r}(L)$ be the graph whose vertices are the iterated \mpt{}s of~$L$ with signature~$k_1<\dots<k_r$, and whose edges are the pairs of iterated \mpt{}s linked by a (complete or partial) flip.

\begin{theorem}\label{mpt:theo:iteratedflip}
The graph of flips~$G^{k_1,\dots,k_r}(L)$ is connected.
\end{theorem}

To prove this proposition, we need the following lemma:

\begin{lemma}\label{mpt:lem:exist}
Any intersection point~$v$ in the \kkernel{k} of a pseudoline arrangement is a contact point in a \pt{k} of it.
\end{lemma}

\begin{proof}
The result holds when~$k=1$. We obtain the general case by iteration.
\end{proof}

\begin{proof}[Proof of Theorem~\ref{mpt:theo:iteratedflip}]
We prove the result by induction on~$r$ ($L$~is fixed). When~$r=1$, we already know that the flip graph is connected. Now, let~$A_-$~and~$A_+$ be two iterated \mpt{}s of~$L$ with signature~$k_1<\dots<k_r$, that we want to join by flips. Let~$B_-$~and~$B_+$ be iterated \mpt{}s of~$L$ with signature~$k_1<\dots<k_{r-1}$, and~$\Lambda_-$~and~$\Lambda_+$ be \pt{k_r}s of~$L$ such that~$A_-=B_-,\Lambda_-$ and~$A_+=B_+,\Lambda_+$.

By induction,~$G^{k_1,\dots,k_{r-1}}(L)$ is connected: let~${B_-=B_1,B_2,\dots,B_{p-1},B_p=B_+}$ be a path from~$B_-$ to~$B_+$ in~$G^{k_1,\dots,k_{r-1}}(L)$. For all~$j$, let~$v_j$ be such that~$B_{j+1}$ is obtained from~$B_j$ by flipping~$v_j$  and let~$w_j$ be such that~$B_j$ is obtained from~$B_{j+1}$ by flipping~$w_j$. Let~$\Lambda_j$ be a \pt{k_r} of~$L$ containing the contact points of~$B_j$ plus~$w_j$ (it exists by Lemma~\ref{mpt:lem:exist}), and let~$C_j=B_j,\Lambda_j$. Let~$D_j$ be the iterated \mpt of~$L$ obtained from the iterated pseudotriangulation~$C_j$ by a partial flip of~$v_j$. Finally, since~$G^{k_r}(B_j)$ is connected, there is a path of complete flips from~$D_{j-1}$ to~$C_j$.

Merging all these paths, we obtain a global path from~$A_-$ to~$A_+$: we transform~$A_-$ into~$C_1$ via a path of complete flips; then~$C_1$ into~$D_1$ by the partial flip of~$v_1$; then~$D_1$ into~$C_2$ via a path of complete flips; then~$C_2$ into~$D_2$ by the partial flip of~$v_2$; and so on.
\end{proof}


\section{Further topics on \mpt{}s}\label{mpt:sec:furthertopics}

To finish this chapter, we discuss the extensions in the context of \mpt{}s of two known results on pseudotriangulations:
\begin{enumerate}
\item The first one concerns the relationship between the greedy pseudotriangulation of a point set and its horizon trees.
\item The second one extends to arrangements of double pseudolines the definition and properties of \mpt{}s.
\end{enumerate}


\subsection{Greedy \mpt{}s and horizon graphs}\label{mpt:subsec:furthertopics:horizon}

We have seen in previous sections that the greedy \pt{k} of a pseudoline arrangement~$L$ can be seen as:
\begin{enumerate}
\item the unique source of the graph of increasing flips on \pt{k}s of~$L$;
\item a greedy choice of crossing points given by a sorting network;
\item a greedy choice of contact points;
\item an iteration of greedy \pt{1}s.
\end{enumerate}
In this section, we provide a ``pattern avoiding'' characterization of the crossing points of the greedy \pt{k} of~$L$.

\svs
Let~$L$ be a pseudoline arrangement, and~$\chi$ be a cut of~$L$. We index by~$\ell_1,\dots,\ell_n$ the pseudolines of~$L$ in the order in which they cross~$\chi$ (it is well defined, up to a complete inversion).

We define the \defn{$k$-upper $\chi$-horizon set}\index{horizon!--- set} of~$L$ to be the set~$\UU^k_\chi(L)$ of crossing points~$\ell_\alpha\wedge\ell_\beta$, with~$1\le \alpha<\beta\le n$, such that there is no~$\gamma_1,\dots,\gamma_k$ satisfying~$\alpha<\gamma_1<\dots<\gamma_k$ and $\ell_\alpha\wedge\ell_{\gamma_i}\cle_\chi\ell_\alpha\wedge\ell_\beta$ for all~$i\in[k]$. In other words, on each pseudoline~$\ell_\alpha$ of~$L$, the set~$\UU^k_\chi(L)$ consists of the smallest~$k$ crossing points of the form~$\ell_\alpha\wedge\ell_\beta$, with~$\alpha<\beta$.

Similarly, define the \defn{$k$-lower $\chi$-horizon set} of~$L$ to be the set~$\LL^k_\chi(L)$ of crossing points~${\ell_\alpha\wedge\ell_\beta}$, with ${1\le \alpha<\beta\le n}$, such that there is no~$\delta_1,\dots,\delta_k$ satisfying~$\delta_1<\dots<\delta_k<\beta$ and ${\ell_\beta\wedge\ell_{\delta_i}\cle_\chi\ell_\alpha\wedge\ell_\beta}$ for all~$i\in[k]$. On each pseudoline~$\ell_\beta$ of~$L$, the set~$\LL^k_\chi(L)$ consists of the smallest~$k$ crossing points of the form~$\ell_\alpha\wedge\ell_\beta$, with~$\alpha<\beta$.

Finally, we define the set~$\GG^k_\chi(L)$ to be the set of crossing points~$\ell_\alpha\wedge\ell_\beta$, with $1\le \alpha<\beta\le n$, such that there is no~$\gamma_1,\dots,\gamma_k$ and~$\delta_1,\dots,\delta_k$ satisfying:
\begin{enumerate}[(i)]
\item $\alpha<\gamma_1<\dots<\gamma_k$, $\delta_1<\dots<\delta_k<\beta$, and~$\delta_k<\gamma_1$;~and
\item $\ell_\alpha\wedge\ell_{\gamma_i}\cle_\chi\ell_\alpha\wedge\ell_\beta$ and~$\ell_\beta\wedge\ell_{\delta_i}\cle_\chi\ell_\alpha\wedge\ell_\beta$ for all $i\in[k]$.
\end{enumerate}

Obviously, the sets~$\UU^k_\chi(L)$ and~$\LL^k_\chi(L)$ are both contained in~$\GG^k_\chi(L)$.

\begin{example}
In \fref{mpt:fig:greedyhorizon}, we have labeled the vertices of the pseudoline arrangement~$L$ of \fref{mpt:fig:dual}(b) with different geometric tags according to their status: 
\begin{itemize}
\item[($\triangle$)] crossing points of the $k$-upper $\chi$-horizon set~$\UU_\chi^k(L)$ are represented by up triangles~$\triangle$;
\item[($\vartriangledown$)] crossing points of the $k$-lower $\chi$-horizon set~$\LL_\chi^k(L)$ are represented by down triangles~$\vartriangledown$;
\item[({\Large $\davidsstar$})] crossing points in both~$\UU_\chi^k(L)$ and~$\LL_\chi^k(L)$) are represented by up and down triangles~{\Large $\davidsstar$};
\item[($\Box$)] crossing points of~$\GG_\chi^k(L)$ but neither in~$\UU_\chi^k(L)$, nor in~$\LL_\chi^k(L)$ are represented by squares~$\Box$.
\end{itemize}
Observe that the remaining vertices are exactly the crossing points of the \greedy{\chi} \pt{k} of~$L$.

\begin{figure}[h]
	\capstart
	\centerline{\includegraphics[scale=1]{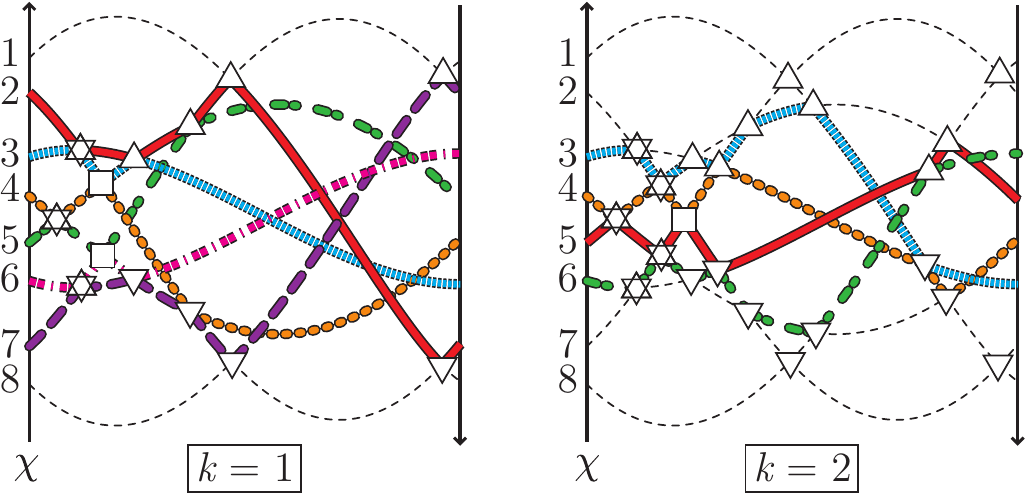}}
	\caption[The sets~$\UU^k_\chi(L)$,~$\LL^k_\chi(L)$, and~$\GG^k_\chi(L)$]{The sets $\UU^k_\chi(L)$, $\LL^k_\chi(L)$, and~$\GG^k_\chi(L)$ for the pseudoline arrangement of \fref{mpt:fig:dual}(b) and~$k\in\{1,2\}$. The underlying \pt{k} is the greedy \pt{k} of~$L$.}
	\label{mpt:fig:greedyhorizon}
\end{figure}

\end{example}

\begin{example}\label{mpt:exm:convex}
We consider the arrangement~$V_n^*$ of~$n$ pseudolines in convex position. Let~$z$ be a vertex on the upper hull of its support,~$F \eqdef \ens{z'}{z\cle z'}$ denote the filter generated by~$z$, and~$\chi$ denote the corresponding cut (see \fref{mpt:fig:convex}). It is easy to check that:
\begin{enumerate}[(i)]
\item $\UU^k_\chi(V_n^*)=\ens{\ell_\alpha\wedge\ell_\beta}{1\le \alpha\le n\text{ and }\alpha<\beta\le \alpha+k}$;
\item $\LL^k_\chi(V_n^*)=\ens{\ell_\alpha\wedge\ell_\beta}{1\le \alpha\le k\text{ and }\alpha< j\le n}$;~and
\item $\GG^k_\chi(V_n^*)=\UU^k_\chi(\cC_n)\cup\LL^k_\chi(\cC_n)$.
\end{enumerate}
Observe again that the remaining vertices are exactly the crossing points of the \greedy{\chi} \pt{k} of~$V_n^*$.

\begin{figure}
	\capstart
	\centerline{\includegraphics[scale=1]{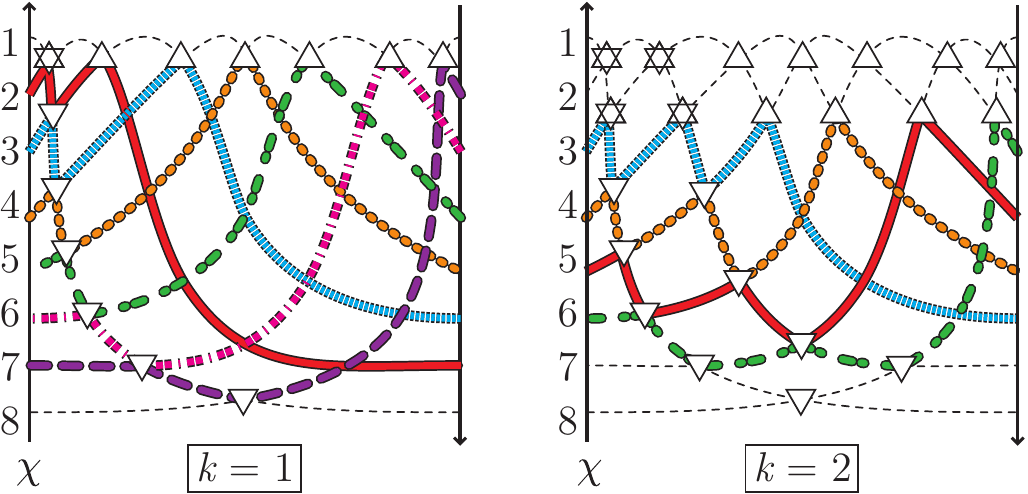}}
	\caption[The sets $\UU^k_\chi(V_8^*)$, $\LL^k_\chi(V_8^*)$, and $\GG^k_\chi(V_8^*)$]{The sets $\UU^k_\chi(V_8^*)$, $\LL^k_\chi(V_8^*)$, and $\GG^k_\chi(V_8^*)$ for the arrangement~$V_8^*$ of~$8$ pseudolines in convex position, and~$k\in\{1,2\}$. The underlying \pt{k} is the greedy \pt{k} of~$V_8^*$.}
	\label{mpt:fig:convex}
\end{figure}

\end{example}

\enlargethispage{.1cm}
Theorem~\ref{mpt:theo:horizon} extends this observation to all pseudoline arrangements, using convex position as a starting point for a proof by mutation.

\begin{theorem}\label{mpt:theo:horizon}
For any pseudoline arrangement~$L$ with no contact point, and any cut~$\chi$ of~$L$, the sets~$V(\Gamma_\chi^k(L))$ and~$\GG^k_\chi(L)$ coincide.
\end{theorem}

\begin{figure}[b]
	\capstart
	\centerline{\includegraphics[scale=1]{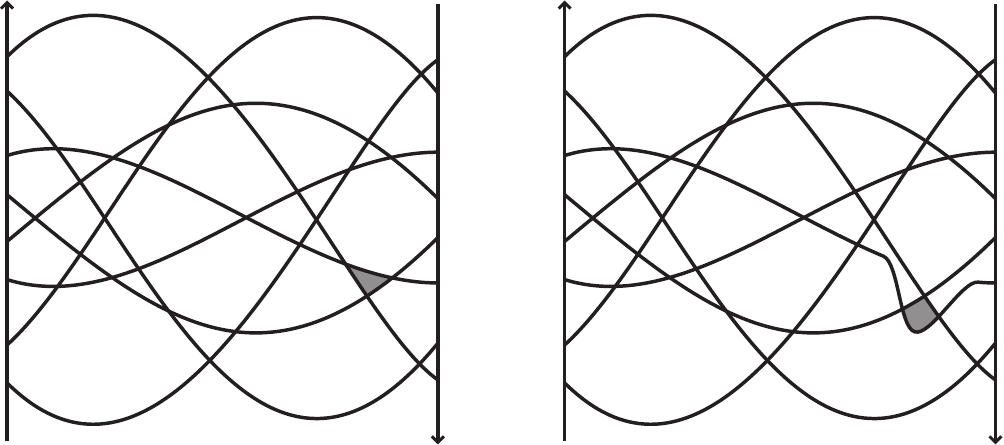}}
	\caption[A mutation in a pseudoline arrangement]{A mutation in the arrangement of \fref{mpt:fig:dual}(b).}
	\label{mpt:fig:mutation}
\end{figure}

The proof of this theorem works by mutation. A \defn{mutation}\index{mutation} is a local transformation of an arrangement~$L$ that only inverts one triangular face of~$L$. More precisely, it is a homotopy of arrangements during which only one curve~$\ell\in L$ moves, sweeping a single vertex of the remaining arrangement~$L\ssm\{\ell\}$ (see \fref{mpt:fig:mutation} and Appendix~\ref{app:sec:dpl} for more details).

If~$P$ is a point set of a topological plane, mutating an empty triangle~$p^*q^*r^*$ of~$P^*$ by sweeping the vertex~$q^*\wedge r^*$ with the pseudoline~$p^*$ corresponds in the primal to moving a little bit~$p$ such that only the orientation of the triangle~$pqr$ changes.

The graph of mutations on pseudoline arrangements is known to be connected: any two pseudoline arrangements (with no contact points and the same number of pseudolines) are homotopic via a finite sequence of mutations (followed by a homeomorphism). In fact, one can even avoid mutations of triangles that cross a given cut of~$L$:

\begin{proposition}\label{mpt:prop:mutation}
Let~$L$ and~$L'$ be two pseudoline arrangements of~$\cM$ (with no contact points and the same number of pseudolines) and~$\chi$ be a cut of both of~$L$ and~$L'$. There is a finite sequence of mutations of triangles disjoint from~$\chi$ that transforms~$L$ into~$L'$.
\end{proposition}

\begin{proof}
We prove that any arrangement~$L$ of~$n$ pseudolines can be transformed into the arrangement~$V_n^*$ of $n$ pseudolines in convex position (see \fref{mpt:fig:convex}).

Let~$\ell_1,\dots,\ell_n$ denote the pseudolines of~$L$ (ordered by their crossings with~$\chi$). Let~$\Delta$ denote the triangle formed by~$\chi$, $\ell_1$~and~$\ell_2$. If there is a vertex of the arrangement~$L\ssm\{\ell_1,\ell_2\}$ inside~$\Delta$, then there is a triangle of the arrangement~$L$ inside~$\Delta$ and adjacent to~$\ell_1$ or~$\ell_2$. Mutating this triangle reduces the number of vertices of~$L\ssm\{\ell_1,\ell_2\}$ inside~$\Delta$ such that after some mutations, there is no more vertex inside~$\Delta$. If~$\Delta$ is intersected by pseudolines of~$L\ssm\{\ell_1,\ell_2\}$, then there is a triangle inside~$\Delta$ formed by~$\ell_1$,~$\ell_2$~and one of these intersecting pseudolines (the one closest to~$\ell_1\wedge\ell_2$). Mutating this triangle reduces the number of pseudolines intersecting~$\Delta$. Thus, after some mutations,~$\Delta$~is a triangle of the arrangement~$L$.

Repeating these arguments, we can affirm that, for all~$i\in\{2,\dots,n-1\}$ and after some mutations, $\ell_i$, $\ell_1$, $\ell_{i+1}$ and~$\chi$ delimit a face of the arrangement~$L$. Thus, one of the two topological disk delimited by~$\chi$ and~$\ell_1$ contains no more vertex of~$L$, and the proof is then straightforward by induction.
\end{proof}

\begin{figure}
	\capstart
	\centerline{\includegraphics[scale=1]{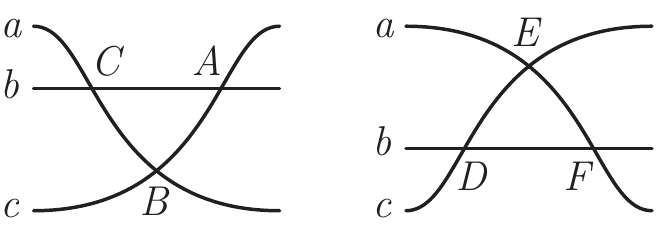}}
	\caption[A local image of a mutation]{A local image of a mutation.}
	\label{mpt:fig:mutationlocal}
\end{figure}

Let~$\nabla$ be a triangle of $L$ not intersecting~$\chi$. Let~$L'$ denote the pseudoline arrangement obtained from~$L$ by mutating the triangle~$\nabla$ into the inverted triangle~$\Delta$.
Let~$a<b<c$ denote the indices of the pseudolines~$\ell_a,\ell_b$ and~$\ell_c$ that form~$\nabla$ and~$\Delta$. In~$\nabla$, we denote~$A=\ell_b\wedge\ell_c$, $B=\ell_a\wedge\ell_c$ and~$C=\ell_a\wedge\ell_b$; similarily, in~$\Delta$, we denote~$D=\ell_b\wedge\ell_c$, $E=\ell_a\wedge\ell_c$ and~$F=\ell_a\wedge\ell_b$ (see \fref{mpt:fig:mutationlocal}).

\begin{lemma}\label{mpt:lem:mutation}
With these notations, the following properties hold:
\begin{enumerate}[(i)]
\item $B\in\UU_\chi^k(L)\Leftrightarrow C\in\UU_\chi^{k-1}(L)\Leftrightarrow E\in\UU_\chi^{k-1}(L')\Leftrightarrow F\in\UU_\chi^k(L')$\\
$A\in\UU_\chi^k(L)\Leftrightarrow D\in\UU_\chi^k(L')$\\
$E\in\LL_\chi^k(L')\Leftrightarrow D\in\LL_\chi^{k-1}(L')\Leftrightarrow B\in\LL_\chi^{k-1}(L)\Leftrightarrow A\in\LL_\chi^k(L)$\\
$F\in\LL_\chi^k(L')\Leftrightarrow C\in\LL_\chi^k(L)$
\item $C\in\UU_\chi^k(L)\Rightarrow A\in\UU_\chi^k(L)$\\
$E\in\UU_\chi^k(L')\Rightarrow D\in\UU_\chi^k(L')$\\
$D\in\LL_\chi^k(L')\Rightarrow F\in\LL_\chi^k(L')$\\
$B\in\LL_\chi^k(L)\Rightarrow C\in\LL_\chi^k(L)$\\
$B\in\GG_\chi^k(L)\Rightarrow C\in\GG_\chi^k(L)\wedge D\in\GG_\chi^k(L')\wedge F\in\GG_\chi^k(L')$\\
$E\in\GG_\chi^k(L')\Rightarrow D\in\GG_\chi^k(L')\wedge C\in\GG_\chi^k(L)\wedge A\in\GG_\chi^k(L)$
\item $C\in\GG_\chi^k(L)\Rightarrow A\in\UU_\chi^k(L)\vee C\in\LL_\chi^k(L)$\\
$D\in\GG_\chi^k(L')\Rightarrow F\in\LL_\chi^k(L')\vee D\in\UU_\chi^k(L')$\\
$A\in\GG_\chi^k(L)\Rightarrow A\in\UU_\chi^k(L)\vee C\in\LL_\chi^{k-1}(L)$\\
$F\in\GG_\chi^k(L')\Rightarrow F\in\LL_\chi^k(L')\vee D\in\UU_\chi^{k-1}(L')$
\item $C\in\GG_\chi^k(L)\wedge E\notin\GG_\chi^k(L')\Rightarrow A\notin\GG_\chi^k(L)$\\
$D\in\GG_\chi^k(L')\wedge B\notin\GG_\chi^k(L)\Rightarrow F\notin\GG_\chi^k(L')$
\end{enumerate}
\end{lemma}

\begin{proof}
By symmetry, it is enough to prove the first line of each of the four points of the lemma.

Properties of point~(i) directly come from the definitions. For example, all the assertions of the first line are false if and only if there exist~$\gamma_1,\dots,\gamma_{k-1}$ with~$a<\gamma_1<\dots<\gamma_{k-1}$ and, for all~$i\in[k-1]$, $\ell_a\wedge\ell_{\gamma_i}\cle_\chi C$ (or equivalently~$\ell_a\wedge\ell_{\gamma_i}\cle_\chi E$).

We derive point~(ii) from the following observation: if~$\gamma>b$ and if~$\ell_b\wedge\ell_\gamma\cle_\chi C$, then~$\gamma>a$ and~$\ell_a\wedge\ell_\gamma\cle_\chi B$.

For point~(iii), assume that~$A\notin\UU_\chi^k(L)$ and~$C\notin\LL_\chi^k(L)$. Then there exist~$\gamma_1,\dots,\gamma_k$ and~$\delta_1,\dots,\delta_k$ such that~$\delta_1<\dots<\delta_k<b<\gamma_1<\dots<\gamma_k$ and, for all $i\in[k]$, 
$\ell_b\wedge\ell_{\gamma_i}\cle_\chi A$ (and therefore $\ell_a\wedge\ell_{\gamma_i}\cle_\chi C$) and $\ell_b\wedge\ell_{\delta_i}\cle_\chi C$. Thus~$C\notin\GG_\chi^k(L)$.

Finally, assume that~$C\in\GG_\chi^k(L)$ and~$E\notin\GG_\chi^k(L')$. Then, there exist~$\gamma_1,\dots,\gamma_k$ and $\delta_1,\dots,\delta_k$ such that $a<\gamma_1<\dots<\gamma_k$, $\delta_1<\dots<\delta_k<c$, $\delta_k<\gamma_1$, and for all $i\in[k]$, $\ell_a\wedge\ell_{\gamma_i}\cle_\chi E$ and $\ell_c\wedge\ell_{\delta_i}\cle_\chi E$. Since $C\in\GG_\chi^k(L)$, we have $\delta_k>b$. Thus $b<\gamma_1<\dots<\gamma_k$ and for all $i\in[k]$, $\ell_b\wedge\ell_{\gamma_i}\cle_\chi A$ and $\ell_c\wedge\ell_{\delta_i}\cle_\chi A$. This implies that~$A\notin\GG_C^k(L)$.
\end{proof}

\begin{figure}
	\capstart
	\centerline{\includegraphics[scale=1]{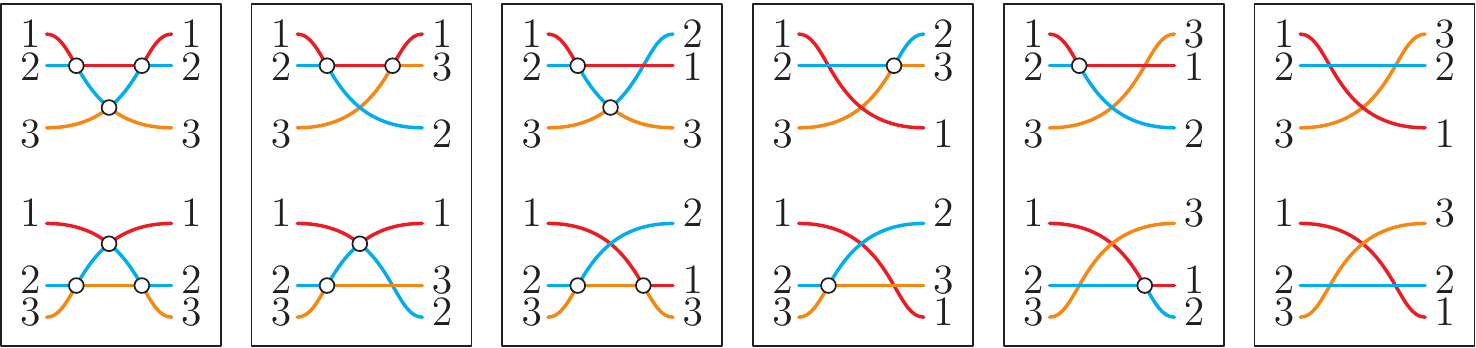}}
	\caption[The six possible cases for the mutation of the greedy \pt{k}]{The six possible cases for the mutation of the greedy \pt{k}.}
	\label{mpt:fig:mutationlocal6sol}
\end{figure}

We are now ready to establish the proof of Theorem~\ref{mpt:theo:horizon}:

\begin{proof}[Proof of Theorem~\ref{mpt:theo:horizon}]
The proof works by mutation. We already observed the result when the pseudoline arrangement is in convex position (see Example~\ref{mpt:exm:convex} and \fref{mpt:fig:convex}). Proposition~\ref{mpt:prop:mutation} ensures that any pseudoline arrangement can be reached from this convex configuration by mutations of triangles not intersecting~$\chi$. Thus, it is sufficient to prove that such a mutation preserves the property.

Assume that~$L$ is a pseudoline arrangement and~$\chi$ is a cut of~$L$, for which the result holds. Let~$\nabla$ be a triangle of~$L$ not intersecting~$\chi$. Let~$L'$ denote the pseudoline arrangement obtained from~$L$ by mutating the triangle~$\nabla$ into the inverted triangle~$\Delta$. Let~$A,B,C$~and~$D,E,F$ denote the vertices of~$\nabla$~and~$\Delta$ as indicated in \fref{mpt:fig:mutationlocal}. 

If~$v$ is a vertex of the arrangement~$L'$ different from~$D,E,F$, then:
$$v\in V(\Gamma^k_\chi(L'))\Leftrightarrow v\in V(\Gamma^k_\chi(L))\Leftrightarrow v\in\GG^k_\chi(L)\Leftrightarrow v\in\GG^k_\chi(L').$$

Thus, we only have to prove the equivalence when~$v\in\{D,E,F\}$. The proof is a (computational) case analysis: using the properties of Lemma~\ref{mpt:lem:mutation} as boolean equalities relating the boolean variables~``$X\in\mathbb{Y}_\chi^p(L)$'' (where~$X\in\{A,B,C,D,E,F\}$, $\mathbb{Y}\in\{\UU,\LL,\GG\}$, and ${p\in\{k-1,k\}}$), we have written a short boolean satisfiability program which affirms that:
\begin{enumerate}[(i)]
\item either $\{A,B,C\}\subset\GG_\chi^k(L)$ and $\{D,E,F\}\subset\GG_\chi^k(L')$;
\item or $\{A,B,C\}\cap\GG_\chi^k(L)=\{A,C\}$ and $\{D,E,F\}\cap\GG_\chi^k(L')=\{D,E\}$;
\item or $\{A,B,C\}\cap\GG_\chi^k(L)=\{B,C\}$ and $\{D,E,F\}\cap\GG_\chi^k(L')=\{D,F\}$;
\item or $\{A,B,C\}\cap\GG_\chi^k(L)=\{A\}$ and $\{D,E,F\}\cap\GG_\chi^k(L')=\{D\}$;
\item or $\{A,B,C\}\cap\GG_\chi^k(L)=\{C\}$ and $\{D,E,F\}\cap\GG_\chi^k(L')=\{F\}$;
\item or $\{A,B,C\}\cap\GG_\chi^k(L)=\emptyset$ and $\{D,E,F\}\cap\GG_\chi^k(L')=\emptyset$.
\end{enumerate}

It is easy to check that these six cases correspond to sorting the six possible permutations of~$\{1,2,3\}$ on~$\nabla$ and~$\Delta$ (see \fref{mpt:fig:mutationlocal6sol}). Thus, if~$V(\Gamma_\chi^k(L))\cap\{A,B,C\}=\GG_\chi^k(L)\cap\{A,B,C\}$, then~$V(\Gamma_\chi^k(L'))\cap\{D,E,F\}=\GG_\chi^k(L')\cap\{D,E,F\}$, which finishes the proof.
\end{proof}

\begin{figure}[h]
	\capstart
	\centerline{\includegraphics[scale=1]{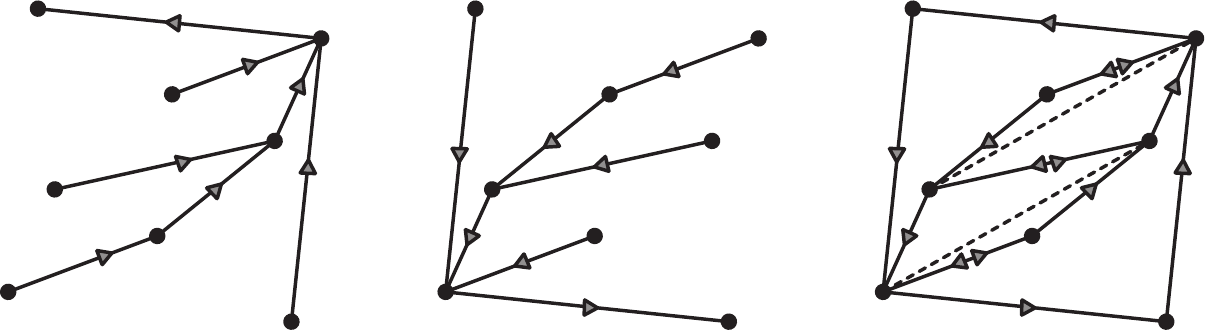}}
	\caption[Horizon trees]{The upper horizon tree (left), the lower horizon tree (middle), and the greedy pseudotriangulation (right) of the point set of \fref{mpt:fig:dual}(a). The dashed edges in the greedy pseudotriangulation are not in the horizon trees.}
	\label{mpt:fig:horizon}
\end{figure}

Let us finish this discussion by recalling the interpretation of the horizon sets when~${k=1}$.
Let~$P$ be a finite point set. For any~$p\in P$, let~$u(p)$ denote the point~$q$ that minimizes the angle~$(Ox,pq)$ among all points of~$P$ with~$y_p<y_q$ (by convention, for the higher point~$p$ of~$P$,~$u(p)=p$). The \defn{upper horizon tree}\index{horizon!--- tree} of~$P$ is the set~$\UU(P)=\ens{pu(p)}{p\in P}$. The \defn{lower horizon tree}~$\LL(P)$ of $P$ is defined symetrically. See \fref{mpt:fig:horizon}.

Choosing a cut~$\chi$ of~$P^*$ corresponding to the point at infinity~$(-\infty,0)$ makes coincide primal and dual definitions of horizon sets: $\UU^1_\chi(P^*)=\UU(P)^*$~and~$\LL^1_\chi(P^*)=\LL(P)^*$.

In~\cite{p-htvpt-97}, Michel~Pocchiola observed that the set~$\UU(P)\cup\LL(P)$ of edges can be completed into a pseudotriangulation of~$P$ just by adding the sources of the faces of~$P^*$ intersected by the cut~$\chi$. The obtained pseudotriangulation is our \greedy{\chi} \pt{1}~$\Gamma^1_\chi(P^*)$.



\subsection{Multipseudotriangulations of double pseudoline arrangements}\label{mpt:subsec:furthertopics:dpl}

In this section, we deal with double pseudoline arrangements, \ie duals of sets of disjoint convex bodies. Definitions and properties of \mpt{}s naturally extend to these objects.

\subsubsection{Definitions}

\begin{figure}[b]
	\capstart
	\centerline{\includegraphics[scale=1]{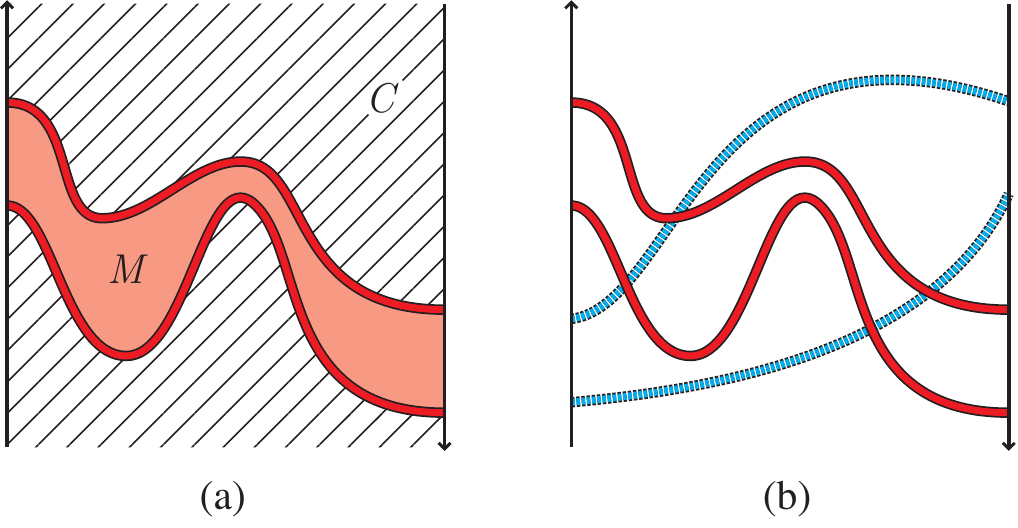}}
	\caption[A double pseudoline and a double pseudoline arrangement in the M\"obius strip]{(a) A double pseudoline; (b) An arrangement of~$2$ double pseudolines.}
	\label{mpt:fig:dpl}
\end{figure}

A simple closed curve in the M\"obius strip can be:
\begin{enumerate}[(i)]
\item either contractible (homotopic to a point);
\item or non separating, or equivalently homotopic to a generator of the fundamental group of~$\cM$: it is a pseudoline;
\item or separating and non-contractible, or equivalently homotopic to the double of a generator of the fundamental group of~$\cM$: it is called a \defn{double pseudoline}\index{double pseudoline} (see \fref{mpt:fig:dpl}(a)).
\end{enumerate}

The complement of a double pseudoline~$\ell$ has two connected components: the bounded one is a M\"obius strip~$M_\ell$ and the undounded one is an open cylinder~$C_\ell$ (see \fref{mpt:fig:dpl}(a)).

\begin{figure}
	\capstart
	\centerline{\includegraphics[scale=1]{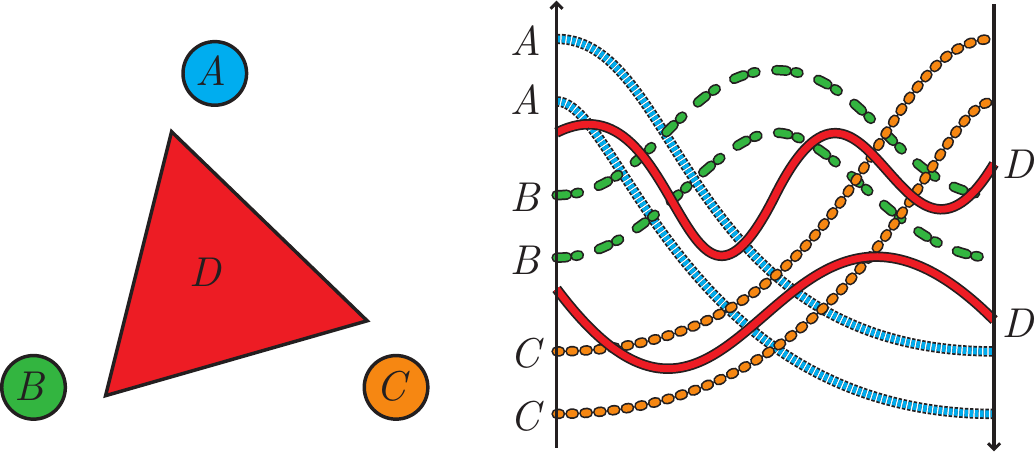}}
	\caption[A configuration of~$4$ disjoint convex bodies and its dual double pseudoline arrangement]{A configuration of four disjoint convex bodies and its dual double pseudoline arrangement.}
	\label{mpt:fig:dplactrexm1}
\end{figure}

\begin{figure}[b]
	\capstart
	\centerline{\includegraphics[scale=1]{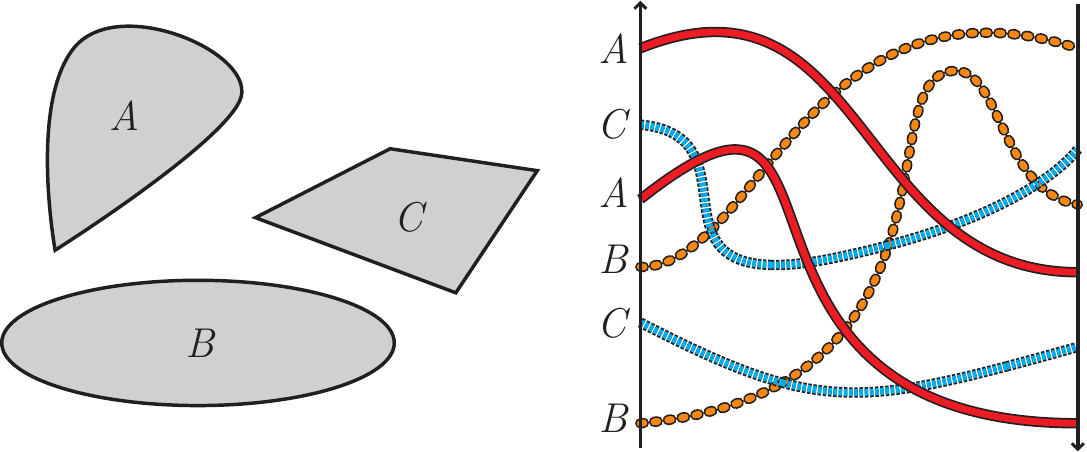}}
	\caption[A configuration of~$3$ disjoint convex bodies and its dual double pseudoline arrangement]{A configuration of three disjoint convex bodies and its dual double pseudoline arrangement.}
	\label{mpt:fig:dpla}
\end{figure}

\index{dual!--- of a convex body|hbf}
\index{dual!--- of a set of convex bodies|hbf}
The canonical example of a double pseudoline is the set~$C^*$ of tangents to a convex body~$C$ of the plane. Observe also that the $p$th level of a pseudoline arrangement is a double pseudoline. If~$C$ is a convex body of the plane, then the M\"obius strip~$M_{C^*}$ corresponds to lines that pierce~$C$, while~$C_{C^*}$ corresponds to lines disjoint from~$C$. If~$C$ and~$C'$ are two disjoint convex bodies, the two corresponding double pseudolines~$C^*$ and~$C'^*$ cross four times (see \fref{mpt:fig:dplactrexm1} and~\fref{mpt:fig:dpla}). Each of these four crossings corresponds to one of the four \defn{bitangents}\index{bitangent} (or \defn{common tangents}) between~$C$~and~$C'$.

\begin{definition}[\cite{hp-adp-08}]
\index{double pseudoline!--- arrangement}
A \defn{double pseudoline arrangement} is a finite set of double pseudolines such that any two of them have exactly four intersection points, cross transversally at these points, and induce a cell decomposition of the M\"obius strip (\ie the complement of their union is a union of topological disks, together with the external cell).
\end{definition}


Given a set~$Q$ of disjoint convex bodies in the plane (or in any topological plane), its dual $Q^* \eqdef \ens{C^*}{C\in Q}$ is an arrangement of double pseudolines (see \fref{mpt:fig:dplactrexm1} and \fref{mpt:fig:dpla}). Furthermore, as for pseudoline arrangements, any double pseudoline arrangement can be represented by (\ie is the dual of) a set of disjoint convex bodies in a topological plane~\cite{hp-adp-08}. For more details on double pseudoline arrangements, we refer to our Appendix~\ref{app:sec:dpl}, where we discuss an algorithm to enumerate all arrangements with few double pseudolines. 

\begin{figure}
	\capstart
	\centerline{\includegraphics[scale=1]{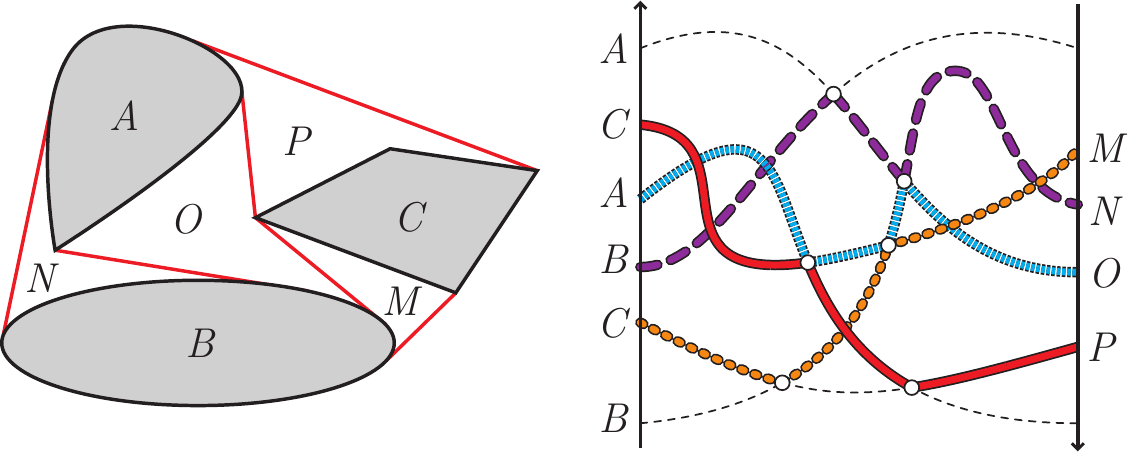}}
	\caption[A pseudotriangulation of a set of disjoint convex bodies]{A pseudotriangulation of the set of disjoint convex bodies of \fref{mpt:fig:dpla}.}
	\label{mpt:fig:pseudotriangulation}
\end{figure}

In this chapter, we only consider \defn{simple} arrangements of double pseudolines. Defining the \defn{support}\index{support (of an arrangement)}, the \defn{levels}\index{levels (of an arrangement)}, and the \defn{kernels}\index{kernel@\kkernel{k} (of an arrangement)} of double pseudoline arrangements as for pseudoline arrangements, we can extend \mpt{}s to double pseudoline arrangements (see \fref{mpt:fig:pseudotriangulation}):

\begin{definition}
\index{multipseudotriangulation!--- of a double pseudoline arrangement}
\index{pseudotriangulation@\pt{k}!--- of a double pseudoline arrangement}
A \pt{k} of a double pseudoline arrangement~$L$ is a pseudoline arrangement supported by the \kkernel{k} of~$L$.
\end{definition}

All the properties related to flips developed in Section~\ref{mpt:sec:enumeration} apply in this context. In the end of this section, we only revisit the properties of the primal of a multipseudotriangulation of a double pseudoline arrangement.

\subsubsection{Elementary properties}

Let~$Q$ be a set disjoint convex bodies in general position in the plane and~$Q^*$ be its dual arrangement. Let~$\Lambda$ be a \pt{k} of~$Q^*$, $V(\Lambda)$~denote all crossing points of~$Q^*$ that are not crossing points of~$\Lambda$, and~$E$ denote the corresponding set of bitangents of~$Q$. As in Subsection~\ref{mpt:subsec:mpt:pointedcrossing}, we discuss the properties of the primal configuration~$E$:

\begin{lemma}
The set~$E$ has~$4|Q|k-|Q|-2k^2-k$ edges.
\end{lemma}

\proof
The number of edges of~$E$ is the number of crossing points of~$Q^*$ minus the number of crossing points of~$\Lambda$, \ie
$$|E|=4{|Q^*| \choose 2}-{|\Lambda| \choose 2}=4{|Q| \choose 2}-{2|Q|-2k \choose 2}=4|Q|k-|Q|-2k^2-k.$$

\vspace{-1cm}\qed
\vspace{.5cm}

We now discuss pointedness. For any convex body~$C$ of~$Q$, we arbitrarily choose a point~$p_C$ in the interior of~$C$, and we consider the set~$X_C$ of all segments between~$p_C$ and a sharp boundary point of~$C$. We denote by~$X \eqdef \bigcup_{C\in Q} X_C$ the set of all these segments.

\begin{lemma}
The set $E\cup X$ cannot contain a \kalter{k}.\index{alternation@\kalter{k}}
\end{lemma}

\begin{figure}
	\capstart
	\centerline{\includegraphics[scale=1]{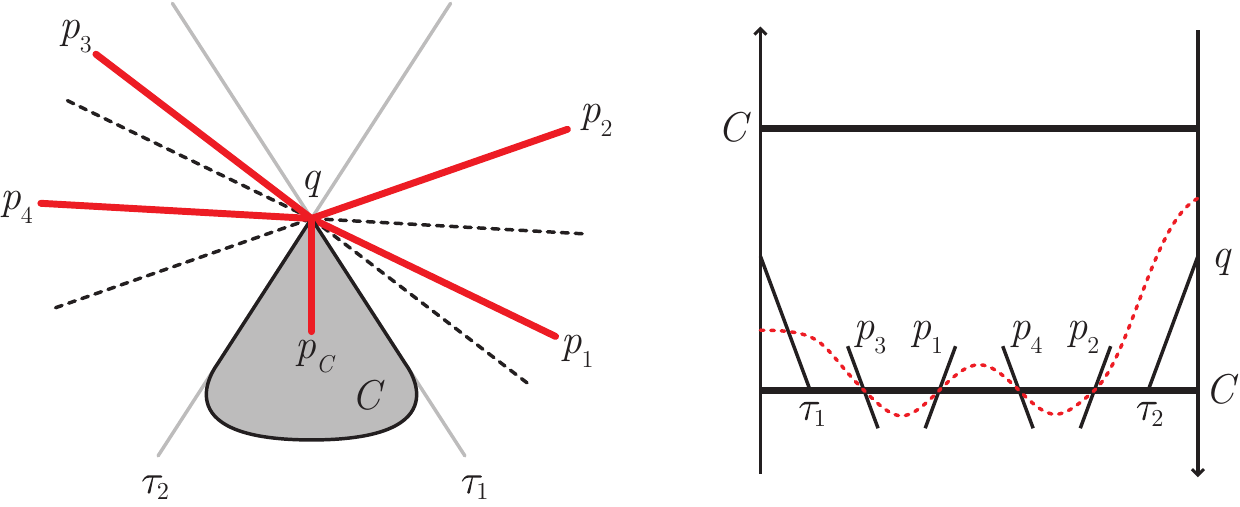}}
	\caption[A \kalter{k} at a sharp point]{A \kalter{k} at a sharp point.}
	\label{mpt:fig:sharp}
\end{figure}

\begin{proof}
Let~$C$ be a convex body of~$Q$, $q$~be a sharp point of~$C$ and~$F \eqdef \ens{[p_i,q]}{i\in[2k]}$ be a set of edges incident to~$q$ such that $\{[p_C,q]\}\cup F$ is a \kalter{k}. We prove that $F$ is not contained in~$E$. Indeed, the dual pseudolines~$\ens{p_i^*}{i\in[2k]}$ intersect alternately the double pseudoline~$C^*$ between the tangents to~$C$ at~$q$ (see \fref{mpt:fig:sharp}). This ensures the existence of a witness pseudoline~$\ell$ that separates all the contact points~$p_i^*\wedge C^*$, while crossing~$C^*$ exactly~$2k$ times and the other double pseudolines of~$Q^*$ exactly has~$q^*$ does. (As usual, we obtain it by a perturbation of the pseudoline~$q^*$.) Counting the crossings of~$\ell$ with~$Q^*$ and~$\Lambda$ respectively, we obtain:
\begin{enumerate}[(i)]
\item $\ell$~crosses~$Q^*$ exactly~$2|(Q\ssm\{C\})^*|+2k=2|Q|+2k-2$ times;
\item $\ell$~crosses~$\Lambda$ at least~$|\Lambda|=2|Q|-2k$ times;
\item for each of the points~$p_i\wedge q^*$, replacing the crossing point by a contact point removes two crossings with~$\ell$.\qedhere
\end{enumerate}
\end{proof}

\subsubsection{Stars}

If~$\lambda$ is a pseudoline of~$\Lambda$, we call \defn{star}\index{star} the envelope~$S(\lambda)$ of the primal lines of the points of~$\lambda$. The star~$S(\lambda)$ contains:
\begin{enumerate}[(i)]
\item all bitangents~$\tau$ between two convex bodies of~$Q$ such that~$\tau$ is a contact point of~$\lambda$; and
\item all convex arcs formed by the tangent points of the lines covered by~$\lambda$ with the convex bodies of~$Q$.
\end{enumerate}
This star is a (non-necessarily simple) closed curve. We again have bounds on the number of \defn{corners}\index{corner} (\ie convex internal angles) of~$S(\lambda)$:

\begin{proposition}\label{mpt:prop:cornersconv}
The number of corners of~$S(\lambda)$ is odd and between~$2k-1$ and~$4k|Q|-2k-1$.
\end{proposition}

\begin{proof}
In the case of double pseudoline arrangements, corners are even easier to characterize: a bitangent~$\tau$ between two convex~$C$ and~$C'$ of~$Q$ always defines two corners, one at each extremity. These corners are contained in one of the two stars adjacent to~$\tau$. Let~$\lambda$ be a pseudoline with a contact point at~$\tau$. In a neighbourhood of~$\tau$, the pseudoline~$\lambda$ can be contained either in~$M_{C^*}$ or in~$C_{C^*}$. In the first case, the star~$S(\lambda)$ contains the corner formed by the bitangent~$\tau$ and the convex~$C$ (or possibly, by the bitangent~$\tau$ and another tangent to~$C$); while in the second case, it does not. (The same observation holds for~$C'$.)

\begin{figure}
	\capstart
	\centerline{\includegraphics[scale=1]{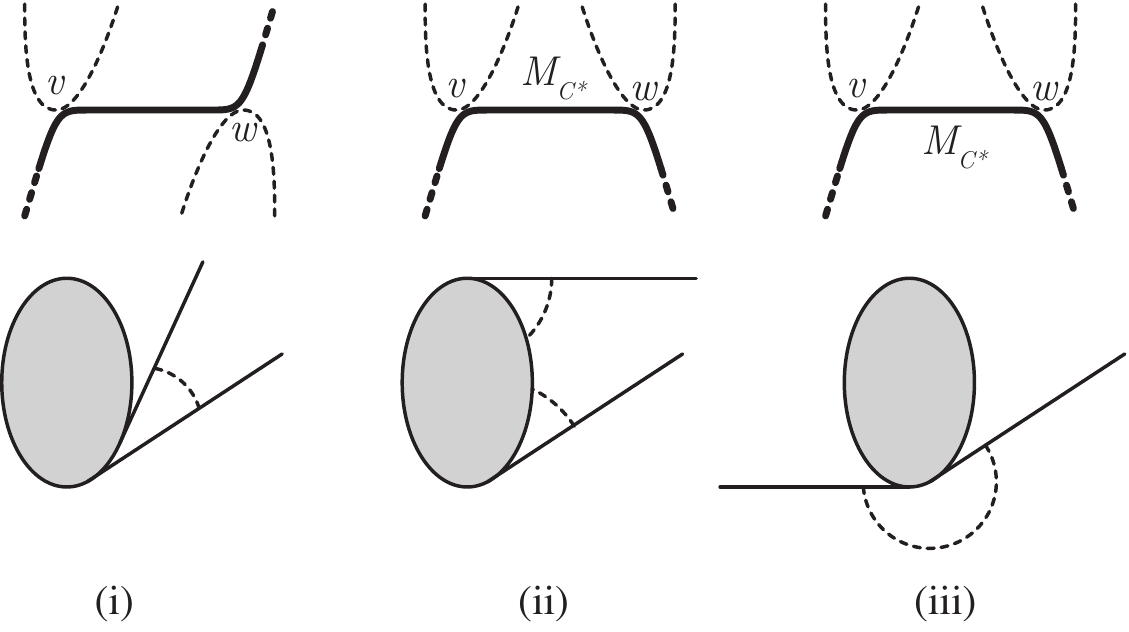}}
	\caption[Corners of a star]{The three possible situations for two consecutive contact points on~$\ell$.}
	\label{mpt:fig:corners}
\end{figure}

In other words, if the double pseudoline~$C^*$ supports a pseudoline~$\lambda$ between two contact points~$v$ and~$w$, then one of the three following situations occurs (see \fref{mpt:fig:corners}):
\begin{enumerate}[(i)]
\item either~$v$ and~$w$ lie on opposite sides of~$\lambda$; then exactly one of these contact points lies in~$M_{C^*}$, and~$S(\lambda)$ has one corner at~$C$.
\item or~$v$ and~$w$ both lie on~$M_{C^*}$, and~$S(\lambda)$ has two corners at~$C$.
\item or~$v$ and~$w$ both lie on~$C_{C^*}$, and~$S(\lambda)$ has no corners at~$C$.
\end{enumerate}
In particular, the number~$c=c(\lambda)$ of corners of~$\lambda$ is the number of situations~(i) plus twice the number of situations~(ii). This proves that~$c$ is odd and (at least) bigger than the number of opposite consecutive contact points of the pseudoline~$\lambda$.

In order to get a lower bound on this number, we construct (as in the proof of Proposition~\ref{mpt:prop:corners}) a witness pseudoline~$\mu$ that crosses~$\lambda$ between each pair of opposite contact points and passes on the opposite side of each contact point. It crosses~$\lambda$ at most~$c$ times and~$\Lambda\ssm\{\lambda\}$ exactly~$|\Lambda|-1$ times. Moreover, if~$\alpha$ is a pseudoline and~$\beta$ is a double pseudoline of~$\cM$, then either~$\alpha$ is contained in~$M_\beta$ and has no crossing with~$\beta$, or~$\alpha$ and~$\beta$ have an even number of crossings. Since~$\mu$ is a pseudoline and can be contained in at most one M\"obius strip~$M_C^*$ (for~$C\in Q$), the number of crossings of~$\mu$ with~$Q^*$ is at least~$2(|Q|-1)$. Thus, we obtain the lower bound $2(|Q|-1)\le2|Q|-2k-1+c$, \ie $c\ge 2k-1$.

From this lower bound, we obtain the upper bound: the total number of corners is at most twice the number of bitangents:
$$2(4k|Q|-|Q|-2k^2-k)\ge\sum_{\mu\in\Lambda} c(\mu)\ge c(\lambda)+(2|Q|-2k-1)(2k-1),$$
and we get,~$c\le4|Q|k-2k-1$.
\end{proof}

When~$k=1$, we can even prove that all stars are pseudotriangles. Indeed, since any star has at least~$3$ corners, the upper bound calculus gives~$2(3|Q|-3)\ge c+3(2|Q|-3)$, \ie $c\le 3$.

\begin{figure}
	\capstart
	\centerline{\includegraphics[scale=1]{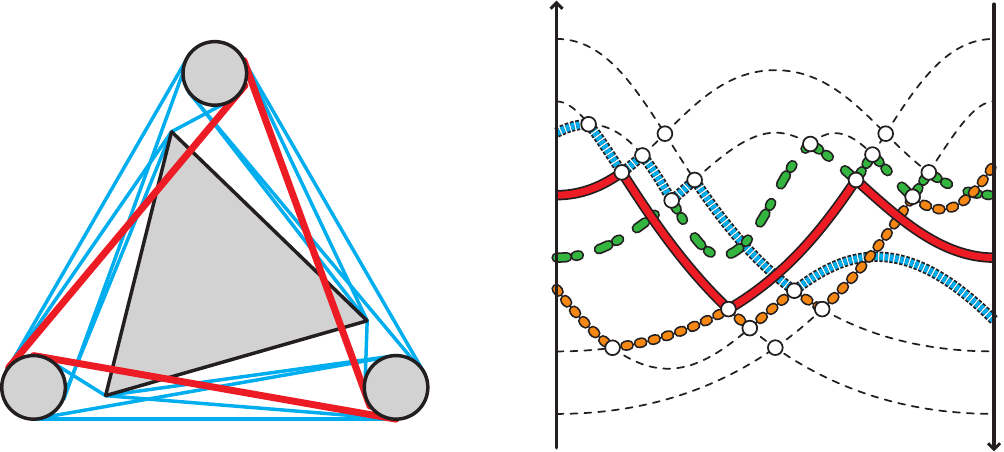}}
	\caption[A \pt{2} of the double pseudoline arrangement of \fref{mpt:fig:dplactrexm1}]{A \pt{2} of the double pseudoline arrangement of \fref{mpt:fig:dplactrexm1}. Observe that the bolded red star has only~$3$ corners.}
	\label{mpt:fig:dplactrexm2}
\end{figure}

For general~$k$, observe that contrary to the case of pseudoline arrangements, a star may have~$2k-1$ corners (see \fref{mpt:fig:dplactrexm2}).

\svs
Let us now give an analogue of Proposition~\ref{mpt:prop:decomposition}. For any point~$q$ in the plane, we denote by~$\eta^k(q)$ the number of crossings between~$q^*$ and the support of~$Q^*$ minus its first~$k$ levels. Let~$\delta^k(q) \eqdef \eta^k(q)/2-|Q|+k$. For any~$\lambda\in\Lambda(U)$ and any point~$q$ in the plane, we still denote by~$\sigma_\lambda(q)$ the \defn{winding number}\index{winding number (of a star)} of~$S(\lambda)$ around~$q$.

\begin{proposition}\label{mpt:prop:decompositionconv}
For any point~$q$ of the plane~$\delta^k(q)=\sum_{\lambda\in\Lambda} \sigma_\lambda(q)$.
\end{proposition}

\begin{proof}
Remember that if~$\tau_\lambda(q)$ denotes the number of intersection points between~$q^*$ and~$\lambda$, then~$\sigma_\lambda(q)=(\tau_\lambda(q)-1)/2$. Thus, we have
$$\eta^k(q)=\sum_{\lambda\in\Lambda} \tau_\lambda(q)=|\Lambda|+2\sum_{\lambda\in\Lambda}\sigma_\lambda(q),$$
\nobreak and we get the result since~$|\Lambda|=2|Q|-2k$.
\end{proof}

When~$k=1$, it is easy to see that~$\delta^1(q)$ is~$1$ if~$q$ is inside the \defn{free space} of the convex hull of~$Q$ (\ie in the convex hull of~$Q$, but not in~$Q$), and~$0$ otherwise. Remember that a \defn{pseudotriangulation}\index{pseudotriangulation} of~$Q$ is a pointed set of bitangents that decomposes the free space of the convex hull of~$Q$ into pseudotriangles~\cite{pv-ptta-96}. Propositions~\ref{mpt:prop:cornersconv} and ~\ref{mpt:prop:decompositionconv} provide, when~$k=1$, the following analogue of Theorem~\ref{mpt:theo:dualitypt}:

\begin{theorem}\label{mpt:theo:dualityptconv}
Let~$Q$ be a set of disjoint convex bodies (in general position) and~$Q^*$ denote its dual arrangement. Then:
\begin{enumerate}[(i)]
\item The dual pseudoline arrangement~$T^* \eqdef \ens{\Delta^*}{\Delta\text{ pseudotriangle of }T}$ of a pseudotriangulation~$T$ of~$Q$ is a \pt{1} of~$Q^*$.
\item  primal set of edges~$E \eqdef \ens{[p,q]}{p,q\in P,\; p^*\wedge q^*\text{ is not a crossing point of } \Lambda}$ of a \pt{1}~$\Lambda$ of~$Q^*$ is a pseudotriangulation of~$Q$.
\end{enumerate}
\end{theorem}

Observe that at least two other arguments are possible to prove~(ii):
\begin{enumerate}
\item either comparing the degrees of the flip graphs as in our first proof of Theorem~\ref{mpt:theo:dualitypt};
\item or checking that all forbidden configurations of the primal (two crossing bitangents, a non-pointed sharp vertex, a non-free bitangent) may not appear in the dual, as in our second proof of Theorem~\ref{mpt:theo:dualitypt}.
\end{enumerate}

\begin{figure}
	\capstart
	\centerline{\includegraphics[scale=1]{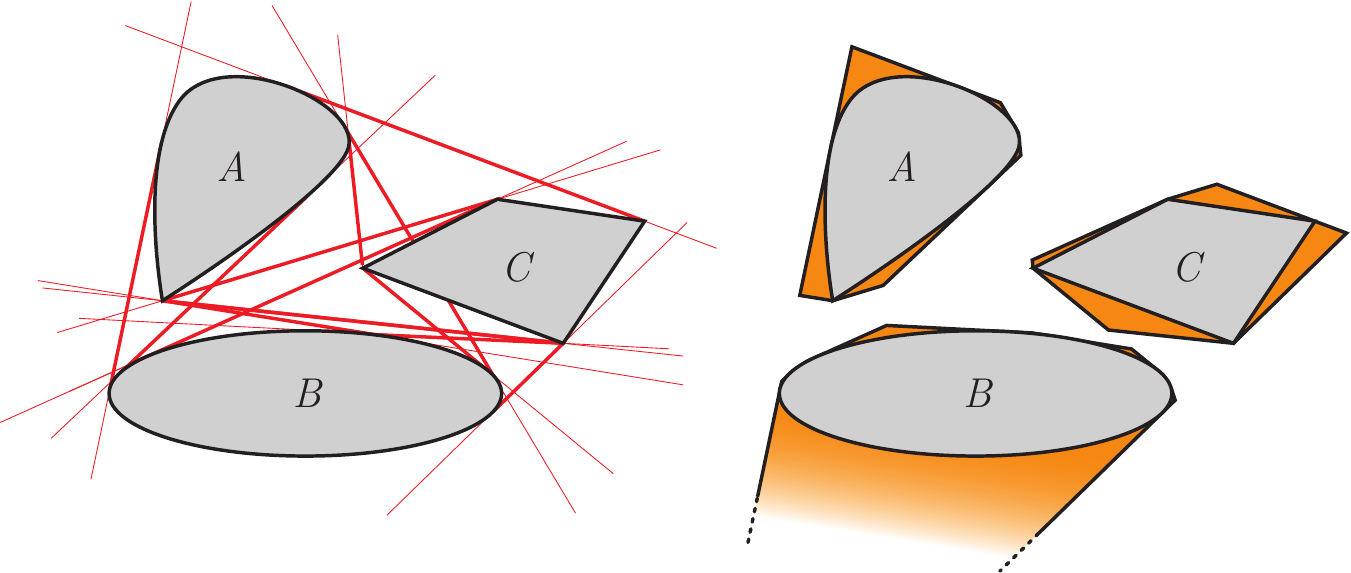}}
	\caption[The set of all bitangents to the arrangement of convex bodies of \fref{mpt:fig:dpla} and the corresponding maximal convex bodies]{The set of all bitangents to the arrangement of convex bodies of \fref{mpt:fig:dpla} and the corresponding maximal convex bodies.}
	\label{mpt:fig:maxconv}
\end{figure}

\svs
\index{depth@\kdepth{k}}
Let us finish this discussion about stars by interpreting the number~$\delta^k(q)$ for general~$k$ and for ``almost every'' point~$q$. For any convex body~$C$ of~$Q$, let~$\nabla C$ denote the intersection of all closed half-planes delimited by a bitangent between two convex bodies of~$Q$, and containing~$C$ (see \fref{mpt:fig:maxconv}). By definition, the bitangents between two convex bodies~$C$ and~$C'$ of~$Q$ coincide with the bitangent of~$\nabla C$ and~$\nabla C'$. Furthermore, the convex bodies~$\nabla C$~($C\in Q$) are maximal for this property. We denote~$\nabla Q \eqdef \ens{\nabla C}{C\in Q}$ the set of maximal convex bodies of~$Q$.

For a point~$q$ outside~$\nabla Q$, the interpretation of~$\delta^k(q)$ is similar to the case of points. We call \defn{level} of a bitangent~$\tau$ the level of the corresponding crossing point in~$Q^*$. Given a point~$q$ outside~$\nabla Q$, the number~$\delta^k(q)$ is the number of bitangents of level~$k$ crossed by any (generic continuous) path from~$q$ to the external face (in the complement of~$\nabla Q$), counted positively when passing from the ``big'' side to the ``small side'', and negatively otherwise.

	\addtocontents{toc}{\vspace*{.2cm}}
	\chapter{Three open problems}\label{chap:multiassociahedron}

In this chapter, we discuss three further topics on multitriangulations:
\begin{enumerate}
\item \defn{Dyck multipaths}: \ktri{k}s of the \gon{n} are counted by the same Hankel determinant of Catalan numbers which counts certain families of non-crossing Dyck paths~\cite{j-gt-03,j-gtdfssp-05}. We investigate explicit bijections between these two combinatorial sets.
\item \defn{Rigidity}: The number of edges in a \ktri{k} is exactly that of a minimally generically rigid graph in a \dimensional{2k} space. We conjecture that \ktri{k}s are minimally rigid graphs in dimension~$2k$, and prove it for~$k=2$.
\item \defn{Multiassociahedron}: The simplicial complex of all \kcross{(k+1)}-free sets of \krel{k} edges of the \gon{n} is known to be a combinatorial sphere~\cite{j-gt-03} and conjectured to be polytopal (as happens for~$k=1$ where it is the polar of the associahedron~\cite{l-atg-89,bfs-ccsp-90,gkz-drmd-94,l-rsp-04,hl-rac-07}). We realize the first non-trivial case (namely,~$k=2$ and~$n=8$) and discuss possible generalizations of two constructions of the associahedron (the first one is the secondary polytope of the regular \gon{n}~\cite{gkz-drmd-94} and the second one is Loday's construction~\cite{l-rsp-04}).
\end{enumerate}
Throughout the chapter, the goal is to present natural and promising ideas based on stars in multitriangulations, although they provide at the present time only partial results to these problems.


\section{Catalan numbers and Dyck multipaths}\label{ft:sec:Dyckpaths}

\index{Catalan number|hbf}
Let~$C_m \eqdef \frac{1}{m+1}{ 2m \choose m}$ denote the $m$th \defn{Catalan number}. Catalan numbers count several combinatorial families (thus called \defn{Catalan families}) such as triangulations, rooted binary trees, Dyck paths, \etc~---~see~\cite{s-ec-99,s-add} for an exhaustive list. These families all share the same recursive structure: a Catalan object of size~$p$ can be decomposed into an object of size~$q$ and another one of size~$p-q-1$, for some~$1\le q\le p-1$. This recursive structure is apparent in the equation~$f(x)=1+xf(x)^2$ on the generating function~$f(x) \eqdef \sum_{m\in\N} C_mx^m$ of a Catalan family, and it inductively defines natural bijections between Catalan families.

The goal of this section is to discuss enumerative and bijective combinatorics of multitriangulations. We give in particular a brief overview of the existing literature on this quesiton.


\subsection{Enumeration of multitriangulations and Dyck multipaths}\label{ft:subsec:Dyckpaths:enumeration}

As mentioned in the introduction, Jakob~Jonsson proved in~\cite{j-gt-03,j-gtdfssp-05} that:

\begin{theorem}[\cite{j-gt-03,j-gtdfssp-05}]\label{ft:theo:enumeration}
The number~$\theta(n,k)$ of \ktri{k}s of the \gon{n} equals:
$$\theta(n,k)=\det(C_{n-i-j})_{1\le i,j\le k}=\det\begin{pmatrix} C_{n-2} & C_{n-3} & \edots & \edots & C_{n-k-1} \\ C_{n-3} & \edots & \edots & C_{n-k-1} & \edots \\ \edots & \edots & \edots & \edots & \edots \\ \edots & C_{n-k-1} & \edots & \edots & C_{n-2k+1} \\ C_{n-k-1} & \edots & \edots & C_{n-2k+1} & C_{n-2k} \end{pmatrix}.$$
\end{theorem}
\vspace*{-.8cm}\qed

\begin{remark}\label{ft:rem:enumeration}
This Hankel determinant can be expressed by the explicit product formula~\cite{v-ciqda-00}:
$$\theta(n,k)=\prod_{1\le i<j\le n-2k-1} \frac{i+j+2k}{i+j}.$$
For the asymptotics of it, observe that the recurrence~$C_n=\frac{4n-2}{n+1}C_{n-1}$ makes each entry of the matrix equal to~$C_n$ times a rational function of degree at most~$2k$ in~$n$. Since~$C_n \in \Theta(4^n n^{-3/2})$ we conclude that, for fixed~$k$, the number of~\ktri{k}s of an \gon{n} grows as~$4^{nk}$, modulo a rational function of degree~$O(k^2)$ in~$n$.
\end{remark}

\begin{figure}[!h]
	\capstart
	\centerline{\includegraphics[scale=1]{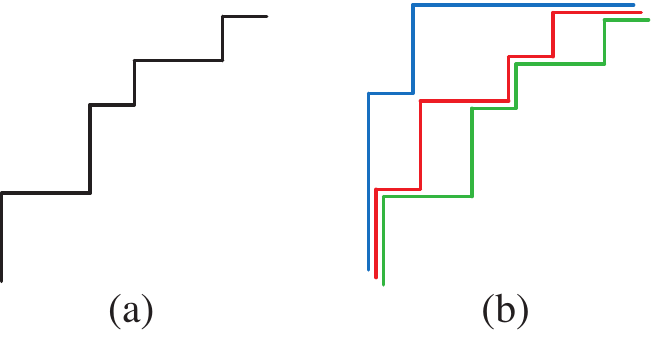}}
	\caption[Dyck paths and Dyck multipaths]{(a) A Dyck path $NNEENNENEENE$ and (b) a Dyck \kpath{3} of semilength~$6$.}
	\label{ft:dyck}
\end{figure}

On the other hand, the determinant in Theorem~\ref{ft:theo:enumeration} also counts another class of objects:

\begin{definition}
\index{Dyck path|hbf}
\index{Dyck path@Dyck \kpath{k}}
A \defn{Dyck path} of \defn{semilength}~$\ell$ is a lattice path using north steps~$N \eqdef (0,1)$ and east steps~$E \eqdef (1,0)$ starting from~$(0,0)$ and ending at~$(\ell,\ell)$, and such that it never goes below the diagonal~$y=x$~---~see \fref{ft:dyck}(a). We call \defn{Dyck \kpath{k}} of \defn{semilength}~$\ell$ any \tuple{k}~$(D_1,\dots,D_k)$ of Dyck paths of semilength~$\ell$ such that each $D_i$ never goes above $D_{i-1}$, for~$2\le i\le k$~---~see \fref{ft:dyck}(b).
\end{definition}

That the Catalan determinant of Theorem~\ref{ft:theo:enumeration} also counts the number of Dyck $k$-paths of semilength $n-2k$ is an almost direct application of Lindstr\"om-Gessel-Viennot technique~\cite{l-vrim-73,gv-bdphlf-85,dcv-ecytbh-86} to count non-intersecting lattice paths. The equality of their cardinality raises the question to find an explicit bijection between the set~$T_{n,k}$ of \ktri{k}s of the \gon{n} and the set~$D_{n,k}$ of Dyck \kpath{k}s of semilength~$n-2k$. This question was stated by Jakob Jonsson in~\cite{j-gt-03,j-gtdfssp-05} and studied in different articles~\cite{e-btdp-07,n-abtdp-09,k-gdidcffs-06,r-idsfmp-07}. Finding an explicit bijection between these two families would provide a simpler proof of the enumeration formula of Theorem~\ref{ft:theo:enumeration}, and may also reveal properties hidden on one family but apparent on the other one.


\subsection{Overview of existing proofs}\label{ft:subsec:Dyckpaths:existingproofs}

In this section, we give a short overview of different approaches to Theorem~\ref{ft:theo:enumeration}. This theorem, conjectured by Andreas~Dress, Jacobus~Koolen, Vincent~Moulton and Volkmar~Welker, was first proved by Jakob Jonsson in~\cite{j-gt-03} with a direct but ``complicated'' method (according to its author). This unpublished manuscript also contains a discussion on reasonable properties that should have an explicit bijection between~$T_{n,k}$ and~$D_{n,k}$ (we develop this discussion later in Section~\ref{ft:subsec:Dyckpaths:failedbijection}). Afterwards, the topic followed two different directions: on the one hand, it was replaced in the more general context of the study of fillings of polyominoes with restricted length of increasing and decreasing chains~\cite{j-gtdfssp-05, k-gdidcffs-06, r-idsfmp-07}; on the other hand, two different bijections between \ktri{2}s and \kpath{2}s recently appeared~\cite{e-btdp-07,n-abtdp-09}.

\subsubsection{Fillings of polyominoes}

\begin{figure}[!h]
	\capstart
	\centerline{\includegraphics[scale=1]{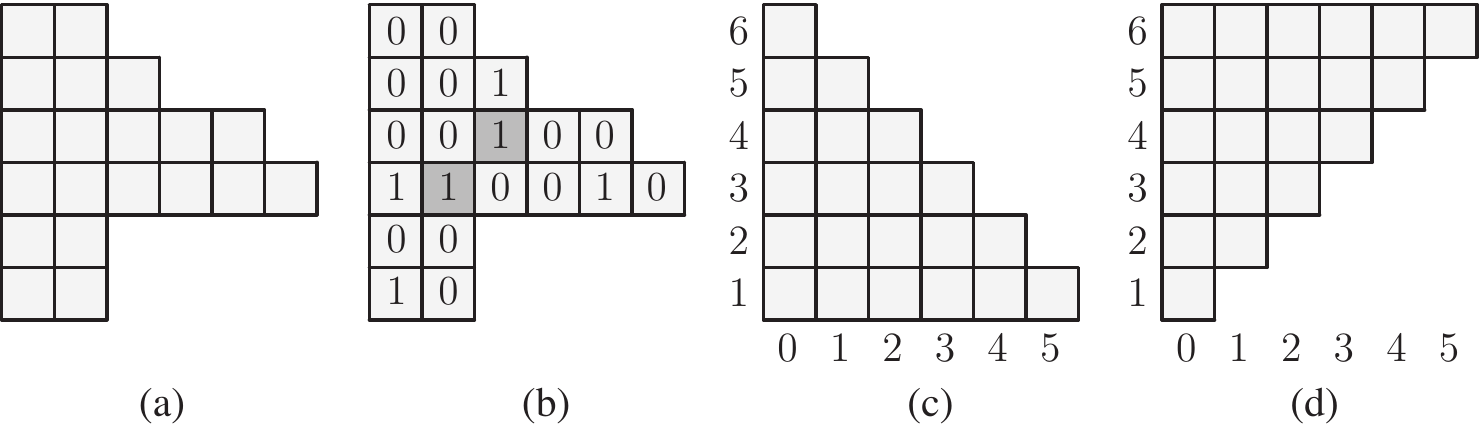}}
	\caption[Stack polyominoes]{(a) A stack polyomino of shape~$(2,2,6,5,3,2)$; (b) a filling of this polyomino, with a shaded \kdiag{2}; (c) the polyomino~$A_7$; and (d) the polyomino~$B_7$.}
	\label{ft:fig:polyominoes1}
\end{figure}

\index{polyomino|hbf}
In~\cite{j-gtdfssp-05}, Jakob Jonsson gave an alternative simpler (but this time indirect) proof of Theorem~\ref{ft:theo:enumeration}, using the following polyominoes (see \fref{ft:fig:polyominoes1} and~\cite{j-gtdfssp-05, k-gdidcffs-06, r-idsfmp-07} for~\mbox{details}):

\begin{definition}
A \defn{stack polyomino} is a subset~$\ens{(i,j)}{0\le i\le s_j,\; j\in[r]}$ of~$\Z^2$, where there exists~$t\in[r]$ such that~$0\le s_1 \le\dots\le s_t$ and~$s_t\ge\dots\ge s_r\ge0$. The \tuple{r}~$(s_1,\dots,s_r)$ is the \defn{shape} of the polyomino, while the multiset~$\{s_1,\dots,s_r\}$ is its \defn{content}. A \defn{filling} of a polyomino is an assignment of~$0$ and~$1$ to its boxes. For~$\ell\in\N$, an~\defn{\kdiag{\ell}} of a filling of a polyomino~$\Lambda$ is a chain~$(\alpha_1,\beta_1),\dots,(\alpha_\ell,\beta_\ell)$ of~$\ell$ boxes of~$\Lambda$ filled with~$1$ and such that:
\begin{itemize}
\item $\alpha_1<\alpha_2<\dots<\alpha_\ell$ and~$\beta_1<\beta_2<\dots<\beta_\ell$, and
\item the rectangle~$\ens{(\alpha,\beta)}{\alpha_1\le\alpha\le\alpha_\ell \text{ and } \beta_1\le\beta\le\beta_\ell}$ is a subset of~$\Lambda$.
\end{itemize}
\end{definition}

Two special polyominoes are particularly interesting for us:
\begin{enumerate}[(i)]
\item It was proved in~\cite{ht-gbmdpi-92} that the Hankel determinant of Catalan numbers of Theorem~\ref{ft:theo:enumeration} also counts the number of maximal \kdiag{(k+1)}-free fillings of the triangular polyomino $A_n \eqdef \ens{(i,j)}{0\le i\le n-j-1,\; j\in[n-1]}$ (see \fref{ft:fig:polyominoes1}(c)).
\item The adjacency matrix of any subset~$E$ of~$E_n$ defines a filling of the triangular polyomino $B_n \eqdef \ens{(i,j)}{0\le i<j\le n-1}$ (see \fref{ft:fig:polyominoes1}(d)) in which \kdiag{\ell}s correspond to \kcross{\ell}s of~$E$. In particular, this gives a straightforward bijection between \ktri{k}s of an \gon{n} and maximal \kdiag{(k+1)}-free fillings of~$B_n$.
\end{enumerate}
According to these two examples, the proof of Theorem~\ref{ft:theo:enumeration} is a direct consequence of:

\begin{theorem}[\cite{j-gtdfssp-05}]\label{ft:theo:jonsson}
The number of maximal \kdiag{(k+1)}-free fillings of a stack polyomino does not depend on its shape, but only on its content.\qed
\end{theorem}

Enumeration properties of fillings of polyominoes with restrictions on the length of the longest increasing and decreasing chains have been studied in other related articles. In particular, Christian Krattenthaler~\cite{k-gdidcffs-06} used growth diagrams to study symmetries of crossings and nestings in fillings of Ferrers shapes (generalizing results of~\cite{cddsy-cnmp-07}). Following this method, Martin Rubey~\cite{r-idsfmp-07} extended Theorem~\ref{ft:theo:jonsson} to more general moon polyominoes (using a very different method from the original proof). He also emphasized that ``the problem of finding a completely bijective proof [...] remains open. However, it appears that this problem is difficult''.

\subsubsection{Two inductive bijections when~$k=2$}

In recent papers, Sergi Elizalde~\cite{e-btdp-07} and Carlos Nicolas~\cite{n-abtdp-09,nicolas-phd} both obtained explicit bijections between \ktri{2}s and Dyck \kpath{2}s. Although the details are quite different, these two papers are both based on induction. The idea is to construct two isomorphic generating trees for the sets~$T_{n,2}$ and~$D_{n,2}$, that is to say, two parallel constructions which generate on the one hand~$T_{n+1,2}$ from~$T_{n,2}$, and on the other hand~$D_{n+1,2}$ from~$D_{n,2}$. Of course, in both articles, the generating tree for \ktri{2}s is based on the flattening \& inflating operation developed in Section~\ref{stars:sec:flat-infl} (this operation is very implicit in~\cite{e-btdp-07}, but completely explicit in~\cite{n-abtdp-09,nicolas-phd}). The two bijections essentially differ in the way they use this flattening \& inflating operation:
\begin{enumerate}
\item In~\cite{e-btdp-07}, Sergi Elizalde starts from a \ktri{2} and inductively flattens the \kstar{2} which contains the last \kear{2}. During this inductive process, the \krel{2} edges of the \ktri{2} are colored with two different colors. Then, each color is used to define one of the Dyck paths (see Section~\ref{ft:subsec:Dyckpaths:failedbijection}).
\item In~\cite{n-abtdp-09,nicolas-phd}, Carlos Nicolas defines a generalized notion of \defn{returns} in a Dyck multipath. There is a parallelism between returns in Dyck \kpath{k}s and \kcross{k}s incident to $\{0,\dots,k-1\}$ in \ktri{k}s: similarly to the inflating operation of Section~\ref{stars:sec:flat-infl}, a return in a Dyck \kpath{k} of semilength~$\ell$ can be inflated to obtain a Dyck \kpath{k} of semilength~$\ell+1$. This parallelism is exploited to obtain two isomorphic generating trees for \ktri{2}s and Dyck \kpath{2}s. By construction, the resulting bijection has the nice property to identify returns with \kcross{k}s incident to~$\{0,\dots,k-1\}$.
\end{enumerate}
Unfortunately, it seems that none of these two bijections can be directly extended to general~$k$.


\subsection{Edge colorings and indegree sequences}\label{ft:subsec:Dyckpaths:failedbijection}

We now present a relevant attempt to define an explicit and direct (meaning non-inductive) bijection between \ktri{k}s and Dyck \kpath{k}s. The starting point of our definition is the following well-known bijection between triangulations and Dyck paths:

\begin{proposition}\label{ft:prop:indegreeTriangulations}
The application~$d:T\mapsto d(T) \eqdef N^{\delta_0(T)}EN^{\delta_1(T)}E\dots N^{\delta_{n-3}(T)}E$, where~$\delta_i(T)$ is the number of triangles of~$T$ whose first vertex is~$i$, is a bijection between triangulations of the \gon{n} and Dyck paths of semilength~$n-2$.\qed
\end{proposition}

In other words, the sequence of north powers of the Dyck path~$d(T)$ is the sequence of indegrees of the vertices of the \gon{n} in the oriented graph whose edges are the interior edges of~$T$ together with the top edge~$[0,n-1]$, all oriented towards their smallest vertex. 

Observe that this bijection
\begin{enumerate}[(i)]
\item respects triangles: each triangle of~$T$ corresponds to a pair of $NE$ steps in~$d(T)$. Namely, the Dyck path~$d(T)$ has a north step at the first vertex of each triangle and an east step at the second vertex of each triangle.
\item is compatible with the natural partial orders: if~$T$ and~$T'$ are two triangulations of the \gon{n}, with $T<T'$ (\ie such that~$T'$ can be obtained from~$T$ by a sequence of slope-increasing flips~---~see Section~\ref{stars:sec:flips}), then~$d(T)$ is above~$d(T')$.
\end{enumerate}

\svs
Generalizing this simple bijection, Jakob Jonsson studied in~\cite{j-gt-03} the repartition of indegree sequences in multitriangulations. He defined the two following parallel notions:

\begin{definition}\label{ft:def:signature}
Given a \ktri{k}~$T$ of the \gon{n}, consider the oriented graph whose edges are the \krel{k} edges of~$T$ together with the the top \kbound{k} edges of the \gon{n} $[0,n-k],$ $[1,n-k+1],\dots,$ $[k-1,n-1]$, all oriented towards their smallest vertex. The \defn{indegree sequence} of~$T$ is the sequence of indegrees of the first~$n-k-1$ vertices of the \gon{n} in this oriented graph.

The \defn{signature} of a Dyck \kpath{k}~$(N^{\nu_0^1}EN^{\nu_1^1}E\dots N^{\nu_{\ell-1}^1}E,\dots,N^{\nu_0^k}EN^{\nu_1^k}E\dots N^{\nu_{\ell-1}^k}E)$ of semi\-length~$\ell$ is the sequence~$(\sum_{i=1}^k \nu_{j-i}^i)_{j\in[\ell+k-1]}$ (with the convention~$\nu_p^q \eqdef 0$ if~$p<0$ or~$p\ge \ell$).
\end{definition}

\begin{figure}
	\capstart
	\centerline{\includegraphics[scale=1]{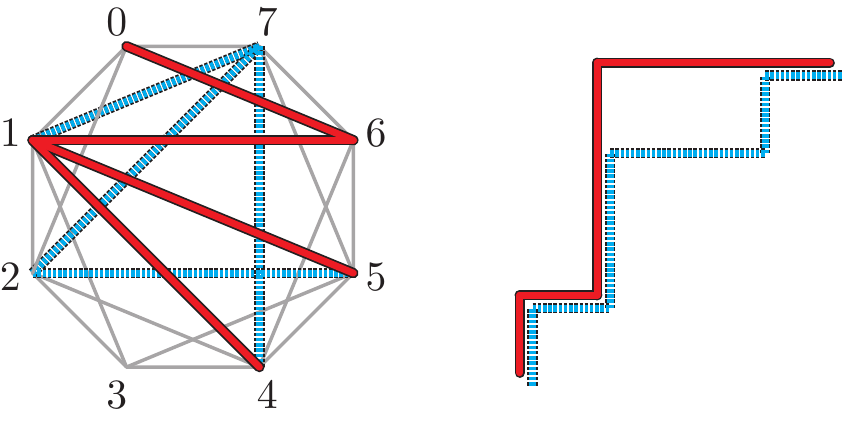}}
	\caption[The Dyck \kpath{2}~$D(T)$ associated to the \ktri{2} of \fref{intro:fig:2triang8points}]{The Dyck \kpath{2}~$D(T)$ associated to the \ktri{2}~$T$ of \fref{intro:fig:2triang8points}.}
	\label{ft:fig:2triang8pointsbijection}
\end{figure}

\begin{example}
In \fref{ft:fig:2triang8pointsbijection}, the indegree sequence of the \ktri{2} and the signature of the Dyck \kpath{2} both equal~$(1,4,2,0,1)$.
\end{example}

The following theorem refines Theorem~\ref{ft:theo:enumeration}, comparing indegree sequences in multitriangulations \vs signatures in Dyck multipaths:

\begin{theorem}[\cite{j-gt-03}]\label{ft:theo:indegree-signature}
For any~$k\ge 1$,~$n\ge 2k+1$ and any sequence~$\delta \eqdef (\delta_1,\dots,\delta_{n-k-1})$, the number of \ktri{k}s of the \gon{n} with indegree sequence~$\delta$ equals the number of Dyck \kpath{k}s of semilength~$n-2k$ with signature~$\delta$.\qed
\end{theorem}

It is reasonable to look for a bijection between multitriangulations and Dyck multipaths which respects these parameters. In other words, the indegree sequence of a \ktri{k} already tells us the signature that should have its associated Dyck \kpath{k}, and it ``only remains'' to separate this signature into~$k$ distinct north power sequences.

\begin{figure}[!h]
	\capstart
	\centerline{\includegraphics[scale=1]{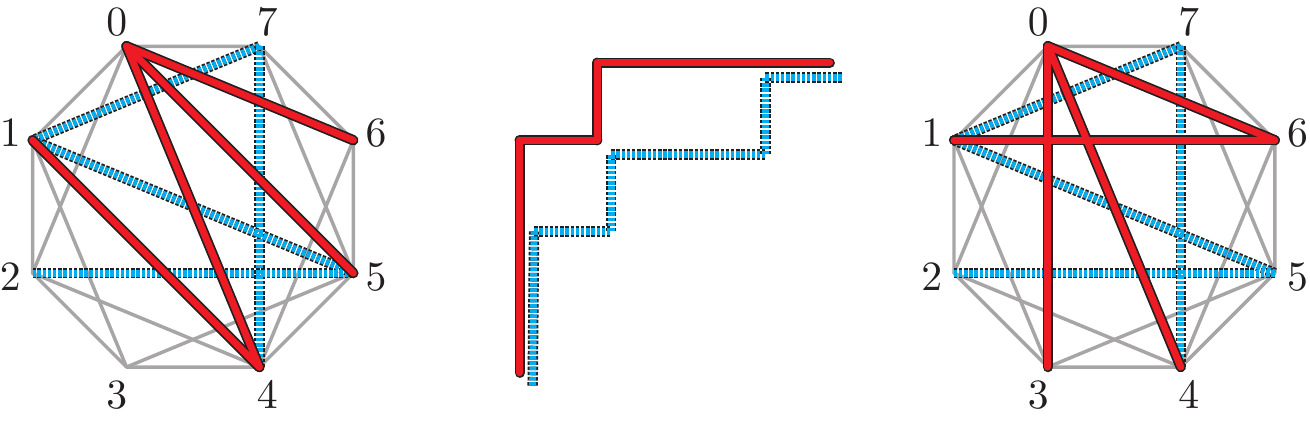}}
	\caption[Two \ktri{2}s of the octagon with the same associated Dyck \kpath{2}]{Two \ktri{2}s of the octagon with the same associated Dyck \kpath{2}.}
	\label{ft:fig:bijectionproblem}
\end{figure}

\svs
Using this heuristic, we attempt to generalize the function~$d$ of Proposition~\ref{ft:prop:indegreeTriangulations} in a seemingly natural fashion. Unfortunately, the generalized function $D$ that we construct is not bijective; but it is surprisingly close to a bijection for small cases (see \fref{ft:fig:bijectionproblem} and the discussion above), which makes it plausible that a clever modification of it makes it a bijection.

We consider the application~$D$, which associates to a \ktri{k} of the \gon{n} the \tuple{k} $D(T) \eqdef (D_1(T),\dots,D_k(T))$ of lattice paths defined by
$$D_j(T) \eqdef N^{\delta_{j-1}^j(T)}EN^{\delta_{j}^j(T)}E\dots N^{\delta_{j+n-2k-2}^j(T)}E,$$
where $\delta_i^j(T)$ denotes the number of \kstar{k}s of~$T$ whose~$j$th vertex is the vertex~$i$. See \fref{ft:fig:2triang8pointsbijection} for an example.

This application can also be interpreted in terms of colored indegrees of the vertices of the \gon{n}. For~$j\in[k]$, and for any \kstar{k}~$S$ of~$T$, with vertices~$0\cle s_1\cl s_2\cl\cdots\cl s_{2k+1}\cle n-1$, color the edge~$[s_j,s_{k+j+1}]$ with color~$j$. This coloring decomposes the set of all \krel{k} edges of~$T$ together with the top \kbound{k} edges~$[0,n-k],\dots,[k-1,n-1]$ into~$k$ disjoint sets~$\tau_1,\dots,\tau_k$ each containing~$n-2k$ edges. The sequence of north powers of the $j$th Dyck path~$D_j(T)$ is the sequence of indegrees of the vertices of the \gon{n} in the oriented graph~$\tau_j$ (where each edge is oriented towards its smallest vertex), after $j-1$ shifts.

\svs
The two following lemmas tend to show that the application~$D$ is a reasonable approximation of the bijection we are looking for.

\begin{lemma}
For any \ktri{k} of the \gon{n}, the \tuple{k} of lattice paths~$D(T)$ is a Dyck \kpath{k} of semilength~$n-2k$.
\end{lemma}

\begin{proof}
We first prove that each~$D_j(T)$ is a Dyck path of semilength~$n-2k$. Observe that for each~$i\in[n-2k]$, there exists a unique \kstar{k}~$S_i$ whose $(k+1)$th vertex is~$k-1+i$ (indeed, it is the unique \kstar{k} bisected by the line passing through the vertex~$k-1+i$ and through the midpoint of~$[0,n-1]$). The $j$th vertex of this \kstar{k}~$S_i$ is certainly smaller or equal~$i+j-2$. Thus, the number of north steps in~$D_j(T)$ before the~$i$th east step (which is the number of \kstar{k}s of~$T$ whose~$j$th vertex is smaller or equal~$i+j-2$) is at least~$i$. Since the total numbers of north steps and of east steps both equal the number~$n-2k$ of \kstar{k}s of~$T$, the lattice path~$D_j(T)$ is indeed a Dyck path of semilength~$n-2k$.

Finally, in each \kstar{k} of~$T$, the $(j-1)$th vertex is certainly before the $j$th vertex, which implies that for all~$i$, the number of \kstar{k}s of~$T$ whose~$(j-1)$th vertex is smaller or equal~$i+j-3$ is at least the number the number of \kstar{k}s of~$T$ whose~$j$th vertex is smaller or equal~$i+j-2$. Consequently, the Dyck path~$D_j(T)$ never goes above the Dyck path~$D_{j-1}(T)$ for all~$j$, and~$D(T)$ is indeed a Dyck \kpath{k}.
\end{proof}

By construction, the function~$D$ respects \kstar{k}s and preserves the parameters of Theorem~\ref{ft:theo:indegree-signature}. Furthermore, it is compatible with the following partial orders:
\begin{enumerate}[(i)]
\item For two \ktri{k}s~$T$ and~$T'$ we write $T<T'$ if~$T'$ can be obtained from~$T$ by a sequence of slope-increasing flips~---~see Section~\ref{stars:sec:flips}.
\item For two Dyck \kpath{k}s~$D \eqdef (D_1,\dots,D_k)$ and~$D' \eqdef (D_1',\dots,D_k')$ we write~$D<D'$ if~$D_1$ is above~$D_1'$, or $D_1=D_1'$ and~$D_2$ is above~$D_2'$, or \etc{} (lexicographic order for the ``above'' relation).
\end{enumerate}

\begin{lemma}
The application~$D$ is increasing with respect to these partial orders: for any two \ktri{k}s~$T$ and~$T'$ of the \gon{n}, $T<T'$ implies~${D(T)<D(T')}$.
\end{lemma}

\begin{proof}
It is sufficient to prove that a single slope-increasing flip in a multitriangulation also increases the associated Dyck multipath (for the order~$<$ on Dyck multipaths).

Consider a slope-increasing flip from a \ktri{k}~$T$ to another one~$T'$. Let~$a<c$ denote the vertices of the edge in~$T\ssm T'$ and~$b<d$ the vertices of the edge in~$T'\ssm T$. Since the flip is slope-increasing, we have~$0\cle a\cl b\cl c\cl d\cle n-1$. Let~$R$ (resp.~$S$) be the \kstar{k} of~$T$ adjacent to~$[a,c]$ and containing the vertex~$b$ (resp.~$d$). Similarly, let~$X$ (resp.~$Y$) be the \kstar{k} of~$T'$ adjacent to~$[b,d]$ and containing the vertex~$a$ (resp.~$c$). See \fref{ft:fig:flipincr} for an example.

\begin{figure}[b]
	\capstart
	\centerline{\includegraphics[scale=1]{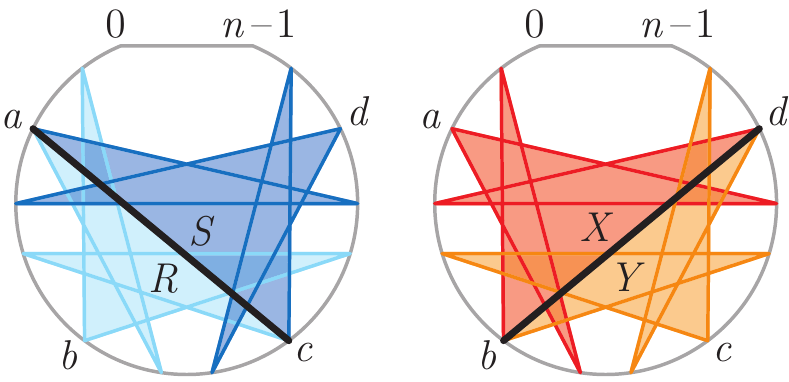}}
	\caption[An increasing flip]{An increasing flip.}
	\label{ft:fig:flipincr}
\end{figure}

It is easy to check that, before vertex~$a$, the vertices of~$X$ are exactly the vertices of~$R$ while the vertices of~$Y$ are exactly the vertices of~$S$. Furthermore, we have more vertices of~$R$ than vertices of~$S$ before~$a$. Consequently, if~$a$ is the $j$th vertex of~$S$, then:
\begin{enumerate}[(i)]
\item for all~$i<j$ we have~$D_i(T)=D_i(T')$, since the $i$th vertex of~$R$ is also the $i$th vertex of~$X$ and the $i$th vertex of~$S$ is also the $i$th vertex of~$Y$; and
\item $D_j(T)$ is strictly above~$D_j(T')$, since the $j$th vertex of~$R$ is also the $j$th vertex of~$X$, but the $j$th vertex of~$S$ (\ie $a$) is strictly before the $j$th vertex of~$Y$ (because $a\notin Y$).
\end{enumerate}
This proves that~$D(T)<D(T')$, by definition of the lexicographic order.
\end{proof}

Despite these promising properties, it turns out that~$D$ is not a bijection between~$T_{n,k}$ and~$D_{n,k}$, but it is surprisingly close to be:
\begin{enumerate}[(i)]
\item For~$k=2$ and~$n\le 7$, the function~$D$ is a bijection.
\item For~$k=2$ and~$n=8$ it is not, but there is a unique pair of \ktri{2}s producing the same Dyck \kpath{2}. They are shown in \fref{ft:fig:bijectionproblem}.
\end{enumerate}

We have tried other possible functions from~$T_{n,k}$ to~$D_{n,k}$, all obtained as a generalization of a well-known bijection between~$T_{n,1}$ and~$D_{n,1}$, and compatible with \kstar{k}s. It turns out that all of them are equivalent to~$D$ (in the sense that for~$k=2$ and~$n=8$, we always obtain twice the Dyck \kpath{2} of \fref{ft:fig:bijectionproblem}).

We hope however that \kstar{k}s will help to find a direct general bijection.


\subsection{Beam arrangement}\label{ft:subsec:Dyckpaths:laser}

To finish this section, we construct another lattice path configuration associated to a multitriangulation, namely its dual pseudoline arrangement presented in Chapter~\ref{chap:mpt} that we embed on the integer grid. For this, we exploit again the adjacency matrix of a multitriangulation, drawing lattice paths on its rows and columns according to its entries.

\begin{figure}[!h]
	\capstart
	\centerline{\includegraphics[scale=1]{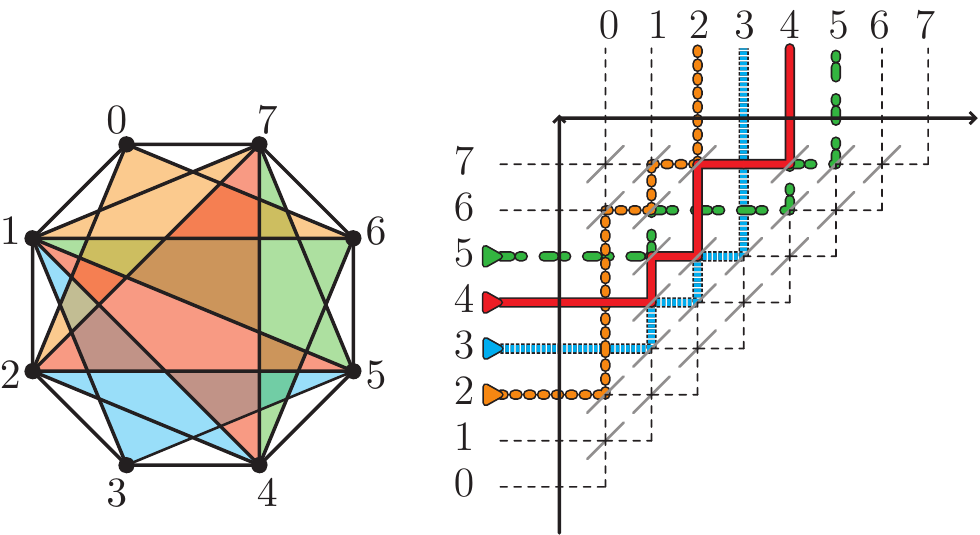}}
	\caption[The beam arrangement of the triangulation of~\fref{intro:fig:2triang8points}]{The beam arrangement of the triangulation of~\fref{intro:fig:2triang8points}.}
	\label{ft:fig:mirrorslasers}
\end{figure}

Let~$T$ be a \ktri{k} of the \gon{n}. For each edge~$(u,v)$ of~$T$ (with~$u<v$), we place a \defn{mirror} at the grid point of coordinates~$(u,v)$. This mirror is a double faced mirror parallel to the diagonal~$x=y$ so that it reflects a ray coming from~$(-\infty,v)$ to a ray going to~$(u,+\infty)$, and a ray coming from~$(u,-\infty)$ to a ray going to~$(+\infty,v)$. Furthermore, for~$1\le i\le n-2k$, we place a \defn{light beam} at~$(-\infty,k-1+i)$ pointing horizontally. We obtain~$n-2k$ beams which reflect on the  mirrors of~$T$~---~see \eg \fref{ft:fig:mirrorslasers}. All throughout this chapter, we will call this ray configuration the \defn{beam arrangement}\index{beam!--- arrangement} associated to~$T$.

Our study of pseudoline arrangements in Chapter~\ref{chap:mpt} implies the following properties:

\begin{proposition}\label{ft:prop:laser}
\begin{enumerate}[(i)]
\item All beams are $x$- and $y$-monotone lattice paths.
\item The $i$th beam comes from~$(-\infty,k-1+i)$ and goes to~$(k-1+i,+\infty)$.
\item Each beam reflects exactly~$2k+1$ times, and thus, has~$k$ vertical segments (plus one vertical half-line) and~$k$ horizontal segments (plus one horizontal half-line).
\item The beams form a pseudoline arrangement: any two of them cross exactly once.\qed
\end{enumerate}
\end{proposition}

Let us recall that the $i$th beam~$B_i$ ``corresponds'' via duality to the \kstar{k}~$S_i$ of~$T$ whose $(k+1)$th vertex is the vertex~$k-1+i$ (in other words, the \kstar{k} bisected by the line passing through the vertex~$k-1+i$ and through the midpoint of~$[0,n-1]$). Indeed:
\begin{enumerate}[(i)]
\item the beam~$B_i$ is (by duality) the set of all bisectors of~$S_i$;
\item the mirrors which reflect~$B_i$ are the edges of~$S_i$;~and
\item the intersection of two beams~$B_i$ and~$B_j$ is the common bisector of~$S_i$ and~$S_j$.
\end{enumerate}
Observe that instead of $k$~Dyck paths of semi-length~$n-2k$, the beam arrangement of a \ktri{k} has~$n-2k$ beams which all have~$k$ horizontal steps. Despite these similarities, we have not found any bijection between Dyck \kpath{k}s and beam arrangements.


\section{Rigidity}\label{ft:sec:rigidity}

Triangulations are rigid graphs in the plane: the only continuous motions of their vertices which preserve their edges' lengths are the isometries of the plane. In this section, we investigate~rigidity properties of multitriangulations. We recall some rigidity notions, and refer to~\cite{g-cf-01,f-gga-04} for very nice introductions to combinatorial rigidity, and to~\cite{gss-cr-93} for a more technical one.


\subsection{Combinatorial rigidity, sparsity and arboricity}\index{rigidity|hbf}

A \dimensional{d} \defn{framework}\index{frameworks@(rigid and flexible) frameworks} is a finite graph together with an embedding of its vertices in~$\R^d$. We think of the straight edges of a framework as rigid bars connected by flexible joints, and we want to distinguish \defn{flexible} frameworks (those that can be deformed preserving their edges' lengths) from \defn{rigid} ones.

\begin{figure}
	\capstart
	\centerline{\includegraphics[width=\textwidth]{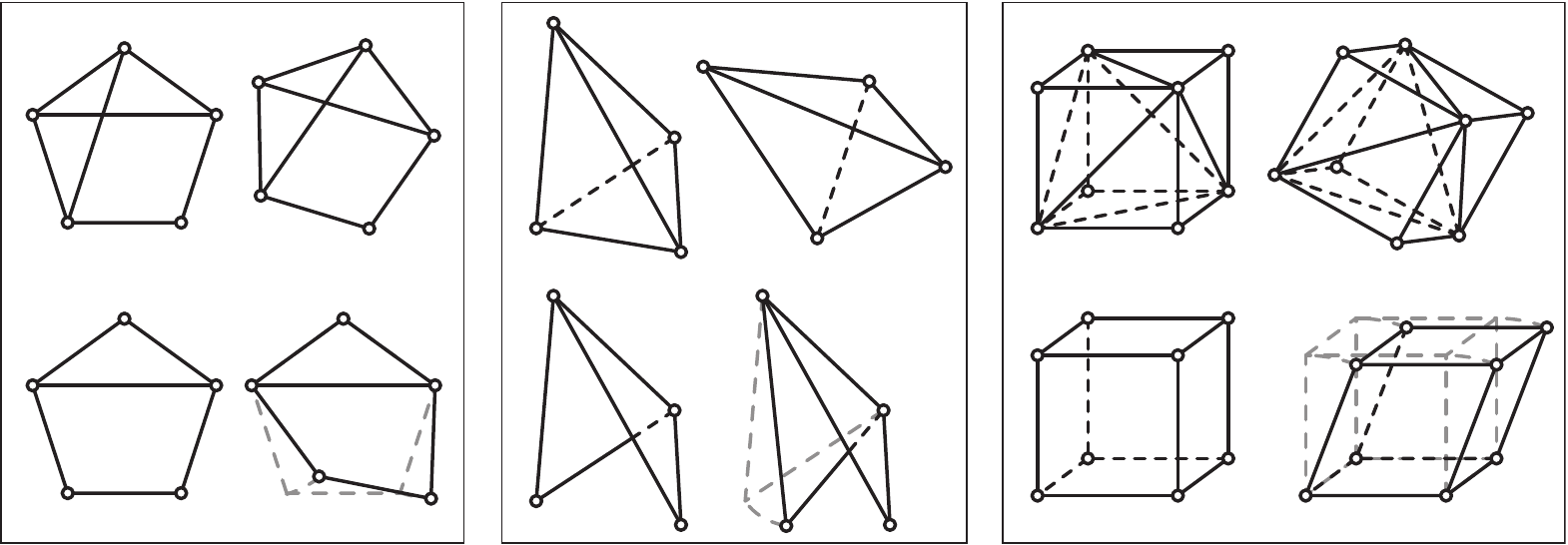}}
	\caption[Motions of frameworks]{Motions of frameworks. Each of the three groups represents a rigid graph (top left) with a rigid motion (top right) and a flexible subgraph (bottom left) with a non-rigid motion (bottom right). The first group lies in the plane while the other two lie in \dimensional{3} space.}
	\label{ft:fig:motions}
\end{figure}

To be more precise, define a \defn{motion}\index{motion} of a framework to be a continuous motion of its vertices which preserves the lengths of all its edges (see \fref{ft:fig:motions}). Of course, any framework admits \defn{rigid motions}\index{motion!rigid ---} which correspond to isometries of~$\R^d$ (\ie in which not only the edges' lengths are preserved, but also all distances between all pairs of vertices). \defn{Rigid} frameworks are those that admit only rigid motions.

In other words, a framework is rigid if the system of quadratic equations whose variables are the coordinates of its vertices and whose equations encode the lengths of its edges has locally no other non-congruent solution. Each edge of the graph is thus encoded by a quadratic equation which normally reduces the dimension of the solution by one (intuitively, one edge removes one ``degree of freedom''). Since the space of geometric embeddings of a graph with~$n$ vertices in~$\R^d$ has dimension~$dn$ while the space of congruent frameworks to a given framework has dimension~${d+1 \choose 2}$, it is reasonable to expect a framework to be rigid as soon as it has $dn-{d+1 \choose 2}$ edges ``properly placed'' (in the sense that no edge is redundant). The reality is a bit more complicated as the following examples show:

\begin{figure}[!h]
	\capstart
	\centerline{\includegraphics[scale=1]{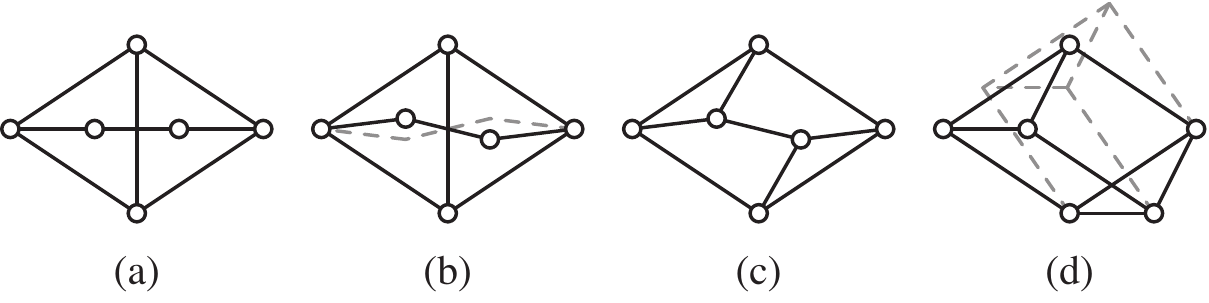}}
	\caption[Four \dimensional{2} frameworks with~$6$ vertices]{Four \dimensional{2} frameworks with~$6$ vertices.}
	\label{ft:fig:nongeneric}
\end{figure}

We would expect~$9$ edges to be required to make a \dimensional{2} framework with~$6$ vertices rigid. However, the behavior of~(a) and~(d) is not the one expected: (a)~is rigid (because the middle horizontal chain is tight) although it has only~$8$ edges and (d) is flexible (because it has three parallel edges) although it has~$9$ edges. These frameworks~(a) and~(d) are \defn{singular} in the sense that there exist arbitrarily close frameworks (for example~(b) and~(c)) which behave differently with respect to rigidity: just break the colinearity in~(a) (like in~(b)) or the parallelism in~(d) (like in (c)). In fact, the \dimensional{d} embeddings of a given finite graph are either almost all rigid or almost all flexible. In the first case, the graph is said to be \defn{generically rigid}\index{rigidity!generic ---} in dimension~$d$. Generic rigidity is a combinatorial property: it only depends on the graph and not of a specific embedding. In the plane, generic rigidity is completely characterized:

\begin{theorem}[\cite{l-grpss-70,h-gsss-11,r-ntarss-84}]\label{ft:theo:rigidity2d}
The following conditions on a graph~$G \eqdef (V,E)$ are equivalent:
\begin{enumerate}
\item $G$~is \defn{minimally generically rigid} in the plane, \ie $G$~is generically rigid and removing any edge from~$G$ produces a generically flexible graph.
\item \defn{Laman's Condition}\index{Laman condition}~\cite{l-grpss-70}: $|E|=2|V|-3$ and any subgraph of~$G$ on~$v$ vertices has at most~$2v-3$ edges.
\item \defn{Henneberg's Constructions}~\cite{h-gsss-11}: $G$~can be constructed from a single edge by a sequence of Henneberg's additions: there exists a sequence $G_0 \eqdef (\{v,w\},\{(v,w)\}),\dots,$ $G_i \eqdef (V_i,E_i),\dots,$ $G_p=G$ such that for any~$i\in[p]$:
\begin{enumerate}
\item either there exist three vertices~$v\in V_i\ssm V_{i-1}$ and~$x,y\in V_{i-1}$ such that $V_i=V_{i-1}\cup\{v\}$ and~$E_i=E_{i-1}\cup\{(v,x),(v,y)\}$;
\item or there exist four vertices~$v\in V_i\ssm V_{i-1}$ and~$x,y,z\in V_{i-1}$ with $(x,y)\in E_{i-1}$ such that~$V_i=V_{i-1}\cup\{v\}$ and~$E_i=E_{i-1}\diffsym\{(v,x),(v,y),(v,z),(x,y)\}$.
\end{enumerate}
\item \defn{Recski's Theorem}~\cite{r-ntarss-84}: For any two distinct vertices~$v,w\in V$, the (multi)graph ${(V,E\cup\{(v,w)\})}$ is the union of two spanning trees.\qed
\end{enumerate}
\end{theorem}

In dimension~$3$ and higher, combinatorial rigidity is far less understood. Generic rigid graphs still require sufficiently many non-redundant edges: 

\begin{proposition}[Laman's Condition]
A minimally generically rigid graph~$G \eqdef (V,E)$ in dimension~$d$ is a \defn{$d$-Laman graph}: $|E|=d|V|-{d+1 \choose 2}$ and any subgraph of~$G$ on~$v$ vertices has at most~$dv-{d+1 \choose 2}$ edges.\qed
\end{proposition}

\begin{figure}[b]
	\capstart
	\centerline{\includegraphics[scale=1.05]{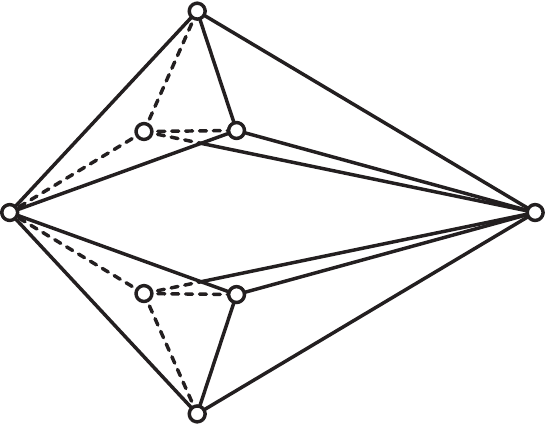}}
	\caption[The double banana]{The double banana: a flexible but Laman \dimensional{3} framework.}
	\label{ft:fig:doublebanana}
\end{figure}

However, this condition is not sufficient to be generically minimally rigid. The classical example is that of the \defn{double banana} (see \fref{ft:fig:doublebanana}) which is obtained by gluing two bananas (which are just two glued tetrahedra) at their endpoints. It has~$8$ vertices and $18$~edges, and any induced subgraph on~$v$ vertices has less than~$3v-6$ edges. However, it is flexible since the bananas can rotate independently around the axis passing through their endpoints.

Laman graphs are particular examples of the following class of graphs:

\begin{definition}
\index{sparse graph}
\index{tight graph|hbf}
A (multi)graph~$G \eqdef (V,E)$ is \defn{\spars{(p,q)}} if any subgraph of~$G$ on~$v$ vertices has at most~$pv-q$ vertices. If furthermore~$|E|=p|V|-q$, then~$G$ is said to be \defn{\tight{(p,q)}}.
\end{definition}

Sparsity for several families of parameters~$(p,q)$ is relevant in rigidity theory, matroid theory and pebble games~---~see~\cite{ls-pgasg-08} and the references therein. The equivalence between properties~(2), (3)~and~(4) of Theorem~\ref{ft:theo:rigidity2d} is a specialization of the following theorem to the case~$(p,q)=(2,3)$.

\begin{theorem}[\cite{h-cag-02}]
\index{arborescence@\arbores{(p,q)}}
Let $p,q$ be to integers such that $0\le p\le q<2p$. The following conditions on a graph~$G \eqdef (V,E)$ are equivalent:
\begin{enumerate}
\item $G$~is \defn{\tight{(p,q)}}: $|E|=p|V|-q$ and any subgraph of~$G$ on~$v$ vertices has at most~$pv-q$ edges.
\item $G$~is a \defn{\arbores{(p,q-p)}}: Adding any~$q-p$ edges to~$G$ (including multiple edges) results in a (multi)graph decomposable into~$p$ edge-disjoint spanning trees.
\end{enumerate}
If furthermore $q\le 3p/2$, then these two conditions are equivalent to the following:
\begin{enumerate}
\setcounter{enumi}{2}
\item $G$~can be constructed by a \defn{Henneberg construction}, starting from the graph with~$2$ vertices and~$2p-q$ (multiple) edges, and repeating the following operation: remove a subset~$F$ of~$E$ with $0\le |F|\le p$, and add a new vertex~$v$ to~$V$ adjacent to the~$2|F|$ endpoints of~$F$ (with possibly multiple edges) and to~$p-|F|$ additional arbitrary vertices of~$V$, in such a way that no edge has multiplicity greater than~$2p-q$.\qed
\end{enumerate}
\end{theorem}

Observe that \tight{(p,p)} graphs are exactly unions of $p$ edge-disjoint spanning trees.


\subsection{Sparsity in multitriangulations}

The goal of this section is to present where rigidity, sparsity and arboricity show up in the context of multitriangulations. The motivation comes from the application of rigidity properties to the construction of the polytope of pseudotriangulations (see Remarks~\ref{ft:rem:rigiditypseudotriangulationspolytope} and~\ref{ft:rem:rigiditypolytope}). The starting observation is the following consequence of Corollary~\ref{stars:coro:starsenumeration}:

\begin{corollary}
\ktri{k}s are \tight{\left(2k,{2k+1 \choose 2}\right)} graphs.\qed
\end{corollary}

As mentioned previously, being \tight{\left(2k,{2k+1 \choose 2}\right)} is a necessary, but not sufficient, condition for a graph to be minimally generically rigid in dimension~$2k$. This suggests the conjecture:

\begin{conjecture}
Every \ktri{k} is minimally generically rigid in dimension~$2k$.
\end{conjecture}

This conjecture is true when~$k=1$. More interestingly, we can prove it for~$k=2$:

\begin{theorem}\label{ft:theo:rigid2d}
Every \ktri{2} is a generically minimally rigid graph in dimension~$4$.
\end{theorem}

\begin{proof}
We prove by induction on~$n$ that \ktri{2}s of the \gon{n} are generically rigid in dimension~$4$. Minimality then follows from the fact that a generically rigid graph in dimension~$d$ needs to have at least~$dn-{d+1 \choose 2}$ edges.

Induction begins with the unique \ktri{2} on five points, that is,  the complete graph~$K_5$, which is generically rigid in dimension~$4$.

For the inductive step, let us recall the following graph-theoretic construction called ``vertex split'' (see \fref{ft:fig:split} and~\cite{w-vsif-90} for more details). Let~$G \eqdef (V,E)$ be a graph,~$u\in V$ be a vertex of~$G$ and~$U \eqdef \ens{u' \in V}{(u,u')\in E}$ denote the vertices adjacent to~$u$. Let~$U = X\sqcup Y\sqcup Z$ be a partition of~$U$. The \defn{vertex split} of~$u$ on~$Y$ is the graph whose vertex set is~$V\cup\{v\}$ (where~$v$ is a new vertex) and whose edge set is obtained from $E$ by removing $\ens{[u,z]}{z\in Z}$ and adding $\{[u,v]\}\cup\ens{[v,y]}{y\in Y}\cup\ens{[v,z]}{z\in Z}$.

\begin{figure}[!h]
	\capstart
	\centerline{\includegraphics[scale=1]{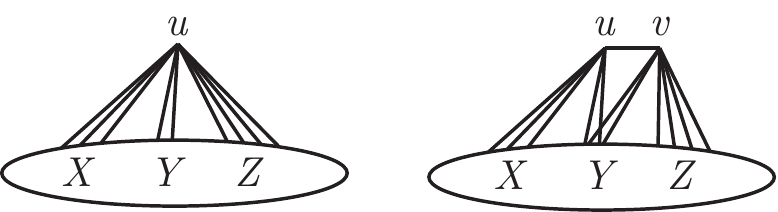}}
	\caption[A vertex split on two edges]{A vertex split on two edges.}
	\label{ft:fig:split}
\end{figure}

The key result that we need is that a vertex split on~$d-1$ edges in a generically rigid graph in dimension~$d$ is also generically rigid in dimension~$d$~\cite{w-vsif-90}.

So, let~$n\ge 5$ and assume that we have already proved that every \ktri{2} of the \gon{n} is rigid. Let~$T$ be a \ktri{2} of the \gon{(n+1)}. Let~$S$ be a \kstar{2} of~$T$ with at least two \kbound{2} edges (such a \kstar{2} exists since it appears in the ``outer'' side of any \kear{2}). It is easy to check that the inverse transformation of the flattening of~$S$ is exactly a vertex split on~$3$ edges. Thus, the result follows.
\end{proof}

\begin{figure}
	\capstart
	\centerline{\includegraphics[scale=1]{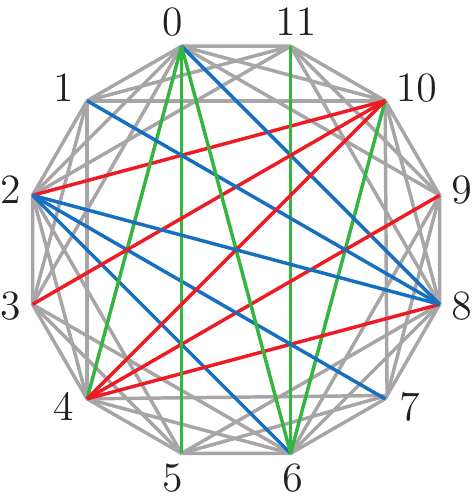}}
	\caption[A \ktri{3} without~$2$ consecutives \kear{3}s]{A \ktri{3} without~$2$ consecutives \kear{3}s.}
	\label{ft:fig:3triang3ear}
\end{figure}

\enlargethispage{.3cm}
Generalizing this proof, observe that if~$T$ is a \ktri{k} and~$S$ a \kstar{k} of~$T$ with~$k$ \kbound{k} edges (or equivalently~$k-1$ consecutive \kear{k}s), then the inverse transformation of the flattening of~$S$ is exactly a vertex split on~$2k-1$ edges. However, our proof of Theorem~\ref{ft:theo:rigid2d} cannot be directly applied since there exist \ktri{k}s with no \kstar{k} containing~$k$ \kbound{k} edges, or equivalently, without~$k-1$ consecutive \kear{k}s, for~$k\ge 3$ (see \fref{ft:fig:3triang3ear}).

\mvs
Sparsity also shows up on dual graphs of multitriangulations. Remember from Section~\ref{stars:subsec:surfaces:dual} that we defined the dual graph~$T^\dual$ of a \ktri{k}~$T$ to be the multigraph whose vertices are the \kstar{k}s of~$T$ and with one edge between two \kstar{k}s for each of their common edges. As well as dual graphs of triangulations are trees, duals of multitriangulations are arborescences:

\begin{proposition}
The dual graph of a \ktri{k} is the union of $k$ edge-disjoint spanning trees. In other words, dual graphs of \ktri{k}s are \tight{(k,k)}.
\end{proposition}

\begin{proof}
The proof works again by induction on~$n$. When~$n=2k+1$, the dual graph is a single vertex and there is nothing to prove. For the inductive step, consider a \ktri{k}~$T$ of the \gon{n}, and a \kcross{k}~$X$ of~$T$ which is external, \ie adjacent to~$k$ consecutive vertices. We want to derive from a decomposition of the dual graph of~$T$ into~$k$ edge-disjoint spanning trees a similar decomposition of the dual graph of the inflating~$\inflating{T}{X}$ of~$X$ in~$T$.

Observe first that the graph~$(\inflating{T}{X})^\dual$ can be easily described from~$T^\dual$: Let $Y$ denote the multiset of vertices of~$T^\dual$ corresponding to all \kstar{k}s of~$T$ adjacent to an edge of~$X$ (a \kstar{k} appears in~$Y$ as many times as its number of edges in~$X$) and let~$Z$ denote the set of edges of~$T^\dual$ corresponding to \krel{k} edges of~$X$. Then~$(\inflating{T}{X})^\dual$ is obtained from~$T^\dual$ by removing the edges of~$Z$ and adding a new vertex~$v$ adjacent to all vertices of~$Y$ (with multiple edges if necessary). Observe that the vertex~$v$ has degree~$|Y|=k+|Z|$.

Now  consider a decomposition~$\tau_1,\dots,\tau_k$ of~$T^\dual$ into~$k$ edge-disjoint spanning trees. For each~$1\le i\le k$, let ~$\alpha_i \eqdef |\tau_i\cap Z|$ be the number of edges of the \kcross{k}~$X$ whose dual edge is in~$\tau_i$. Then~$\tau_i\ssm Z$ is a forest with~$1+\alpha_i$ trees. For all~$i$ such that $\alpha_i\ne0$, we transform this forest into a tree~$\tilde\tau_i$ by choosing~$1+\alpha_i$ edges between~$v$ and the~$2\alpha_i$ endpoints of~$\tau_i\cap Z$. Then, the vertex~$v$ still has~$k+|Z|-\sum_{i\,|\, \alpha_i\ne 0} (1+\alpha_i) = |\ens{i}{\alpha_i=0}|$ uncolored edges. We arbitrarily affect one such edge to connect~$v$ to each tree $\tau_i$ with $\alpha_i=0$. The resulting trees~$\tilde\tau_1,\dots,\tilde\tau_k$ form a decomposition of~$(\inflating{T}{X})^\dual$. This finishes the proof since any \ktri{k} of the \gon{(n+1)} can be obtained from a \ktri{k} of the \gon{n} by inflating an external \kcross{k}.
\end{proof}

\begin{example}
\kcolorable{k} \ktri{k}s are exactly the \ktri{k}s whose dual graph is the union of $k$ edge-disjoint paths. Indeed, we already implicitly proved in Section~\ref{stars:sec:ears} that the dual graph of a \kcolorable{k} \ktri{k} is a union of~$k$ disjoint paths (remember the algorithm which computes the~$k$ disjoint \kaccordion{k}s of a \kcolorable{k} \ktri{k} in the proof of Theorem~\ref{stars:theo:kcolorable}). Reciprocally, if the dual graph of a \ktri{k} can be decomposed into~$k$ paths, then no \kstar{k} of~$T$ can contain more than~$2k$ \krel{k} edges. Thus all \kstar{k}s of~$T$ are external, which characterizes \kcolorable{k} \ktri{k}s (Theorem~\ref{stars:theo:kcolorable}).
\end{example}

\begin{figure}
	\capstart
	\centerline{\includegraphics[scale=1]{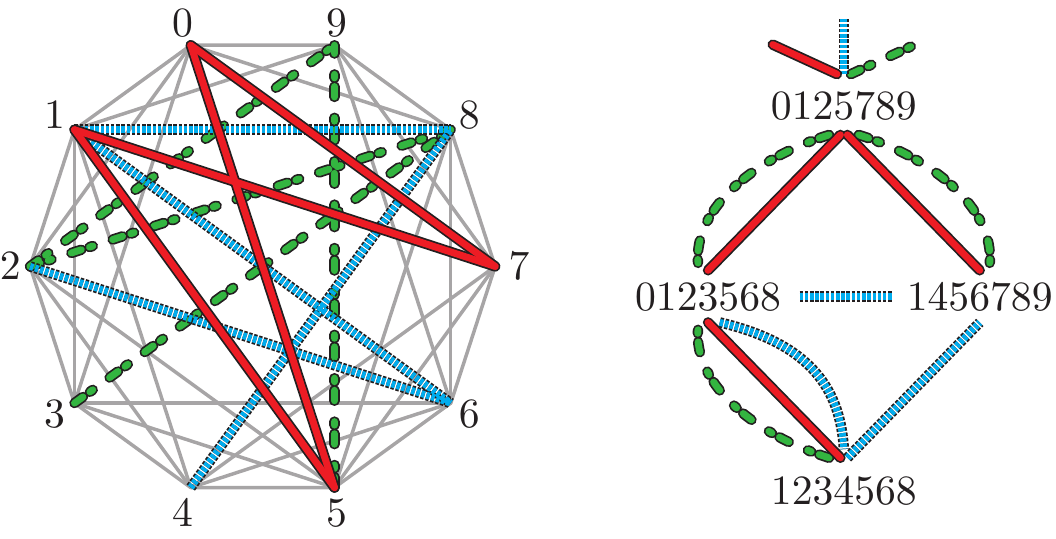}}
	\caption[A counter-example]{A \ktri{3} of the \gon{10} in which coloring in each \kstar{3} the first top edge in red, the second in blue and the third in green does not decompose the dual graph into three edge-disjoint spanning trees.}
	\label{ft:fig:bad3coloring}
\end{figure}

\begin{example}
In Section~\ref{ft:subsec:Dyckpaths:failedbijection}, we have defined a $k$-coloring of the \krel{k} edges of a \ktri{k}~$T$:  for any \kstar{k} of~$T$ with vertices~$0\cle s_1\cl s_2\cl\cdots\cl s_{2k+1}\cle n-1$, color the first top edge~$[s_1,s_{k+2}]$ with color~$1$, the second top edge~$[s_2,s_{k+3}]$ with color~$2$, and so on until the $k$th top edge~$[s_k,s_{2k+1}]$ with color~$k$.

When~$k=2$, this $2$-coloring of the \krel{2} edges of a \ktri{2} defines a decomposition of the dual graph~$T^\dual$ into two edge-disjoint spanning trees. Indeed, it is easy to check that if a \kstar{2}~$S$ (with vertices~$0\cle s_1\cl\cdots\cl s_5\cle n-1$) shares its first top edge~$[s_1,s_4]$ with another \kstar{2}~$R$ (with vertices~$0\cle r_1\cl\cdots\cl r_5\cle n-1$), then the first top edge~$[r_1,r_4]$ of~$R$ is smaller than that of~$S$ for the order defined by $[a,b]<[c,d] \Leftrightarrow (a<c) \text{ or } (a=c \text{ and } b>d)$. Thus, the subgraph of~$T^\dual$ corresponding to the edges colored with~$1$ has no cycle (its vertices are partially ordered). Since it has one edge less than the number of vertices, it is a spanning tree. The proof is symmetric for the graph produced by the edges colored with~$2$.

However, when~$k>2$, the resulting coloring is not necessarily a decomposition of~$T^\dual$ into edge-disjoint spanning trees. See for example \fref{ft:fig:bad3coloring}, where one color has a cycle.
\end{example}

\section{Multiassociahedron}\label{ft:sec:multiassociahedron}

In this section, we consider the following complex:

\begin{definition}\label{ft:def:deltank}
We denote by~$\Delta_{n,k}$ the simplicial complex of all \kcross{(k+1)}-free subsets of \krel{k} edges of~$E_n$.
\end{definition}

It is indeed a simplicial complex since \kcross{(k+1)}-freeness is a hereditary condition: if a set of edges has no \kcross{(k+1)}, then any of its subsets is \kcross{(k+1)}-free.

Observe that we do not include the non-\krel{k} edges of the \gon{n} in the sets of~$\Delta_{n,k}$. Indeed, since it cannot appear in a \kcross{(k+1)}, a non-\krel{k} edge would be a cone point of the complex: the simplicial complex of all \kcross{(k+1)}-free subsets of~$E_n$ is a~$kn$-fold cone over~$\Delta_{n,k}$. However, in order to keep notations and explanations simple, we sometimes implicitly include all the non-\krel{k} edges of the \gon{n} in the elements of~$\Delta_{n,k}$, as for example in the next paragraph.

Our study of \ktri{k}s of the \gon{n} already gave some combinatorial information about the simplicial complex~$\Delta_{n,k}$:
\begin{enumerate}[(i)]
\item The maximal elements of~$\Delta_{n,k}$ are the \ktri{k}s of the \gon{n}. Since all \ktri{k}s of the \gon{n} have~$k(n-2k-1)$ \krel{k} egdes (Corollary~\ref{stars:coro:starsenumeration}), the simplicial complex~$\Delta_{n,k}$ is \defn{pure} of dimension~$k(n-2k-1)-1$.
\item The codimension~$1$ elements of~$\Delta_{n,k}$ are flips between \ktri{k}s: they are obtained by removing one \krel{k} edge from a \ktri{k}, and we know that such a set of edges is contained in exactly two \ktri{k}s. Thus, every codimension~$1$ face of~$\Delta_{n,k}$ is contained in exactly two facets of~$\Delta_{n,k}$ (in other words, $\Delta_{n,k}$~is a pseudo-manifold). The graph of flips that we studied in Section~\ref{stars:sec:flips} is the ridge-graph of~$\Delta_{n,k}$. 
\item Finally, a subset of at most~$k$ edges never contains a \kcross{(k+1)}. Thus, the simplicial complex~$\Delta_{n,k}$ is~\defn{\neighborly{k}}: it contains all sets of at most~$k$ elements.
\end{enumerate}

\mvs
Before going further, we revisit the examples of Section~\ref{stars:subsec:notations:examples}:

\begin{example}[$k=1$]\label{ft:exm:associahedron}
The simplicial complex~$\Delta_{n,1}$ of all crossing-free sets of internal diagonals of the \gon{n} is isomorphic to the boundary complex of the polar of the associahedron~\cite{l-atg-89,bfs-ccsp-90,gkz-drmd-94,l-rsp-04,hl-rac-07}. We recall two constructions of this polytope in Section~\ref{ft:subsec:multiassociahedron:associahedron}.
\end{example}

\begin{example}[$n=2k+\varepsilon$ with ${\varepsilon\in[3]}$]\label{ft:exm:n=2k+eps}
Since the \gon{(2k+1)} has no \krel{k} edge, and thus a unique \ktri{k}, the simplicial complex~$\Delta_{2k+1,k}=\{\emptyset\}$ is trivial.

More interestingly, the \krel{k} edges of the \gon{(2k+2)} are its $k+1$ long diagonals, which form the unique \kcross{(k+1)}. Any set of at most~$k$ of them is \kcross{(k+1)}-free, and the simplicial complex~$\Delta_{2k+2,k}$ is isomorphic to the boundary complex of the \simp{k}~$\simplex_k$ (see \fref{ft:fig:2k+2associahedron}).

\begin{figure}[!h]
	\capstart
	\centerline{\includegraphics[scale=.8]{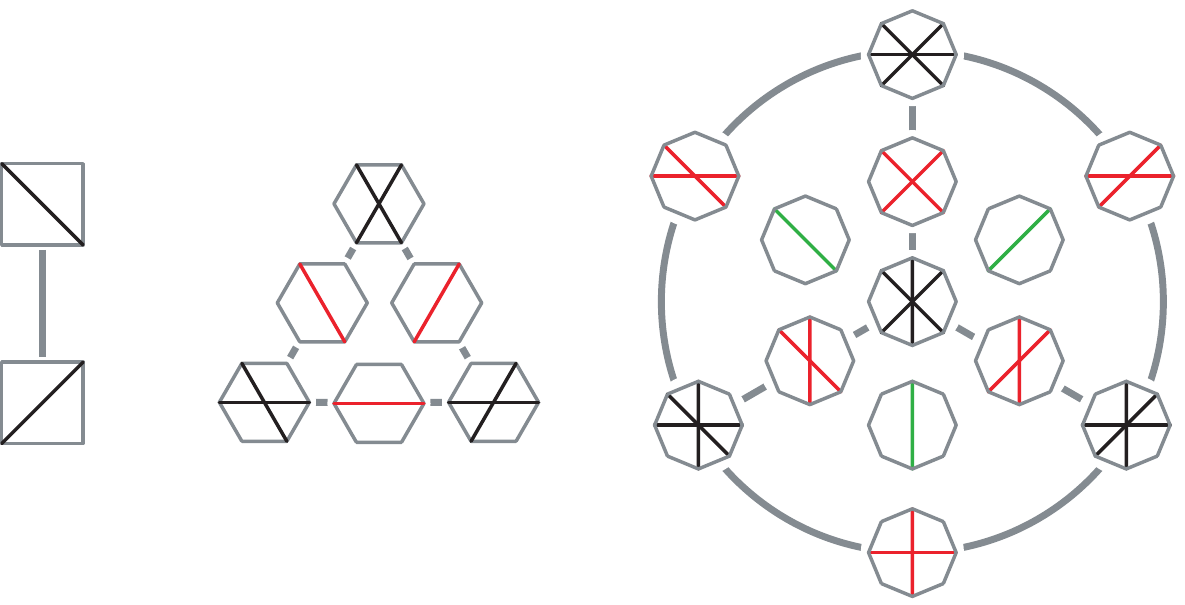}}
	\caption[The simplicial complex~$\Delta_{2k+2,k}$]{The simplicial complex~$\Delta_{2k+2,k}$ for~$k\in[3]$ (in the rightmost picture, the external face corresponds to the horizontal edge).}
	\label{ft:fig:2k+2associahedron}
\end{figure}

Finally, the last small case completely described in Section~\ref{stars:subsec:notations:examples} is when $n=2k+3$. The corresponding simplicial complex has also a very nice description:

\begin{lemma}\label{ft:lem:2k+3cyclicpolytope}
The simplicial complex~$\Delta_{2k+3,k}$ is the boundary complex of the cyclic polytope of dimension~$2k$ with~$2k+3$ vertices.
\end{lemma}

\begin{proof}
Remember from Example~\ref{stars:exm:n=2k+3} that the \ktri{k}s of the \gon{(2k+3)} are exactly the disjoint unions of~$k$ pairs of consecutive long diagonals of~$E_{2k+3}$. This translates, by Gale's evenness criterion (see Proposition~\ref{intro:prop:gale}), into the facets of the cyclic polytope.
\end{proof}

\end{example}

According to these examples, it is reasonable to wonder whether there exists a polytope generalizing the associahedron for all values of~$k$:

\begin{question}\label{ft:qu:polytope}
\index{boundary complex}
\index{realization!polytopal ---}
Is~$\Delta_{n,k}$ the boundary complex of a polytope for any~$n$ and~$k$?
\end{question}

If so, we will say that~$\Delta_{n,k}$ is \defn{realizable}. Observe that any \defn{realization} would certainly be a simplicial \poly{k(n-2k-1)}tope.

An important step towards the answer of this question is to know whether the simplicial complex~$\Delta_{n,k}$ has the topology of a sphere. This was proved by Jakob Jonsson~\cite[Theorem~1]{j-gt-03} in the following sense:

\begin{theorem}[\cite{j-gt-03}]
The simplicial complex~$\Delta_{n,k}$ is a vertex-decomposable piece-wise linear sphere of dimension~$k(n-2k-1)-1$.\qed
\end{theorem}

Being a topological sphere is a necessary~---~but not sufficient~---~condition to be realizable.

Despite our efforts to solve it, Question~\ref{ft:qu:polytope} remains open in general. However, we present here two interesting steps in its direction:
\begin{enumerate}
\item In Section~\ref{ft:subsec:multiassociahedron:delta82}, we prove that the simplicial complex~$\Delta_{8,2}$ is realizable, which solves the first non-trivial case (remember Examples~\ref{ft:exm:associahedron} and~\ref{ft:exm:n=2k+eps}). In fact, we do even more than just providing a single realization of~$\Delta_{8,2}$, since we completely describe its space of symmetric realizations.
\item In Section~\ref{ft:subsec:multiassociahedron:loday}, we propose a natural generalization of a construction of the associahedron due to Jean-Louis Loday~\cite{l-rsp-04}. This results in the construction of a polytope with a rich combinatorial structure: its graph is the graph of flips restricted to certain \ktri{k}s, and it has a very simple facet description. \textit{A priori} this polytope could have been a projection of a polytope realizing~$\Delta_{n,k}$ (we will see why it cannot be such a projection).
\end{enumerate}


\subsection{Two constructions of the associahedron}\label{ft:subsec:multiassociahedron:associahedron}\index{associahedron|hbf}

In this section, we recall two constructions of the associahedron of which we will discuss the generalization later. In both constructions, a vector is associated to each triangulation of the \gon{n}, and the associahedron is obtained as the convex hull of the vectors associated to all triangulations of the \gon{n}. However, to associate a vector to a triangulation, the point of view is very different from one construction to the other: the first construction is based on the vertices of the \gon{n} (the vector has one coordinate for each vertex), while the second construction focusses on the triangles of a triangulation (the vector has one coordinate per triangle).

\subsubsection{The secondary polytope of the \gon{n}}\label{ft:subsubsec:secondarypolytope}

This first construction is due to Israel Gel'fand, Mikhail Kapranov and Andrei Zelevinski. In this construction, the vector associated to each triangulation is defined as follows:

\begin{definition}[\cite{gkz-drmd-94}]\label{ft:def:secondarypolytope}
\index{secondary polytope}
\index{polytope!secondary ---}
The \defn{area vector} of a triangulation~$T$ of the \gon{n} is the vector~$\phi(T)$ of~$\R^{n}$, whose $v$th coordinate (for~$v\in V_n$) is the sum of the areas of the triangles of~$T$ containing the vertex~$v$. The \defn{secondary polytope} $\Sigma_n$ is the convex hull of the area vectors of all triangulations of the \gon{n}.
\end{definition}

\begin{remark}\label{ft:rem:secondarypolytopedimension}
We obtain affine relations among the coordinates of the area vectors of the triangulations of the \gon{n} by using the triangles to decompose the total area and the center of mass of the \gon{n}. Indeed, the areas of all the triangles of a triangulation~$T$ sum to the area of the \gon{n}:
$$\cA(\conv(V_n))=\sum_{\simplex\in T} \cA(\simplex) = \dotprod{\one}{\phi(T)}.$$
Furthermore, the center of mass~$\cm(\conv(V_n))$ of the \gon{n} is the barycenter, with coefficients~$\cA(\simplex)$, of the centers of mass~$\cm(\simplex)$ of the triangles of~$T$. Since the center of mass of a triangle is also its isobarycenter, we obtain:
$$\cm(\conv(V_n))=\sum_{\simplex\in T} \cA(\simplex)\cm(\simplex)=\sum_{\simplex\in T}\left(\cA(\simplex)\sum_{v\in\simplex} \frac{v}{3}\right)=\frac{1}{3}\sum_{v\in V_n} \left(v\sum_{v\in\simplex\in T} \cA(\simplex)\right).$$
Consequently, the secondary polytope~$\Sigma_n$ lies in an affine subspace of codimension~$3$.
\end{remark}

\begin{theorem}[\cite{gkz-drmd-94}]\label{ft:theo:secondarypolytope}
The boundary complex of the polar~$\Sigma_n^\polar$ of the secondary polytope is isomorphic to the simplicial complex~$\Delta_{n,1}$ of crossing-free sets of diagonals of the \gon{n}.
\end{theorem}

\begin{proof}
Consider a vector~$\omega\in\R^n$. For any triangulation~$T$ of the \gon{n},
$$\dotprod{\omega}{\phi(T)} = \sum_{v\in V_n} \left(\omega_v \sum_{v\in\simplex\in T} \cA(\simplex)\right) = \sum_{\simplex\in T} \left(\cA(\simplex)\sum_{v\in\simplex} \omega_v\right) = 3 \sum_{\simplex\in T} \cV(\prism(\simplex,\omega)),$$
where~$\cV(\prism(\simplex,\omega))$ denotes the volume of the prism with vertices~$(v,0)_{v\in\simplex}$ and~$(v,\omega_v)_{v\in\simplex}$. The last sum equals the volume below the piecewise linear function~${f(\omega,T):\conv(V_n)\to\R}$ which is linear on each triangle of~$T$ and which evaluates to~$\omega$ on~$V_n$ (\ie whose value at any vertex~$v\in V_n$ is $\omega_v$). Consequently, the scalar product~$\dotprod{\omega}{\phi(T)}$ is minimized precisely for the triangulations~$T$ for which~$f(\omega,T)$ coincides with the lower envelope of the lifted point set $\ens{(v,\omega_v)}{v\in V_n}\subset \R^3$, that is, for the triangulations which refine the polygonal subdivision of the \gon{n} obtained as the projection of the lower envelope of this lifted point set.

\index{subdivision!regular ---}
Any crossing-free set~$S$ of diagonals of the \gon{n} is \defn{regular}: there exists a lifting vector~$\omega$ such that the projection of the lower envelope of the lifted point set coincides with the polygonal subdivision induced by~$S$. According to the discussion above, $\omega$ is a normal vector of a face of~$\Sigma_n$ which contains exactly the area vectors of the triangulations of the \gon{n} containing~$S$.
\end{proof}

For example, the secondary polytopes~$\Sigma_3$ (of a triangle),~$\Sigma_4$ (of a square), and~$\Sigma_5$ (of a regular pentagon) are respectively a single point, a segment, and a regular pentagon. More interesting is the secondary polytope of an hexagon:

\begin{example}[The \dimensional{3} associahedron]
\fref{ft:fig:associahedron} displays the secondary polytope~$\Sigma_6$ of a regular hexagon and its polar polytope~$\Sigma_6^\polar$. As mentioned previously, a vertex (resp.~an edge, resp.~a facet) of~$\Sigma_6$ corresponds to a triangulation (resp.~a flip, resp.~a diagonal) of the hexagon (and the reverse holds for the polar polytope).

\begin{figure}
	\capstart
	\centerline{\includegraphics[scale=1]{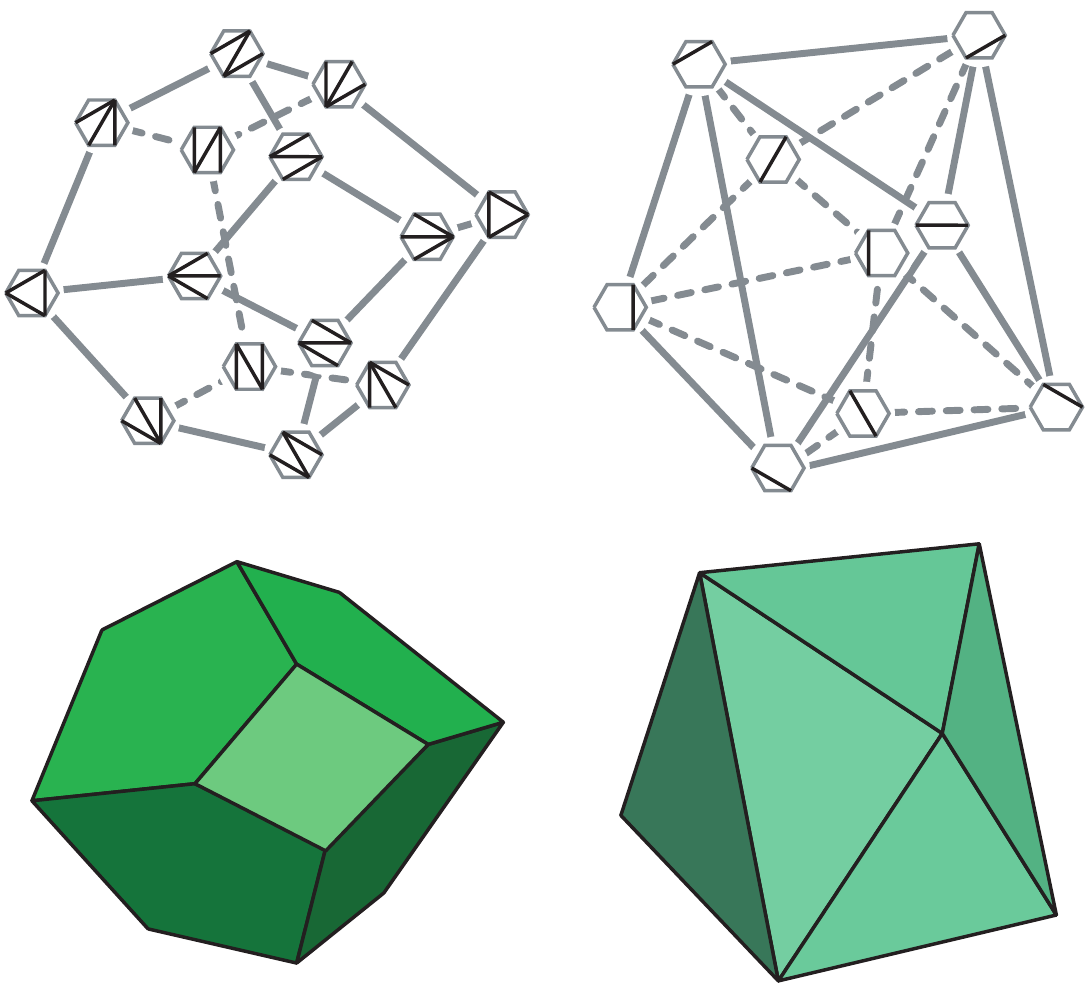}}
	\caption[The \dimensional{3} associahedron and its polar]{The \dimensional{3} associahedron (left) and its polar (right).}
	\label{ft:fig:associahedron}
\end{figure}

It is easy to observe on the picture that the facets of the associahedron are products of two smaller associahedra. Here the facets of~$\Sigma_6$ are either pentagons (\ie $\Sigma_3\times\Sigma_5$) or squares (\ie $\Sigma_4\times\Sigma_4$). This is a general fact since an internal diagonal separates the \gon{n} into two smaller polygons, one of size~$\ell+1$ and the other one of size~$n-\ell+1$ (where~$\ell$ is the length of the diagonal). More generally, all codimension~$p$ faces of a \dimensional{d} associahedron are products of~$p+1$ smaller associahedra whose dimensions sum to~$d-p$.
\end{example}

It is interesting to observe that the full isometry group of the regular \gon{n} defines an isometry action on the secondary polytope:

\begin{lemma}\label{ft:lem:secondarypolytopesymmetry}
The dihedral group~$\D_n$ of isometries of the regular \gon{n} acts on its secondary polytope~$\Sigma_n$ by isometry.
\end{lemma}

\begin{proof}
An isometry of the regular \gon{n} translates on the secondary polytope into a permutation of the coordinates.
\end{proof}

\begin{example}[The \dimensional{3} associahedron continued]

\begin{figure}
	\capstart
	\centerline{\includegraphics[scale=1]{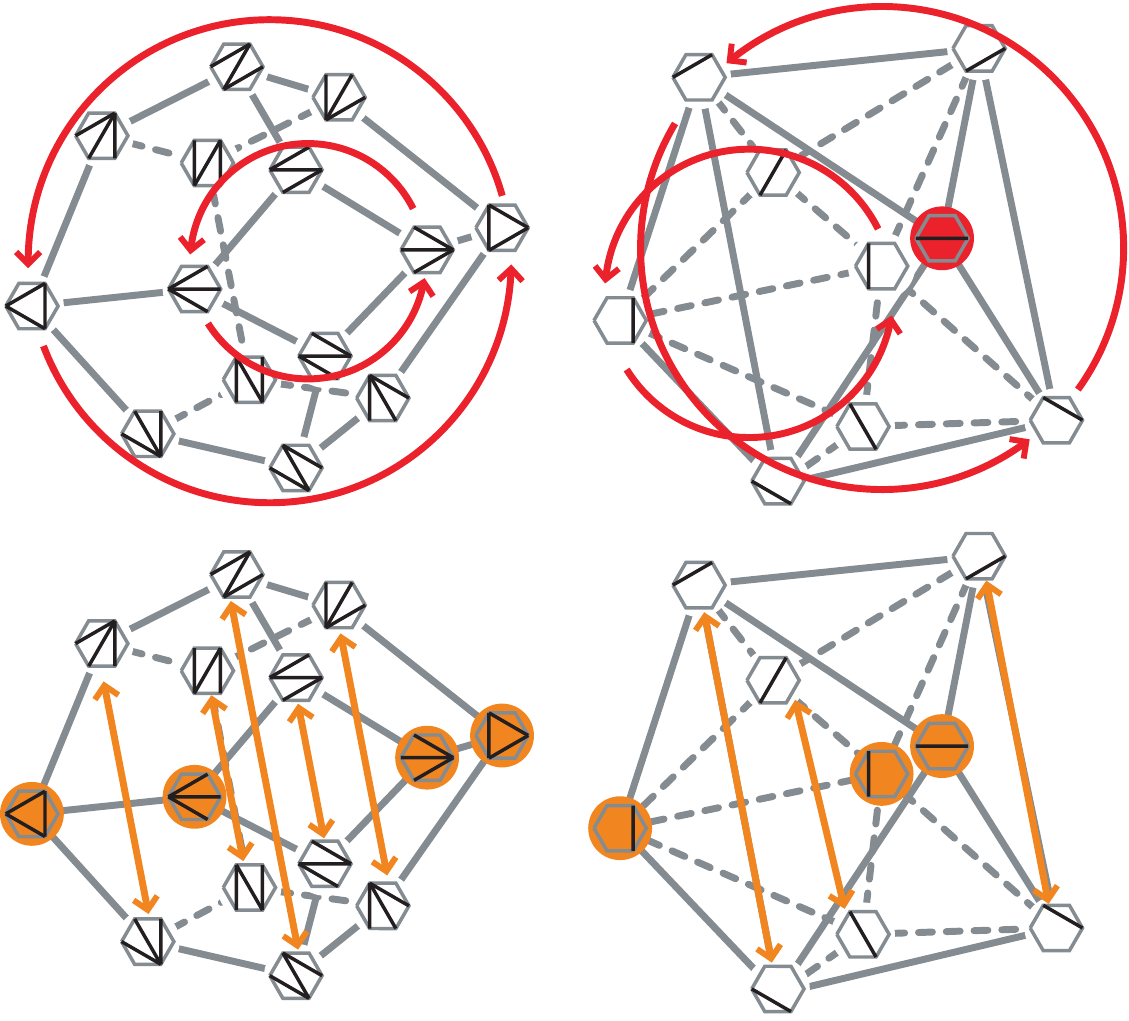}}
	\caption[Two symmetries on the \dimensional{3} associahedron]{Two symmetries on the \dimensional{3} associahedron.}
	\label{ft:fig:associahedronSymmetry}
\end{figure}

\fref{ft:fig:associahedronSymmetry} displays the action of two isometries of the regular hexagon on the \dimensional{3} associahedron~$\Sigma_6$ and its dual~$\Sigma_6^\polar$. The first one is the reflection~$\rho_v$ of the hexagon with respect to the vertical axis, which translates on the secondary polytope~$\Sigma_6$ into a rotation (whose axis contains the point corresponding to the horizontal edge of the hexagon). The second one is the reflection~$\rho_h$ of the hexagon with respect to the horizontal axis, which translates on the secondary polytope~$\Sigma_6$ into a reflection (with respect to the plane containing the fixed triangulations). Even if our two isometries are both indirect on the hexagon (they are both reflections), the corresponding isometries on the secondary polytope can be direct or indirect.

Observe that our two reflections generate the complete dihedral group~$\D_6$ of the regular hexagon. Thus, the action of any isometry of the hexagon on the secondary polytope can be deduced from \fref{ft:fig:associahedronSymmetry}. It is however not always easy to see: imagine for example the action on~$\Sigma_6$ of the rotation of angle~$\pi/3$ of the hexagon.
\end{example}

\begin{remark}
This construction of the secondary polytope is not restricted to the very special situation of convex point configurations in the plane. If~$P$ is an arbitrary point set in~$\R^d$, the secondary polytope of~$P$ is defined as the convex hull of the volume vectors of all triangulations of~$P$~---~where the volume vector of a triangulation~$T$ of~$P$ is the vector of~$\R^{|P|}$ whose $p$th coordinate (for $p\in P$) is the sum of the volumes of the simplices of~$T$ adjacent to the vertex~$p$. Then the face lattice of the secondary polytope is isomorphic to the refinement lattice on regular subdivisions of~$P$. It has dimension~$|P|-d-1$, and its vertices (resp.~edges) correspond to regular triangulations of~$P$ (resp.~flips between them). We refer to~\cite{bfs-ccsp-90,gkz-drmd-94,lrs-tri} for further details, and restrict to convex configurations in the plane.
\end{remark}

\subsubsection{Loday's construction of the associahedron}\label{ft:subsubsec:lodayconstruction}

The second construction of the associahedron that we try to generalize is due to Jean-Louis Loday~\cite{l-rsp-04}. The resulting polytope is a bit less symmetric, but has other advantages: in particular, it has integer coordinates, and its facets have simple equations. Here, we follow the original presentation~\cite{l-rsp-04} which first defines a polytope associated to binary trees. We will see later in Section~\ref{ft:subsec:multiassociahedron:loday} an equivalent definition in terms of beam arrangements and its generalization to multitriangulations.

We call \defn{\dfslabeling{}} of a rooted binary tree the labeling obtained by a depth-first search of the tree: we start from the root, and for each internal node, we first completely search its left child, then label our node with the smallest free label, and finally visit its right child. In other words, the \dfslabeling{} is the unique labeling of the tree with the property that the label of any node is bigger than any label in its left child and smaller than any label in its right child. 

The construction associates a vector to each binary tree as follows:

\begin{definition}[\cite{l-rsp-04}]\label{ft:def:loday}
The \defn{\dfsvector{}} of a rooted binary tree~$\tau$ on $m$ nodes is the vector $\psi(\tau)\in\R^m$ whose~$i$th coordinate (for~$i\in[m]$) is the product of the number of leaves in the left child by the number of leaves in the right child of the node labeled by~$i$ in the \dfslabeling{}~of~$\tau$. The \defn{\dfspolytope{}}~$\Omega_m$ is the convex hull of the \dfsvector{}s of all rooted binary trees on~$m$~nodes.
\end{definition}

Remember that triangulations of the \gon{n} and binary trees on~$n-2$ nodes are in bijective correspondence (see Section~\ref{intro:sec:triangulations}). Thus, we define by extension the \defn{\dfsvector{}} of a triangulation of the \gon{n} to be the \dfsvector{} of its dual binary tree on~$n-2$ nodes.

\begin{theorem}[\cite{l-rsp-04}]\label{ft:theo:loday}
The boundary complex of the polar~$\Omega_{n-2}^\polar$ of the \dfspolytope{} is isomorphic to the simplicial complex~$\Delta_{n,1}$ of crossing-free sets of diagonals of the \gon{n}. Indeed,
\begin{enumerate}[(i)]
\item $\Omega_{n-2}$~is included in the affine hyperplane $\ens{x\in\R^{n-2}}{\dotprod{\one}{x}={n-1 \choose 2}}$;~and
\item for any $0\le u <v\le n-1$, the equation $\sum_{p=u+1}^{v-1} x_p = {v-u \choose 2}$ defines a facet of~$\Omega_{n-2}$ which contains all the \dfsvector{}s of the triangulations of the \gon{n} containing the edge~$[u,v]$, and only these ones.\qed
\end{enumerate}
\end{theorem}

We do not repeat here the proof of~\cite{l-rsp-04} since our interpretation of this construction in Section~\ref{ft:subsec:multiassociahedron:loday} in terms of beam arrangements will provide a visual and simple proof.

\begin{example}
\fref{ft:fig:associahedronLoday} shows the \dimensional{3} associahedron as realized by Loday's construction.
\begin{figure}
	\capstart
	\centerline{\includegraphics[scale=1]{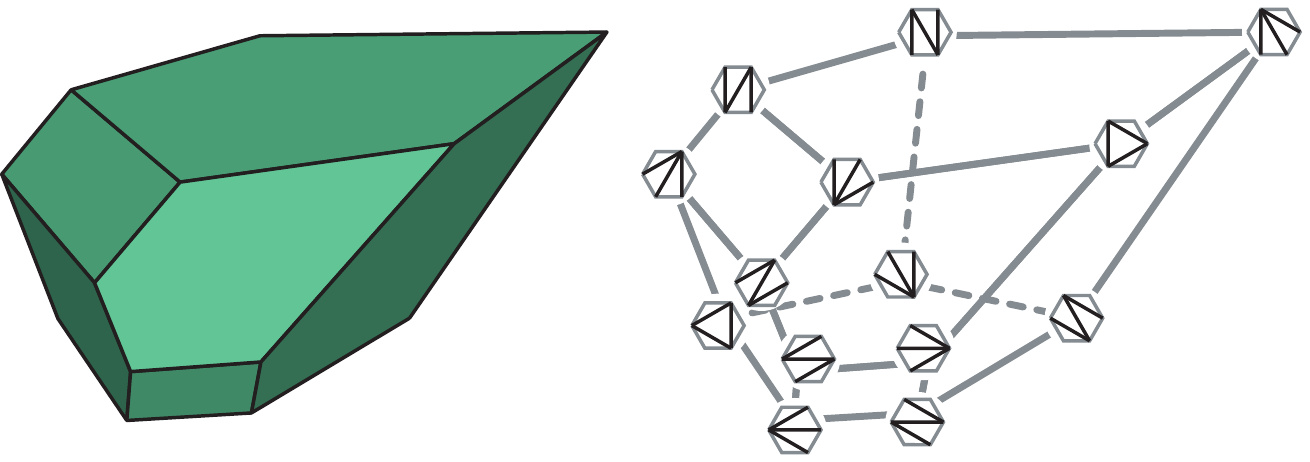}}
	\caption[Loday's construction of the \dimensional{3} associahedron]{Loday's construction of the \dimensional{3} associahedron.}
	\label{ft:fig:associahedronLoday}
\end{figure}
\end{example}

Observe that since this polytope is defined on binary trees, its symmetry group is automatically not the full dihedral group of isometries of the regular \gon{n}. However, the construction conserves the symmetry of the trees:

\begin{lemma}\label{ft:lem:lodaysymmetry}
The polytope~$\Omega_m$ is symmetric under reversing coordinates.
\end{lemma}

\begin{proof}
Reversing the coordinates of the \dfsvector{} corresponds to inverting left and right in the depth-first search.
\end{proof}

This symmetry corresponds to the reflection of the \gon{n} with respect to the diameter passing through the middle of the edge~$[0,n-1]$. In \fref{ft:fig:associahedronLoday}, it is a rotation around an axis passing through the center of the preserved square and the center of the preserved edge.


\subsection{The space of symmetric realizations of $\Delta_{8,2}$}\label{ft:subsec:multiassociahedron:delta82}

In this section, we focus on symmetric realizations of the simplicial complex~$\Delta_{n,k}$, in the following sense:

\begin{definition}\label{ft:def:symmetricrealization}
\index{realization!polytopal ---}
\index{realization!symmetric polytopal ---}
A polytope~$P\subset\R^d$ is a \defn{realization} of a simplicial complex~$\Delta$ if its boundary complex~$\partial P$ is isomorphic to~$\Delta$, \ie if there is a bijection~$\gamma:\Delta\to\partial P$ which respects inclusion: $S\subset T\Leftrightarrow \gamma(S)\subset\gamma(T)$,~for all~$S,T\in\Delta$. We say that~$P$ is a \defn{symmetric realization} under a group~$G$ acting on~$\Delta$ if for any~$g\in G$, the action~$\partial P\to \partial P: y\mapsto gy \eqdef \gamma(g\gamma^{-1}(y))$ of~$g$ on~$\partial P$ is an isometry.
\end{definition}

When~$G$ is the complete symmetry group of the simplicial complex~$\Delta$, we sometimes do not mention~$G$ explicitely, and we just say that~$P$ is a symmetric realization of~$\Delta$. As far as the simplicial complex~$\Delta_{n,k}$ is concerned, we already know examples of symmetric realizations:
\begin{itemize}
\item[~~~~\fbox{$k=1$}] The associahedron~$\Sigma_n$ constructed as the secondary polytope of the regular \gon{n} is a symmetric realization of~$\Delta_{n,1}$ under the dihedral group~$\D_n$.
\item[~~~~\fbox{$n=2k+1$}] A point is a symmetric realization of~$\Delta_{2k+1,k}$.
\item[~~~~\fbox{$n=2k+2$}] The regular \simp{k} is symmetric under any permutation of its vertices, and thus, is a symmetric realization of~$\Delta_{2k+2,k}$ under~$\sym_{k+1}$.
\item[~~~~\fbox{$n=2k+3$}] Define the \dimensional{2k} cyclic polytope with~$2k+3$~vertices embedded on the \defn{Caratheodory curve} $\nu_{2k}:\theta\mapsto(\cos(\ell\theta),\sin(\ell\theta))_{\ell\in[k]}$ as
$$C_{2k}(2k+3) \eqdef \conv\ens{\nu_{2k}\left(e^{\frac{2i\pi m}{2k+3}}\right)}{m\in\Z_{2k+3}}\subset\R^{2k}.$$
Then~$C_{2k}(2k+3)$ is a symmetric realization of~$\Delta_{2k+3,k}$ under the dihedral group~$\D_{2k+3}$ (see for example~\cite{kw-agcp-10}).
\end{itemize}

\mvs
The main result of this section concerns the first non-trivial case: we study symmetric realizations of the simplicial complex~$\Delta_{8,2}$ of \kcross{3}-free sets of \krel{2} edges of the octagon. We use computer enumeration to obtain a complete description of the space of symmetric realizations of~$\Delta_{8,2}$. This enumeration works in two steps:
\begin{enumerate}[(i)]
\item We first enumerate all possible ``combinatorial'' configurations (matroid polytopes) which could realize the simplicial complex~$\Delta_{8,2}$.
\item Then, from the knowledge of these combinatorial configurations, we deduce all ``geometric'' polytopes realizing~$\Delta_{8,2}$ and symmetric under the dihedral group.
\end{enumerate}
Before presenting this result, we briefly remind the notions we need in the world of oriented matroids (which will make precise what we mean by ``combinatorial configurations''), and we expose the method on the basic example of~$\Delta_{6,1}$.

\subsubsection{Chirotopes and symmetric matroid realizations}\label{ft:subsubsec:orientedmatroids}

For our purpose, \defn{oriented matroids}\index{oriented matroids} are combinatorial abstractions of point configurations: they keep track of the arrangement in a point set, by encoding relative positions of its points. For example, we have already met pseudoline arrangements in Chapter~\ref{chap:mpt}, which correspond to point configurations in a (topological) plane. We refer to~\cite{bvswz-om-99,b-com-06,rgz-om-97} for expositions on the different aspects of oriented matroids, and we focus here on a short overview of what we will really use in the next sections.

Among other possible definitions of oriented matroids, we use the concept of chirotopes:

\begin{definition}\label{ft:def:chirotopepointconfig}
\index{chirotope}
The \defn{chirotope} of a point set~$P \eqdef \ens{p_i}{i\in I}\subset\R^d$ indexed by a set~$I$ is the application~$\chi_P:I^{d+1}\to\{-1,0,1\}$ which associates to each \tuple{(d+1)}~$(i_0,\dots,i_d)$ of elements of~$I$ the \defn{orientation} of the simplex formed by the points of~$P$ indexed by~$i_0,\dots,i_d$, that~is:
$$\chi_P(i_0,\dots,i_d) \eqdef \sign\det\begin{pmatrix} 1 & 1 & \dots & 1 \\ p_{i_0} & p_{i_1} & \dots & p_{i_d} \end{pmatrix}.$$
\end{definition}

\begin{figure}[!h]
	\capstart
	\centerline{\includegraphics[scale=1.2]{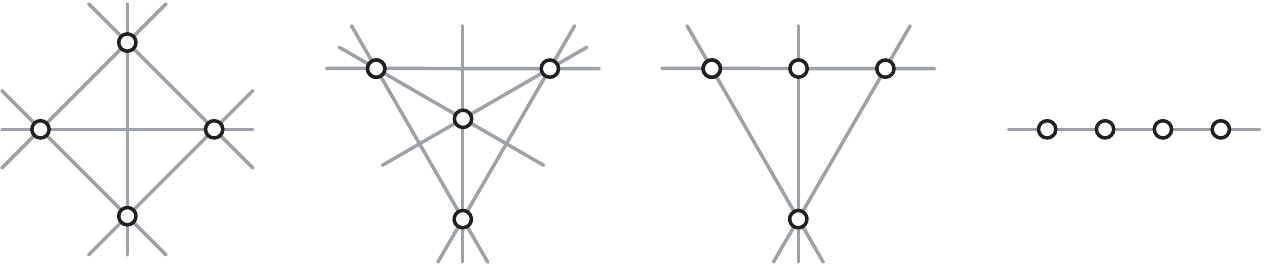}}
	\caption[The four possible order types on configurations of four (distinct) points]{The four possible order types on configurations of four (distinct) points.}
	\label{ft:fig:4points}
\end{figure}

Two points configurations~$P$ and~$Q$ indexed by the same set~$I$ have the same \defn{order type} if there exists a bijection~$\alpha:I\to I$ which makes coincide the chirotopes of~$P$~and~$Q$, \ie such that $\chi_P(i_0,\dots,i_d)=\chi_Q(\alpha(i_0),\dots,\alpha(i_d))$ for any \tuple{(d+1)}~$(i_0,\dots,i_d)\in I^{d+1}$. For example, there are only four different order types of configurations of four points in the plane, represented in \fref{ft:fig:4points}.

From the geometry of point configurations, we derive the following combinatorial properties of chirotopes:

\begin{lemma}\label{ft:lem:propertieschirotope}
The chirotope $\chi:I^{d+1}\to\{-1,0,1\}$ of a \dimensional{d} point configuration indexed by~$I$ satisfies the following relations:
\begin{enumerate}[(i)]
\item \defn{Alternating relations}:
for any \tuple{(d+1)}~$(i_0,\dots,i_d)\in I^{d+1}$ and any permutation~$\pi$ of~$\{0,\dots,d\}$ of signature~$\sigma$, 
$$\chi(i_{\pi(0)},\dots,i_{\pi(d)})=\sigma\chi(i_0,\dots,i_d).$$
\item \defn{Matroid property}: the set~$\chi^{-1}(\{-1,1\})$ of non-degenerate \tuple{(d+1)}s of~$I$ satisfies the \defn{Steinitz exchange axiom}: for any non-degenerate \tuple{(d+1)}s $X,Y\in\chi^{-1}(\{-1,1\})$ and any element~$x\in X\ssm Y$, there exists~$y\in Y\ssm X$ such that $X\diffsym\{x,y\}\in\chi^{-1}(\{-1,1\})$.
\item \defn{\GP relations}\index{Grassman@\GP relations|hbf}:
for any~$i_0,\dots,i_{d-2},j_1,j_2,j_3,j_4\in I$, the set 
\begin{eqnarray*}
\{ & \chi(i_0,\dots,i_{d-2},j_1,j_2)\,\chi(i_0,\dots,i_{d-2},j_3,j_4), \\
 & -\chi(i_0,\dots,i_{d-2},j_1,j_3)\,\chi(i_0,\dots,i_{d-2},j_2,j_4), \\
 & \chi(i_0,\dots,i_{d-2},j_1,j_4)\,\chi(i_0,\dots,i_{d-2},j_2,j_3) & \}
\end{eqnarray*}
either contains~$\{-1,1\}$ or is contained in~$\{0\}$.
\end{enumerate}
\end{lemma}

\begin{proof}
Alternating relations come from alternating properties of the determinant. The set of non-degenerate \tuple{(d+1)}s corresponds to the set of bases of our point set, and thus satisfies the matroid property. The \GP relations are derived from the $3$-terms \GP equality on determinants: for any matrix~$M$ with~$d+1$~rows and~$d+3$~columns,
$$\det(M_{\hat1\hat2})\det(M_{\hat3\hat4})-\det(M_{\hat1\hat3})\det(M_{\hat2\hat4})+\det(M_{\hat1\hat4})\det(M_{\hat2\hat3})=0,$$
where~$M_{\hat p\hat q}$ denotes the submatrix of~$M$ obtained by deleting its $p$th and $q$th columns. In turn, this equality is a simple calculation when multiplying~$M$ by the inverse of the bottom square submatrix~$M_{\hat1\hat2}$.
\end{proof}

These properties are the basement of the combinatorial definition of chirotopes: we now completely forget point configurations and define chirotopes as follows.

\begin{definition}\label{ft:def:chirotopes}
\index{chirotope}
A \defn{chirotope} of rank~$r$ on a set~$I$ is any application~$\chi:I^r\to\{-1,0,1\}$ which is not identically zero and satisfies the three properties of Lemma~\ref{ft:lem:propertieschirotope}.
\end{definition}

By definition, the chirotope of a (non-degenerate) \dimensional{d} point configuration is a chirotope of rank~$d+1$. However, for the same reason that some pseudoline arrangements are not realizable, there exist chirotopes which do not arise as chirotopes of point configurations. These two sentences together give a method as well as its limits to look for a (small) point configuration satisfying certain combinatorial properties: we can first enumerate all chirotopes which satisfy the constraints, and then try to realize these chirotopes as euclidean point configurations.

\svs
We will apply this method to obtain all symmetric realizations of a simplicial complex.We first express the combinatorial constraints we are focussing on:

\begin{lemma}\label{ft:lem:chirotoperealization}
Let~$\Delta\subset 2^I$ be a simplicial complex (whose vertices are indexed by~$I$), and~$G$ be a group acting on~$\Delta$.
Let~$P \eqdef \ens{p_i}{i\in I}\subset\R^d$ be a \dimensional{d} point configuration (indexed by~$I$) whose convex hull is a symmetric realization of~$\Delta$ under~$G$. Then the chirotope~$\chi_P$ of~$P$ satisfies the following properties:
\begin{enumerate}[(i)]
\item \defn{Alternating relations}.
\item \defn{Matroid property}.
\item \defn{\GP relations}.
\item \defn{Necessary Simplex Orientations}:
if~$i_0,\dots,i_d\in I$ are such that both~$\{i_0,\dots,i_{d-2},i_{d-1}\}$ and~$\{i_0,\dots,i_{d-2},i_{d}\}$ are facets of~$\Delta$, then for any~$j\in I\ssm\{i_0,\dots,i_d\}$,
$$\chi_P(i_0,\dots,i_{d-2},i_{d-1},i_d)=\chi_P(i_0,\dots,i_{d-2},i_{d-1},j)=\chi_P(i_0,\dots,i_{d-2},j,i_d).$$
\item \defn{Symmetry}:
there exists a morphism~$\tau: G\to \{\pm1\}$ such that for any~$i_0,\dots,i_d\in I$, and any~$g\in G$,
$$\chi_P(gi_0,\dots,gi_d)=\tau(g)\chi_P(i_0,\dots,i_d).$$
\end{enumerate}
\vspace{-.6cm}\qed
\end{lemma}

As before, these properties lead to the following combinatorial abstraction:

\begin{definition}\label{ft:def:chirotoperealization}
\index{realization!symmetric matroid ---|hbf}
A \defn{symmetric matroid realization} of~$\Delta$ under~$G$ is any application which associates to a \tuple{(d+1)} of~$I$ a sign in~$\{-1,0,1\}$ and satisfies the five properties of Lemma~\ref{ft:lem:chirotoperealization}.
\end{definition}

In the following paragraphs, we investigate symmetric matroid realizations of~$\Delta_{n,k}$ as the first step to understand its symmetric geometric realizations. This enables us to describe the space of symmetric realizations of little cases such as the first non-trivial case~$\Delta_{8,2}$. To expose the method, we start with a simpler example. 

\subsubsection{A basic example: the space of symmetric realizations of~$\Delta_{6,1}$}\label{ft:subsubsec:delta61}

As our first example of the use of oriented matroids, we study symmetric realizations of the polar of the \dimensional{3} associahedron. It is convenient to use letters $\ra$, $\rb$, $\rc$, $\rd$, $\re$, $\rf$, $\rG$, $\rH$ and $\rI$ to denote the edges $[0,2]$, $[1,3]$, $[2,4]$, $[3,5]$, $[4,0]$, $[5,1]$, $[0,3]$, $[1,4]$ and $[2,5]$ of the hexagon respectively (the three capital letters correspond to the long diagonals). We use these letters to index the vertices of~$\Delta_{6,1}$.

Since the polar~$\Sigma_6^\polar$ of the secondary polytope of the regular hexagon is a symmetric realization of~$\Delta_{6,1}$, its chirotope~$\chi_{\Sigma_6^\polar}$ provides a symmetric matroid realization. To describe this chirotope, instead of listing the signs of all possible \tuple{4}s of~$\{\ra,\rb,\dots,\rI\}$, we give the sign of one representative for each orbit under permutation of the coordinates and under the action~of~$\D_6$:
\begin{center}
\begin{tabular}{c@{$\,=\,$}c@{\quad}c@{$\,=\,$}c@{\quad}c@{$\,=\,$}c@{\quad}c@{$\,=\,$}c@{\quad}c@{$\,=\,$}c@{\quad}c@{$\,=\,$}c}
$|\ra\rb\rc\rd|$&$1$ & $|\ra\rb\rc\re|$&$-1$ & $|\ra\rb\rc\rG|$&$-1$ & $|\ra\rb\rc\rI|$&$1$ & $|\ra\rb\rd\re|$&$0$ & $|\ra\rb\rd\rG|$&$1$ \\
$|\ra\rb\rd\rH|$&$-1$ & $|\ra\rb\rd\rI|$&$-1$ & $|\ra\rb\rG\rH|$&$1$ & $|\ra\rb\rH\rI|$&$1$ & $|\ra\rc\re\rG|$&$-1$ & $|\ra\rc\rG\rH|$&$0$ \\
$|\ra\rc\rG\rI|$&$-1$ & $|\ra\rd\rG\rH|$&$1$ & $|\ra\rd\rG\rI|$&$1$ & $|\ra\rG\rH\rI|$&$1$
\end{tabular}
\end{center}
where~$|\rx\ry\rz\rt|$ is an abbreviation for~$\chi_{\Sigma_6^\polar}(\rx,\ry,\rz,\rt)$. You can recover all the chirotope~$\chi_{\Sigma_6^\polar}$ by remembering that the reflection~$\rho_v$ (resp.~$\rho_h$) of the hexagon, with respect to the vertical (resp.~horizontal) axis, corresponds to a direct (resp.~indirect) isometry of~$\Sigma_6^\polar$. In other words,~$\tau(\rho_v)=1$ and~$\tau(\rho_h)=-1$.

Using our computer enumeration (see Appendix~\ref{app:sec:enumerationmatroidpolytopes} for a presentation of the \haskell implementation), we obtain that the chirotope given by the secondary polytope is in fact the unique solution:

\begin{proposition}
The chirotope~$\chi_{\Sigma_6^\polar}$ is the unique symmetric matroid realization of~$\Delta_{6,1}$ (up to the inversion of all the signs).\qed
\end{proposition}

It remains to understand all possible symmetric geometric realizations of this chirotope (for the moment, we know only one such realization: the polar of the secondary polytope itself). Given such a realization, consider the matrix~$M$ with 4 rows and 9 columns formed by the homogeneous coordinates of its vertices. Let~$N$ be the square submatrix of~$M$ formed by the columns~$\ra$, $\rG$, $\rH$ and $\rI$. Since~$|\ra\rG\rH\rI|=1$, we know that~$N$ is invertible, and we denote
$$M'  \eqdef  N^{-1}M  \eqdef  \begin{pmatrix} 1 & b_0 & c_0 & d_0 & e_0 & f_0 & 0 & 0 & 0 \\ 0 & b_1 & c_1 & d_1 & e_1 & f_1 & 1 & 0 & 0 \\ 0 & b_2 & c_2 & d_2 & e_2 & f_2 & 0 & 1 & 0 \\ 0 & b_3 & c_3 & d_3 & e_3 & f_3 & 0 & 0 & 1 \end{pmatrix}.$$
We determine the unknown coefficients of this matrix, using the fact that not only the signs of the determinants of the $(4\times4)$-submatrices of~$M'$ (\ie the chirotope) are symmetric under the action of~$\D_6$, but also the determinants themselves. For example, applying the rotation of the hexagon of angle~$\pi/3$ (which acts as an indirect isometry on the chirotope), we know that
$$[\ra\rG\rH\rI]=-[\rb\rG\rH\rI]=[\rc\rG\rH\rI]=-[\rd\rG\rH\rI]=[\re\rG\rH\rI]=-[\rf\rG\rH\rI],$$
(where~$[\rx\ry\rz\rt]$ denotes the determinant of the square submatrix of~$M'$ formed by the columns $\rx$, $\ry$, $\rz$ and $\rt$). Thus, we derive that $1=-b_0=c_0=-d_0=e_0=-f_0$. Similarly, we obtain the following relations between the coefficients of the matrix~$M'$:
\begin{align*}
[\ra\rb\rH\rI]=[\ra\rf\rG\rH]\quad & \Longrightarrow\quad b_1=f_3; \\
[\ra\rb\rG\rH]=-[\ra\rb\rG\rI]=[\ra\rf\rG\rI]=-[\ra\rf\rH\rI]=-[\rb\rc\rH\rI]\quad & \Longrightarrow\quad b_3=b_2=f_2=f_1=b_1+c_1; \\
[\ra\rc\rG\rI]=[\ra\rc\rH\rI]=[\ra\rf\rG\rH]=[\ra\rf\rG\rI]\quad & \Longrightarrow\quad -c_2=c_1=e_3=-e_2; \\
[\ra\rc\rG\rH]=[\ra\re\rH\rI]=[\rb\rd\rH\rI]=0\quad & \Longrightarrow\quad c_3=e_1=b_1-d_1=0; \\
\text{and}\qquad[\ra\rd\rG\rH]=[\ra\rd\rH\rI]\quad & \Longrightarrow\quad d_3=d_1.
\end{align*}
Furthermore, in each column of~$M'$, the coordinates sum to~$1$ (since we have homogeneous coordinates when multipling by the matrix~$N$). This implies that~$b_1=2-2b_2$ and~$d_2=2-2d_1$. Thus, the matrix~$M'$ can be written
$$M' = \begin{pmatrix} 1 & -1 & 1 & -1 & 1 & -1 & 0 & 0 & 0 \\ 0 & -2t+2 & 3t-2 & -2t+2 & 0 & t & 1 & 0 & 0 \\ 0 & t & -3t+2 & 4t-2 & -3t+2 & t & 0 & 1 & 0 \\ 0 & t & 0 & -2t+2 & 3t-2 & -2t+2 & 0 & 0 & 1 \end{pmatrix}.$$
Finally, from our knowledge on the chirotope, we know that~$|\ra\rb\rG\rH|=|\ra\rG\rH\rI|=1$ and that $|\ra\rd\rG\rI|=|\ra\rG\rH\rI|=-1$, which implies that~$0<t<1/2$.

\svs
In order to complete our understanding of the space of symmetric realizations of~$\Delta_{6,1}$, it only remains to study the possible values of the matrix~$N$. By symmetry, the triangle formed by the vertices labeled by~$\rG$, $\rH$ and $\rI$ is equilateral. Thus, since a dilation (the composition of an homothecy by an isometry) does not destruct the symmetry, we can assume without loss of generality that this triangle is formed by the vectors of the canonical basis in~$\R^3$, \ie that the matrix~$N$ is of the form
$$N \eqdef \begin{pmatrix} 1 & 1 & 1 & 1 \\ x & 1 & 0 & 0 \\ y & 0 & 1 & 0 \\ z & 0 & 0 & 1 \end{pmatrix}.$$
Since~$\ra$ is equidistant from~$\rG$ and~$\rI$, we have $x=z$, and since~$\rb$ is equidistant from~$\rG$ and~$\rH$, we have~$y=x+3t-2$. Finally, since~$N$ is invertible, its determinant $\det N=3-3x-3t$ does not vanish, and~$x\ne 1-t$. 

\svs
Knowing the matrices~$M'$ and~$N$, we have a complete description of the matrix~$M$:

\begin{proposition}
Any symmetric realization of the simplicial complex~$\Delta_{6,1}$ is a dilation of the set of column vectors of the matrix
$$\begin{pmatrix} x & -x-2t+2 & x+3t-2 & -x-2t+2 & x & -x+t & 1 & 0 & 0 \\ x+3t-2 & -x-2t+2 & x & -x+t & x & -x-2t+2 & 0 & 1 & 0 \\ x & -x+t & x & -x-2t+2 & x+3t-2 & -x-2t+2 & 0 & 0 & 1\end{pmatrix}$$
for a couple $(t,x)\in\R^2$, with $0<t<1/2$ and~$x\ne1-t$. Reciprocally, any such dilation is a symmetric realization of~$\Delta_{6,1}$.\qed
\end{proposition}

The reciprocal part of this proposition is easy to check (for example with a computer algebra system): we only have to verify that for any admissible values of~$t$ and~$x$:
\begin{enumerate}[(i)]
\item The boundary complex of the convex hull of the column vectors of~$M$ is indeed isomorphic to~$\Delta_{6,1}$. For this, we check that the chirotope of the column vectors of~$M'$ coincides with the chirotope~$\chi_{\Sigma_6^\polar}$ (and we only have to check equality for one representative of each of the~$16$ orbits under~$\D_6$).
\item The dihedral group~$\D_6$ indeed acts by isometry on the columns vectors of~$M$ (and we only have to check the action of the two reflections~$\rho_v$ and~$\rho_h$, since they generate~$\D_6$).
\end{enumerate}

\begin{corollary}
Up to dilation, the space of symmetric realizations of the simplicial complex~$\Delta_{6,1}$ has dimension~$2$.\qed
\end{corollary}

\begin{figure}[!h]
	\capstart
	\centerline{\includegraphics[scale=1.2]{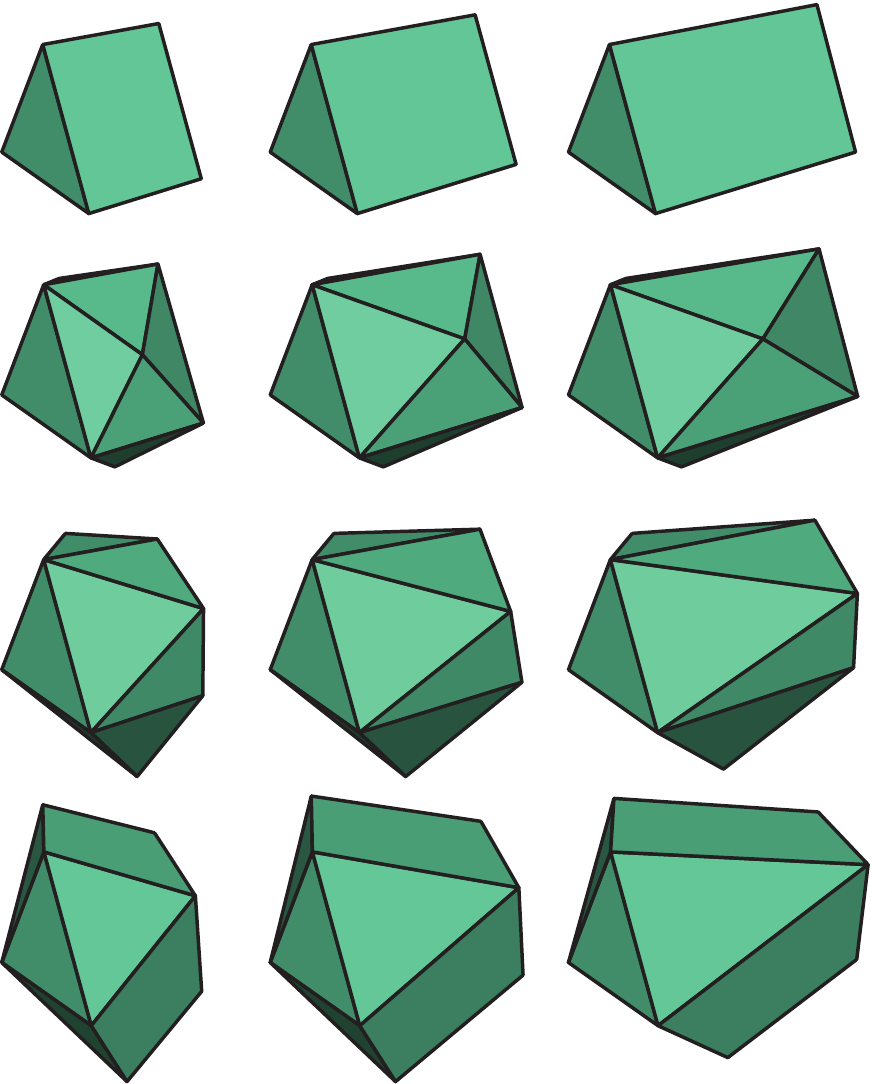}}
	\caption[The space of symmetric realizations of~$\Delta_{6,1}$]{The space of symmetric realizations of~$\Delta_{6,1}$.}
	\label{ft:fig:associahedronparameters}
\end{figure}

To illustrate this result, \fref{ft:fig:associahedronparameters} shows the convex hull of the column vectors of~$M$ for different values of~$t$ and~$x$. The two dimensions of the space of symmetric realizations are apparent in this picture: if you see the associahedron as three egyptian pyramids stacked on the square facets of a triangular prism, you can variate the height of the prism and the height of the pyramids. The bounds on the possible values of~$t$ correspond to the moment when the pyramids disappear ($t=0$) and when they merge ($t=1/2$). See \fref{ft:fig:associahedronparameters}.

\subsubsection{The first non-trivial case: the space of symmetric realizations of~$\Delta_{8,2}$}\label{ft:subsubsec:delta82}

Using the same method, we now study the symmetric realizations of the simplicial complex~$\Delta_{8,2}$ of \kcross{3}-free sets of \krel{2} edges of the octagon. The polytope we want to construct would be a \poly{6}tope, with~$12$~vertices (the \krel{2} edges of the octagon), $66$~edges (because $\Delta_{8,2}$ is \neighborly{2}),~$192$~\face{2}s,~$306$~\face{3}s, $252$~ridges, and~$84$~facets (the \ktri{2}s of the octagon). As before, it is convenient to label the \krel{2} edges of the octagon with letters: let $\ra \eqdef [0,3]$, $\rb \eqdef [1,4]$, $\rc \eqdef [2,5]$, $\rd \eqdef [3,6]$, $\re \eqdef [4,7]$, $\rf \eqdef [5,0]$, $\rg \eqdef [6,1]$, $\rh \eqdef [7,2]$, $\rI \eqdef [0,4]$, $\rJ \eqdef [1,5]$, $\rK \eqdef [2,6]$, and $\rL \eqdef [3,7]$ (the capital letters denote the long diagonals).

As before, a computer enumeration (see Appendix~\ref{app:sec:enumerationmatroidpolytopes}) provides the complete list of all possible symmetric matroid realizations of our simplicial complex: 

\begin{proposition}\label{ft:prop:delta82matroid}
There are~$15$ symmetric matroid realizations of~$\Delta_{8,2}$ (up to the inversion of all the signs).\qed
\end{proposition}

To present these solutions, we only give the sign of one representative for each of the~$62$~orbits of~$\{\ra,\rb,\dots,\rL\}^7$ under permutation and symmetry. First, all solutions have the following~$59$ common signs:
\begin{center}
\begin{tabular}{@{}c@{$\,=\,$}c@{\quad}c@{$\,=\,$}c@{\quad}c@{$\,=\,$}c@{\quad}c@{$\,=\,$}c@{\quad}c@{$\,=\,$}c@{}}
$|\ra\rb\rc\rd\re\rf\rg|$&$0$ & $|\ra\rb\rc\rd\re\rf\rI|$&$-1$ & $|\ra\rb\rc\rd\re\rf\rJ|$&$1$ & $|\ra\rb\rc\rd\re\rf\rK|$&$-1$ & $|\ra\rb\rc\rd\re\rg\rI|$&$1$ \\
$|\ra\rb\rc\rd\re\rg\rJ|$&$-1$ & $|\ra\rb\rc\rd\re\rI\rK|$&$-1$ & $|\ra\rb\rc\rd\re\rI\rL|$&$1$ & $|\ra\rb\rc\rd\re\rJ\rK|$&$1$ & $|\ra\rb\rc\rd\rf\rg\rI|$&$-1$ \\
$|\ra\rb\rc\rd\rf\rg\rJ|$&$1$ & $|\ra\rb\rc\rd\rf\rg\rL|$&$1$ & $|\ra\rb\rc\rd\rf\rI\rJ|$&$1$ & $|\ra\rb\rc\rd\rf\rI\rK|$&$1$ & $|\ra\rb\rc\rd\rf\rI\rL|$&$-1$ \\
$|\ra\rb\rc\rd\rf\rJ\rK|$&$-1$ & $|\ra\rb\rc\rd\rf\rJ\rL|$&$1$ & $|\ra\rb\rc\rd\rf\rK\rL|$&$-1$ & $|\ra\rb\rc\rd\rI\rJ\rK|$&$-1$ & $|\ra\rb\rc\rd\rI\rJ\rL|$&$1$ \\
$|\ra\rb\rc\rd\rI\rK\rL|$&$-1$ & $|\ra\rb\rc\re\rf\rg\rI|$&$0$ & $|\ra\rb\rc\re\rf\rg\rK|$&$0$ & $|\ra\rb\rc\re\rf\rI\rJ|$&$-1$ & $|\ra\rb\rc\re\rf\rI\rK|$&$1$ \\
$|\ra\rb\rc\re\rf\rI\rL|$&$1$ & $|\ra\rb\rc\re\rf\rJ\rK|$&$-1$ & $|\ra\rb\rc\re\rf\rJ\rL|$&$-1$ & $|\ra\rb\rc\re\rf\rK\rL|$&$1$ & $|\ra\rb\rc\re\rg\rI\rJ|$&$1$ \\
$|\ra\rb\rc\re\rg\rI\rK|$&$-1$ & $|\ra\rb\rc\re\rg\rI\rL|$&$-1$ & $|\ra\rb\rc\re\rg\rK\rL|$&$-1$ & $|\ra\rb\rc\re\rI\rJ\rL|$&$-1$ & $|\ra\rb\rc\re\rI\rK\rL|$&$1$ \\
$|\ra\rb\rc\re\rJ\rK\rL|$&$-1$ & $|\ra\rb\rc\rf\rI\rJ\rK|$&$-1$ & $|\ra\rb\rc\rf\rI\rK\rL|$&$-1$ & $|\ra\rb\rc\rI\rJ\rK\rL|$&$1$ & $|\ra\rb\rd\re\rg\rI\rJ|$&$1$ \\
$|\ra\rb\rd\re\rg\rI\rK|$&$1$ & $|\ra\rb\rd\re\rg\rI\rL|$&$1$ & $|\ra\rb\rd\re\rg\rJ\rK|$&$-1$ & $|\ra\rb\rd\re\rI\rJ\rK|$&$1$ & $|\ra\rb\rd\re\rI\rJ\rL|$&$-1$ \\
$|\ra\rb\rd\rf\rI\rJ\rL|$&$1$ & $|\ra\rb\rd\rf\rI\rK\rL|$&$1$ & $|\ra\rb\rd\rf\rJ\rK\rL|$&$-1$ & $|\ra\rb\rd\rg\rI\rJ\rK|$&$-1$ & $|\ra\rb\rd\rg\rI\rJ\rL|$&$-1$ \\
$|\ra\rb\rd\rg\rJ\rK\rL|$&$1$ & $|\ra\rb\rd\rI\rJ\rK\rL|$&$-1$ & $|\ra\rb\re\rf\rI\rJ\rK|$&$-1$ & $|\ra\rb\re\rf\rI\rJ\rL|$&$-1$ & $|\ra\rb\re\rf\rJ\rK\rL|$&$1$ \\
$|\ra\rb\re\rI\rJ\rK\rL|$&$1$ & $|\ra\rc\re\rg\rI\rJ\rK|$&$-1$ & $|\ra\rc\re\rI\rJ\rK\rL|$&$-1$ & $|\ra\rc\rf\rI\rJ\rK\rL|$&$1$
\end{tabular}
\end{center}

\svs
The three remaining orbits are those of the tuples $(\ra,\rb,\rc,\rd,\re,\rI,\rJ)$, $(\ra,\rb,\rc,\re,\rI,\rJ,\rK)$ and $(\ra,\rb,\rd,\rf,\rI,\rJ,\rK)$. The following table summarizes the possible signs for these three orbits (only~$15$ of the~$27$ possibilities are admissible):
\begin{center}
\begin{tabular}{|c|c|c|c|c|c|c|c|c|c|c|c|c|c|c|c|}
\hline
 & $A$ & $B$ & $C$ & $D$ & $E$ & $F$ & $G$ & $H$ & $I$ & $J$ & $K$ & $L$ & $M$ & $N$ & $O$ \\
\hline
$|\ra\rb\rc\rd\re\rI\rJ|$ & $-1$ & $-1$ & $-1$ & $0$ & $0$ & $0$ & $1$ & $1$ & $1$ & $1$ & $1$ & $1$ & $1$ & $1$ & $1$  \\
\hline
$|\ra\rb\rc\re\rI\rJ\rK|$ & $1$ & $1$ & $1$ & $1$ & $1$ & $1$ & $-1$ & $-1$ & $-1$ & $0$ & $0$ & $0$ & $1$ & $1$ & $1$ \\
\hline
$|\ra\rb\rd\rf\rI\rJ\rK|$ & $-1$ & $0$ & $1$ & $-1$ & $0$ & $1$ & $-1$ & $0$ & $1$ & $-1$ & $0$ & $1$ & $-1$ & $0$ & $1$ \\
\hline
\end{tabular}
\end{center}

\svs
The second step is to realize geometrically these matroid polytopes. Given such a realization, we define again the matrix~$M$ with 7 rows and 12 columns formed by the homogeneous coordinates of its vertices. Let~$N$ be the square submatrix of~$M$ formed by the columns~$\ra$, $\rb$, $\rc$, $\rI$, $\rJ$, $\rK$ and $\rL$. Since~$|\ra\rb\rc\rI\rJ\rK\rL|=1$, we know that~$N$ is invertible, and we denote
$$M'  \eqdef  N^{-1}M  \eqdef  \begin{pmatrix} 1 & 0 & 0 & d_0 & e_0 & f_0 & g_0 & h_0 & 0 & 0 & 0 & 0 \\ 0 & 1 & 0 & d_1 & e_1 & f_1 & g_1 & h_1 & 0 & 0 & 0 & 0 \\ 0 & 0 & 1 & d_2 & e_2 & f_2 & g_2 & h_2 & 0 & 0 & 0 & 0 \\ 0 & 0 & 0 & d_3 & e_3 & f_3 & g_3 & h_3 & 1 & 0 & 0 & 0 \\ 0 & 0 & 0 & d_4 & e_4 & f_4 & g_4 & h_4 & 0 & 1 & 0 & 0 \\ 0 & 0 & 0 & d_5 & e_5 & f_5 & g_5 & h_5 & 0 & 0 & 1 & 0 \\ 0 & 0 & 0 & d_6 & e_6 & f_6 & g_6 & h_6 & 0 & 0 & 0 & 1 \end{pmatrix} \eqdef \begin{pmatrix} I_3 & T & 0_{3\times4} \\ 0_{4\times3} & B & I_4 \end{pmatrix},$$
where $T\in \R^{3\times 5}$ and $B\in \R^{4\times5}$ denote respectively the ``top'' and ``bottom'' unknown submatrices. We use the symmetry of the determinants of the $(7\times 7)$-submatrices of $M'$ under the action of the dihedral group to  determine the unknown coefficients in the matrices~$T$ and~$B$.

\begin{lemma}\label{ft:lem:T}
The matrix~$T$ equals
$$\begin{pmatrix}
-1 &  1+\sqrt{2} & -2-\sqrt{2} & 2+\sqrt{2} & -1-\sqrt{2} \\
-1-\sqrt{2} & 2+2\sqrt{2} & -3-2\sqrt{2} & 2+2\sqrt{2} & -1-\sqrt{2} \\
-1-\sqrt{2} & 2+\sqrt{2} & -2-\sqrt{2} & 1+\sqrt{2}  & -1
\end{pmatrix}.$$
\end{lemma}

\begin{proof}
From symmetry, we derive the following equalities between the determinants of the submatrices of~$M'$:
\begin{align*}
[\ra\rb\rc\rI\rJ\rK\rL]=-[\rb\rc\rd\rI\rJ\rK\rL]=-[\ra\rb\rh\rI\rJ\rK\rL]\quad & \Longrightarrow\quad 1=-d_0=-h_2; \\
[\ra\rb\rd\rI\rJ\rK\rL]=-[\ra\rb\rg\rI\rJ\rK\rL]=-[\ra\rc\rd\rI\rJ\rK\rL]=-[\ra\rc\rh\rI\rJ\rK\rL] & =-[\rb\rc\re\rI\rJ\rK\rL]=[\rb\rc\rh\rI\rJ\rK\rL] \\
\Longrightarrow\quad d_2= & -g_2=d_1=h_1=-e_0=h_0 \eqfed \alpha; \\
[\ra\rb\re\rI\rJ\rK\rL]=-[\ra\rb\rf\rI\rJ\rK\rL]=-[\rb\rc\rf\rI\rJ\rK\rL]=[\rb\rc\rg\rI\rJ\rK\rL]\quad & \Longrightarrow\quad e_2=-f_2=-f_0=g_0 \eqfed \beta; \\
\text{and}\qquad[\ra\rc\re\rI\rJ\rK\rL]=[\ra\rc\rg\rI\rJ\rK\rL]\quad & \Longrightarrow\quad -e_1=-g_1 \eqfed -\gamma.
\end{align*}

Thus, we can already write the matrix~$T$ as follows:
$$T=\begin{pmatrix} -1 & -\alpha & -\beta & \beta & \alpha \\  \alpha & \gamma & \delta & \gamma & \alpha \\ \alpha & \beta & -\beta & -\alpha & -1 \end{pmatrix}.$$

Furthermore, 
\begin{align*}
[\rb\rd\rh\rI\rJ\rK\rL]=-[\ra\rc\re\rI\rJ\rK\rL]\quad & \Longrightarrow\quad \alpha^2-1=\gamma; \\
\text{and}\qquad[\ra\rd\rh\rI\rJ\rK\rL]=-[\ra\rb\re\rI\rJ\rK\rL]\quad & \Longrightarrow\quad -\alpha^2-\alpha=-\beta.
\end{align*}

Thus, the matrix~$T$ can be written:
$$T=\begin{pmatrix} -1 & -\alpha & -\alpha^2-\alpha & \alpha^2+\alpha & \alpha \\  \alpha & \alpha^2-1 & \delta & \alpha^2-1 & \alpha \\ \alpha & \alpha^2+\alpha & -\alpha^2-\alpha & -\alpha & -1 \end{pmatrix}.$$

Finally,
\begin{align*}
[\rd\re\rf\rI\rJ\rK\rL]=-[\ra\rb\rc\rI\rJ\rK\rL]\quad & \Longrightarrow\quad -\alpha^4-2\alpha^3-2\alpha^2-\alpha+\delta\alpha=-1 \\
& \Longrightarrow\quad \delta=\frac{1}{\alpha}(\alpha^4+2\alpha^3+2\alpha^2+\alpha-1); \\
\text{and}\qquad[\rd\re\rh\rI\rJ\rK\rL]=-[\ra\rb\re\rI\rJ\rK\rL]\quad & \Longrightarrow\quad \alpha^3+2\alpha^2-1=-\alpha^2-\alpha \\
& \Longrightarrow\quad\alpha\in\{-1,-1+\sqrt{2},-1-\sqrt{2}\}.
\end{align*}

According to our description of all symmetric matroid realizations, we have~$|\ra\rd\rg\rI\rJ\rK\rL|\ne 0$ and~$|\ra\rb\rd\rI\rJ\rK\rL|<0$. Thus,~$\alpha\ne -1$ and~$\alpha<0$. Consequently,~$\alpha=-1-\sqrt{2}$, which completes the proof of the lemma.
\end{proof}

\begin{lemma}\label{ft:lem:B}
There exists~$u\in(-1-2\sqrt{2},-1-3\sqrt{2}/2)$ such that the matrix~$B$ equals
{\small
$$\begin{pmatrix}
1+\frac{\sqrt{2}}{2} & u & -(1+\sqrt{2})(1+u) & (1+\sqrt{2})(1+u)+1 & -\frac{\sqrt{2}}{2}-u \\
-\frac{\sqrt{2}}{2}-u & (1+\sqrt{2})(1+u)+1 & -(1+\sqrt{2})(1+u) & u & 1+\frac{\sqrt{2}}{2} \\
1+\frac{\sqrt{2}}{2} & -2-2\sqrt{2}-u & (1+\sqrt{2})(3+u)+2 & -(1+\sqrt{2})(3+u)-1 & 2+\frac{3\sqrt{2}}{2}+u \\
2+\frac{3\sqrt{2}}{2}+u & -(1+\sqrt{2})(3+u)-1 & (1+\sqrt{2})(3+u)+2 & -2-2\sqrt{2}-u & 1+\frac{\sqrt{2}}{2} 
\end{pmatrix}.$$
}
\end{lemma}

\begin{proof}
From the symmetry, we derive the following relations:
\begin{align*}
[\ra\rb\rc\rd\rI\rJ\rK]=-[\ra\rb\rc\rh\rI\rJ\rL]\quad & \Longrightarrow\quad -d_6=-h_5; \\
[\ra\rb\rc\rd\rI\rJ\rL]=[\ra\rb\rc\rd\rJ\rK\rL]=-[\ra\rb\rc\rh\rI\rJ\rK]=- & [\ra\rb\rc\rh\rI\rK\rL]\quad\Longrightarrow\quad d_5=d_3=h_6=h_4; \\
[\ra\rb\rc\rd\rI\rK\rL]=-[\ra\rb\rc\rh\rJ\rK\rL]\quad & \Longrightarrow\quad -d_4=-h_3; \\
[\ra\rb\rc\re\rI\rJ\rK]=-[\ra\rb\rc\rg\rI\rJ\rL]\quad & \Longrightarrow\quad -e_6=-g_5; \\
[\ra\rb\rc\re\rI\rJ\rL]=-[\ra\rb\rc\rg\rI\rJ\rK]\quad & \Longrightarrow\quad e_5=g_6; \\
[\ra\rb\rc\re\rI\rK\rL]=-[\ra\rb\rc\rg\rJ\rK\rL]\quad & \Longrightarrow\quad -e_4=-g_3; \\
[\ra\rb\rc\re\rJ\rK\rL]=-[\ra\rb\rc\rg\rI\rK\rL]\quad & \Longrightarrow\quad e_3=g_4; \\
[\ra\rb\rc\rf\rI\rJ\rK]=-[\ra\rb\rc\rf\rI\rJ\rL]\quad & \Longrightarrow\quad -f_6=-f_5; \\
[\ra\rb\rc\rf\rI\rK\rL]=-[\ra\rb\rc\rf\rJ\rK\rL]\quad & \Longrightarrow\quad -f_4=-f_3.
\end{align*}

With these equalities, the matrix~$B$ can be written:
$$B=\begin{pmatrix} d_3 & e_3 & f_3 & e_4 & d_4 \\ d_4 & e_4 & f_3 & e_3 & d_3 \\ d_3 & e_5 & f_5 & e_6 & d_6 \\ d_6 & e_6 & f_5 & e_5 & d_3 \end{pmatrix}.$$

Furthermore, for any symmetric matroid realization,~$d_3=|\ra\rb\rc\rd\rJ\rK\rL|\ne 0$. Thus:
\begin{align*}
[\ra\rb\rc\rd\re\rK\rL]=-[\ra\rb\rc\rd\rh\rI\rL]\quad & \Longrightarrow\quad d_3e_4-d_4e_3=d_3^2-d_4d_6 \quad \Longrightarrow\quad e_4=\frac{1}{d_3}(d_3^2-d_4d_6+d_4e_3); \\
[\ra\rb\rc\rd\re\rJ\rL]=[\ra\rb\rc\rd\rh\rI\rK]\quad & \Longrightarrow\quad d_3(e_3-e_5)=d_3(d_6-d_4) \quad \Longrightarrow\quad e_5=e_3+d_4-d_6; \\
[\ra\rb\rc\rd\re\rJ\rK]=[\ra\rb\rc\rd\rh\rI\rJ]\quad & \Longrightarrow\quad d_3e_6-d_6e_3=d_3^2-d_6^2 \quad \Longrightarrow\quad e_6=\frac{1}{d_3}(d_3^2-d_6^2+d_6e_3).
\end{align*}

We also know that~$|\ra\rb\rc\rd\rh\rK\rL|\ne 0$ and~$|\ra\rb\rc\rd\rh\rI\rJ|\ne 0$. Thus,~$d_3\ne d_4$ and~$d_3\ne d_6$, and
\begin{align*}
[\ra\rb\rc\rd\rf\rK\rL]=-[\ra\rb\rc\rd\rg\rI\rL]\quad & \Longrightarrow\quad f_3(d_3-d_4)=d_3e_3-d_4e_6 \\
& \Longrightarrow\quad f_3=\frac{1}{d_3(d_3-d_4)}(d_3^2e_3-d_4d_3^2+d_4d_6^2-d_4d_6e_3); \\
[\ra\rb\rc\rd\rf\rI\rJ]=[\ra\rb\rc\rd\rg\rJ\rK]\quad & \Longrightarrow\quad f_5(d_3-d_6)=d_3e_5-d_6e_4 \\
& \Longrightarrow\quad f_5=\frac{1}{d_3(d_3-d_6)}(d_3^2e_3+d_3^2d_4-2d_3^2d_6+d_4d_6^2-d_4d_6e_3).
\end{align*}

From our knowledge of the symmetric matroid realizations, we derive that
$$0=[\ra\rb\rc\rd\re\rg\rh]=(d_4-d_6)^2(2d_3+d_4+d_6)(2d_3-d_4-d_6).$$
Since~$d_4\ne d_6$, this implies that~$d_4=\pm 2d_3-d_6$. Assume first that~$d_4=-2d_3-d_6$. Then
$$0=[\ra\rb\rc\rd\re\rf\rg]=\frac{32(d_3+d_6)^4(d_3+d_6-e_3)^2}{(3d_3+d_6)(d_3-d_6)},$$
and we obtain~$d_6=-d_3+e_3$ or~$d_6=-d_3$. But this is impossible since~$d_3=[\ra\rb\rc\rd\rJ\rK\rL]>0$, $e_3=[\ra\rb\rc\re\rJ\rK\rL]<0$ and~$d_6=-[\ra\rb\rc\rd\rI\rJ\rK]>0$. Consequently,~$d_4=2d_3-d_6$.

\svs
Similarly, the knowledge of the chirotope ensures that
$$0=[\ra\rb\rc\re\rf\rg\rI]=\frac{4}{d_3^2}(d_3-d_6)^2(d_6-e_3-d_3)(d_6-(1+\sqrt{2})d_3-e_3)(d_6-(1-\sqrt{2})d_3-e_3).$$
Since~$d_3\ne d_6$, this implies that~$d_6=d_3+e_3$ or~$d_6=(1\pm \sqrt{2})d_3+e_3$. The solution ~$d_6=d_3+e_3$ is eliminated since $0\ne[\ra\rb\rc\rd\re\rg\rJ]=(d3-d6)^2(e3+d3-d6)$. The solution $d_6=(1-\sqrt{2})d_3+e_3$ is eliminated since~$d_3>0$, $e_3<0$ and $d_6>0$. Consequently, we are sure that $d_6=(1+\sqrt{2})d_3+e_3$.

\svs
Furthermore, in each column of the matrix~$M'$, the coordinates sum to~$1$ (since we have homogeneous coordinates when multipling by the matrix~$N$). This implies that~$d_3=1+\sqrt{2}/2$ and the matrix~$B$ can be written as stated in the lemma. Finally, the bounds on the possible value of~$u$ are derived from the signs of the chirotope: we have~${|\ra\rb\rc\rd\rf\rI\rK|=1}$, ${|\ra\rb\rc\rd\rf\rI\rJ|=1}$ and~${|\ra\rb\re\rf\rI\rJ\rL|=-1}$, which implies that~${u+1+\sqrt{2}<0}$, ${(u+1+\sqrt{2})(u+1+2\sqrt{2})<0}$ and ${2u+2+3\sqrt{2}<0}$ respectively. Thus, we have~$-1-2\sqrt{2}<u<-1-\frac{3\sqrt{2}}{2}$.
\end{proof}

These lemmas imply that only three of the~$15$ symmetric matroid realizations of~$\Delta_{8,2}$ are realizable geometrically by a polytope with symmetric determinants: if~${u<-2-\sqrt{2}}$ we obtain the chirotope~$G$, if~${u=-2-\sqrt{2}}$ we obtain~$K$, and if~${-2-\sqrt{2}<u}$ we obtain~$O$.

In order to complete our understanding of the space of symmetric realizations of~$\Delta_{8,2}$, it only remains to study the possible values of the matrix~$N$. To determine~$N$, we again use symmetry, but this time on the length of the edges of~$P$. For example, we know that the vertices $\rI$, $\rJ$, $\rK$, and~$\rL$ span a \simp{3} with~$|\rI\rJ|=|\rJ\rK|=|\rK\rL|=|\rI\rL|$ and~$|\rI\rK|=|\rJ\rL|$ (where~$|\rx\ry|$ denotes the euclidean distance from~$\rx$ to~$\ry$). Thus, since a dilation does not destruct the symmetry, we can assume that the matrix~$N$ is of the form:
$$N \eqdef \begin{pmatrix}
1 & 1 & 1 & 1 & 1 & 1 & 1 \\
x_1 & 0 & 0 & 0 & 0 & 0 & 0 \\
x_2 & y_2 & 0 & 0 & 0 & 0 & 0 \\
x_3 & y_3 & z_3 & 0 & 0 & 0 & 0 \\
x_4 & y_4 & z_4 & 1 & -1 & 1 & -1 \\
x_5 & y_5 & z_5 & v & 0 & -v & 0 \\
x_6 & y_6 & z_6 & 0 & v & 0 & -v
\end{pmatrix},$$
with~$x_1> 0$,~$y_2> 0$,~$z_3> 0$, and~$v>0$.

\mvs
Using the remaining equations given by the symmetries, we obtain (with the help of a computer algebra system) the following constraints:

\begin{lemma}
The coefficients of the matrix $N$ satisfy:
$$x_1=\sqrt{\frac{(2+\sqrt{2})(2y_3+z_3\sqrt{2})(y_3+z_3)z_3}{y_3-z_3}},\qquad x_2=\frac{y_3(2+\sqrt{2})(y_3+z_3)}{\sqrt{z_3^2-y_3^2}},$$
$$y_2=\sqrt{z_3^2-y_3^2}, \qquad x_3=-(2+\sqrt{2})(y_3+z_3\frac{\sqrt{2}}{2}),\qquad x_4=y_4=z_4=0,$$
$$x_5=-x_6=y_5=y_6=-z_5=z_6=-\frac{1}{2}\sqrt{2}v(1+\sqrt{2}+u),$$
with~$-z_3<y_3<-\frac{z_3}{\sqrt{2}}$.\qed
\end{lemma}

Reciprocally, it is easy to check (again with a computer algebra system) that under these conditions, the convex hull of the column vectors of~$M=NM'$ is a symmetric realization of~$\Delta_{8,2}$: we only have to verify that for any admissible values of~$u$, $v$, $y_3$ and~$z_3$:
\begin{enumerate}[(i)]
\item The boundary complex of the convex hull of the column vectors of~$M$ is indeed isomorphic to~$\Delta_{6,1}$. This is encoded in the chirotope (and we only have to check the signs for one representative of each orbit under the dihedral group~$\D_8$).
\item The dihedral group~$\D_8$ indeed acts by isometry on the columns vectors of~$M$ (and we only have to check it on two generators of~$\D_8$).
\end{enumerate}
Thus, we obtain the main result of this section:

\begin{proposition}\label{ft:prop:delta82}
Up to dilation, the space of symmetric realizations of the simplicial complex~$\Delta_{8,2}$ has dimension~$4$.\qed
\end{proposition}

\subsubsection{A general construction?}\label{ft:subsubsec:generalconstruction}

Even if our study of symmetric realizations of the simplicial complex~$\Delta_{8,2}$ cannot be extended to general~$n$ and~$k$, we consider that it provides new evidence and motivation for the general investigation. In particular, it seems reasonable to attempt to construct symmetric realizations of the multiassociahedron, in much the same way the secondary polytope of the \gon{n} is constructed: 
\begin{enumerate}[(i)]
\item First, associate to each \ktri{k}~$T$ of the \gon{n} a vector~$\phi(T)\in\R^{nk}$, with~$k$ coordinates for each vertex~$v$ of the \gon{n}. These coordinates would be computed according to the \kstar{k}s of~$T$ which contain~$v$.
\item Then, consider the convex hull of the vectors~$\phi(T)$ of all \ktri{k}s of the \gon{n}, which should be of codimension~$k(2k+1)$.
\end{enumerate}

\begin{figure}
	\capstart
	\centerline{\includegraphics[scale=1]{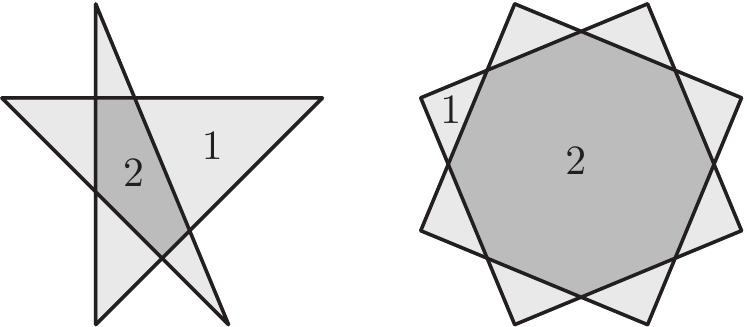}}
	\caption[The area of a \kstar{2} and the \kdepth{2} of the octagon]{The area of a \kstar{2} (each region is counted with multiplicity equal to the winding number of the \kstar{2} around it) and the \kdepth{2} of the octagon.}
	\label{ft:fig:areadepth}
\end{figure}

We don't have a full candidate for the vector~$\phi(T)$, but it seems at least natural to choose, for the first coordinate of~$\phi(T)$ associated to a vertex~$v$ of the \gon{n}, the sum of the areas of the \kstar{k}s of~$T$ containing~$v$. (Remember that the area of a \kstar{k}~$S$ of the \gon{n} is the integral over~$\R^2$ of the winding number of~$S$; that is to say, a point of the plane is counted with multiplicity equal to the number of revolutions of~$S$ around it). 
However, it turns out that the only relation between these coordinates comes from the total area covered by the \kstar{k}s of~$T$:
$$\cA_k(V_n)=\sum_{S\in T} \cA(S)=\frac{1}{2k+1}\sum_{v\in V_n} \sum_{v\in S\in T} \cA(S),$$
where~$\cA_k(V_n)$ is the integral over~$\R^2$ of the \kdepth{k} with respect to~$V_n$ (in other words, the sum of the winding numbers of the \kbound{k} components of the \gon{n}). See \fref{ft:fig:areadepth}. In particular, there is no more relations between the centers of mass of the \kstar{k}s of a \ktri{k}. For example, the $84$~area vectors of~$\R^8$ corresponding to the \ktri{2}s of the octagon (the $v$th coordinate of the vector associated to~$T$ is the sum of the areas of the \kstar{2}s of~$T$ containing~$v$) span an affine space of dimension~$7$ which is already to big to be completed into a \poly{6}tope.

\subsection{A generalization of Loday's associahedron}\label{ft:subsec:multiassociahedron:loday}

In this section, we present another attempt to construct a multiassociahedron. We reinterpret \mbox{Loday's} construction~\ref{ft:def:loday} of the associahedron in terms of beam arrangements. This interpretation has two main advantages: on the one hand, the construction becomes much more visual, which leads to a simple proof of Theorem~\ref{ft:theo:loday}; on the other hand our definition naturally extends to multitriangulations and gives rise to a polytope with nice combinatorial properties. We start with a short reminder on the properties of the vector configuration formed by the columns of the incidence matrix of a directed graph.

\subsubsection{The incidence cone of a directed graph}

Let~$G$ be a directed (multi)graph on~$p$ vertices, whose underlying undirected graph is connected. An oriented edge with origin~$i$ and endpoint~$j$ is denoted~$(i,j)$. Let~$(e_1,\dots,e_p)$ be the canonical basis of~$\R^p$.

\begin{definition}
The \defn{incidence configuration} of the directed (multi)graph~$G$ is the vector configuration~$I(G) \eqdef \ens{e_i-e_j}{(i,j)\in G}$. The \defn{incidence cone} of~$G$ is the cone~$C(G)$ generated by~$I(G)$, \ie its positive span.
\end{definition}

In other words, the incidence configuration of a graph consists of the column vectors of its incidence matrix. We will use the following relations between the graph properties of~$G$ and the orientation properties of~$I(G)$. We refer to~\cite[Section~1.1]{bvswz-om-99} for details.

\begin{observation}
Consider a subgraph~$H$ of~$G$. Then the vectors of~$I(H)$:
\begin{enumerate}[(i)]
\item are independent if and only if $H$~has no (non-necessarily oriented) cycle;
\item form a basis of the hyperplane~$\sum x_i = 0$ if and only if $H$~is a spanning tree;
\item form a circuit if and only if $H$~is a (non-necessarily oriented) cycle; the positive and negative parts of the circuit correspond to the subsets of edges oriented in one or the other direction along the cycle; in particular,~$I(H)$ is a positive circuit if and only if $H$~is an oriented cycle;
\item form a cocircuit if and only if $H$~is a minimal (non-necessarily oriented) cut; the positive and negative parts of the cocircuit correspond to the edges in one or the other direction in this cut; in particular,~$I(H)$ is a positive cocircuit if and only if $H$~is an oriented cut.
\end{enumerate}
\end{observation}


\begin{observation}\label{ft:obs:incidencecone}
\begin{enumerate}[(i)]
\item The incidence cone~$C(G)$ is pointed if and only if~$G$ is an acyclic directed graph. That is, if it induces a partial order on its set of nodes.
\item In this case, the rays of~$C(G)$ correspond to the edges of the Hasse diagram of~$G$. In particular, the cone is simplicial if and only if the Hasse diagram of~$G$ is a tree.
\item The facets of~$C(G)$ correspond to the complement of the minimal directed cuts of~$G$.
\end{enumerate}
\end{observation}

With this tool in mind, we can tackle serenely the topic of this section.

\subsubsection{Loday's construction revisited}\label{ft:subsubsec:lodaybis}

\begin{figure}[b]
	\capstart
	\centerline{\includegraphics[scale=1]{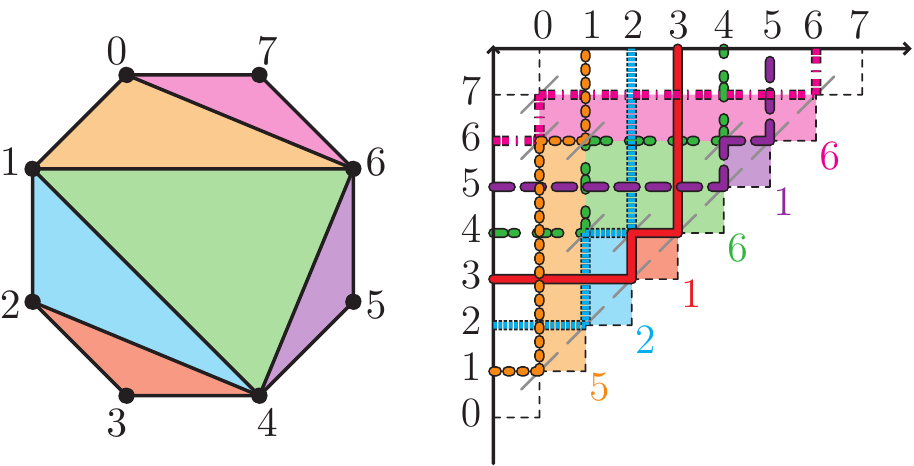}}
	\caption[Loday's \dfsvector{} of a triangulation interpreted on its beam arrangement]{Loday's \dfsvector{} of a triangulation interpreted on its beam arrangement.}
	\label{ft:fig:mirrorslaserstriangulation}
\end{figure}

Let~$T$ be a triangulation of the \gon{n}. Consider its beam arrangement (\ie its dual pseudoline arrangement drawn on the integer grid)~---~remember if necessary Section~\ref{ft:subsec:Dyckpaths:laser} and see \fref{ft:fig:mirrorslaserstriangulation}. As stated in Proposition~\ref{ft:prop:laser}, the $i$th beam~$B_i$ is an $x$- and $y$-monotone lattice path, coming from~$(-\infty,i)$ and going to~$(i,+\infty)$, which reflects in exactly three mirrors. For such a beam~$B_i$, we define its \defn{area}~$\cA_i$ to be the area inside the box~$[0,i]\times[i,n-1]$ and below~$B_i$.

\begin{lemma}
Loday's \dfsvector{}~$\psi(T)$ coincides with the \defn{beam vector}~$(\cA_1,\dots,\cA_{n-2})$.
\end{lemma}

\begin{proof}
Let~$T$ be a triangulation of the \gon{n} and~$\tau$ denote its dual binary tree on~$n-2$ nodes (rooted at the triangle of~$T$ containing the edge~$[0,n-1]$). The node~$x_i$ labeled by~$i$ in the \dfslabeling{} of~$\tau$ corresponds to the triangle whose second vertex is $i$, that is, to the $i$th beam. Furthermore, if this triangle has vertices~$h<i<j$, then the number of leaves in the left (resp.~right) child of~$x_i$ is exactly~$i-h$ (resp.~$j-i$). Thus, the $i$th coordinate of~$\psi(T)$, which is by definition the product of the numbers of leaves in both children of~$x_i$, equals the area $\cA_i=(i-h)(j-i)$ of the square located between the beam~$B_i$ and the lines~$x=i$ and~$y=i$.
\end{proof}

Consequently, we see the \dfspolytope{}~$\Omega_{n-2}$ as the convex hull of the beam vectors of all triangulations of the \gon{n}. With this interpretation, we present a simple proof of Loday's Theorem~\ref{ft:theo:loday}, based on the behavior of the beam vectors with respect to flips:

\begin{lemma}
Let~$T$ and~$T'$ be two triangulations of the \gon{n} related by a flip which exchanges the diagonals in the quadrangle formed by gluing their $i$th and $j$th triangles. Then the difference~$\psi(T)-\psi(T')$ of their beam vectors is parallel to~$e_i-e_j$.
\end{lemma}

\begin{proof}
Only the beams~$B_i$ and~$B_j$ are perturbed, and the area lost by~$B_i$ is transferred to~$B_j$.
\end{proof}

\begin{corollary}\label{ft:coro:lodayrelation}
$\Omega_{n-2}$~is contained in the affine hyperplane $\ens{x\in\R^{n-2}}{\dotprod{\one}{x}={n-1 \choose 2}}$.
\end{corollary}

\begin{proof}
The beam vector of the minimal triangulation~$T_{n,1}^{\min}$ belongs to this hyperplane. Since the difference between the beam vectors of two triangulations related by a flip is orthogonal to the vector~$\one$, the result follows from the connectedness of the flip graph.
\end{proof}

Let~$\Box_{u,v}$ denote the unit grid square~$\{(u,v),(u-1,v),(u,v+1),(u-1,v+1)\}$ whose bottom right vertex is~$(u,v)$. Another way to prove the previous lemma is to observe that, for any triangulation~$T$ of the \gon{n} and any~$1\le u\le v\le n-2$, the square~$\Box_{u,v}$ is covered exactly once by the beam arrangement of~$T$, meaning that there exists a unique~$1\le i\le n-2$ such that~$\Box_{u,v}$ lies between the beam~$B_i$ and the lines~$x=i$ and~$y=i$. We refer to \fref{ft:fig:mirrorslaserstriangulation} for an illustration of this property.


\begin{proposition}\label{ft:prop:lodayfacets}
For any internal edge~$[u,v]$ of the \gon{n}, the equation ${\sum_{p=u+1}^{v-1} x_p = {v-u \choose 2}}$ defines a facet~$F_{u,v}$ of~$\Omega_{n-2}$. This facet contains precisely the beam vectors of all triangulations of the \gon{n} containing the edge~$[u,v]$.
\end{proposition}

\begin{proof}
Any square~$\Box_{x,y}$ is covered by a beam~$B_i$ which satisfies~$x\le i \le y$. Consequently, the beams~$\ens{B_i}{u<i<v}$ cover at least the zone~$Z_{u,v}$ located to the right of the line~$x=u$ and below the line~$y=v$. This zone has area~${v-u \choose 2}$ which proves that the inequality is valid.

\enlargethispage{.3cm}
Consider now a triangulation~$T$ of the \gon{n} containing the edge~$[u,v]$. For any~$i,j$ such that~$u<i<v<j$, the edge~$[i,j]$ cannot be in~$T$. Consequently,~$B_i$ is straight from~$(i,v)$ to~$(i,+\infty)$. Similarly,~$B_i$ is straight from~$(-\infty,i)$ to~$(u,i)$. Thus, the beams $\ens{B_i}{u<i<v}$ are the only beams passing in the region~$Z_{u,v}$, and all their crossing points are located in~$Z_{u,v}$. Consequently, we can interpret them as a beam arrangement of the \gon{(v-u+1)}, which implies (by the previous lemma) that $\sum_{p=u+1}^{v-1} x_p = {v-u \choose 2}$.

Reciprocally, if the edge~$[u,v]$ is not in~$T$, then by maximality there is another edge~$[x,y]$ preventing it to be in~$T$. Assuming that~$0\le x<y\le n-1$, we have either~$u<x<v<y$, or~$x<u<y<v$. In the first (resp.~second) case, the beam~$B_x$ (resp.~$B_y$) cannot be straight from~$(x,v)$ to~$(x,+\infty)$ (resp.~from~$(-\infty,y)$ to~$(u,y)$), and thus,~$B_x$ (resp.~$B_y$) covers at least one square outside~$Z_{u,v}$. Consequently, the beams~$\ens{B_i}{u<i<v}$ cover at least the zone~$Z_{u,v}$ plus one square, which proves that the inequality is strict if~$[u,v]$ is not is~$T$.

Finally, we shall prove that~$F_{u,v}$ is indeed a facet, and not a face of smaller dimension. Given a triangulation~$T$ containing~$[u,v]$, we stay in~$F_{u,v}$ when we flip any edge different from~$[u,v]$. This gives rise to $n-4$~linearly independent directions in~$F_{u,v}$. According to Corollary~\ref{ft:coro:lodayrelation}, this ensures that~$F_{u,v}$ is a facet and that the polytope~$\Omega_{n-2}$ has dimension~$n-3$.
\end{proof}

\begin{proposition}\label{ft:prop:coneLoday}
Let~$T$ be a triangulation of the \gon{n}, and~$\tau$ denote its dual binary tree, where the edges are oriented towards the root. Then the cone~$C(T)$ of~$\Omega_{n-2}$ at~$\psi(T)$ equals the incidence~cone~$C(\tau)$.
\end{proposition}

\begin{proof}
The cone~$C(T)$ certainly contains the incidence cone~$C(\tau)$: for any~$(i,j)\in\tau$, the ray $\R^+(e_i-e_j)$ of~$C(\tau)$ is obtained in~$\Omega_{n-2}$ by flipping the edge common to the $i$th and $j$th triangles of~$T$. Reciprocally, the facet of~$C(\tau)$ obtained as the convex hull of all rays except that of the edge~$(i,j)\in\tau$ corresponds to the facet~$F_{i,j}$ of~$\Omega_{n-2}$.
\end{proof}

\begin{corollary}
Every facet of~$\Omega_{n-2}$ is of the form~$F_{u,v}$ for some relevant edge~$[u,v]$ of the \gon{n}.
\end{corollary}

\begin{proof}
For any triangulation~$T$, all the facets of the cone at~$\psi(T)$ are of~this~form.
\end{proof}

These results prove that~$\Omega_{n-2}^\polar$ realizes~$\Delta_{n,1}$: the boundary complex of~$\Omega_{n-2}^\polar$ and the simplicial complex~$\Delta_{n,1}$ have exactly the same vertex-facet incidences, and these incidences characterize the complete complex. To summarize our results, we have shown that:
\begin{enumerate}
\item All beam vectors lie in a hyperplane orthogonal to~$\one$, and~$\Omega_{n-2}$ has dimension~$n-3$.
\item For each relevant edge~$[u,v]$ of the \gon{n}, the vector $\sum_{p=u+1}^{v-1} e_i$ is the normal vector of a facet of~$\Omega_{n-2}$, whose vertices are precisely the beam vectors of all triangulations containing~$[u,v]$. Furthermore all facets of~$\Omega_{n-2}$ are of this form.
\item For each triangulation~$T$, the cone of~$\Omega_{n-2}$ at the beam vector~$\psi(T)$ equals the incidence cone of the dual tree of~$T$.
\item $\Omega_{n-2}^\polar$~realizes the simplicial complex~$\Delta_{n,1}$.
\end{enumerate}

\subsubsection{The beam polytope}\label{ft:subsubsec:beampoly}

\begin{figure}
	\capstart
	\centerline{\includegraphics[scale=1]{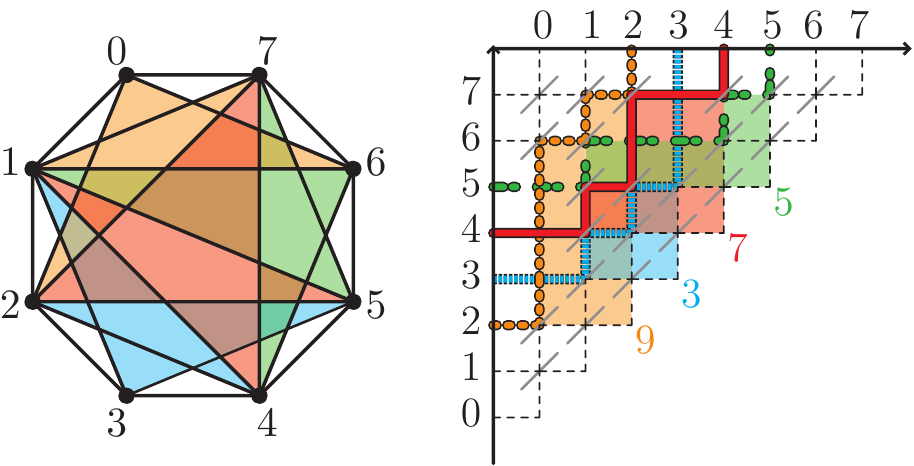}}
	\caption[The beam vector of the \ktri{2} of~\fref{intro:fig:2triang8points}]{The beam vector of the \ktri{2} of~\fref{intro:fig:2triang8points}.}
	\label{ft:fig:mirrorslasersarea}
\end{figure}

The main advantage of our interpretation of Loday's construction via the beam arrangement is its natural extension to multitriangulations. Remember that a \ktri{k} of the \gon{n} can also be seen by duality as a beam arrangement (\ie a pseudoline arrangement drawn on the integer grid)~---~see~\fref{ft:fig:mirrorslasersarea}. The $i$th beam~$B_i$ is an $x$- and $y$-monotone lattice path, coming from~$(-\infty,k-1+i)$ and going to~$(k-1+i,+\infty)$, which reflects in exactly~$2k+1$ mirrors.

\begin{definition}
\index{beam!--- vector}
\index{beam!--- polytope}
The \defn{beam vector} of a \ktri{k}~$T$ of the \gon{n} is the vector~$\psi(T)\in\R^{n-2k}$ whose $i$th coordinate is the area~$\cA_i$ of the region located below the beam~$B_i$ and inside the box~$[0,k-1+i]\times[k-1+i,n-1]$. The~\defn{beam polytope}~$\Omega_{n,k}$ is the convex hull of the beam vectors of all \ktri{k}s of the \gon{n}.
\end{definition}

\begin{remark}
We could have defined directly the beam vector~$\psi(T)$ of a \ktri{k}~$T$ with a formula: its $i$th coordinate is given by
$$\psi(T)_i \eqdef \sum_{j=1}^k(s_{j+1}-s_j)(s_{k+1+j}-s_{k+1}),$$
where $0\cle s_1\cl\dots\cl s_{2k+1}\cle n-1$ are the vertices of the \kstar{k} of~$T$ whose $(k+1)$th vertex~$s_{k+1}$ is the vertex~$k-1+i$.
\end{remark}

As for triangulations, the interesting property of the beam vectors is their behavior with respect to flips:

\begin{lemma}
Let~$T$ and~$T'$ be two \ktri{k}s of the \gon{n} related by a flip which involves their $i$th and $j$th \kstar{k}s. Then the difference~$\psi(T)-\psi(T')$ of their beam vectors is parallel to~$e_i-e_j$.
\end{lemma}

\begin{proof}
Only the beams~$B_i$ and~$B_j$ are perturbed, and the area lost by~$B_i$ is transferred to~$B_j$.
\end{proof}

\begin{corollary}
The polytope~$\Omega_{n,k}$ is contained in the affine hyperplane defined by
$$\ens{x\in\R^{n-2k}}{\dotprod{\one}{x}=\frac{k(n-2k)(n-k)}{2}}.$$
\end{corollary}

\proof
The proof is similar to that of Corollary~\ref{ft:coro:lodayrelation}: the lemma is derived from the connectedness of the flip graph and from the calculation:

\vspace*{.6cm}\qed
\vspace*{-1cm}
$$\dotprod{\one}{\psi(T_{n,k}^{\min})}=\sum_{i=1}^{n-2k} \left(\frac{k(k+1)}{2}+k(n-2k-i)\right)=\frac{k(n-2k)(n-k)}{2}.$$

\begin{example}
The beam vectors of the \ktri{2}s of the heptagon are:
\[
\begin{bmatrix} 3 \\ 5 \\ 7 \end{bmatrix},
\begin{bmatrix} 7 \\ 3 \\ 5 \end{bmatrix},
\begin{bmatrix} 5 \\ 7 \\ 3 \end{bmatrix},
\begin{bmatrix} 3 \\ 7 \\ 5 \end{bmatrix},
\begin{bmatrix} 5 \\ 3 \\ 7 \end{bmatrix},
\begin{bmatrix} 7 \\ 5 \\ 3 \end{bmatrix},
\begin{bmatrix} 5 \\ 5 \\ 5 \end{bmatrix},
\begin{bmatrix} 5 \\ 4 \\ 6 \end{bmatrix},
\begin{bmatrix} 3 \\ 8 \\ 4 \end{bmatrix},
\begin{bmatrix} 6 \\ 3 \\ 6 \end{bmatrix},
\begin{bmatrix} 4 \\ 8 \\ 3 \end{bmatrix},
\begin{bmatrix} 6 \\ 4 \\ 5 \end{bmatrix},
\begin{bmatrix} 4 \\ 4 \\ 7 \end{bmatrix},
\begin{bmatrix} 7 \\ 4 \\ 4 \end{bmatrix}.
\]
The sum of their coordinates always equals~$15$. Projecting on the hyperplane~$\dotprod{(1,1,1)}{x}=15$, we obtain the point set in~\fref{ft:fig:loday72}. Observe that we have chosen the projection in such a way that the symmetry of this point set (obtained by reversing the coordinates of the beam vectors in~$\R^3$) is the orthogonal reflection with respect to the horizontal axis. For convenience, we have labeled the \krel{2} edges by letters ($\ra \eqdef [0,3]$, $\rb \eqdef [1,4]$, \etc), and the \ktri{2}s by their set of \krel{2} edges. For example,~$\ra\rb\re\rf$ is the \ktri{2} whose set of \krel{2} edges is~$\{[0,3],[1,4],[0,4],[1,5]\}$.

\begin{figure}
	\capstart
	\centerline{\includegraphics[scale=1]{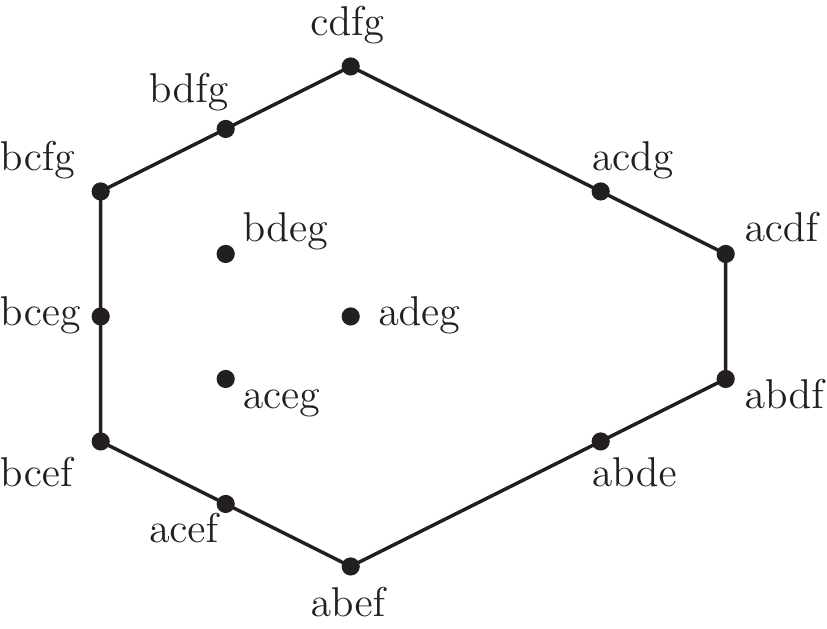}}
	\caption[The \poly{2}tope $\Omega_{7,2}$]{The \poly{2}tope $\Omega_{7,2}$.}
	\label{ft:fig:loday72}
\end{figure}

The $f$-vectors of the polytopes~$\Omega_{8,2}$, $\Omega_{9,2}$~and~$\Omega_{10,2}$ are~$(22,33,13)$, $(92,185,118,25)$~and $(420,1062,945,346,45)$ respectively. We have represented~$\Omega_{8,2}$ and~$\Omega_{9,2}$ in Figures~\ref{ft:fig:loday82} and~\ref{ft:fig:loday92}. The polytope~$\Omega_{8,2}$ is simple while the polytope~$\Omega_{9,2}$ has two non-simple vertices (which are contained in the projection facet of the Schlegel diagram on the right of \fref{ft:fig:loday92}) and the polytope~$\Omega_{10,2}$ has~$24$ non-simple vertices. We will characterize later non-simple vertices in the beam polytope.

\begin{figure}
	\capstart
	\centerline{\includegraphics[scale=1.15]{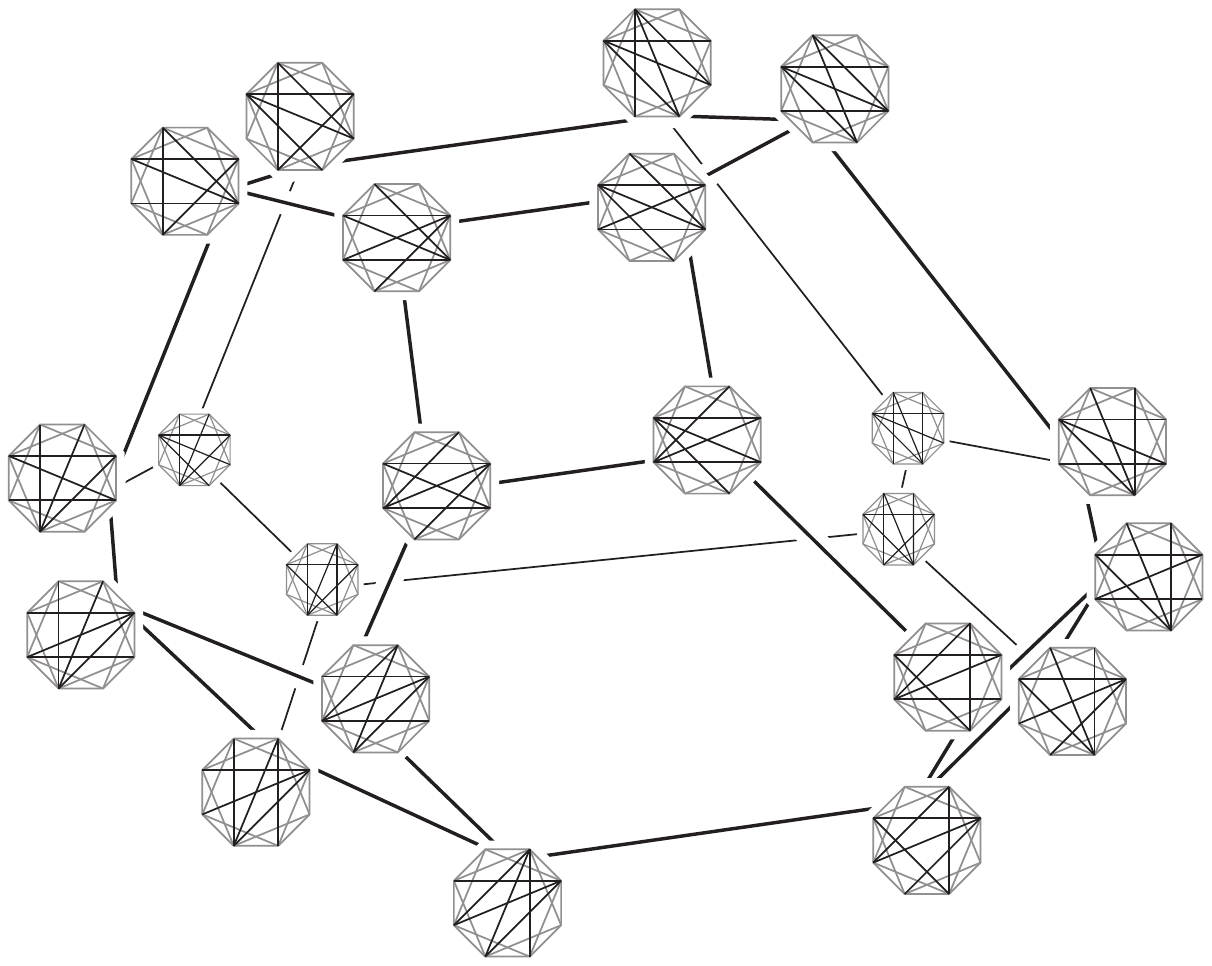}}
	\caption[The \poly{3}tope $\Omega_{8,2}$]{The \poly{3}tope $\Omega_{8,2}$.}
	\label{ft:fig:loday82}
\end{figure}

\begin{figure}
	\capstart
	\centerline{\includegraphics[scale=.64]{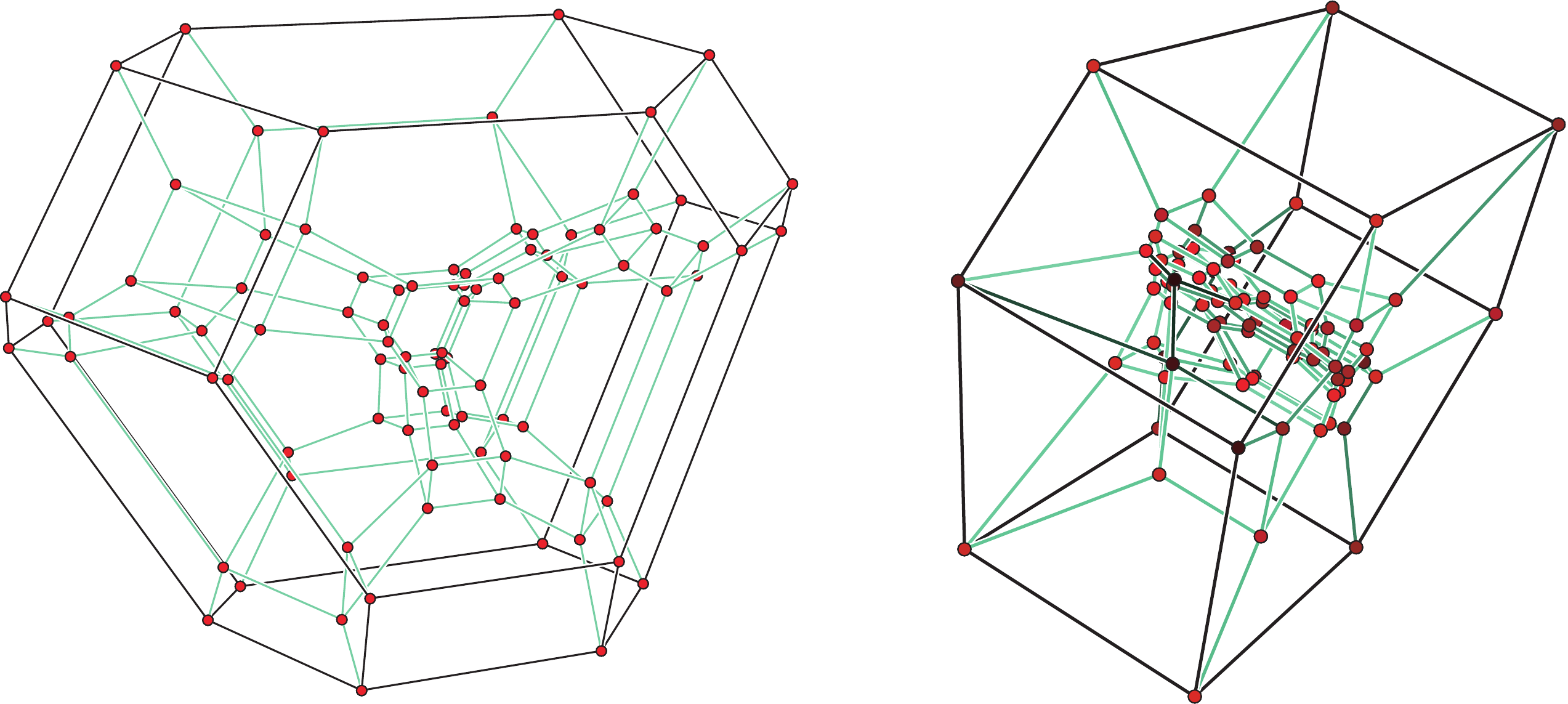}}
	\caption[The \poly{4}tope $\Omega_{9,2}$]{Two Schlegel diagrams of the \poly{4}tope $\Omega_{9,2}$. On the second one, the two leftmost vertices of the projection facet are non-simple vertices.}
	\label{ft:fig:loday92}
\end{figure}

\end{example}

The dimension of the beam polytope~$\Omega_{n,k}$ is~$n-2k-1$, and thus differs from the dimension~$k(n-2k-1)$ of the $k$-associahedron we would like to construct. However, it is reasonable to think that the beam polytope~$\Omega_{n,k}$ could be a projection of a $k$-associahedron. In this case, we could use the beam polytope~$\Omega_{n,k}$ as a starting point to construct this $k$-associahedron, splitting each of the coordinates of~$\Omega_{n,k}$ into~$k$ coordinates. With this motivation, we study in the next section the facets of the beam polytope and check that they are compatible with~$\Delta_{n,k}$: we prove that the \ktri{k}s whose beam vector belongs to a given facet of the beam polytope form a face of~$\Delta_{n,k}$.


\subsubsection{\kvalid{k} \zoseq{}s and facet description}\label{ft:subsubsec:facetdescription}

To emphasize the rich combinatorial structure of our beam polytope~$\Omega_{n,k}$, we give a simple description of its facet-defining equations and of its vertex-facet incidences, which generalizes Proposition~\ref{ft:prop:lodayfacets}. The main characters of this description are the following sequences:

\begin{definition}
\index{valid sequence@\kvalid{k} sequence}
We call \defn{\kvalid{p} \zoseq{}} of size~$q$ any sequence of~$\{0,1\}^q$ which is neither~$(0,0,\dots,0)$, nor~$(1,1,\dots,1)$ and does not contain a subsequence~$10^r1$ for~$r\ge p$. In other words, all subsequences of~$p$ consecutive zeros appear before the first one or after the last one.
\end{definition}

\begin{example}
\kvalid{1} \zoseq{}s are sequences whose only zeros are located at the beginning or at the end, \ie sequences of the form~$0^x1^y0^z$ with~$x+y+z=q$ and~$y\ne0\ne x+z$. Observe that these sequences are exactly the normal vectors of Loday's \dfspolytope{}~$\Omega_q$ (see Proposition~\ref{ft:prop:lodayfacets}).
\end{example}

\begin{remark}
Let~$\VS_{p,q}$ denote the number of \kvalid{p} \zoseq{}s of size~$q$. It is easy to see that~$\VS_{1,q}=\frac{1}{2}(q-1)(q+2)$ (the number of internal diagonals of a \gon{(q+2)}!) and that $\VS_{2,q}=F_{q+4}-(q+4)$, where~$F_n$ denotes the $n$th Fibonacci number (indeed, if $u_q$ denotes the number of \zoseq{}s of length~$q$ which start with a~$1$ and do not contain a subsequence~$00$ before their last~$1$, then $u_q=u_{q-1}+u_{q-2}+1$ and ~$\VS_{2,q}=\VS_{2,q-1}+u_q$). To compute~$\VS_{p,q}$ in general, consider the non-ambiguous rational expression~$0^*(1,10,100,\dots,10^{p-1})^*10^*$. The corresponding rational language consists of all \kvalid{p} \zoseq{}s plus all non-empty sequences of~$1$. Thus, the generating function of \kvalid{p} \zoseq{}s is:
$$\sum_{q\in\N} \VS_{p,q}x^q=\frac{1}{1-x}\,\frac{1}{1-\sum_{i=1}^{p} x^i}\,x\,\frac{1}{1-x}-\frac{x}{1-x}=\frac{x^2(2-x^p)}{(1-2x+x^{p+1})(1-x)}.$$
\end{remark}

We use these sequences to describe the facets of the beam polytope. For this, we need one more definition.
Let~$\sigma$ be a \zoseq{} of length~$n-2k$. Let~$|\sigma|_0$ denote the number of zeros in~$\sigma$. For all~$i\le |\sigma|_0+2k$, we denote by~$\zeta_i(\sigma)$ the position of the $i$th zero in the \zoseq{}~$0^k\sigma0^k$, obtained from~$\sigma$ by appending a prefix and a suffix of $k$ consecutive zeros. We associate to~$\sigma$ the set of edges:
$$D(\sigma) \eqdef \ens{[\zeta_i(\sigma),\zeta_{i+k}(\sigma)]}{i\in[|\sigma|_0+k]}.$$
Observe that~$D(\sigma)$ can contain some \kbound{k} edges of the \gon{n}.

\begin{theorem}\label{ft:theo:validsequencesfacetdescription}
Each \kvalid{k} \zoseq{}~$\sigma$ of size~$n-2k$ is the normal vector of a facet~$F_\sigma$ of the beam polytope $\Omega_{n,k}$, which contains precisely the beam vectors of the \ktri{k}s of the \gon{n} containing~$D(\sigma)$. 
\end{theorem}

\begin{proof}
We denote by~$Z(\sigma) \eqdef \ens{(x,y)\in\R^2}{\exists\, [u,v]\in D(\sigma), x\ge u \text{ and } y\le v}$ the zone of all points of the plane which see~$D(\sigma)$ in their north-west quadran.

We first prove that the value of~$\dotprod{\sigma}{\psi(T)}$ is the same for all \ktri{k}s of the \gon{n} which contain~$D(\sigma)$. Let~$T$ be such a \ktri{k}. For any integer~$i$ with~$\sigma_i=1$, any mirror located out of the zone~$Z(\sigma)$ and on one of the lines~$x=i$ or~$y=i$ creates a \kcross{(k+1)} with~$D(\sigma)$. Thus, the beam~$B_i$ arrives straight from~$(-\infty,i+k-1)$ until it enters~$Z(\sigma)$, and goes straight to~$(i+k-1,+\infty)$ as soon as it goes out of~$Z(\sigma)$. In other words, the beams~$B_i$ with~$\sigma_i=1$ are the only beams passing through a vertex in~$Z(\sigma)$ not in~$D(\sigma)$. Consequently, the only contact points between a beam~$B_i$ with~$\sigma_i=1$ and a beam~$B_j$ with $\sigma_j=0$ are those of~$D(\sigma)$. Thus, the flip of any edge not in~$D(\sigma)$ does not change the value of~$\dotprod{\sigma}{\psi(T)}$. This value is constant on all \ktri{k}s of the \gon{n} which contain~$D(\sigma)$ by connectedness of the flip graph.

Consider now a \ktri{k}~$T$ of the \gon{n} which does not contain the set~$D(\sigma)$. For any $u\le i+k-1\le v$, the beam~$B_i$ covers the square~$\Box_{u,v}$ except if it passes in a level between the points~$(u,v)$ and $(i+k-1,i+k-1)$. Consequently, each square of~$Z(\sigma)$ is covered by the beams~$(B_i)_{\sigma_i=1}$ of~$T$ at least as many times as it is covered by those of any \ktri{k} containing~$D(\sigma)$. Furthermore, if~$[u,v]$ is an edge of~$D(\sigma)$ which is not in~$T$, then at least one beam~$B_i$ of~$T$ with~$u\le i+k-1\le v$ and~$\sigma_i=1$ covers a square out of~$Z(\sigma)$. This ensures that~$\dotprod{\sigma}{.}$ is strictly bigger on~$\psi(T)$ than on the beam vector of a \ktri{k} containing~$D(\sigma)$.

We still have to prove that~$F_\sigma$ is indeed a facet, and not a face of smaller dimension. This is a direct consequence of the study of the vertex cones of~$\Omega_{n,k}$ to come next.
\end{proof}

In the next section, we prove furthermore that all the facets of the beam polytope are derived from \kvalid{k} \zoseq{}s as those of Theorem~\ref{ft:theo:validsequencesfacetdescription}.

\subsubsection{Directed dual multigraphs, cones, and vertex characterization}


We consider the dual multigraph~$T^\dual$ of~$T$. Remember (see the definition in Section~\ref{stars:subsec:surfaces:dual}) that~$T^\dual$ has one vertex for each \kstar{k} of~$T$ and one edge between its two vertices corresponding to the \kstar{k}s~$R$ and~$S$ of~$T$ for each common \krel{k} edge of~$R$ and~$S$. We label by~$i$ the vertex of~$T^\dual$ which corresponds to the $i$th \kstar{k} of~$T$, that is, to the $i$th beam~$B_i$. For the edges, we need two notations: on the one hand, we denote by~$e^\dual$ the edge of~$T^\dual$ which is dual to the \krel{k} edge~$e$ of~$T$; on the other hand, we write~$(i,j)$ for an edge of~$T^\dual$ which relates its vertices~$i$ and~$j$ (we can have more than one edge~$(i,j)$). Observe that~$T^\#$ comes naturally with its edges oriented: each edge~$(i,j)$ of~$T^\dual$, corresponding to a mirror separating the two beams~$B_i$ and~$B_j$ of the beam arrangement of~$T$, is oriented from~$i$ to~$j$ if~$B_i$ lies above~$B_j$ at this mirror, and in the other direction otherwise.

We first need to understand the directed cuts of this oriented multigraph:

\begin{lemma}
Let~$T$ be a \ktri{k} of the \gon{n}. Consider a set of \krel{k} edges~$E\subset T$ and its dual set~$E^\dual \eqdef \ens{e^\dual}{e\in E}\subset T^\dual$. Assume that~$E^\dual$ is a directed cut of~$T^\dual$: there is a partition~$U\sqcup V=[n-2k]$ of the vertices of~$T^\dual$ such that $E^\dual=T^\dual\cap\ens{(u,v)}{u\in U, v\in V}$ and $T^\dual\cap\ens{(v,u)}{u\in U, v\in V}=\emptyset$. Let~$\sigma$ be the caracteristic \zoseq{} of~$V$, \ie the sequence of length~$n-2k$ with~$\sigma_i=1$ if~$i\in V$ and~$0$ otherwise. Then:
\begin{enumerate}[(i)]
\item $E$ is precisely the set of \krel{k} edges in $D(\sigma)$.
\item If the directed cut~$E^\dual$ is minimal, then~$\sigma$ is \kvalid{k}.
\end{enumerate}
\end{lemma}

\begin{proof}
We define again the zone~$Z(\sigma) \eqdef \ens{(x,y)\in\R^2}{\exists\, [u,v]\in D(\sigma), x\ge u \text{ and } y\le v}$ of all points of the plane which see~$D(\sigma)$ in their north-west quadran.

We claim that all the beams~$B_i$ with~$i\in V$ are straight out of~$Z(\sigma)$: otherwise, we would have in~$T^\dual$ an edge from~$V$ to~$U$. Consequently, the beams~$B_i$ with~$i\in V$ precisely cover the zone~$Z(\sigma)$ plus all the vertical an horizontal straight lines~$x=i$ and~$y=i$ for~$i\in V$. This implies that the contact points between a beam~$B_i$ with~$i\in V$ and a beam~$B_j$ with~$j\in U$ are precisely the \krel{k} edges in $D(\sigma)$.

For the second point, assume that~$\sigma$ is not \kvalid{k}. Then there exists a sequence of at least~$k$ zeros between two ones. Let~$V_1$ (resp.~$V_2$) denote the set of all positions of ones before (resp.~after) this sequence of zeros. Then the beams~$B_i$ with~$i\in V_1$ have no contact points with the beams~$B_j$ with~$j\in V_2$. Consequently, the sets~$V_1$ and~$U\cup V_2$ are separated by a directed cut of~$T^\dual$ which is strictly contained in~$E^\dual$. This proves that~$E^\dual$ was not minimal.
\end{proof}

\begin{theorem}\label{ft:theo:conebeampolytope}
Let~$T$ be a \ktri{k} of the \gon{n} and~$T^\dual$ denote its oriented dual multigraph. The cone~$C(T)$ of~$\Omega_{n,k}$ at~$\psi(T)$ equals the incidence cone~$C(T^\dual)$.
\end{theorem}

\begin{proof}
The proof is similar to that of Proposition~\ref{ft:prop:coneLoday}. Observe first that the cone~$C(T)$ certainly contains the incidence cone~$C(T^\dual)$: for any~$(i,j)\in T^\dual$, the ray~$\R^+(e_i-e_j)$ of~$C(T^\dual)$ is obtained in~$C(T)$ by flipping an edge common to the $i$th and $j$th stars of~$T$.

To prove the reciprocal inclusion, we show that each facet of~$C(T^\dual)$ is also a facet of~$C(T)$. Let~$F$ be a facet of~$C(T^\dual)$. According to Observation~\ref{ft:obs:incidencecone}(iii), there is a minimal directed cut~$E^\dual$ of~$T^\dual$ such that~$F$ is the cone over~$\ens{e_i-e_j}{(i,j)\in T^\dual\ssm E^\dual}$. The previous lemma then ensures that~$E^\dual$ is the dual of the set of \krel{k} edges in~$D(\sigma)$, where~$\sigma$ is a \kvalid{k} \zoseq{}. Consequently, the facet~$F$ of~$C(T^\dual)$ is also a facet~$F_\sigma$ of~$C(T)$.
\end{proof}

Observe that this proof also finishes the proof of Theorem~\ref{ft:theo:validsequencesfacetdescription}: any~$F_\sigma$ is indeed a facet of the cones at its vertices. It ensures moreover that:

\begin{corollary}
Any facet of~$\Omega_{n,k}$ is of the form~$F_\sigma$ for a \kvalid{k} \zoseq{}~$\sigma$ of length~$n-2k$.
\end{corollary}

\begin{proof}
All the facets of the cone at any vertex of~$\Omega_{n,k}$ are of this form.
\end{proof}

\begin{example}
For example, the numbers of facets of the beam polytopes~$\Omega_{7,2}$, $\Omega_{8,2}$, $\Omega_{9,2}$ and~$\Omega_{10,2}$ are respectively~$6$, $13$, $25$ and~$45$. They equal the numbers~$\VS_{2,n-4}=F_n-n$ of \kvalid{2} \zoseq{}s of length~$n-4$ (where~$F_n$ denotes the $n$th Fibonacci number).
\end{example}

Furthermore, Theorem~\ref{ft:theo:conebeampolytope} together with Observation~\ref{ft:obs:incidencecone}(i) provide a criterion to characterize the vertices of the beam polytope:

\begin{corollary}
Let~$T$ be a \ktri{k} of the \gon{n}.
\begin{enumerate}[(i)]
\item The beam vector~$\psi(T)$ is a vertex of the beam polytope~$\Omega_{n,k}$ if and only if the oriented dual multigraph~$T^\dual$ of~$T$ is acyclic.
\item When~$T^\dual$ is acyclic, it defines a partially ordered set on its vertices, and the vertex~$\psi(T)$ is simple if and only if the Hasse diagram of this poset is a tree.
\end{enumerate}
\end{corollary}

\begin{example}
The beam polytope~$\Omega_{9,2}$ has~$92$ vertices, two of which are not simple. For one of them, we have represented in \fref{ft:fig:hassediagram} the corresponding \ktri{2} of the \gon{9}, its dual multigraph, and its Hasse diagram.

\begin{figure}
	\capstart
	\centerline{\includegraphics[scale=1]{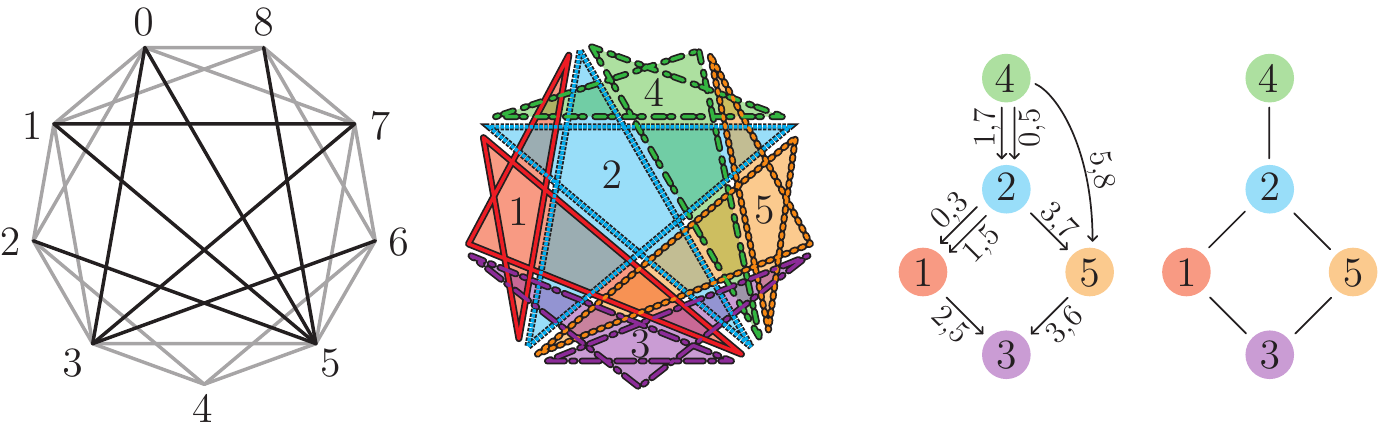}}
	\caption[A \ktri{2} of the \gon{9} whose oriented dual multigraph is acyclic, but whose Hasse diagram is not a tree]{A \ktri{2}~$T$ of the \gon{9} whose oriented dual multigraph~$T^\dual$ is acyclic, but whose Hasse diagram is not a tree. Its beam vector~$\psi(T)$ is a non-simple vertex of the beam polytope~$\Omega_{9,2}$. We have represented the \ktri{2}~$T$, its decomposition into \kstar{2}s, its oriented dual multigraph~$T^\dual$ and the associated Hasse diagram. We have labeled by~$i$ the node of~$T^\dual$ corresponding to the $i$th \kstar{2} of~$T$, and by~$u,v$ the dual edge of~$[u,v]\in T$.}
	\label{ft:fig:hassediagram}
\end{figure}
\end{example}

Let us summarize our combinatorial description of the beam polytope:
\begin{enumerate}
\item All beam vectors lie in a common hyperplane orthogonal to~$\one$, and the beam polytope has dimension~$n-2k-1$. 
\item Each \kvalid{k} \zoseq{}~$\sigma$ of length~$n-2k$ is the normal vector of a facet~$F_{\sigma}$ of~$\Omega_{n,k}$, and the set of \ktri{k}s whose beam vector is contained in~$F_{\sigma}$ is a face of the simplicial complex~$\Delta_{n,k}$. Furthermore, all facets of~$\Omega_{n,k}$ are of this form.
\item The cone of~$\Omega_{n,k}$ at the beam vector~$\psi(T)$ of a \ktri{k}~$T$ of the \gon{n} is the cone~$C(T^\dual)$ of the oriented dual multigraph~$T^\dual$ of~$T$.
\item The vertices of the beam polytope~$\Omega_{n,k}$ are the beam vectors of the \ktri{k}s whose oriented dual multigraph is acyclic. A vertex~$\psi(T)$ is simple if and only if the Hasse diagram of~$T^\dual$ is a tree.
\end{enumerate}

\subsubsection{Projection obstruction}\label{ft:subsubsec:projectionobstruction}

Despite these encouraging properties, it turns out that the beam polytope can in fact not be the projection of a realization of~$\Delta_{n,k}$. We finish this section with a proof that~$\Omega_{7,2}$ is not the projection of a realization of~$\Delta_{7,2}$.

\begin{figure}
	\capstart
	\centerline{\includegraphics{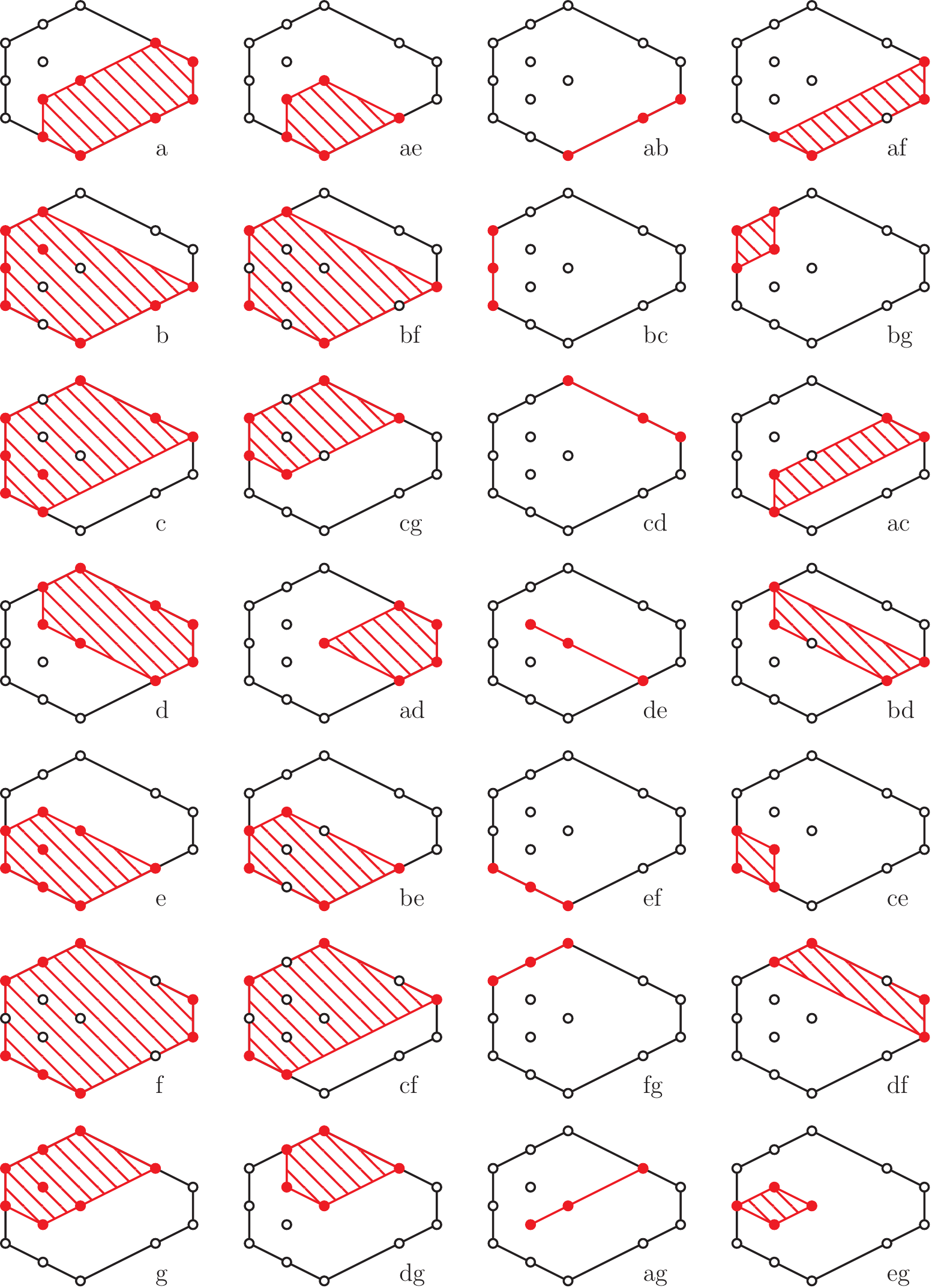}}
	\caption[Projections of the $1$- and $2$-codimensional faces of~$\Delta_{7,2}^\polar$]{Projections of the $1$- and $2$-codimensional faces of~$\Delta_{7,2}^\polar$.}
	\label{ft:fig:loday72facets}
\end{figure}

\begin{proposition}
Let~$P$ be a \poly{4}tope with~$14$ vertices labeled by the~$14$ \ktri{2}s of the heptagon such that the faces of~$P$ are faces of~$\Delta_{7,2}$. Then it is not possible to project~$P$ down to the plane such that, for each \ktri{2} of the heptagon, the vertex of~$P$ labeled by~$T$ is sent to~$\psi(T)$.
\end{proposition}

\begin{proof}
We prove the result by contradiction. Assume that there exists a \dimensional{4} point set~$V \eqdef \ens{v_T}{T\; 2\text{-triangulation of the heptagon}}\subset\R^4$ labeled by the \ktri{2}s of the heptagon such that:
\begin{enumerate}[(i)]
\item the faces of~$\conv(V)$ are faces of~$\Delta_{n,k}$; and
\item the orthogonal projection on the plane~$\R x+\R y$ sends~$v_T$ to~$\psi(T)$.
\end{enumerate}
Then, projecting the point set~$V$ only partially on~$\R x+\R y+\R z$ yields a \dimensional{3} point set~$W \eqdef \ens{w_T}{T\; 2\text{-triangulation of the heptagon}}\subset\R^3$ with the same properties~(i) and~(ii). We prove that such a point configuration~$W$ cannot exist.

The argument is based on a careful study of the projections of the faces of~$\Delta_{7,2}^\polar$. For any face~$F$ of~$\Delta_{7,2}^\polar$ we denote by~$\psi(F) \eqdef \ens{\psi(T)}{T\in F}$ its projection on the plane. For example, we have represented in \fref{ft:fig:loday72facets} the projections of all $1$- and $2$-codimensional faces of~$\Delta_{7,2}^\polar$: the left column represents, for each \krel{2} edge~$e$ of the heptagon, the set of \ktri{2}s which contain~$e$; the other three columns represent, for each pair of \krel{2} edges~$\{e,f\}$ of the heptagon, the set of \ktri{2}s which contain~$\{e,f\}$.

We denote by~$\cF^+$ the set of faces~$F$ of~$\Delta_{7,2}^\polar$ such that~$\conv\ens{w_T}{T\in F}$ is an upper facet of~$\conv W$ (with respect ot the last coordinate~$z$). Similarly, let~$\cF^-$ denote the set of faces of~$\Delta_{7,2}^\polar$ corresponding to lower facets of~$\conv W$. The projections of the upper facets of~$\conv W$ define a tiling of the polytope~$\Omega_{7,2}$, \ie they satisfy:
\begin{enumerate}[(i)]
\item \defn{covering property}: $\Omega_{7,2}=\bigcup_{F\in\cF^+}\conv\psi(F)$; and
\item \defn{intersection property}: for all~$F,F'\in\cF^+$, the intersection $\conv\psi(F)\cap\conv\psi(F')$ is either empty, or an edge or a vertex of both~$\conv\psi(F)$ and~$\conv\psi(F')$, and~$F$ and~$F'$ coincide on this interesection:
$$\conv\psi(F)\cap\conv\psi(F')\cap\psi(F)=\conv\psi(F)\cap\conv\psi(F')\cap\psi(F').$$
\end{enumerate}
A case analysis proves that we have only ten solutions for this tiling, represented in \fref{ft:fig:loday72facettilings}. Thus, we already know all possibilities for the upper and lower parts~$\cF^+$ and~$\cF^-$.

\begin{figure}[t]
	\capstart
	\centerline{\includegraphics{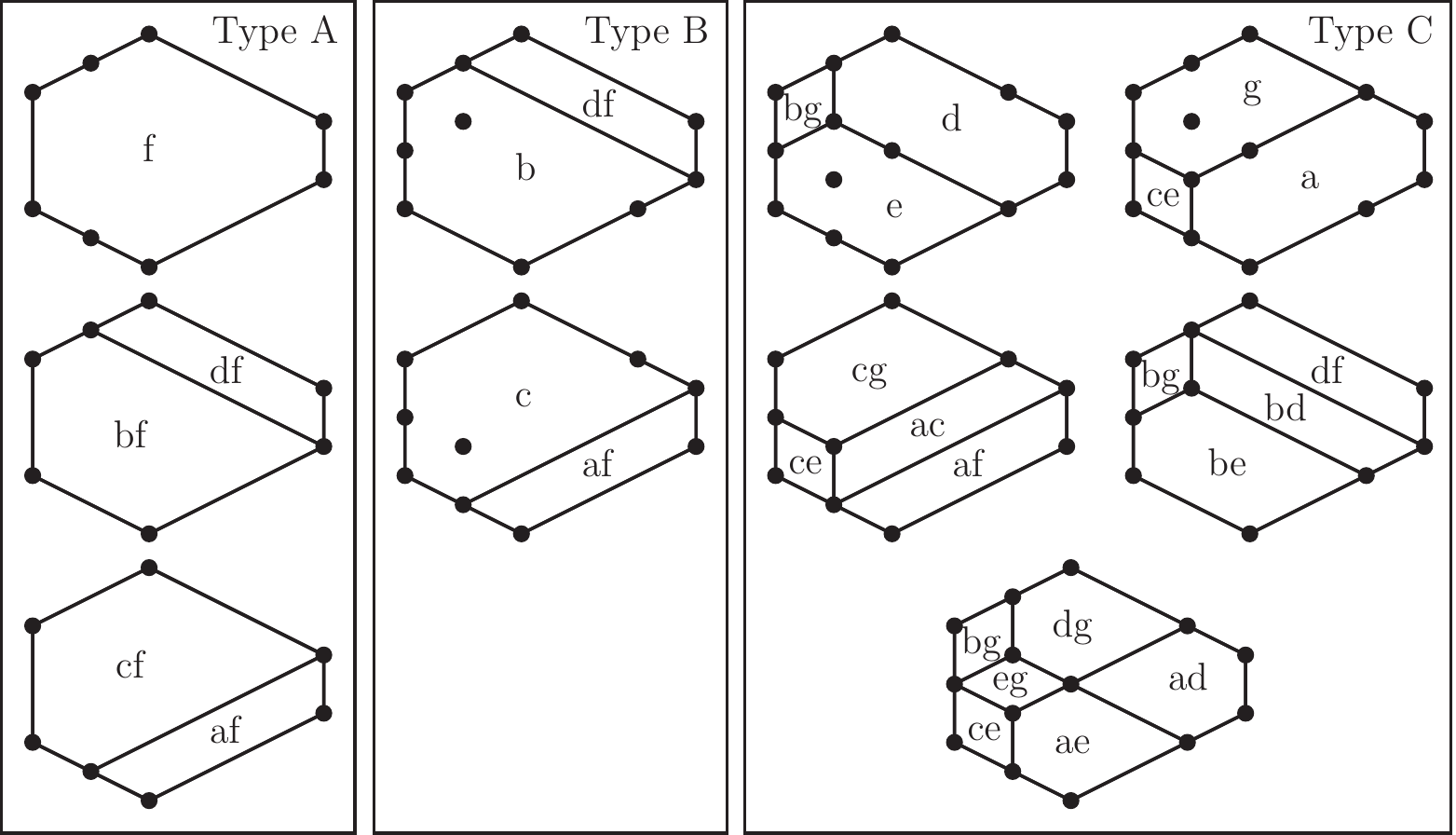}}
	\caption[The ten possible solutions for the upper and lower hull]{The ten possible solutions for the upper and lower hull of~$W$.}
	\label{ft:fig:loday72facettilings}
\end{figure}

Now, the upper and lower hulls of~$W$ have to be glued together, using if necessary some vertical facets. We have five possible vertical facets, which appear in the third column of \fref{ft:fig:loday72facets}: any triple of points of~$W$ which projects to an edge of~$\Omega_{7,2}$ either is colinear in~$W$, or forms a vertical triangular facet of~$\conv W$. If it is colinear, this triple belongs to one facet of~$\cF^+$ and to one facet of~$\cF^-$. If we have a vertical facet, then the two endpoints belong to two different facets in~$\cF^+$ (both containing the middle point), and to one facet in~$\cF^-$ which does not contain the middle point, or \viceversa. For example, the rightmost edge separates the ten possible tilings of \fref{ft:fig:loday72facettilings} into three groups: call type~A the tilings~$\rf$, $\rb\rf+\rd\rf$ and~$\rc\rf+\ra\rf$ in which the rightmost edge contains only two vertices; call type~B the tilings~$\rb+\rd\rf$ and~$\rc+\ra\rf$ in which the rightmost edge contains three vertices and is contained in a face; and call type~C the five remaining tilings, in which the rightmost edge contains three vertices and is subdivided between two faces. Then either both~$\cF^+$ and~$\cF^-$ are of type~B, or~$\cF^+$ is of type~A and~$\cF^-$ is of type~C (or \viceversa). It is easy to check that in fact there exists no pair of tilings in \fref{ft:fig:loday72facettilings} which are compatible for all the boundary edges. This implies that~$W$ cannot exist, and thus, finishes the proof of the proposition.
\end{proof}


\subsection{Pseudoline arrangements with the same support}\label{ft:subsec:multiassociahedron:pseudolinearrangements}

To finish this section, we discuss the extension of the polytopality question to a more general family of simplicial complexes arising from the framework developed in Chapter~\ref{chap:mpt}:

\begin{definition}
Given the support~$\cS$ of a pseudoline arrangement, let~$\Delta(\cS)$ denote the simplicial complex whose maximal simplices are the sets of contact points of pseudoline arrangements supported by~$\cS$.
\end{definition}

For any integers~$n$ and~$k$, the simplicial complex~$\Delta_{n,k}$ is isomorphic to~$\Delta(V_n^{*k})$, where~$V_n^{*k}$ denotes the \kkernel{k} of the dual pseudoline arrangement of~$V_n$. We present in the next section the other motivating example, whose graph is the flip graph on pseudotriangulations of~a~point~set.

\subsubsection{The polytope of pseudotriangulations via rigidity properties}

Let~$P$ be a point set in general position in the plane and~$P^{*1}$ denote the first kernel of the dual arrangement of~$P$. According to the results of Chapter~\ref{chap:mpt}, the complex~$\Delta(P^{*1})$ is isomorphic to the simplicial complex of all non-crossing and pointed sets of internal edges of~$E$. The maximal elements of this complex are in bijection with the pseudotriangulations of~$P$, and its ridge graph is the graph of flips~$G(P)$. We have seen in Chapter~\ref{chap:mpt} that this graph is connected and regular (of degree~$2|P|-3-h$, where~$h$ denotes the number of points of~$P$ on its convex hull). 

\begin{figure}
	\capstart
	\centerline{\includegraphics[scale=.95]{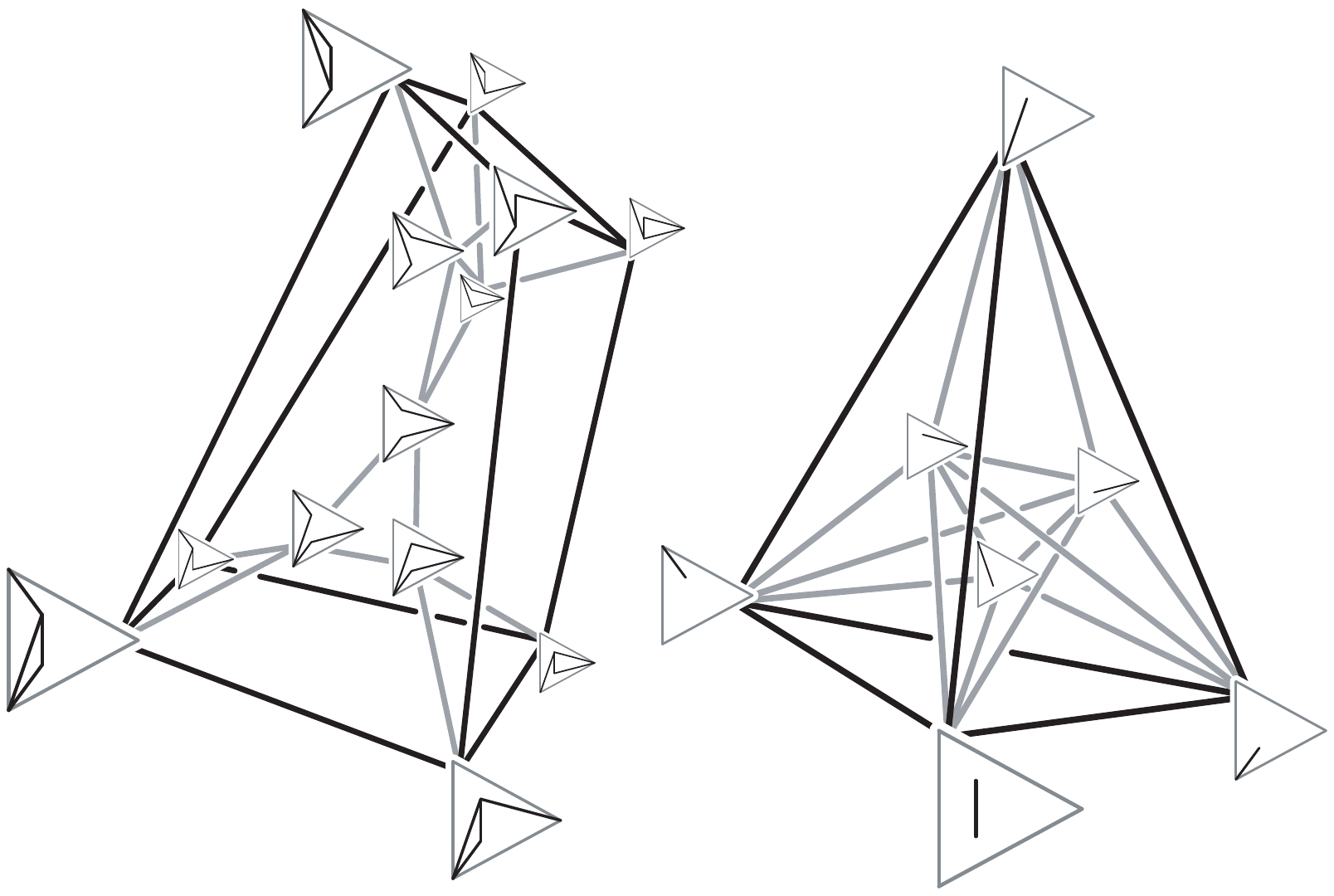}}
	\caption[The polytope of (pointed) pseudotriangulations of a configuration of five points]{The polytope of (pointed) pseudotriangulations of a configuration~$P$ of five points (left). Its dual polytope (right) realizes the simplicial complex~$\Delta(P^{*1})$. The only missing edge in~$\Delta(P^{*1})$ corresponds to the only pair of crossing edges in~$P$.}
	\label{ft:fig:pseudotriangulationPolytope}
\end{figure}

\svs
\index{polytope!--- of pseudotriangulations}
\index{pseudotriangulation!polytope of ---s}
In~\cite{rss-empppt-03}, G\"unter Rote, Francisco Santos and Ileana Streinu proved that when~$P$ is an Euclidean point set, the simplicial complex~$\Delta(P^{*1})$ is in fact polytopal (see \fref{ft:fig:pseudotriangulationPolytope}):

\begin{theorem}[\cite{rss-empppt-03}]
Let~$P$ be a point set in general position in the Euclidean plane, with~$h$ points on the convex hull. Then the simplicial complex~$\Delta(P^{*1})$ is the boundary complex of a simplicial \poly{(2|P|-3-h)}tope.\qed
\end{theorem}

\begin{remark}\label{ft:rem:rigiditypseudotriangulationspolytope}
The proof in~\cite{rss-empppt-03} relies on rigidity properties of pseudotriangulations~\cite{s-ptrmp-05}:
\begin{enumerate}[(i)]
\item As a consequence of Lemma~\ref{mpt:lem:streinu} and Theorem~\ref{ft:theo:rigidity2d}, all pseudotriangulations are generically minimally rigid in the plane. In fact, every planar minimally generically rigid graph can be embedded as a pseudotriangulation~\cite{horsssssw-pmrgpt-05}.
\item The framework obtained by removing one boundary edge from a pseudotriangulation is an expansive mechanism\footnote{An \defn{infinitesimal motion}\index{motion!infinitesimal ---} of a framework~$(P,E)$ is an assignement~$v:P\to\R^2$ of a velocity to each point such that $\dotprod{p-q}{v(p)-v(q)}=0$ for any edge~$(p,q)\in E$ of the framework. Infinitesimal motions correspond to the first derivatives of the motions of~$(P,E)$. We quotient the space of infinitesimal motions by the subspace of rigid ones, that is, by the first derivatives of rigid motions of the whole plane.

An infinitesimal motion is said to be~\defn{expansive}\index{motion!expansive ---} if all distances increase: $\dotprod{p-q}{v(p)-v(q)}\ge 0$ for all~$p,q\in P$. A \defn{mechanism} is a framework whose space of infinitesimal motions is \dimensional{1}. A mechanism is \defn{expansive} if its (unique up to a scalar) infinitesimal motion is expansive.}.
\end{enumerate}

The polytope of pseudotriangulations is obtained as a perturbation of the cone
$$\cE(P) \eqdef \ens{v:P\to\R^2}{\dotprod{p-q}{v(p)-v(q)}\ge 0 \text{ for all } p,q\in P}$$
of infinitesimal expansive motions of~$P$. An extreme ray of the \defn{expansion cone}~$\cE(P)$ corresponds to an expansive mechanism, and thus, to many ``pseudotriangulations minus one boundary edge'' (which have the same rigid components). A suitable choice of perturbation of the inequalities of the expansion cone~$\cE(P)$ is used to ``separate'' all pseudotriangulations which correspond to the same expansive mechanism. The vertices of the resulting polyhedron are in bijection with the pseudotriangulations of~$P$ and are all contained in a bounded face, whose boundary complex is the dual complex of~$\Delta(P^{*1})$.
\end{remark}

\begin{remark}\label{ft:rem:rigiditypolytope}
In much the same way as for pseudotriangulations, rigidity properties of multitriangulations would provide insight to the question of the polytopality of~$\Delta_{n,k}$. We shortly present the ideas of the possible construction via the rigidity matrix (see~\cite{g-cf-01,f-gga-04} for rigidity notions and~\cite{rss-empppt-03} for the construction).

Consider a generic embedding~$P:[n]\to\R^d$ of~$n$ points in a \dimensional{d} space. We associate to each pair~$\{u,v\}\in{[n] \choose 2}$ its \defn{rigidity vector} $\phi_{u,v} \eqdef (f_1,f_2,\dots,f_n)\in(\R^d)^n=\R^{dn}$ defined by~$f_u \eqdef P(u)-P(v)$, $f_v \eqdef P(v)-P(u)$ and~$f_i \eqdef 0\in\R^d$ for all~$i\in[n]\ssm\{u,v\}$. This vector configuration naturally lives in the orthogonal complement~$X$ (of dimension~$dn-{d+1 \choose 2}$) of the space of rigid infinitesimal motions.
The rigidity properties of a graph~$G$ on~$n$ vertices can be read on its \defn{rigidity configuration}~$\phi(G) \eqdef \ens{\phi_{u,v}}{\{u,v\}\in G}$:
\begin{enumerate}[(i)]
\item $G$ is infinitesimally rigid if and only if $\phi(G)$ spans~$X$; and
\item $G$ is \defn{stress-free} (\ie no edge of~$G$ can be removed preserving the same rigid components) if and only if~$\phi(G)$ is independent.
\end{enumerate}
Consequently, the minimally generically rigid graphs on~$n$ vertices in dimension~$d$ are precisely the graphs whose rigidity configuration~$\phi(G)$ forms a basis of~$X$.

Assume now that we can prove that every \ktri{k} of the \gon{n} is generically rigid in dimension~$2k$. Then the cone hull~$C(T)$ of the rigidity configuration~$\phi(T)$ of each \ktri{k}~$T$ of the \gon{n} is a full-dimensional simplicial cone in the vector space~$X$. Since each edge in a \ktri{k} can be uniquely flipped, the complex~$\Theta_{n,k}$ formed by all cones~$C(T)$ associated to all \ktri{k}s of the \gon{n} is a pure pseudo-manifold (every facet of every full-dimensional simplicial cone is a facet of precisely two cones).

Assume furthermore that there exists a generic embedding~$P:[n]\to\R^{2kn}$ such that for every \ktri{k}s~$T$ and~$T'$ of the \gon{n} related by a flip, in the expansive mechanism~$(P,T\cap T')$, the edge of~$T\ssm T'$ expands when the edge of~$T'\ssm T$ contracts. Then in the complex~$\Theta_{n,k}$, two cones which share a facet are separated by this facet. This would imply that the simplicial complex~$\Delta_{n,k}$ is a complete simplicial fan. And this simplicial fan would hopefully be regular, so that it would be the face fan of a simplicial polytope.

This was our initial motivation to study rigidity properties of multitriangulations.
\end{remark}

\subsubsection{The polytope of pseudoline arrangements with acyclic dual graphs}

In this section, we observe that our construction of the beam polytope of Section~\ref{ft:subsec:multiassociahedron:loday} can be extended to the framework developed in Chapter~\ref{chap:mpt}. We skip the proofs of this extension (because they are essentially the same as those of Section~\ref{ft:subsec:multiassociahedron:loday}, but also because we prove in the next section that multitriangulations are sufficiently universal to cover the general setting).

Let~$\cS$ be the support of a pseudoline arrangement in the M\"obius strip, and~$\chi$ be a cut of~$\cS$. When we cut the M\"obius strip along~$\chi$, we obtain a representation of~$\cS$ on the band~${[0,\pi)\times\R}$. Consider now a pseudoline arrangement~$\Lambda$ supported by~$\cS$, and denote by~$\lambda_1,\dots,\lambda_n$ the pseudolines of~$\Lambda$, such that their points on the vertical axis are ordered from top to bottom. We associate to~$\Lambda$ the vector~$\psi(\Lambda)\in\R^n$ whose $i$th coordinate is the number of faces of the arrangement located below the pseudoline~$\lambda_i$. Finally, we define the polytope~$\Omega(\cS,\chi)$ as the convex hull of the vectors~$\psi(\Lambda)$ associated to all pseudoline arrangements~$\Lambda$ supported by~$\cS$.

The crucial property of the beam polytope remains valid in this setting: if~$\Lambda$ and~$\Lambda'$ are two pseudoline arrangements with the same support related by a flip which involve their $i$th and $j$th pseudolines, then the difference~$\psi(\Lambda)-\psi(\Lambda')$ is parallel to~$e_i-e_j$. This implies in particular that~$\Omega(\cS,\chi)$ is contained in an hyperplane orthogonal to~$\one$.

To describe the vertices of~$\Omega(\cS,\chi)$, we need to consider the oriented dual multigraph~$\Lambda^\dual$ of a pseudoline arrangement~$\Lambda$ supported by~$\cS$. This graph has one vertex for each pseudoline of~$\Lambda$ and an edge between two vertices for each contact point between the corresponding pseudolines of~$\Lambda$. Each edge is oriented from the pseudoline passing above the contact point to the pseudoline passing below it. With this definition, the cone of~$\Omega(\cS,\chi)$ at the vertex~$\psi(\Lambda)$ equals the incidence cone~$C(\Lambda^\dual)$ of the oriented dual multigraph~$\Lambda^\dual$. In particular, the vector~$\psi(\Lambda)$ is a vertex of~$\Omega(\cS,\chi)$ if and only if the oriented graph~$\Lambda^\dual$ is acyclic. Moreover, this vertex is simple when the Hasse diagram of~$\Lambda^\dual$ is a tree.

Finally, the normal vectors of the facets of~$\Omega(\cS,\chi)$ are precisely the normal vectors of all the cones~$C(\Lambda^\dual)$, that is, all the characteristic vectors of the sink sets of the minimal directed cuts of the oriented dual graphs~$\Lambda^\dual$.

\subsubsection{Universality of the multiassociahedron}

We finish our discussion concerning the multiassociahedron with the proof of a universality result. The goal of this section is to observe that the family of multiassociahedra~$(\Delta_{n,k})_{n\ge 2k+1}$ together with all their faces already contains full generality: any simplicial complex of the form~$\Delta(\cS)$ is a coface of a multiassociahedron~$\Delta_{n,k}$ for some~$n,k\in\N$. More precisely:

\begin{theorem}[Universality of the multiassociahedron]\label{ft:theo:universality}
If~$\cS$ is the support of an arrangement of~$n$ pseudolines with~$m$ intersection points (crossing or contact points), then the simplicial complex~$\Delta(\cS)$ is isomorphic to a coface of~$\Delta_{n+2m-2,m-1}$.
\end{theorem}

\begin{figure}
	\capstart
	\centerline{\includegraphics[scale=1]{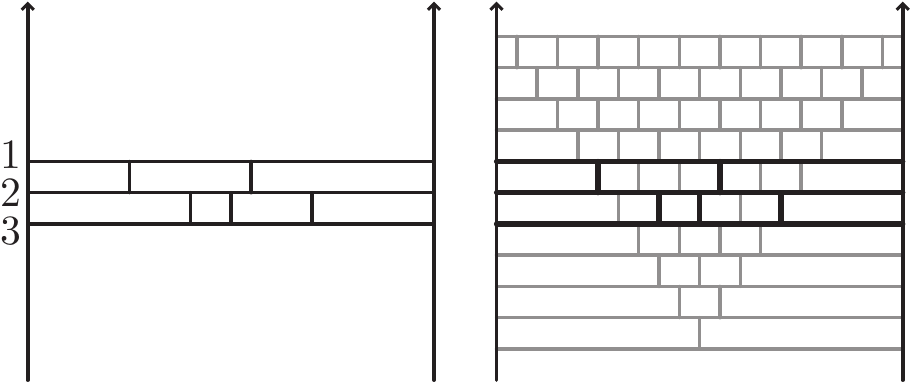}}
	\caption[Universality of multitriangulations]{Universality of multitriangulations: the wiring diagram of the left pseudoline arrangement is a subdiagram of the wiring diagram in convex position.}
	\label{ft:fig:universality}
\end{figure}

\begin{proof}
We choose an arbitrary sweep of the support~$\cS$ and order the vertices of~$\cS$ according to it. We represent the result as a \defn{wiring diagram} of~$\cS$ (see~\cite[Chapter~6]{f-gga-04} for the definition and \fref{ft:fig:universality}): a vertex of~$\cS$ is represented by a (vertical) comparator between two consecutive (horizontal) supporting lines. Index the supporting lines with~$[n]$ from top to bottom, and for any~$i\in[m]$, let~$i^\square$ be such that the $i$th comparator join the horizontal lines~$i^\square$ and~$i^\square+1$. Define the set~$W \eqdef \ens{[i-1,i+i^\square+m-2]}{i\in[m]}$.

We claim that the simplicial complex~$\Delta(\cS)$ is isomorphic to the simplicial complex of all \kcross{m}-free subsets of \krel{(m-1)} edges of the \gon{(n+2m-2)} which contain the complement of~$W$. Indeed, by duality, the \ktri{(m-1)}s of this complex correspond to the pseudoline arrangements supported by~$V_{n+2m-2}^{*(m-1)}$ (\ie the \kkernel{(m-1)} of the  arrangement of~$n+2m-2$ pseudolines in convex position) whose sets of contact points contain~$\ens{u^*\wedge v^*}{[u,v]\notin W}$. Thus, our claim follows from the observation that, by construction, the remaining contact points~$\ens{u^*\wedge v^*}{[u,v]\in W}$ are positioned exactly as those of~$\cS$ (see \fref{ft:fig:universality}).
\end{proof}

\begin{corollary}
Any pseudoline arrangement is (isomorphic to) the dual pseudoline arrangement of a multitriangulation.\qed
\end{corollary}

\begin{figure}[h]
	\capstart
	\centerline{\includegraphics[scale=1]{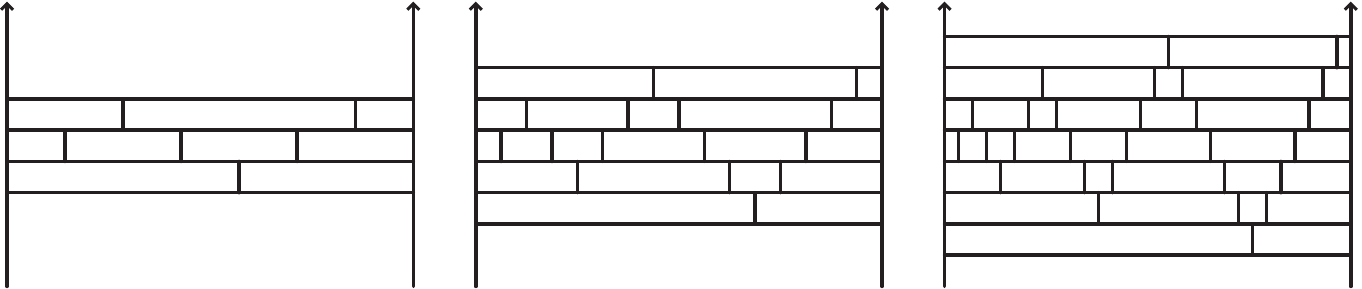}}
	\caption[A family of pseudoline arrangements whose representing multitriangulation is necessarily quadratic]{A family of pseudoline arrangements whose representing multitriangulation is necessarily quadratic.}
	\label{ft:fig:waves}
\end{figure}

\begin{remark}
Theorem~\ref{ft:theo:universality} ensures that an arrangement of~$n$ pseudolines with no contact point is a multitriangulation of an \gon{N}, with~$N\le n^2+2n-2$. The family of pseudoline arrangements suggested by \fref{ft:fig:waves} is an example for which~$N$ is necessarily quadratic in~$n$.
\end{remark}

By universality, all the results and the conjectures on ``multitriangulations of a convex polygon'' translate to \defn{equivalent} results and conjectures on ``pseudoline arrangements with a given support''. Let us give three examples to close our discussion:

\begin{corollary}
The simplicial complex~$\Delta(\cS)$ associated to the support~$\cS$ of a pseudoline arrangement is a combinatorial sphere.\qed
\end{corollary}

\begin{conjecture}
The shortest path in the flip graph~$G(\cS)$ between two pseudoline arrangements~$\Lambda$ and~$\Lambda'$ supported by~$\cS$ never flips a common contact point of~$\Lambda$ and~$\Lambda'$.
\end{conjecture}

\begin{conjecture}
Given any support~$\cS$ of a pseudoline arrangement, the simplicial complex~$\Delta(\cS)$ is the boundary complex of a polytope.
\end{conjecture}



\begin{appendices} 

\chapter{Two enumeration algorithms}\label{app:implementations}

This appendix is devoted to the presentation of the implementation of two enumeration algorithms that we mentioned in Chapters~\ref{chap:mpt} and~\ref{chap:multiassociahedron}: the first one enumerates arrangements of pseudolines and double pseudolines, and the second one deals with symmetric matroid polytopes realizing the simplicial complex of \kcross{(k+1)}-free sets of \krel{k} edges of the \gon{n}.


\section{Double pseudoline arrangements}\label{app:sec:dpl}

\index{enumeration algorithm!--- for double pseudoline arrangements}
Pseudoline arrangements (or in higher dimension, pseudohyperplane arrangements) have been extensively studied in the last decades as a useful combinatorial abstraction of configurations of points~\cite{bvswz-om-99,b-com-06,k-ah-92,g-pa-97}. Recently, Luc~Habert and Michel~Pocchiola~\cite{hp-adp-08}, motivated by the question of efficient visibility graph algorithms for convex shapes, introduced double pseudoline arrangements as a combinatorial abstraction of configurations of disjoint convex bodies. The main structural properties of pseudoline arrangements~---~connectedness of their mutation graph, embedability in projective or affine geometries, axiomatic characterization in terms of chirotopes~---~extend to double pseudoline arrangements. Definitions and results of~\cite{hp-adp-08} are recalled in Section~\ref{app:subsec:dpl:preliminaries}.

\begin{figure}
	\capstart
	\centerline{\includegraphics[scale=1]{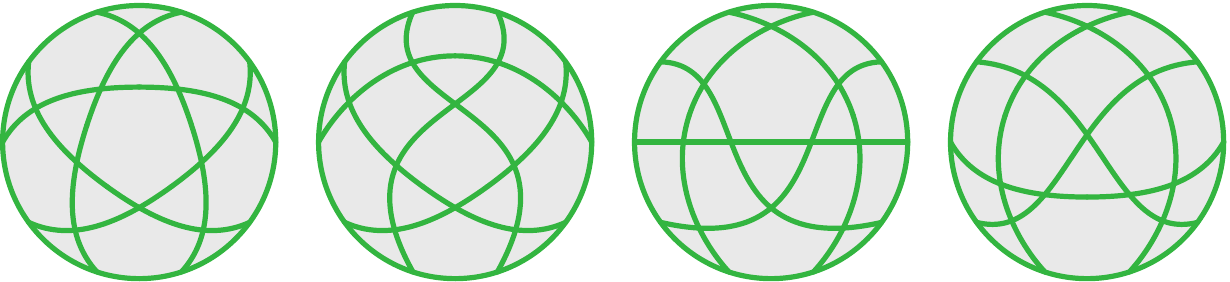}}
	\caption[Simple projective arrangements of six pseudolines]{The four isomorphism classes of simple projective arrangements of six pseudolines.}
	\label{app:fig:6pseudolines}
\end{figure}

To help our understanding of double pseudoline arrangements, and in order to carry out computer experiments, it is interesting to develop algorithms to enumerate arrangements with few double pseudolines. For pseudoline arrangements, or more generally for pseudohyperplane arrangements, different enumeration algorithms have been proposed and implemented~\cite{bg-gom-00,ff-gomgt-02,aak-eotsp-02}. The number $p_n$ of isomorphism classes of simple pseudoline arrangements in the projective plane is given by the following the table:
\begin{center}
\begin{tabular}{c|ccccccccccc}
$n$ & $1$ & $2$ & $3$ & $4$ & $5$ & $6$ & $7$ & $8$ & $9$ & $10$ & $11$ \\
\hline
$p_n$ & $1$ & $1$ & $1$ & $1$ & $1$ & $4$  & $11$ & $135$ & $4\,382$ & $312\,356$ & $41\,848\,591$
\end{tabular}
\end{center}
For example, \fref{app:fig:6pseudolines} gives one representative for each of the four isomorphism classes of simple projective arrangements of six pseudolines.

In this section, we present the implementation of an incremental algorithm to enumerate simple projective arrangements of double pseudolines. Before the description of the algorithm in Section~\ref{app:subsec:dpl:algo}, we shortly recall some preliminaries on arrangements in Section~\ref{app:subsec:dpl:preliminaries}, and prove in Section~\ref{app:subsec:dpl:extensions} the connectedness property on which our incremental algorithm relies. Finally in Section~\ref{app:subsec:dpl:variations}, we briefly discuss three possible variations of this algorithm: to mixed arrangements (with both pseudolines and double pseudolines), to arrangements in the M\"obius strip, and to non-simple arrangements (where three pseudolines can meet at the same vertex).


\subsection{Preliminaries}\label{app:subsec:dpl:preliminaries}

\subsubsection{Pseudoline and double pseudoline arrangements}

\index{projective!--- plane|hbf}
\index{plane!projective ---|hbf}
We denote the \defn{projective plane} by~$\cP$ and represent it as a disk with antipodal boundary points identified.

\begin{figure}[!h]
	\capstart
	\centerline{\includegraphics[scale=1]{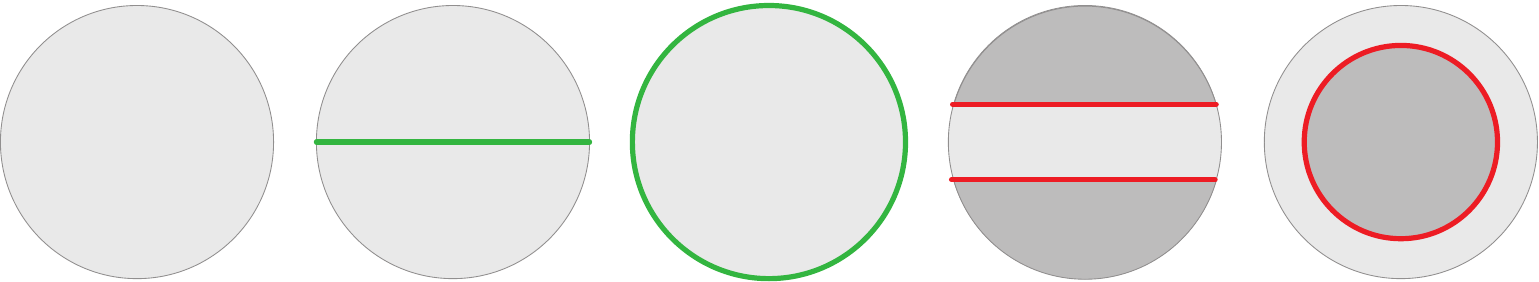}}
	\caption[A projective plane, two pseudolines, and two double pseudolines]{A projective plane, two pseudolines, and two double pseudolines (their enclosed topological disk is represented darker).}
	\label{app:fig:pldpl}
\end{figure}

A simple closed curve in~$\cP$ is a \defn{pseudoline}\index{pseudoline} if it is non-separating (or equivalently, non-contractible), and a \defn{double pseudoline}\index{double pseudoline|hbf} otherwise (see \fref{app:fig:pldpl}). The complement of a pseudoline is a topological disk. The complement of a double pseudoline~$\ell$ has two connected components: a M\"obius strip~$M_\ell$ and a topological disk~$D_\ell$.

\begin{figure}[!h]
	\capstart
	\centerline{\includegraphics[scale=1]{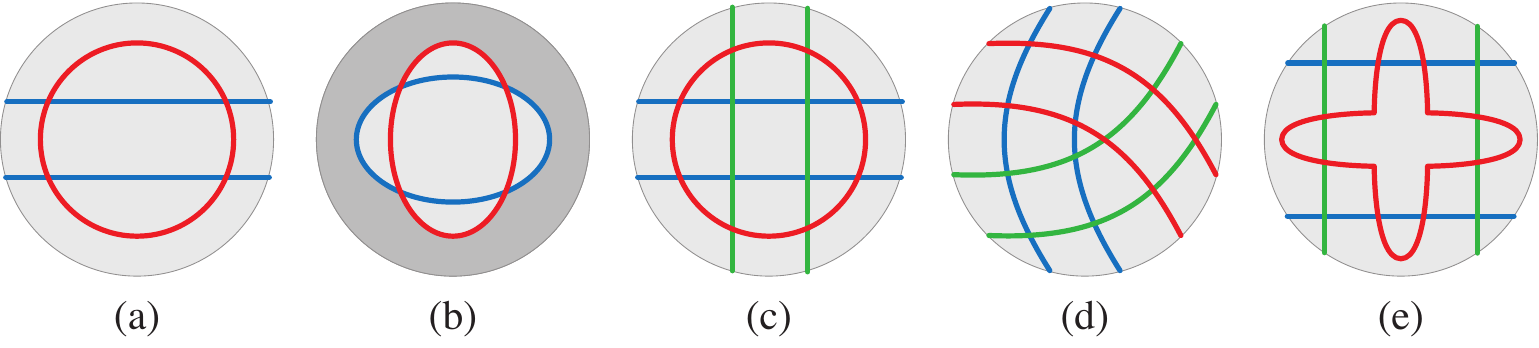}}
	\caption[Double pseudoline arrangements]{(a) An arrangement of $2$~double pseudolines; (b) Two double pseudolines which do not form an arrangement; (c-d-e) Arrangements of $3$~double pseudolines.}
	\label{app:fig:dplarrangements}
\end{figure}

 A \defn{pseudoline arrangement}\index{pseudoline!--- arrangement} is a finite set of pseudolines such that any two of them intersect exactly once (see \fref{app:fig:6pseudolines}). A \defn{double pseudoline arrangement}\index{double pseudoline!--- arrangement|hbf} is a finite set of double pseudolines such that any two of them have exactly four intersection points, cross transversally at these points, and induce a cell decomposition of~$\cP$ (see \fref{app:fig:dplarrangements}). To simplify the presentation, we first focus only on \defn{simple} arrangements, that is, where no three curves meet at the same point. Non-simple arrangements will be discussed in Section~\ref{app:subsec:dpl:variations}.

Two arrangements~$L$ and~$L'$ are \defn{isomorphic} if there is an homeomorphism of the projective plane that sends~$L$ on~$L'$ (or equivalently, if there is an isotopy joining~$L$ to~$L'$). For example, the arrangements (c)~and~(d) in \fref{app:fig:dplarrangements} are isomorphic while (e)~is not isomorphic to them. As mentioned previously, different algorithms have been implemented to enumerate isomorphism classes of pseudoline arrangements, and their number is know up to $11$~pseudolines. In this section, we develop an algorithm to enumerate isomorphism classes of double pseudolines arrangements.

\subsubsection{Mutations}

\begin{figure}
	\capstart
	\centerline{\includegraphics[scale=1]{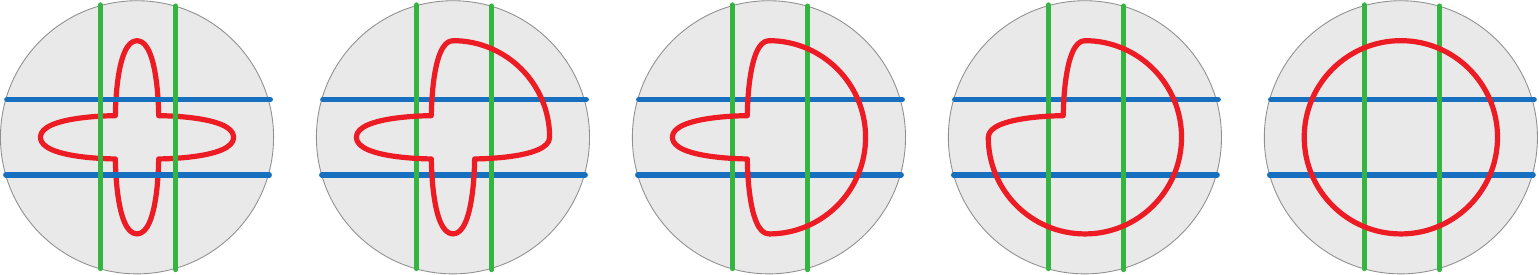}}
	\caption[Mutations]{A sequence of mutations from the double pseudoline arrangement of \fref{app:fig:dplarrangements}(e) to the double pseudoline arrangement of \fref{app:fig:dplarrangements}(c). All the vertices in the enclosed M\"obius strip of the ``flower'' double pseudoline of the left arrangement are pumped out.}
	\label{app:fig:mutations}
\end{figure}

A \defn{mutation}\index{mutation|hbf} is a local transformation of an arrangement~$L$ that only inverts a triangular face of~$L$ (see \fref{app:fig:mutations}). More precisely, it is a homotopy of arrangements in which only one curve~$\ell$ moves, sweeping a single vertex of the remaining arrangement~$L\smallsetminus\{\ell\}$. Iterating this elementary transformation is in fact sufficient to obtain all possible arrangements:

\begin{theorem}[\cite{hp-adp-08}]\label{app:theo:connectednessArrangements}
The graph of mutations on arrangements is connected:  any two arran\-gements with the same numbers of (double) pseudolines are homotopic via a finite sequence of mutations followed by an isotopy.
\end{theorem}

This result is known as Ringel's Homotopy Theorem for arrangements of pseudolines. For double pseudolines, it is proved in~\cite{hp-adp-08} by reduction to the case of pseudoline arrangements. The argument consists in mutating each double pseudoline of the arrangement until it becomes \defn{thin}\index{double pseudoline!--- arrangement! thin --- ---} (that is, until there remains no point inside its enclosed M\"obius strip), and to observe that arrangements of thin double pseudolines behave exactly as pseudoline arrangements. The fact that an arrangement can be mutated until all its double pseudolines become thin is ensured by the following crucial lemma (see \fref{app:fig:mutations}):

\begin{lemma}[Pumping Lemma~\cite{hp-adp-08}]\label{app:lem:pumping1}
In an arrangement~$L$, if a double pseudoline~$\ell$ encloses at least one vertex of~$L$ in its enclosed M\"obius strip~$M_\ell$, then there is a triangular face of~$L$ included in~$M_\ell$ and supported by~$\ell$.
\end{lemma}

From Theorem~\ref{app:theo:connectednessArrangements}, we may derive a simple enumeration algorithm which consists in traversing the graph of mutations starting from any given arrangement. This naive algorithm is sufficient for the enumeration of small cases but already fails (because of RAM memory limitations) for arrangements of five double pseudolines. In order to go a little bit further (and particularly, to enumerate arrangements of five double pseudolines), we use an incremental version of this algorithm~---~still based on mutations~---~presented in Sections~\ref{app:subsec:dpl:extensions} and~\ref{app:subsec:dpl:algo}.

\subsubsection{Geometric Representation Theorem}

A \defn{projective geometry}\index{projective!--- geometry} is a couple~$(\PP,\LL)$, where~$\PP$ is a projective plane and~$\LL$ is a set of pseudolines of~$\PP$, such that any two pseudolines of~$\LL$ have exactly one intersection point (which depends continuously on the two lines) and through any pair of points of~$\PP$ passes a unique pseudoline of~$\LL$ (which depends continuously on the two points). In a projective geometry, the \defn{dual}\index{dual!--- of a point} of a point~$p$ (of~$\PP$) is its set~$p^*$ of incident pseudolines (of~$\LL$), and the \defn{dual}\index{dual!--- of a point set} of a point set~$P$ is the set~$P^* \eqdef \ens{p^*}{p\in P}$ of duals of points of~$P$. Similarly, the \defn{dual}\index{dual!--- of a convex body} of a convex body~$C$ (of~$\PP$) is its set~$C^*$ of tangent pseudolines (of~$\LL$), and the \defn{dual}\index{dual!--- of a set of convex bodies} of a set~$S$ of convex bodies is the set~$S^* \eqdef \ens{C^*}{C\in S}$ of duals of bodies of~$S$. It is known that:
\begin{enumerate}[(i)]
\item The line space~$\LL$ is a projective plane.
\item The dual of a finite point set is a pseudoline arrangement (which is simple if the point set is in general position).
\item The dual of a finite set of pairwise disjoint convex bodies is a double pseudoline arrangement (which is simple if no three convex share a common tangent).
\end{enumerate}
The \defn{Geometric Representation Theorem}\index{Geometric Representation Theorem} asserts that the converse is also true:

\begin{theorem}[Geometric Representation Theorem]
Any arrangement of pseudolines (resp.~of double pseudolines) is isomorphic to the dual arrangement of a finite point set (resp.~of a finite set of pairwise disjoint convex bodies) of a projective geometry.
\end{theorem}

For pseudoline arrangements, this theorem is a consequence of duality in projective geometries combined with the embedabbility of any pseudoline arrangement in the line space of a projective geometry~\cite{gpwz-atp-94}. It was extended in~\cite{hp-adp-08} for double pseudoline arrangements. From this result, it is easy to derive the following extension of the Pumping Lemma~\ref{app:lem:pumping1}, which will be useful in our algorithm:

\begin{lemma}[Improved Pumping Lemma]\label{app:lem:pumping2}
Let~$L$ be an arrangement, $\ell$~be a double pseudoline of~$L$ and $x$~be a point on~$\ell$. If the M\"obius strip~$M_\ell$ enclosed by~$\ell$ contains at least one vertex of~$L$, then there is a triangular face of~$L$ included in~$M_\ell$, supported by~$\ell$, and avoiding the point~$x$.
\end{lemma}

\subsubsection{Chirotopes}

An \defn{indexed oriented arrangement} is an arrangement whose curves are oriented and one-to-one indexed with an indexing set. The \defn{chirotope}\index{chirotope} of an indexed oriented arrangement is the application that assigns to each triple~$\tau$ of indices the isomorphism class of the subarrangement indexed by~$\tau$. Since there are only two indexed and oriented simple arrangements of three pseudolines, it is convenient to call them~$+$ and~$-$ and to define a chirotope as an application from triples of indices to $\{+,-\}$ (this relates to the other definition of chirotopes given in the context of oriented matroids in Chapter~\ref{chap:multiassociahedron} and used in Appendix~\ref{app:sec:enumerationmatroidpolytopes}). For double pseudoline arrangements, there are much more indexed and oriented simple arrangements of three double pseudolines. However, appart from the number, chirotopes of pseudoline arrangements and double pseudoline arrangements behave exactly the same. On the one hand, the chirotope of an arrangement characterizes its isomorphism class:

\begin{theorem}[\cite{hp-adp-08}]
The isomorphism class of an indexed oriented arrangement only depends on its chirotope.
\end{theorem}

On the other hand, chirotopes are characterized by their restrictions to subconfigurations of size at most five:

\begin{theorem}[\cite{hp-adp-08}]\label{app:theo:chirotopes}
Given an application~$\chi$ that assigns to each triple~$\tau$ of indices an isomorphism class of an oriented arrangement indexed by~$\tau$, the following properties are equivalent:
\begin{enumerate}[(i)]
\item $\chi$~is the chirotope of an indexed oriented arrangement; and
\item the restriction of~$\chi$ to the set of triples of any subset of at most five indices is the chirotope of an indexed oriented arrangement.
\end{enumerate}
\end{theorem}

This result provides a strong motivation for enumerating arrangements of at most five simple and double pseudolines: the complete list of these arrangements gives an axiomatization of the chirotopes of all arrangements.


\subsection{Connectedness of the one-extension space}\label{app:subsec:dpl:extensions}

Our enumeration algorithm is incremental, and relies on the following objects:

\begin{definition}
A \defn{one-extension}\index{one-extension (of an arrangement)} of an arrangement~$L$ of~$n$ double pseudolines is an arrangement of~$n+1$ double pseudolines of which~$L$ is a subarrangement. The added double pseudoline is called the \defn{one-extension curve}.
\end{definition}

A mutation of a one-extension of~$L$ is certainly still a one-extension of~$L$ when the moving curve is the one-extension curve. We call such a mutation a \defn{one-extension mutation}. It is crucial for our algorithm that these restricted mutations still generate the complete one-extension space:

\begin{theorem}
The graph of one-extension mutations is connected: any two one-extensions of an arrangement are homotopic via a finite sequence of mutations (followed by an isotopy) during which the only moving curve is the one-extension curve.
\end{theorem}

\proof
Let~$L$ be an arrangement,~$\lambda$ be a distinguished double pseudoline of~$L$, and~$M \eqdef L\cup\{\mu\}$ and~$N \eqdef L\cup\{\nu\}$ be two one-extensions of~$L$. Using continuous motions~---~thanks to the duality principle~---~we can assume the following facts on~$\mu$:
\begin{enumerate}[(i)]
\item $\mu$ is a \defn{thin} double pseudoline in~$M$ (\ie its enclosed M\"obius strip contains no vertex of~$M$);
\item $\lambda$~and~$\mu$ are \defn{touching}: a cell~$\sigma$ of size two of the subarrangement~$\{\lambda,\mu\}$ is also a cell of the whole arrangement~$M$;
\item the M\"obius strip~$M_\mu$ enclosed by $\mu$ is a tubular neighborhood of a pseudoline~$\mu_*$ which intersects any double pseudoline of~$L$ in exactly two points, the intersection being~\mbox{transversal}.
\end{enumerate}
Similarly, we can assume that~$\nu$ is a thin double pseudoline, which touches~$\lambda$, and whose enclosed M\"obius strip defines a tubular neighborhood of a pseudoline~$\nu_*$ which intersects (transversally) any double pseudoline of~$L$ in two points. Furthermore, we can restrict the analysis to the following case:
\begin{enumerate}[(i)]
\setcounter{enumi}{3}
\item $\mu$~and~$\nu$ coincide in the topological disk~$D_\lambda$;
\item $\mu_*$~and~$\nu_*$ coincide in~$D_{\lambda}$, have finitely many intersections in~$M_\lambda$, and intersect transversally at these points;
\item $M_\mu\cup M_\nu$ forms a tubular neighborhood of~$\mu_*\cup\nu_*$.
\end{enumerate}

We consider the arrangements~$M_* \eqdef L\cup\{\mu_*\}$ and~$N_* \eqdef L\cup\{\nu_*\}$. As in the proof of Theorem~\ref{app:theo:connectednessArrangements}, the proof of this theorem boils down to show that the mixed arrangements~$M_*$ and~$N_*$ are homotopic via a finite sequence of mutations during which the only moving curves are the pseudolines~$\mu_*$ and~$\nu_*$. We prove this by induction on the number of vertices in the M\"obius strip~$M_\lambda$ enclosed by~$\lambda$:
\begin{itemize}
\item if~$\lambda$ is a thin double pseudoline, then we are done up to an isotopy;
\item otherwise, the Improved Pumping Lemma~\ref{app:lem:pumping2} ensures that there exists a triangular face~$\Delta$ of the arrangement supported by~$\lambda$ and not adjacent to the ``touching cell''~$\sigma$. We can assume that~$\mu_*$ and~$\nu_*$ do not intersect~$\Delta$, modulo a finite sequence of mutations in~$M_*$~and~$N_*$, during which the only moving curves are the pseudolines~$\mu_*$~and~$\nu_*$ (these mutations are possible since for each connected component~$X$ of the trace of a double pseudoline of~$L$ on~$M_\lambda$, both $\mu_*$~and~$\nu_*$ intersect~$X$ at most once). We can then perform the mutation of~$\Delta$: let $\tilde{M}_*$ and~$\tilde{N}_*$ denote the mutated arrangements. By induction hypothesis, there is a finite sequence of mutations which transforms~$\tilde{M}_*$ to~$\tilde{N}_*$ only moving the pseudolines~$\mu_*$ and~$\nu_*$. But this sequence also transforms~$M_*$ to~$N_*$, which finishes the proof.\qed
\end{itemize}


\subsection{The incremental algorithm}\label{app:subsec:dpl:algo}

\subsubsection{Description}

Let~$\cA_n$ be the set of isomorphism classes of arrangements of~$n$ double pseudolines and~$q_n \eqdef |\cA_n|$.
Our algorithm enumerates~$\cA_n$ from~$\cA_{n-1}$ by mutating an added distinguished double pseudoline.

A \defn{pointed} arrangement is an arrangement with a distinguished double pseudoline. For any integer~$n$, let~$\cA_n^\bullet$ denote the set of isomorphism classes of pointed arrangements of~$n$ double pseudolines (two pointed arrangements are isomorphic if there is an isomorphism of the projective plane which sends one arrangement to the other, respecting their distinguished double pseudoline). We always use the notation~$A^\bullet$ for a pointed arrangement and~$A$ for its non-pointed version (and similarly for sets of pointed arrangements). We also use the notation $\subarr(A)$ for the set of subarrangements of an arrangement~$A$.

For each~$i\in\{1,\ldots,q_{n-1}\}$, the algorithm enumerates the set~$S_i^\bullet$ of arrangements of~$\cA_n^\bullet$ containing the $i$th arrangement~$a_i$ of~$\cA_{n-1}$ as a subarrangement, by mutations of a distinguished added double pseudoline. From the set~$S_i$, it selects the subset~$R_i$ of arrangements with no subarrangements in~$\{a_1,\ldots,a_{i-1}\}$. In other words,~$R_i$ is the set of arrangements of~$\cA_n$ whose first subarrangement, among the set $\cA_{n-1} \eqdef \{a_1,\ldots,a_{q_{n-1}}\}$, is the $i$th arrangement~$a_i$. Thus,~$\cA_n$ is the disjoint union of the sets~$R_i$.


\svs
\noindent
\framebox[\textwidth]{
\begin{minipage}{.95\textwidth}
\svs
\noindent\textsc{Incremental enumeration}
\svs
\ligne
\svs

\begin{algorithmic}

\REQUIRE $\cA_{n-1}=\{a_1,\ldots,a_{q_{n-1}}\}$.
\ENSURE $\cA_n$.

\medskip
\FOR{$i$ from $1$ to $q_{n-1}$}
	\STATE $A^\bullet\leftarrow\;$add a pointed double pseudoline to $a_i$.
	\STATE $Q^\bullet\leftarrow[A^\bullet]$. $S^\bullet\leftarrow\{A^\bullet\}$. $R\leftarrow\{\}$.
	\IF{$\subarr(A)\cap \{a_1,\ldots, a_{i-1}\}=\emptyset$}
		\STATE write $A$. $R\leftarrow\{A\}$.
	\ENDIF
	\WHILE{$Q^\bullet\ne\emptyset$}
		\STATE $A^\bullet\leftarrow \text{pop } Q^\bullet$.
		\STATE $T\leftarrow\;$list the triangles of $A^\bullet$ adjacent to its pointed double pseudoline.
		\FOR{$t\in T$}
			\STATE $B^\bullet\leftarrow\;$mutate the triangle $t$ in $A^\bullet$. 
			\IF{$B^\bullet\notin S^\bullet$} 
				\IF{$\subarr(B)\cap \{a_1,\ldots, a_{i-1}\}=\emptyset$ and $B\notin R$}
					\STATE write $B$. $R\leftarrow R\cup\{B\}$.
				\ENDIF
				\STATE $ \text{push}(B^\bullet,Q^\bullet)$. $S^\bullet\leftarrow S^\bullet\cup\{B^\bullet\}$.
			\ENDIF	
		\ENDFOR
	\ENDWHILE
\ENDFOR
\end{algorithmic}
\svs
\end{minipage}
}


\svs
An alternative approach\footnote{We thank Luc Habert for this suggestion.} for counting arrangements is to enumerate the subsets~$S_i^\bullet$ and to compute, for each arrangement~$A^\bullet$ of~$S_i^\bullet$, the number~$\sigma(A^\bullet)$ of double pseudolines~$\alpha$ of~$A$ such that~$A$ pointed at~$\alpha$ is isomorphic to~$A^\bullet$. Then:
$$q_n=\frac{1}{n}\sum_{i=1}^{q_{n-1}} \sum_{A^\bullet\in S_i^\bullet} \sigma(A^\bullet).$$
The disadvantage of this version is that it only counts~$q_n$ and cannot provide a data base for~$\cA_n$.

\subsubsection{Encoding an arrangement}

In order to manipulate arrangements, one can represent them in several different ways. We use two different encodings:
\begin{enumerate}
\item \defn{Flag representation}:
An arrangement defines a map on a surface (the projective plane), and thus can be manipulated with the corresponding tools. A \defn{flag} of an arrangement~$L$ is any triple~$(v,e,f)$ consisting of one vertex~$v$, one edge~$e$ and one face~$f$ (of the cell complex defined by~$L$), such that~$v\in e\subset f$. The geometric information of the arrangement is contained in the three involutions~$\sigma_0$,~$\sigma_1$ and~$\sigma_2$ which respectively change the vertex, the edge and the face of a flag (see \fref{app:fig:sigma}). This representation is convenient to perform all the necessary elementary operations we need to apply on arrangements (such as mutations, extensions, \etc).

\begin{figure}[h]
	\capstart
	\centerline{\includegraphics[scale=1]{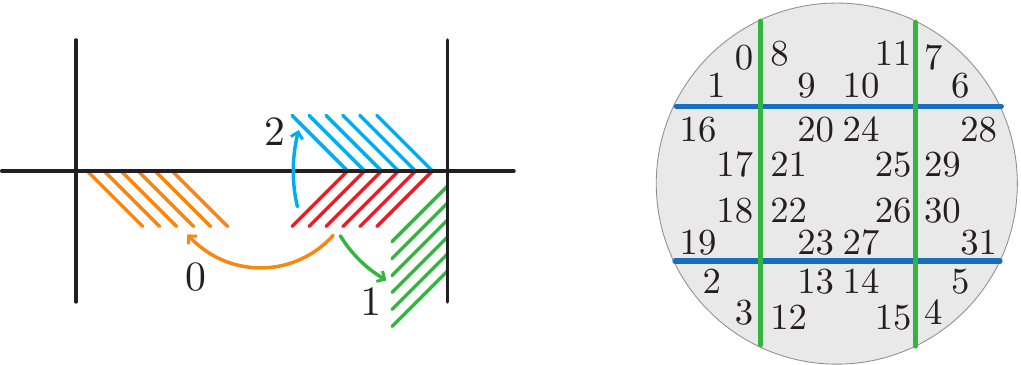}}
	\caption[The three involutions $\sigma_0$, $\sigma_1$ and $\sigma_2$]{The three involutions $\sigma_0$, $\sigma_1$ and $\sigma_2$ (left). An example of flag representation (right).}
	\label{app:fig:sigma}
\end{figure}

\item \defn{Encoding}: 
In the enumeration process, once we have finished to manipulate an arrangement~$L$, we still have to remember it and to test whether the new arrangements we find afterwards are or not isomorphic to~$L$. For this, we compute another representation of~$L$, which is shorter than the flag representation, and which allows a quick isomorphism test.

We first associate to each flag~$\phi=(v,e,f)$ of the arrangement~$L$ a word $w_\phi$ constructed as follows. Let~$\ell_1$ be the simple or double pseudoline containing~$e$, and, for all~${2\le p\le n}$, let~$\ell_p$ be the~$p$th double pseudoline crossed by~$\ell_1$ on a walk starting at~$\phi$ and oriented by~$\sigma_0\sigma_1\sigma_2\sigma_1$. This walk also defines a starting flag~$\phi_i$ for each~$\ell_i$. We then walk successively on~$\ell_1,\dots,\ell_n$, starting from~$\phi_i$ and in the direction given by~$\sigma_0\sigma_1\sigma_2\sigma_1$, and index the vertices by~$1,2,\ldots,V$ in the order of appearance. For all~$i$, let~$w_i$ denote the word formed by reading the indices of the vertices of~$\ell_i$ starting from~$\phi_i$. The word~$w_\phi$ is the concatenation of~$w_1,w_2,\ldots,w_n$. Finally, we associate to the arrangement~$L$ the lexicographically smallest word among all the~$w_\phi$ where~$\phi$ ranges over all the flags of~$L$. It is easy to check that two arrangements are isomorphic if and only if they get the same encoding.
\end{enumerate}

\subsubsection{Adding a double pseudoline}

\begin{figure}
	\capstart
	\centerline{\includegraphics[scale=1]{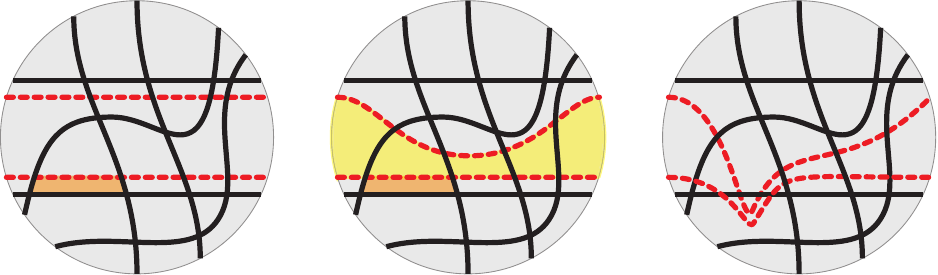}}
	\caption[Three steps to insert a double pseudoline in an arrangement]{Three steps to insert a double pseudoline in an arrangement: duplicate a double pseudoline (left), flatten it (middle) and add four crossings (right).}
	\label{app:fig:addADPL}
\end{figure}

One of the important steps of the incremental method is to add a double pseudoline to an initial arrangement. It is easy to achieve when the initial arrangement contains at least one pseudoline, but it is more involved when we have only double pseudolines. Our method uses three steps (see \fref{app:fig:addADPL}):
\begin{enumerate}[(a)]
\item \defn{Duplicate a double pseudoline}:
We choose one arbitrary double pseudoline~$\ell$ and duplicate it, drawing a new double pseudoline~$\ell'$ completely included in the M\"obius strip~$M_\ell$ enclosed by~$\ell$. We choose an arbitrary rectangle~$R$ delimited by~$\ell$ and~$\ell'$.
\item \defn{Flatten}: We pump the double pseudoline~$\ell'$ (using Lemma~\ref{app:lem:pumping2}) such that no vertex of the arrangement lies in the M\"obius strip~$M_{\ell'}$. During the process, we do not touch the rectangle~$R$.
\item \defn{Add four crossings}: We replace the rectangle~$R$ by four crossings between~$\ell$ and~$\ell'$.
\end{enumerate}

If we think of our double pseudoline arrangement as the dual of a configuration of convex bodies, this method corresponds to: (a) choosing one convex~$C$ and drawing a new convex~$C'$ inside~$C$; (b) flattening the convex~$C'$ until it becomes almost a single point, maintaining it almost in contact with the boundary of $C$; and (c) moving~$C'$ outside~$C$.


\subsection{Results}

We have implemented this algorithm in~\verb$C++$ programming language. The documentation (as well as the source code) of this implementation is available upon request by email.

\begin{figure}
	\capstart
	\centerline{\includegraphics[scale=1]{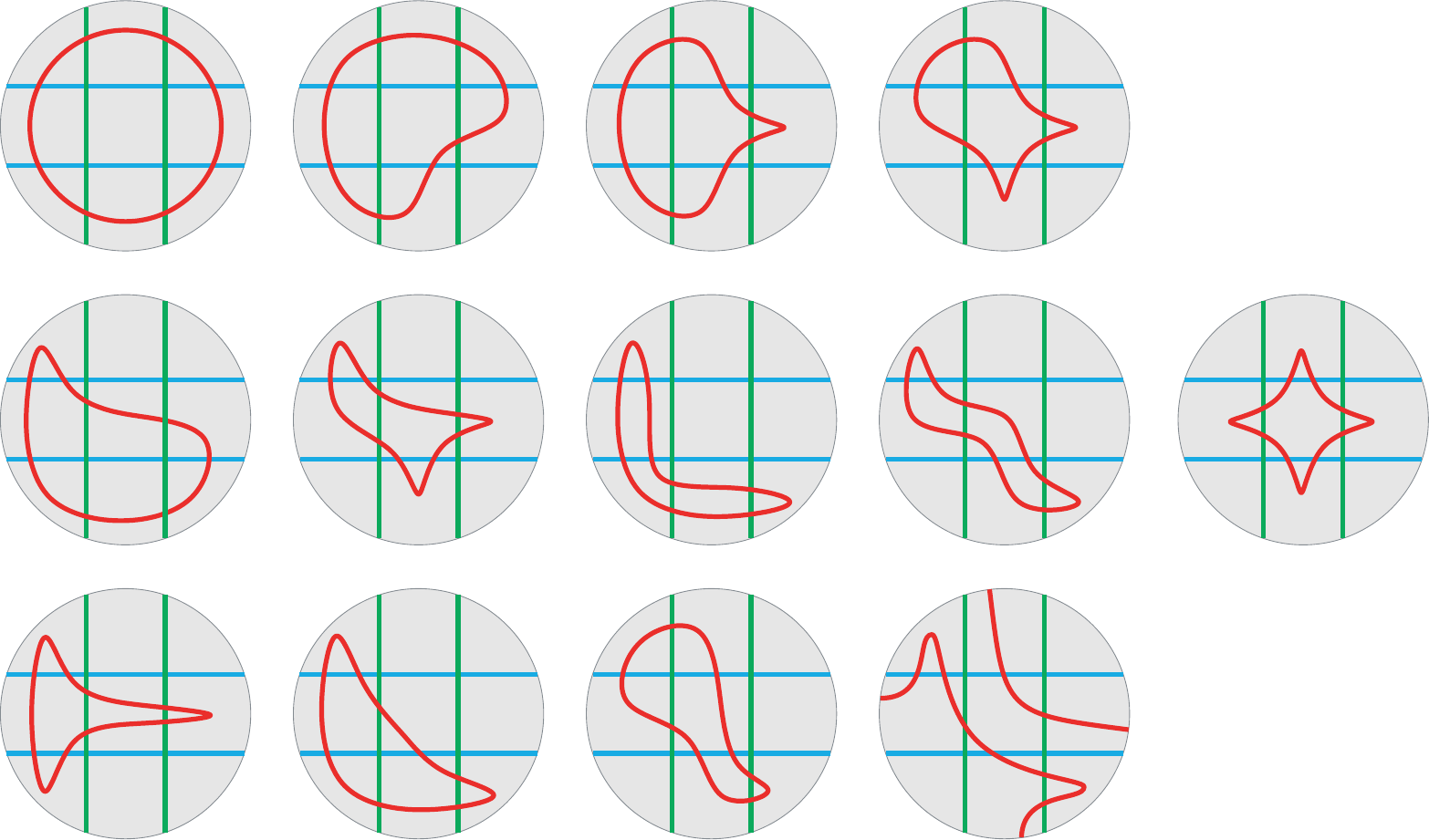}}
	\caption[Simple projective arrangements of three double pseudolines]{The thirteen isomorphism classes of simple projective arrangements of three double pseudolines.}
	\label{app:fig:3doublepseudolines}
\end{figure}

%
This implementation provided us with the following values of the number~$q_n$ of isomorphism classes of simple arrangements of~$n$ double pseudolines:
\begin{center}
\begin{tabular}{c|ccccc}
$n$ & $1$ & $2$ & $3$ & $4$ & $5$ \\
\hline
$q_n$ & $1$ & $1$ & $13$ & $6570$ & $181\,403\,533$
\end{tabular}
\end{center}
For example, \fref{app:fig:3doublepseudolines} gives one representative for each of the thirteen isomorphism classes of simple projective arrangements of three double pseudolines.

\svs
Let us briefly comment on running time and working space. Observe first that our algorithm can be parallelized very easily (separating each enumeration of~$S_i$ and~$R_i$, for~$i\in\{1,\ldots,q_{n-1}\}$). In order to obtain the number~$q_5$ of simple arrangements of five double pseudolines, we used four $2$GHz processors for almost three weeks. The working space is bounded by~${\max|S_i^\bullet|=279\,882}$ (times the space of the encoding of a single configuration, \ie about~$80$ characters).

\begin{figure}
	\capstart
	\centerline{\includegraphics[scale=.99]{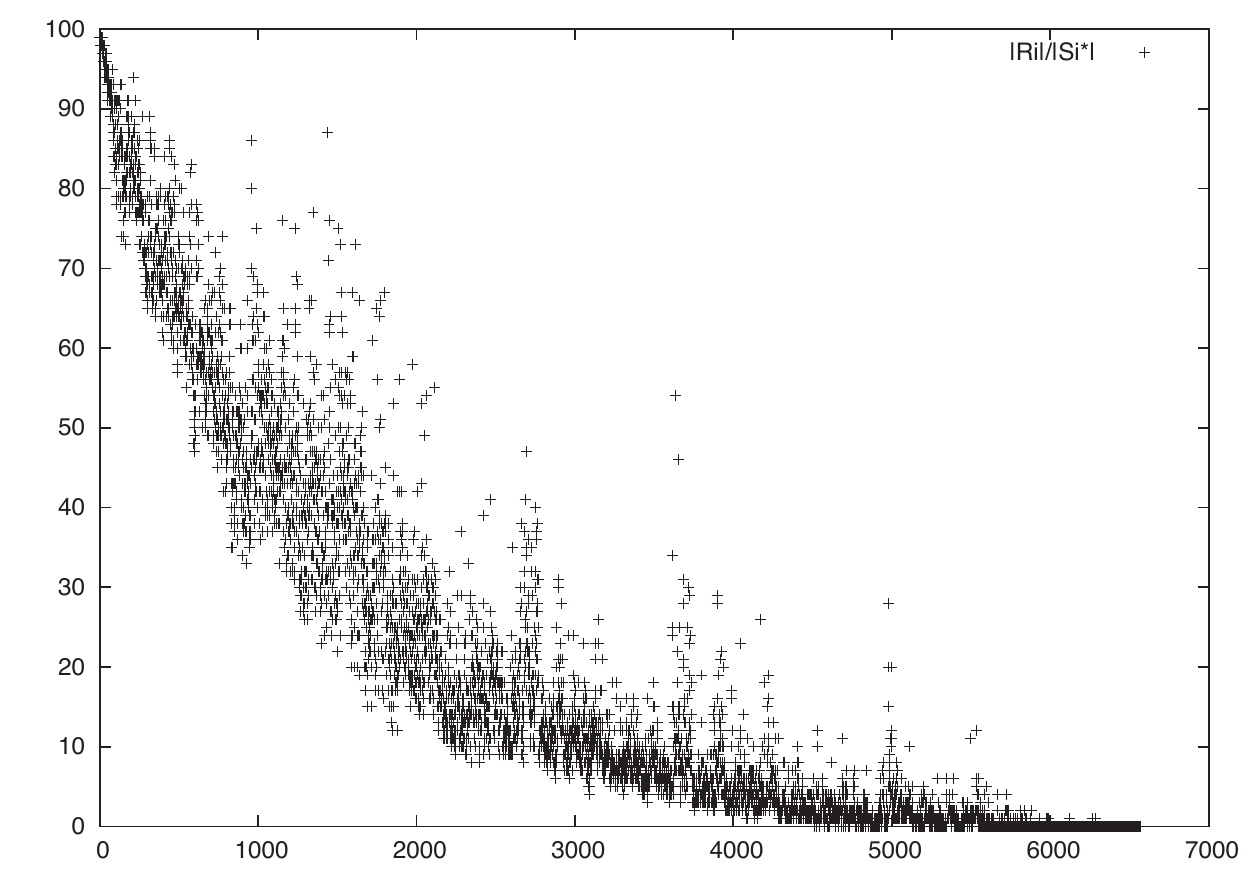}}
	\caption[Percentage of new configurations during the enumeration process]{Percentage~$|R_i|/|S_i^\bullet|$ of new configurations among all generated configurations during the recursive enumeration of arrangements of five double pseudolines.}
	\label{app:fig:statistics}
\end{figure}

\fref{app:fig:statistics} shows the evolution of the ratio between the sizes of the sets~$R_i$ and~$S_i^\bullet$, during the enumeration of arrangements of five double pseudolines. It is also interesting to observed that $\sum |S_i^\bullet| / \sum |R_i|\simeq 5$, which confirms that hardly any configurations of five convex bodies have symmetries.

\subsection{Three variations}\label{app:subsec:dpl:variations}


\subsubsection{Mixed arrangements}

\begin{figure}[!h]
	\capstart
	\centerline{\includegraphics[scale=1]{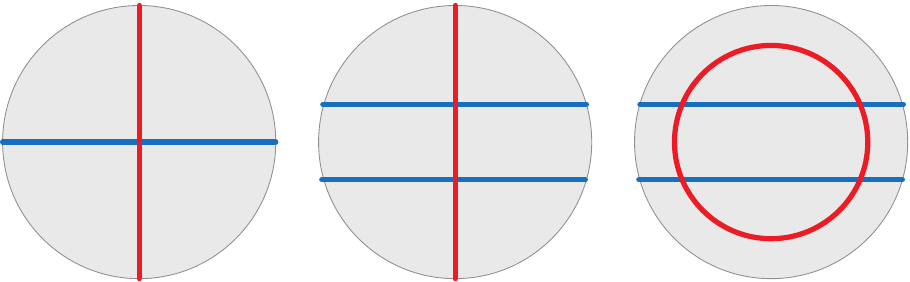}}
	\caption[Mixed arrangements of size~$2$]{The three mixed arrangements of size~$2$.}
	\label{app:fig:mixedarrangements}
\end{figure}

An \defn{arrangement of simple and double pseudolines} (or \defn{mixed arrangement}) is a finite set of pseudolines and double pseudolines such that:
\begin{enumerate}[(i)]
\item any two pseudolines have a unique intersection point;
\item a pseudoline and a double pseudoline have exactly two intersection points and cross transversally at these points; and \item any two double pseudolines have exactly four intersection points, cross transversally at these points, and induce a cell decomposition of~$\cP$.
\end{enumerate}
\fref{app:fig:mixedarrangements} shows these three conditions, \ie the three admissible mixed arrangements of size~$2$.

All the properties of arrangements~---~connectedness of their mutation graph, embedability in projective  geometries, axiomatic characterization in terms of chirotopes~---~remain valid in this context.

As already used in the proofs, there is a correspondence between the isomorphism classes of arrangements with~$n$ pseudolines and~$m$ double pseudolines, and isomorphism classes of arrangements of~$n+m$ double pseudolines, $n$ of which are thin (\ie do not contain any vertex of the arrangement in there enclosed M\"obius strip).
%

Mixed arrangements can be enumerated by the same algorithm as before. The numbers of mixed (simple projective) arrangements with~$n$ pseudolines and~$m$ double pseudolines are given by the following table:
\begin{center}
\begin{tabular}{c|cccccc}
\begin{picture}(15,15)(0,0)
	\put(9,9){$m$}
	\put(0,0){$n$}
	\put(4,4){$\ssm$}
\end{picture}
  & $0$ & $1$ & $2$ & $3$ & $4$ & $5$  \\
\hline
$0$ &  &  & $1$ & $13$ & $6\,570$ & $181\,403\,533$\\
$1$ & & $1$ & $4$ & $626$ & $4\,822\,394$ \\
$2$ & $1$ & $2$ & $48$ & $86\,715$ \\
$3$ & $1$ & $5$ & $1\,329$ \\
$4$ & $1$ & $25$ & $80\,253$ \\
$5$ & $1$ & $302$ \\
$6$ & $4$ & $9\,194$ \\
$7$ & $11$ & $556\,298$ \\
$8$ & $135$ \\
$9$ & $4\,382$ \\
$10$ & $312\,356$
\end{tabular}
\end{center}

\subsubsection{M\"obius arrangements}

If we consider configurations of points and convex bodies in a plane (rather than projective) geometry, then we need to enumerate mixed arrangements in the M\"obius strip (rather than the projective plane). Since a M\"obius strip is a projective plane minus one point (at infinity), we enumerate mixed arrangements in the M\"obius strip with essentially the same method. The only difference is that we have to maintain one additional information: we mark the external face, which contains the point at infinity. This information is used to deal with the two following constraints: first, the external face can never be mutated; second, two M\"obius arrangements are isomorphic if there is an isomorphism sending one to the other, and conserving the  external face. The numbers of mixed simple M\"obius arrangements with~$n$ pseudolines and~$m$ double pseudolines are given by the following table:
\begin{center}
\begin{tabular}{c|cccccc}
\begin{picture}(15,15)(0,0)
	\put(9,9){$m$}
	\put(0,0){$n$}
	\put(4,4){$\ssm$}
\end{picture}
  & $0$ & $1$ & $2$ & $3$ & $4$ & $5$  \\
\hline
$0$ &  &  & $1$ & $16$ & $11\,502$ & $238\,834\,187$\\
$1$ & & $1$ & $7$ & $1\,499$ & $9\,186\,477$ \\
$2$ & $1$ & $3$ & $140$ & $245\,222$ \\
$3$ & $1$ & $13$ & $5\,589$ \\
$4$ & $2$ & $122$ & $416\,569$ \\
$5$ & $3$ & $2\,445$ \\
$6$ & $16$ & $102\,413$ \\
$7$ & $135$ & $7\,862\,130$ \\
$8$ & $3\,315$ \\
$9$ & $158\,830$ \\
\end{tabular}
\end{center}

\subsubsection{Non-simple arrangements}

\begin{figure}[!h]
	\capstart
	\centerline{\includegraphics[scale=1]{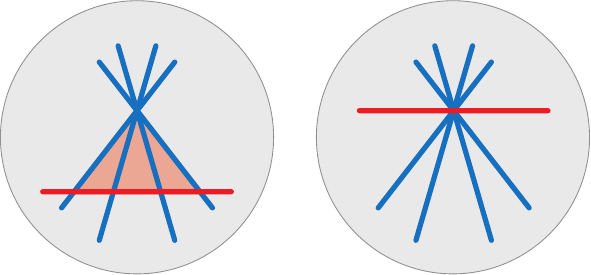}}
	\caption[Half mutations in a non-simple arrangement]{Half mutations in a non-simple arrangement.}
	\label{app:fig:halfmutations}
\end{figure}

Finally, our algorithm and our implementation can easily be extended to non-simple arrangements (\ie arrangements in which more than two curves can cross at a given vertex. The main difference lies in the possibility of the two following inverse \defn{half mutations} (see \fref{app:fig:halfmutations}):
\begin{enumerate}
\item \defn{Leaving half mutation}: A curve~$\ell$ passing through a non-simple vertex~$v$ is pushed out of~$v$ and creates a fan with apex~$v$ and based on~$\ell$.
\item \defn{Reaching half mutation}: A fan with apex~$v$ based on a curve~$\ell$ is destroyed, the curve~$\ell$ being pushed to be incident to the vertex~$v$.
\end{enumerate}

The numbers of mixed (non-necessary simple) arrangements with~$n$ pseudolines and~$m$ double pseudolines are given by the following table:
\begin{center}
\begin{tabular}{c|ccccc}
\begin{picture}(15,15)(0,0)
	\put(9,9){$m$}
	\put(0,0){$n$}
	\put(4,4){$\ssm$}
\end{picture}
  & $0$ & $1$ & $2$ & $3$ & $4$  \\
\hline
$0$ &  &  & $1$ & $46$ & $153\,528$ \\
$1$ & & $1$ & $9$ & $6\,998$ \\
$2$ & $1$ & $3$ & $265$ \\
$3$ & $1$ & $16$ & $18\,532$ \\
$4$ & $2$ & $159$ \\
$5$ & $4$ & $4\,671$ \\
$6$ & $17$ & $342\,294$ \\
$7$ & $143$ \\
$8$ & $4\,890$ \\
$9$ & $461\,053$ \\
\end{tabular}
\end{center}


\subsection{Further developments}\label{app:subsec:dpl:furtherDevelopments}

To close this section, we mention three possible further developments related to this algorithm:
\begin{enumerate}
\item \defn{Drawing an arrangement}: We have seen a method to add a pseudoline in a mixed arrangement. Combined with a planar-graph-drawing algorithm, this provides an algorithm to draw an arrangement in the unit disk. For example, the number of isomorphism classes of mixed arrangements with~$1$ pseudoline and~$n$ double pseudolines can be interpreted as the number of drawings of the arrangements of $\cA_n$. As an illustration, we have represented in \fref{app:fig:drawings} the four isomorphism classes of arrangements with~$1$ pseudoline and~$3$ double pseudolines, together with the corresponding planar drawings.

\begin{figure}[!h]
	\capstart
	\centerline{\includegraphics[width=\textwidth]{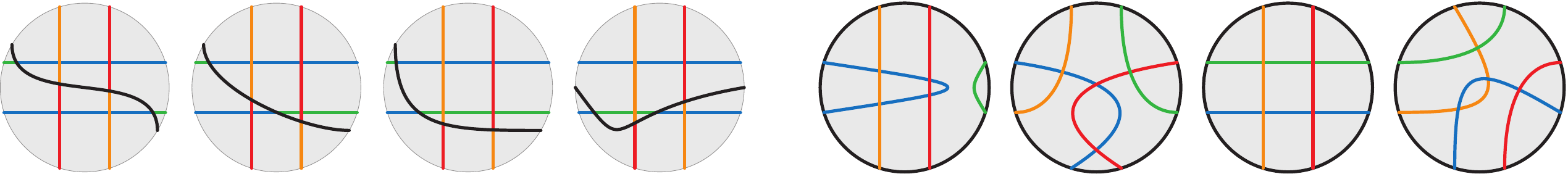}}
	\caption[Drawings of arrangements]{Drawings of arrangement: the four isomorphism classes of arrangements with~$1$ pseudoline and~$3$ double pseudolines (left) and the corresponding planar drawings obtained by cutting the projective plane along the pseudoline of each arrangement (right).}
	\label{app:fig:drawings}
\end{figure}

\item \defn{Axiomatization}: Pseudoline arrangements admit simple axiomatizations~\cite{bvswz-om-99,k-ah-92}, with few axioms dealing with configurations of at most five pseudolines. Theorem~\ref{app:theo:chirotopes} affirms that the complete list of arrangements of at most five simple and double pseudolines is an axiomatization of mixed arrangements. Is there any simpler axiomatization? Is it possible to algorithmically reduce this axiomatization?
\item \defn{Realizability}: It is well-known that certain pseudoline arrangements are not realizable in the Euclidean plane. Inflating pseudolines into thin double pseudolines in such an arrangement gives rise to non-realizable double pseudoline arrangements. Are there smaller examples? Are all arrangements of at most five curves realizable? 
\end{enumerate}


\newpage
\section{Symmetric matroid polytopes}\label{app:sec:enumerationmatroidpolytopes}

\index{enumeration algorithm!--- for symmetric matroid realizations}
\index{realization!symmetric matroid ---}
Our second algorithm~\cite{bp-srmt-09} enumerates the symmetric matroid realizations of the simplicial complex~$\Delta_{n,k}$ of \kcross{(k+1)}-free sets of \krel{k} edges of the \gon{n}. The matroid \mbox{realizations} it provides when~$n=8$ and~$k=2$ (Proposition~\ref{ft:prop:delta82matroid}) are the starting point of our study of symmetric geometric realizations of~$\Delta_{8,2}$ (Section~\ref{ft:subsubsec:delta82}). Since it was written with the functional language \haskell, the code is sufficiently concise and simple to be presented in this document. We start with a very brief presentation of the language features of \haskell.


\subsection{About Haskell}

The programming language \haskell, officially born in 1987, is named in honor of the mathematician and logician Haskell~B.~Curry\footnote{Haskell~B.~Curry worked on the theory of formal languages and processes. His name is in particular associated with the interpretation of a function with multiple arguments as a chain of functions with a single argument. For example, for a two-arguments function,~$f:(a,b)\mapsto c \;\simeq\; \curry(f):a\mapsto(b\mapsto c)$.}. Its particularities are the following:
\begin{itemize}
\svs
\item[~~~~\fbox{\defn{functions}}] \haskell is a \defn{purely functional programming language}, which means that its basic method of computation is the evaluation of functions on arguments. Contrarily to imperative languages, a functional language avoids state and mutable data. For example, to compute the factorial of an integer~\verb$n$, when an imperative language multiplies iteratively the~\verb$n$ first integers, \haskell defines \defn{recursively}:
\begin{center}\verb$factorial 0 = 1$ and \verb$factorial (n+1) = (n+1)*factorial n$.\end{center}

\svs
\item[~~~~\fbox{\defn{types}}] All functions and arguments in \haskell have a type (\verb$Int$, \verb$Bool$, \verb$Int -> Bool$, \etc), which can be explicitly mentioned by the user after the symbol~\verb$::$. If no type is mentioned, \haskell will anyway deduce automatically the type, checking that types are coherent. It also provides \defn{polymorphic type} (a type that can be parametrized by another type, \eg the type~\verb$[a]$ of a list containing elements of  type~\verb$a$), and \defn{type classes} (types that share common operations, \eg the~\verb$Eq$ class of all types supporting the operations~\verb$==$~and~\verb$/=$).

\svs
\item[~~~~\fbox{\defn{list comprehension}}] One of the main tools that make \haskell code readable for a non-programmer reader is the possibility to create lists based on existing lists, exactly in the same way we construct sets from other sets in mathematics. For example, given a function \verb$f :: a -> b$, a condition~\verb$cond :: a -> Bool$ and a list\footnote{In \haskell, capital letters are reserved for constructors. For the sake of readability, we have decided however to use capital letters in our presentation, in particular for lists.}~\verb$L :: [a]$, we can construct a new list 
\begin{center}\verb$[ f x | x <- L, cond x ] :: [b]$\end{center}
 which contains the images under~\verb$f$ of all elements of~\verb$L$ which satisfy the condition~\verb$cond$. The list~\verb$L$ is the generator, and the condition~\verb$cond$ is a guard. We can furthermore have multiple generators and/or guards. An interesting example is the~\verb$map$ function (whose type is~\verb$(a -> b) -> [a] -> [b]$) which applies a function~\verb$f$ to all elements of a list~\verb$L$ and can be defined as~\verb$map f L = [ f x | x <- L ]$.

\svs
\item[~~~~\fbox{\defn{lazy evaluation}}] \haskell is a lazy programming language in the sense that it delays the evaluation of an expression until the result is  required. It avoids unnecessary computation, but also enables the user to perform constructions that would not be possible otherwise, such as infinite structures. For example, finding the first integer which satisfies a condition~\verb$cond :: Int -> Bool$ can be achieved by evaluating~\verb$head [ n | n <- [1..], cond n ]$.

\svs
\item[~~~~\fbox{\defn{guarded equations \& pattern matching}}] The result of a function often depends on a condition. In \haskell, instead of writing a conditional control structure with~\verb$if$~--~\verb$then$~--~\verb$else$ like in many other programming languages, we can simply use guarded equations indicated by a~\verb$|$ symbol. For example, the code:
\begin{verbatim}
																f args | cond1  = out1
																	| cond2  = out2
																	| otherwise  = out3
\end{verbatim}
means that the function~\verb$f$, applied on the arguments~\verb$args$, returns~\verb$out1$ if~\verb$cond1$ is fulfilled, otherwise~\verb$out2$ if~\verb$cond2$ is fulfilled, and finally,~\verb$out3$ otherwise.

Finally, as many functional languages, \haskell also provides the possibility of matching the arguments of a function with a given pattern. This was the case for example when we wrote \haskell's definition of the factorial: we match its argument either with~\verb$0$ or with~\verb$n+1$. It is also useful when we program with lists.
\end{itemize}

\svs
Lists in \haskell are managed by the~\verb$ListPrelude$ and~\verb$List$ modules. The functions we use on a list~\verb$L$ are the usual functions:

\begin{center}
\begin{tabular}{|cl|cl|}
\hline
\verb$head L$ & first element of~\verb$L$ & \verb$take p L$ & takes \verb$p$~elements from the front of~\verb$L$ \\
\verb$tail L$ & all but the first element of~\verb$L$ & \verb$drop p L$ & drops \verb$p$ elements from the front of~\verb$L$ \\
\verb$length L$ & length of the list & \verb$L!!p$ & \verb$p$th~element of~\verb$L$ (\verb$0$$\,\le\,$\verb$p$$\,<\,$\verb$length L$)\\
\verb$L++M$ & concatenation of~\verb$L$ and~\verb$M$ & \verb$L\\M$ & \verb$L$~setminus~\verb$M$ \\
\verb$x:L$ & adds~\verb$x$ to the front of~\verb$L$ & \verb$delete x L$ & deletes~\verb$x$ from~\verb$L$ \\
\verb$elem x L$ & is~\verb$x$ an element of~\verb$L$? & \verb$elemIndices x L$\hspace*{-.25cm} & lists positions where~\verb$x$ occurs in~\verb$L$ \\
\verb$sort L$ & sorts~\verb$L$ in increasing order & \verb$nub L$ & eliminates repetitions in~\verb$L$ \\
\verb$and L$ & conjonction of booleans & \verb$or L$ & disjonction of booleans \\
\hline
\end{tabular}
\end{center}

To illustrate certain of these functions, and to give a first example of \haskell code, we describe in detail the function~\verb$tuples$ which, for a given integer~\verb$p$ and a list~\verb$L$, lists all ordered \verb$p$-tuples of~\verb$L$ (\ie all tuples which appear sorted as in the list):

\begin{lignes}
\begin{verbatim}
-- tuples p L = all ordered p-tuples of the list L
tuples :: Int -> [a] -> [[a]]
tuples 0 _ = [[]]
tuples _ [] = [] 
tuples p (h:t) = [ h:tu | tu <- tuples (p-1) t ] ++ tuples p t
\end{verbatim}
\end{lignes}

Let us briefly comment this code since all our functions follow the same pattern. In the first line of all our functions (starting with~\verb$--$), we give a brief description of the function. The second line declares the type of the function. The other lines define the function.

Here, the function~\verb$tuples$ is recursive and should be read as follows: an ordered \verb$p$-tuple of~\verb$L$ either starts with the first element of~\verb$L$, followed by an ordered~\verb$(p-1)$-tuple of the remaining list, or is an ordered~\verb$p$-tuple of the tail of~\verb$L$. Observe the use of pattern matching for initializations. 

This function, applied on the arguments~\verb$3 [1,4,2,5,3]$ returns the following list:

\begin{lignes}
\begin{verbatim}
> tuples 3 [1,4,2,5,3]
[[1,4,2], [1,4,5], [1,4,3], [1,2,5], [1,2,3], [1,5,3], [4,2,5], [4,2,3], [4,5,3], [2,5,3]]
\end{verbatim}
\end{lignes}

We believe that this short introduction is sufficient to understand our enumeration algorithm, and we refer to~\cite{h-ph-07,b-com-06} for further reading on \haskell programming language.


\subsection{Enumeration of all \ktri{k}s of the \gon{n}}

We first need to generate all facets of~$\Delta_{n,k}$, that is, all \ktri{k}s of the \gon{n} (remember that we only consider the \krel{k} edges of~$E_n$ since any \ktri{k} contains all non-\krel{k} edges). According to Corollary~\ref{stars:coro:starsenumeration}, one possibility for this is to enumerate all \tuple{k(n-2k-1)}s of \krel{k} edges of~$E_n$ and to keep only \kcross{(k+1)}-free such tuples. The corresponding \haskell code is very concise using list comprehension:

\begin{lignes}
\begin{verbatim}
-- Fix the parameters n and k
n, k :: Int
n = 8; k = 2

-- kRelevantEdges = all k-relevant edges of the n-gon
kRelevantEdges :: [[Int]]
kRelevantEdges = nub [ sort [i, mod (i+j) n] | j <- [k+1..n-k-1], i <- [0..n-1] ]

-- crossing e f = do e and f cross each other?
crossing :: [Int] -> [Int] -> Bool
crossing [a,b] [c,d] = and [a<c, c<b, b<d] || and [c<a, a<d, d<b]

-- kCrossingFree E = is the set of edges E (k+1)-crossing-free?
kCrossingFree :: [[Int]] -> Bool
kCrossingFree E = and [ or [ not (crossing e f) | [e,f] <- tuples 2 F ] | F <- tuples (k+1) E ]

-- kTriangulations = all k-triangulations of the n-gon
kTriangulations :: [[[Int]]]
kTriangulations = [ T | T <- tuples (k*(n-2*k-1)) kRelevantEdges, kCrossingFree T ]
\end{verbatim}
\end{lignes}

Even if this method is really na\"ive, it is sufficient for our purpose (after all, we only have to enumerate the \ktri{2}s of the octagon). However, for completeness, and to demonstrate the use of flips as local transformations, we present another way to enumerate \ktri{k}s, based on the graph of flips.

\newpage
The first step is to compute the two \kstar{k}s adjacent to a given \krel{k} edge of a \ktri{k}. According to Lemma~\ref{stars:lem:angles}, a \kstar{k} can be traced following its angles in the \ktri{k}:

\begin{lignes}
\begin{verbatim}
-- angle T [b,c] = the vertex a such that (a,b,c)={[a,b],[b,c]} is an angle of T
angle :: [[Int]] -> [Int] -> Int
angle T [b,c] = head [ a | a <- map (\i -> mod (c+i) n) [1..], mod (a+k) n == b || elem (sort [a,b]) T ]

-- completeStar T e = the k-star of T adjacent to the (k-relevant) edge e
completeStar :: [[Int]] -> [Int] -> [Int]
completeStar T L
	| length L == 2*k+1  = L
	| otherwise  = completeStar T ((angle T (take 2 L)):L)
\end{verbatim}
\end{lignes}

Observe the use of lazy evaluation in the definition of our function~\verb$angle$: \haskell will only compute the first vertex~\verb$a$ after~\verb$c$ such that the edge~\verb$[a,b]$ is an edge of the \ktri{k} (which might be a \kbound{k} edge).

The next step is the definition of the flip operation. Remember that this operation is local: given a \krel{k} edge~$e$ of a \ktri{k}~$T$, the flip of~$e$ in~$T$ consists in exchanging in~$T$ the edge~$e$ with the common bisector~$f$ of the two \kstar{k}s~$R$ and~$S$ of~$T$ containing~$e$. The previous function~\verb$completeStar$ computes the vertices of the \kstar{k}s~$R \eqdef \{r_0,\dots,r_{2k}\}$ and $S \eqdef \{s_0,\dots,r_{2k}\}$ in star order (and with $e=[r_{2k-1},r_{2k}]=[s_{2k},s_{2k-1}]$), and we just have to find their common bisector~$f$. One way to compute it would be to try, for any~$p,q\in\Z_{2k+1}$, whether $[r_p,s_q]$ bisects both~$R$ and~$S$, \ie whether~$r_{p+1}\cl s_q\cl r_{p-1}$ and~$s_{q+1}\cl r_p\cl s_{q-1}$. But according to Lemma~\ref{stars:lem:commonedge}, we even have to consider only the edges~$[r_p,s_p]$ (\ie the case~$p=q$).

\begin{lignes}
\begin{verbatim}
-- edgeFlip T e = flips the (k-relevant) edge e in the k-triangulation T
edgeFlip :: [[Int]] -> [Int] -> [[Int]]
edgeFlip T e = sort (f:(delete e T))
 where R = completeStar T e
			S = completeStar T (reverse e)
			cyclic :: [Int] -> Bool
			cyclic [a,b,c] = elem (sort (nub [a,b,c])) [[a,b,c], [b,c,a], [c,a,b]]
			f = head [ sort [R!!p, S!!p] | p <- [0..], 
													   cyclic [R!!(mod (p+1) (2*k+1)), S!!p, R!!(mod (p-1) (2*k+1))]
											     && cyclic [S!!(mod (p+1) (2*k+1)), R!!p, S!!(mod (p-1) (2*k+1))] ]
\end{verbatim}
\end{lignes}

Inside the code of a function, \haskell allows to locally define expressions and functions, either with the syntax~\verb$let$~--~\verb$in$ (expression style), or with the word~\verb$where$ (declaration style). For example, if~\verb$g :: a -> Int$, the three following declarations of~\verb$f :: a -> Int$ are equivalent:
\begin{center}
\verb$f x = (g x)^2 + g x	            f x = let y = g x in y^2 + y            f x = y^2 + y where y = g x$
\end{center}
The advantage of the two last formulations is that~\verb$g x$ is written and computed only once.

\enlargethispage{.2cm}
Concerning the code of the function~\verb$edgeFlip$, let us also mention that the function~\verb$reverse$ changes the orientation of the edge~\verb$e$: we compute one star on each side of~\verb$e$.

We are now ready to traverse the graph of flips to enumerate all \ktri{k}s of the \gon{n}. According to Lemma~\ref{stars:lem:increase}, we can start from the minimal \ktri{k}~$T_{n,k}^{\min}$ (called~\verb$initial$ in the following code) and flip only \krel{k} edges incident to one of the first~$k$ vertices of the~\gon{n}:

\begin{lignes}
\begin{verbatim}
-- kTriangulations = all k-triangulations of the n-gon
kTriangulations :: [[[Int]]]
kTriangulations = searchFlipGraph [initial] [ (initial, e) | e <- initial ]
 where initial = [ [i, mod (i+j) n] | i <- [0..k-1], j <- [k+1..n-k-1] ]
			searchFlipGraph :: [[[Int]]] -> [([[Int]],[Int])] -> [[[Int]]]
			searchFlipGraph L [] = L
			searchFlipGraph L (h:t)
				| elem T L   = searchFlipGraph L t
				| otherwise   = searchFlipGraph (T:L) (t ++ [ (T,e) | e <- T, head e < k ])
			 where T = edgeFlip (fst h) (snd h)
\end{verbatim}
\end{lignes}

Let us now simplify notations: we index the \krel{k} edges of~$E_n$ by numbers (their position in the list~\verb$kRelevantEdges$) and we transform any \ktri{k} in~\verb$kTriangulations$ into a list of indices:

\begin{lignes}
\begin{verbatim}
-- edgeIndex e = the index of the edge e in the list of k-relevant edges of the n-gon
edgeIndex :: [Int] -> Int
edgeIndex e = head (elemIndices e kRelevantEdges)

-- facets = all facets of Delta_nk
facets :: [[Int]]
facets = map (\T -> sort (map edgeIndex T)) kTriangulations
\end{verbatim}
\end{lignes}


We finish this section by checking that we obtain the expected number of \ktri{2}s of the octagon:

\begin{lignes}
\begin{verbatim}
> length facets
84
\end{verbatim}
\end{lignes}


\subsection{Necessary simplex orientations}

Let~$m \eqdef {n \choose 2}-kn$ denote the number of \krel{k} edges of~$E_n$, let~$I \eqdef \{0,\dots,m-1\}$, and let~$d \eqdef k(n-2k-1)$ be the number of \krel{k} edges in a \ktri{k} of the \gon{n} (which is the dimension of the polytopes we want to construct). We first fix these parameters in the code ($d$ is denoted~\verb$dim$ to avoid collision later):

\begin{lignes}
\begin{verbatim}
-- Fix the parameters m and dim
m, dim :: Int
m = div (n*(n-2*k-1)) 2; dim = k*(n-2*k-1)
\end{verbatim}
\end{lignes}

We now start computing all symmetric matroid realizations of~$\Delta_{n,k}$. Remember that each solution is a function~$\chi:I^{d+1}\to\{-1,0,1\}$ which satisfies the five properties of Lemma~\ref{ft:lem:chirotoperealization}:
\begin{enumerate}[(i)]
\item \defn{Alternating relations}:
for any \tuple{(d+1)}~$(i_0,\dots,i_d)\in I^{d+1}$ and any permutation~$\pi$ of~$\{0,\dots,d\}$ of signature~$\sigma$, 
$$\chi(i_{\pi(0)},\dots,i_{\pi(d)})=\sigma\chi(i_0,\dots,i_d).$$
\item \defn{Matroid property}: the set~$\chi^{-1}(\{-1,1\})$ of non-degenerate \tuple{(d+1)}s of~$I$ satisfies the \defn{Steinitz exchange axiom}: for any non-degenerate \tuple{(d+1)}s $X,Y\in\chi^{-1}(\{-1,1\})$ and any element~$x\in X\ssm Y$, there exists~$y\in Y\ssm X$ such that $X\diffsym\{x,y\}\in\chi^{-1}(\{-1,1\})$.
\item \defn{\GP relations}:
for any~$i_0,\dots,i_{d-2},j_1,j_2,j_3,j_4\in I$, the set 
\begin{eqnarray*}
\{ & \chi(i_0,\dots,i_{d-2},j_1,j_2)\,\chi(i_0,\dots,i_{d-2},j_3,j_4), \\
 & -\chi(i_0,\dots,i_{d-2},j_1,j_3)\,\chi(i_0,\dots,i_{d-2},j_2,j_4), \\
 & \chi(i_0,\dots,i_{d-2},j_1,j_4)\,\chi(i_0,\dots,i_{d-2},j_2,j_3) & \}
\end{eqnarray*}
either contains~$\{-1,1\}$ or is contained in~$\{0\}$.
\item \defn{Necessary Simplex Orientations}:
if~$i_0,\dots,i_d\in I$ are such that both~$\{i_0,\dots,i_{d-2},i_{d-1}\}$ and~$\{i_0,\dots,i_{d-2},i_{d}\}$ are facets of~$\Delta$, then for any~$j\in I\ssm\{i_0,\dots,i_d\}$,
$$\chi_P(i_0,\dots,i_{d-2},i_{d-1},i_d)=\chi_P(i_0,\dots,i_{d-2},i_{d-1},j)=\chi_P(i_0,\dots,i_{d-2},j,i_d).$$
\item \defn{Symmetry}:
there exists a morphism~$\tau: G\to \{\pm1\}$ such that for any~$i_0,\dots,i_d\in I$, and any~$g\in G$,
$$\chi_P(gi_0,\dots,gi_d)=\tau(g)\chi_P(i_0,\dots,i_d).$$
\end{enumerate}

\svs
We see such a function as a set of pairs~$(\tu,s)$, that we call \defn{oriented bases}, where~$\tu$ is a \tuple{(d+1)} of elements of~$I$, and~$s=\chi(\tu)\in\{-1,0,1\}$. To limit the number of oriented bases we are considering, we only remember those corresponding to ordered \tuple{(d+1)}s of~$I$, which we also call \defn{normalized oriented base}. Of course, we need a function to normalize an arbitrary oriented base:

\begin{lignes}
\begin{verbatim}
-- norm ob = normalizes the oriented base ob
norm :: ([Int],Int) -> ([Int],Int)
norm ([],s) = ([],s)
norm (h:t, s) = (take i tu ++ [h] ++ drop i tu, r*(-1)^(mod i 2))
 where (tu,r) = norm (t,s)
			i = length [ x | x <- tu, h > x ]
\end{verbatim}
\end{lignes}

Thus, a symmetric matroid realization is collection of ${m \choose d+1}$ normalized oriented bases which satisfies the four properties~(ii),~(iii),~(iv) and~(v). To construct such a collection of signs, we start with the necessary simplex orientations, and try to deduce from these initial informations the signs of all normalized oriented bases.

To obtain the necessary simplex orientations, we first need to compute all pairs~$\{\subf,\{x,y\}\}$ (where~$\subf$ is a subfacet of~$\Delta_{n,k}$ and~$x$ and~$y$ are two vertices disjoint from~$\subf$) such that both $\subf\cup\{x\}$ and~$\subf\cup\{y\}$ are facets of~$\Delta_{n,k}$. We name such pairs \defn{subfacet bounds}.

\begin{lignes}
\begin{verbatim}
-- subfacets = all subfacets of Delta_nk
subfacets :: [[Int]]
subfacets = nub (concat (map (\f -> [ delete j f | j <- f ]) facets))

-- subfacetBounds = all subfacet bounds of Delta_nk
subfacetBounds :: [([Int],[Int])]
subfacetBounds = [ (sf, [x,y]) | sf <- subfacets, [x,y] <- tuples 2 ([0..m-1]\\sf),
												 elem (sort (x:sf)) facets, elem (sort (y:sf)) facets ]
\end{verbatim}
\end{lignes}

Now that we have all subfacet bounds, we arbitrarily fix the orientation of an arbitrary subfacet bound (say the first subfacet bound in the list~\verb$subfacetBounds$ gets oriented positively), and deduce from this initial oriented base all necessary simplex orientations. For this, we have a function~\verb$fan$ which deduces, from the orientation of a subfacet bound, all implied orientations (see Property~(iv) of all symmetric matroid realizations):

\begin{lignes}
\begin{verbatim}
-- fan sfbd s = all orientations implied by the orientation of a subfacet bound sfbd
fan :: ([Int],[Int]) -> Int -> [([Int],Int)]
fan (sf, [x,y]) s = [ norm (sf ++ [x,z], s) | z <- vertices ] ++ [ norm (sf ++ [z,y], s) | z <- vertices ]
 where vertices = [0..m-1]\\(sf ++ [x,y])

-- necessarySimplexOrientations = necessary simplex orientations of Delta_nk
necessarySimplexOrientations :: [([Int],Int)]
necessarySimplexOrientations = generateNSO subfacetBounds [norm (sf0 ++ p0, 1)]
 where (sf0,p0) = head subfacetBounds
			generateNSO :: [([Int],[Int])] -> [([Int],Int)] -> [([Int],Int)]
			generateNSO [] OM = OM
			generateNSO (sfbd@(sf,p):t) OM
				| elem (norm (sf ++ p, 1)) OM  = generateNSO t (nub (OM ++ fan sfbd 1))
				| elem (norm (sf ++ p, -1)) OM  = generateNSO t (nub (OM ++ fan sfbd (-1)))
				| otherwise  = generateNSO (t ++ [sfbd]) OM
\end{verbatim}
\end{lignes}

By construction, at the end of this process, we have all necessary simplex orientations (and in particular, all subfacet bounds have been oriented). It is interesting to observe that we already know~$304$ orientations of the~${12 \choose 7} = 792$ normalized oriented bases:

\begin{lignes}
\begin{verbatim}
> length necessarySimplexOrientations
304
\end{verbatim}
\end{lignes}

In order to get more information, we need to use properties~(iii) and~(v) of the symmetric matroid realizations of~$\Delta_{n,k}$.


\subsection{\GP relations}

\index{Grassman@\GP relations}
We have now an embryo of oriented matroid, that we want to complete into a full oriented matroid (with all signs of the~${m \choose d+1}$ normalized oriented bases). For this, we will use extensively \GP relations.

First of all, we need to manipulate signs of our oriented bases. Given a \tuple{(d+1)} of~$I$, the following function tells us its sign ($-1$, $0$, or $1$) in our embryo of oriented matroid. For this, it normalizes this tuple and looks for its normalized version in our set of normalized oriented bases. If our embryo has still no orientation for this tuple, the function will answer~$\pm2$.

\begin{lignes}
\begin{verbatim}
-- sign tu OM = what is the sign of the tuple tu in the oriented matroid OM
sign :: [Int] -> [([Int],Int)] -> Int
sign tu OM = ns*signNormalized ntu OM
 where (ntu, ns) = norm (tu, 1)
			signNormalized :: [Int] -> [([Int],Int)] -> Int
			signNormalized _ [] = 2
			signNormalized tu1 ((tu2,s):t)	
				| tu1 == tu2  = s
				| otherwise  = sign tu1 t
\end{verbatim}
\end{lignes}

Concerning \GP relations, we need two different functions: given a set~\verb$OM$ of oriented bases, we need:
\begin{enumerate}[(i)]
\item a \defn{verification} function~\verb$contrGP$ which tells whether the set of oriented bases~\verb$OM$ contradicts \GP relations;~and
\item a \defn{deduction} function~\verb$consGP$ which computes all consequences of \GP relations from the signs in~\verb$OM$.
\end{enumerate}

\svs
For the verificiation function, we check that for any $w,x,y,z\in I$, and for any hyperline~$\hl$ disjoint from~$\{w,x,y,z\}$, the set
$$\{\chi(\hl,w,x)\chi(\hl,y,z),-\chi(\hl,w,y)\chi(\hl,x,z),\chi(\hl,w,z)\chi(\hl,x,y)\}$$
either contains~$\{-1,1\}$ or is contained in~$\{0\}$:

\begin{lignes}
\begin{verbatim}
-- contradiction hl [w,x,y,z] OM = does (hl,[w,x,y,z]) contradict Grassmann-Plücker relations?
contradiction :: [Int] -> [Int] -> [([Int],Int)] -> Bool
contradiction hl [w,x,y,z] OM = elem (sort (nub ([a*b, -c*d, e*f]))) [[1], [-1], [0,1], [-1,0]]
 where a = sign (hl ++ [w,x]) OM; b = sign (hl ++ [y,z]) OM
		 	c = sign (hl ++ [w,y]) OM; d = sign (hl ++ [x,z]) OM
			e = sign (hl ++ [w,z]) OM; f = sign (hl ++ [x,y]) OM

-- contrGP OM = does OM contradict Grassmann-Plücker relations?
contrGP :: [([Int],Int)] -> Bool
contrGP OM = or [ contradiction hl q OM | q <- tuples 4 [0..m-1], hl <- tuples (dim-1) ([0..m-1]\\q) ]
\end{verbatim}
\end{lignes}

The deduction function, is very similar. For example, given a \tuple{4}~$(w,x,y,z)\in I^4$ and a hyperline~$\hl$ disjoint form~$\{w,x,y,z\}$:
\begin{enumerate}[(i)]
\item if $-\chi(\hl,w,y)\chi(\hl,x,z)=\chi(\hl,w,z)\chi(\hl,x,y)=0$ and~$\chi(\hl,y,z)\ne 0$, then~$\chi(\hl,w,x)$ necessarily equals $0$;
\item if $\{-\chi(\hl,w,y)\chi(\hl,x,z),\chi(\hl,w,z)\chi(\hl,x,y)\}\cap\{-1,1\}=\{\varepsilon\}$, with~$\varepsilon\in\{-1,1\}$, then~$\chi(\hl,w,x)$ necessarily equals $-\varepsilon\chi(\hl,y,z)$;
\end{enumerate}
Similarly, we have equations relating~$\chi(\hl,w,y)$ with~$\chi(\hl,x,z)$ and~$\chi(\hl,w,z)$ with~$\chi(\hl,x,y)$. These equations can be translated into \haskell as follows:

\begin{lignes}
\begin{verbatim}
-- consequence hl [w,x,y,z] OM = consequence of the Grassmann-Plücker relations
consequence :: [Int] -> [Int] -> [([Int],Int)] -> [([Int],Int)]
consequence hl [w,x,y,z] OM
	| abs a == 2 && abs b == 1 && elem [-c*d, e*f] set  = [norm (hl ++ [w,x], b*signum (c*d-e*f))]
	| abs b == 2 && abs a == 1 && elem [-c*d, e*f] set  = [norm (hl ++ [y,z], a*signum (c*d-e*f))]
	| abs c == 2 && abs d == 1 && elem [a*b, e*f] set  = [norm (hl ++ [w,y], d*signum (a*b+e*f))]
	| abs d == 2 && abs c == 1 && elem [a*b, e*f] set  = [norm (hl ++ [x,z], c*signum (a*b+e*f))]
	| abs e == 2 && abs f == 1 && elem [a*b, -c*d] set  = [norm (hl ++ [w,z], f*signum (-a*b+c*d))]
	| abs f == 2 && abs e == 1 && elem [a*b, -c*d] set  = [norm (hl ++ [x,y], e*signum (-a*b+c*d))]
	| otherwise  = []
 where a = sign (hl ++ [w,x]) OM; b = sign (hl ++ [y,z]) OM
		 	c = sign (hl ++ [w,y]) OM; d = sign (hl ++ [x,z]) OM
			e = sign (hl ++ [w,z]) OM; f = sign (hl ++ [x,y]) OM
			set = [[-1,-1], [0,0], [1,1], [0,-1], [-1,0], [0,1], [1,0]]
\end{verbatim}
\end{lignes}

In the previous function, the function~\verb$signum :: Int -> Int$ gives the sign ($-1$,~$0$, or~$1$) of an integer, and should not be confused with the function~\verb$sign$ (which gives the sign of a tuple in a set of oriented bases). 
The reader can also observe that for a given \tuple{4} and a given disjoint hyperline, we can have either no or one consequence of the \GP relations. This is the reason why we return a list of consequences (which can either be empty or have a single element).

We can finally write the function~\verb$consGP$ which iteratively deduces all consequences of \GP relations in a set~\verb$OM$ of oriented bases. At each step, it first computes all immediate consequences of \GP relations, applying the function~\verb$consequence$ to all \tuple{4}s and disjoint hyperlines. If we have no new consequence, we return~\verb$OM$. If we have new consequences, but if they are not coherent (\ie if two consequences contradict themself), then we return the empty set: no oriented matroid can be constructed from the set~\verb$OM$ of oriented bases. Finally, if we have new coherent consequences, then we add them to the set of oriented bases and we run~\verb$consGP$ again.

\begin{lignes}
\begin{verbatim}
-- coherent OM = is OM coherent?
coherent :: [([Int],Int)] -> Bool
coherent OM = and [ fst (sortedOM!!(i-1)) /= fst (sortedOM!!i) | i<-[1..(length sortedOM)-1] ]
 where sortedOM = sort OM
\end{verbatim}
\end{lignes}

\begin{lignes}
\begin{verbatim}
-- consGP OM = consequences of Grassmann-Plücker relations in OM
consGP :: [([Int],Int)] -> [[([Int],Int)]]
consGP OM
	| newSigns == []  = [OM]
	| coherent newSigns  = consGP (nub (OM ++ newSigns))
	| otherwise  = []
 where newSigns = nub (concat [ consequence hl q OM | q <- tuples 4 [0..m-1],
 																						   hl <- tuples (dim-1) ([0..m-1]\\q) ])
\end{verbatim}
\end{lignes}





\subsection{Symmetry under the dihedral group}

We now focus on the action of the dihedral group on our oriented matroids. We first define the action on the oriented matroid of the rotation~$\rho:v\mapsto v+1$ and the reflection~$\sigma:v\mapsto -v$:

\begin{lignes}
\begin{verbatim}
-- rotation i = the rotation v -> v + i
rotation :: Int -> [Int] -> [Int]
rotation i T = map (\e -> edgeIndex (sort (map (\v -> mod (v+i) n) (kRelevantEdges!!e)))) T

-- reflection = the reflexion v -> -v
reflection :: [Int] -> [Int]
reflection T = map (\e -> edgeIndex (sort (map (\v -> mod (-v) n) (kRelevantEdges!!e)))) T
\end{verbatim}
\end{lignes}

These transformations generate the complete dihedral group:
$$\D_n = \ens{\rho^i}{i\in[n]}\cup\ens{\rho^i\circ\sigma}{i\in[n]}.$$
In order to deduce from the action of the dihedral group the signs of new oriented bases, we need to know the signs of the action of~$\rho$ and~$\sigma$ (\ie the signs $\tau(\rho)$ and~$\tau(\sigma)$). Since the dihedral group already acts on the necessary simplex orientations, these signs are already determined:

\begin{lignes}
\begin{verbatim}
-- signRot = the sign of the action of the rotation v -> v + 1 on the oriented matroid
signRot :: Int
signRot = s*(sign (rotation 1 tu) necessarySimplexOrientations)
 where (tu,s) = head necessarySimplexOrientations

-- signRef = the sign of the action of the reflectino v -> -v on the oriented matroid
signRef :: Int
signRef = s*(sign (reflection tu) necessarySimplexOrientations)
 where (tu,s) = head necessarySimplexOrientations
\end{verbatim}
\end{lignes}

Now, given a set of oriented bases, we can compute the orbit of these signs under the action of the dihedral group:

\begin{lignes}
\begin{verbatim}
-- orbitDG OM = the orbit under the dihedral group of the set of oriented bases OM
orbitDG :: [([Int],Int)] -> [([Int],Int)]
orbitDG OM = sort (nub (concat (map (\(tu,s) -> ([ norm (rotation i tu, (signRot^i)*s) | i <- [0..n-1] ]
						++ [ norm (rotation i (reflection tu), (signRot ^i)*signRef*s) | i <- [0..n-1] ])) OM)))
\end{verbatim}
\end{lignes}

Applying the function~\verb$orbitDG$ to a set of oriented bases yields a set of oriented bases closed under the action of~$\D_n$. It is convenient to present such a set of oriented bases keeping only one representative per orbit under~$\D_n$:

\begin{lignes}
\begin{verbatim}
-- reduceOrbitDG OM = one representative by orbit under the dihedral group
reduceOrbitDG :: [([Int],Int)] -> [([Int],Int)]
reduceOrbitDG [] = []
reduceOrbitDG (h:t) = h:(reduceOrbitDG (t\\orbit))
 where orbit = orbitDG [h]
\end{verbatim}
\end{lignes}

Using this function, we can count the number of orbits in our necessary simplex orientations (observe that we are sure that they are closed under the action of~$\D_n$ since they were deduced from a symmetric set of facets):

\begin{lignes}
\begin{verbatim}
> length (reduceOrbitDG necessarySimplexOrientations)
25
\end{verbatim}
\end{lignes}



To finish this section on the dihedral group, we want to observe that imposing a realization of a simplicial complex to be symmetric sometimes forces the determinant of certain tuples to vanish. Indeed, if~$\chi:I^{d+1}\to\{-1,0,1\}$ is a symmetric matroid realization of a simplicial complex~$\Delta$ under a group~$G$, then for any~$i_0,\dots,i_d\in I$, and any $g\in G$,
$$\chi(i_0,\dots,i_d)=\tau(g)\chi(gi_0,\dots,gi_d)=\tau(g)\chi(i_{\pi(0)},\dots,i_{\pi(d)})= \varepsilon\tau(g)\chi(i_0,\dots,i_d),$$
where~$\pi$ is the permutation of~$\{0,\dots,d\}$ defined by~$i_{\pi(j)}=gi_j$, and~$\varepsilon$ is the signature of~$\pi$. In particular, if~$\varepsilon=-\tau(g)$, then it forces~$\chi(i_0,\dots,i_d)$ to vanish. We say that~$\chi(i_0,\dots,i_d)$ is a \defn{necessary zero}. These necessary zeros are easy to compute: they are the \tuple{(d+1)}s of~$I$ whose orbit under the dihedral group is not coherent:

\begin{lignes}
\begin{verbatim}
-- necessaryZeros = tuples whose determinant is forced to vanish by the dihedral action
necessaryZeros :: [([Int],Int)]
necessaryZeros = [ (b,0) | b <- tuples (dim+1) [0..m-1], not (coherent (orbitDG [(b,1)])) ]
\end{verbatim}
\end{lignes}

Observe that when we add these necessary zeros to our necessary simplex orientations, we obtain new \GP consequences:

\begin{lignes}
\begin{verbatim}
> length necessaryZeros
16
> length (consGP (sort (necessarySimplexOrientations ++ necessaryZeros)))
456
\end{verbatim}
\end{lignes}

Thus, we already now~$456$ normalized oriented bases (of the~$792$ in total). In order to complete this embryo of oriented matroid, we need to start ``guessing signs''.


\subsection{Guessing signs}

We now have to guess new signs to complete our current set of oriented bases into a full oriented matroid. The method consists in picking one \tuple{(d+1)} of~$I$ whose sign remains unknown, and to try the three possibilities~$-1$,~$0$ and~$1$. Applying symmetry and \GP relations, we may obtain:
\begin{enumerate}[(1)]
\item either a contradiction, and we eliminate this choice of sign for this \tuple{(d+1)};
\item or a complete symmetric matroid realization of~$\Delta_{8,2}$;
\item or a certain number of additional signs, but not all, and we have to iterate the guessing process until we arrive at situation (1) or (2).
\end{enumerate}
The point of this section is to make a good choice for the \tuple{(d+1)} for which we would like to guess the sign, that is, a choice which will provide as much information as possible.

\svs
Let~$\hl$ be an hyperline and~$w,x,y,z$ be four indices not in~$\hl$. Consider the six \tuple{(d+1)}s $\hl\cup\{w,x\}$, $\hl\cup\{y,z\}$, $\hl\cup\{w,y\}$, $\hl\cup\{x,z\}$, $\hl\cup\{w,z\}$, and $\hl\cup\{x,y\}$ which appear in the \GP relation~$\{\hl,[w,x,y,z]\}$. In our current set of oriented bases, the signs of some of these six \tuple{(d+1)}s are known, while the others remain unknown. Assume for example that both $\chi(\hl,w,x)\chi(\hl,y,z)$ and $-\chi(\hl,w,y)\chi(\hl,x,z)$ are known and equal, while both $\chi(\hl,w,z)$ and $\chi(\hl,x,y)$ remain unknown. Then finding $\chi(\hl,w,z)$ would automatically provide $\chi(\hl,x,y)$ and \viceversa. We call \defn{key base} (of~$\{\hl,[w,x,y,z]\}$) any such base, \ie any base whose sign, if known, would allow us to derive the sign of another base. The function~\verb$keyBases$ lists all key bases of a given \GP relation:

\begin{lignes}
\begin{verbatim}
-- keyBases OM GPrelation = key bases of the Grassmann-Plücker relation
keyBases :: [([Int],Int)] -> ([Int],[Int]) -> [[([Int],Int)]]
keyBases OM (hl,q) = [ [(b,s)] | b <- map (\p -> sort (hl ++ p)) (tuples 2 q), s <- [-1,1],
												  sign b OM == 2, consequence hl q ((b,s):OM) /= [] ]
\end{verbatim}
\end{lignes}

A good choice for the \tuple{(d+1)} for which we would like to guess the sign should be, as well as its images under the dihedral group, a key base of many \GP relations. The following function~\verb$bestTuple$ computes our best possible choice: it enumerates the key bases of all \GP relations, and chooses the \tuple{(d+1)} whose orbit under~$\D_n$ is the most represented among all these key bases. (The routine~\verb$mostRep$ chooses in a list the most represented element: for example,~\verb$mostRep [1,4,2,1,2,2]$ returns~\verb$[2]$.)

\begin{lignes}
\begin{verbatim}
-- mostRep L = the most represented element in L
mostRep :: Ord a => [a] -> [a]
mostRep L
	| L == []  = []
	| otherwise  = [head (head (sortBy (\x y -> compare (length y) (length x)) (group (sort L))))]

-- bestTuple OM = the best choice for the tuple for which we would like to guess the sign
bestTuple :: [([Int],Int)] -> [Int]
bestTuple OM
	| b == []  = head ((tuples (dim+1) [0..m-1])\\(map fst OM))
	| otherwise  = fst (head (head b))
 where GPrels = [ (hl,q) | hl <- tuples (dim-1) [0..m-1], q <- tuples 4 ([0..m-1]\\hl) ]
            b = mostRep (map orbitDG (concat (map (keyBases  OM) GPrels)))
\end{verbatim}
\end{lignes}

With this rule for choosing the best possible tuple, we can now write the~\verb$run$ function, which starts from a set of oriented bases and iteratively guess new signs, until it founds a full oriented matroid:

\begin{lignes}
\begin{verbatim}
-- run = the iterative guessing process
run :: [[([Int],Int)]] -> [[([Int],Int)]]
run [] = []
run (OM:t)
	| length OM == length (tuples (dim+1) [0..m-1])  = OM:(run t)
	| otherwise  = run ([ [newOM] | newOM <- newOMs, not (contrGP newOM) ] ++ t)
 where b = bestTuple OM
            newOMs = concat [ consGP ((orbitDG [(b,s)]) ++ OM) | s <- [1,0,-1] ]
\end{verbatim}
\end{lignes}


\subsection{Results}

We finally obtain our symmetric realizations of~$\Delta_{8,2}$:

\begin{lignes}
\begin{verbatim}
-- symmetricMatroidRealizations = the list of symmetric matroid realizations of Delta_nk
symmetricMatroidRealizations :: [[([Int],Int)]]
symmetricMatroidRealizations = run [necessarySimplexOrientations ++ necessaryZeros]

-- reducedSMR = the list of symmetric matroid realizations, with one representative per orbit 
reducedSMR :: [[([Int],Int)]]
reducedSMR = map reduceOrbitDG symmetricMatroidRealizations
\end{verbatim}
\end{lignes}

We obtain Proposition~\ref{ft:prop:delta82matroid} with the following command:

\begin{lignes}
\begin{verbatim}
> length symmetricMatroidRealizations
15
\end{verbatim}
\end{lignes}

In order to describe these~$15$ symmetric matroid realizations, we reduce each of them, keeping only one representative for each of the~$62$ orbits under the dihedral group: we obtain the list~\verb$reducedSMR$. These~$15$ realizations share the following~$59$ orbits:

\begin{lignes}
\begin{verbatim}
> simplify (sort (interesection reducedSMR))
(abcdefg, 0), (abcdefI, 1), (abcdefJ, -1), (abcdefK, 1), (abcdegI, -1), (abcdegJ, 1), (abcdeIK, 1), 
(abcdeIL, -1), (abcdeJK, -1), (abcdfgI, 1), (abcdfgJ, -1), (abcdfgL, -1), (abcdfIJ, -1), (abcdfIK, -1), 
(abcdfIL, 1), (abcdfJK, 1), (abcdfJL, -1), (abcdfKL, 1), (abcdIJK, 1), (abcdIJL, -1), (abcdIKL, 1), 
(abcefgI, 0), (abcefgK, 0), (abcefIJ, 1), (abcefIK, -1), (abcefIL, -1), (abcefJK, 1), (abcefJL, 1), 
(abcefKL, -1), (abcegIJ, -1), (abcegIK, 1), (abcegIL, 1), (abcegKL, 1), (abceIJL, 1), (abceIKL, -1), 
(abceJKL, 1), (abcfIJK, 1), (abcfIKL, 1), (abcIJKL, -1), (abdegIJ, -1), (abdegIK, -1), (abdegIL, -1), 
(abdegJK, 1), (abdeIJK, -1), (abdeIJL, 1), (abdfIJL, -1), (abdfIKL, -1), (abdfJKL, 1), (abdgIJK, 1), 
(abdgIJL, 1), (abdgJKL, -1), (abdIJKL, 1), (abefIJK, 1), (abefIJL, 1), (abefJKL, -1), (abeIJKL, -1), 
(acegIJK, 1), (aceIJKL, 1), (acfIJKL, -1)
\end{verbatim}
\end{lignes}

Here, the function~\verb$intersection$ intersects all realizations, and the function \verb$simplify$ transforms the tuple~\verb$[0,1,2,3,4,5,6]$ into the word~\verb$abcdefg$. Finally, the three remaining orbits are those of~\verb$abcdeIJ$,~\verb$abceIJK$ and~\verb$abdfIJK$. The following command list their respective signs in the~$15$ realizations:

\begin{lignes}
\begin{verbatim}
> map (map snd) (map (\x -> sort (x\\(intersection reducedSMR))) reducedSMR)
[[1,-1,1], [1,-1,0], [1,-1,-1], [0,-1,1], [0,-1,0], [0,-1,-1], [-1,1,1], [-1,1,0], [-1,1,-1], [-1,0,1], [-1,0,0],
[-1,0,-1], [-1,-1,1], [-1,-1,0], [-1,-1,-1]]
\end{verbatim}
\end{lignes}

We finish by a short remark. To obtain all possible symmetric matroid realizations of~$\Delta_{8,2}$, we started from the necessary simplex orientations (and the necessary zeros), and we iteratively guessed signs, applying at each step symmetry and \GP relations. Thus, by construction, the sets of oriented bases we obtain satisfy points (i), (iii), (iv) and (v) of the definition of symmetric matroid realization. To conclude, we still have to check point (ii), namely that the set of non-degenerate \tuple{(d+1)}s forms a matroid. This is easily done by the following function:

\begin{lignes}
\begin{verbatim}
-- matroid M = checks that the set of (sorted) bases M forms a matroid
matroid :: [[Int]] -> Bool
matroid M = and [ or [ elem (sort (y:(delete x X))) M | y <- Y\\X ] | (X,Y) <- tuples 2 M, x <- X\\Y ]
\end{verbatim}
\end{lignes}

\begin{lignes}
\begin{verbatim}
> and [ map (\OM -> matroid [ tu | (tu,s) <- OM, s /= 0 ]) symmetricMatroidRealizations ]
True
\end{verbatim}
\end{lignes}

\end{appendices}


\renewcommand{\partfigure}{squareddominopoly}
\renewcommand{\namepart}{Polytopality of products}
\part{Polytopality of products}
	\chapter{Introduction}\label{chap:introduction_polytopes}

The second part of this dissertation focusses on \defn{polytopality}\index{polytopality} questions: given an abstract polytopal \complex{k}, we ask \defn{``whether it is the \skeleton{k} of a polytope?''}, \defn{``what dimension can a realizing polytope have?''} and \defn{``how does a realizing polytope look like?''}. This kind of questions can be studied for complexes:
\begin{enumerate}[(i)]
\item either derived from combinatorial structures: for example, the simplicial complex of all non-crossing diagonals of a convex polygon, or the graph of flips on multitriangulations of a point set (remember our discussion in Section~\ref{ft:sec:multiassociahedron}), or the graph whose vertices are permutations of~$[n]$ and where two permutations are related by an edge if they differ by a single adjacent transposition, \etc{} See~\cite[Lecture~9]{z-lp-95} for more examples.
\item or obtained by certain operations on other complexes: for example, the next two chapters deal with complexes obtained as products of other complexes.
\end{enumerate}

The most prominent result on realizability problems is Steinitz' Theorem (see Theorem~\ref{nonpolytopal:theo:steinitz}) which characterizes graphs of \poly{3}topes and describes their realization space. No similar statement can be expected in higher dimensions, where realizability is a difficult question in general. In fact, even \poly{4}topes can have arbitrarily wild and complicated realization spaces (this is known as the \defn{Universality Theorems} for polytopes~\cite{rg-rsp-96}).

In the second part of this thesis, we consider a specific realizability question: it concerns \complex{k}es obtained as products of smaller complexes. More precisely:
\begin{itemize}
\item In Chapter~\ref{chap:nonpolytopal}, we explore polytopality questions on \defn{Cartesian products of graphs}. This product is defined in such a way that the graph of a product of polytopes is the product of their graphs; in particular, products of polytopal graphs are automatically polytopal. Our question addresses the reciprocal statement: are the two factors of a polytopal product necessarily polytopal?
\item In Chapter~\ref{chap:psn}, we study \defn{prodsimplicial-neighborly polytopes}, whose \skeleton{k} is combinatorially equivalent to that of a given product of simplices. The goal is to construct \psn polytopes which minimize the possible dimension among all \psn polytopes.
\end{itemize}

Before summarizing our results in Section~\ref{intro:sec:results}, we fix in the next section some notation and terminology. Since we will use them extensively in the coming chapters, we also recall the definition and elementary properties of the famous cyclic polytopes as an illustration of dimensional ambiguous polytopes. The reader familiar with polytope theory is invited to jump directly to Section~\ref{intro:sec:results}.


\section{Generalities on polytopes}\label{intro:sec:generalities}


\subsection{Notation and terminology}\label{intro:subsec:generalities:notations}

A \defn{polytope}\index{polytope|hbf} is the convex hull of a finite point set of~$\R^n$, or equivalently a bounded intersection of finitely many half-spaces in~$\R^n$. We write~$P \eqdef \conv(X) \eqdef \ens{x\in\R^n}{Ax\le b}$, where~$X\subset\R^n$ is finite, $A\in\R^{n\times f}$ and~$b\in\R^f$. The \defn{dimension} of a polytope is the dimension of its affine span (we abbreviate ``\dimensional{d} polytope'' into ``\poly{d}tope'').

Dating back to antiquity, the study of polytopes was initially limited to dimensions~$2$ and~$3$, and focussed on their (sym)metric properties: famous examples are regular polygons and Platonic solids. During the last centuries, the study of polytopes was extended to higher dimensions, and the interest on polytopes shifted from their geometric properties to more combinatorial aspects. We now pay particular attention to their face structure.

A \defn{face}\index{face} of a polytope~$P$ is an intersection~$F \eqdef P\cap H$ of~$P$ with a supporting hyperplane~$H$ of~$P$ (the empty set is also considered as a face of~$P$). If $H=\ens{x\in\R^n}{\dotprod{h}{x}=c}$ and ${\dotprod{h}{x}<c}$ for all~$x\in P$, then we say that~$h$ is an (outer) \defn{normal vector}\index{normal vector} of the face~$F$. A \face{0} (resp.~\face{1}, resp.~\face{(d-1)}) of a \poly{d}tope is called a \defn{vertex} (resp.~an \defn{edge}, resp.~a \defn{facet}\index{facet}). The \defn{\fvector{}} of a \poly{d}tope $P$ is the vector $f(P) \eqdef (f_{-1}(P),f_0(P),\dots,f_{d-1}(P))$, where $f_k(P)$ denotes the number of \face{k}s of $P$ (by convention, the empty face has dimension~$-1$). For~${k\in\N}$, the \defn{\skeleton{k}}\index{skeleton@\skeleton{k}} of a polytope is the polytopal \complex{k} formed by all its faces of dimension at most~$k$. In particular, the \skeleton{1} of a polytope is its \defn{graph} and the \skeleton{(d-1)} of a \poly{d}tope is its \defn{boundary complex}\index{boundary complex|hbf}.

As mentioned previously, we are interested in combinatorial properties of polytopes: two polytopes are considered as equivalent if the incidence relations between their faces behave similarly in both polytopes. More precisely, two complexes~$\cC$ and~$\cC'$ are \defn{combinatorially equivalent}\index{combinatorially equivalent} if there is an inclusion-preserving bijection from one to the other: $\phi:\cC\to\cC'$ with $A\subset B\Leftrightarrow \phi(A)\subset\phi(B)$.

\begin{example}
A \simp{d} is the convex hull of~$d+1$ affinely independent points: for example, $\simplex_d \eqdef \conv\ens{e_k}{k\in[d+1]}$, where~$(e_1,\dots,e_{d+1})$ denotes the canonical orthogonal basis of~$\R^{d+1}$. The convex hull of any subset of its vertices forms a face of the simplex. In particular, its graph is the complete graph, its \fvector{} is given by~$f_k(\simplex_d)={d+1\choose k+1}$, and all simplices are combinatorially equivalent.
\end{example}

To finish, let us recall that a polytope is \defn{simplicial}\index{polytope!simplicial ---} if all its facets are simplices, and \defn{simple}\index{polytope!simple ---} if all its vertex figures are simplices (the vertex figure of a vertex~$v$ in a polytope~$P$ is the intersection of~$P$ with an hyperplane that only cuts off vertex~$v$). In other words, the vertices of a simplicial polytope are in general position (no~$d+1$ of them lie in an hyperplane) while the facet-defining inequalities of a simple polytope are in general position. For example, a simplex is both simple and simplicial.

All along the text of the next two chapters, we will recall classical results of polytope theory when we need them. For a detailed presentation on polytopes, we refer to the excellent books of Branko Gr\"unbaum~\cite{g-cp-03} and G\"unter Ziegler~\cite{z-lp-95}, as well as the expository chapters~\cite[Chapter~5]{m-ldg-02} and~\cite{hrgz-bpcp-97,bb-fnpc-97,k-psp-97}.


\subsection{Polytopality and ambiguity}\label{intro:subsec:generalities:ambiguity}

The main property we will be interested in (in particular for graphs) is the following:

\begin{definition}
\index{polytopal}
A \complex{k}~$\cC$ is \defn{polytopal} if it is combinatorially equivalent to the \skeleton{k} of some polytope~$P$. If~$P$ is \dimensional{d}, we say that~$\cC$ is \defn{\poly{d}topal}.
\end{definition}

As soon as a \complex{k} is polytopal, we want to know which (combinatorial types of) polytopes realize it. Certain complexes determine the combinatorial type of their possible realizations (for example, graphs of \poly{3}topes~---~see Theorem~\ref{nonpolytopal:theo:steinitz}); some others are far from that: even the dimension of the realization is sometimes not fixed by the \complex{k}.

\begin{definition}
\index{polytopality!--- range}
\index{polytopality!--- dimension}
A \complex{k} is \defn{dimensionally ambiguous} if there exist two polytopes realizing it in two different dimensions. We define the \defn{polytopality range} of a \complex{k}~$\cC$ to be the set of integers~$d$ for which~$\cC$ is \poly{d}topal, and the \defn{polytopality dimension} to be the minimal dimension of a realizing polytope.
\end{definition}

As an illustration, we recall the definition and basic properties of cyclic polytopes, which are fundamental examples in polytope theory because of their extremal properties:

\begin{definition}\label{intro:def:cyclic}
\index{polytope!cyclic ---|hbf}
Let~$t\mapsto \mu_d(t) \eqdef (t,t^2,\dots,t^d)^T$ denote the \defn{moment curve} in~$\mathbb{R}^d$, and choose $n$ arbitrary distinct real numbers $t_1<t_2<\dots<t_n$. We denote by~$C_d(n) \eqdef \conv\ens{\mu_d(t_i)}{i\in[n]}$ ``the'' \defn{cyclic polytope} of dimension~$d$ with~$n$ vertices (the convex hull of any~$n$ distinct points on the moment curve always leads to the same combinatorial polytope).
\end{definition}

\begin{figure}[h]
	\capstart
	\centerline{\includegraphics[scale=1]{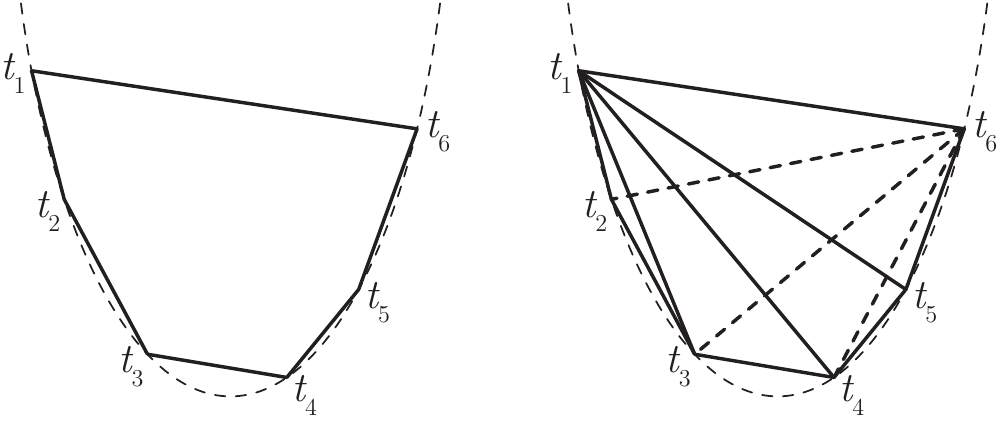}}
	\caption[Cyclic polytopes]{The cyclic polytopes $C_2(6)$ (with a unique upper facet $\{1,6\}$) and $C_3(6)$ (with its upper facets $\{i,i+1,6\}$, $i\in[4]$).}
	\label{intro:fig:cyclicPolytope}
\end{figure}

\begin{proposition}\label{intro:prop:cyclic}
\index{polytope!neighborly ---}\index{neighborly polytope}
The cyclic polytope~$C_d(n)$ is \defn{\neighborly{\Fracfloor{d}{2}}}: any subset of at most $\Fracfloor{d}{2}$ vertices forms a face of~$C_d(n)$.
\end{proposition}

In particular, for any~$4\le d< n$, the graph of the cyclic polytope~$C_d(n)$ is the complete graph. In other words, for~$n\ge 5$, the polytopality range of the complete graph~$K_n$ contains, and in fact equals~$\{4,\dots,n-1\}$.

The proof of Proposition~\ref{intro:prop:cyclic} is classic but we recall it since we will use similar ideas in Section~\ref{psn:sec:cyclic}. We naturally identify a face $\conv\ens{\mu(t_i)}{i\in F}$ of~$C_d(n)$ with the subset~$F$ of~$[n]$. Gale's evenness criterion (Proposition~\ref{intro:prop:gale}) characterizes which subsets of~$[n]$ are facets of~$C_d(n)$. To state it, define the \defn{blocks} of a subset~$F$ of~$[n]$ to be the maximal subsets of consecutive elements of~$F$. The \defn{initial} (resp.~\defn{final}) block is the block containing~$1$ (resp.~$n$)~---~when it exists. Other blocks are called \defn{inner} blocks. For example, $\{1,2,3,6,7,9,10,11\}\subset [11]$ has~$3$ blocks: $\{1,2,3\}$ (initial), $\{6,7\}$ (inner) and $\{9,10,11\}$ (final).

\begin{proposition}[Gale's evenness criterion~\cite{g-ncp-63}]\label{intro:prop:gale}
A subset~$F$ of~$[n]$ is a facet of~$C_d(n)$ if and only if~$|F|=d$ and all inner blocks of~$F$ have even size. Furthermore, such a facet~$F$ is supported by the hyperplane $H_F \eqdef \ens{z\in\R^d}{\dotprod{(\gamma_i(F))_{i\in[d]}}{z}=-\gamma_0(F)}$, where~$\gamma_0(F),\dots,\gamma_d(F)$ are defined as the coefficients of~$t^0,\dots,t^d$ in the polynomial:
$$\Pi_F(t) \eqdef \prod_{i\in F}(t-t_i) \eqdef \sum_{i=0}^d \gamma_i(F)t^i.$$
Finally, according to the parity of the size~$\ell$ of the final block of~$F$:
\begin{enumerate}[(i)]
\item If~$\ell$ is odd, then the whole cyclic polytope lies below~$H_F$ (with respect to the last coordinate). We say that~$F$ is an \defn{upper facet}\index{facet!upper and lower ---s}. Its normal vector is~$(\gamma_i(F))_{i\in [d]}$.
\item If~$\ell$ is even, then the whole cyclic polytope lies above~$H_F$. We say that~$F$ is a \defn{lower facet}.  Its normal vector is~$(-\gamma_i(F))_{i\in [d]}$.
\end{enumerate}
\end{proposition}

\begin{proof}
Observe first that:
\begin{enumerate}[(i)]
\item Vandermonde determinant
$$\det\begin{pmatrix}1 & 1 & \dots & 1 \\ \mu(x_0) & \mu(x_1) & \dots & \mu(x_d) \end{pmatrix}=\prod_{0\le i<j\le d} (x_j-x_i)$$
ensures that any~$d+1$ points on the \dimensional{d} moment curve are affinely independent, and thus, that the cyclic polytope is simplicial.
\item For any~$t\in\R$,
$$\dotprod{(\gamma_i(F))_{i\in[d]}}{\mu_d(t)}+\gamma_0(F)=\sum_{i=0}^d \gamma_i(F)t^i=\Pi_F(t).$$
\item the coefficient~$\gamma_d(F)$ is always~$1$, and thus, the vector~$(\gamma_i(F))_{i\in[d]}$ points upwards (with respect to the last coordinate).
\end{enumerate}

Let~$F$ be a subset of~$[n]$ of size~$d$. Then:
\begin{enumerate}[(i)]
\item for all~$j\in F$,~$\Pi_F(t_j)=0$; thus,~$H_F$ is the affine hyperplane spanned by~$F$;
\item for all~$j\notin F$, the sign of~$\Pi_F(t_j)$ is~$(-1)^{|F\cap\{j+1,\dots,n\}|}$.
\end{enumerate}
In particular, if~$F$ has an odd inner block~$\{a,a+1,\dots,b\}$, then~$\Pi_F(t_{a-1})$ and~$\Pi_F(t_{b+1})$ have different signs, and~$F$ is not a facet. Reciprocally, if all inner blocks have even size, then the sign of all~$\Pi_F(t_j)$ is~$(-1)^{\ell}$, where~$\ell$ is the size of the final block. Thus,~$F$ is an upper facet when~$\ell$ is odd, and a lower facet when~$\ell$ is even.
\end{proof}


\section{Polytopality of products: our results}\label{intro:sec:results}

The \defn{Cartesian product}\index{Cartesian product!--- of polytopes|hbf} of two polytopes~$P$, $Q$ is the polytope~$P\times Q \eqdef \ens{(p,q)}{p\in P, q\in Q}$. The combinatorial structure of the product is completely understood from the structure of its factors:
\begin{enumerate}[(i)]
\item The dimension of~$P\times Q$ is the sum of the dimensions of~$P$ and~$Q$.
\item The non-empty faces of~$P\times Q$ are precisely the products of non-empty faces of~$P$ by non-empty faces of~$Q$.
\item The \fvector{} of~$P\times Q$ is given by~$f_k(P\times Q)=\sum_{i+j=k} f_i(P)f_j(Q)$.
\end{enumerate}

\begin{figure}[h]
	\capstart
	\centerline{\includegraphics[scale=.8]{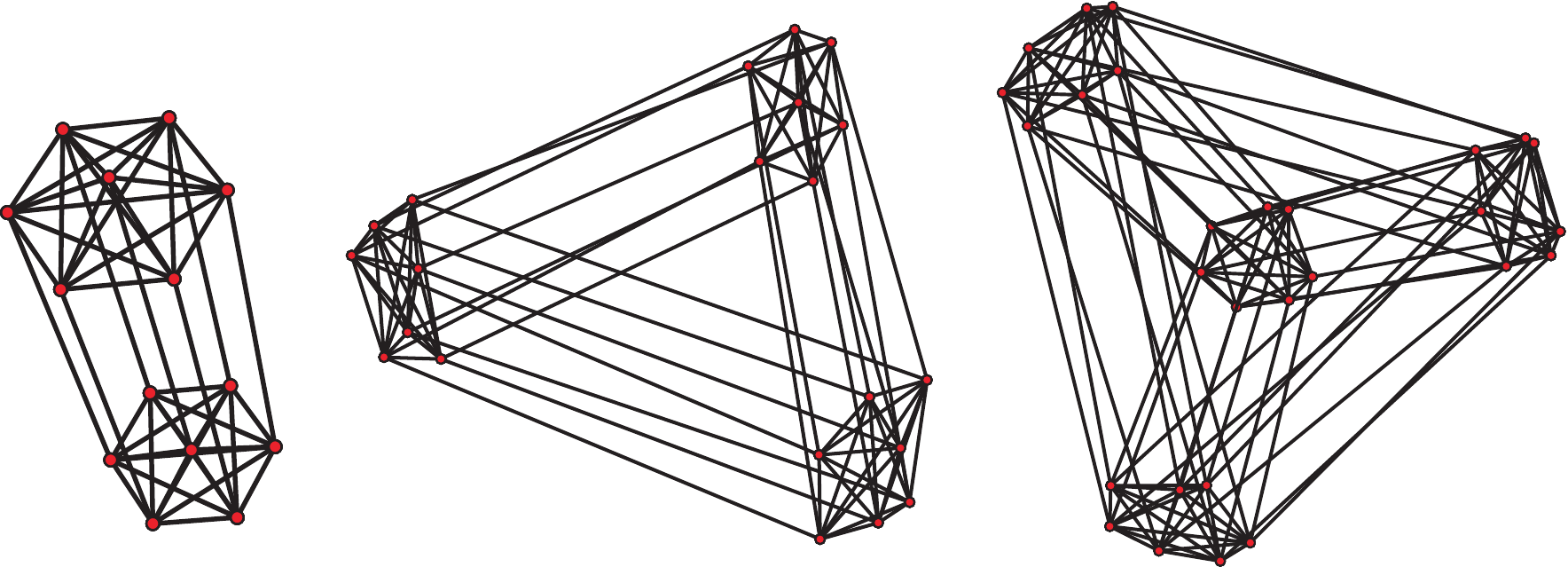}}
	\caption[The graphs of the products~$\simplex_{(i,6)} = \simplex_i\times\simplex_6$, for~$i\in\{1,2,3\}$]{The graphs of the products~$\simplex_{(i,6)} \eqdef \simplex_i\times\simplex_6$, for~$i\in[3]$.}
	\label{intro:fig:prodSimplGraphs}
\end{figure}

\begin{example}
Let~$\sub{n} \eqdef (n_1\dots,n_r)$ be a tuple of positive integers. The product of simplices $\simplex_{\sub{n}} \eqdef \simplex_{n_1}\times\dots\times\simplex_{n_r}$ has dimension~$\sum n_i$ and its \fvector{} is given by
$$f_k(\simplex_{\sub{n}})=\sum_{\substack{0\le k_i\le n_i \\ k_1+\dots+k_r=k}} \prod_{i\in[r]} {n_i+1\choose k_i+1}.$$
In particular,~$f_0(\simplex_{\sub{n}})=\prod n_i$ and~$f_1(\simplex_{\sub{n}})=\sum_{i\in[r]} {n_i+1\choose 2}\prod_{j\ne i} (n_j+1)$. The graphs of~$\simplex_{(i,6)}$ ($i\in[3]$) are represented in \fref{intro:fig:prodSimplGraphs}.
\end{example}

\subsection{Polytopality of products of non-polytopal graphs}\label{intro:subsec:results:nonpolytopal}

In Chapter~\ref{chap:nonpolytopal}, we set out to study the polytopality of products of graphs. The \defn{Cartesian product}\index{Cartesian product!--- of graphs} of two graphs~$G$ and~$H$ is the graph~$G\times H$ with vertex set~$V(G\times H) \eqdef V(G)\times V(H)$ and edge set~$E(G\times H) \eqdef \big(V(G)\times E(H)\big)\cup\big(E(G)\times V(H)\big)$. In other words, for~$a,c\in V(G)$ and~$b,d\in V(H)$, the vertices~$(a,b)$ and~$(c,d)$ of~$G\times H$ are adjacent if either~$a=c$ and $\{b,d\}\in E(H)$, or~$b=d$ and~$\{a,c\}\in E(G)$. This product is consistent with the product of polytopes: the graph of a product of two polytopes is the product of their graphs. In particular, the product of two polytopal graphs is automatically polytopal. In this chapter, we consider the reciprocal question: given two graphs~$G$ and~$H$, does the polytopality of the product~$G\times H$ imply the polytopality of the factors~$G$ and~$H$?

The product of a triangle by a path of length~$2$ is a simple counter-example to this question: although the path is not polytopal, the resulting graph is the graph of a \poly{3}tope obtained by gluing two triangular prisms along a triangular face. We neutralize such simple examples by requiring furthermore both factors to be regular graphs. Let~$d$ and~$e$ denote the regularity degree of~$G$ and~$H$ respectively. In this case, the product~$G\times H$ is \regular{(d+e)} and it is natural to wonder whether it is the graph of a simple \poly{(d+e)}tope. The answer is given by the following theorem:

\begin{theorem}\label{intro:theo:simpleproduct}
The product~$G\times H$ is the graph of a simple polytope if and only if both~$G$ and~$H$ are graphs of simple polytopes. In this case, there is a unique (combinatorial type of) simple polytope whose graph is~$G\times H$: it is precisely the product of the unique (combinatorial types of) simple polytopes whose graphs are~$G$ and~$H$ respectively.\qed
\end{theorem}

In this theorem, the uniqueness of the simple polytope realizing~$G\times H$ is a direct application of the fact that a simple polytope is uniquely determined by its graph~\cite{bm-ppi-87,k-swtsp-88}. These results rely on the following basic property of simple polytopes: any~$k+1$ edges adjacent to a vertex of a simple polytope~$P$ define a \face{k} of~$P$.

As an application of Theorem~\ref{intro:theo:simpleproduct} we obtain a large family of non-polytopal \regular{4} graphs: the product of a non-polytopal \regular{3} graph by a segment is non-polytopal~and~\regular{4}.

\begin{figure}[b]
	\capstart
	\centerline{\includegraphics[width=.9\textwidth]{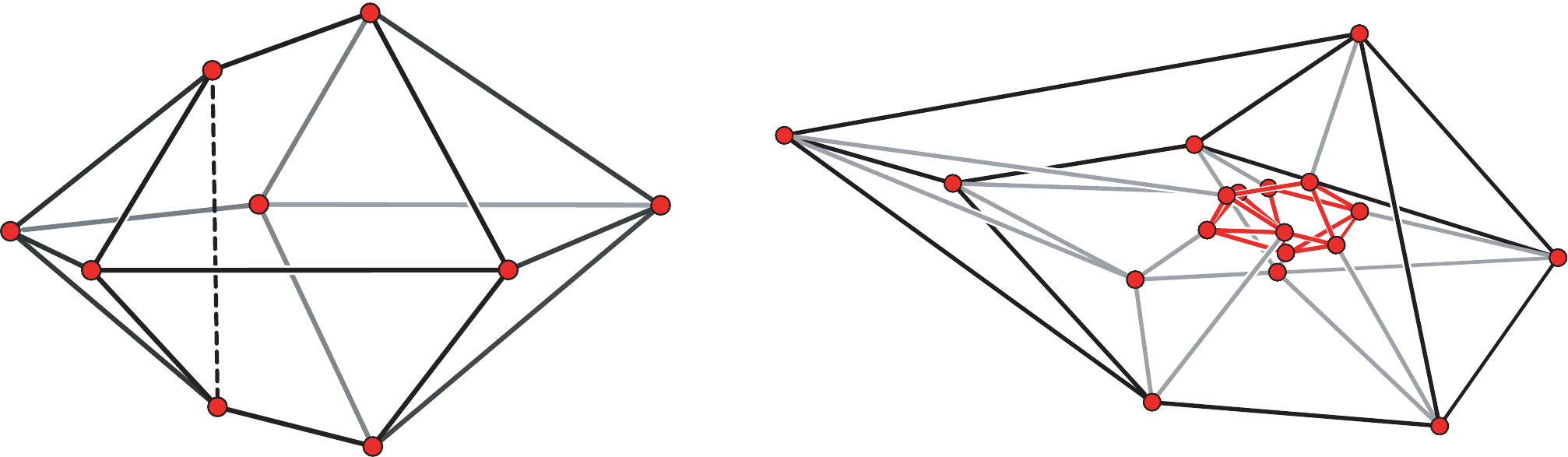}}
	\caption[A non-polytopal \regular{4} graph~$H$ and the Schlegel diagram of a \poly{4}tope whose graph is the product of~$H$ by a segment]{A non-polytopal \regular{4} graph~$H$ which is the graph of a regular subdivision of a \poly{3}tope (left) and the Schlegel diagram of a \poly{4}tope whose graph is the product of~$H$ by a segment (right).}
	\label{intro:fig:truncatedoctahedron}
\end{figure}

We then wonder whether the product of two non-polytopal regular graphs can be polytopal in a dimension smaller than its degree. The following examples partially answer this question:

\begin{theorem}
\begin{enumerate}[(i)]
\item For~$n\ge 3$, the product~$K_{n,n}\times K_2$ of a complete bipartite graph by a segment is not polytopal.
\item The product of a \poly{d}topal graph by the graph of a regular subdivision of an \poly{e}tope is \poly{(d+e)}topal. 
This provides polytopal products of non-polytopal regular graphs (see for example \fref{intro:fig:truncatedoctahedron}).\qed
\end{enumerate}
\end{theorem}

\subsection{Prodsimplicial-neighborly polytopes}\label{intro:subsec:results:psn}

In Chapter~\ref{chap:psn}, we consider polytopes whose skeletons are that of a product of simplices:

\begin{definition}
\index{polytope!psn@\xpsn{(k,\sub{n})} ---}
Let~$k\ge0$ and~$\sub{n} \eqdef (n_1,\dots,n_r)$, with~$r\ge1$ and~$n_i\ge 1$ for all~$i$. A polytope is \defn{$(k,\sub{n})$-prodsimplicial-neighborly}~---~or \defn{\xpsn{(k,\sub{n})}} for short~---~if its \skeleton{k} is combinatorially equivalent to that of the product of simplices~$\simplex_{\sub{n}} \eqdef \simplex_{n_1}\times\cdots\times\simplex_{n_r}$.
\end{definition}

This definition is essentially motivated by two particular classes of \psn polytopes:
\index{polytope!neighborly ---}
\index{neighborly polytope}
\index{polytope!neighborly cubical ---}
\index{neighborly cubical polytope}
\begin{enumerate}[(i)]
\item \defn{neighborly} polytopes arise when $r=1$;
\item \defn{neighborly cubical} polytopes~\cite{jz-ncp-00,js-ncps-07,sz-capdp} arise when~$\sub{n}=(1,1,\dots,1)$.
\end{enumerate}

\begin{remark}
In the literature, a polytope is \defn{\neighborly{k}} if any subset of at most~$k$ of its vertices forms a face. Observe that such a polytope is \xpsn{(k-1,n)} with our notation.
\end{remark}

Obviously, the product~$\simplex_{\sub{n}}$ itself is a \xpsn{(k,\sub{n})} polytope of dimension~$\sum n_i$. We are naturally interested in finding \xpsn{(k,\sub{n})} polytopes in smaller dimensions. For example, the cyclic polytope~$C_{2k+2}(n+1)$ is a \xpsn{(k,n)} polytope in dimension~$2k+2$. We denote by~$\delta(k,\sub{n})$ the smallest possible dimension that a \xpsn{(k,\sub{n})} polytope can have.

\psn polytopes can be obtained by projecting the product~$\simplex_{\sub{n}}$, or a combinatorially equivalent polytope, onto a smaller subspace. For example, the cyclic polytope~$C_{2k+2}(n+1)$ (just like any polytope with~$n+1$ vertices) can be seen as a projection of the simplex~$\simplex_n$ to~$\R^{2k+2}$.

\begin{definition}
\index{polytope!ppsn@\xppsn{(k,\sub{n})} ---}
We say that a \xpsn{(k,\sub{n})} polytope is \defn{$(k,\sub{n})$-projected-prodsimplicial-neighborly} ---~or \defn{\xppsn{(k,\sub{n})}} for short~---~if it is a projection of a polytope combinatorially equivalent to~$\simplex_{\sub{n}}$.
\end{definition}

We denote by~$\delta_{pr}(k,\sub{n})$ the smallest possible dimension of a \xppsn{(k,\sub{n})} polytope.

\mvs
Chapter~\ref{chap:psn} may be naturally divided into two parts. In the first part (Sections~\ref{psn:sec:cyclic} and~\ref{psn:sec:deformedproducts}), we present three methods for constructing low-dimensional \ppsn polytopes:
\begin{enumerate}[(i)]
\item Reflections of cyclic polytopes;
\item Minkowski sums of cyclic polytopes;
\item Deformed product constructions in the spirit of Raman Sanyal and G\"unter Ziegler~\cite{z-ppp-04,sz-capdp}.
\end{enumerate}
The second part (Section~\ref{psn:sec:topologicalObstruction}) derives topological obstructions for the existence of such objects, using techniques developed by Raman Sanyal in~\cite{s-tovnms-09} to bound the number of vertices of Minkowski sums. In view of these obstructions, our constructions in the first part turn out to be optimal for a wide range of parameters.

\paragraph{Constructions.}
Our first non-trivial example is a \xpsn{(k,(1,n))} polytope in dimension~$2k+2$, obtained in Section~\ref{psn:subsec:cyclic:reflection} by reflecting the cyclic polytope~$C_{2k+2}(n+1)$ in a well-chosen~hyperplane:

\begin{proposition}
For any~$k\ge0$,~$n\ge2k+2$ and~$\lambda\in\R$ sufficiently large, the polytope
$$\conv\left(\ens{(t_i,\dots,t_i^{2k+2})^T}{i\in[n+1]} \cup \ens{(t_i,\dots,t_i^{2k+1},\lambda-t_i^{2k+2})^T}{i\in[n+1]}\right)$$
is a \xpsn{(k,(1,n))} polytope of dimension~$2k+2$.\qed
\end{proposition}

For example, this provides us with a \poly{4}tope whose graph is~$K_2\times K_n$, for any~$n\ge3$.

\svs
Next, we form in Section~\ref{psn:subsec:cyclic:minkowskiSumCyclicPolytopes} well-chosen Minkowski sums of cyclic polytopes to obtain explicit coordinates for \xppsn{(k,\sub{n})} polytopes:

\begin{theorem}\label{intro:theo:UBminkowskiCyclic}
Let~$k\ge0$ and~$\sub{n} \eqdef (n_1,\dots,n_r)$ with~$r\ge1$ and~$n_i\ge1$ for all~$i$. There exist index sets~$I_1,\dots,I_r\subset\R$, with~$|I_i|=n_i$ for all~$i$, such that the polytope
$$\conv\ens{w_{a_1,\dots,a_r}}{(a_1,\dots,a_r)\in I_1\times\cdots\times I_r} \subset \R^{2k+r+1}$$
is \xppsn{(k,\sub{n})}, where $w_{a_1,\dots,a_r} \eqdef \big(a_1,\dots,a_r,\sum_{i\in[r]}a_i^2,\dots,\sum_{i\in[r]} a_i^{2k+2}\big)^T$. Consequently:
$$\delta(k,\sub{n}) \le \delta_{pr}(k,\sub{n}) \le 2k+r+1.$$
\end{theorem}
\vspace{-.7cm}\qed
\vspace{.4cm}

For $r=1$ we recover neighborly polytopes. 

\svs
Finally, we extend in Section~\ref{psn:sec:deformedproducts} Raman Sanyal and G\"unter Ziegler's technique of ``projecting deformed products of polygons''~\cite{z-ppp-04,sz-capdp} to products of arbitrary simple polytopes: we  suitably project a suitable polytope combinatorially equivalent to a given product of simple polytopes in such a way as to preserve its complete \skeleton{k}. More concretely, we describe how to use colorings of the graphs of the polar polytopes of the factors in the product to raise the dimension of the preserved skeleton. The basic version of this technique yields the following result:

\begin{proposition}
Let~$P_1,\dots,P_r$ be simple polytopes. For each polytope~$P_i$, denote by~$n_i$ its dimension, by~$m_i$ its number of facets, and by~$\chi_i \eqdef \chi(\gr(P_i^\polar))$ the chromatic number of the graph of its polar polytope~$P_i^\polar$. For a fixed integer~$d\le n$, let~$t$ be maximal such that~${\sum_{i=1}^t n_i\le d}$. Then there exists a \poly{d}tope whose \skeleton{k} is combinatorially equivalent to that of the product $P_1\times\cdots\times P_r$ as soon as
$$0 \le k \le \sum_{i=1}^r (n_i-m_i) + \sum_{i=1}^t (m_i-\chi_i) + \Floor{\frac{1}{2}\left(d-1+\sum_{i=1}^t(\chi_i-n_i)\right)}.$$
\end{proposition} 
\vspace{-1cm}\qed
\vspace{.8cm}

Among polytopes that minimize the last summand are products of \defn{even polytopes} (all \face{2}s have an even number of vertices). See Example~\ref{psn:exm:even} for the details, and the end of Section~\ref{psn:subsec:deformedproducts:general} for extensions of this technique.

\svs
Specializing the factors to simplices provides another construction of \ppsn~polytopes. When some of these simplices are small compared to~$k$, this technique in fact yields our best examples of \ppsn polytopes:

\begin{theorem}
For any~$k\ge0$ and~$\sub{n} \eqdef (n_1,\dots,n_r)$ with~${1=n_1=\cdots=n_s<n_{s+1}\le\cdots\le n_r}$,
$$\delta_{pr}(k,\sub{n}) \le
\begin{cases}
     2(k+r)-s-t & \text{if } 3s \le 2k+2r, \\
     2(k+r-s)+1 & \text{if } 3s = 2k+2r+1, \\
     2(k+r-s+1) & \text{if } 3s \ge 2k+2r+2,
\end{cases}$$
where~$t\in\{s,\dots,r\}$ is maximal such that~$3s+\sum_{i=s+1}^{t}(n_i+1) \le 2k+2r$.\qed
\end{theorem}

If~$n_i=1$ for all~$i$, we recover the neighborly cubical polytopes of~\cite{sz-capdp}.

\paragraph{Obstructions.}
In order to derive lower bounds on the minimal dimension~$\delta_{pr}(k,\sub{n})$ that a \xppsn{(k,\sub{n})} polytope can have, we apply in Section~\ref{psn:sec:topologicalObstruction} a method due to Raman Sanyal~\cite{s-tovnms-09}. For any projection which preserves the \skeleton{k} of~$\simplex_{\sub{n}}$, we construct via Gale duality a simplicial complex guaranted to be embeddable in a space of a certain dimension. The argument is then a topological obstruction based on Sarkaria's criterion for the embeddability of a simplicial complex in terms of colorings of Kneser graphs~\cite{m-ubut-03}. We obtain the following result:

\begin{theorem}\label{intro:theo:topObstr}
Let~$\sub{n} \eqdef (n_1,\dots,n_r)$ with~$1=n_1=\cdots=n_s<n_{s+1}\le\cdots\le
n_r$. Then:
\begin{enumerate}
\item If 
$$0 \le k \le \sum_{i=s+1}^r \Fracfloor{n_i-2}{2} + \max\left\{0,\Fracfloor{s-1}{2}\right\},$$
then~$\delta_{pr}(k,\sub{n}) \ge 2k+r-s+1$.
\item If~$k\ge \Floor{\frac{1}{2} \sum_i n_i}$ then~$\delta_{pr}(k,\sub{n}) \ge \sum_i n_i$.\qed
\end{enumerate}
\end{theorem}

In particular, the upper and lower bounds provided by Theorems~\ref{intro:theo:UBminkowskiCyclic} and~\ref{intro:theo:topObstr} match over a wide range of parameters:

\begin{theorem}\label{intro:theo:mainResult}
For any~$\sub{n} \eqdef (n_1,\dots,n_r)$ with~$r\ge1$ and~$n_i\ge2$ for all~$i$, and for any~$k$ such that~$0\le k\le \sum_{i\in [r]}\Fracfloor{n_i-2}{2}$, the smallest \xppsn{(k,\sub{n})} polytope has dimension exactly~$2k+r+1$. In other words:
$$\delta_{pr}(k,\sub{n}) = 2k+r+1.$$
\end{theorem}
\vspace{-.7cm}\qed
\vspace{.2cm}

\begin{remark}
The last two sections of Chapter~\ref{chap:psn} consist in applying (and partially extending) methods and results on projections of polytopes developed by Raman Sanyal and G\"unter Ziegler~\cite{z-ppp-04,sz-capdp,s-tovnms-09}. We have decided to develop them in this dissertation since their application to products of simplices yields new results which complete our study on polytopality of products. Moreover, after the completion of our work on obstructions for projections of products of simplices, we learned that Thilo R\"orig and Raman Sanyal obtained similar results in a recent work~\cite{rs-npps} (see also~\cite{rorig-phd,sanyal-phd}).
\end{remark}


\section{Sources of material}

The work presented in the second part of this thesis was initiated during my stay at the \defn{Centre de Recerca Matem\`atica} in Barcelona for the \defn{i-\textsc{math} Winter School DocCourse Combinatorics and Geometry 2009: Discrete and Computational Geometry}~\cite{crm}, under the tutoring of Julian Pfeifle. The results presented in the coming chapters have been submitted in two papers in collaboration with participants of this event:
\begin{enumerate}
\item The motivation for studying polytopality of products of graphs came from the course of G\"unter Ziegler~\cite{crm} who presented the prototype example of this question: ``is the product of two Petersen graphs polytopal?'' The (partial) answers that we propose to this question are developed in a preprint written in collaboration with Julian Pfeifle and Francisco Santos~\cite{pps-ppg}.
\item The study of \psn polytopes was proposed by Julian Pfeifle as a research project during this DocCourse~\cite{crm}. Our results are exposed in a preprint in collaboration with Benjamin Matschke and Julian Pfeifle~\cite{mpp-psnp}.
\end{enumerate}
	\chapter{Cartesian products of non-polytopal graphs}\label{chap:nonpolytopal}

This chapter is devoted to a discussion on polytopality of graphs. In Section~\ref{nonpolytopal:sec:polytopalitygraphs}, we shortly review and discuss the current knowledge concerning polytopality of general graphs. We attach a particular attention to a variety of examples (some of which are well-known while some others are originals), which we will use afterwards. In Section~\ref{nonpolytopal:sec:product}, we focus on Cartesian products of graphs. This product is defined in such a way that the graph of a product of polytopes is the product of their graphs; in particular, products of polytopal graphs are automatically polytopal. This Chapter concerns the reciprocal statement: are the two factors of a polytopal product necessarily polytopal?


\section{Polytopality of graphs}\label{nonpolytopal:sec:polytopalitygraphs}

\begin{definition}
\index{polytopal}
A graph~$G$ is \defn{polytopal} if it is isomorphic to the \skeleton{1} of some polytope~$P$. If~$P$ is \dimensional{d}, we say that~$G$ is \poly{d}topal.
\end{definition}

In small dimension, polytopality is easy to deal with. For example, \poly{2}topal graphs are exactly cycles. The first interesting question is \poly{3}topality, which is characterized by Steinitz' ``fundamental theorem of convex types'':

\begin{theorem}[Steinitz~\cite{s-pr-22}]\label{nonpolytopal:theo:steinitz}\index{Steinitz' Theorem|hbf}
A graph~$G$ is the graph of a \poly{3}tope~$P$ if and only if $G$~is planar and \connected{3}. Moreover, the combinatorial type of~$P$ is uniquely determined by~$G$.\qed
\end{theorem}

We refer to~\cite{g-cp-03,z-lp-95} for a discussion on three possible approaches for the proof of this fundamental theorem.

Throughout this chapter, we will notice that a first step to realize a graph~$G$ is to understand the possible face lattice of a polytope realizing~$G$. For example, it is often difficult to decide what cycles of~$G$ can be the graphs of \face{2}s of a \poly{d}tope realizing~$G$. In dimension~$3$, graphs of \face{2}s are characterized by the following separation condition:

\begin{theorem}[Whitney~\cite{w-nspg-32}]\label{nonpolytopal:theo:whitney}
Let~$G$ be the graph of a \poly{3}tope~$P$. The graphs of the faces of~$P$ are precisely the induced cycles in~$G$ that do not separate~$G$.\qed
\end{theorem}

In contrast to the easy $2$- and \dimensional{3} worlds, \poly{d}topality becomes much more involved as soon as~$d\ge 4$. As an illustration, the existence of neighborly polytopes (see the cyclic polytopes in Proposition~\ref{intro:prop:cyclic}) proves that all possible edges can be present in the graph of a \poly{4}tope. Starting from a neighborly polytope, and stacking vertices on undesired edges, one can even observe the following:

\begin{observation}[Perles]
Every graph is an induced subgraph of the graph of a \poly{4}tope.
\end{observation}

It is a long-standing question of polytope theory how to determine whether a graph is~\poly{d}to\-pal or not. In the next section, we recall some general necessary conditions and apply them to discuss polytopality of small examples.


\subsection{Necessary conditions for polytopality}\label{nonpolytopal:subsec:polytopalitygraphs:necessaryconditions}

\begin{proposition}\label{nonpolytopal:prop:necessaryconditions}
A \poly{d}topal graph~$G$ satisfies the following properties:
\begin{enumerate}
\item \defn{Balinski's Theorem}: $G$~is \connected{d}~\cite{b-gscps-61}.

\item \defn{Principal Subdivision Property} ($d$-PSP): Every vertex of~$G$ is the principal vertex of a principal subdivision of~$K_{d+1}$. Here, a \defn{subdivision} of~$K_{d+1}$ is obtained by replacing edges by paths, and a \defn{principal subdivision}\index{subdivision!principal ---} of~$K_{d+1}$ is a subdivision in which all edges incident to a distinguished \defn{principal vertex} are not subdivided~\cite{b-ncp-67}.

\item \defn{Separation Property}: The maximal number of components into which~$G$ may be separated by removing~$n>d$ vertices equals~$f_{d-1}\big(C_d(n)\big)$, the maximum number of facets of a \poly{d}tope with $n$~vertices~\cite{k-ppg-64}.\qed
\end{enumerate}
\end{proposition}

\begin{remark}\label{nonpolytopal:remark:3poly}
The principal subdivision property together with Steinitz' Theorem ensure that no graph of a \poly{3}tope is \poly{d}topal for~$d\ne3$. In other words, any \poly{3}tope is the unique polytopal realization of its graph. This property is also obviously true in dimension~$0$, $1$ or $2$. In contrast, it is strongly wrong in dimension~$4$ and higher.
\end{remark}

Before providing examples of application of Proposition~\ref{nonpolytopal:prop:necessaryconditions}, let us insist on the fact that these necessary conditions are not sufficient (see also Examples~\ref{nonpolytopal:exm:marc+antonio} and~\ref{nonpolytopal:exm:diamond}):

\begin{example}[Non-polytopality of the complete bipartite graph~\cite{b-ncp-67}]
For any two integers $m,n\ge 3$, the complete bipartite graph~$K_{m,n}$ is not polytopal, although~$K_{n,n}$ satisfies all properties of Proposition~\ref{nonpolytopal:prop:necessaryconditions} to be \poly{4}topal as soon as~$n\ge 7$.

Indeed, assume that~$K_{n,m}$ is the graph of a \poly{d}tope~$P$. Then~$d\ge4$ because~$K_{n,m}$ is non-planar. Consider the induced subgraph~$H$ of~$K_{n,m}$ corresponding to some \face{3}~$F$ of~$P$. Because~$H$ is induced and has minimum degree at least~$3$, it contains a $K_{3,3}$ minor, so $F$~was not a \face{3} after all.
\end{example}

\begin{example}[Circulant graphs]
\index{circulant graph}
Let~$n$ be an integer and~$S$~be a subset of~$\left\{1,\dots,\Fracfloor{n}{2}\right\}$. The \defn{circulant} graph~$\Gamma_n(S)$ is the graph whose vertex set is~$\Z_n$ and whose edge set is the set of pairs of vertices whose difference lies in~$S\cup(-S)$. Observe that the degree of~$\Gamma_n(S)$ is precisely $|S\cup(-S)|$ (in particular, the degree is odd only if~$n$ is even and~$S$ contains~$n/2$) and that $\Gamma_n(S)$ is connected if and only if~$S\cup\{n\}$ is relatively prime. For example, \fref{nonpolytopal:fig:circulant} represents all connected circulant graphs on at most~$8$ vertices.

\begin{figure}
	\capstart
	\centerline{\includegraphics[width=\textwidth]{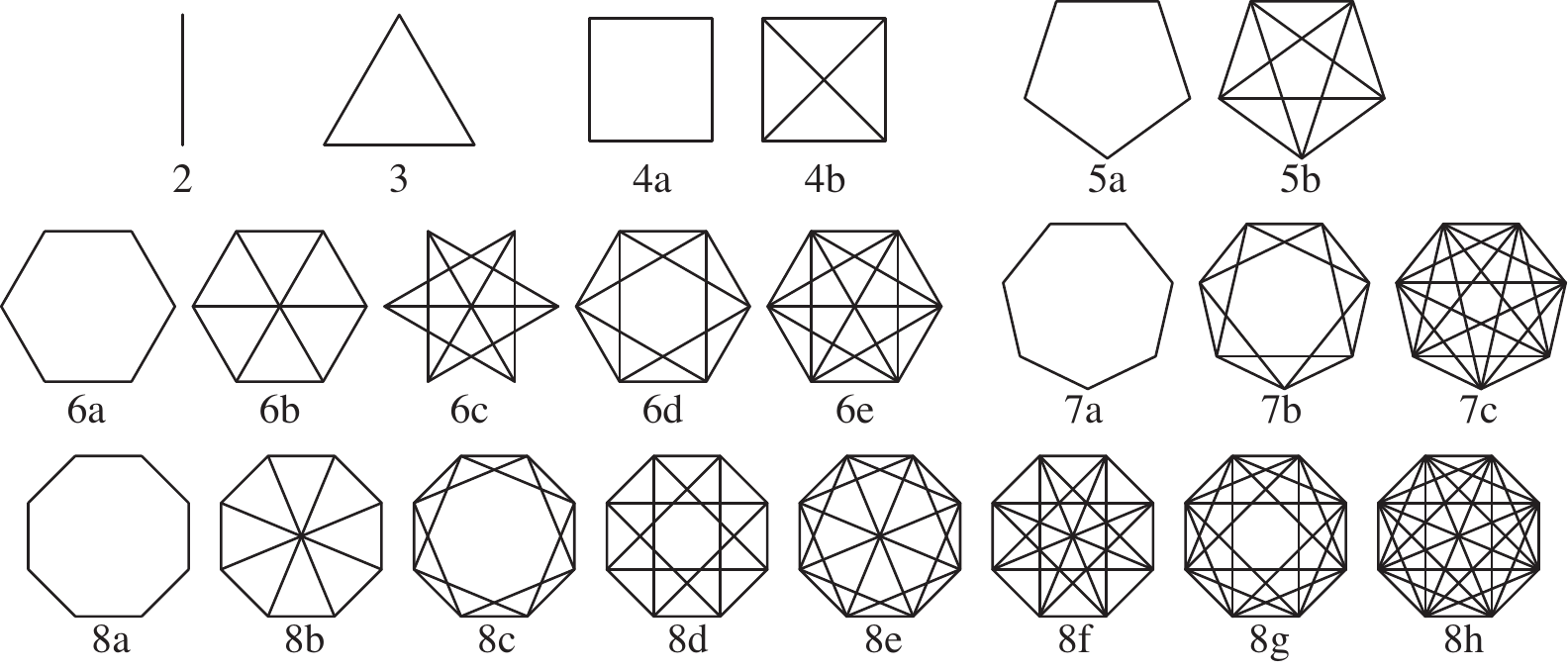}}
	\caption[Connected circulant graphs with at most~$8$ vertices]{Connected circulant graphs with at most~$8$ vertices.}
	\label{nonpolytopal:fig:circulant}
\end{figure}

\newpage
Using Proposition~\ref{nonpolytopal:prop:necessaryconditions} we can determine the polytopality of various circulant graphs:
\begin{itemize}
\svs
\item[~~~~\fbox{Degree~$2$}] A connected graph of degree~$2$ is a cycle, and thus the graph of a polygon.

\svs
\item[~~~~\fbox{Degree~$3$}] Up to isomorphism, the only connected circulant graphs of degree~$3$ are~$\Gamma_{2m}(1,m)$ and~$\Gamma_{4m+2}(2,2m+1)$. When~$m\ge 3$, the first one is not planar, and thus not polytopal. The second one is the graph of a prism over a \gon{(2m+1)}.

\svs
\item[~~~~\fbox{Degree~$4$}] As soon as we reach degree~$4$, we cannot provide a complete description of polytopal circulant graphs, but we can discuss special cases, namely the circulant graphs~$\Gamma_n(1,s)$ for~$s\in\{2,3,4\}$:
\begin{enumerate}[(a)]
\item $s=2$: For any~$m\ge 2$, the graph~$\Gamma_{2m}(1,2)$ is the graph of an antiprism over an \gon{m}. In contrast, for any~$m\ge 3$, the graph~$\Gamma_{2m+1}(1,2)$ is not polytopal: it is not planar and does not satisfy the principal subdivision property for dimension~$4$.
\item $s\in\{3,4\}$: For any~$n\ge7$, the graph~$\Gamma_n(1,3)$ is not polytopal. Indeed, the \cycle{4}s induced by the vertices~$\{1,2,3,4\}$ and~$\{2,3,4,5\}$ should define \face{2}s of any realization (because of Theorem~\ref{nonpolytopal:theo:whitney} in dimension~$3$ and of Proposition~\ref{nonpolytopal:prop:cycles} in dimension~$4$), but they intersect improperly. Similarly, for any~$n\ge9$, the graph~$\Gamma_n(1,4)$ is not polytopal.
\end{enumerate}

\svs
\item[~~~~\fbox{Degree~$n-2$}] The graph~$\Gamma_{2m}(1,2,\dots,m-1)$ is the only circulant graph with two vertices more than its degree. It is not planar when~$m\ge 4$ and it is not \poly{(2m-2)}topal since it does not satisfy the principal subdivision property in this dimension. However, it is always the graph of the \dimensional{m} cross-polytope, and when~$m$ is even, it is also the graph of the join of two \dimensional{(m/2)} cross-polytopes.

\svs
\item[~~~~\fbox{Degree~$n-1$}] The complete graph on~$n$ vertices is the graph of any neighborly polytope, and its polytopality range is~$\{4,\dots,n-1\}$ (as soon as~$n\ge 5$).
\end{itemize}

\svs
The sporadic cases developed above are sufficient to determine the polytopality range of all circulant graphs on at most~$8$ vertices, except the graphs~8e and~8f of \fref{nonpolytopal:fig:circulant} that we treat separately now. None of them can be \poly{3}topal since they are not planar. We prove that they are not \poly{4}topal by discussing what could be the \face{3}s of a possible realization:
\begin{itemize}
\item We start with the graph~$\Gamma_8(1,3,4)$ represented in \fref{nonpolytopal:fig:circulant}(8f). Consider any subgraph of~$\Gamma_8(1,3,4)$ induced by~$6$ vertices. If the distance between the two missing vertices is odd (resp.~even), then the subgraph is not planar (resp.~not \connected{3}). Consequently, any subgraph of~$\Gamma_8(1,3,4)$ induced by~$7$ vertices is not planar, while any subgraph of~$\Gamma_8(1,3,4)$ induced by~$5$ vertices is not \connected{3}. Thus, the only possible \face{3}s are tetrahedra, but~$\Gamma_8(1,3,4)$ contains only~$4$ induced~$K_4$. Thus,~$\Gamma_8(1,3,4)$ is not polytopal.
\item The case of~$\Gamma_8(1,2,4)$ is more involved. Up to rotation, its only \poly{3}topal induced subgraphs are represented in \fref{nonpolytopal:fig:C8124}. Assume that the subgraph induced by~$\{0,1,2,3,4,5\}$ defines a \face{3}~$F$ in a realization~$P$ of~$\Gamma_8(1,2,4)$. Then the triangle~$123$ is a \face{2} of~$P$ and thus should be contained in another \face{3} of~$P$. But any \face{3} which contains~$123$ also contains either~$0$ or~$4$, and thus intersects improperly with~$F$. Consequently, the subgraph induced by~$\{0,1,2,3,4,5,6\}$ cannot define a \face{3} of a realization of~$\Gamma_8(1,2,4)$. For the same reason, the subgraphs induced by~$\{0,1,2,3,4,6\}$ and~$\{0,1,2,3,4\}$ cannot define \face{3}s. Assume now that the subgraph induced by $\{0,1,2,3,5,6\}$ forms a \face{3} in a realization~$P$ of~$\Gamma_8(1,2,4)$. Then the triangle~$123$ is a \face{2} of~$P$ and should be contained in another \face{3} of~$P$. The only possibility is the subgraph induced by~$\{1,2,3,4,6,7\}$ which intersects improperly~$F$. Finally, the only possible \face{3}s are the two tetrahedra induced respectively by the odd and the even vertices, and thus,~$\Gamma_8(1,2,4)$ is not polytopal.
\end{itemize}

\begin{figure}
	\capstart
	\centerline{\includegraphics[scale=1]{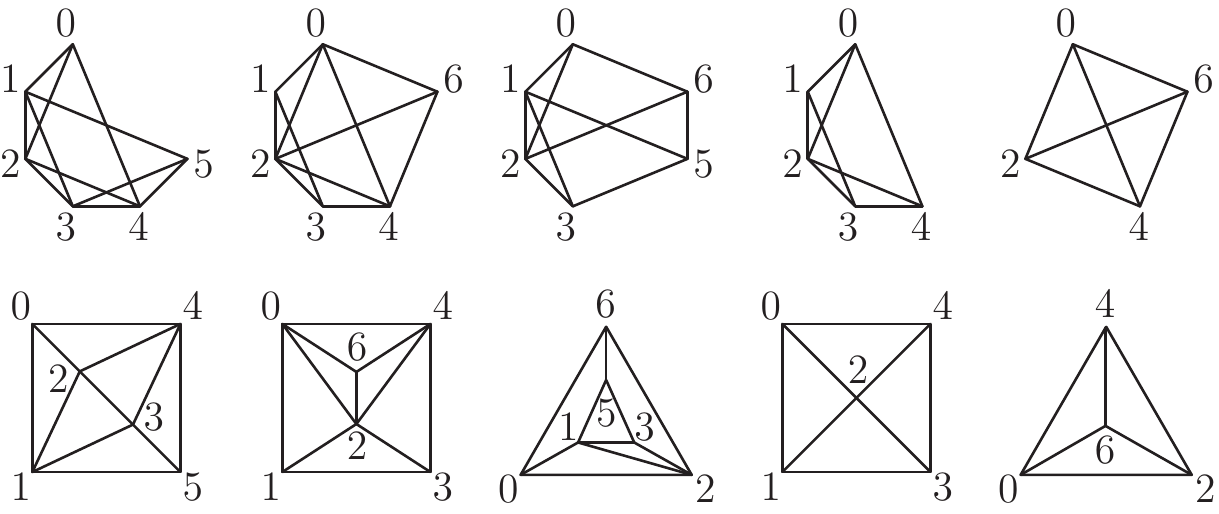}}
	\caption[The \poly{3}topal induced subgraphs in the circulant graph~$\Gamma_8(1,2,4)$]{The \poly{3}topal induced subgraphs of the circulant graph~$\Gamma_8(1,2,4)$. The faces of the planar drawing below each of these subgraphs are the \face{2}s of the corresponding \poly{3}tope.}
	\label{nonpolytopal:fig:C8124}
\end{figure}

\svs
The following table summarizes the polytopality range of all the graphs of \fref{nonpolytopal:fig:circulant}:
\begin{center}
\begin{tabular}{|c|c|cc|cc|ccccc|}
\hline
2 & 3 & 4a & 4b & 5a & 5b & 6a & 6b & 6c & 6d & 6e \\
$\{1\}$ & $\{2\}$ & $\{2\}$ & $\{3\}$ & $\{2\}$ & $\{4\}$ & $\{2\}$ & $\emptyset$ & $\{3\}$ & $\{3\}$ & $\{4,5\}$ \\
\hline
\end{tabular}

\enlargethispage{.3cm}
\svs
\begin{tabular}{|ccc|cccccccc|}
\hline
7a & 7b & 7c & 8a & 8b & 8c & 8d & 8e & 8f & 8g & 8h \\
$\{2\}$ & $\emptyset$ & $\{4,5,6\}$ & $\{2\}$ & $\emptyset$ & $\{3\}$ & $\emptyset$ & $\emptyset$ & $\emptyset$ & $\{4,5\}$ & $\{4,5,6,7\}$ \\
\hline
\end{tabular}
\end{center}
\end{example}

\begin{example}[A graph whose polytopality range is~$\{d\}$~\cite{k-ppg-64}]
An interesting application of the separation property of Proposition~\ref{nonpolytopal:prop:necessaryconditions} is the possibility to construct, for any integer~$d$, a polytope whose polytopality range is exactly the singleton~$\{d\}$. The construction, proposed by Victor Klee~\cite{k-ppg-64}, consists in stacking a vertex on all facets of the cyclic polytope~$C_d(n)$ (for example on all facets of a simplex). The graph of the resulting polytope can be separated in~$f_{d-1}(C_d(n))$ isolated points by removing the~$n$ initial vertices, and thus is not \poly{d'}topal for~$d'<d$, by the separation property. It can not be \poly{d'}topal for $d'>d$ either, since the stacked vertices have degree~$d$ (because the cyclic polytope is simplicial). Thus, the dimension of the resulting graph is not ambiguous.
\end{example}


\subsection{Simple polytopes}

\index{polytope!simple ---|hbf}
A \poly{d}tope is \defn{simple} if its vertex figures are simplices. In other words, its facet-defining hyperplanes are in general position, so that a vertex is contained in exactly $d$~facets, and also in exactly $d$~edges (and thus the graph of a simple \poly{d}tope is \regular{d}). Surprisingly, a \regular{d} graph can be realized by at most one simple polytope:

\begin{theorem}[\cite{bm-ppi-87,k-swtsp-88}]\label{nonpolytopal:theo:kalai}
Two simple polytopes are combinatorially equivalent if and only if they have the same graph.\qed
\end{theorem}

This property, conjectured by Micha~Perles, was first proved by Roswitha Blind and Peter Mani~\cite{bm-ppi-87}. Gil Kalai~\cite{k-swtsp-88} then gave a very simple way of reconstructing the face lattice from the graph, and Eric Friedman~\cite{f-fsppt-09} showed that this can even be done in polynomial time.

As mentioned previously, the first step to find a polytopal realization of a graph is often to understand how can the face lattice of this realization look like. Theorem~\ref{nonpolytopal:theo:kalai} ensures that if the realization is simple, there is only one choice. This motivates to  temporarily restrict the study of realization of regular graphs only to simple polytopes:

\begin{definition}
\index{polytopal!simply ---|hbf}
A graph is \defn{simply \poly{d}topal} if it is the graph of a simple \poly{d}tope.
\end{definition}

We can exploit properties of simple polytopes to obtain results on the simple polytopality of graphs. For us, the key property turns out to be that any \tuple{k} of edges incident to a vertex of a simple polytope is contained in a \face{k}. For example, this implies the following result:

\begin{proposition}\label{nonpolytopal:prop:cycles}
All induced cycles of length $3$, $4$ and $5$ in the graph of a simple \poly{d}tope~$P$ are graphs of \face{2}s of~$P$.
\end{proposition}

\begin{proof}
For \cycle{3}s, the result is immediate: any two adjacent edges of a \cycle{3} induce a \face{2}, which must be a triangle because the graph is induced.

Next, let~$\{a,b,c,d\}$ be consecutive vertices of a \cycle{4} in the graph of a simple polytope~$P$. Any pair of edges emanating from a vertex lies in a \face{2} of~$P$. Let~$C_a$ be the \face{2} of~$P$ that contains the edges~$\conv\{a,b\}$ and~$\conv\{a,d\}$. Similarly, let~$C_c$ be the \face{2} of~$P$ that contains~$\conv\{b,c\}$ and~$\conv\{c,d\}$. If~$C_a$ and~$C_c$ were distinct, they would intersect improperly, at least in the two vertices~$b$~and~$d$. Thus,~$C_a=C_c=\conv\{a,b,c,d\}$ is a \face{2} of~$P$.

The case of \cycle{5}s is a little more involved. We first show it for \poly{3}topes. If a \cycle{5}~$C$ in the graph~$G$ of a simple \poly{3}tope does not define a \face{2}, it separates~$G$ into two nonempty subgraphs~$A$~and~$B$ (Theorem~\ref{nonpolytopal:theo:whitney}). Since~$G$~is \connected{3}, both~$A$~and~$B$ are connected to~$C$ by at least three edges. But the endpoints of these six edges must be distributed among the five vertices of~$C$, so one vertex of~$C$ receives two additional edges, and this contradicts simplicity.

For the general case, we show that any \cycle{5}~$C$ in a simple polytope is contained in some \face{3}, and apply the previous argument (a face of a simple polytope is simple). First observe that any three consecutive edges in the graph of a simple polytope lie in a common \face{3}. This is true because any two adjacent edges define a \face{2}, and a \face{2} together with another adjacent edge defines a \face{3}. Thus, four of the vertices of~$C$ are already contained in a \face{3}~$F$. If the fifth vertex~$w$ of~$C$ lies outside~$F$, then the \face{2} defined by the two edges of~$C$ incident to~$w$ intersects improperly with~$F$.
\end{proof}

\begin{remark}
Observe that there is an induced \cycle{6} in the graph of the cube (resp.~an induced \cycle{p} in the graph of a double pyramid over a \cycle{p}, for~$p\ge 3$) which is not the graph of a \face{2}. It is also interesting to notice that contrarily to dimension~$3$ (Theorem~\ref{nonpolytopal:theo:whitney}), the \face{2}s of a \poly{4}tope are not characterized by a separation property: a pyramid over a cube has a non-separating induced \cycle{6} which does not define a \face{2}.
\end{remark}

\begin{corollary}\label{nonpolytopal:coro:cycles}
A simply polytopal graph cannot:
\begin{enumerate}[(i)]
\item be separated by an induced cycle of length~$3$, $4$ or~$5$.
\item contain two induced cycles of length $4$ or~$5$ which share~$3$ vertices.\qed
\end{enumerate}
\end{corollary}

\begin{example}[An infinite family of non-polytopal graphs for non-trivial reasons~\cite{gn-pc-09}]\label{nonpolytopal:exm:marc+antonio}
Consider the family of graphs suggested in \fref{nonpolytopal:fig:marc+antonio}. The $n$th graph of this family is the graph~$G_n$ whose vertex set is~$\Z_{2n+3}\times\Z_2$ and where the vertex~$(x,y)$ is related with the vertices~$(x+y+1,y)$, $(x+y,y+1)$, $(x-y-1,y)$ and~$(x+y-1,y+1)$.

\begin{figure}[h]
	\capstart
	\centerline{\includegraphics[scale=1]{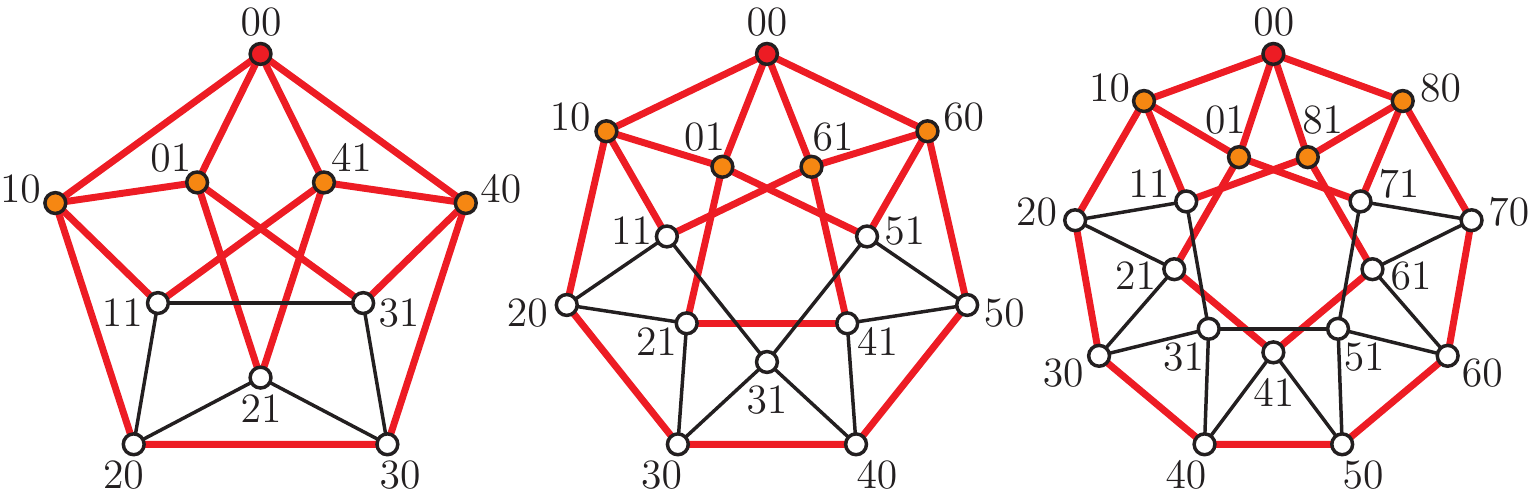}}
	\caption[An infinite family of non-polytopal graphs (for non-trivial reasons)]{An infinite family of non-polytopal graphs (for non-trivial reasons). The vertex~$00$ is the principal vertex of a principal subdivision of~$K_5$ whose edges are colored in red.}
	\label{nonpolytopal:fig:marc+antonio}
\end{figure}

Observe first that the graphs of this family satisfy all necessary conditions of Proposition~\ref{nonpolytopal:prop:necessaryconditions}:
\begin{enumerate}
\item They are \connected{4}: when we remove~$3$ vertices, either the external cycle~$\ens{i0}{i\in\Z_{2n+1}}$ or the internal cycle~$\ens{i1}{i\in\Z_{2n+1}}$ remains a path, to which all the vertices are connected. 
\item They satisfy the principal subdivision property for dimension~$4$: the edges of a principal subdivision of~$K_5$ with principal vertex~$00$ are colored in \fref{nonpolytopal:fig:marc+antonio}.
\item They satisfy the separation property: the cyclic \poly{4}tope on~$m$ vertices has~$\frac{m(m-3)}{2}$ facets while removing~$m$ vertices from~$G_n$ cannot create more than~$m$ connected components.
\end{enumerate}

Consider the first graph~$G_1$ of this family (on the left in \fref{nonpolytopal:fig:marc+antonio}). Since the \cycle{5}s induced by~$\{00,10,20,30,40\}$ and~$\{00,10,20,21,41\}$ share two edges, $G_1$ is not polytopal (because of Theorem~\ref{nonpolytopal:theo:whitney} in dimension~$3$ and of Proposition~\ref{nonpolytopal:prop:cycles} in dimension~$4$). In fact, Proposition~\ref{nonpolytopal:prop:cycles} even excludes all graphs of the family:

\begin{lemma}
None of the graphs of the infinite family suggested in \fref{nonpolytopal:fig:marc+antonio} is polytopal.
\end{lemma}

\begin{proof}
Since they contain a subdivision of~$K_5$, they are not \poly{3}topal.

Denote by~$e_i$ the edge of the external cycle from vertex~$i0$ to vertex~$(i+1)0$. If the graph~$G_n$ were \poly{4}topal, then all $3$- and \cycle{4}s would define \face{2}s. Now consider two consecutive angles $e_1e_2$ and $e_2e_3$ of the external cycle. Each of them defines a \face{2} by simplicity. These two \face{2}s must in fact coincide, since $e_2$~is already contained in a square and a triangular \face{2}, none of which contain the angles~$e_1e_2$ and~$e_2e_3$. By iterating this argument, we obtain that the entire external cycle forms a \face{2}.

Consider a \face{3}~$F$ containing the external cycle. The edge $e_1$ must also be contained in either the adjacent square or the adjacent triangle; without loss of generality, let it be the square. Then the triangle adjacent to the next edge~$e_2$ must also be in~$F$ (because~$F$ already contains two of its edges). By the same reasoning, the square adjacent to~$e_3$ is also contained in~$F$, and iterating this argument (and using that $n$ is odd) shows that in fact~$F$~contains all the squares and triangles of~$G_n$, contradiction.
\end{proof}
\end{example}


\subsection{Truncation and star-clique operation}

We consider the polytope~$\tau_v(P)$ obtained by cutting off a single vertex~$v$ in a polytope~$P$: the set of inequalities defining~$\tau_v(P)$ is that of~$P$ together with a new inequality satisfied by all the vertices of~$P$ except~$v$. The faces of~$\tau_v(P)$ are: (i) all the faces of~$P$ which do not contain~$v$; (ii) the truncations~$\tau_v(F)$ of all faces~$F$ of~$P$ containing~$v$; and (iii) the vertex figure of~$v$ in~$P$ together with all its faces. In particular, if~$v$ is a simple vertex in~$P$, then the truncation of~$v$ in~$P$ replaces~$v$ by a simplex. On the graph of~$P$, it translates into the following transformation:

\begin{definition}
\index{star-clique operation}
Let~$G$ be a graph and~$v$ be a vertex of degree~$d$ of~$G$. The \defn{star-clique operation} (at~$v$) replaces vertex~$v$ by a \clique{d}~$K$, and assigns one edge incident to~$v$ to each vertex of~$K$. The resulting graph~$\sigma_v(G)$ has~$d-1$ more vertices and~${d \choose 2}$ more edges.
\end{definition}

\begin{example}
The truncation of a vertex in a \simp{d} creates (combinatorially) a prism over a \simp{(d-1)}. Its graph is the product~$K_2\times K_d$ and thus can be realized as the product of a segment by any neighborly polytope~$P$ on~$d$ vertices. The next chapter will furthermore prove that it is also \poly{4}topal.
Observe that the complete graph is also realized by any neighborly polytope but that the result of the truncation of a vertex in such a polytope is not always a star-clique operation on the graph (it is the case only when the vertex figure of the cutted vertex in the realization is neighborly).
\end{example}

\begin{proposition}\label{nonpolytopal:prop:starclique}
Let~$v$ be a vertex of degree~$d$ in a graph~$G$. Then~$\sigma_v(G)$ is \poly{d}topal if and only if~$G$ is \poly{d}topal.
\end{proposition}

\begin{proof}
If a \poly{d}tope~$P$ realizes~$G$, then the truncation~$\tau_v(P)$ realizes~$\sigma_v(G)$. For the other direction, consider a \poly{d}tope~$Q$ which realizes~$\sigma_v(G)$. We use simplicity to assert that the \clique{d} replacing~$v$ forms a facet~$F$ of~$Q$. Up to a projective transformation, we can assume that the~$d$ facets of~$Q$ adjacent to~$F$ intersect behind~$F$. Then, removing the inequality defining~$F$ from the facet description of~$Q$ creates a polytope which realizes~$G$.
\end{proof}

\begin{remark}
In degree~$3$, star-clique operations, usually called $\Delta Y$-transformations, are used to prove Steinitz' Theorem~\ref{nonpolytopal:theo:steinitz}. The argument is that any \connected{3} and planar graph can be reduce to the complete graph~$K_4$ by a sequence of such transformations (see~\cite{z-lp-95} for details).
\end{remark}

We can exploit Proposition~\ref{nonpolytopal:prop:starclique} to construct several families of non-polytopal graphs:

\begin{corollary}\label{nonpolytopal:coro:starclique}
Any graph obtained from the graph of a \regular{4} \poly{3}tope by a finite nonempty sequence of star-clique operations is non-polytopal.
\end{corollary}

\begin{proof}
No such graph can be \poly{3}topal since it is not planar. If the resulting graph were \poly{4}topal, Proposition~\ref{nonpolytopal:prop:starclique} would assert that the original graph was also \poly{4}topal, which would contradict Remark~\ref{nonpolytopal:remark:3poly}.
\end{proof}

\begin{example}\label{nonpolytopal:exm:diamond}
For~$n\ge3$, consider the family of graphs suggested by \fref{nonpolytopal:fig:diamond}. They are constructed as follows: place a regular \gon{2n}~$C_{2n}$ into the plane, centered at the origin. Draw a copy~$C_{2n}'$ of~$C_{2n}$ scaled by~$\tfrac{1}{2}$ and rotated by~$\frac{\pi}{2n}$, and lift the vertices of~$C_{2n}'$ alternately to heights~$1$ and~$-1$ into the third dimension. The graph~$\Diamond_n$ is the graph of the convex hull of the result.

\begin{figure}[b]
	\capstart
	\centerline{\includegraphics[scale=1]{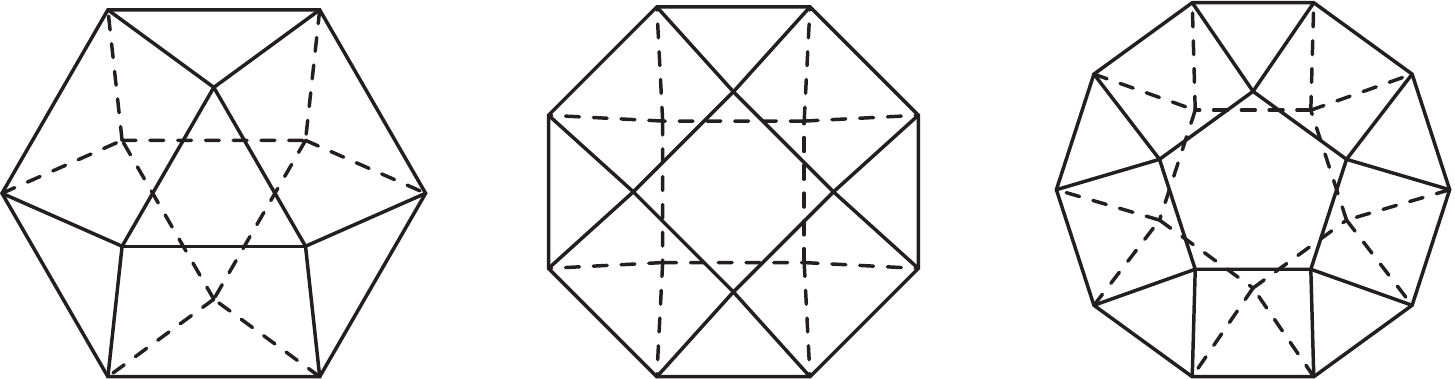}}
	\caption[The graphs $\Diamond_n$ for $n\in\{3,4,5\}$]{The graphs $\Diamond_n$ for $n\in\{3,4,5\}$.}
	\label{nonpolytopal:fig:diamond}
\end{figure}

Let~$\Diamond_n^\star$~be the result of successively applying the star-clique operation to all vertices on the external cycle~$C_{2n}$. Corollary~\ref{nonpolytopal:coro:starclique} ensures that~$\Diamond_n^\star$ is not polytopal, although it satisfies all necessary conditions to be \poly{4}topal (we skip this discussion which is similar to that in Example~\ref{nonpolytopal:exm:marc+antonio}).
\end{example}


\section{Polytopality of products of graphs}\label{nonpolytopal:sec:product}

\index{Cartesian product!--- of graphs|hbf}
Define the \defn{Cartesian product}~$G\times H$ of two graphs~$G$~and~$H$ to be the graph with vertex set $V(G\times H) \eqdef V(G)\times V(H)$, and edge set $E(G\times H) \eqdef \big(V(G)\times E(H)\big)\cup\big(E(G)\times V(H)\big)$. In other words, for~$a,c\in V(G)$ and~$b,d\in V(H)$, the vertices~$(a,b)$ and~$(c,d)$ of~$G\times H$ are adjacent if either~$a=c$ and~$\{b,d\}\in E(H)$, or~$b=d$ and~$\{a,c\}\in E(G)$. Notice that this product is usually denoted by~$G\Box H$ in graph theory. We choose to use the notation~$G\times H$ to be consistent with the Cartesian product of polytopes: if~$G$ and~$H$ are the graphs of the polytopes~$P$ and~$Q$ respectively, then the product~$G\times H$ is the graph of the product~$P\times Q$. In this section, we focus on the polytopality of products of non-polytopal graphs.

\begin{figure}[b]
	\capstart
	\centerline{\includegraphics[scale=1]{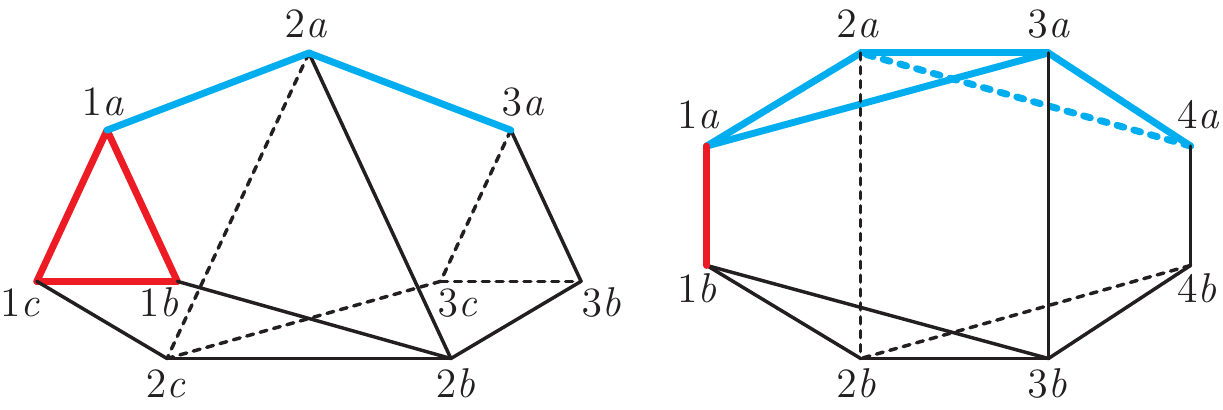}}
	\caption[Polytopal products of non-polytopal graphs]{Polytopal products of non-polytopal graphs: the product of a triangle~$abc$ by a path~$123$ (left) and the product of a segment~$ab$ by two glued triangles~$123$ and~$234$ (right).}
	\label{nonpolytopal:fig:ctrm}
\end{figure}

As already mentioned, the factors of a polytopal product are not necessarily polytopal: consider for example the product of a triangle by a path, or the product of a segment by two glued triangles (see \fref{nonpolytopal:fig:ctrm} and more generally Proposition~\ref{nonpolytopal:prop:subdiv}). We neutralize these elementary examples by further requiring the product~$G\times H$, or equivalently the factors~$G$ and~$H$, to be regular (the degree of a vertex~$(v,w)$ of~$G\times H$ is the sum of the degrees of the vertices~$v$ of~$G$ and~$w$ of~$H$). In this case, it is natural to investigate when such regular products can be simply polytopal. The answer is given by Theorem~\ref{nonpolytopal:theo:simplypolytopalproducts}.

Our study of polytopality of products of graphs was inspired by G\"unter~Ziegler's prototype question: ``is the Cartesian product of two Petersen graphs polytopal?'' In fact, observe that we incidentally already answered this question in the case of dimension~$6$ in Corollary~\ref{nonpolytopal:coro:cycles}: the product of two Petersen graphs cannot be simply polytopal since it contains two induced \cycle{5}s which share three vertices (more generally, Theorem~\ref{nonpolytopal:theo:simplypolytopalproducts} characterizes simple polytopality of products). However, we have no answer for dimensions~$4$ and~$5$: the combinatorial approach we use is not sufficient to prove non-polytopality since there exists a pseudo-manifold whose graph is the product of two Petersen graphs (see Proposition~\ref{nonpolytopal:prop:pseudomanifold}).

Before starting, let us observe that the necessary conditions of Proposition~\ref{nonpolytopal:prop:necessaryconditions} are  preserved under Cartesian products in the following sense:

\begin{proposition}
If two graphs~$G$ and~$H$ are respectively $d$- and \connected{e}, and respectively satisfy $d$- and $e$-PSP, then their product~$G\times H$ is \connected{(d+e)} and satisfies $(d+e)$-PSP. 
\end{proposition}

\begin{proof}
The connectivity of a Cartesian product of graphs was studied in~\cite{cs-ccpg-99}. In fact, it is even proved in~\cite{s-ccpg-08} that $\kappa(G\times H)=\min(\kappa(G)|H|,\kappa(H)|G|,\delta(G)+\delta(H))\ge\kappa(G)+\kappa(H)$,
where~$\kappa(G)$ and~$\delta(G)$ respectively denote the connectivity and the minimum degree of a graph~$G$.

For the principal subdivision property, consider a vertex~$(v,w)$ of~$G\times H$. Choose a principal subdivision of~$K_{d+1}$ in~$G$ with principal vertex~$v$ and neighbors~$N_v$, and a principal subdivision of~$K_{e+1}$ in~$H$ with principal vertex~$w$ and neighbors~$N_w$. This gives rise to a principal subdivision of~$K_{d+e+1}$ in $G\times H$ with principal vertex~$(v,w)$ and neighbors~$(N_v\times\{w\})\cup(\{v\}\times N_w)$. Indeed, for~$x,x'\in N_v$, the vertices~$(x,w)$ and~$(x',w)$ are connected by a path in~$G\times w$ by construction; similarly, for~$y,y'\in N_w$, the vertices~$(v,y)$ and~$(v,y')$ are connected by a path in~$v\times H$. Finally, for each~$x\in N_v$ and~$y\in N_w$, connect~$(x,w)$ to~$(v,y)$ via the path of length~$2$ that passes through~$(x,y)$. All these paths are disjoint by construction.
\end{proof}


\subsection{Simply polytopal products}

\index{polytopal!simply ---}
A product of simply polytopal graphs is automatically simply polytopal. We prove that the reciprocal statement is also true:

\begin{theorem}\label{nonpolytopal:theo:simplypolytopalproducts}
A product of graphs is simply polytopal if and only if its factors are.
\end{theorem}

Applying Theorem~\ref{nonpolytopal:theo:kalai}, we obtain a strong characterization of the simply polytopal~\mbox{products}:

\begin{corollary}\label{nonpolytopal:coro:simplypolytopalproducts}
Products of simple polytopes are the only simple polytopes whose graph is a product.\qed
\end{corollary}

Let~$G$ and~$H$ be two connected regular graphs of degree~$d$~and~$e$ respectively, and assume that the graph~$G\times H$~is the graph of a simple \poly{(d+e)}tope~$P$. By Proposition~\ref{nonpolytopal:prop:cycles}, for all edges~$a$ of~$G$ and~$b$ of~$H$, the \cycle{4} $a\times b$ is the graph of a \face{2} of~$P$.

\begin{observation}\label{nonpolytopal:obs:propagate}
Let $F$~be any facet of $P$, let $v$~be a vertex of~$G$, and let~$\{x,y\}$~be an edge of~$H$ such that $(v,x)\in F$ and $(v,y)\notin F$. Then~$G\times\{x\}\subset F$ and~$G\times\{y\}\cap F=\emptyset$.
\end{observation}

\begin{proof}
Since the polytope is simple, all neighbors of~$(v,x)$ except~$(v,y)$ are connected to~$(v,x)$ by an edge of~$F$. Let~$v'$ be a neighbor of~$v$ in~$G$, and let~$C$ be the \face{2}~$\conv\{v,v'\}\times\conv\{x,y\}$ of~$P$. If~$(v',y)$ were a vertex of~$F$, the intersection~$C\cap F$ would consist of exactly three vertices (because~$(v,y)\notin F$), a contradiction. In summary,~$(v',x)\in F$ and~${(v',y)\notin F}$, for all neighbors~$v'$ of~$v$. Repeating this argument and using the fact that~$G$~is connected yields~$G\times\{x\}\subset F$ and~$G\times\{y\}\cap F=\emptyset$.   
\end{proof}

\begin{lemma}\label{nonpolytopal:lem:facets}
The graph of any facet of~$P$ is either of the form~$G'\times H$ for a \regular{(d-1)} induced subgraph~$G'$ of~$G$, or of the form~$G\times H'$ for an \regular{(e-1)} induced subgraph~$H'$~of~$H$.
\end{lemma}

\begin{proof}
Assume that the graph of a facet~$F$ is not of the form~$G'\times H$.  Then there exists a vertex~$v$ of~$G$ and an edge~$\{x,y\}$ of~$H$ such that~$(v,x)\in F$ and~$(v,y)\notin F$. By Observation~\ref{nonpolytopal:obs:propagate}, the subgraph~$H'$ of~$H$ induced by the vertices~$y\in H$ such that~$G\times \{y\}\subset F$ is nonempty. We now prove that the graph~$\gr(F)$ of~$F$ is exactly~$G\times H'$. 

The inclusion~$G\times H'\subset\gr(F)$ is clear: by definition,~$G\times\{y\}$ is a subgraph of~$\gr(F)$ for any vertex~$y\in H'$. For any edge~$\{x,y\}$ of~$H'$ and any vertex~$v\in G$, the two vertices~$(v,x)$ and~$(v,y)$ are contained in~$F$, so the edge between them is an edge of~$F$; if not, we would have an improper intersection between~$F$ and this edge. 

For the other inclusion, let~$H'' \eqdef \ens{y\in H}{G\times\{y\}\cap F = \emptyset}$ and~$H''' \eqdef H \ssm(H'\cup H'')$. If $H'''\ne\emptyset$, the fact that~$H$~is connected ensures that there is an edge between some vertex of~$H'''$ and either a vertex of~$H'$ or~$H''$. This contradicts Observation~\ref{nonpolytopal:obs:propagate}.

We have proved that~$G\times H'=\gr(F)$. The fact that~$F$ is a simple \poly{(d+e-1)}tope and the \regular{d}ity of~$G$ together ensure that~$H'$~is \regular{(e-1)}.
\end{proof}

\begin{proof}[Proof of Theorem~\ref{nonpolytopal:theo:simplypolytopalproducts}]
One direction is clear. For the other direction, proceed by induction on~$d+e$, the cases~$d=0$ and~$e=0$ being trivial.  Now assume that~$d,e\ge1$, that~$G\times H=\gr(P)$, and that $G$~is not the graph of a \poly{d}tope. By Lemma~\ref{nonpolytopal:lem:facets}, all facets of~$P$ are of the form~$G'\times H$ or~$G\times H'$, where~$G'$ (resp.~$H'$) is an induced \regular{(d-1)} (resp.~\regular{(e-1)}) subgraph of~$G$ (resp.~$H$). By induction, the second case does not arise. We fix a vertex~$w$ of~$H$. Then induction tell us that~$F_w \eqdef G'\times\{w\}$ is a face of~$P$, and~$G'\times H$ is the only facet of~$P$ that contains~$F_w$ by Lemma~\ref{nonpolytopal:lem:facets}. This cannot occur unless~$F_w$~is a facet, but this only happens in the base case~$H=\{w\}$.
\end{proof}



\begin{corollary}
Consider a graph~$G$ \regular{d}, \connected{d}, satisfying $d$-PSP, but not simply \poly{d}topal. Then, any product of~$G$ by a simply \poly{e}topal graph is \regular{(d+e)}, \connected{(d+e)}, satisfies $(d+e)$-PSP, but is not simply \poly{(d+e)}topal.
\end{corollary}

\begin{example}
The product of the circulant graph~$C_8(1,4)$ by the graph of the \dimensional{d} cube is a non simply polytopal graph for non-trivial reasons. For any~$m\ge 4$, the product of the circulant graph~$C_{2m}(1,m)$ by a segment is non-polytopal for non-trivial reasons.
\end{example}


\subsection{Polytopal products of non-polytopal graphs}

\index{subdivision!regular ---}
In this section, we give a general construction to obtain polytopal products starting from a polytopal graph~$G$ and a non-polytopal one~$H$. We need the graph~$H$ to be the graph of a \defn{regular subdivision} of a polytope~$Q$, that is, the graph of the upper\footnote{The unusual convention to define a subdivision as the projection of the upper facets of the lifting simplifies the presentation of the construction.} envelope (the set of all upper facets with respect to the last coordinate) of the convex hull of the point set~$\ens{(q,\omega(q))}{q\in V(Q)}\subset\R^{e+1}$ obtained by lifting the vertices of~$Q\subset\R^e$ according to a \defn{lifting function}~$\omega:V(Q)\to\R$.

\begin{proposition}\label{nonpolytopal:prop:subdiv}
If~$G$ is the graph of a \poly{d}tope~$P$, and~$H$ is the graph of a regular subdivision of an \poly{e}tope~$Q$, then~$G\times H$~is \poly{(d+e)}topal. In the case~$d>1$, the regular subdivision of~$Q$ can even have internal vertices.
\end{proposition}

\begin{proof}
Let $\omega:V(Q)\to\R_{>0}$ be a lifting function that induces a regular subdivision of~$Q$ with graph~$H$. Assume without loss of generality that the origin of~$\R^d$ lies in the interior of~$P$. For each~$p\in V(P)$ and~$q\in V(Q)$, we define the point~$\rho(p,q) \eqdef (\omega(q)p,q)\in\R^{d+e}$. Consider
$$R \eqdef \conv\ens{\rho(p,q)}{p\in V(P), q\in V(Q)}.$$
Let~$g$ be a facet of~$Q$ defined by the linear inequality~$\dotprod{\psi}{y}\le 1$. Then~$\dotprod{(0,\psi)}{(x,y)}\le 1$ defines a facet of~$R$, with vertex set~$\ens{\rho(p,q)}{p\in P, q\in g}$, and isomorphic to~$P\times g$.

Let~$f$ be a facet of~$P$ defined by the linear inequality~$\dotprod{\phi}{x}\le 1$. Let~$c$ be a cell of the subdivision of~$Q$, and let~$\psi_0 h + \dotprod{\psi}{y}\le 1$ be the linear inequality that defines the upper facet corresponding to~$c$ in the lifting. Then we claim that the linear inequality
$$\chi(x,y) = \psi_0\dotprod{\phi}{x}+\dotprod{\psi}{y} \le 1$$
selects a facet of~$R$ with vertex set~$\ens{\rho(p,q)}{p\in f, q\in c}$ that is isomorphic to~$f\times c$. Indeed, 
$$\chi\big(\rho(p,q)\big) = \chi(\omega(q)p,q) = \psi_0\omega(q)\dotprod{\phi}{p}+\dotprod{\psi}{q} \le 1,$$
where equality holds if and only if~$\dotprod{\phi}{p}=1$ and~$\psi_0\omega(q)+\dotprod{\psi}{q}=1$, so that~$p\in f$ and~$q\in c$.

The above set~$\cF$ of facets of~$R$ in fact contains all facets: indeed, any \face{(d+e-2)} of a facet in~$\cF$ is contained in precisely two facets in~$\cF$. Since the union of the edge sets of the facets in~$\cF$ is precisely~$G\times H$, it follows that the graph of~$R$ equals~$G\times H$.

A similar argument proves the same statement in the case when~$d>1$ and~$H$ is a regular subdivision of~$Q$ with internal vertices (meaning that not only the vertices of~$Q$ are lifted, but also a finite number of interior points).
\end{proof}

We already mentioned two examples obtained by such a construction in the beginning of this section (see \fref{nonpolytopal:fig:ctrm}): the product of a polytopal graph by a path and the product of a segment by a subdivision of an \gon{n} with no internal vertex. Proposition~\ref{nonpolytopal:prop:subdiv} even produces examples of regular polytopal products which are not simply polytopal:

\begin{example}\label{nonpolytopal:exm:diamond-subdiv}

\begin{figure}[b]
	\capstart
	\centerline{\includegraphics[width=.85\textwidth]{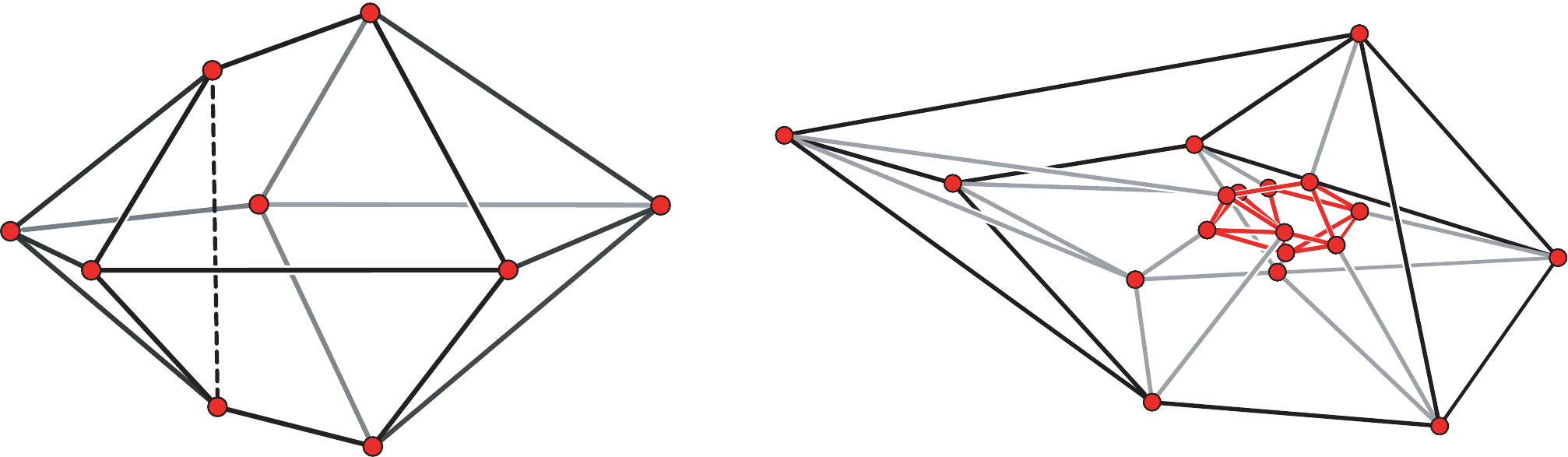}}
	\caption[A non-polytopal \regular{4} graph~$H$ and the Schlegel diagram of a \poly{4}tope whose graph is the product of~$H$ by a segment]{A non-polytopal \regular{4} graph~$H$ which is the graph of a regular subdivision of a \poly{3}tope (left) and the Schlegel diagram of a \poly{4}tope whose graph is the product of~$H$ by a segment (right).}
	\label{nonpolytopal:fig:truncatedoctahedron}
\end{figure}

Let~$H$ be the graph obtained by a star-clique operation from the graph of an octahedron. It is non-polytopal (Corollary~\ref{nonpolytopal:coro:starclique}), but it is the graph of a regular subdivision of a \poly{3}tope (see \fref{nonpolytopal:fig:truncatedoctahedron}). Consequently, the product of~$H$ by any regular polytopal graph is polytopal. Thus, there exist regular polytopal products which are not simply polytopal.
\end{example}

Finally, Proposition~\ref{nonpolytopal:prop:subdiv} also produces polytopal products of two non-polytopal graphs:

\begin{example}[Product of dominos]\index{domino graph}
Define the \defn{\domino{p} graph}~$D_p$ to be the product of a path~$P_p$ of length~$p$ by a segment. Let~$p,q\ge 2$. Observe that~$D_p$ and~$D_q$ are not polytopal and that~$D_p\times P_q$ is a regular subdivision of a \poly{3}tope. Consequently, the product of dominos ${D_p\times D_q}$ is a \poly{4}topal product of two non-polytopal graphs (see Figure~\ref{nonpolytopal:fig:domino}). 

Finally, let us observe that the product~$D_p\times D_q=P_p\times P_q\times (K_2)^2$ can be decomposed in different ways into a product of two graphs. However, in any such decomposition, at least one of the factors is non-polytopal.

\begin{figure}[h]
	\capstart
	\centerline{\includegraphics[scale=.6]{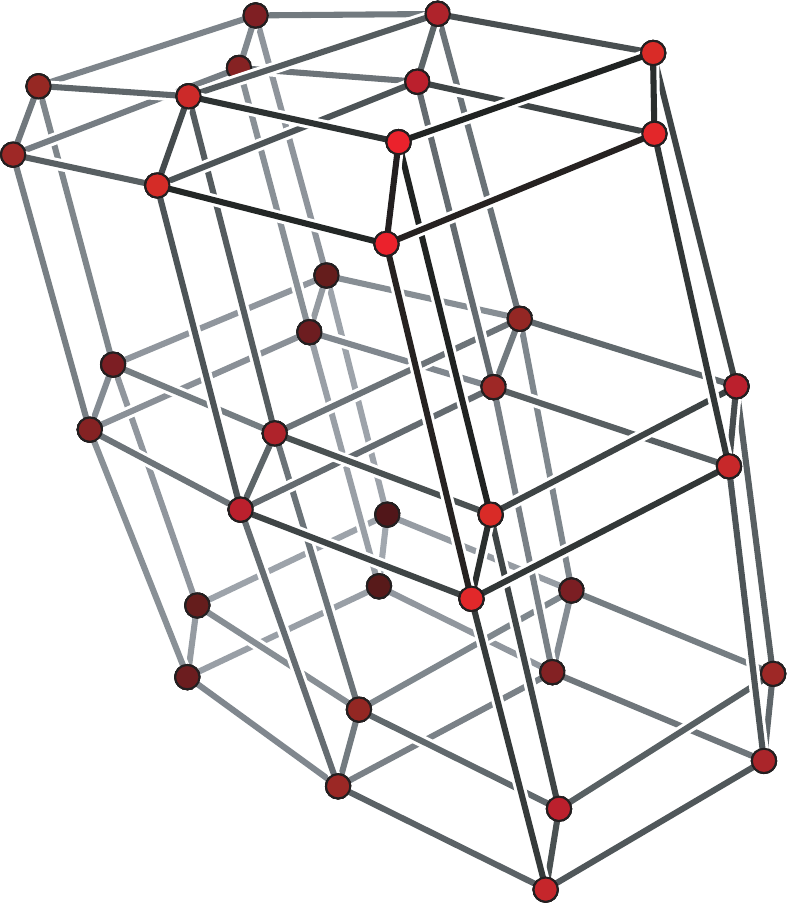} \qquad \includegraphics[scale=.5]{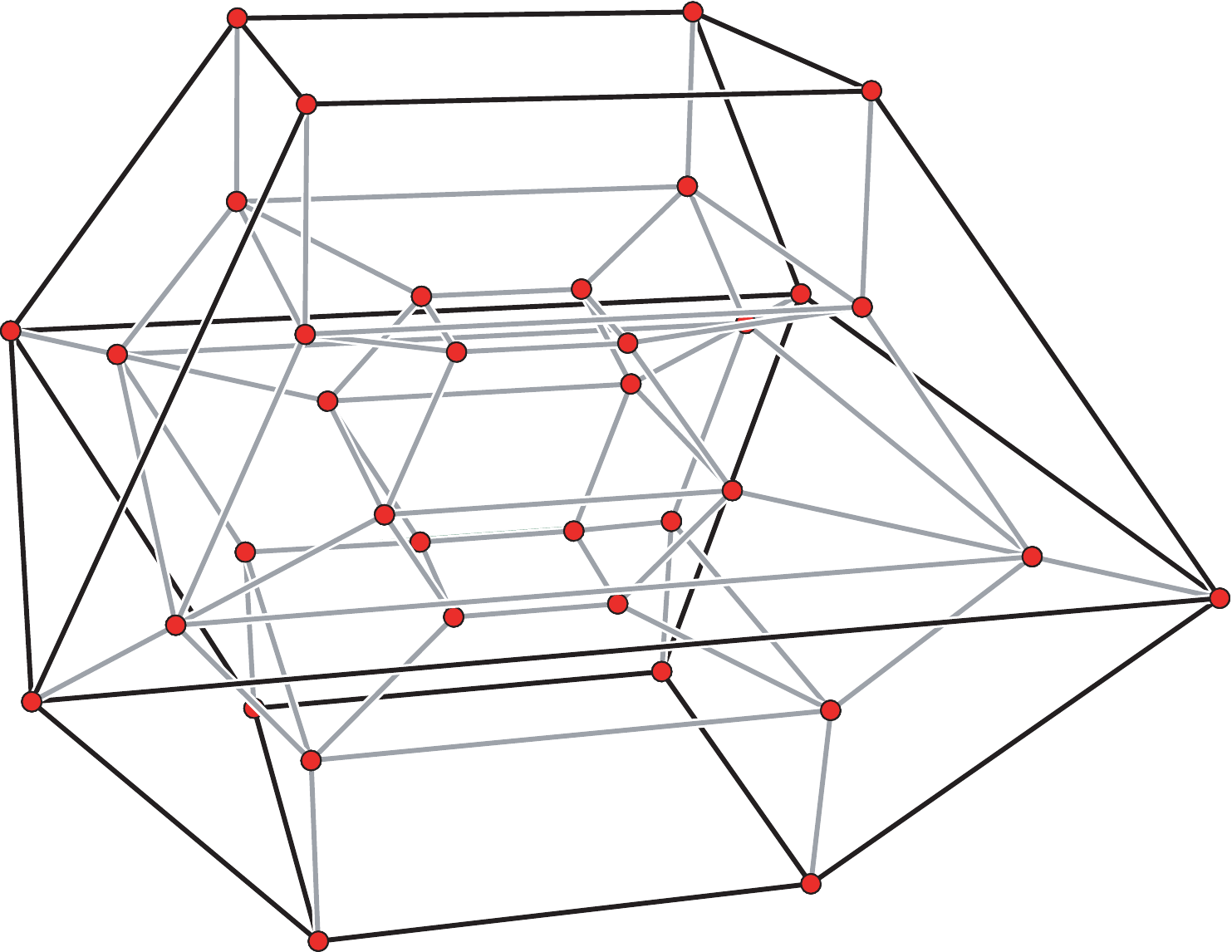}}
	\caption[The graph of the product of two dominos and the Schlegel diagram of a realizing \poly{4}tope]{The graph of the product of two \domino{2}s (left) and the Schlegel diagram of a realizing \poly{4}tope (right).}
	\label{nonpolytopal:fig:domino}
\end{figure}
\end{example}


\subsection{Product with a segment}

In this section, we complete our list of examples of products of a segment by a regular graph~$H$. The goal is to illustrate all possible behaviors of such a product regarding polytopality:
\begin{enumerate}
\item If~$H$ is polytopal, then~$K_2\times H$ is polytopal. However, some ambiguities can appear:
\begin{enumerate}
\item either the dimension is ambiguous. For example,~$K_2\times K_n$ is realized by the product of a segment by any neighborly polytope. See Chapter~\ref{chap:psn} for a discussion on dimensional ambiguity of products of complete graphs.
\item or the dimension is unambiguous, but there exist realizations in the same dimension which are not combinatorially equivalent. In this case, $H$~is not simply polytopal (Theorem~\ref{nonpolytopal:theo:kalai}). In Proposition~\ref{nonpolytopal:prop:prismoctahedron}, we determine all possible realizations of the graph of a prism over an octahedron.
\item or there is no ambiguity at all. It is the case for example if~$H$ is simply \poly{3}topal.
\end{enumerate}
\item If~$H$ is not polytopal, then~$K_2\times H$ is not simply polytopal (Theorem~\ref{nonpolytopal:theo:simplypolytopalproducts}). However:
\begin{enumerate}
\item either~$K_2\times H$ is polytopal in smaller dimension (Example~\ref{nonpolytopal:exm:diamond-subdiv}).
\item or~$K_2\times H$ is not polytopal at all. It is the case when~$H$ is the complete graph~$K_{n,n}$ (Proposition~\ref{nonpolytopal:prop:knnxk2}) or when~$H$ is non-polytopal and \regular{3} (Proposition~\ref{nonpolytopal:prop:GxK2}).
\end{enumerate}
\end{enumerate}

\begin{proposition}\label{nonpolytopal:prop:knnxk2}
For~$n\ge3$, the graph~$K_2\times K_{n,n}$ is not polytopal.
\end{proposition}

To prove this proposition, we will need the following well known lemma:

\begin{lemma}\label{nonpolytopal:lem:8verts}
A \poly{3}tope with no triangular facet has at least~$8$~vertices.
\end{lemma}

\begin{proof}
Let~$P$ be a \poly{3}tope. For~$k\ge3$, denote by~$v_k$ the number of vertices of degree~$k$ and by~$p_k$ the number of \face{2}s with~$k$~vertices. By double counting and Euler's Formula (see \cite[Chapter~13]{g-cp-03} for details),
$$v_3+p_3 = 8+\sum_{k\ge5} (k-4)(v_k+p_k).$$
The lemma immediately follows.
\end{proof}

\begin{proof}[Proof of Proposition~\ref{nonpolytopal:prop:knnxk2}]
Observe that~$K_2\times K_{n,n}$ is not \poly{d}topal for~$d\le3$ because it contains a~$K_{3,3}$-minor, and for~$d=n+1$ by Theorem~\ref{nonpolytopal:theo:simplypolytopalproducts}.

The proof proceeds by contradiction. Suppose that~$K_2\times K_{n,n}$ is the graph of a \poly{d}tope~$P$, for some~$d$ with~$3\le d\le n$, and consider a \face{3}~$F$ of~$P$. Since~$K_2\times K_{n,n}$~contains no triangle, Lemma~\ref{nonpolytopal:lem:8verts} says that $F$~has at least $8$~vertices. Denote by~$A$~and~$B$ the two maximal independent sets in~$K_{n,n}$, and by~$A_0,B_0,A_1,B_1$ their corresponding copies in the cartesian product~$K_2\times K_{n,n}$. We discuss the possible repartition of the vertices of~$F$ in these sets.

Assume first that~$F$~has at least three vertices in~$A_0$; let~$x,y,z$ be three of them. Then it cannot have more than two vertices in~$B_0$, because otherwise its graph would contain a copy of~$K_{3,3}$. In fact, there must be exactly two vertices~$u,v$ in~$B_0$: since any vertex of~$F$ has degree at least~$3$, and each vertex in~$A_0$ can only be connected to vertices in~$B_0$ or to its corresponding neighbor in~$A_1$, each vertex of~$F$ in~$A_0$ must have at least, and thus exactly, two neighbors in~$B_0$ and one in~$A_1$. Thus,~$F$ also has at least three vertices in~$A_1$, and by the same reasoning, there must be exactly two vertices in~$B_1$; call one of them~$w$. But now~$\{x,y,z\}$ and~$\{u,v,w\}$ are the two maximal independent sets of a subdivision of~$K_{3,3}$ included in~$F$.

By symmetry and Lemma~\ref{nonpolytopal:lem:8verts}, $F$~has exactly two vertices in each of the sets $A_0,B_0,A_1,B_1$. Since all these vertices have degree~$3$, we have proved that~$P$'s only $3$-faces are combinatorial cubes whose graphs are cartesian products of~$K_2$ with \cycle{4} in~$A_0\cup B_0=K_{n,n}$. However, this \cycle{4} is not contained in any other \face{3}, which is an obstruction to the existence of~$P$.
\end{proof}

\begin{proposition}\label{nonpolytopal:prop:GxK2}
If~$H$ is a non-polytopal and \regular{3} graph, then~$K_2\times H$ is non-polytopal.
\end{proposition}

\proof
We distinguish two cases:
\begin{enumerate}[(i)]
\item If~$H$ contains a~$K_4$-minor, then~$K_2\times H$ is not \poly{3}topal because it contains a~$K_5$-minor, and it is not \poly{4}topal by Theorem~\ref{nonpolytopal:theo:simplypolytopalproducts}. Since~$K_2\times H$ is \regular{4}, these are the only possibilities.
\item Otherwise,~$H$ is a series-parallel graph. Thus, it can be obtained from~$K_2$ by a sequence of \defn{series} and \defn{parallel} extensions, \ie subdividing or duplicating an edge. Since duplicating an edge creates a double edge, and subdividing an edge yields a vertex of degree two, $H$~is either not simple or not \regular{3}.\qed
\end{enumerate}

To complete our collection of examples of products with a segment, we examine the possible realizations of the graph of the prism over the octahedron:

\begin{proposition}\label{nonpolytopal:prop:prismoctahedron}
The graph of the prism over the octahedron is realized by exactly four combinatorially different polytopes.
\end{proposition}

In order to exhibit four different realizations, we need to remember the context and the proof of Proposition~\ref{nonpolytopal:prop:subdiv}. Given the graph~$G$ of a \poly{d}tope~$P$ and the graph~$H$ of a regular subdivision of an \poly{e}tope~$Q$ defined by a lifting function~$\omega:V(Q)\to\R$, we construct a \poly{(d+e)}tope with graph~$G\times H$ as follows: we start from the product~$P\times Q$ and we deform each face~$\{p\}\times Q$ according to the lifting function~$\omega$. This subdivides the face~$\{p\}\times Q$ and makes appear the subgraph~$\{p\}\times H$ of the product~$G\times H$. We want to observe now that the deformation can be different at each vertex of~$P$: we can use at each vertex of~$P$ a different lifting function, which will produce combinatorially different polytopes.

To illustrate this, we come back to the prism over the octahedron. Denote by~$H$ the graph of the octahedron. Observe that the octahedron has four regular subdivisions with no additional edges: the octahedron itself (for a constant lifting function) and the three subdivisions of the octahedron into two egyptian pyramids glued along their square face (for a lifting function positive on one vertex, negative on the opposite one, and vanishing on the other four vertices). This leads to four combinatorially different realizations of~$K_2\times H$: in our previous construction, we can choose either the octahedron at both ends of the segment (thus obtaining the prism over the octahedron), or the octahedron at one end and the glued square pyramid at the other, or the glued square pyramids for both ends of the segment (and this leads to two possibilities according on whether we choose the same square or two orthogonal squares to subdivide the octahedron).

\begin{figure}[h]
	\capstart
	\centerline{\includegraphics[scale=.55]{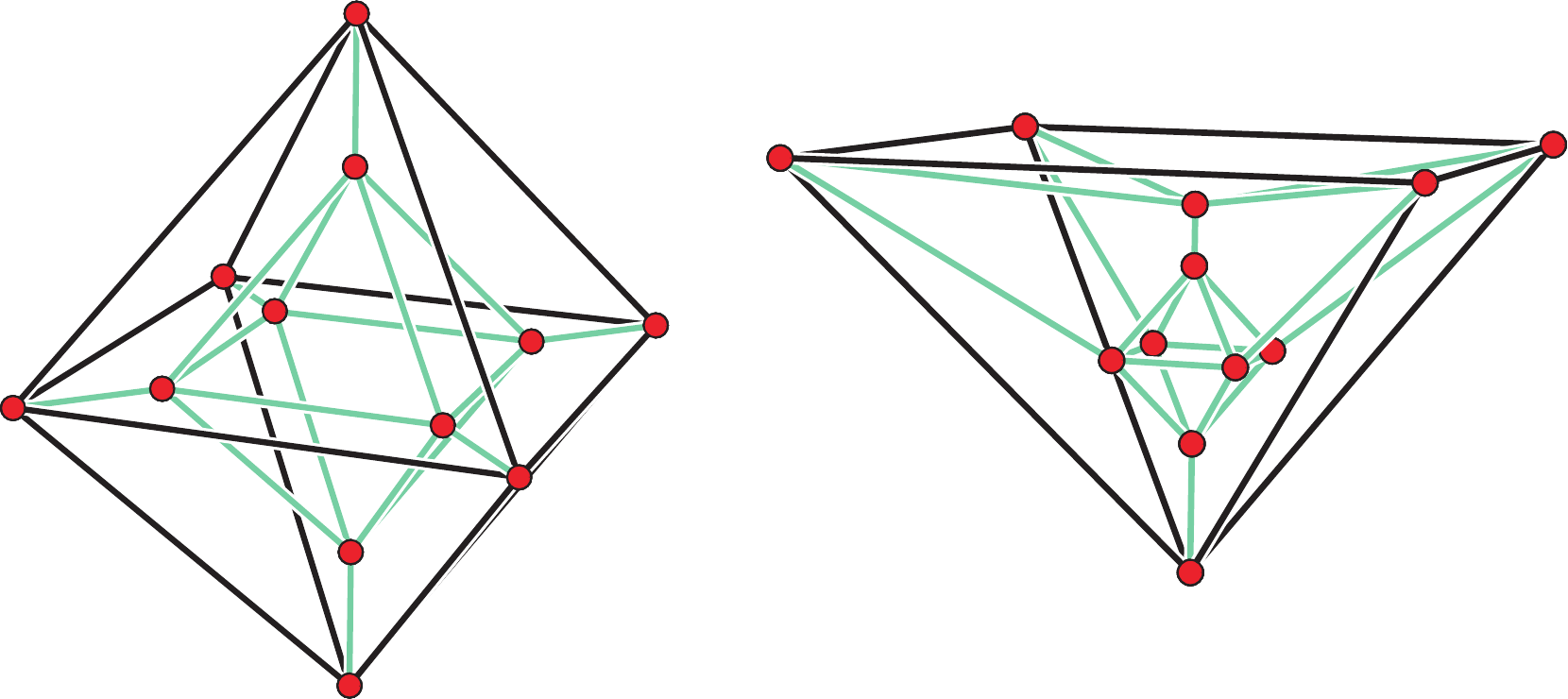}}
	\caption[The prism over the octahedron and a combinatorially different polytope with the same graph]{The prism over the octahedron (left) and a combinatorially different polytope with the same graph (right).}
	\label{nonpolytopal:fig:prismoctahedron}
\end{figure}

In fact, by the same argument, we can even slightly improve Proposition~\ref{nonpolytopal:prop:subdiv}:

\begin{observation}
Let~$G$ be the graph of a \poly{d}tope~$P$ and~$H$ be the graph of an \poly{e}tope~$Q$. For all~$v\in G$, choose a lifting function~$\omega_v:V(Q)\to\R$, and denote by~$H_v$ the graph of the corresponding regular subdivision of~$Q$. Then the graph obtained by replacing in~$G\times H$ the subgraph~$\{v\}\times H$ by~$\{v\}\times H_v$ is polytopal.

The result remains true if~$d>1$ and all subdivisions of~$Q$ have the same interior points.
\end{observation}

It remains to prove that any realization of~$K_2\times H$ is combinatorially equivalent to one of the four described above.  First, the dimension is unambiguous:~$K_2\times H$ can only be \poly{4}topal (by Remark~\ref{nonpolytopal:remark:3poly} and Theorem~\ref{nonpolytopal:theo:simplypolytopalproducts}). In particular, any realization has the following property:

\begin{definition}
\index{polytope!almost simple ---}
A \poly{d}tope is \defn{almost simple} if its graph is \regular{(d+1)}.
\end{definition}

The vertex figures of a simple polytope are all simplices, which implies that any two incident edges in a simple polytope lie in a common \face{2}. For almost simple polytopes, the vertex figures are almost as restricted: they are \defn{$(d-1)$-circuits}, that is, \poly{(d-1)}topes with $d+1$ vertices. This implies the following property:

\begin{proposition}\label{nonpolytopal:prop:almostsimple}
Let~$\{v,w\}$ be an edge of an almost simple \poly{d}tope $P$. Then:
\begin{enumerate}[(a)]
\item either~$\{v,w\}$ together with any other edge incident to~$v$ forms a \face{2};
\item or there exists exactly one more edge incident to~$v$ which does not form a \face{2} with~$\{v,w\}$, but any two \face{2}s both incident to~$\{v,w\}$ lie in a \face{3}.
  \end{enumerate}
\end{proposition}

\begin{proof}
Consider the vertex figure~$F_v$ of~$v$. It is a \poly{(d-1)}tope with~$d+1$ vertices, one of which, say~$\overline{w}$, corresponds to the edge~$\{v,w\}$. This vertex~$\overline{w}$ can be adjacent to either~$d$ or~$d-1$ vertices of~$F_v$. The first case corresponds to statement~(a). In the second case, $\overline{w}$ has exactly one missing edge in~$F_v$ (corresponding to a missing \face{2} in~$P$), but the edge figure of~$\{v,w\}$ is a \simp{(d-2)}. This implies statement~(b).
\end{proof}

With this in mind, we can finally prove Proposition~\ref{nonpolytopal:prop:prismoctahedron}:

\begin{proof}[Proof of Proposition~\ref{nonpolytopal:prop:prismoctahedron}]
We first introduce some notations: let~$V \eqdef \{1,\bar1,2,\bar2,3,\bar3\}$ denote the $6$ vertices of the octahedron such that~$\{1,\bar1\}$, $\{2,\bar2\}$ and $\{3,\bar3\}$ are the three missing edges, and let $a$ and $b$ denote the two endpoints of the segment factor. We denote the vertices of $K_2\times H$ by~$\{1a,1b,\bar1a,\dots,\bar3b\}$. We call \defn{horizontal edges} the edges of the form~$\{ia,ib\}$, for $i\in V$, and \defn{vertical edges} the edges of the form~$\{ix,jx\}$, for~$i\in V$, $j\in V\ssm\{i,\bar i\}$ and~$x\in\{a,b\}$.

We first use ``almost simplicity'' to study the possible \face{2}s of a realization~$P$ of~$K_2\times H$. Assume that there exists a \face{2} $F$ which is neither a triangle nor a square. It has to contain an angle between a horizontal edge and a vertical edge, say without loss of generality~$\{1a,1b\}$ and~$\{1a,2a\}$. By inducedness, the next edges of~$F$ are necessarily~$\{2a,\bar1a\}$ and~$\{\bar1a,\bar1b\}$. Since the edges~$\{1a,1b\}$ and~$\{1a,2a\}$ form an angle of~$F$, the two edges~$\{2a,2b\}$ and~$\{1b,2b\}$ cannot form an angle: otherwise the \cycle{4}~$(1a,1b,2b,2a)$ would form a square face which improperly intersects~$F$. Similarly, since the edges~$\{2a,\bar1a\}$ and~$\{\bar1a,\bar1b\}$ form an angle, the edges~$\{2a,2b\}$ and~$\{2b,\bar1b\}$ cannot form an angle. Thus,~$\{2a,2b\}$ is adjacent to two missing angles, which is impossible by Proposition~\ref{nonpolytopal:prop:almostsimple}. We conclude that the \face{2}s of any realization of~$K_2\times H$ can only be squares and triangles. 

We now use this information on the \face{2}s to understand the possible \face{3}s of~$P$. Assume first that none of the angles of the \cycle{4}s~$(1a,2a,\bar1a,\bar2a)$, $(1a,3a,\bar1a,\bar3a)$, and~$(2a,3a,\bar2a,\bar3a)$ forms a \face{2}. Then for each~$i\in V$, the vertex~$ia$ has already two missing angles. Consequently, the remaining angles necessarily form a \face{2} of~$P$ by Proposition~\ref{nonpolytopal:prop:almostsimple}. By inducedness, we obtain all the triangles of the $a$-copy of $H$, and any two adjacent of these triangles are contained in a common \face{3}. This \face{3} is necessarily an octahedron.

Assume now that one of the angles of the \cycle{4}s~$(1a,2a,\bar1a,\bar2a)$, $(1a,3a,\bar1a,\bar3a)$, and $(2a,3a,\bar2a,\bar3a)$ forms a \face{2}. By symmetry, we can suppose that it is the angle defined by the edges~$\{1a,2a\}$ and~$\{2a,\bar1a\}$. Let~$F$ denote the corresponding \face{2} of~$P$. By inducedness, the last vertex of~$F$ cannot be either~$3a$ or~$\bar3a$, and~$F$ is necessarily the square~$(1a,2a,\bar1a,\bar2a)$. It is now easy to see that none of the angles of the \cycle{4}~$(1a,3a,\bar1a,\bar3a)$ (resp.~$(2a,3a,\bar2a,\bar3a)$) can be an angle of a \face{2} of $P$: otherwise, this \cycle{4} would be a \face{2} of~$P$ (by a symmetric argument), which would intersect improperly with~$F$. All together, this implies that the vertices~$3a$ and~$\bar3a$ both have already two missing angles, and thus, that all the other angles form \face{2}s by Proposition~\ref{nonpolytopal:prop:almostsimple}. Furthermore, any two \face{2}s adjacent to an edge~$\{3a,ia\}$, with~$i\in\{1,\bar1,2,\bar2\}$, form a \face{3}. This implies that all angles adjacent to a vertex~$ia$, except the angles of the \cycle{4}s~$(1a,3a,\bar1a,\bar3a)$ and~$(2a,3a,\bar2a,\bar3a)$ form a \face{2}.

Since the two above cases can occur independently at both ends of the segment~$(a,b)$, we obtain the claimed result.
\end{proof}


\subsection{Topological products}

To finish, we come back to G\"unter Ziegler's motivating question: ``is the product of two Petersen graphs polytopal?'' We proved in Theorem~\ref{nonpolytopal:theo:simplypolytopalproducts} that it is not~$6$ polytopal, but the question remains open in dimension~$4$ and~$5$.

\begin{proposition}\label{nonpolytopal:prop:pseudomanifold}
The product of two Petersen graphs is the graph of a \dimensional{4} pseudo-manifold of Euler characteristic~$1$.
\end{proposition}

\begin{figure}[h]
	\capstart
	\centerline{\includegraphics[scale=1]{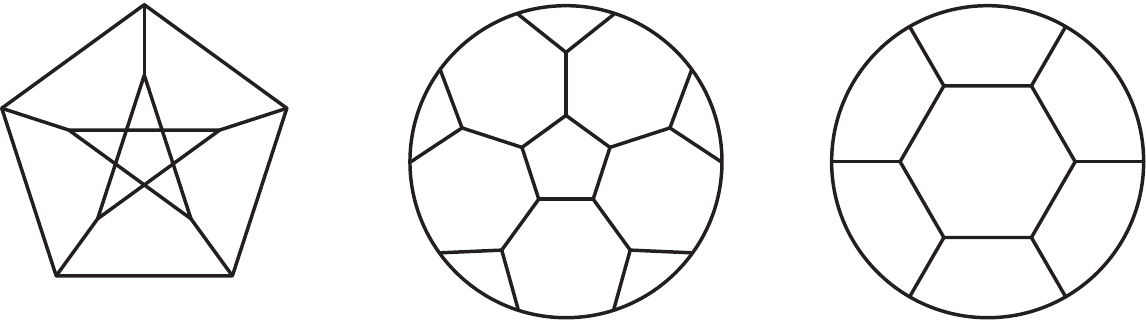}}
	\caption[The Petersen graph and~$K_{3,3}$ embedded on the projective plane]{The Petersen graph (left), its embedding on the projective plane (middle) and the embedding of~$K_{3,3}$ on the  projective plane (right). Antipodal points on the circle are identified.}
	\label{nonpolytopal:fig:projectiveembeddings}
\end{figure}

\begin{proof}
\index{projective!--- plane}
\index{plane!projective ---}
The Petersen graph is the graph of a cellular decomposition of the projective plane~$\cP$ with~$6$ pentagons (see \fref{nonpolytopal:fig:projectiveembeddings}). Consequently, the product of two Petersen graphs is the graph of a cellular decomposition of~$\cP\times\cP$, whose Euler characteristic is~$1$ (the Euler characteristic of a product is the product of the Euler characteristics of the factors). The maximal cells of this decomposition are~$36$ products of two pentagons.
\end{proof}

This proposition tells that understanding the possible \face{4}s of a realization, and their possible incidence relations (as we did for example in Proposition~\ref{nonpolytopal:prop:knnxk2}) is not enough to disprove the polytopality of the product of two Petersen graphs. Observe that the same remark holds for the product of any graphs of cellular decompositions of manifolds: for example, the product of a triangle by the Petersen graph is the graph of a \dimensional{3} pseudo-manifold of Euler characteristic~$0$.

Another interesting example is the product of~$K_{3,3}$ by a triangle. Indeed, in contrast with $K_{3,3}\times K_2$, the graph~$K_{3,3}\times K_3$ is the graph of a \dimensional{3} pseudo-manifold of Euler characteristic~$0$. To see this, embed~$K_{3,3}$ in the projective plane~$\cP$ as in \fref{nonpolytopal:fig:projectiveembeddings}, and multiply this embedding by a triangle. The only difference with Proposition~\ref{nonpolytopal:prop:pseudomanifold} is that the central hexagon in the embedding of Figure~\ref{nonpolytopal:fig:projectiveembeddings} intersects improperly the three external squares. Consequently, in the product with the triangle, each of the three hexagonal prisms intersects improperly three cubes. This can be solved by replacing the chain of three hexagonal prisms by a chain of six triangular prisms with the same boundary. We leave the details to the reader.

	\chapter{Prodsimplicial neighborly polytopes}\label{chap:psn}

\index{Cartesian product!--- of polytopes}
In this dissertation, we focus on polytopal realizations of Cartesian products of graphs (or more generally of complexes). In Chapter~\ref{chap:nonpolytopal}, we mainly studied the existence of a realization. In this chapter, we investigate the possible dimension of a realization of a product. For this purpose, the simplest products are already interesting (and challenging) enough:

\begin{definition}
\index{polytope!psn@\xpsn{(k,\sub{n})} ---|hbf}
Let~$k\ge0$ and~$\sub{n} \eqdef (n_1,\dots,n_r)$, with~$r\ge1$ and~$n_i\ge 1$ for all~$i$. A polytope is \defn{$(k,\sub{n})$-prodsimplicial-neighborly}~---~or \defn{\xpsn{(k,\sub{n})}} for short~---~if its \skeleton{k} is combinatorially equivalent to that of the product of simplices~$\simplex_{\sub{n}} \eqdef \simplex_{n_1}\times\cdots\times\simplex_{n_r}$.
\end{definition}

\begin{remark}
This definition is essentially motivated by two particular classes of \psn polytopes:
\begin{enumerate}[(i)]
\item \defn{neighborly} polytopes arise when $r=1$. In the literature, a polytope is \defn{\neighborly{k}} if any subset of at most~$k$ of its vertices forms a face. Note that such a polytope is \xpsn{(k-1,n)} with our notation.
\item \defn{neighborly cubical} polytopes~\cite{jz-ncp-00,js-ncps-07,sz-capdp} arise when~$\sub{n}=(1,1,\dots,1)$.
\end{enumerate}
\index{polytope!neighborly ---|hbf}
\index{neighborly polytope|hbf}
\index{polytope!neighborly cubical ---|hbf}
\index{neighborly cubical polytope|hbf}
\end{remark}

In this chapter, we investigate the possible dimension of a \xpsn{(k,\sub{n})} polytope. Since the product~$\simplex_{\sub{n}}$ itself is \xpsn{(k,\sub{n})}, the highest dimension is obviously~$\sum n_i$. On the other side, we denote by~$\delta(k,\sub{n})$ the smallest possible dimension that a \xpsn{(k,\sub{n})} polytope can have.

In order to construct \xpsn{(k,\sub{n})} polytopes, a natural approach is to project the product~$\simplex_{\sub{n}}$ (or a combinatorially equivalent polytope) onto a smaller subspace in such a way to conserve the \skeleton{k}. This yields the following fundamental subclass of \psn polytopes:

\begin{definition}
\index{polytope!ppsn@\xppsn{(k,\sub{n})} ---|hbf}
We say that a \xpsn{(k,\sub{n})} polytope is \defn{$(k,\sub{n})$-projected-prodsimplicial-neighborly} ---~or \defn{\xppsn{(k,\sub{n})}} for short~---~if it is a projection of a polytope combinatorially equivalent to~$\simplex_{\sub{n}}$.
\end{definition}

We denote by~$\delta_{pr}(k,\sub{n})$ the smallest possible dimension of a \xppsn{(k,\sub{n})} polytope.

\svs
In this chapter, we present upper and lower bounds on this minimal possible dimension. For the upper bounds, we present two methods to construct \psn polytopes:
\begin{enumerate}
\item In Section~\ref{psn:sec:cyclic}, we give \defn{explicit} constructions based on Minkowski sums of \mbox{cyclic polytopes}.
\item In Section~\ref{psn:sec:deformedproducts}, we extend Raman Sanyal and G\"unter Ziegler's technique of \defn{projecting deformed products of polygons}~\cite{z-ppp-04,sz-capdp} to products of arbitrary simple polytopes, and apply it to products of simplices.
\end{enumerate}
To obtain lower bounds, we apply in Section~\ref{psn:sec:topologicalObstruction} the \defn{topological obstruction method} developed by Raman Sanyal to bound the number of vertices in a Minkowski sum.  In view of these obstructions, our constructions in the first part turn out to be optimal for a wide range of parameters.


\section{Constructions from cyclic polytopes}\label{psn:sec:cyclic}

\index{polytope!cyclic ---}
Remember that~$C_d(n) \eqdef \conv\ens{(t_i,t_i^2,\dots,t_i^d)^T}{i\in[n]}$ denotes the cyclic \poly{d}tope with~$n$ vertices. Cyclic polytopes and their products yield our first examples of \psn polytopes:

\begin{example}
For any integers~$k\ge0$ and~$n\ge2k+2$, the cyclic polytope~$C_{2k+2}(n+1)$ is \xppsn{(k,n)}: it is \psn by Proposition~\ref{intro:prop:cyclic}, and \ppsn because any polytope with~$n+1$ vertices is a projection of the \simp{n}.
\end{example}

\begin{example}\label{psn:exm:prodcyclic}
For any~$k\ge0$ and~$\sub{n} \eqdef (n_1,\dots,n_r)$ with~$r\ge1$ and~$n_i\ge1$ for all~$i$, define $I \eqdef \ens{i\in[r]}{n_i\ge2k+3}$. Then the product
$$\prod_{i\in I}C_{2k+2}(n_i+1)\times\prod_{i\notin I}\simplex_{n_i}$$
is a \xppsn{(k,\sub{n})} polytope of dimension~$(2k+2)|I|+\sum_{i\notin I} n_i$ (which is smaller than~$\sum_{i\in[r]} n_i$ when~$I$ is not empty). Consequently:
$$\delta(k,\sub{n}) \le \delta_{pr}(k,\sub{n}) \le (2k+2)|I|+\sum_{i\notin I} n_i.$$
\end{example}


\subsection{Reflections of cyclic polytopes}\label{psn:subsec:cyclic:reflection}

Our next example deals with the special case of the product~$\simplex_1\times\simplex_n$ of a segment by a simplex. Using products of cyclic polytopes as in Example~\ref{psn:exm:prodcyclic}, we can realize the \skeleton{k} of this polytope in dimension~$2k+3$. We can lower this dimension by~$1$ by reflecting the cyclic polytope~$C_{2k+2}(n+1)$ in a well-chosen hyperplane:

\begin{proposition}\label{psn:prop:reflect}
For any~$k\ge0$,~$n\ge2k+2$ and~$\lambda\in\R$ sufficiently large, the polytope
$$P \eqdef \conv\big(\ens{(t_i,\dots,t_i^{2k+2})^T}{i\in[n+1]}\cup\ens{(t_i,\dots,t_i^{2k+1},\lambda-t_i^{2k+2})^T}{i\in[n+1]}\big)$$
is a \xpsn{(k,(1,n))} polytope of dimension~$2k+2$.
\end{proposition}

\begin{proof}
The polytope~$P$ is obtained as the convex hull of two copies of the cyclic polytope $C_{2k+2}(n+1)$. The first one~$Q \eqdef \conv\ens{\mu_{2k+2}(t_i)}{i\in[n+1]}$ lies on the moment curve~$\mu_{2k+2}$, while the second one is obtained as a reflection of ~$Q$ with respect to a hyperplane that is orthogonal to the last coordinate vector~$e_{2k+2}$ and sufficiently far away. During this process:
\begin{enumerate}[(i)]
\item We destroy all the faces of~$Q$ only contained in upper facets of~$Q$.
\item We create prisms over faces of~$Q$ that lie in at least one upper and one lower facet of~$Q$. In other words (see Definition~\ref{psn:def:spf}), we create prisms over the faces of~$Q$ strictly preserved under the orthogonal projection~$\pi:\R^{2k+2}\to\R^{2k+1}$ with kernel~$\R e_{2k+2}$.
\end{enumerate}

The projected polytope~$\pi(Q)$ is nothing but the cyclic polytope~$C_{2k+1}(n+1)$. Since this polytope is \neighborly{k}, any \face{(\le k-1)} of~$Q$ is strictly preserved by~$\pi$, and thus, we take a prism over all \face{(\le k-1)}s of~$Q$.

Thus, in order to complete the proof that the \skeleton{k} of~$P$ is that of~$\simplex_1\times\simplex_n$, it is enough to show that any \face{k} of~$Q$ remains in~$P$. This is obviously the case if this \face{k} is also a \face{k} of~$C_{2k+1}(n+1)$, and follows from the next combinatorial lemma otherwise.
\end{proof}

\begin{lemma}
A \face{k} of~$C_{2k+2}(n+1)$ which is not a \face{k} of~$C_{2k+1}(n+1)$ is only contained in lower facets of~$C_{2k+2}(n+1)$.
\end{lemma}

\begin{proof}
Let~$F\subset[n+1]$ be a \face{k} of~$C_{2k+2}(n+1)$ that is contained in at least one upper facet~$G\subset[n+1]$ of~$C_{2k+2}(n+1)$. By Proposition~\ref{intro:prop:gale}, we have $n+1\in G$. Thus:
\begin{enumerate}[(i)]
\item If~$n+1\notin F$, then~$G\ssm\{n+1\}$ is a facet of~$C_{2k+1}(n+1)$ containing~$F$.
\item If~$n+1\in F$, and~$F' \eqdef F\ssm\{n+1\}$ has only~$k$ elements. Thus,~$F'$~is a face of~$C_{2k}(n)$, and can be completed to a facet of~$C_{2k}(n)$. Adding the index~$n+1$ back to this facet, we obtain a facet of~$C_{2k+1}(n+1)$ containing~$F$.
\end{enumerate}
In both cases, we have shown that~$F$ is a \face{k} of~$C_{2k+1}(n+1)$.
\end{proof}

\begin{remark}
Observe that the previous construction is similar to our construction of polytopal products of non-polytopal graphs of Proposition~\ref{nonpolytopal:prop:subdiv}, except that we now care not only about the graph, but about a bigger skeleton. The crucial property (presented in the previous lemma) is that there exists a regular subdivision of the cyclic polytope~$C_{2k+1}(n+1)$ whose \skeleton{k} is combinatorially equivalent to that of the \simp{n}. In particular, this method can be used to construct \xpsn{(k,\sub{n})} polytopes in dimension $(2k+1)|I|+\max\left(1,\sum_{i\notin I} n_i\right)$ where ${I \eqdef \ens{i\in[r]}{n_i\ge2k+3}}$. We do not explicitly present the construction here since our next section provides even better results as soon as~$\sub{n}\ne(1,n)$.
\end{remark}


\subsection{Minkowski sums of cyclic polytopes}\label{psn:subsec:cyclic:minkowskiSumCyclicPolytopes}

Our next examples are Minkowski sums of cyclic polytopes. We first describe an easy construction that avoids all technicalities, but only yields \xppsn{(k,\sub{n})} polytopes in dimension~$2k+2r$. After that, we show how to reduce the dimension to~$2k+r+1$, which according to Corollary~\ref{psn:coro:topObstr} is best possible for large~$n_i$'s.

\begin{proposition}
Let~$k\ge0$ and~$\sub{n} \eqdef (n_1,\dots,n_r)$ with~$r\ge1$ and~$n_i\ge1$ for all~$i$. For any pairwise disjoint index sets~$I_1,\dots, I_r\subset\R$, with~$|I_i|=n_i$ for all~$i$, the polytope 
$$P \eqdef \conv\ens{v_{a_1,\dots,a_r}}{(a_1,\dots,a_r)\in I_1\times\cdots\times I_r} \subset \R^{2k+2r}$$
is \xppsn{(k,\sub{n})}, where~$v_{a_1,\dots,a_r}  \eqdef  \left(\sum_{i\in[r]} a_i,\sum_{i\in[r]} a_i^2,\dots,\sum_{i\in[r]} a_i^{2k+2r}\right)^T\in\R^{2k+2r}$.
\end{proposition}

\begin{proof}
The vertex set of~$\simplex_{\sub{n}}$ is indexed by~$I_1\times\cdots\times I_r$. Let~$A \eqdef A_1\times\cdots\times A_r\subset I_1\times\cdots\times I_r$ define a \face{k} of $\simplex_{\sub{n}}$. Consider the polynomial
$$f(t)  \eqdef  \prod_{i\in[r]}\prod_{a\in A_i} (t-a)^2  \eqdef  \sum_{j=0}^{2k+2r} c_j t^j.$$
Since~$A$~indexes a \face{k} of~$\simplex_{\sub{n}}$, we know that~$\sum |A_i|=k+r$, so that the degree of~$f(t)$ is indeed~$2k+2r$. Since~$f(t)\ge0$, and equality holds if and only if~$t\in\bigcup_{i\in[r]} A_i$, the inner product~$\dotprod{(c_1,\dots,c_{2k+2r})}{v_{a_1,\dots,a_r}}$ equals
$$(c_1,\dots, c_{2k+2r})\begin{pmatrix}\sum_{i\in[r]} a_i\\ \vdots \\ \sum_{i\in[r]} a_i^{2k+2r}\end{pmatrix} = \sum_{i\in[r]}\sum_{j=1}^{2k+2r} c_j a_i^j = \sum_{i\in[r]} \big(f(a_i)-c_0\big) \ge -rc_0,$$
with equality if and only if~$(a_1,\dots,a_r)\in A$. Thus,~$A$~indexes a face of~$P$ defined by the linear inequality~$\sum_{i\in[r]} c_i x_i \ge -rc_0$.
\end{proof}

To realize the \skeleton{k} of~$\simplex_{n_1}\times\cdots\times\simplex_{n_r}$ even in dimension~$2k+r+1$, we slightly modify this construction in the following way:

\begin{theorem}\label{psn:theo:UBminkowskiCyclic}
Let~$k\ge0$ and~$\sub{n} \eqdef (n_1,\dots,n_r)$ with~$r\ge1$ and~$n_i\ge1$ for all~$i$. There exist index sets~$I_1,\dots, I_r\subset\R$, with~$|I_i|=n_i$ for all~$i$, such that the polytope
$$P  \eqdef  \conv\ens{w_{a_1,\dots,a_r}}{(a_1,\dots,a_r)\in I_1\times\cdots\times I_r} \subset \R^{2k+r+1}$$
is \xppsn{(k,\sub{n})}, where $w_{a_1,\dots,a_r}  \eqdef  \left(a_1,\dots,a_r,\sum_{i\in[r]} a_i^2,\dots,\sum_{i\in[r]} a_i^{2k+2}\right)^T\in\R^{2k+r+1}$.
\end{theorem}

\begin{proof}
We will choose the index sets~$I_i$ to be sufficiently separated in a sense that will be made explicit later in the proof. This choice will enable us, for each \face{k}~$F$ of~$\simplex_{\sub{n}}$ indexed by $A_1\times\cdots\times A_r\subset I_1\times\cdots\times I_r$, to construct a \defn{monic} polynomial (\ie a polynomial with leading coefficient equal to~$1$)
$$f_F(t)  \eqdef  \sum_{j=0}^{2k+2} c_jt^j$$
that has, for all~$i\in[r]$, the form
$$f_F(t) = Q_i(t)\prod_{a\in A_i}(t-a)^2  +s_it+r_i,$$
with polynomials~$Q_i(t)$ everywhere positive and certain reals~$r_i$ and~$s_i$. From the coefficients of this polynomial, we build the vector
$$n_F  \eqdef  (s_1-c_1,\dots,s_r-c_1,-c_2,-c_3,\dots,-c_{2k+2}) \in \R^{2k+r+1}.$$
To prove that~$n_F$ is a face-defining normal vector for~$F$, take an arbitrary \tuple{r}~$(a_1,\dots,a_r)$ in ${I_1\times\cdots\times I_r}$, and consider the following inequality for the inner product:
\begin{eqnarray*}
\dotprod{n_F}{w_{a_1,\dots,a_r}}
&=& \sum_{i\in[r]}\left(s_i a_i-\sum_{j=1}^{2k+2} c_ja_i^j\right) = \sum_{i\in[r]}\left(s_ia_i+c_0-f_F(a_i)\right) \\
&=& \sum_{i\in[r]}\left(c_0-Q_i(a_i)\prod_{a\in A_i}(a_i-a)^2-r_i\right) \le rc_0-\sum_{i\in[r]} r_i.
\end{eqnarray*}
Equality holds if and only if~$(a_1,\dots,a_r)\in A_1\times\cdots\times A_r$. Given the existence of a polynomial~$f_F$ with the claimed properties, this proves that~$A_1\times\cdots\times A_r$ indexes all~$w_{a_1,\dots,a_r}$'s that lie on a face~$F'$ in~$P$, and they of course span~$F'$ by definition of~$P$. To prove that~$F'$ is combinatorially equivalent to~$F$ it suffices to show that each~$w_{a_1,\dots,a_r}\in F'$ is in fact a vertex of~$P$, since~$P$ is a projection of~$\simplex_{\sub{n}}$. This can be shown with the normal vector~$(2a_1,\dots,2a_r,-1,0,\dots,0)$, using the same calculation as before.

Before showing how to choose the index sets~$I_i$ that enable us to construct the polynomials~$f_F$ in general, we would like to make a brief aside to show the smallest example.
\end{proof}

\begin{example}
For~$k=1$ and~$r=2$, choose the index sets~$I_1$,~$I_2\subset\R$ arbitrarily, but separated in the sense that the largest element of~$I_1$ be smaller than the smallest element of~$I_2$. For any \face{1}~$F$ of~$P$ indexed by~$\{a,b\}\times\{c\}\subset I_1\times I_2$, consider the polynomial~$f_F$ of degree~$2k+2=4$:
$$f_F(t)  \eqdef  (t-a)^2(t-b)^2 = (t^2+\alpha t+\beta)(t-c)^2 + s_2t + r_2,$$
where
\vspace{-.4cm}
$$\begin{array}{ccl}
\alpha & \!\!\!\!\!\eqdef\!\!\!\!\! & 2(-a-b+c), \\ 
r_2    & \!\!\!\!\!\eqdef\!\!\!\!\! &  a^2b^2 - \beta c^2,
\end{array}
\qquad
\begin{array}{ccl}
\beta  & \!\!\!\!\!\eqdef\!\!\!\!\! & a^2 + b^2 + 3c^2 + 4ab - 4ac - 4bc, \\
s_2    & \!\!\!\!\!\eqdef\!\!\!\!\! &  -2a^2b - 2ab^2 - \alpha c^2 + 2\beta c.
\end{array}
$$
Since the index sets~$I_1$,~$I_2$ are separated, the discriminant~$\alpha^2-4\beta = -8(c-a)(c-b)$ is negative, which implies that the polynomial~$Q_2(t)=t^2+\alpha t+\beta$ is positive for all values of~$t$.
\end{example}

\begin{proof}[Proof of Theorem~\ref{psn:theo:UBminkowskiCyclic}, continued]
We still need to show how to choose the index sets~$I_i$ that enable us to construct the polynomials~$f_F$ in general. Once we have chosen these index sets, finding~$f_F$ is equivalent to the task of finding polynomials~$Q_i(t)$ such that:

\begin{enumerate}[(i)]
\item $Q_i(t)$~is monic of degree~$2k+2-2|A_i|$.
\item The~$r$ polynomials~$f_i(t) \eqdef Q_i(t)\prod_{a\in A_i}(t-a)^2$ are equal up to possibly the coefficients in front of~$t^0$ and~$t^1$.
\item $Q_i(t)>0$~for all~$t\in\R$.
\end{enumerate}

The first two items form a linear equation system on the coefficients of the~$Q_i(t)$'s which has the same number of equations as variables. We will show that it has a unique solution if one chooses the right index sets~$I_i$ (the third item will be dealt with at the end). To do this, choose pairwise distinct reals~$\bar a_1,\dots,\bar a_r\in\R$ and look at the similar equation system:

\enlargethispage{.75cm}
\begin{enumerate}[(i)]
\item $\bar Q_i(t)$~are monic polynomials of degree~$2k+2-2|A_i|$.
\item The~$r$ polynomials~$\bar f_i(t) \eqdef \bar Q_i(t)(t-\bar a_i)^{2|A_i|}$ are equal up to possibly the coefficients in front of~$t^0$ and~$t^1$.
\end{enumerate}

The first equation system moves into the second when we deform the points of the sets~$A_i$ continuously to~$\bar a_i$, respectively. If we show that the second equation system has a unique solution then so has the first equation system as long as we have chosen the sets~$I_i$ close enough to the~$\bar a_i$'s, by continuity of the determinant (note that in the end, we can fulfill all these closeness conditions required for all \face{k}s of~$\simplex_{\sub{n}}$ since there are only finitely many \face{k}s).

Note that a polynomial~$\bar f_i(t)$ of degree~$2k+2$ has the form 
\begin{equation}\label{psn:eq:fstarhasform}
\bar Q_i(t)(t-\bar a_i)^{2|A_i|}+ s_it+r_i
\end{equation}
for a monic polynomial~$\bar Q_i$ and some reals~$s_i$ and~$r_i$, if and only if~$\bar f''_i(t)$ has the form
\begin{equation}\label{psn:eq:fprimeprimestarhasform}
R_i(t)(t-\bar a_i)^{2(|A_i|-1)}
\end{equation}
for a polynomial~$R_i(t)$ with leading coefficient~$(2k+2)(2k+1)$. The backward direction can be settled by assuming without loss of generality~$\bar a_i=0$ (otherwise just make a variable shift $(t-\bar a_i)\mapsto t$) and then integrating \eqref{psn:eq:fprimeprimestarhasform} twice with integration constants zero to obtain \eqref{psn:eq:fstarhasform}.

Therefore the second equation system is equivalent to the following third one:

\begin{enumerate}[(i)]
\item $R_i(t)$~are polynomials of degree~$2k-2(|A_i|-1)$ with leading coefficient~$(2k+2)(2k+1)$.
\item The~$r$ polynomials~$g_i(t) \eqdef R_i(t)(t-\bar a_i)^{2(|A_i|-1)}$ all equal the same polynomial, say~$g(t)$.
\end{enumerate}

Since~$\sum_i 2(|A_i|-1)=2k$, it has the unique solution
$$R_i(t) = (2k+2)(2k+1)\prod_{j\neq i} (t-\bar a_j)^{2(|A_j|-1)},$$
with
$$g(t) = (2k+2)(2k+1)\prod_{j\in[r]} (t-\bar a_j)^{2(|A_j|-1)}.$$

Therefore the second system also has a unique solution, where the~$\bar f_i(t)$ are obtained by integrating~$g_i(t)$ twice with some specific integration constants. For a fixed~$i$ we can again assume~$\bar a_i=0$. Then both integration constants have been zero for this~$i$, hence~$\bar f_i(0)=0$ and~$\bar f'_i(0)=0$. Since~$g_i$ is non-negative and zero only at isolated points,~$\bar f_i$ is strictly convex, hence non-negative and zero only at~$t=0$. Therefore~$\bar Q_i(t)$ is positive for~$t\neq 0$. Since we chose~$\bar a_i=0$, we can quickly compute the correspondence between the coefficients of ${\bar Q_i(t)=\sum_j \bar q_{i,j}t^j}$ and of~$R_i(t)=\sum_j r_{i,j}t^j$:
$$r_{i,j} = \big(2|A_i|(2|A_i|-1)+4j|A_i|+j(j-1)\big)\bar q_{i,j}.$$
In particular 
$$\bar Q_i(0)=\bar q_{i,0}=\frac{r_{i,0}}{2|A_i|(2|A_i|-1)}=\frac{R_i(0)}{2|A_i|(2|A_i|-1)}>0,$$
therefore~$\bar Q_i(t)$ is everywhere positive, hence so is~$Q_i(t)$ if one chooses~$I_i$ possibly even closer to~$\bar a_i$, since the solutions of linear equation systems move continuously when one deforms the entries of the equation system by a homotopy (as long as the determinant stays non-zero), since the determinant and taking the adjoint matrix are continuous maps. The positivity of~$Q_i(t)$ finishes the proof.
\end{proof}


\section{Projections of deformed products}\label{psn:sec:deformedproducts}

In the previous section, we saw an explicit construction of polytopes whose \skeleton{k} is equivalent to that of a product of simplices. In this section, we provide another construction of \xppsn{(k,\sub{n})} polytopes, using Raman Sanyal and G\"unter Ziegler's technique of \defn{projecting deformed products of polygons}~\cite{z-ppp-04,sz-capdp} and generalizing it to products of arbitrary simple polytopes. This generalized technique consists in suitably projecting suitable polytopes combinatorially equivalent to a given product of simple polytopes in such a way as to preserve its complete \skeleton{k}. The case of products of simplices then yields \xppsn{(k,\sub{n})} polytopes.


\subsection{Projection of deformed product of simple polytopes}\label{psn:subsec:deformedproducts:general}

We first discuss the general setting: for any given product~$P \eqdef P_1\times\cdots\times P_r$ of simple polytopes, we construct 
a polytope $\defP$ combinatorially equivalent to~$P$ and whose \skeleton{k} is preserved under the projection on the first~$d$ coordinates.

\paragraph{Deformed products of simple polytopes.}

\index{deformed products of simple polytopes}
Let~$P_1,\dots,P_r$ be \defn{simple} polytopes of respective dimensions~$n_1,\dots, n_r$ and facet descriptions~$P_i = \ens{x\in\R^{n_i}}{A_i x\le b_i}$, where each real matrix $A_i\in\R^{m_i\times n_i}$~has one row for each of the~$m_i$~facets
of~$P_i$ and~$n_i=\dim P_i$ many columns, and~$b_i$~is a right-hand side vector in~$\R^{m_i}$. The product~$P \eqdef P_1\times\cdots\times P_r$ then has dimension $n \eqdef \sum_{i\in[r]} n_i$, and is given by the~$m \eqdef \sum_{i\in[r]} m_i$ inequalities:
$$\begin{pmatrix} A_1 & & \\ & \ddots & \\ & & A_r \end{pmatrix} x \le \begin{pmatrix} b_1 \\ \vdots \\ b_r \end{pmatrix}.$$
The left hand~$m\times n$ matrix shall be denoted by~$A$. It is proved in~\cite{az-dpmsp-99} that for any matrix~$\defA$ obtained from~$A$ by \defn{arbitrarily} changing the zero entries above the diagonal blocks, there exists a right-hand side~$\defb$ such that the deformed polytope~$\defP$ defined by the inequality system ${\defA x\le \defb}$ is combinatorially equivalent to~$P$. The equivalence is the obvious one: it maps the facet defined by the~$i$th row of~$A$ to the one given by the $i$th row of~$\defA$, for all~$i$. Following the ideas of~\cite{sz-capdp}, we will use this ``deformed product'' construction in such a way that the projection of~$\defP$ to the first~$d$~coordinates preserves its \skeleton{k} in the following sense.

\paragraph{Preserved faces and the Projection Lemma.}

\begin{figure}
	\capstart
	\centerline{\includegraphics[scale=1]{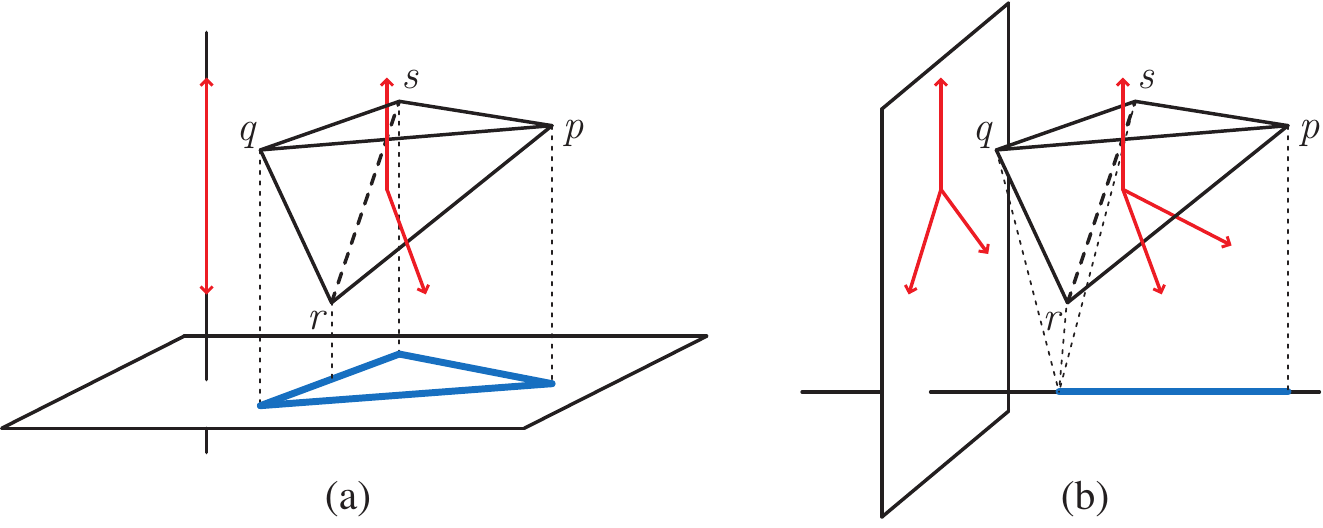}}
	\caption[Preserved faces under a projection]{(a)~Projection of a tetrahedron onto~$\R^2$: the edge~$pq$ is strictly preserved, while neither the edge~$qr$, nor the face~$\textit{qrs}$, nor the edge~$\textit{qs}$ are (because of conditions~(i), (ii)~and~(iii) respectively). (b)~Projection of a tetrahedron to~$\R$: only the vertex~$p$ is strictly preserved.}
	\label{psn:fig:projection}
\end{figure}

For integers~$n>d$, let~$\pi:\R^n\to\R^d$ denote the orthogonal projection to the first~$d$ coordinates, and~$\tau:\R^n\to\R^{n-d}$ denote the dual orthogonal projection to the last~$n-d$ coordinates. Let~$P$ be a full-dimensional simple polytope
in~$\R^n$, with~$0$ in its interior. We consider the following notion of preserved faces (see \fref{psn:fig:projection}):

\begin{definition}[\cite{z-ppp-04}]\label{psn:def:spf}
\index{strictly preserved face}
A proper face~$F$ of a polytope~$P$ is \defn{strictly preserved} under~$\pi$ if:
\begin{enumerate}[(i)]
\item $\pi(F)$~is a face of~$\pi(P)$,
\item $F$~and~$\pi(F)$ are combinatorially isomorphic, and
\item $\pi^{-1}(\pi(F))$~equals~$F$.
\end{enumerate}
\end{definition}

The characterization of strictly preserved faces of~$P$ uses the normal vectors of the facets of~$P$. Let~$F_1,\dots,F_m$ denote the facets of~$P$ and for all~$i\in[m]$, let~$f_i$ denote the normal vector of~$F_i$, and~$g_i \eqdef \tau(f_i)$. For any face~$F$ of~$P$, let~$\varphi(F)$ denote the set of indices of the facets of~$P$ containing~$F$, \ie such that~$F=\bigcap_{i\in \varphi(F)} F_i$.

\begin{lemma}[Projection Lemma~\cite{az-dpmsp-99,z-ppp-04}]\label{psn:lem:projection}
\index{Projection Lemma}
A face~$F$ of the polytope~$P$ is strictly preserved under the projection~$\pi$ if and only if~$\ens{g_i}{i\in \varphi(F)}$ is positively spanning.\qed
\end{lemma}

\paragraph{A first construction.}

\begin{figure}
	\capstart
	\centerline{\includegraphics[scale=1]{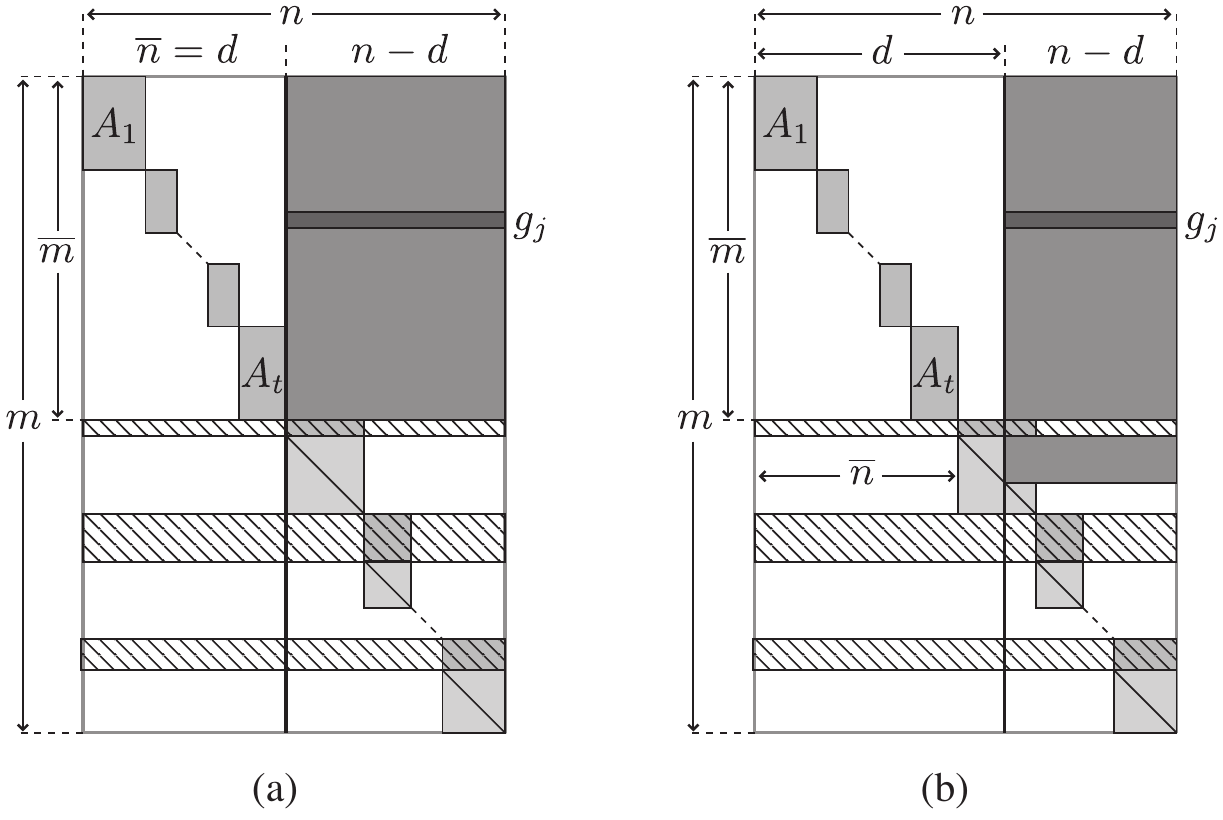}}
	\caption[The deformed matrix~$\defA$, basic case]{The deformed matrix~$\defA$~(a)~when the projection does not slice any block~($\bn=d$), and~(b)~when the block~$A_{t+1}$ is sliced~($\bn<d$).  Horizontal hatched boxes denote bad row vectors. The top right solid block is formed by the vectors~$g_1,\dots,g_\bm$. }
	\label{psn:fig:defA}
\end{figure}

Let~$t\in\{0,1,\dots,r\}$ be maximal such that the matrices~$A_1,\dots, A_t$ are entirely contained in the first~$d$~columns
of~$A$. Let~$\bm \eqdef \sum_{i=1}^t m_i$ and~$\bn \eqdef \sum_{i=1}^t n_i$. By changing bases appropriately, we can assume that the bottom~$n_i\times n_i$ block of~$A_i$ is the identity matrix, for each~$i\ge t+1$. In order to simplify the exposition, we also assume first that~$\bn=d$, \ie that the projection on the first~$d$ coordinates separates the first~$t$ block matrices from the last~$r-t$. See \fref{psn:fig:defA}(a).

Let~$\{g_1,\dots, g_\bm\}\subset\R^{n-d}$ be a set of vectors such that~$G \eqdef \{e_1,\dots,e_{n-d}\}\cup\{g_1,\dots,g_\bm\}$~is the Gale transform of a full-dimensional simplicial neighborly polytope~$Q$~---~see~\cite{z-lp-95,m-ldg-02} for definition and properties of Gale duality. By elementary properties of the Gale transform,~$Q$ has~${\bm+n-d}$ vertices, and~$\dim Q=(\bm+n-d)-(n-d)-1=\bm-1$. In particular, every subset of~$\Fracfloor{\bm-1}{2}$~vertices spans a face of~$Q$, so every subset of~$\bm+n-d-\Fracfloor{\bm-1}{2} \eqfed \alpha$ elements of~$G$ is positively spanning.

We deform the matrix~$A$ into the matrix~$\defA$ of \fref{psn:fig:defA}(a), using the vectors~$g_1,\dots,g_\bm$ to deform the top~$\bm$ rows. We denote by~$\defP$ the corresponding deformed product. We say that a facet of~$\defP$ is \defn{good} if the right part of the corresponding row of~$\defA$ is covered by a vector of $G$, and \defn{bad} otherwise. Bad facets are hatched in \fref{psn:fig:defA}(a). Observe that there are~$\beta \eqdef m-\bm-n+d$ bad facets in total.

Let~$F$ be a \face{k} of~$\defP$. Since~$\defP$ is a simple \poly{n}tope, $F$~is the intersection of~$n-k$~facets, among which at least~$\gamma \eqdef n - k - \beta$ are good facets.  If the corresponding elements of~$G$ are positively spanning, then~$F$~is strictly preserved under projection on the first~$d$~coordinates. Since we have seen that any subset of~$\alpha$~vectors of~$G$ is positively spanning,~$F$~will surely be preserved if~$\alpha\le\gamma$, which is equivalent to:
$$k \le n - m  + \Fracfloor{\bm-1}{2}.$$
Thus, under this assumption, we obtain a \poly{d}tope whose \skeleton{k} is combinatorially equivalent to the one of~$P=P_1\times\cdots\times P_r$.

\paragraph{When the projection slices a block.}

We now discuss the case when~$\bn < d$, for which the method is very similar.  We consider~$\bm + d - \bn$ vectors~$g_1,\dots,g_{\bm + d - \bn}$ such that the set of vectors ${G \eqdef \{e_1,\dots,e_{n-d}\}\cup\{g_1,\dots,g_{\bm + d -\bn}\}}$ is the Gale dual of a neighborly polytope. We deform the matrix~$A$ into the matrix~$\defA$ shown in \fref{psn:fig:defA}(b), using again the vectors~$g_1,\dots,g_\bm$ to deform the top~$\bm$ rows and the vectors~$g_{\bm+1}\dots,g_{\bm + d -\bn}$ to deform the top~$d-\bn$ rows of the~$n_{t+1}\times n_{t+1}$ bottom identity submatrix of~$A_{t+1}$. This is indeed a valid deformation since we can prescribe the~$n_{t+1}\times n_{t+1}$ bottom submatrix of~$A_{t+1}$ to be any upper triangular matrix, up to changing the basis appropriately.  For the same reasons as before:
\begin{enumerate}[(i)]
\item any subset of at least~$\alpha \eqdef \bm + n - \bn - \Fracfloor{\bm + d - \bn - 1}{2}$ elements of~$G$ is positively spanning;
\item the number of bad facets is~$\beta  \eqdef  m-\bm-n+\bn$, and thus any \face{k} of~$\defP$ is contained in at least~$\gamma \eqdef n-k-\beta$ good facets.
\end{enumerate}
Thus, the condition~$\alpha\le\gamma$ translates to:
$$k \le n - m  + \Fracfloor{\bm+d-\bn-1}{2},$$
and we obtain the following proposition:

\begin{proposition}\label{psn:prop:defp1}
Let~$P_1,\dots,P_r$ be simple polytopes of respective dimension~$n_i$, and with~$m_i$ many facets. For a fixed integer~$d\le n$, let~$t$ be maximal such that~$\sum_{i=1}^t n_i\le d$. Then there exists a \poly{d}tope whose \skeleton{k} is combinatorially equivalent to that of the product~$P_1\times\cdots\times P_r$ as soon as
$$0 \le k \le \sum_{i=1}^r (n_i-m_i) + \Floor{\frac{1}{2}\left(d-1+\sum_{i=1}^t(m_i-n_i)\right)}.$$
\end{proposition}
\vspace{-1.1cm}\qed
\vspace{1cm}

In the next two paragraphs, we present two improvements on the bound of this proposition. Both use colorings of the graphs of the polar polytopes~$P_i^\polar$, in order to weaken the condition $\alpha\le\gamma$, in two different directions:
\begin{enumerate}[(i)]
\item the first improvement decreases the number of required vectors in the Gale transform~$G$, which in turn, decreases the value of~$\alpha$;
\item the second one decreases the number of bad facets, and thus, increases the value of~$\gamma$.
\end{enumerate}

\paragraph{Multiple vectors.}

In order to raise our bound on~$k$, we can save vectors of~$G$ by repeating some of them several times. Namely, any two facets that have no \face{k} in common can share the same vector~$g_j$. Since any two facets of a simple polytope containing a common \face{k} share a ridge, this condition can be expressed in terms of incidences in the graph of the polar polytope: facets not connected by an edge in this graph can use the same vector~$g_j$. For a graph~$H$, we denote as usual its chromatic number by~$\chi(H)$. Then, each~$P_i$ with~$i\le t$ only contributes~$\chi_i \eqdef \chi(\gr(P_i^\polar))$ different vectors in~$G$, instead of~$m_i$ of them. Thus, we only need in total~$\bchi \eqdef \sum_{i=1}^t \chi_i$ different vectors~$g_j$. This improvement replaces~$\bm$ by~$\bchi$ in the formula of~$\alpha$, while~$\beta$~and~$\gamma$ do not change, and the condition~$\alpha\le\gamma$ is equivalent to
$$k \le n - m + \bm - \bchi + \Fracfloor{\bchi-d-\bn-1}{2}.$$
Thus, we obtain the following improved proposition:

\vspace{3.6cm}\qed
\vspace{-4.1cm}
\begin{proposition}\label{psn:prop:defp2}
Let~$P_1,\dots,P_r$ be simple polytopes. For each polytope~$P_i$, denote by~$n_i$ its dimension, by~$m_i$ its number of facets, and by~$\chi_i \eqdef \chi(\gr(P_i^\polar))$ the chromatic number of the graph of its polar polytope~$P_i^\polar$. For a fixed integer~$d\le n$, let~$t$ be maximal such that~${\sum_{i=1}^t n_i\le d}$. Then there exists a \poly{d}tope whose \skeleton{k} is combinatorially equivalent to that of the product~$P_1\times\cdots\times P_r$ as soon as
$$0 \le k \le \sum_{i=1}^r (n_i-m_i) + \sum_{i=1}^t (m_i-\chi_i)+\Floor{\frac12\left(d-1+\sum_{i=1}^t(\chi_i-n_i)\right)}.$$
\end{proposition}

\begin{example}\label{psn:exm:even}
\index{polytope!even ---}
Since polars of simple polytopes are simplicial,~$\chi_i\ge n_i$ is an obvious lower bound for the chromatic number of the dual graph of~$P_i$. Polytopes that attain this lower bound with equality are characterized by the property that all their \face{2}s have an even number of vertices, and are called \defn{even} polytopes. 

If all~$P_i$ are even polytopes, then~$\bn=\bchi$, and we obtain a \poly{d}tope with the same \skeleton{k} as~$P_1\times\cdots\times P_r$ as soon as
$$k \le n - m + \bm - \bn + \Fracfloor{d-1}{2}.$$
In order to maximize~$k$, we should maximize~$\bm-\bn$, subject to the condition~$\bn\le d$. For example, if all~$n_i$ are equal, this amounts to ordering the~$P_i$ by decreasing number of facets.
\end{example}

\paragraph{Scaling blocks.}

\begin{figure}
	\capstart
	\centerline{\includegraphics[scale=1]{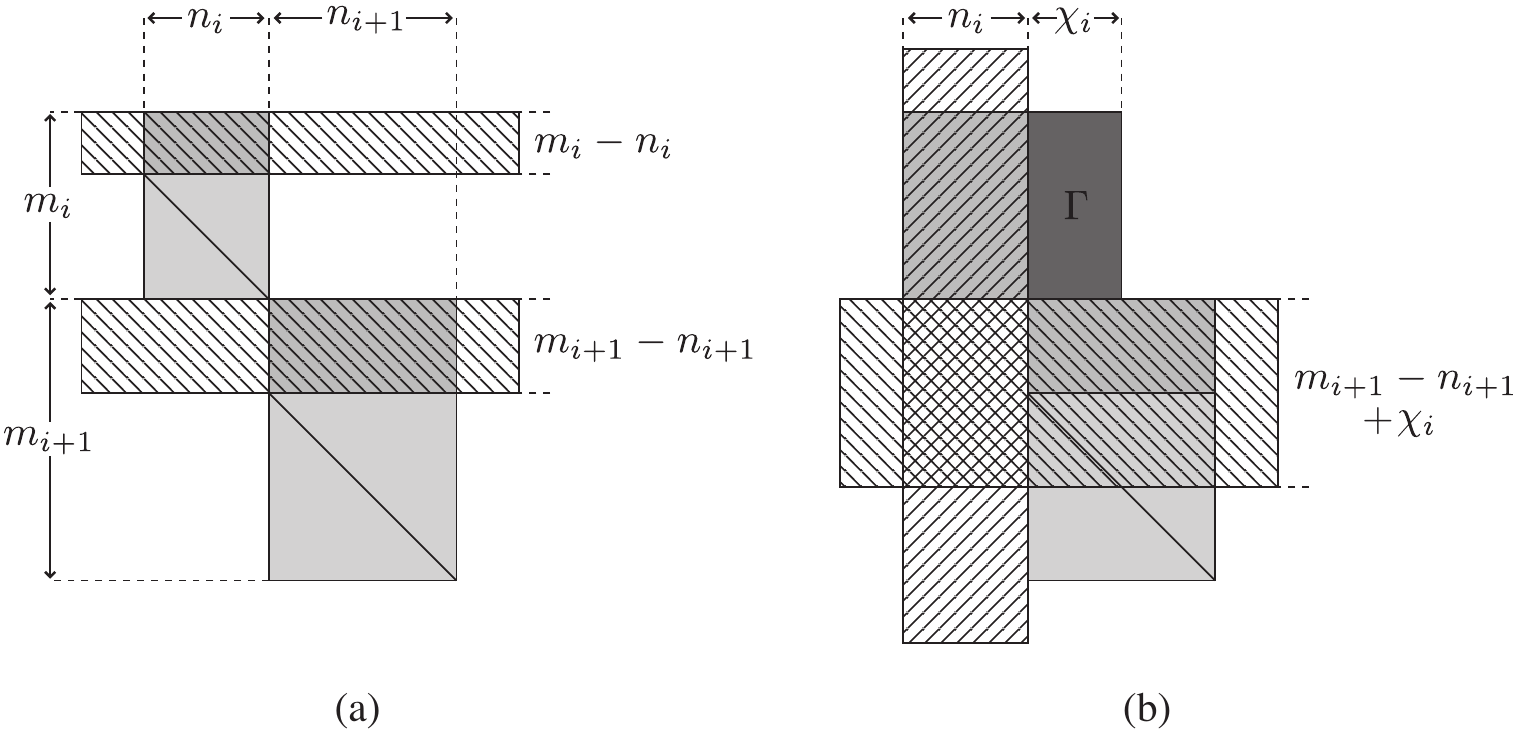}}
	\caption[The deformed matrix~$\defA$, colorings]{How to raise the dimension of the preserved skeleton by inserting the incidence matrix~$\Gamma$ of a coloring of the graph of the polar polytope~$P_i^\polar$. Part~(a) shows the situation before the insertion of~$\Gamma$, and Part~(b) the changes that have occurred.  Bad row vectors and unnecessary columns are hatched. The entries in the matrix to the left of~$\Gamma$ must be rescaled to retain a valid inequality description of~$P$.}
	\label{psn:fig:chi}
\end{figure}

We can also apply colorings to the blocks~$A_i$ with~$i\ge t+1$, by filling in the area below~$G$ and above the diagonal blocks. To explain this, assume for the moment that~$\chi_i\le n_{i+1}$ for a certain fixed~$i\ge t+2$. Assume that the rows of~$A_i$ are colored with~$\chi_i$~colors using a valid coloring~$c:[m_i]\to[\chi_i]$ of the graph of the polar polytope~$P_i^\polar$. Let~$\Gamma$ be the incidence matrix of~$c$, defined by~$\Gamma_{j,k}=1$ if~$c(j)=k$, and~$\Gamma_{j,k}=0$ otherwise. Thus,~$\Gamma$ is a matrix of size~$m_i\times \chi_i$. We put this matrix to the right of~$A_i$ and above~$A_{i+1}$ as in \fref{psn:fig:chi}(b), so that we append the same unit vector to each row of~$A_i$ in the same color class. Moreover, we scale all entries of the block~$A_i$ by a sufficiently small constant~$\varepsilon>0$.

In this setting, the situation is slightly different:
\begin{enumerate}[(i)]
\item In the Gale dual~$G$, we do not need anymore the~$n_i$ basis vectors of~$\R^{n-d}$ hatched in \fref{psn:fig:chi}(b). Let~$a \eqdef \sum_{j<i} n_j$ denote the index of the last column vector of~$A_{i-1}$ and ${b \eqdef 1+\sum_{j\le i} n_j}$ denote the index of the first column vector of~$A_{i+1}$. We define the vector set $G \eqdef \{e_1,\dots,e_{a-d},e_{b-d},\dots,e_{n-d}\}\cup\{g_1,\dots,g_\bm\}$ to be the Gale transform of a simplicial neighborly polytope~$Q$ of dimension~$\bm-1-n_i$. As before, any subset of at least $\alpha \eqdef \bm+n-\bn-n_i-\Fracfloor{\bm+d-\bn-n_i-1}{2}$ vectors of~$G$ positively spans~$\R^{n-d}$.

\item \defn{Bad} facets are defined as before, except that the top~$m_i-n_i$ rows of~$A_i$ are not bad anymore, but all of the first~$m_{i+1}-n_{i+1}+\chi_i$~rows of~$A_{i+1}$ are now bad. Thus, the net change in the number of bad rows is~$\chi_i-m_i+n_i$, so that any \face{k} is contained in at least ${\gamma \eqdef 2n-k-m+\bm-\bn+m_i-n_i-\chi_i}$ good rows. Up to~$\varepsilon$-entry elements, the last $n-d$ coordinates of these rows correspond to pairwise distinct elements of~$G$.
\end{enumerate}

Applying the same reasoning as above, the \skeleton{k} of~$\defP$ is strictly preserved under projection to the first~$d$ coordinates as soon as~$\alpha\le\gamma$, which is equivalent to:
$$k \le n - m + m_i -\chi_i + \Fracfloor{\bm+d-\bn-n_i-1}{2}.$$

Thus, we improve our bound on~$k$ as soon as
$$\Delta \eqdef m_i-\chi_i+\Fracfloor{\bm+d-\bn-n_i-1}{2} - \Fracfloor{\bm+d-\bn-1}{2}>0.$$
For example, this difference~$\Delta$ is big for polytopes whose polars have many vertices but a small chromatic number.

\mvs

\begin{figure}
	\capstart
	\centerline{\includegraphics[scale=1]{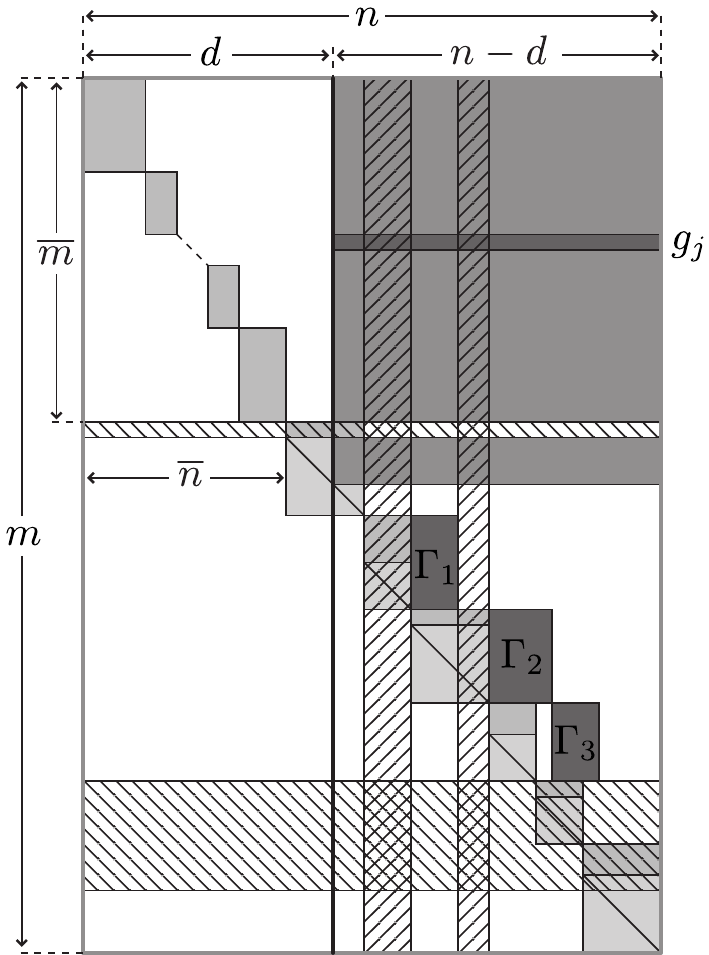}}
	\caption[The deformed matrix~$\defA$, general case]{How to reduce the number of vectors in the Gale transform using various coloring matrices of polar polytopes. Situations where~$\chi_i > n_{i+1}$ can be accommodated for as illustrated by the matrix~$\Gamma_2$ in the picture.}
	\label{psn:fig:matrixfinal}
\end{figure}

Finally, observe that one can apply this ``scaling'' improvement even if~$\chi_i>n_{i+1}$ (except that it will perturb more
than the two blocks~$A_i$ and~$A_{i+1}$) and to more than one matrix~$A_i$. See the example in \fref{psn:fig:matrixfinal}. In this picture, the~$\Gamma$~blocks are incidence matrices of colorings of the graphs of the polar
polytopes. Call ``diagonal entries'' all entries on the diagonal of the~$n_i \times n_i$ bottom submatrix of a factor~$A_i$. A column is unnecessary (hatched in the picture) if its diagonal entry has a~$\Gamma$ block on the right and no~$\Gamma$ block above. Good rows are those covered by a vector~$g_j$ or a~$\Gamma$ block, together with the basis vectors whose diagonal entry has no~$\Gamma$ block above (bad rows are hatched in the picture).

\begin{example}
\begin{enumerate}
\item If~$P_i$ is a segment, then~$n_i=1$, $m_i=2$~and~$\chi_i=1$, so that~$\Delta=1$ if~$\bm$~is even and~$0$~otherwise. Iterating this, if~$P_i$~is an $s$-cube, then~$\Delta\simeq\frac{s}{2}$. This yields \defn{neighborly cubical polytopes}~---~see~\cite{jz-ncp-00,js-ncps-07,sz-capdp}.
\item If~$P_i$ is an even cycle, then~$n_i=2$, $m_i=2p$~and~$\chi_i=2$, so that~$\Delta=2p-3$. This yields \defn{projected products of polygons}~---~see \cite{z-ppp-04,sz-capdp}.
\end{enumerate}
\index{polytope!neighborly cubical ---}
\index{neighborly cubical polytope}
\end{example}

In general, it is difficult to give the explicit ordering of the factors and choice of deformation that will yield the largest possible value of~$k$ attainable by a concrete product~$P_1\times\dots\times P_r$ of simple polytopes, and consequently to summarize this improvement by a precise proposition as we did for our first improvement. However, this best value can clearly be found by optimizing over the finite set of all possible orderings and types of deformation. Furthermore, we can be much more explicit for products of simplices, as we detail in the next section.


\subsection{Projection of deformed product of simplices}

We are now ready to apply this general construction to the particular case of products of simplices. For this, we represent the simplex~$\simplex_{n_i}$ by the inequality system~$A_i x \le b_i$, where
$$A_i  \eqdef  \begin{pmatrix} -1 & \dots & -1 \\ 1 & & \\ & \ddots & \\ & & 1 \end{pmatrix}$$
and~$b_i$ is a suitable right-hand side. We express the results of the construction with a case distinction according to the number~$s \eqdef |\ens{i\in[r]}{n_i=1}|$ of segments in the product~$\simplex_{\sub{n}}$.

\begin{proposition}\label{psn:prop:defp-ppsn}
Let~$\sub{n} \eqdef (n_1,\dots,n_r)$ with~$1=n_1=\cdots=n_s<n_{s+1}\le\cdots\le n_r$. Then
\begin{enumerate}[(1)]
\item for any~$0\le d\le s-1$, there exists a \dimensional{d} \xppsn{(k,\sub{n})} polytope as soon as
$$k \le \Fracfloor{d}{2}-r+s-1.$$
\item for any~$s\le d\le n$, there exists a \dimensional{d} \xppsn{(k,\sub{n})} polytope as soon as
$$k \le \Fracfloor{d+t-s}{2}-r+s,$$
where~$t\in\{s,\dots,r\}$ denotes the maximal integer such that~$\sum_{i=1}^{t} n_i\le d$.
\end{enumerate}
\end{proposition}

\begin{figure}[h]
	\capstart
	\centerline{\includegraphics[scale=1]{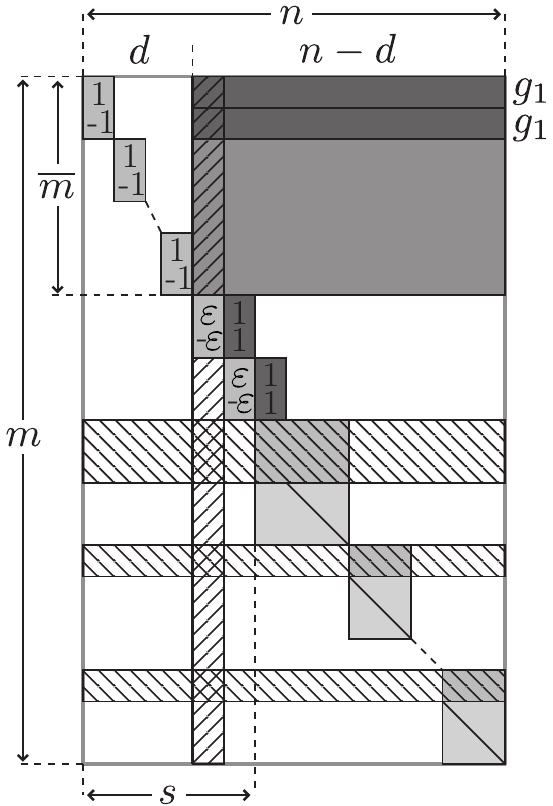}}
	\caption[\ppsn polytopes via deformed products, case~$s>d$]{How to obtain \ppsn polytopes from a deformed product construction, when the number~$s$ of segment factors exceeds the target dimension~$d$ of the projection.}
	\label{psn:fig:defp-ppsn1}
\end{figure}

\begin{proof}[Proof of (1)]
This is a special case of the results obtainable with the methods of Section~\ref{psn:subsec:deformedproducts:general}. The best construction is obtained using the matrix in \fref{psn:fig:defp-ppsn1}, from which we read off that:
\begin{enumerate}[(i)]
\item any subset of at least~$\alpha \eqdef n-\Fracfloor{d}{2}$ vectors in~$G$ is positively spanning;
\item the number of bad facets is~$\beta \eqdef r-s+1$, and therefore any \face{k} of~$\defP$ is contained in at least~$\gamma \eqdef n-k-r+s-1$ good facets.
\end{enumerate}
From this, the claim follows.
\end{proof}

\begin{figure}
	\capstart
	\centerline{\includegraphics[scale=1]{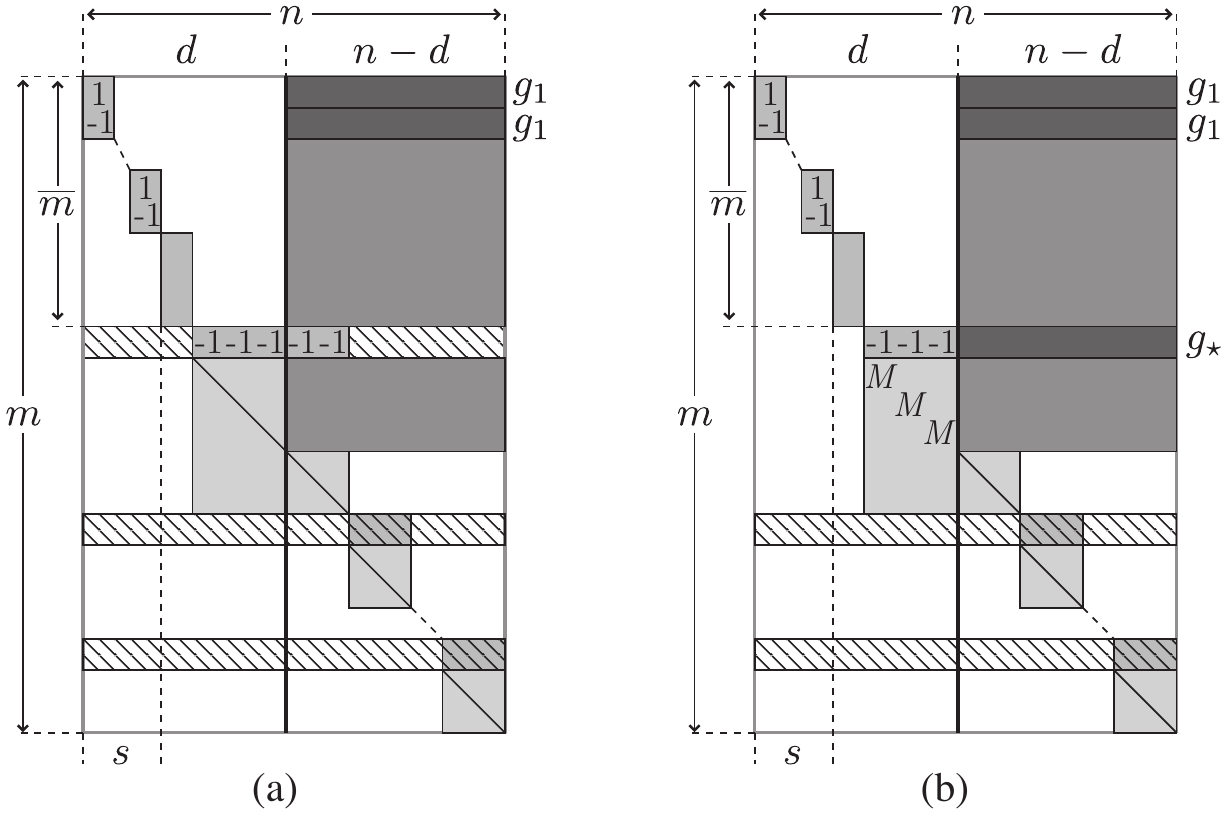}}
	\caption[\ppsn polytopes via deformed products, case~$s\le d$]{Obtaining \ppsn polytopes from a deformed product construction, when few of the factors are segments. Part~(a) shows the technique used so far, and Part~(b) an additional optimization that exchanges a bad facet for a new vector in the Gale transform.}
	\label{psn:fig:defp-ppsn2}
\end{figure}

\vspace{-.5cm}
\proof[Proof of (2)]
Consider the deformed product of \fref{psn:fig:defp-ppsn2}(a). Using similar calculations as before, we deduce that:
\begin{enumerate}[(i)]
\item any subset of at least~$\alpha \eqdef t-s+n-\Fracfloor{d+t-s-1}{2}$ vectors in~$G$ is positively spanning;
\item the number of bad facets is~$\beta \eqdef r-t$, and therefore any \face{k} of~$\defP$ is contained in at least~$\gamma \eqdef n-k-r+t$ good facets.
\end{enumerate}
This yields a bound of
$$k \le \Fracfloor{d+t-s-1}{2}-r+s.$$

We optimize the final~`$-1$' away by suitably deforming the matrix~$A_{t+1}$ as in \fref{psn:fig:defp-ppsn2}(b). This amounts to adding one more vector~$g_\star$ to the Gale diagram, so that the first row of~$A_{t+1}$ ceases to be a bad facet. This deformation is valid because:
\begin{enumerate}[(i)]
\item the matrix 
$$\begin{pmatrix}
      -1 & \dots & -1 & \star & \dots &  \star\\
      M \\
      & \ddots\\
      && M\\
      &&& 1\\
      &&&& \ddots\\
      &&&&& 1
\end{pmatrix}$$
still defines a simplex, as long as the~`$\star$' entries are negative and~$M\gg0$ is chosen sufficiently large. 
\item we can in fact choose the new vector~$g_\star$ to have only negative entries, by imposing an additional restriction on the Gale diagram~$G=\{e_1,\dots, e_{n-d}$, $g_1, \dots, g_{d+t},g_\star\}$ of~$Q$. Namely, we require the vertices of the simplicial \poly{(d+t)}tope~$Q$ that correspond to the Gale vectors~$g_1,\dots,g_{d+t}$ to lie on a facet. This forces the remaining vectors~$e_1,\dots,e_{n-d},g_\star$ to be positively spanning, so that~$g_\star$ has only negative entries.\qed
\end{enumerate}

\bvs
Finally, we reformulate Proposition~\ref{psn:prop:defp-ppsn} to express, in terms of~$k$ and~$\sub{n} \eqdef (n_1,\dots,n_r)$, what dimensions a \xppsn{(k,\sub{n})} polytope can have. This yields upper bounds on~$\delta_{pr}(k,\sub{n})$.

\begin{theorem}\label{psn:theo:defp-ppsn}
For any~$k\ge0$ and~$\sub{n} \eqdef (n_1,\dots,n_r)$ with~${1=n_1=\cdots=n_s<n_{s+1}\le\cdots\le n_r}$,
$$\delta_{pr}(k,\sub{n}) \le
   \begin{cases}
     2(k+r)-s-t   & \text{if } 3s \le 2k+2r, \\
     2(k+r-s)+1   & \text{if } 3s = 2k+2r+1, \\
     2(k+r-s+1)   & \text{if } 3s \ge 2k+2r+2,
   \end{cases}
$$
where~$t\in\{s,\dots,r\}$ is maximal such that 
$$3s+\sum_{i=s+1}^{t}(n_i+1) \le 2k+2r.$$
\end{theorem}

\begin{proof}
Apply Part~(1) of Proposition~\ref{psn:prop:defp-ppsn} when~$3s\ge2k+2r+2$ and Part~(2) otherwise.
\end{proof}

\begin{remark}
When all the~$n_i$'s are large compared to~$k$, the dimension of the \xppsn{(k,\sub{n})} polytope provided by this theorem is bigger than the dimension~$2k+r+1$ of the \xppsn{(k,\sub{n})} polytope obtained by the Minkowski sum of cyclic polytopes of Theorem~\ref{psn:theo:UBminkowskiCyclic}. However, if we have many segments (neighborly cubical polytopes), or more generally if many~$n_i$'s are small compared to~$k$, this construction provides our best examples of \ppsn polytopes.
\end{remark}


\section{Topological obstructions}\label{psn:sec:topologicalObstruction}

In this section, we give lower bounds on the minimal dimension~$\delta_{pr}(k,\sub{n})$ of a \xppsn{(k,\sub{n})} polytope, applying a method developed by Raman Sanyal~\cite{s-tovnms-09} to bound the number of vertices of Minkowski sums of polytopes. This method provides lower bounds on the target dimension of any linear projection that preserves a given set of faces of a polytope.  It uses Gale duality to associate a certain simplicial complex~$\cK$ to the set of faces that are preserved under the projection. Then lower bounds on the embeddability dimension of~$\cK$ transfer to lower bounds on the target dimension of the projection. In turn, the embeddability dimension is bounded via colorings of the Kneser graph of the system of minimal non-faces of~$\cK$, via Sarkaria's Embeddability~Theorem. 

After the completion of this work, we learned that Thilo R\"orig and Raman Sanyal~\cite{rs-npps,rorig-phd,sanyal-phd} also applied Sanyal's topological obstruction method to derive lower bounds on the target dimension of a projection preserving skeleta of different kind of products (products of polygons, products of simplices, and wedge products of polytopes). In particular, for a product~$\simplex_n\times\cdots\times\simplex_n$ of~$r$ identical simplices, $r\ge 2$, they obtain our Theorem~\ref{psn:theo:topObstr-small-k} and a result~\cite[Theorem~4.5]{rs-npps} that is only slightly weaker than Theorem~\ref{psn:theo:topObstr-all-k} in this setting. We decided to present this part in this dissertation since it completes our study on products~of~\mbox{polytopes}.

For the convenience of the reader, we first quickly recall Sarkaria's embeddability criterion. We then provide a brief overview of Raman Sanyal's method before applying it to obtain lower bounds on the dimension of \xppsn{(k,\sub{n})} polytopes. As mentioned in the introduction, these bounds match the upper bounds obtained from our different constructions for a wide range of parameters, and thus give the exact value of the minimal dimension of a \ppsn polytope.


\subsection{Sarkaria's embeddability criterion}

\subsubsection{Kneser graphs}\label{subsubsec:kneser}

Recall that a \defn{\kcoloring{k}} of a graph~$G \eqdef (V,E)$ is a map~$c:V\to[k]$ such that~$c(u)\ne c(v)$ for~$(u,v)\in E$. As usual, we denote~$\chi(G)$ the \defn{chromatic number} of~$G$ (\ie the minimal~$k$ such that~$G$ admits a \kcoloring{k}). We are interested in the chromatic number of Kneser graphs.

Let~$\cZ$ be a subset of the power set~$2^{[n]}$ of~$[n]$. The \defn{Kneser graph}\index{Kneser graph} on~$\cZ$, denoted~$\KG(\cZ)$, is the graph with vertex set~$\cZ$, where two vertices~$X,Y\in\cZ$ are related if and only if~$X\cap Y=\emptyset$:
$$\KG(\cZ) \eqdef \left(\cZ,\ens{(X,Y)\in\cZ^2}{X\cap Y=\emptyset}\right).$$
Let~$\KG_{n}^{k} \eqdef \KG\big({[n] \choose k}\big)$ denote the Kneser graph on the set of \tuple{k}s of~$[n]$.  For example, the graph~$\KG_{n}^{1}$ is the complete graph~$K_n$ (of chromatic number~$n$) and the graph~$\KG_{5}^{2}$ is the Petersen graph (of chromatic number~$3$). 


\begin{remark}
\begin{enumerate}[(i)]
\item If~$n\le 2k-1$, then any two \tuple{k}s of~$[n]$ intersect and the Kneser graph~$\KG_{n}^{k}$ is independent (\ie it has no edge). Thus its chromatic number is~$\chi(\KG_{n}^{k})=1$.
\item If~$n\ge 2k-1$, then~$\chi(\KG_{n}^{k})\le n-2k+2$. Indeed, the map~$c:{[n] \choose k}\to[n-2k+2]$ defined by~$c(F) \eqdef \min(F\cup\{n-2k+2\})$ is a \kcoloring{(n-2k+2)} of~$\KG_{n}^{k}$.
\end{enumerate}
\end{remark}

In fact, it turns out that this upper bound is the exact chromatic number of the Kneser graph:~$\chi(\KG_{n}^{k})=\max\{1,n-2k+2\}$. This result has been conjectured by Martin Kneser and proved by L\'aszl\'o Lov\'asz applying the Borsuk-Ulam Theorem~---~see~\cite{m-ubut-03} for more details. However, we will only need the upper bound for the topological obstruction.

\subsubsection{Sarkaria's Theorem}\label{psn:subsubsec:sarkaria}

\index{Sarkaria's Theorem}
Our lower bounds on the dimension of \xppsn{(k,\sub{n})} polytopes rely on lower bounds for the dimension in which certain simplicial complexes can be embedded. Among other possible methods (see~\cite{m-ubut-03}), we use Sarkaria's Coloring and Embedding Theorem.

We associate to any simplicial complex~$\cK$ the set system~$\cZ$ of \defn{minimal non-faces} of~$\cK$, that is, the inclusion-minimal sets of~$2^{V(\cK)}\ssm\cK$. For example, the system of minimal non-faces of the \skeleton{k} of the \simp{n} is ${[n+1] \choose k+2}$. Sarkaria's Theorem provides a lower bound on the dimension into which~$\cK$ can be embedded, in terms of the chromatic number of the~Kneser~graph~of~$\cZ$.

\begin{theorem}[Sarkaria's Theorem]\label{psn:theo:sarkaria}
Let~$\cK$ be a simplicial complex embeddable in~$\R^d$,~$\cZ$~be the system of minimal non-faces of~$\cK$, and~$\KG(\cZ)$ be the Kneser graph on~$\cZ$. Then
$$d \ge |V(\cK)|-\chi(\KG(\cZ))-1.$$
\end{theorem}
\vspace{-1cm}\qed
\vspace{.4cm}

In other words, we get large lower bounds on the possible embedding dimension of~$\cK$ when we obtain colorings with few colors of the Kneser graph on the system of minimal non-faces of~$\cK$. We refer to the excellent treatment in~\cite{m-ubut-03} for further details.


\subsection{Sanyal's topological obstruction method}

For given integers~$n>d$, we consider the orthogonal projection~$\pi:\R^n\to\R^d$ to the first~$d$ coordinates, and its dual projection~$\tau:\R^n\to\R^{n-d}$ to the last~$n-d$ coordinates. Let~$P$ be a full-dimensional simple polytope in~$\R^n$, with~$0$ in its interior, and assume that its vertices are strictly preserved under~$\pi$. Let~$F_1,\dots,F_m$ denote the facets of~$P$, and for all~$i\in[m]$, let~$f_i$ denote the normal vector of~$F_i$, and~$g_i \eqdef \tau(f_i)$. For any face~$F$ of~$P$, let~$\varphi(F)$ denote the set of indices of the facets of~$P$ containing~$F$, \ie such that~$F=\bigcap_{i\in \varphi(F)} F_i$.

\begin{lemma}[Sanyal \cite{s-tovnms-09}]\label{psn:lem:gale_transform}
The vector configuration~$G \eqdef \ens{g_i}{i\in[m]}\subset\R^{n-d}$ is the Gale transform of the vertex set of a (full-dimensional) polytope~$Q$ of~$\R^{m-n+d-1}$. Up to a slight perturbation of the facets of~$P$, we can even assume~$Q$ to be simplicial.\qed
\end{lemma}

\index{Sanyal's projection polytope}
\index{polytope!Sanyal's projection ---}
We will refer to the polytope~$Q$ as \defn{Sanyal's projection polytope}. We denote its vertices by~$a_1,\dots,a_m$ such that~$a_i$ is transformed into~$g_i$ by Gale duality. The faces of this polytope capture the key notion of strictly preserved faces of~$P$~---~remember Definition~\ref{psn:def:spf}. Indeed, the Projection Lemma~\ref{psn:lem:projection} ensures that for any face~$F$ of~$P$ strictly preserved by the projection~$\pi$, the set~$\ens{g_i}{i\in\varphi(F)}$ is positively spanning, which implies by Gale duality that the set of vertices~$\ens{a_i}{i\in[m]\ssm\varphi(F)}$ forms a face of~$Q$.

\begin{figure}[b]
	\capstart
	\centerline{\includegraphics[scale=1]{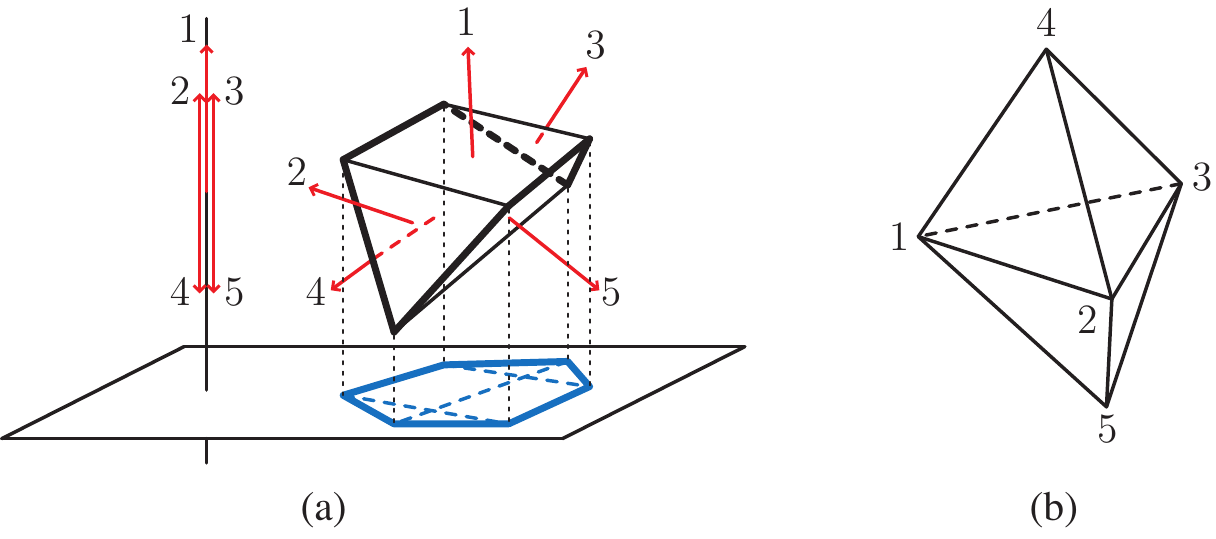}}
	\caption[A projection of a triangular prism and its associated projection polytope~$Q$]{(a)~A projection of a triangular prism and (b)~its associated projection polytope~$Q$. The six faces of~$Q$ corresponding to the six edges of~$P$ preserved under projection (bold) make up the entire boundary complex of~$Q$. }
	\label{psn:fig:entireQ}
\end{figure}

\begin{example}\label{ex:entireQ}
Let~$P$ be a triangular prism in $3$-space that projects to a hexagon as in \fref{psn:fig:entireQ}(a), so that~$n=3$, $d=2$~and~$m=5$. The vector configuration~$G\subset\R^1$ obtained by projecting~$P$'s normal vectors consists of three vectors pointing up and two pointing down, so that Sanyal's projection polytope~$Q$ is a bipyramid over a triangle. An edge~$F_i\cap F_j$ of~$P$ that is preserved under projection corresponds to the face~$[5]\ssm\{i,j\}$ of~$Q$. Notice that the six faces of~$Q$ corresponding to the six edges of~$P$ that are preserved under projection (in bold in \fref{psn:fig:entireQ}(a)) make up the entire boundary complex of the bipyramid~$Q$.
\end{example}

\bvs
Let~$\cF$ be a subset of the set of all strictly preserved faces of~$P$ under~$\pi$. Define~$\cK$ to be the simplicial complex induced by~$\ens{[m]\ssm\varphi(F)}{F\in \cF}$. 

\begin{remark}
Notice that not all non-empty faces of~$\cK$ correspond to non-empty faces in~$\cF$: in Example~\ref{ex:entireQ}, if~$\cF$ consists of all strictly preserved edges, then~$\cK$ is the entire boundary complex of Sanyal's projection polytope~$Q$, so that it contains the edge~$\{2,3\}$. But then the complementary intersection of facets,~$F_1\cap F_4\cap F_5$, does not correspond to any non-empty face of~$P$ (and \textit{a fortiori} of~$\cF$).
\end{remark}

Since the set of vertices~$\ens{a_i}{i\in[m]\ssm\varphi(F)}$ forms a face of~$Q$ for any face~$F\in\cF$, and since $Q$~is
simplicial, $\cK$~is a subcomplex of the face complex of~$Q\subset\R^{m-n+d-1}$. In particular, when~$\cK$~is not the entire boundary complex of~$Q$, it embeds into~$\R^{m-n+d-2}$ by stereographic projection (otherwise, it only embeds into~$\R^{m-n+d-1}$, as happens in Example~\ref{ex:entireQ}).

Thus, given the simple polytope~$P\subset\R^n$, and a set~$\cF$ of faces of~$P$ that we want to preserve under projection, the study of the embeddability of the corresponding abstract simplicial complex~$\cK$ provides lower bounds on the dimension~$d$ in which we can project~$P$. We proceed in the following way:
\begin{enumerate}
\item we first choose our subset~$\cF$ of strictly preserved faces sufficiently simple to be understandable, and sufficiently large to provide an obstruction;
\item we then understand the system~$\cZ$ of minimal non-faces of the simplicial complex~$\cK$;
\item finally, we find a suitable coloring of the Kneser graph on~$\cZ$ and apply Sarkaria's Theorem~\ref{psn:theo:sarkaria} to bound the dimension in which~$\cK$~can be embedded: a \kcoloring{t} of~$\KG(\cZ)$ ensures that~$\cK$ is not embeddable into~$|V(\cK)|-t-2=m-t-2$, which by the previous paragraph bounds the dimension~$d$ from below as follows:
\end{enumerate}

\begin{theorem}[Sanyal~\cite{s-tovnms-09}]\label{psn:theo:sanyal}
Let~$P$ be a simple polytope in~$\R^n$ whose facets are in general position, and let~$\pi:\R^n\to\R^d$ be a projection. Let~$\cF$ be a subset of the set of all strictly preserved faces of~$P$ under~$\pi$, let~$\cK$~be the simplicial complex induced by~$\ens{[m]\ssm\varphi(F)}{F\in \cF}$, and let~$\cZ$~be its system of minimal non-faces. If the Kneser graph~$\KG(\cZ)$ is
\kcolorable{t}, then:
\begin{enumerate}[(1)]
\item if~$\cK$ is not the entire boundary complex of the Sanyal polytope~$Q$, then~$d\ge n-t+1$;
\item otherwise,~$d\ge n-t$.\qed
\end{enumerate}
\end{theorem}

In the remainder of this section, we apply Sanyal's topological obstruction to our problem. The hope was initially to extend it to bound the target dimension of a projection preserving the \skeleton{k} of an arbitrary product of simple polytopes. However, the combinatorics involved to deal with this general question turn out to be too complicated and we restrict our attention to products of simplices. We obtain in this manner bounds on the minimal dimension~$\delta_{pr}(k,\sub{n})$ of a \xppsn{(k,\sub{n})} polytope.


\subsection{Preserving the \skeleton{k} of a product of simplices}

In this section, we understand the abstract simplicial complex~$\cK$ corresponding to our problem, and describe its system
of minimal non-faces.

The facets of~$\simplex_{\sub{n}}$ are exactly the products
$$\psi_{i,j} \eqdef \simplex_{n_1}\times\cdots\times\simplex_{n_{i-1}}\times(\simplex_{n_i}\ssm\{j\})\times\simplex_{n_{i+1}}\times\cdots\times\simplex_{n_r},$$
for~$i\in[r]$ and~$j\in[n_i+1]$.  We identify the facet~$\psi_{i,j}$ with the element~$j\in[n_i+1]$ of the disjoint union~$[n_1+1]\uplus[n_2+1]\uplus\cdots\uplus[n_r+1]$.

Let~$F \eqdef F_1\times\cdots\times F_r$ be a \face{k} of~$\simplex_{\sub{n}}$. Then~$F$ is contained in a facet~$\psi_{i,j}$ of~$\simplex_{\sub{n}}$ if and only if~$j\notin F_i$. Thus, the set of facets of~$\simplex_{\sub{n}}$ that do not contain~$F$ is exactly~${F_1\uplus\cdots\uplus F_r}$. Consequently, if we want to preserve the \skeleton{k} of~$\simplex_{\sub{n}}$, then the abstract simplicial complex~$\cK$ we are interested in is induced by
\begin{equation}\label{psn:eq:complexK}
\bigg\{F_1\uplus\cdots\uplus F_r \ \bigg| \ \emptyset\ne F_i\subset [n_i+1] \text{ for all } i\in[r], \text{ and } \sum_{i\in[r]} (|F_i|-1)=k\bigg\}.
\end{equation}
 
\begin{remark}\label{psn:rem:entireQ}
In contrast to the general case, when we want to preserve the \defn{complete} \skeleton{k} of a product of simplices, the complex~$\cK$ cannot be the entire boundary complex of the Sanyal polytope~$Q$. In consequence, the better lower bound from Part~(1) of Sanyal's Theorem~\ref{psn:theo:sanyal} always holds, and we always use it from now on without further notice.

To prove that~$\cK$ cannot cover the entire boundary complex of~$Q$, observe that
$$\dim Q = m-n+d-1 = \sum (n_i+1) - \sum n_i + d - 1 = r+d-1,$$
while~$\dim\cK=r+k-1$ by~\eqref{psn:eq:complexK}. A necessary condition for~$\cK$ to be the entire boundary complex of~$Q$ is that~$\dim\cK=\dim Q-1$, which translates to~$d=k+1$. Now suppose that the entire \skeleton{k} of~$\simplex_{\sub{n}}$ is preserved under projection to dimension~$k+1$. Then the projections of those \face{k}s are facets of~$\pi(\simplex_{\sub{n}})$. Since any ridge of the projected polytope is contained in exactly two facets, and the \defn{entire} \skeleton{k} of~$\simplex_{\sub{n}}$ is preserved, we know that any \face{(k-1)} of~$\simplex_{\sub{n}}$ is also contained in exactly two \face{k}s. But this can only happen if~$k=n-1$, which means~$n=d$.

Observe again that~$\cK$ can be the entire boundary complex of~$Q$ if we do not preserve \defn{all} \face{k}s of~$\simplex_{\sub{n}}$~---~see Example~\ref{ex:entireQ}.
\end{remark}

The following lemma gives a description of the minimal non-faces of~$\cK$:

\begin{lemma}
The system of minimal non-faces of $\cK$ is
$$\cZ = \bigg\{G_1\uplus\cdots\uplus G_r \ \bigg| \ |G_i|\ne 1 \text{ for all }i\in[r], \text{ and } \sum_{i\,|\,G_i\ne\emptyset} (|G_i|-1)=k+1\bigg\}.$$
\end{lemma}

\begin{proof}
A subset~$G \eqdef G_1\uplus\cdots\uplus G_r$ of~$[n_1+1]\uplus[n_2+1]\uplus\cdots\uplus[n_r+1]$ is a face of~$\cK$ when it can be extended to a subset~$F_1\uplus\cdots\uplus F_r$ with~$\sum (|F_i|-1)=k$ and~$\emptyset\ne F_i\subset [n_i+1]$ for all~$i\in[r]$, that is, when
$$k \ge \left|\ens{i\in[r]}{G_i=\emptyset}\right| + \sum_{i\in[r]} (|G_i|-1) = \sum_{i\,|\,  G_i\ne\emptyset} (|G_i|-1).$$
Thus,~$G$ is a non-face if and only if 
$$\sum_{i\, | \,  G_i\ne\emptyset} (|G_i|-1) \ge k+1.$$

If~$\sum_{i\, | \,  G_i\ne\emptyset} (|G_i|-1)> k+1$, then removing any element provides a smaller non-face. If there is an~$i$ such that~$|G_i|=1$, then removing the unique element of~$G_i$ provides a smaller non-face. Thus, if~$G$ is a minimal non-face, then~$\sum_{i\, | \,  G_i\ne\emptyset} (|G_i|-1)=k+1$, and~$|G_i|\ne1$ for all~$i\in[r]$.

Reciprocally, if~$G$ is a non-minimal non-face, then it is possible to remove one element keeping a non-face. Let~$i\in[r]$ be such that we can remove one element from~$G_i$, keeping a non-face. Then, either~$|G_i|=1$, or 
$$\sum_{j\, | \,  G_j\ne\emptyset} (|G_j|-1) \ge 1+(|G_i|-2)+\sum_{j\ne i\, | \,  G_j\ne\emptyset} (|G_j|-1)\ge k+2,$$
since we keep a non-face.
\end{proof}


\subsection{Colorings of~$\KG(\cZ)$}

The next step consists in providing a suitable coloring for the Kneser graph on the system~$\cZ$ of minimal non-faces of~$\cK$. Let~$S \eqdef \ens{i\in[r]}{n_i=1}$ denote the set of indices corresponding to the segments, and~$R \eqdef \ens{i\in[r]}{n_i\ge2}$ the set of indices corresponding to the non-segments in the product~$\simplex_{\sub{n}}$. We first provide a coloring for two extremal situations.

\begin{theorem}[Topological obstruction for low-dimensional skeleta]\label{psn:theo:topObstr-small-k}
If~$k\le \sum_{i\in R} \Fracfloor{n_i-2}{2}$, then the dimension of any \xppsn{(k,\sub{n})} polytope cannot be smaller than~$2k+|R|+1$: 
$$\delta_{pr}(k,\sub{n}) \ge 2k+|R|+1.$$
\end{theorem}

\proof
Let~$k_1,\dots,k_r\in\N$ be such that 
$$\sum_{i\in[r]} k_i = k \quad\text{and}\quad 
\begin{cases}
  k_i=0 & \text{for } i\in S;\\
  0\le k_i\le \frac{n_i-2}{2} & \text{for } i\in R.
\end{cases}
$$

Observe that:
\begin{enumerate}
\item Such a tuple exists since~$k\le \sum_{i\in R} \Fracfloor{n_i-2}{2}$.
\item For any minimal non-face~$G \eqdef G_1\uplus\cdots\uplus G_r$ of~$\cZ$, there exists~$i\in[r]$ such that~$|G_i|\ge k_i+2$. Indeed, if~$|G_i|\le k_i+1$ for all~$i\in[r]$, then
$$k+1 = \sum_{i\, | \,  G_i\ne \emptyset} (|G_i|-1) \le \sum_{i\, | \,  G_i\ne \emptyset} k_i \le \sum_{i\in[r]} k_i = k,$$
which is impossible.
\end{enumerate}

For all~$i\in[r]$, we fix a proper coloring~$\gamma_i:{[n_i+1]\choose[k_i+2]}\to[\chi_i]$ of the Kneser graph~$\KG_{n_i+1}^{k_i+2}$, with~$\chi_i=1$ color if~$i\in S$ and~$\chi_i=n_i-2k_i-1$ colors if~$i\in R$~---~see Section~\ref{subsubsec:kneser}. Then, we define a coloring~$\gamma:\cZ\to[\chi_1]\uplus\cdots\uplus[\chi_r]$ of the Kneser graph on~$\cZ$ as follows. Let~$G \eqdef G_1\uplus\cdots\uplus G_r$ be a given minimal non-face of~$\cZ$. We choose arbitrarily an~$i\in[r]$ such that~$|G_i|\ge k_i+2$, and again arbitrarily a subset~$g$ of~$G_i$ with~$k_i+2$ elements. We color~$G$ with the color of~$g$ in~$\KG_{n_i+1}^{k_i+2}$, that is, we define~$\gamma(G) \eqdef \gamma_i(g)$.

The coloring~$\gamma$ is a proper coloring of the Kneser graph~$\KG(\cZ)$. Indeed, let~$G \eqdef G_1\uplus\cdots\uplus G_r$ and~$H \eqdef H_1\uplus\cdots\uplus H_r$ be two minimal non-faces of~$\cZ$ related by an edge in~$\KG(\cZ)$, which means that they do not intersect. Let~$i\in[r]$ and~$g\subset G_i$ be such that we have colored~$G$ with~$\gamma_i(g)$, and similarly~$j\in[r]$ and~$h\subset G_j$ be such that we have colored~$H$ with~$\gamma_j(h)$. Since the color sets of~$\gamma_1,\dots,\gamma_r$ are disjoint, the non-faces~$G$~and~$H$ can receive the same color~$\gamma_i(G)=\gamma_j(H)$ only if~$i=j$ and~$g$~and~$h$ are not related by an edge in~$\KG_{n_i+1}^{k_i+2}$, which implies that~$g\cap h\ne\emptyset$. But this cannot happen, because~$g\cap h\subset G_i\cap H_i$, which is empty by assumption.  Thus,~$G$~and~$H$ get different colors.

This provides a proper coloring of~$\KG(\cZ)$ with~$\sum\chi_i$ colors. According to Theorem~\ref{psn:theo:sanyal} and Remark~\ref{psn:rem:entireQ}, we know that the dimension~$d$ of the projection is at least
$$\sum_{i\in[r]} n_i - \sum_{i\in[r]} \chi_i +1 = 2k+|R|+1.$$

\vspace*{-1.4cm}\qed
\vspace{.7cm}

\begin{theorem}[Topological obstruction for high-dimensional skeleta]\label{psn:theo:topObstr-large-k}
If~$k\ge \Floor{\frac12\sum_i n_i}$, then any \xppsn{(k,\sub{n})} polytope is combinatorially equivalent to~$\simplex_{\sub{n}}$:
$$\delta_{pr}(k,\sub{n}) \ge \sum n_i.$$
\end{theorem}

\begin{proof}
Let~$G \eqdef G_1\uplus\cdots\uplus G_r$ and~$H \eqdef H_1\uplus\cdots\uplus H_r$ be two minimal non-faces of~$\cZ$. Let $A \eqdef \ens{i\in[r]}{G_i\ne\emptyset \text{ or } H_i\ne\emptyset}$. Then
\begin{eqnarray*}
\sum_{i\in A} (|G_i|+|H_i|) & \ge & \sum_{G_i\ne\emptyset} (|G_i|-1) + \sum_{H_i\ne\emptyset} (|H_i|-1) + |A| \\ 
  & = & 2k+2+|A| > \sum_{i\in[r]} n_i+|A| \ge \sum_{i\in A} (n_i+1).
\end{eqnarray*}
Thus, there exists~$i\in A$ such that~$|G_i|+|H_i|>n_i+1$, which implies that~$G_i\cap H_i\ne\emptyset$, and proves that~$G\cap H\ne\emptyset$.

Consequently, the Kneser graph~$\KG(\cZ)$ is independent (and we can color it with only one color). We obtain that the dimension~$d$ of the projection is at least~$\sum n_i$. In other words, in this extremal case, there is no better \xpsn{(k,\sub{n})} polytope than the product~$\simplex_{\sub{n}}$ itself.
\end{proof}

\begin{remark}\label{remark:betterColoring}
Theorem \ref{psn:theo:topObstr-large-k} can sometimes be strengthened a little: if~$k = \frac{1}{2}\sum n_i-1$, and~$k+1$ is not representable as a sum of a subset of~$\{n_1,\dots,n_r\}$, then~$\delta_{pr}(k,\sub{n})=\sum n_i$. 
\end{remark}

\begin{proof}
As in the previous theorem, we prove that the Kneser graph~$\KG(\cZ)$~is independent. Indeed, assume that~$G \eqdef G_1\uplus\cdots\uplus G_r$ and~$H \eqdef H_1\uplus\cdots\uplus H_r$ are two minimal non-faces of~$\cZ$ related by an edge in~$\KG(\cZ)$. Then,~$G\cap H$ is empty, which implies that for all~$i\in [r]$,
\begin{equation}\label{psn:eq:gihi}
|G_i|+|H_i| \le n_i+1.
\end{equation}
Let~$U \eqdef \ens{i}{G_i\ne\emptyset}$ and~$V \eqdef \ens{i}{H_i\ne\emptyset}$. Then,
\begin{eqnarray*}
\sum_{i\in U\cup V} (|G_i| + |H_i|) & = & \sum_{i\in U} (|G_i|-1) + \sum_{i\in V} (|H_i|-1) + |U| + |V| = 2k+2 + |U| + |V| \\
& = & \sum_{i\in[r]} n_i + |U| + |V| \stackrel{(\star)}{\ge} \sum_{i\in U\cup V} n_i + |U\cup V| = \sum_{i\in U\cup V} (n_i+1).
\end{eqnarray*}
Summing \eqref{psn:eq:gihi} over~$i\in U\cup V$ implies that both the inequality~$(\star)$ and \eqref{psn:eq:gihi} for~$i\in U\cup V$ are in fact equalities. The tightness of~$(\star)$ implies furthermore that~$|U|+|V|=|U\cup V|$, so that~$U\cap V=\emptyset$; in other words,~$H_i$~is empty whenever~$G_i$~is not. The equality in~\eqref{psn:eq:gihi} then asserts that~$|G_i|=n_i+1$ for all~$i\in U$, and therefore
$$k+1 = \sum_{i\in U} (|G_i|-1) = \sum_{i\in U} n_i$$
is representable as a sum of a subset of the~$n_i$, which contradicts the assumption. 
\end{proof}

Finally, to fill the gap in the ranges of~$k$ covered by Theorems~\ref{psn:theo:topObstr-small-k} and \ref{psn:theo:topObstr-large-k}, we merge both coloring ideas as follows.

We partition~$[r]=A\uplus B$ and choose~$k_i\ge0$ for all~$i\in A$ and~$k_B\ge0$ such that
\begin{equation}\label{psn:eq:SumKiAndKB}
\left(\sum_{i\in A} k_i \right) + k_B \le k.
\end{equation}
We will see later what the best choice for~$A$,~$B$,~$k_B$ and the~$k_i$'s is. Let~$n_B \eqdef \sum_{i\in B}n_i$. Color the Kneser graphs~$\KG_{n_i+1}^{k_i+2}$ for~$i\in A$ and~$\KG_{n_B}^{k_B+1}$ with pairwise disjoint color sets with
$$\chi_i=
\begin{cases}
 n_i-2k_i-1 &\textnormal{if }2k_i\le n_i-2, \\
 1 &\textnormal{if }2k_i\ge n_i-2, \\
\end{cases}
\qquad
\text{and}
\qquad
\chi_B=
\begin{cases}
  0 & \textnormal{if } n_B=0, \\
  n_B-2k_B &\textnormal{if }2k_B\le n_B-1, \\
  1 &\textnormal{if }2k_B\ge n_B-1, \\
\end{cases}
$$
colors respectively. 

Observe now that for all minimal non-faces~$G \eqdef G_1\uplus\cdots\uplus G_r$,
\begin{itemize}
\item either there is an~$i\in A$ such that~$|G_i|\ge k_i+2$,
\item or~$\sum_{i\in B\,|\,G_i\neq\emptyset} (|G_i|-1)\ge k_B+1$.
\end{itemize}
Indeed, otherwise
$$k+1 = \sum_{i\,|\,G_i\neq\emptyset}(|G_i|-1) \le \left(\sum_{i\in A}k_i\right) + k_B \le k.$$
This permits us to define a coloring of~$\KG(\cZ)$ in the following way. For each minimal non-face~$G=G_1\uplus\cdots\uplus G_r$, we arbitrarily choose one of the following strategies:
\begin{enumerate}
\item If we can find an~$i\in A$ such that~$|G_i|\ge k_i+2$, we choose an arbitrary subset~$g$ of~$G_i$ with~$k_i+2$ elements, and color~$G$ with the color of~$g$ in~$\KG_{n_i+1}^{k_i+2}$;
\item Otherwise,~$\sum_{i\in B\,|\,G_i\neq\emptyset} (|G_i|-1)\ge k_B+1$, and we choose an arbitrary subset~$g$ of 
$$\biguplus_{i\in B}(G_i\ssm \{n_i+1\}) \subset \biguplus_{i\in B}[n_i]$$
with~$k_B+1$ elements and color~$G$ with the color of~$g$ in~$\KG_{n_B}^{k_B+1}$.
\end{enumerate}
By exactly the same argument as in the proof of Theorem~\ref{psn:theo:topObstr-small-k}, one can verify that this provides a valid coloring of the Kneser graph~$KG(\cZ)$ with
$$\chi  \eqdef  \chi(A,B,\sub{k_i},k_B)  \eqdef  \sum_{i\in A}\chi_i+\chi_B$$
many colors.  Therefore Sanyal's Theorem~\ref{psn:theo:sanyal} and Remark~\ref{psn:rem:entireQ} yield the following lower bound on the dimension~$d$ of any \xppsn{(k,\sub{n})} polytope:
$$d \ge d_k  \eqdef  d_k(A,B,\sub{k_i},k_B)  \eqdef  \sum_i n_i+1-\chi \ge \delta_{pr}(k,\sub{n}).$$

It remains to choose parameters~$A$,~$B$, and~$\ens{k_i}{i\in A}$ and~$k_B$ that maximize this bound. We proceed algorithmically, by first fixing~$A$ and~$B$, and choosing the~$k_i$'s and~$k_B$ to maximize the bound on the dimension~$d_k$. For this, we first start with~$k_i=0$ for all~$i$ and~$k_B=0$, and observe the variation of~$d_k$ as we increase individual~$k_i$'s or~$k_B$. By (\ref{psn:eq:SumKiAndKB}), we are only allowed a total of~$k$ such increases. During this process, we will always maintain the conditions~$2k_i\le n_i-1$ for all~$i\in A_i$ and~$2k_B\le n_B$ (which makes sense by the formulas for~$\chi_i$ and~$\chi_B$).

We start with~$k_i=0$ for all~$i$ and~$k_B=0$. Then
\begin{eqnarray*}
\chi(A,B,\sub{0},0) & = & \sum_{i\in A}(n_i-1)+|S\cap A| + n_B \\
 & = & \sum_{i\in A} n_i-|A|+|S\cap A|+\sum_{i\in B} n_i  = \sum_{i\in[r]} n_i-r+|B\cup S|,
\end{eqnarray*}
and
$$d_k(A,B,\sub{0},0) = 1+r-|B\cup S|,$$
where~$S \eqdef \ens{i\in[r]}{n_i=1}$ denotes the set of segments. 

We now study the variation of~$d_k$ as we increase each of the~$k_i$'s and~$k_B$ by one. For~$i\in A$, increasing~$k_i$ by one decreases~$\chi_i$ by
$$
\begin{cases}
 2,&\textnormal{ if }2k_i\le n_i-4,\\
 1,&\textnormal{ if }2k_i= n_i-3,\\
 0,&\textnormal{ if }2k_i\ge n_i-2,
\end{cases}
$$
and hence increases~$d_k$ by the same amount. Observe in particular that~$d_k$ remains invariant if we increase~$k_i$ for some segment~$i\in S$ (because~$n_i=1$ for segments). Thus, it makes sense to choose~$B$ to contain all segments. Similarly, increasing~$k_B$ by one decreases~$\chi_B$ by 
$$
\begin{cases}
 2,&\textnormal{ if }2k_B\le n_B-3,\\
 1,&\textnormal{ if }2k_B= n_B-2,\\
 0,&\textnormal{ if }2k_B\ge n_B-1,
\end{cases}
$$
and increases~$d_k$ by the same amount.

Recall that we are allowed at most~$k$ increases of~$k_i$'s or~$k_B$ by~\eqref{psn:eq:SumKiAndKB}. Heuristically, it seems reasonable to first increase the~$k_i$'s or~$k_B$ that increase~$d_k$ by two, and then these that increase~$d_k$ by one. Hence we get a case distinction on~$k$, which  also depends on~$A$~and~$B$:

\begin{theorem}[Topological obstruction, general case]\label{psn:theo:topObstr-all-k}
Let~$k\ge0$ and~$\sub{n} \eqdef (n_1,\dots,n_r)$ with $r\ge1$ and~$n_i\ge1$ for all~$i$. Let~$[r]=A\uplus B$ be a partition of~$[r]$ with~$B\supset S \eqdef \ens{i\in[r]}{n_i=1}$.
Define
\begin{eqnarray*}
K_1 &  \!\!\!\! \eqdef \!\!\!\!  & K_1(A,B)  \eqdef  \sum_{i\in A}\Fracfloor{n_i-2}{2} + \max\left\{0,\Fracfloor{n_B-1}{2}\right\}, \\ 
K_2 &  \!\!\!\! \eqdef \!\!\!\!  & K_2(A,B)  \eqdef  \left|\ens{i\in A}{n_i\text{ is odd}}\right| +
 \begin{cases}
   1 & \text{if } n_B\text{ is even and non-zero}, \\
   0 & \text{otherwise}.
 \end{cases}
\end{eqnarray*}
Then the following lower bounds hold for the dimension of a \xppsn{(k,\sub{n})} polytope:
\begin{center}
\renewcommand{\arraystretch}{1.4}
\begin{tabular}{ll}
(1) If~$0\le k\le K_1$, then & $\delta_{pr}(k,\sub{n}) \ge r+1-|B|+2k$; \\
(2) If~$K_1\le k\le K_1+K_2$, then & $\delta_{pr}(k,\sub{n}) \ge r+1-|B|+K_1+k$; \\
(3) If~$K_1+K_2\le k$, then & $\delta_{pr}(k,\sub{n}) \ge r+1-|B|+2K_1+K_2.$
\end{tabular}
\end{center}
\end{theorem}

This theorem enables us to recover Theorem~\ref{psn:theo:topObstr-small-k} and Theorem~\ref{psn:theo:topObstr-large-k}:

\begin{corollary}\label{psn:coro:topObstr}
Let~$k\ge0$ and~$\sub{n} \eqdef (n_1,\dots,n_r)$ with~$r\ge1$ and~$n_i\ge1$ for all~$i$, and define~$S \eqdef \ens{i\in[r]}{n_i=1}$  and~$R \eqdef \ens{i\in[r]}{n_i\ge2}$. Then:
\begin{enumerate}
\item If 
$$0 \le k \le \sum_{i\in R}\Fracfloor{n_i-2}{2} + \max\left\{0,\Fracfloor{|S|-1}{2}\right\},$$
then~$\delta_{pr}(k,\sub{n}) \ge 2k+|R|+1$.
\item If~$k\ge \Floor{\frac12\sum n_i}$ then~$\delta_{pr}(k,\sub{n})\ge \sum_i n_i$.
\end{enumerate}
\end{corollary}

\begin{proof}
Take~$A=R$ and~$B=S$ for~(1), and~$A=\emptyset$ and~$B=[r]$ for~(2).
\end{proof}

To conclude, observe that the upper bounds on~$\delta_{pr}(k,\sub{n})$ provided by the constructions of Section~\ref{psn:sec:cyclic} and~\ref{psn:sec:deformedproducts} match the lower bound of Corollary~\ref{psn:coro:topObstr} for a wide range of parameters. In particular:

\begin{theorem}\label{psn:theo:exactvaluedelta}
For any~$\sub{n} \eqdef (n_1,\dots,n_r)$ with~$r\ge1$ and~$n_i\ge2$ for all~$i$, and for any~$k$ such that~$0\le k\le \sum_{i\in [r]}\Fracfloor{n_i-2}{2}$, the smallest \xppsn{(k,\sub{n})} polytope has dimension exactly~$2k+r+1$. In other words:
$$\delta_{pr}(k,\sub{n}) = 2k+r+1.$$
\end{theorem}



\backmatter 
\setcounter{compteurlof}{0}

\pagestyle{hautpageintro}
\addtocontents{toc}{ \vspace{.5cm} \demiligne \vspace{.1cm} }

\renewcommand{\namepart}{Bibliography}
\bibliographystyle{alpha}
\bibliography{biblio.bib}

\cleardoublepage

\renewcommand{\namepart}{Index}
\printindex

\pagestyle{empty}

\cleardoublepage

~

\newpage

\vspace*{-2.4cm}
\enlargethispage{2.1cm}

Cette thèse aborde deux sujets particuliers de géométrie discrète, les \defn{multitriangulations} et les \defn{réalisations polytopales de produits}, dont la problématique commune est la recherche de réalisations polytopales de structures combinatoires.

Une \ktri{k} est un ensemble maximal de cordes du \gon{n}e convexe ne contenant pas de sous-ensemble de~$k+1$ cordes qui se croisent deux à deux.
Nous proposons une étude combinatoire et géométrique des multitriangulations basée sur leurs étoiles, qui jouent ici le même rôle que les triangles des triangulations. Cette étude nous amène à interpréter les multitriangulations par dualité comme des \defn{arrangements de pseudodroites avec points de contact} de support donné. Nous exploitons finalement ces résultats pour discuter quelques problèmes ouverts sur les multitriangulations, en particulier la question de la réalisation polytopale de leurs graphes~de~flip.

Nous étudions dans un deuxième temps la polytopalité des produits cartésiens. Nous nous interrogeons d'abord sur l'existence de réalisations polytopales d'un produit cartesien de graphes, puis nous recherchons la dimension minimale que peut avoir un polytope dont le $k$-sque\-lette est celui d'un produit de simplexes.

\svs
\noindent\textbf{Mots-clés~:} triangulation, multitriangulation, pseudotriangulation, graphe de flip, associahèdre, nombre de Catalan, arrangement de pseudodroites, rigidité, polytope, produit cartésien.

\svs
\centerline{\rule[3pt]{5cm}{.5pt}}
\svs

Esta memoria versa sobre dos temas particulares de geometría discreta, las \defn{multitriangulaciones} y las \defn{realizaciones politopales de productos}, cuya problemática común es la búsqueda de realizaciones politopales de estructuras combinatorias.

Una $k$-triangulación es un conjunto maximal de cuerdas del \gon{n}o convexo que no contiene ningún subconjunto de~$k+1$ cuerdas que se cruzan dos a dos. Proponemos un estudio combinatorio y geométrico de las multitriangulaciones basado en sus estrellas, que desempeñan el mismo papel que los triángulos de los triangulaciones. Este estudio nos lleva a interpretar las multitriangulaciones por dualidad como \defn{arreglos de pseudorectas con puntos de contacto} de soporte dado. Explotamos finalmente estos resultados para discutir algunos problemas abiertos sobre multitriangulaciones, en particular la pregunta de la realización politopal de sus~grafos~de~flip.

Estudiamos en segundo lugar la politopalidad de productos cartesianos. Nos preguntamos primero sobre la existencia de realizaciones politopales de un producto cartesiano de grafos, y buscamos después la dimensión mínima que puede tener un politopo cuyo $k$-esque\-leto es el de un producto de símplices.

\svs
\noindent\textbf{Palabras-clave:} triangulación, multitriangulación, pseudotriangulación, grafo de flip, associaedro, número de Catalan, arreglo de pseudorectas, rigidez, politopo, producto cartesiano.

\svs
\centerline{\rule[3pt]{5cm}{.5pt}}
\svs

This thesis explores two specific topics of discrete geometry, the \defn{multitriangulations} and the \defn{polytopal realizations of products}, whose connection is the problem of finding polytopal realizations of a given combinatorial structure.

A \ktri{k} is a maximal set of chords of the convex \gon{n} such that no~$k+1$ of them mutually cross. We propose a combinatorial and geometric study of multitriangulations based on their stars, which play the same role as triangles of triangulations. This study leads to interpret multitriangulations by duality as \defn{pseudoline arrangements with contact points} covering a given support. We exploit finally these results to discuss some open problems on multitriangulations, in particular the question of the polytopal realization of their flip graphs.

We study secondly the polytopality of Cartesian products. We investigate the existence of polytopal realizations of cartesian products of graphs, and we study the minimal dimension that can have a polytope whose \skeleton{k} is that of a product of simplices.

\svs
\noindent\textbf{Keywords:} triangulation, multitriangulation, pseudotriangulation, graph of flip, associahedron, Catalan number, pseudoline arrangement, rigidity, polytope, Cartesian product.

\end{document}